\documentclass[10pt]{amsbook}
\usepackage{amsmath,amssymb, amsfonts}
\usepackage[all]{xy}
\usepackage{hyperref}
\usepackage{color}
\usepackage{verbatim}

\usepackage{enumerate}
\usepackage{tikz-cd}
\usepackage{MnSymbol}
\newcommand{\heart}{\ensuremath\heartsuit}

\DeclareMathAlphabet{\mathpzc}{OT1}{pzc}{m}{it}
\newcommand{\adj}[4]{#1\negmedspace: #2\rightleftarrows #3:\negmedspace #4}
\newtheorem{thm}{Theorem}[chapter]
\newtheorem{cor}[thm]{Corollary}
\newtheorem{example}[thm]{Example}
\newtheorem{lem}[thm]{Lemma}
\newtheorem{obs}[thm]{Observation}
\newtheorem{ass}[thm]{Assumptions}
\newtheorem{prop}[thm]{Proposition}
\theoremstyle{definition}
\newtheorem{defn}[thm]{Definition}
\newtheorem{notation}[thm]{Notation}

\newtheorem{rem}[thm]{Remark}

\DeclareFontFamily{U}{rsf}{} \DeclareFontShape{U}{rsf}{m}{n}{
  <5> <6> rsfs5 <7> <8> <9> rsfs7 <10->  rsfs10}{}
\DeclareMathAlphabet{\mathscr}{U}{rsf}{m}{n}

\newcommand*{\defeq}{\mathrel{\vcenter{\baselineskip0.5ex \lineskiplimit0pt
                     \hbox{\scriptsize.}\hbox{\scriptsize.}}}%
                     =}

\newcommand{\F}{\textrm{F}}

\renewcommand{\imath}{\sqrt{-1}}

\DeclareMathOperator{\QCoh}{QCoh}

\DeclareMathOperator{\C}{C} 
  
\DeclareMathOperator{\h}{h}

\DeclareMathOperator{\Hom}{Hom}

\DeclareMathOperator{\id}{id}
\DeclareMathOperator{\Rlim}{Rlim}
\DeclareMathOperator{\colim}{colim}

\DeclareMathOperator{\hocolim}{hocolim}

\DeclareMathOperator{\spec}{spec}

\DeclareMathOperator{\coker}{coker}
\DeclareMathOperator{\Ban}{Ban}
\DeclareMathOperator{\xNorm}{xNorm}

\DeclareMathOperator{\Norm}{Norm}
\DeclareMathOperator{\Born}{Born}

\DeclareMathOperator{\Map}{Map}

\DeclareMathOperator{\Comm}{Comm} 
 
\DeclareMathOperator{\Ho}{Ho}
\DeclareMathOperator{\Ch}{Ch}

\newcommand{\ootimes}{\overline{\otimes}}
\newcommand{\wotimes}{\widehat{\otimes}}

\def\limpro{\mathop{\lim\limits_{\displaystyle\leftarrow}}}
\def\limind{\mathop{\lim\limits_{\displaystyle\rightarrow}}}

\newcommand{\Mod}{{\text{\bfseries\sf{Mod}}}}

\makeatletter\@addtoreset{thm}{chapter}\makeatother

\usepackage{amsrefs}

\usepackage{minitoc}
\begin{document}

\title[]{A Perspective on the Foundations of Derived Analytic Geometry}
\author{Oren Ben-Bassat, Jack Kelly, Kobi Kremnizer}\thanks{}%
\address{Oren Ben-Bassat, ben-bassat@math.haifa.ac.il, Department of Mathematics, University of Haifa, Haifa 3498838, Israel}
\address{Jack Kelly, Jack.Kelly@lincoln.ox.ac.uk,  Mathematical Institute,
University of Oxford,
Andrew Wiles Building,
Radcliffe Observatory Quarter,
Woodstock Road,
Oxford,
OX2 6GG, England
}
\address{Kobi Kremnizer,Yakov.Kremnitzer@maths.ox.ac.uk, Mathematical Institute,
University of Oxford,
Andrew Wiles Building,
Radcliffe Observatory Quarter,
Woodstock Road,
Oxford,
OX2 6GG, England}

\dedicatory{}
\thanks{ The first author acknowledges the support of the European Commission under the Marie Curie Programme for the IEF grant which enabled this research to take place. The contents of this article reflect the views of the authors and not the views of the European Commission. The second author acknowledges the support of the Simons Foundation through its `Targeted Grants to Institutes' program.}

\begin{abstract}
We show how one can do algebraic geometry with respect to the category of simplicial objects in an exact category. As a biproduct, we get a theory of derived analytic geometry.
\noindent 
\end{abstract}
\maketitle 
\tableofcontents

\chapter{Introduction, Notation, and Conventions}

\section{Introduction}

The perspective on derived analytic geometry which we present in this work is based on the following maxim.
\begin{quote}
$$\textrm{geometry }=\textrm{ algebra }+\textrm{ topology }$$
\end{quote}
Classically this is evidenced through \textit{ringed spaces}. A ringed space is a pair $(X,\mathcal{O}_{X})$ where $X$ is a topological space (topology) and $\mathcal{O}_{X}$ is a sheaf of rings on $X$ (algebra). Consider for example the case $X=\mathbb{C}$ with its Euclidean topology. The analytic geometry on $\mathbb{C}$ is realised through the sheaf $\mathcal{O}^{an}_{\mathbb{C}}$ of analytic functions on $\mathbb{C}$, and the smooth geometry through the sheaf $\mathcal{O}^{sm}_{\mathbb{C}}$ of smooth functions on $\mathbb{C}\cong\mathbb{R}^{2}$. 

Let $\mathrm{Spaces}$ denote the category of ringed spaces. The assigment $(X,\mathcal{O}_{X})\rightarrow\mathcal{O}_{X}(X)$ determines a contravariant functor 
$$\mathcal{O}(-):\mathrm{Spaces}\rightarrow\mathrm{Rings}$$
 Grothendieck observed that there is a full subcategory $\mathrm{Aff}$ of \textit{affine schemes} such that this functor restricts to an equivalence. In other words this functor has a fully faithful left adjoint.
 $$\mathrm{Spec}:\mathrm{Rings}^{op}\rightarrow\mathrm{Spaces}$$
 One then defines a `space' (precisely a prestack) to be a contravariant functor $\mathcal{X}:\mathrm{Aff}^{op}\rightarrow\mathrm{Grpds}$
  from the category of affine schemes to the category of groupoids. Schemes live inside the category of prestacks via their functor of points, $X\mapsto\mathrm{Hom}_{\mathrm{Spaces}}(\mathrm{Spec}(A),X)$. Using various Grothendieck topologies on $\mathrm{Aff}$ (Zariski, \'{e}tale, Nisnevich, smooth, fpqc, fppf$\ldots$), one can identify, in a purely formal and category-theoretic way, subcategories of `geometric stacks' (e.g. schemes, algebraic spaces, Artin stacks, Deligne-Mumford stacks) which are glued together from affines using the topology. This gives another interpretation of our maxim - geometry is the combination of the category of unital commutative rings (algebra) and a (Grothendieck) topology. In this picture, the algebra describes local pieces (coordinate charts), and the Grothendieck topology tells us how these pieces can be glued together to obtain global objects.
 
 This categorical approach provides a stunningly powerful dictionary between algebra and geometry. Moreover, this formalism works over any base ring, not just fields. This leads to a `universal' geometry over the integers, unifying number theory with complex geometry through, for example, the Weil conjectures. This `functor of points' framework can be formulated in any presentably symmetric monoidal category $\mathrm{C}$ by the following procedure
\begin{enumerate}
\item
 \textit{Define} the category of affine schemes over $\mathrm{C}$ to be the category $\mathrm{Aff_{C}}\defeq\mathrm{Comm}(\mathrm{C})^{op}$ opposite to the category of commutative monoids internal to $\mathrm{C}$. 
 \item
 \textit{Define} the category $\mathrm{PreStk_{C}}$ of prestacks over $\mathrm{C}$ to be the category of contravariant functors
 $$\mathrm{Aff_{C}}^{op}\rightarrow\mathrm{Grpd}$$
 \item
 Identify a Grothendieck topology on $\mathrm{Aff_{C}}$, and use this to identify geometric stacks.
 \end{enumerate}
 Moreover there is a ready-baked definition of the category of quasi-coherent sheaves on a pre-stack. Using the natural equivalence
 $$\mathrm{QCoh}(\mathrm{Spec}(A))\cong{}_{A}\mathrm{Mod}$$
 from classical algebraic geometry, for an object $A\in\mathrm{Comm(C)}$ one simply defines
   $$\mathrm{QCoh}(\mathrm{Spec}(A))\defeq{}_{A}\mathrm{Mod}(\mathrm{C})$$
   to be the category of $A$-modules internal to $\mathrm{C}$ (here $\mathrm{Spec}(A)$ is just notation, it refers to $\mathrm{A}$ regarded as an object of the opposite category $\mathrm{Comm(C)}^{op}$). This defines a $(2-)$functor
   $$\mathrm{QCoh}:\mathrm{Aff_{C}}\rightarrow\mathrm{Cat}^{R}$$
   to the $(2-)$category of presentable categories with morphisms being right adjoint functors. By Kan extension this extends to a functor
      $$\mathrm{QCoh}:\mathrm{PreStk_{C}}\rightarrow\mathrm{Cat}^{R}$$
   In particular any morphisms of pre-stacks $f:\mathcal{X}\rightarrow\mathcal{Y}$ automatically determines an adjunction
   $$\adj{f^{*}}{\mathrm{QCoh}(\mathcal{Y})}{\mathrm{QCoh}(\mathcal{X})}{f_{*}}$$
 This relative approach to geometry developed for instance in \cite{deligne2007categories} and \cite{banerjee2017noetherian}, also lends itself well to theories of derived geometry through the work of To\"{e}n and Vezzosi \cite{toen2008homotopical}. Essentially, one replaces $\mathrm{C}$ with a combinatorial symmetric monoidal model category satisfying some technical assumptions, the category of prestacks by the category of simplicial prestacks, and the Grothendieck topology by a homotopy Grothendieck topology. For derived algebraic geometry over a ring $R$, $\mathrm{C}$ will be the category of simplicial $R$-modules, and there are derived enhancements of the usual topologies.


The goal of this monograph is to simultaneously give a relative algebraic formulation of analytic and smooth geometry, and to develop the foundations of derived versions of these geometries. The idea is as follows. Let $X$ be a complex analytic space. The algebra of functions $\mathcal{O}(X)$ has a canonical complete locally convex topology compatible with the algebraic structure. Denoting by $\hat{\mathcal{T}}_{c}$ the category of complete locally convex spaces over $\mathbb{C}$, the assigment
$$\mathcal{O}(-):\mathrm{Stein Spaces}\rightarrow\mathrm{Comm}(\hat{\mathcal{T}_{c}})$$
is fully faithful, where $\hat{\mathcal{T}_{c}}$ is equipped with the completed projective tensor product $\hat{\otimes}_{\pi}$. Moreover for $X$ and $Y$ Stein spaces we have $\mathcal{O}(X\times Y)\cong\mathcal{O}(X)\hat{\otimes}_{\pi}\mathcal{O}(Y)$, where $\hat{\otimes}_{\pi}$ is the complete projective tensor product of locally convex spaces. As we explain, this story can also be repeated for rigid geometry over a non-Archimedean Banach field $k$ by replacing Stein spaces with affinoids, and also for smooth geometry over $\mathbb{R}$. Since Stein spaces are in essence the affine objects of complex geometry, as a first approximation, one might then consider algebraic geometry as geometry relative to the category $\hat{\mathcal{T}}_{c}$ of complete locally convex topological vector spaces over $\mathbb{C}$. There are numerous problems with the approach. The first is that locally convex vector spaces are not an abelian category. This is easily surmountable: they have a Quillen exact structure, and if completeness is not imposed the category is even quasi-abelian. More fatally though they do not form a closed symmetric monoidal category with the projective tensor product, they do not form a presentable category, and there are not enough projectives. 

For the purposes of geometry though we are not interested in all locally convex spaces. Those which appear as algebras of functions are typically Fr\'{e}chet spaces or Banach spaces. There are several approaches to extending the category of  Fr\'{e}chet spaces to a well-behaved closed symmetric monoidal abelian category. Indeed condensed mathematics  due to Clausen and Scholze \cite{scholze2019lectures}, \cite{scholze20192lectures}, \cite{scholzecomplexlectures} - specifically the liquid theory in the Archimedean setting over $\mathbb{R}$ or $\mathbb{C}$, and the solid theory in the non-Archimedean setting - is one such extremely successful approach. Condensed mathematics is mathematics relative to a certain compact projectively generated topos, namely the topos of condensed sets. This is the free sifted cocompletion of a certain subcategory of compact Hausdorff spaces. 

Our method is to instead consider bornological structures on vector spaces, which regards boundedness rather than openness as primitive. Given a locally convex topological vector space $X$ over $\mathbb{C}$, we say that a subset $B$ is bounded if all continuous semi-norms on $X$ are bounded functions when restricted to $B$. A map $f:X\rightarrow Y$ of locally convex vector spaces is said to be \textit{bounded} if it sends bounded sets to bounded sets. Modulo some details one can then define a bornological vector space over $\mathbb{C}$ to be a $\mathbb{C}$-vector $V$ space equipped with a collection of \textit{bounded} subsets of $V$ such that addition and scalar multiplication are bounded maps - i.e. they send bounded subsets to bounded subsets. Bornological vector spaces arrange into a category in which the morphisms are bounded maps. We denote this category by $\mathrm{Born}_{\mathbb{C}}$. The construction above actually yields a functor
$$\hat{\mathcal{T}}_{c}\rightarrow\mathrm{Born}_{\mathbb{C}}$$
This functor is fully faithful when restricted to Fr\'{e}chet spaces. Moreover, the image of this functor in fact lands in the extremely well-behaved category $\mathrm{CBorn}_{\mathbb{C}}$ of complete bornological $\mathbb{C}$-vector spaces of convex type. These are bornological spaces which can be described as increasing unions of Banach spaces. Formally there is an equivalence $\mathrm{CBorn}_{\mathbb{C}}\cong\mathrm{Ind^{m}(Ban}_{\mathbb{C}}\mathrm{)}$ with the free cocompletion by monomorphic filtered coimits of $\mathrm{Ban}_{\mathbb{C}}$. For technical reasons it is convenient also to work with the category $\mathrm{Ind(Ban}_{\mathbb{C}}\mathrm{)}$, which is the free cocompletion of $\mathrm{Ban}_{\mathbb{C}}$ by all filtered colimits. Both $\mathrm{CBorn}_{\mathbb{C}}$ and $\mathrm{Ind(Ban}_{\mathbb{C}}\mathrm{)}$ are \textit{quasi-abelian} categories. Consequently their derived categories have both a left and right $t$-structure. The heart of the left $t$-structure is denoted $\mathrm{LH}$, and we have equivalences
$$\mathrm{LH}(\mathrm{CBorn}_{\mathbb{C}})\cong\mathrm{LH}(\mathrm{Ind(Ban}_{\mathbb{C}}\mathrm{)})$$
The category $\mathrm{LH}(\mathrm{Ind(Ban}_{\mathbb{C}}\mathrm{)})$ is a monoidal elementary abelian category in the sense of \cite{kelly2016homotopy}, and therefore is a particularly convenient setting for homotopical algebra and derived geometry. Moreover the composite functor
$$\hat{\mathcal{T}}_{c}\rightarrow\mathrm{CBorn}_{\mathbb{C}}\rightarrow\mathrm{LH}(\mathrm{Ind(Ban}_{\mathbb{C}}\mathrm{)})$$
is not only fully faithful when restricted to Fr\'{e}chet spaces, but even strongly monoidal when restricted to nuclear Fr\'{e}chet spaces.

In summary we get a contravariant functor
$$\mathcal{O}(-):\mathrm{Analytic\;Spaces}\rightarrow\mathrm{Complete\;Bornological\;Rings}$$
which \textit{is fully faithful} when restricted to Stein spaces. Moreover any $X$ has an atlas consisting of Stein spaces. The hope then is that we may realise analytic geometry as geometry relative to the symmetric monoidal category $\mathrm{CBorn}_{\mathbb{C}}$ of complete bornological $\mathbb{C}$-vector spaces, with analytic spaces being `schemes' relative to some topology on $\mathrm{Aff}_{\mathrm{CBorn}_{\mathbb{C}}}$. In this work we will realise this hope.  Finally, the category $\mathrm{CBorn}_{\mathbb{C}}$ is an exact (in fact quasi-abelian) category, and this homological structure gives rise to a (combinatorial monoidal) model structure on the category $\mathrm{sCBorn}_{\mathbb{C}}$ of simplicial complete bornological $\mathbb{C}$-modules. Moreover there is a Grothendieck topology on $\mathrm{Aff}_{\mathrm{sCBorn}_{\mathbb{C}}}$ with respect to which complex analytic spaces are schemes. This will be the basis of our theory of derived analytic geometry. There are several compelling reasons for desiring this model.
\begin{enumerate}
\item
It gives us a large toolset to define analytic moduli spaces using such as mapping stacks, quotient stacks.
\item
It gives us a definition of the (derived) category of quasi-coherent sheaves on an analytic space \textit{which it turns out} satisfies (derived) descent.
\item
Algebraic geometry translates results from (homological) commutative algebra into results about (derived) local/ affine geometry. For example Koszul duality between Lie algebras and (co)commutative (co)algebras tells us about deformation theory. The Hochschild-Kostant-Rosenberg Theorem relates the shifted tangent complex of a stack to the loops space of a stack. As proved in \cite{kelly2019koszul} and \cite{kelly2021analytic}, this is also possible for analytic geometry.
\item
One can define the category $\mathrm{CBorn}_{R}$ of complete bornological $R$-modules over any Banach ring $R$. $\mathrm{sCBorn}_{R}$ is still a combinatorial monoidal model category, and our Grothendieck topology can be define in this generality. Thus one obtains a theory of analytic geometry relative to any Banach ring. This includes the Banach ring $\mathbb{Z}_{an}$ of integers equipped with the absolute value norm, giving rise to a theory of `universal analytic geometry' which specialises to both complex analytic geometry and rigid analytic geometry.
\item
A derived enhancement of analytic geometry will furnish an arsenal of new tools from intersection, deformation theory, and obstruction theory.
\end{enumerate}

Although point (4) has also proved extremely powerful in algebraic geometry, one can of course very successfully work in the algebraic context without ever worrying about the derived world. In contrast, our relative model of analytic geometry  \textit{must be derived from the outset}. The reason is quite simple: in the topologies of interest in algebraic geometry all maps are flat. This is emphatically not true in the analytic setting. As we will show in the complex analytic case (and has as already been shown in the rigid case in \cite{ben2020fr}), even the restriction map $\mathcal{O}(\mathbb{C})\rightarrow\mathcal{O}(B(0,r))$ from analytic functions on $\mathbb{C}$ to analytic functions on an open ball is not flat. In fact it is precisely this fact which is fundamental in Gabber's example explaining why the na\"{i}ve definition of a quasi-coherent sheaf on a complex analytic space fails to satisfy descent, and why the passage to the derived world solves this.

A comparison between bornological and condensed mathematics will not be the subject of this paper. However let us mention that much of the abstractly formulated content of the present work applies immediately to the solid theory of condensed mathematics, at least in characteristic zero. For the liquid theory and general condensed abelian groups one has to be somewhat careful as the tensor product of two projectives is not necessarily projective - however we expect with some minor modifications it should apply there as well. We also expect that the categories of analytic geometric stacks in each geometric context, suitably defined and with appropriate restrictions, should be equivalent, at least when their structure sheaves have coherent higher cohomology. What will differ between the two approaches will be the categories of quasi-coherent sheaves defined on such spaces. There are two noteworthy technical differences between the theories.
The first is the aforementioned issue with the tensor product of projectives. The second arises when in the condensed setting one wishes to introduce various concepts of completion, and to treat Archimedean and non-Archimedean geometry simultaneously. Here one must introduce a slight variation on algebra internal to a symmetric monoidal category in order to incorporate both the solid and liquid theories at once. The affine objects are no longer just commutative algebra objects, but the data of a condensed ring \textit{and} a symmetric monoidal reflective subcategory of its category of modules. This is only a slight additional complication, and it is highly unlikely that it will make proofs significantly more onerous. In fact in the bornological world in order to incorporate Archimedean and non-Archimedean geometry at once we must also do something similar, however for us there are only two necessary reflective subcategories: the entire category of complete bornological abelian groups, and the reflective sub-category of non-Archimedean complete bornological abelian groups. This is because in the complete bornological world from the outset everything is complete. As a result, particularly in the context of Archimedean geometry, the category $\mathrm{CBorn}_{\mathbb{R}}$ seems considerably easier to define  than the category of $p$-liquid $\mathbb{R}$-vector spaces.

One point of similarity to note which is of fundamental importance for $K$-theory, is that we believe it should not be too difficult to show that nuclear objects in the bornological context over $\mathbb{C}$ or $\mathbb{R}$, and the liquid context, coincide. Indeed it is quite obvious that, as in the liquid theory, the category of nuclear bornological spaces is generated under $\aleph_{1}$-filtered colimits by bounded complexes of dual nuclear Fr\'{e}chet (DNF) spaces which are precisely the $\aleph_{1}$-compact nuclear bornological spaces. All that remains to check is that derived mapping spaces between DNF spaces computed in $\mathrm{CBorn}_{\mathbb{R}}$ and $\mathrm{Liq}_{p}$ coincide. 

Let us highlight a more interesting and fundamental point of distinction between the contexts. We wish to stress that from our perspective this does not represent a drawback for either side, it simply gives rise to different interesting mathematics. As in condensed mathematics there is a topos/ $(\infty,1)$-topos inside of which bornological algebra takes place. Roughly speaking it is the category of bornological normed sets. It has a formal description as the $(\infty,1)$-categorical free sifted cocompletion of the category of semi-normed sets (a set equipped with a function to the non-negative real numbers). One can stabilise this category to obtained the category of bornological spectra. Objects in these categories represent `bornological cohomology theories' which include, for example, an enhancement of bounded cohomology. These cohomology theories can in principle be evaluated on an analytic space, to obtain new invariants. A similar procedure works for condensed spectra, and here one has a natural representative for continuous cohomology. It is not at all clear to us how one would naturally represent bounded cohomology as a condensed spectrum, or continuous cohomology as a bornological spectrum (although perhaps it is clear to someone!). Thus the bornological and condensed setups provide different, interesting cohomology theories. 

There is also an approach to derived analytic geometry in the spirit of ringed spaces developed by Porta and Yue Yu (\cite{porta2016higher},\cite{porta2015derived1} \cite{porta2015derived}, \cite{porta2017derived},\cite{porta2018derived}, \cite{porta2017representability}, \cite{porta2018derivedhom}). Spivak has introduced a similar approach to derived smooth geometry in \cite{spivak2010derived}. These are based on Lurie's approach to derived geometry using a higher-categorical generalisation of ringed spaces, namely structured $(\infty,1)$-topoi. This is an $(\infty,1)$-topos $\mathcal{X}$ together with a limit-preserving functor 
$$\mathcal{O}:\mathpzc{G}\rightarrow\mathcal{X}$$
where $\mathpzc{G}$ is a {geometry} - an $(\infty,1)$-category with the structure of a Grothendieck topology and a class of `admissible' morphisms. The functor $\mathcal{O}$ is required to respect these structures in a precise way. For example taking $(\mathpzc{G})^{op}$ to be the $(\infty,1)$-category of smooth schemes gives a reasonable notion of derived algebraic stacks. For analytic geometry Porta and Yue Yu take $\mathcal{G}$ to be the category $\mathcal{T}_{an}$ of smooth analytic spaces (over $\mathbb{C}$ or over a non-Archimedean field $k$), and for smooth geometry Spivak takes as $(\mathpzc{G})^{op}$ the category of smooth manifolds. Thanks to Porta and Yue Yu the field of derived analytic geometry has developed rapidly, leading to GAGA, base-change, and Riemann-Hilbert type theorems. They have also announced a Hochschild-Kostant-Rosenberg theorem. Despite these crucial results there are some technical difficulties. For example, there is no obvious definition of the category of quasi-coherent sheaves on a derived analytic space. Restricting to coherent sheaves makes pullback and push-forward functors tricky to construct, and base-change theorems for coherent sheaves difficult to formulate and prove. Finally the volume of data bound up in the definition of a $\mathcal{T}_{an}$-structured $(\infty,1)$-topos makes it somewhat difficult to check examples. Pridham (\cite{MR4036665}) has shown that at least for Stein spaces, the full pregeometry $\mathcal{T}_{an}$ may be replaced by a discrete pre-geometry arising from a Lawvere theory. This has some advantages, including the existence of concrete model category presentations. However this still does not resolve the difficulty in defining quasi-coherent sheaves, push-forward functors, etc. It also does not encompass Porta and Yue Yu's derived rigid geometry in full. The relative geometry approach will resolve many of the issues of the models of Porta and Yue Yu, and Pridham. Moreover in the complex analytic setting we will show that their categories of ($m$-skeletal) geometric stacks embed fully faithfully in ours. In future work we will also establish this for their non-Archimedean geometric stacks.

More generally in this work we will develop a theory of derived geometry relative to any sufficiently well behaved monoidal abelian/ exact category $\mathpzc{E}$. Specialising to the category of complete bornological spcases will give a theory of derived analytic geometry, and sketch a theory of derived smooth geometry. One difficulty is defining what one means by finitely presented/ smooth/ \'{e}tale in the various geometries we consider. For example standard localisations such as inclusions of polydiscs \textit{do not} correspond to finitely presented maps of algebras internal to $\mathrm{CBorn}$. This is resolved through the use of Lawvere theories which determine the geometry. In summary we will in fact develop derived geometry relative to three pieces of data: a monoidal exact category $\mathpzc{E}$, a possibly multi-sorted Lawvere theory $\mathrm{T}$, and a fully faithful, finite coproduct-preserving functor $\mathrm{F}:\mathrm{T}^{op}\rightarrow\mathrm{Comm}(\mathpzc{E})$ to the category of commutative monoids in $\mathrm{Comm}(\mathpzc{E})$. Including this additional data allows us to define generalised notions of algebras and maps of finite presentation, permitting us to treat geometry relative to various different kinds of algebras - such as Tate algebras, algebras of entire functions, algebras of smooth functions ... - in a unified and rigorous manner.

\section{Outline of the Paper}
Let us give an overview of the structure of the paper and its contents.

\subsection{Chapter 2: Homotopical Algebra}
 In Chapter \ref{HA} we develop the necessary theory of affine/ local geometry relative to a symmetric monoidal $(\infty,1)$-category, i.e. commutative algebra internal to such a category. 
 
 Inspired by work of Raksit \cite{raksit2020hochschild}, in Section \ref{sec:infalgcontext} we develop a unified $(\infty,1)$-categorical setup for simultaneously dealing with geometry relative to simplicial commutative algebras, and geometry relative to $E_{\infty}$-algebras. Essentially, such a category consists of a a presentably monoidal stable $(\infty,1)$-category $\mathbf{C}$ with $t$-structure together with a graded monad $\mathbf{D}$, and a map of monads 
 $$\theta:\mathbf{Comm}(-)\rightarrow\mathbf{D}$$
 satisfying some axioms. Most of the time $\mathbf{C}$ will have the structure of an algebraic context in the sense of \cite{raksit2020hochschild}, and $\mathbf{LSym}$ will be the derived symmetric monad which roughly speaking encodes simplicial commutative rings and their non-connective versions. The only other example which will appear in practice is $\mathbf{D}=\mathbf{Comm}(-)$ with $\theta$ being the identity. This will be useful for Dirac geometry \cite{HessPstra} and its generalisations to other contexts. The main results of this section will concern the cotangent complex and its formal properties. Crucially, many of the important results of \cite{HA} Chapter 7 carry through to this slightly greater generality. Under an additional constraint on connectivity of the monad $\mathbf{D}$ - that the $t$-structure is \textit{Koszul} for $\mathbf{D}$ = the core technical result from which all useful properties follow is 
 
 \begin{thm}\ref{thm:connectivcot}
 Let the $t$-structure be Koszul. Let $f:A\rightarrow B$ be a morphism in $\mathbf{Alg_{D}}(\mathbf{C}_{\ge0})$ and suppose that $\mathrm{cofib}(f)\in\mathbf{C}_{\ge n}$ where $n\ge0$. Then $\mathrm{fib}(\epsilon_{f})\in\mathbf{C}_{\ge n}$ when $n=0$
\end{thm}

In particular it implies the following central results in obstruction theory regarding the Postnikov tower
$$A\rightarrow\cdots \rightarrow \tau_{\le n}A\rightarrow \cdots \rightarrow \tau_{\le n-1}A\rightarrow\cdots\rightarrow\tau_{\le0 }A$$
of a $\mathbf{D}$-algebra.

\begin{thm}\ref{thm:squarezeropostnikov}
Let $A\in\mathbf{Alg}_{\mathbf{D}}^{cn}(\mathbf{C})$. Then for each $n\le -1$ there is a derivation
$$d^{A}_{n}\in\pi_{0}\mathbf{Der}(\tau_{\le n}A,\Sigma^{n+2}\pi_{n+1}(A))$$
such that $\tau_{\le n+1}A$ is the square-zero extension of $\tau_{\le n}A$ along $d^{A}_{n}$. 
\end{thm}
 
In Sections \ref{sec:spapc} and \ref{sec:spac} we introduce generalisations of Raksit's dervived algebraic contexts \cite{raksit2020hochschild}, which encompasses monoidal $(\infty,1)$-categories of modules of differential-graded algebras, and will also in future apply to spectral versions of analytic geometry. These are stable, presentable, symmetric monoidal $(\infty,1)$-categories with $t$-structure, with some compatibilities between the monoial structure and the $t$-structure. The data is written as follows
$$(\mathbf{C},\mathbf{C}_{\ge0},\mathbf{C}_{\le0},\mathbf{C}^{0})$$
with $\mathbf{C}^{0}$ being a specificied set of projective generators for $\mathbf{C}_{\ge0}$. The heart $\mathbf{C}^{\heart}$ is itself a symmetric monoidal (elementary) abelian category. One of the most important technical properties is the following. 

\begin{lem}[Lemma \ref{lem:flattorspec}]
Let $(\mathbf{C},\mathbf{C}_{\ge0},\mathbf{C}_{\le0},\mathbf{C}^{0})$ be a flat spectral algebraic pre-context. Then there is a spectral sequence
$$\mathrm{Tor}_{p}^{\pi_{*}(A)}(\pi_{*}(B),\pi_{*}(C))_{q}\Rightarrow\pi_{p+q}(B\otimes^{\mathbb{L}}_{A}C)$$
\end{lem}

After some discussion of transversality and projectivity, we introduce various classes of maps, including formally $P$-smooth and formally \'{e}tale maps. We also define discrete versions of these classes of maps in $\mathbf{C}^{\heart}$, and obtain the following comparison which is crucial for obstruction theory. 

\begin{prop}[Proposition \ref{prop:discretestrongsmooth}]
\begin{enumerate}
\item
Let $f:A\rightarrow B$ be a derived strong map such that $\pi_{0}(A)\rightarrow\pi_{0}(B)$ is formally $P$-smooth. Then $f$ is formally $P$-smooth. 
\item
If $f$ is $P$-smooth then $f$ then $\pi_{0}(A)\rightarrow\pi_{0}(B)$ is discrete $P$-smooth. 
\item
If $f$ and $\pi_{0}(f)$ are formally smooth then $f$ is derived strong. 
\end{enumerate}
\end{prop}

\subsection{Chapter 3: Model Category Presentations}

Chapter \ref{MCP} provides concrete models for those $(\infty,1)$-algebra contexts which form the basis of this paper, i.e. those arising from so-called elementary exact categories. This will permit explicit computations. In Section \ref{sec:haec} we give a large class of examples of derived algebraic contexts, arising from certain monoidal (one-sided) exact categories.We also explain the link to  To\"{e}n and Vezzosi's formalism \cite{toen2008homotopical}, by showing that all of our examples of interest are in fact presented by HA contexts. Most of the results in this section are of a technical nature. The main result is the following.

\begin{thm}[Theorem \ref{thm:HAcontextsimplicial}]
Let $(\mathcal{P},\mathcal{E}_{\mathcal{P}})$ be a projective class on a locally presentably closed symmetric monoidal additive category $\mathcal{C}$. Suppose that $(\mathcal{C},\mathcal{P},\mathcal{E}_{\mathcal{P}})$ is projectively monoidal, has symmetric projectives, and is $\mathbf{SplitMon}$-elementary. Then $\mathrm{s}\mathcal{C}$ is a HA context.
\end{thm}

We also define nuclear and conuclear objects in exact categories, and prove some general results under which tensor products commute with limits in $\infty$-categories of chain complexes in exact categories.


\subsection{Chapter 4: Lawvere Theories}

Chapter \ref{LT} is on Lawvere Theories. In Section \ref{sec:multi-sort} we begin by recalling some elementary properties of $\Lambda$-sorted Lawvere theories, denoted $\mathrm{T}$, and their algebras. We also introduce multi-sorted Fermat theories, and so-called filtered versions of these as well. After detailing some basic homotopy theory of such algebras, we then in Section \ref{sec:LAWVFunc} move on to understanding how to realise algebras over Lawvere theories as simplicial commutative algebras internal to some symmetric monoidal odel category. Specifically we fix a closed symmetric monoidal category $\mathrm{C}$ with a collection of compact objects $\mathcal{Z}$ defining a model structure on $\mathrm{sC}$ which satisfies some conditions such that the transferred model structure exists on $\mathrm{sComm(C)}$, a filtered Lawvere theory $\mathrm{T}$, and a finite coproduct preserving functor
$$\mathrm{F}:\mathrm{T}^{op}\rightarrow\mathrm{sComm(C)}$$
The right hand side presents an $(\infty,1)$-category $\mathrm{L^{H}}(\mathrm{sComm(C)})$, and the functor $\mathrm{F}$ extends to by sifted colimits to a functor
$$\mathbf{P}_{\Sigma}(\mathrm{T}^{op})\cong\mathbf{sAlg}_{\mathrm{T}}\cong\mathrm{L^{H}}(\mathrm{sAlg_{T}})\rightarrow\mathrm{L^{H}}(\mathrm{sComm(C)})$$
We are particularly interested in functors $\mathrm{F}$ which realise $\mathrm{T}^{op}$ as a Lawvere theory of \textit{homotopy }$\mathrm{C}$-\textit{polynomial type} (see Definition \ref{defn:indhtpypoly}) which roughly means that for any $\underline{\lambda}\in\Lambda^{n}$ there is a homotopy epimorphism
$$\mathrm{Sym}(\mathbb{I}^{n})\rightarrow\mathrm{F}(\mathrm{T}(\underline{\lambda}))$$
i.e., that $\mathrm{F}(\mathrm{T}(\underline{\lambda}))$ are in some sense completions of polynomial algebras. In the case that $\Lambda=\{*\}$ and $\mathrm{F}$ is concretely of homotopy $\mathrm{C}$-polynomial type, i.e. that the map
$$\mathbf{Map}_{\mathbf{sAlg}_{\mathrm{T}}}(\mathrm{T}(1),A)\rightarrow\mathbf{Map}(\mathrm{Sym}(\mathbb{I}),\mathbf{F}(A))$$
is an equivalence, we prove that $\mathbf{F}$ is fully faithful (Lemma \ref{lem:Temb}). We also identify the essential image of this functor. Let $\mathbf{sAlg}_{\mathrm{T}}^{\pi_{0}-fp}$ denote the full subcategory of $\mathbf{sAlg}_{\mathrm{T}}$ consisting of those algebras $A$ such that $\pi_{0}(A)$ is finitely presented as a discrete $\mathrm{T}$-algebra.

\begin{thm}[Theorem \ref{thm:essimageF}]
Suppose that the functors
$$\mathrm{F}:\mathrm{Alg}_{\mathrm{T}}\rightarrow\mathrm{Comm}(\mathbf{C}^{\heart})$$ 
are fully faithful and that for any finitely presented $C\in\mathrm{Alg_{T}}$ the functor
$${}_{C}\mathrm{Mod}^{alg}\rightarrow{}_{\mathrm{F(C)}}\mathrm{Mod}(\mathbf{C}^{\heart})$$
is exact.
Then the essential image of $\mathbf{F}|_{\mathbf{sAlg}_{\mathrm{T}}^{\pi_{0}-fp}}$ consists of all objects $A$ of $\mathbf{DAlg}^{cn}(\mathbf{C})$ such that 
\begin{enumerate}
\item
$\pi_{0}(A)$ is in the essential image of $\mathrm{F}$.
\item
$\pi_{n}(A)\in{}_{\pi_{0}(A)}\mathrm{Mod}^{alg}(\mathbf{C}^{\heart})\subset {}_{\pi_{0}(A)}\mathrm{Mod}(\mathbf{C}^{\heart})$.
\end{enumerate}
\end{thm}

Finally in Section \ref{sec:Tderloc} we discuss derived qutients of $\mathrm{T}$-algebras and define $\mathrm{T}$-algebras of finite presentations internal to $\mathbf{DAlg}^{cn}(\mathbf{C})$ and show that they have perfect relative cotangent complexes. We conclude by defining $\mathrm{T}$-standard \'{e}tale and smooth maps:

\begin{defn}[Definition \ref{defn:Tstandardsmooth}]
 A map $f:A\rightarrow B$ in $\mathbf{DAlg}^{cn}(\mathbf{C}^{\heart})$ is said to be 
 \begin{enumerate}
 \item
 $\mathrm{T}$-\textit{standard smooth of relative dimension }$k$ if there is an isomorphism
$$B\cong A\otimes^{\mathbb{L}}\mathrm{T}(\lambda_{1},\ldots,\lambda_{c+k})\big\slash\big\slash(f_{0},\ldots,f_{c})$$
where $f_{0},\ldots,f_{c}:R\rightarrow\pi_{0}(A)\otimes\mathrm{T}(\lambda_{1},\ldots,\lambda_{n})$ are such that 
$$\mathrm{det}(\frac{\partial f_{i}}{\partial x_{j}})_{1\le i,j\le c}$$
is a unit in $\pi_{0}(B)$
\item
$\mathrm{T}$-\textit{standard \'{e}tale} if it is $\mathrm{T}$-standard smooth of relative dimension $0$.
\end{enumerate}
\end{defn}

\subsection{Chapter 5: Higher Bornological Algebra}

In Chapter 5 we introduce the main example of interest in this paper - the derived algebraic context of and complete bornological modules over a Banach ring $R$. We begin in Section \ref{sec:bornsp} by recalling some basic properties and constructions in the category $\mathrm{Ban}_{R}$ of Banach modules over $R$ including the existence of the projective closed monoidal structure, and the procedure for constructing a projective generating set. The upshot is that the category $\mathrm{Ind(Ban_{R})}$ is a monoidal elementary exact category with symmetric projectives, and thus presents a derived algebraic context.

\begin{cor}[Corollary \ref{cor:IndBangood}]
$\mathrm{Ind}(\Ban_{R})$ is a monoidal elemenetary quasi-abelian category when equipped with the projective tensor product, and has symmetric projectives.
\end{cor}

Moreover we prove that this context is Koszul, so that it has a rich obstruction theory.

We also recall from \cite{bambozzi2015stein} how Fr\'{e}chet spaces over a field embed fully faithfully in $\mathrm{Ind(Ban_{k})}$, and that when restricted to nuclear Fr\'{e}chet spaces this functor is strong monoidal.

In Section \ref{sec:genbornalg} we develop some basic bornological algebra. In particular we recall facts abut nilpotence in bornological algebras from \cite{bambozzi}. We give a definition of bornological Fr\'{e}chet-Stein algebras which works in both the Archimedean and non-Archimedean contexts. We also define derived Fr\'{e}chet-Stein algebras and, following \cite{ben2020fr}, give examples of metrisable and nuclear objects in $\mathrm{Ind(Ban}_{R}\mathrm{)}$. 

In Section \ref{sec:bornlawvere} we introduce the Lawvere theories of interest. They are all generated by algebras in $\mathrm{Comm(}\mathrm{Ind(Ban}_{R}\mathrm{)}\mathrm{)}$. Moreover, apart from $(3)$, they are all of $\mathrm{CBorn}_{R}$-polynomial type.

\begin{enumerate}
\item
The $1$-sorted formal power series Lawvere theory $\mathrm{FPS}_{R}$ generated by $\widehat{\mathrm{Sym}}(R)=\prod_{m=0}^{\infty}\mathrm{Sym}^{m}(R)$.
\item
The $\mathbb{R}_{>0}$-sorted non-Archimedean Tate Lawvere theory generated by 
$$R\Bigr<\frac{x}{r}\Bigr>\defeq R\Bigr<R_{r}\Bigr>$$
where $R_{r}$ denotes the rescaling of $R$ by $r\in\mathbb{R}_{>0}$, and $ R\Bigr<R_{r}\Bigr>$ is the symmetric algebra computed in the non-expanding Banach category $\mathrm{Ban}^{nA,\le1}_{R}$ where we restrict to those morphisms of norm at most $1$.
\item
The $\mathbb{R}_{>0}$-sorted Archimedean Tate Lawvere theory generated by 
$$R\Bigr\{\frac{x}{r}\Bigr\}\defeq R\Bigr\{R_{r}\Bigr\}$$
where $R_{r}$ denotes the rescaling of $R$ by $r\in\mathbb{R}_{>0}$, and $ R\Bigr\{R_{r}\Bigr\}$ is the symmetric algebra computed in the non-expanding Banach category $\mathrm{Ban}^{\le1}_{R}$ where we restrict to those morphisms of norm at most $1$.
\item
The $\mathbb{R}_{>0}$-sorted non-Archimedean dagger affinoid Lawvere theory generated by
$$R\Bigr<\frac{x}{r}\Bigr>^{\dagger}\defeq \underset{\rho>r}\colim R\bigr<\frac{x}{\rho}\Bigr>$$
\item
The $\mathbb{R}_{>0}$-sorted Archimedean dagger affinoid Lawvere theory generated by
$$R\Bigr\{\frac{x}{r}\Bigr\}^{\dagger}\defeq \underset{\rho>r}\colim R\bigr\{\frac{x}{\rho}\Bigr\}$$
\item
The $1$-sorted Archimedean and non-Archimedean entire functional calculus Lawvere theories, $\mathrm{EFC}_{K}$, generated by
$$\lim_{r>0}R\Bigr\{\frac{x}{r}\Bigr\}^{\dagger}$$
\item 
Over $\mathbb{R}$ the smooth Lawvere theory $\mathrm{CartSp}_{\mathrm{smooth}}$ generated by the algebras of smooth functions $C^{\infty}(\mathbb{R})$.
\end{enumerate}

In particular algebras over the $1$-sorted Lawvere theories, including $\mathbf{EFC}_{K}$ for $K$ a non-trivially valued Banach field, and $\mathrm{CartSp}_{\mathrm{smooth}}$ embed fully faithfully in $\mathbf{DAlg}^{cn}(\mathrm{Ind(Ban_{K})})$. In particular thanks to work of Pridham \cite{MR4036665}, the $\mathbf{EFC}_{\mathbb{C}}$ case gives a fully faithful functor from Porta and Yue Yu's category of derived Stein spaces to the category of affines realtive to $\mathbf{Ch}_{\ge0}(\mathrm{Ind(Ban_{K})}$.

\subsection{Chapter 6: Stacks and Presheaves}

In Chapter 6 we introduce our conventions for geometric stacks, and presheaves of categories on sites. This chapter largely consists of setting up useful formalisms. In particular in Subsection \ref{subsubsec:local}, given a class of maps $\mathbf{P}$ and a topology $\tau$ n a category, we define the class of maps $\mathbf{P}^{\tau}$ in $\mathbf{P}$  which are local for the topology. In Section \ref{sec:stacks} we recall material from \cite{kelly2021analytic} concerning relative $(\infty,1)$-geometry tuples, which provide convenient abstract settings for $(\infty,1)$-categorical geometry. In particular, we recall the definition of geometric stacks in such contexts. In Section \ref{sec:sheaves} we provide some basic notation and definitions for presheaves of categories. We also establish some formal base change results. We conclude in Section \ref{sec:monoidalpresheaves} by considering presheaves of monoidal categories and our main example of interest, the presheaf $\mathbf{QCoh}$ relative to a spectral algebraic context.

\subsection{Chapter 7: Descent}

In Chapter 7 we establish some general descent results for modules. We begin with some somewhat standard generalities, including when descent for presheaves of categories implies descent for sub-presheaves, and relating sheaves on pre-stacks to sheaves on their stackifications. Particularly for Stein geometry and smooth geometry, both non-flatness and non-finiteness can cause several issues, not least of which is the fact that these topologies do not seem to satisfy hyper-descent for quasi-coherent sheaves. For finite covers they will still satisfy \v{C}ech descent. However for infinite covers they only satisfy descent upon restriction to certain subcategories of the categories of modules. We are ultimately forced to consider weaker forms of descent not just for a presheaf of categories $\mathbf{Q}$, but a presheaf together with a subpresheaf $\mathbf{N}$ of `well-behaved' objects for which some form of decsent is satisfied. 

Ultimately we consider the following general `kinds' of descent results for modules in Section \ref{sec:descentmod}. We begin by recalling the class of descendable maps introduced by Mathew in \cite{MR3459022}. In particular we simply observe that Mathew's results immediately imply that such maps do define a topology on $\mathbf{Aff}^{cn}_{\mathbf{C}}$ for which quasi-coherent sheaves satisfy \v{C}ech descent. We also consider faithfully flat descent. Following closely Section D.4 of \cite{HA} we show that faithfully flat maps satisfy \v{C}ech descent:

 \begin{cor}[Corollary \ref{cor:ffdesc}]
Let $f:A\rightarrow B$ be faithfully flat. The functor $B\otimes^{\mathbb{L}}_{A}(-):{}_{A}\mathbf{Mod}(\mathbf{C})\rightarrow{}_{B}\mathbf{Mod}(\mathbf{C})$ is comonadic.
 \end{cor}

The final pre-topology we consider is the $\aleph_{1}$-embedding pre-topology.

\begin{defn}
    Fix a class of maps $\mathbf{P}$ in $\mathbf{Aff_{\mathbf{C}}}$. A collection of maps $\{\mathrm{Spec}(B_{i})\rightarrow\mathrm{Spec}(A)\}_{i\in\mathcal{I}}$ is said to be a cover in the $(\aleph_{1},\mathbf{P})$-\textit{embedding pre-topology} if 
    \begin{enumerate}
    \item 
    for each $i$ the map $A\rightarrow B_{i}$ is in $\mathbf{P}$
        \item 
        for each $i$ the map $B_{i}\otimes_{A}^{\mathbb{L}}B_{i}\rightarrow B_{i}$ is an equivalence.
        \item 
        Let ${}_{A}\mathbf{Mod}^{\aleph_{1},\mathbf{P},RR}$ denote the category of $A$-modules $M$ such that 
        \begin{enumerate}
        \item 
        $M$ is connective
            \item 
            each $\pi_{n}(M)$ is transverse to maps in $\pi_{0}(\mathbf{P})$
            \item 
            each $\pi_{n}(M)\otimes_{\pi_{0}(A)}(-)$ commutes with countable products..
        \end{enumerate}
        There exists a countable subset $\mathcal{J}\subseteq\mathcal{I}$ such that if $\alpha:M\rightarrow N$ is a map of formally $\kappa$-filtered $A$-modules such that $B_{j}\otimes^{\mathbb{L}}_{A}\alpha$ is an equivalence for all $j\in\mathcal{J}$, then $\alpha$ was already an equivalence.
    \end{enumerate}
\end{defn}
In the Stein geometry case $\pi_{0}(\mathbf{P})$ will consist of inclusion of open subdomains of Steins, and modules such that each $\pi_{m}(M)$ is coherent will be the real objects of interest in ${}_{A}\mathbf{Mod}^{\aleph_{1},\mathbf{P},RR}$. 

We also consider descent for subsheaves of perfect and coherent complexes, and give a straightforward generalisation of \cite{toen2008homotopical} Corollary 1.3.2.5 guarantees when \v{C}ech descent implies hyperdescent.

\subsection{Chapter 8: Relative Algebraic Geometry}

In Chapter 8 we give a very general formulation of geometry, starting with geometry relative to a spectral algebraic context in Section \ref{sec:geometryrelac}.In particular we record some basic obstruction theory results in Subsection \ref{subsec:cotangent}. 

In Section \ref{sec:geomlawv} we develop geometry relative to a Lavvere theory $\mathrm{T}$ of homotopy polynomial type. In particular we define the associated $G$-topolgy, as well as the \'{e}tale and smooth topologies relative to $\mathrm{T}$. We also consider base-change of such topologies, and define the $\mathrm{T}$-analytification functor. Finally we consider some basic constructions, which will be elaborated on in further work, of mapping stacks and quotient stacks as well as their cotangent complexes.

\subsection{Chapter 9: Bornological Analytic Geometry}

To conclude in Chapter 9 we specialise the material of Chapter 8 to algebraic geometry relative to the derived algebraic context of complete bornological algebras over a Banach ring $R$. We define the $G$-sites and the \'{e}tale sites. In particular we define the open and \'{e}tale sites for derived Stein spaces, derived dagger affinoid spaces, and derived affinoid spaces relative to any Banach ring $R$ (non-Archimedean in the affinoid case). Over $\mathbb{C}$ we prove that the open Stein site gives a reasonable enhancement of the classical Stein site from analytic geometry. Note that the functor
$$\mathcal{O}:\mathrm{Stein}_{\mathbb{C}}^{op}\rightarrow\mathrm{Comm}(\mathrm{Ind(Ban_{\mathbb{C}})})$$
which takes $X$ to the algebra of analytic functions on $X$ is fully faithful. An object of $\mathrm{Comm}(\mathrm{Ind(Ban_{\mathbb{C}})})$ is a Stein algebra if it is in the essential image of this functor. 

\begin{thm}[Theorem \ref{thm:steintopchar}]
    Let $\{\mathrm{Spec}(A_{i})\rightarrow\mathrm{Spec}(A)\}_{i\in\mathcal{I}}$ be a collection of maps of derived Stein algebras, that is $\pi_{0}(A_{i})$ (resp. $\pi_{0}(A)$) are finitely presented Stein algebras, and each $\pi_{m}(A_{i})$ (resp. $\pi_{m}(A)$ ) is a coadmissible (coherent) module over $\pi_{0}(A_{i})$ (resp. $\pi_{0}(A)$. Then 
    $$\{\mathrm{Spec}(A_{i})\rightarrow\mathrm{Spec}(A)\}_{i\in\mathcal{I}}$$
    is a cover in our open Stein site, if and only if
    \begin{enumerate}
        \item 
        $$\{\mathrm{Spec}(\pi_{0}(A_{i}))\rightarrow\mathrm{Spec}(\pi_{0}(A))\}$$ 
        corresponds to a cover of Stein spaces by open Stein subspaces.
        \item 
        for each $i$ and each $n$ the map
        $$\pi_{n}(A)\otimes_{\pi_{0}(A)}^{\mathbb{L}}\pi_{0}(A_{i})\rightarrow\pi_{n}(A_{i})$$
        is an equivalence.
    \end{enumerate}
\end{thm}

In particular we use this to prove that the category $\mathbf{dAn}^{f,<\infty}$ of locally finitiley embedable ($m$-coskeletal) geometric stacks in the sense of Porta and Yue Yu embed fully faithfully in ($m$-coskeletal) geometric stacks in our Stein context.

\begin{thm}[Theorem \ref{thm:pyycompare}]
There is a fully faithful functor
$$\mathbf{dAn}^{f,<\infty}\rightarrow\mathbf{Stk}_{geom}^{<\infty}(\mathbf{Aff}_{\mathbf{Ch}_{\ge0}(\mathrm{Ind(Ban_{\mathbb{C}})}})$$
\end{thm}

Over a general non-trivially valued Banach field $k$ we also get a contravriant embedding of dagger affinoid spaces, and we have the following result in this context.

\begin{thm}[Theorem \ref{lem:afndoverchar}]
    Let $\{\mathrm{Spec}(A_{i})\rightarrow\mathrm{Spec}(A)\}_{i\in\mathcal{I}}$ be a collection of maps of derived dagger affinoid algebras, that is $\pi_{0}(A_{i})$ (resp.$\pi_{0}(A)$) are finitely presented dagger affinoid algebras, and each $\pi_{m}(A_{i})$ (resp. $\pi_{m}(A)$ ) is a finitely presented module over $\pi_{0}(A_{i})$ (resp. $\pi_{0}(A)$. Then 
    $$\{\mathrm{Spec}(A_{i})\rightarrow\mathrm{Spec}(A)\}_{i\in\mathcal{I}}$$
    is a cover in our open dagger affinoid site, if and only if
    \begin{enumerate}
        \item 
        $$\{\mathrm{Spec}(\pi_{0}(A_{i}))\rightarrow\mathrm{Spec}(\pi_{0}(A))\}$$ 
        corresponds to a cover of dagger affinoid spaces by dagger affinoid subdomains.
        \item 
        for each $i$ and each $n$ the map
        $$\pi_{n}(A)\otimes_{\pi_{0}(A)}^{\mathbb{L}}\pi_{0}(A_{i})\rightarrow\pi_{n}(A_{i})$$
        is an equivalence.
    \end{enumerate}
\end{thm}

In addition for $k$ non-Archimedean we get a contravriant embedding of affinoid spaces, and we get the following.

\begin{thm}[Theorem \ref{lem:affinoidcoverchar}]
    Let $\{\mathrm{Spec}(A_{i})\rightarrow\mathrm{Spec}(A)\}_{i\in\mathcal{I}}$ be a collection of maps of derived affinoid algebras, that is $\pi_{0}(A_{i})$ (resp.$\pi_{0}(A)$) are finitely presented affinoid algebras, and each $\pi_{m}(A_{i})$ (resp. $\pi_{m}(A)$ ) is a finitely presented module over $\pi_{0}(A_{i})$ (resp. $\pi_{0}(A))$. Then 
    $$\{\mathrm{Spec}(A_{i})\rightarrow\mathrm{Spec}(A)\}_{i\in\mathcal{I}}$$
    is a cover in our open affinoid site, if and only if
    \begin{enumerate}
        \item 
        $$\{\mathrm{Spec}(\pi_{0}(A_{i}))\rightarrow\mathrm{Spec}(\pi_{0}(A))\}$$ 
        corresponds to a cover of affinoid spaces by affinoid subdomains.
        \item 
        for each $i$ and each $n$ the map
        $$\pi_{n}(A)\otimes_{\pi_{0}(A)}^{\mathbb{L}}\pi_{0}(A_{i})\rightarrow\pi_{n}(A_{i})$$
        is an equivalence.
    \end{enumerate}
\end{thm}

Finally we consider the \'{e}tale site for a non-trivially valued non-Archimedean Banach ring $k$. In this instance we only get an analogue of Theorem \ref{lem:affinoidcoverchar} if we restrict to those $A$ such that $\pi_{0}(A)$ is sous-perfectoid. The reason for this restriction ultimately boils down to a result of Camargo, \cite{camargo2024analytic} Lemma 3.5.13.

\begin{thm}[Theorem \ref{thm:sousperfetale}]
     Let $\{\mathrm{Spec}(A_{i})\rightarrow\mathrm{Spec}(A)\}_{i\in\mathcal{I}}$ be a collection of maps of derived sous-perfectoid algebras, that is $\pi_{0}(A_{i})$ (resp.$\pi_{0}(A)$) are sous-perfectoid algebras, and each  $\pi_{m}(A_{i})$ (resp. $\pi_{m}(A)$ ) is a finitely presented module over $\pi_{0}(A_{i})$ (resp. $\pi_{0}(A))$. 
       Then  $\mathrm{Spec}(B)\rightarrow\mathrm{Spec}(A)$ is a cover in our \'{e}tale topoloy if and only if 
    \begin{enumerate}
        \item 
        $$\{\mathrm{Spec}(\pi_{0}(A_{i}))\rightarrow\mathrm{Spec}(\pi_{0}(A))\}$$ 
        corresponds to an \'{e}tale cover of affinid in the usual sense.
        \item 
        for each $i$ and each $n$ the map
        $$\pi_{n}(A)\otimes_{\pi_{0}(A)}^{\mathbb{L}}\pi_{0}(A_{i})\rightarrow\pi_{n}(A_{i})$$
        is an equivalence.
    \end{enumerate}
\end{thm}

In all of these contexts we establish hyperdescent results for certain categories of perfect and coherent sheaves.

In Subsection \ref{subsec:cinftygeom} we say a few words about $C^{\infty}$-geometry. In particular using recent results of Carchedi, \cite{carchedi2023derived} we show that the category of derived manifolds embeds fully faithfully in affines relative to $\mathbf{Ch}_{\ge0}(\mathrm{Ind(Ban_{\mathbb{R})})}$. We also give an algebraic characteristaion of the standard pre-topology on $\mathbf{DMfld}$.

In addition we define `universal' Stien, dagger affinoid, and affinoid sites over the integers and prove some base change and analytification results. Finally in Section \ref{sec:classicalsheaves} we relate classical categories of sheaves on analytic spaces to our abstract definition of the category $\mathbf{QCoh}$ of quasi-coherent sheaves.

\subsection{Appendix A: Some Technical Details}

In the appendix we record some technical details concerning model category-theoretic presentations of free sifted cocompletions and functors out of them, as well as some facts about fixed points of adjunctions between stable $(\infty,1)$-categories with $t$-structure.

\section{Notation and Conventions}

We conclude by listing some of our notational conventions.

\begin{itemize}
\item 
We will use the mathrm font $\mathrm{A},\mathrm{B},\mathrm{C},\ldots$ to denote $1$-categories.
\item 
We will use the mathpzc font $\mathpzc{A},\mathpzc{B},\mathpzc{C},\ldots$ to denote model categories.
\item
We will use boldface $\mathbf{A},\mathbf{B},\mathbf{C},\ldots$ to denote $(\infty,1)$-categories.
\item 
When working with $t$-structures and complexes we will use homological grading, so that notation matches with \cite{HA}.
\item When working in $(\infty,1)$-categories the mapping space will be denoted $\mathbf{Map}$. In stable $(\infty,1)$-categories this will usually also denote the mapping spectrum functor. On one or two occassions we will have cause to distinguish bettwen the mapping space and mapping spectrum, and in this case the latter will be denoted $\mathbf{Map}^{s}$.
    \item 
    When working in monoidal $(\infty,1)$-categories we will denote the tensor product by $\otimes^{\mathbb{L}}$, and the internal hom by $\underline{\mathbf{Map}}$. The monoidal unit will be denoted $\mathbb{I}$.
    
    \item 
    In the context of a stable monoidal category with compatible $t$-structure, the tensor product on the heart will be denoted $\otimes$, and the internal hom is $\underline{\mathrm{Hom}}$.
\end{itemize}

\subsection*{Acknowledgements}

The authors would like to thank Federico Bambozzi, Sof\'{i}a Marlasca Aparicio, Rhiannon Savage, and Arun Soor for many useful conversations, and for identifying several errors in earlier versions of this work.

\chapter{Homotopical Algebra}\label{HA}

In this first chapter we develop the fundamentals of homotopical algebra relative to very general classes of symmetric monoidal $(\infty,1)$-categories which we will need in this paper and in its sequels. 

\section{$(\infty,1)$-Algebra Contexts}\label{sec:infalgcontext}
Following \cite{toen2008homotopical}, in this section we introduce  $(\infty,1)$-algebra contexts. These are $(\infty,1)$-category versions of the HA contexts of \cite{toen2008homotopical}. There is a lot of data involved in this definition, so we will separate it into several parts, and explain the intuition behind each. Throughout we let $\mathbf{C}$ be a locally presentable closed symmetric monoidal $(\infty,1)$-category with monoidal unit $\mathbb{I}$. 

\subsection{The Definition}

\subsubsection{The Monad}

  We fix a monad $\mathbf{D}:\mathbf{C}\rightarrow\mathbf{C}$ and a map of monads
$$\theta:\mathbf{Comm}(-)\rightarrow \mathbf{D}$$
such that the induced functor
$$\Theta:\mathbf{Alg}_{\mathbf{D}}\rightarrow\mathbf{Comm}(\mathbf{C})$$
commutes with limits and colimits. For $\mathbf{D}$ as above and $A\in\mathbf{Alg_{D}}(\mathbf{C})$ we define
$${}_{A}\mathbf{Mod}\defeq{}_{\Theta(A)}\mathbf{Mod}$$
The idea behind the monad $\mathbf{D}$ and the map $\theta$ is simply to permit us to deal with simplicial commutative algebras and $E_{\infty}$-algebras at the same time. Indeed in this work all of our examples of such monads of interest will have one of the following forms
\begin{enumerate}
\item
$\mathbf{D}=\mathbf{Comm}(-)$ and
$$\theta:\mathbf{Comm}(-)\rightarrow \mathbf{D}$$
is the identity map. 
\item
$\mathbf{C}$ is presented by an additive model category $\mathrm{C}$ satisfying some assumptions (typically a HA context of \cite{toen2008homotopical}), such that the transferred model structure exists on the category $\mathrm{Comm}(\mathrm{C})$ of commutative monoids in $\mathpzc{C}$. The free-forgetful Quillen adjunction for the transferred model structure presents a moandic adjunction of $(\infty,1)$-categories.
\end{enumerate}

\subsubsection{Trivial Square-Zero Extensions}

Recall that for $A\in\mathbf{Comm}(\mathbf{C})$ we may write
$$\mathbf{Stab}({}_{A}\mathbf{Mod})\cong\mathbf{Stab}(\mathbf{Comm}(\mathbf{Stab}(\mathbf{C})_{\big\slash \Sigma^{\infty}A}))$$
The \textit{trivial square-zero extension} functor $\mathrm{sqz}_{A}:{}_{A}\mathbf{Mod}\rightarrow\mathbf{Comm}(\mathbf{C})_{\big\slash A}$ is defined to be the composition
\begin{displaymath}
\xymatrix{
\mathbf{Stab}(\mathbf{Comm}(\mathbf{Stab}(\mathbf{C})_{\big\slash \Sigma^{\infty}A}))\ar[r]^{\Omega^{\infty}} & \mathbf{Comm}(\mathbf{Stab}(\mathbf{C})_{\big\slash \Sigma^{\infty}A})\ar[r]^{\Omega_{\infty}} & \mathbf{Comm}(\mathbf{C})_{\big\slash A}
}
\end{displaymath}
where both functors are infinite loop functors. Let 
$$\theta:\mathbf{Comm}(-)\rightarrow\mathbf{D}$$
be a map of monads as above, and let $A\in\mathbf{Alg_{D}}(\mathbf{C})$. Consider the category
$$\mathbf{Alg_{D}}(\mathbf{C})_{\big\slash A}$$
There is a limit-preserving functor
$$\Theta_{A}:\mathbf{Alg_{D}}(\mathbf{C})_{\big\slash A}\rightarrow\mathbf{Comm}(\mathbf{C})_{\big\slash\Theta(A)}$$
and therefore a functor
\begin{equation}
    \begin{split}
\Theta^{\mathrm{Mod}}_{A}:\mathbf{Stab}(\mathbf{Alg_{D}}(\mathbf{C})_{\big\slash A})\rightarrow & \mathbf{Stab}(\mathbf{Comm}(\mathbf{C})_{\big\slash\Theta(A)})\\
& \cong\mathbf{Stab}({}_{\Theta(A)}\mathbf{Mod})   \\
&\rightarrow{}_{\Theta(A))}\mathbf{Mod}
    \end{split}
\end{equation}
This is a right-adjoint functor, and therefore admits a left-adjoint which we denote by $\mathrm{sqz}_{A}^{\mathbf{D}}$. There is a natural map
$$\theta_{A}:\Theta(\mathrm{sqz}_{A}^{\mathbf{D}})\rightarrow\mathrm{sqz}_{\Theta(A)}$$

We shall assume that $\theta_{A}$ is an \textit{equivalence}. This allows us to take trivial square-zero extensions of $A\in\mathbf{Alg_{D}}$ by modules $M\in\mathbf{Stab}({}_{\Theta(A)}\mathbf{Mod})$, such that the underlying object of this square-zero extension coincides with the underlying object of the trivial square-zero extension of $\Theta(A)$ by $M$. We will also denote this functor by $\mathrm{sqz}_{A}$.

\subsubsection{The Grading}

The monad $\mathbf{Comm}$ is in fact a graded functor, $\mathbf{Comm}\cong\bigoplus_{n=0}^{\infty}\mathbf{Sym}^{n}(-)$ where $\mathbf{Sym}^{n}(-)$ is the $n$th symmetric power functor. We shall assume that the monad $\mathbf{D}$ is also graded, and that the map $\theta:\mathbf{Comm}(-)\rightarrow\mathbf{D}$ is a map of graded monads.

\subsubsection{The `$t$-structure'}

For the remainder of this subsection fix a subcategory $\mathbf{G}_{0}\subset\mathbf{C}$ satisfying the following conditions. 

\begin{ass}[\cite{toen2008homotopical} Assumption 1.1.0.6]\label{assinf:c0}
\begin{enumerate}
\item
$\mathbb{I} \in\mathbf{G}_{0}$.
\item
$\mathbf{G}_{0}\subset\mathbf{C}$ is stable by isomorphisms and small colimits.
\item
$\mathbf{G}_{0}$ is closed under the tensor product.
\end{enumerate}
\end{ass}

As in \cite{toen2008homotopical}, this is supposed to be some weak version of a $t$-structure. Our main examples of non-trivial $\mathbf{G}_{0}$ will in fact arise from $t$-structures. However for any potential future work on complicial geometry as in \cite{toen2008homotopical}, this greater generality is needed.

\begin{defn}[\cite{toen2008homotopical}, Page 52]
Let $(\mathbf{C},\mathbf{D},\theta,\mathbf{G}_{0},\mathbf{S})$ be a stable $(\infty,1)$-algebra context, and let $A\in\mathbf{Alg_{D}}$. An object $M$ of ${}_{A}\mathbf{Mod}$ is said to be
\begin{enumerate}
\item
$\mathbf{G}_{0}$-\textit{connective} if it is in ${}_{A}\mathbf{Mod}^{\mathbf{G}}_{0}$.
\item
for $n>0$, $\mathbf{G}_{-n}$-\textit{connective} if it is of the form $\Omega M'$ for some $\mathbf{G}_{-(n-1)}$-connective $M'$
\item 
for $n>0$, $\mathbf{G}_{n}$-connective if it is of the form $\Sigma M'$ for some $(n-1)$-connective $M'$
\item
$\mathbf{G}$-\textit{connective} if it is $\mathbf{G}_{-n}$-connective for some $n$.
\end{enumerate}
\end{defn}

For example if $\mathbf{G}_{0}=\mathbf{C}$, then all objects are $\mathbf{G}$-connective. 

Denote by $\mathbf{Alg}_{\mathbf{D},n}$ the fibre product
$$\mathbf{Alg}_{\mathbf{D}}\times_{\mathbf{C}}\mathbf{G}_{n}$$
For $A\in\mathbf{Alg}_{\mathbf{D}}$ and $n\ge0$ denote by ${}_{A}\mathbf{Mod}^{\mathbf{G}}_{n}$ the fibre product
$${}_{A}\mathbf{Mod}\times_{\mathbf{C}}\mathbf{G}_{n}$$
Finally, denote by ${}_{A\big\backslash}\mathbf{Alg}_{\mathbf{D},n}$ the fibre product
$${}_{A\big\backslash}\mathbf{Alg}_{\mathbf{D},n}\times_{\mathbf{C}}\mathbf{G}_{n}$$

Consider the $(\infty,1)$-functor 
$$\mathbf{h}_{0}^{-}:{}_{A}\mathbf{Mod}^{op}\rightarrow\mathbf{Fun}({}_{A}\mathbf{Mod}^{\mathbf{G}}_{0},\mathbf{sSet})$$
defined by the composition
$$\mathbf{h}_{A,0}^{-}:{}_{A}\mathbf{Mod}^{op}\rightarrow\mathbf{Fun}({}_{A}\mathbf{Mod},\mathbf{sSet})\rightarrow\mathbf{Fun}({}_{A}\mathbf{Mod}^{\mathbf{G}}_{0},\mathbf{sSet})$$
where the first functor is the Yoneda embedding, and the second the obvious restriction functor. 

\begin{defn}[\cite{toen2008homotopical} Definition 1.1.0.10]
$A\in\mathbf{Alg_{D}}(\mathbf{C})$ is said to be \textit{good with respect to} $\mathbf{G}_{0}$ if the functor $\mathbf{h}_{A,0}^{-}$ is fully faithful. 
\end{defn}

\subsubsection{Putting It All Together}
\begin{defn}\label{defn:infA}
An $(\infty,1)$-\textit{algebra context} is a tuple 
$$(\mathbf{C},\mathbf{D},\theta,\mathbf{G}_{0},\mathbf{S})$$ 
where
\begin{enumerate}
\item
$\mathbf{C}$ is a locally presentable closed symmetric monoidal $(\infty,1)$-category
\item
$\mathbf{D}$ is a graded monad on $\mathbf{C}$, and 
$$\theta:\mathbf{Comm}(-)\rightarrow\mathbf{D}$$
is a map of graded monads, where $\mathbf{Comm}(-)$ is equipped the symmetric powers grading
\item
 $\mathbf{G}_{0}\subset \mathbf{C}$ and $\mathbf{S}\subset \Comm(\mathbf{C})$ are full subcategories which are stable by isomorphisms. 
\end{enumerate}
such that
\begin{enumerate}
\item
$\mathbf{C}$ is an additive $(\infty,1)$-category (i.e. its homotopy category is additive as a $1$-category).
\item
$$\Theta:\mathbf{Alg}_{\mathbf{D}}\rightarrow\mathbf{Comm}(\mathbf{C})$$
commutes with limits and colimits
\item
$\theta_{A}$ is an equivalence for any $A\in\mathbf{Alg}_{\mathbf{D}}(\mathbf{C})$.
\item
$\mathbf{G}_{0}$ satisfies Assumptions \ref{assinf:c0}
\item
Every $A\in\mathbf{S}$ is $\mathbf{G}_{0}$-good.
\end{enumerate}
\end{defn}

Let is give somewhat trivial, but useful, general classes of examples.

\begin{example}
Let $\mathbf{C}$ be a locally presentable closed symmetric monoidal $(\infty,1)$-category
\begin{enumerate}
\item
Let $\mathrm{Id}:\mathbf{Comm}(-)\rightarrow\mathbf{Comm}(-)$ be the identity functor. then
$$(\mathbf{C},\mathbf{Comm}(-),\mathrm{Id},\mathbf{C},\mathbf{Comm}(\mathbf{C}))$$
is an $(\infty,1)$-algebra context.
\item
Let $\theta:\mathbf{D}\rightarrow\mathbf{Comm}(-)$ be a map of graded monads such that the induced functor $\Theta$ commutes with limits and colimits, and $\theta_{A}$ is an equivalence for each $A\in\mathbf{Alg_{D}}(\mathbf{C})$. Let $\mathbf{S}\subset\mathrm{Alg}_{\mathbf{D}}(\mathbf{C})$ be any full subcategory. Then 
$$(\mathbf{C},\mathbf{D},\theta,\mathbf{C},\mathbf{S})$$
is an $(\infty,1)$-algebra context.
\end{enumerate}
\end{example}

%
%
%
\subsection{Stable $(\infty,1)$-Algebra Contexts and $t$-structures}

In this subsection we introduce some conventions and recall some basic properties of stable $(\infty,1)$-categories with $t$-structure. Let $(\mathbf{C},\mathbf{D},\theta,\mathbf{G}_{0},\mathbf{S})$ be an $(\infty,1)$-algebra context. Consider the symmetric monoidal $(\infty,1)$-category $\mathbf{Stab(\mathbf{C})}$. It is not clear (or indeed true in general) that the monad $\mathbf{D}$ can be extended to one on $\mathbf{Stab(\mathbf{C})}$ to produce a `stabilised' context. However it is true in some important examples, for example the derived algebraic contexts of \cite{raksit2020hochschild} which we will return to later.


\subsubsection{$t$-Structures}

If $\mathbf{C}$ is a stable $(\infty,1)$-category we denote a $t$-structure on $\mathbf{C}$ by $(\mathbf{C}_{\ge0},\mathbf{C}_{\le0})$. Note we use homological grading conventions. We then say that  $(\mathbf{C},\mathbf{C}_{\ge0},\mathbf{C}_{\le0})$ is a stable $(\infty,1)$-category with $t$-structure. The homology functors associated to the $t$-structure will be denoted by $\pi_{n}:\mathbf{C}\rightarrow\mathbf{C}^{\heart}$. We will also write $\pi_{*}\defeq\bigoplus_{n\in\mathbb{Z}}\pi_{n}$ and regard it as a graded object of $\mathbf{C}^{\heart}$.

\begin{defn}[\cite{raksit2020hochschild}]
Let $\mathbf{C}$ be a stably monoidal $(\infty,1)$-category with unit $\mathbb{I}$. A $t$-structure $(\mathbf{C}_{\ge0},\mathbf{C}_{\le0})$ on $\mathbf{C}$ is said to be \textit{compatible} if 
\begin{enumerate}
\item
$\mathbf{C}_{\le0}$ is closed under filtered colimits in $\mathbf{C}$.
\item
$\mathbf{C}_{\ge0}$ is closed under the tensor product.
\item
$\mathbb{I}$ is contained in $\mathbf{C}_{\ge0}$. 
\end{enumerate}
\end{defn}

\begin{example}
Let $\mathbf{C}$ be a presentably symmetric monoidal $(\infty,1)$-category in which filtered colimits commute with finite limits. Consider its stabilisation $\mathbf{Stab}(\mathbf{C})$. Then the standard induced $t$-structure is compatible. 
\end{example}

\begin{defn}\label{defn:stablecontext}
Let $(\mathbf{C},\mathbf{D},\theta,\mathbf{G}_{0},\mathbf{S})$ be a stable $(\infty,1)$-algebra context. A $t$-structure $(\mathbf{C}_{\ge0},\mathbf{C}_{\le0})$ on $\mathbf{C}$ is said to be \textit{compatible} if 
\begin{enumerate}
\item
it is left and right complete
\item
the monoidal structure on $\mathbf{C}$ is compatible with the $t$-structure.
\item 
There is an increasing sequence $0\le a(n)\le n$ such that for every $k>1$ and $M\in\mathbf{C}_{\ge n}$, $\mathbf{D}^{k}(M)\in\mathbf{C}_{\ge n+a(n)}$
\item
The functor
$$\mathbf{Alg_{D}}(\mathbf{C})\times_{\mathbf{C}}\mathbf{C}^{\heart}\rightarrow \mathbf{Comm}(\mathbf{C})\times_{\mathbf{C}}\mathbf{C}^{\heart}$$
is an equivalence.
\item
$\tau_{\ge0}(G)\in\mathbf{G}_{0}$ for any $G\in\mathbf{G}_{0}$.
\end{enumerate}

\end{defn}

Let us explain why these conditions are reasonable. The first condition implies that $\mathbf{C}\cong\mathbf{Stab}(\mathbf{C}_{\ge0})$. The second condition means that $\mathbf{Comm}(-)$ satisfies conditions $(2)$ and $(3)$ where in $(3)$ we can take $a(n)=n$. The fourth condition simply means that $\mathbf{D}$-algebras should be thought of as some kind of derived enhancement of commutative algebras in the heart. Let us deduce some useful properties.

The first is obvious from condition $(3)$.

\begin{prop}
    $\mathbf{D}$ restricts to a functor $\mathbf{C}_{\ge0}\rightarrow\mathbf{C}_{\ge0}$.
\end{prop}

This guarantees that $\mathbf{D}$ restricts to a monad $\mathbf{D}:\mathbf{C}_{\ge0}\rightarrow\mathbf{C}_{\ge0}$ so that we have
$$\mathbf{Alg}^{cn}_{\mathbf{D}}(\mathbf{C}_{\ge0})\defeq\mathbf{Alg}_{\mathbf{D}}(\mathbf{C})\times_{\mathbf{C}}\mathbf{C}_{\ge0}\cong\mathbf{Alg}_{\mathbf{D}}(\mathbf{C}_{\ge0})$$

\begin{defn}
    A map $f:M\rightarrow N$ in $\mathbf{C}$ is said to be $n$-\textit{connective} if it induces a surjection on $\pi_{n}$ and an isomorphism on $\pi_{m}$ for $m<n$.
\end{defn}

\begin{lem}[c.f. \cite{HA} Lemma 7.4.3.16]
    Let $f:A\rightarrow B$ be a morphism in $\mathbf{Alg_{D}}(\mathbf{C}_{\ge0})$ and let $M,N\in{}_{B}\mathbf{Mod}(\mathbf{C}_{\ge0})$. If $\mathrm{cofib}(f)$ is $n$-connective then the induced map $\theta:M\otimes_{A}^{\mathbb{L}}N\rightarrow M\otimes_{B}^{\mathbb{L}}N$ is $n$-connective. 
\end{lem}

\begin{proof}
    This works exactly as in \cite{HA} Lemma 7.4.3.16.
\end{proof}

\begin{cor}
For any $n\ge0$ and any $A\in\mathbf{C}_{\ge0}$ the map
$$\tau_{\le n}\mathbf{D}(A)\rightarrow\tau_{\le n}\mathbf{D}(\tau_{\le n}A)$$
is an equivalence.
\end{cor}

\begin{proof}
Let $F$ denote the fibre of $A\rightarrow\tau_{\le n}A_{n}$, which is also $n+1$ connective. We then have an equivalence
$$\mathbf{D}(A)\otimes_{\mathbf{D}(F)}^{\otimes \mathbb{L}}\mathbb{I}\cong\mathbf{D}(A_{\le n})$$
Now 
$$\mathbf{D}(A)\cong\bigoplus_{k=0}^{\infty}\mathbf{D}^{k}(A)\cong\mathbb{I}\oplus \bigoplus_{k=1}^{\infty}\mathbf{D}^{k}(M)$$
By assumption $\bigoplus_{k=1}^{\infty}\mathbf{D}^{k}(M)$ is in $\mathbf{C}_{\ge n+1}$. Thus cofibre of the map $\mathbb{I}\rightarrow \mathbf{D}(F)$ is $n+1$-connective, and the the map 
$$\mathbf{D}(A)\cong\mathbf{D}(A)\otimes_{\mathbb{I}}^{\mathbb{L}}\mathbb{I}\rightarrow\mathbf{D}(A)\otimes_{\mathbf{D}(F)}^{\otimes \mathbb{L}}\mathbb{I}\cong\mathbf{D}(A_{\le n})$$
is $n+1$-connective, as required. 
\end{proof}

This in particular means that for any $A\in\mathbf{C}_{\ge0}$, $\tau_{\le n}(A)$ may be canonically endowed with the structure of a $\mathbf{D}$-algebra. We therefore get a Postnikov tower of $\mathbf{D}$-algebras:

\begin{defn}
Let $(\mathbf{C},\mathbf{D},\theta,\mathbf{G}_{0},\mathbf{S})$ be a stable $(\infty,1)$-algebra context equipped with a compatible $t$-structure $(\mathbf{C}_{\ge0},\mathbf{C}_{\le0})$, and $A\in\mathrm{Comm}(\mathbf{C})$. The \textit{Postnikov tower} of $A$ is
$$A\rightarrow\cdots\rightarrow A_{\le n}\rightarrow A_{\le n-1}\rightarrow\cdots\rightarrow A_{\le0}=\pi_{0}(A)$$
\end{defn}

 The first and final condition guarantee that there is an associated connected $(\infty,1)$-algebraic context, using the following straightforward result.

\begin{prop}
Let $(\mathbf{C},\mathbf{D},\theta,\mathbf{G}_{0},\mathbf{S})$ be a stable $(\infty,1)$-algebra context, and $(\mathbf{C}_{\ge0},\mathbf{C}_{\le0})$ a compatible $t$-structure. Then 
$$(\mathbf{C}_{\ge0},\mathbf{D},\theta,\mathbf{G}_{0}\cap\mathbf{C}_{\ge0},\mathbf{S}\cap\mathbf{C}_{\ge0})$$
 is an $(\infty,1)$-algebra context.
\end{prop}

Define
$$\mathbf{Alg_{D}}(\mathbf{C}^{\heart})\defeq\mathbf{Alg_{D}}(\mathbf{C})\times_{\mathbf{C}}\mathbf{C}^{\heart}$$
We have the following $(\infty,1)$-algebra context associated to $\mathbf{C}$. 

\begin{lem}[c.f. \cite{toen2008homotopical} 2.3.1.1]
$(\mathbf{C},\mathbf{D},\theta,\mathbf{C}_{\ge0},\mathbf{Alg}_{\mathbf{D}}(\mathbf{C}^{\heart}))$ is an $(\infty,1)$-algebra context.
\end{lem}

\begin{proof}
Let $A\in\mathbf{Alg_{D}}(\mathbf{C}^{\heart})$. As in \cite{toen2008homotopical} Lemma 2.3.1.1, since $A$ is in the heart, by replacing $\mathbf{C}$ with ${}_{A}\mathbf{Mod}$ we may assume that $A=\mathbf{D}(0)$ Thus we need to show that the functor
$$\mathbf{C}^{op}\rightarrow\mathbf{Fun}(\mathbf{C}_{\ge0},\mathbf{sSet})$$
is fully faithful. This follows since the $t$-structure on $\mathbf{C}$ is left complete. 
\end{proof}

The existence of a (compatible) $t$-structure $(\mathbf{C}_{\ge0},\mathbf{C}_{\le0})$ on a stable $(\infty,1)$-algebra context $(\mathbf{C},\mathbf{D},\theta,\mathbf{G}_{0},\mathbf{S})$ of course allows us to say when an object of $\mathbf{C}$ is $n$-connective. However particularly for obstruction theory in the geometry of stable algebras, it is useful to have the other notion of connectivity determined by $\mathbf{G}_{0}$.

\begin{example}
    Let $(\mathbf{C},\mathbf{C}_{\ge0},\mathbf{C}_{\le0})$ be a closed symmetric monoidal stable locally presentable $(\infty,1)$-category equipped with a $t$-structure which is compatible with the monoidal structure. Then the $t$-structure is compatible with the stable $(\infty,1)$-algebra context $(\mathbf{C},\mathbf{Comm}(-),\mathrm{Id},\mathbf{C}_{\ge0},\mathbf{Comm}(\mathbf{C}))$.
\end{example}

\begin{proof}
    The only claim which is not entirely trivial is that for any $n\ge0$, any $M\in\mathbf{C}_{\ge n}$ and any $k\ge 2$, $\mathbf{D}^{k}(M)\in\mathbf{C}_{\ge 2n}$. But clearly $M^{\otimes^{\mathbb{L}k}}$ is in $\mathbf{C}_{\ge 2n}$ whenever $k\ge 2$. Since $\mathbf{C}_{\ge 2n}$ is closed under colimits the claim follows.
\end{proof}

\subsection{The Cotangent Complex}

We now introduce cotangent complexes of maps of $\mathbf{D}$-algebras. Because of our assumptions on $\mathbf{D}$, the theory behaves somewhat similarly to the case $\mathbf{D}=\mathbf{Comm}(-)$. In particular much of \cite{HA} Section 7.3 works in this more general context, particular under the Koszul assumption, as we shall see below.

 Let $\mathbb{T}_{\mathbf{Alg_{D}}}$ denote the tangent bundle of the category $\mathbf{Alg_{D}}$, as defined in \cite{HA} Definition 7.3.1.9. An object of $\mathbb{T}_{\mathbf{Alg_{D}}}$ may be thought of as a pair $(A,M)$ where $A\in\mathbf{Alg_{D}}(\mathbf{C})$, and $M$ is a $\Theta(A)$-module. There is a projection functor $\pi_{1}:\mathbb{T}_{\mathbf{Alg_{D}}}\rightarrow\mathbf{Alg_{D}}$ (informally sending $(A,M)$ to $A$) such that for $A\in\mathbf{Alg_{D}}$, we have an equivalence
 $$\mathbb{T}_{\mathbf{Alg_{D}}}\times_{\mathbf{Alg_{\mathbf{D}}}}\{A\}\cong{}_{\Theta(A)}\mathbf{Mod}$$
There is also a `forgetful functor'
$$R:\mathbb{T}_{\mathbf{Alg_{D}}}\rightarrow\mathbf{Fun}(N(\Delta^{1}),\mathbf{Alg_{D}})\rightarrow\mathbf{Fun}(\{0\}),\mathbf{Alg_{D}})$$
The restriction of $R$ to the fibre over $A$ maps to $\mathbf{Alg_{D}}_{\big\slash A}$, and is given by square-zero extension. 

As explained in \cite{HA} 7.3, the functor $R$ has a left adjoint $\mathbb{L}$, such that the composite $\pi_{1}\circ\mathbb{L}$ is naturally equivalent to the identity functor. Thus the value of $\mathbb{L}$ on $A\in\mathbf{Alg_{D}}$, denoted $\mathbb{L}_{A}$, may naturally be regarded as an object of ${}_{\Theta(A)}\mathbf{Mod}$. 

\begin{defn}
$\mathbb{L}_{A}$ is called the \textit{cotangent complex of }$A$. 
\end{defn}

By construction, we have a natural isomorphism of functors ${}_{\Theta(A)}\mathbf{Mod}\rightarrow\mathbf{sSet}$
$$\mathbf{Map}_{{}_{\Theta(A)}\mathbf{Mod}}(\mathbb{L}_{A},(-))\cong\mathbf{Map}_{\mathbf{Alg_{D}}(\mathbf{C})_{\big\slash A}}(A,\mathrm{sqz}_{A}(-))$$
We write
$$\mathbf{Der}_{A}(A,M)\defeq\mathbf{Map}_{\mathbf{Alg_{D}}(\mathbf{C})_{\big\slash A}}(A,\mathrm{sqz}_{A}(M))$$

\begin{example}\label{ex:cotangentfree}
Let $P$ be an object of $\mathbf{C}$, and consider the free algebra $\mathbf{D}(P)$. We claim that  $\mathbb{L}_{\mathbf{D}(P)}\cong\mathbf{D}(P)\otimes^{\mathbb{L}}P$. Indeed we have
\begin{align*}
\mathbf{Map}_{{}_{\mathbf{D}(P)}\mathbf{Mod}}(\mathbf{D}(P)\otimes^{\mathbb{L}}P,N)&\cong\mathbf{Map}_{\mathbf{C}}(P,N)\\
&\cong\mathbf{Map}_{\mathbf{C}_{\big\slash \mathbf{D}(P)}}(P,\mathbf{D}(P)\oplus N)\\
&\cong\mathbf{Map}_{\mathbf{Alg_{D}}_{\big\slash \mathbf{D}(P)}}(\mathbf{D}(P),\mathbf{D}(P)\oplus N)\\
&\cong\mathbf{Map}_{{}_{\mathbf{D}(P)}\mathbf{Mod}}(\mathbb{L}_{\mathbf{D}(P)},N)
\end{align*}
\end{example}

Note that if $A\rightarrow B$ is a map in $\mathrm{Alg}_{\mathbf{D}}(\mathbf{C})$, then there is a natural map $\mathbb{L}_{A}\rightarrow\mathbb{L}_{B}$ and therefore a natural map of $B$-modules
$$B\otimes_{A}^{\mathbb{L}}\mathbb{L}_{A}\rightarrow\mathbb{L}_{B}$$

\begin{defn}
Let $A\rightarrow B$ be a map in $\mathrm{Alg}_{\mathbf{D}}(\mathbf{C})$. The \textit{relative cotangent complex of }$f$, denoted $\mathbb{L}_{B\big\slash A}$, is the cofibre in $\mathbf{Mod}_{B}$ of the map 
$$B\otimes_{A}^{\mathbb{L}}\mathbb{L}_{A}\rightarrow\mathbb{L}_{B}$$
\end{defn}

The object $\mathbb{L}_{B\big\slash A}$ corepresents the functor
$$\mathbf{Der}_{A}(B,-)\defeq\mathbf{Map}_{{}_{A\big\backslash}\mathrm{Alg}_{\mathbf{D}}(\mathbf{C})_{\big\slash B}}(B,\mathrm{sqz}_{B}(-))$$

Let $\mathcal{M}^{\mathbb{T}}(\mathbf{Alg}_{\mathbf{D}}(\mathbf{C}))$ be a tangent correspondence (this is nothing but the category fibered over $N(\Delta^{1})$ encoding the ajdunction $L\dashv R$, see \cite{HA} Definition 7.3.6.9). It contains both $\mathbf{Alg}_{\mathbf{D}}(\mathbf{C})$ and $\mathbb{T}_{\mathbf{Alg}_{\mathbf{D}}(\mathbf{C})}$ as full subcategories, and is equipped with a projection map 
$$\pi:\mathcal{M}^{\mathbb{T}}(\mathbf{Alg}_{\mathbf{D}}(\mathbf{C}))\rightarrow N(\Delta^{1})\times\mathbf{Alg}_{\mathbf{D}}(\mathbf{C})$$ 
Denote by $\underline{\mathbf{Der}}(\mathbf{Alg}_{\mathbf{D}}(\mathbf{C}))$ the fibre product
$$\mathbf{Fun}(N(\Delta^{1}),\mathcal{M}^{\mathbb{T}}(\mathbf{Alg}_{\mathbf{D}}(\mathbf{C})))\times_{\mathbf{Fun}(N(\Delta^{1}),N(\Delta^{1})\times\mathbf{Alg_{D}}(\mathbf{C}))}\mathbf{Alg_{D}}(\mathbf{C})$$
An object of this category may be identified with a derivation $d:A\rightarrow A\oplus M$, where $A\in\mathbf{Alg_{D}}(\mathbf{C})$ and $M\in{}_{\Theta(A)}\mathbf{Mod}$.

By \cite{HA} Section 7.4 there is a functor
$$\Phi:\underline{\mathbf{Der}}(\mathbf{Alg}_{\mathbf{D}}(\mathbf{C}))\rightarrow\mathbf{Fun}(N(\Delta^{1}),\mathbf{Alg}_{\mathbf{D}}(\mathbf{C}))$$
which sends a derivation $d:A\rightarrow A\oplus M$ to $(A\oplus_{d}\Omega M\rightarrow A)$, where $A\oplus_{d}\Omega M$ is given by the pull-back diagram below \begin{displaymath}
\xymatrix{
A\oplus_{d}\Omega M\ar[d]\ar[r] & A\ar[d]^{d_{0}}\\
A\ar[r]^{d} & A\oplus M
}
\end{displaymath}
Here $\eta_{0}:A\rightarrow A\oplus M$ is the derivation classified by the zero map $0:\mathbb{L}_{A}\rightarrow M$. 

Exactly as in \cite{toen2008homotopical} Proposition 1.4.2.5 one can prove the following.

\begin{prop}\label{prop:lifting}
Let $f:A\rightarrow B$ and $g:B\rightarrow C$ be maps in $\mathbf{Alg_{D}}(\mathbf{C})$. 
Let $d:C\rightarrow M$ be a derivation. 
\begin{displaymath}
    \xymatrix{
A\ar[d]^{f}\ar[r]^{p} & C\oplus_{d}\Omega M\ar[d]^{q}\\
B\ar[r]^{g} & C
    }
\end{displaymath}
be a commutative diagram, which we call $x$.
\begin{enumerate}
\item
There exists a natural obstruction
$$\alpha(x)\in\pi_{0}(\mathbf{Map}_{{}_{A}\mathbf{Mod}(\mathbf{C})}(f^{*}\mathbb{L}_{B\big\slash A},M))$$
which vanishes if and only if there exists a lift in the diagram $x$, i.e. a map $g_{d}:B\rightarrow C\oplus_{d}\Omega M$ satisfying $g_{d}\circ f= q$ and $q\circ g_{d}=g$ in $\mathrm{Ho}({}_{A\big\backslash}\mathbf{Alg_{D}})$. 
\item
If $\alpha(x)=0$ then the space of lifts is non-canonically equivalent to 
$$\Omega_{\alpha(x),0}\mathbf{Map}_{{}_{A}\mathbf{Mod}(\mathbf{C})}(f^{*}\mathbb{L}_{B\big\slash A},M)$$
\end{enumerate}
\end{prop}

\subsubsection{$t$-Structures and The Cotangent Complex}

Now suppose that $(\mathbf{C},\mathbf{D},\theta)$ is a stable $(\infty,1)$-algebra context equipped with acompatible $t$-structure $(\mathbf{C}_{\ge0},\mathbf{C}_{\le0})$. 
Our first important result concerns connectivity of the cotangent complex. Let $f:A\rightarrow B$ be a map of $\mathbf{D}$-algebras, and consider the induced map of $B$-modules
$$\eta:\mathbb{L}_{B}\rightarrow\mathbb{L}_{B\big\slash A}$$
Let $B^{\eta}$ denote the corresponding square-zero extension of $B$ by $\mathbb{L}_{B\big\slash A}[1]$. $f$ factors as 
\begin{displaymath}
\xymatrix{
A\ar[r]^{f'} & B^{\eta}\ar[r]^{f''} & B
}
\end{displaymath}
This gives rise to a map of $A$-modules $\mathrm{cofib}(f)\rightarrow\mathrm{cofib}(f'')$. Consider the adjoint map of $B$-modules
$$\epsilon_{f}:B\otimes^{\mathbb{L}}_{A}\mathrm{cofib}(f)\rightarrow\mathrm{cofib}(f'')\cong\mathbb{L}_{B\big\slash A}$$
often called the \textit{Hurewicz map}. Here we analyse this map. Much of this section closely follows \cite{HA} Chapter 7.4, but we would also like to give due credit to forthcoming work of Marlasca Aparicio \cite{SMAUltrasolid}. They are doing similar work in the context of ultrasolid modules in the sense of condensed mathematics which they kindly shared with us, and pointed out an important error in a previous version of the present paper.

\begin{thm}[\cite{HA} Theorem 7.4.3.12]\label{thm:connectivcot}
 Let $f:A\rightarrow B$ be a morphism in $\mathbf{Alg_{D}}(\mathbf{C}_{\ge0})$ and suppose that $\mathrm{cofib}(f)$ is $n$-connective where $n\ge 0$. Then
$\mathrm{fib}(\epsilon_{f})\in\mathbf{C}_{\ge n+a(n)}$.

\end{thm}

We shall prove this in several steps.

\begin{lem}[c.f. \cite{HA} Lemma 7.4.3.15]\label{lem:Madjoin}
    Let $f:A\rightarrow B$ be a morphism in $\mathbf{Alg_{D}}(\mathbf{C}_{\ge0})$ such that $\mathrm{cofib}(f)$ is $n$-connective for some $n\ge0$. Then there exists an object $M\in\mathbf{C}_{\ge n-1}$ and a commutative diagram
    \begin{displaymath}
        \xymatrix{
\mathbf{D}(M)\ar[d]\ar[r]^{\phi} & \mathbb{I}\ar[d]\\
A\ar[dr]^{f}\ar[r] & A'\ar[d]^{f'}\\
& B
        }
    \end{displaymath}
    in $\mathbf{Alg_{D}}(\mathbf{C})$ where the upper square is a pushout $A'\in\mathbf{Alg_{D}}^{cn}(\mathbf{C})$, $\mathrm{cofib}(f')\in\mathbf{C}_{\ge n+1}$ and $\phi$ is adjoint to the zero map $M\rightarrow\mathbb{I}$.
\end{lem}

\begin{proof}
The proof of this works exactly as in \cite{HA} Lemma 7.4.3.15. All that is required here is that whenever $M$ is $n$ connective for $n\ge 0$, then $\mathbf{D}(M)$ is $n$-connective.
\end{proof}

Then exactly as in the proof of \cite{HA} Lemma 7.4.3.12 we have the following,

\begin{cor}\label{cor:algpres}
Let
$$f:A\rightarrow B$$ 
be a map in $\mathbf{Alg_{D}}^{cn}(\mathbf{C})$ such that $\mathrm{cofib}(f)$ is $n$-connective for $n\ge0$. There is a sequence of objects
$$A_{n}\rightarrow A_{n+1}\rightarrow A_{n+2}\rightarrow\ldots$$
 in $\mathbf{Alg_{D}}^{cn}(\mathbf{C})_{\big\slash B}$ with
 \begin{enumerate}
     \item 
     $A_{n}=A$
     \item $\mathrm{cofb}(A_{m}\rightarrow B)\in\mathbf{C}_{\ge m}$.
     \item 
     for each $m\ge n$ there is an object $M\in\mathbf{C}_{\ge m-1}$ and a pushout diagram
     \begin{displaymath}
\xymatrix{
\mathbf{D}(M)\ar[d]\ar[r]^{\phi_{m}} & \mathbb{I}\ar[d]\\
A_{m}\ar[r]^{g_{m,m+1}} & A_{m+1}
}         
     \end{displaymath}
     where $g_{j,k}$ denotes the map $A_{j}\rightarrow A_{k}$ in our sequence, and $\phi_{m}$ is adjoint to the zero map $M\rightarrow 0$ in $\mathbf{C}$.
 \end{enumerate}
 In particular $\colim_{m}A_{m}\rightarrow B$ is an equivalence. 
\end{cor}

\begin{lem}[c.f. \cite{HA} Lemma 7.4.3.17]
    Let $f:A\rightarrow B$ be a morphism in $\mathbf{Alg_{D}}(\mathbf{C}_{\ge0})$. Suppose that for $n\ge0$
 $f$ induces an equivalence $\tau_{\le n}A\rightarrow\tau_{\le n}B$. Then $\tau_{\le n}\mathbb{L}_{B\big\slash A}\cong 0$
 \end{lem}

\begin{proof}
    This works exactly as in \cite{HA} Lemma 7.4.3.17.
\end{proof}

\begin{proof}[proof of Theorem \ref{thm:connectivcot}]
Say that a morphism $f:A\rightarrow B$ in $\mathbf{Alg_{D}}(\mathbf{C})$ is $n$-\textit{good} if $\mathrm{fib}(\epsilon_{f})\in\mathbf{C}_{\ge n+a(n)}$.
\begin{enumerate}
    \item
    Exactly as in the proof of \cite{HA} Theorem 7.4.3.12, part (a), if $f:A\rightarrow B$ and $g:B\rightarrow C$ are $n$-good, and $\mathrm{cofib}(f),\mathrm{cofib}(g)\in\mathbf{C}_{\ge n}$ then $h$ is $n$-good. Indeed as in loc. cit. we consider the commutative diagram
    \begin{displaymath}
    \xymatrix{
        C\otimes_{A}^{\mathbb{L}}\mathrm{cofib}(f)\ar[d]^{\epsilon'}\ar[r] &  C\otimes_{A}^{\mathbb{L}}\mathrm{cofib}(h)\ar[d]^{\epsilon_{h}}\ar[r] &  C\otimes_{A}^{\mathbb{L}}\mathrm{cofib}(g)\ar[d]^{\epsilon''}\\
        C\otimes_{B}^{\mathbb{L}}\mathbb{L}_{B\big\slash A}\ar[r] & \mathbb{L}_{C\big\slash A}\ar[r] & \mathbb{L}_{C\big\slash B}
}
    \end{displaymath}
    where both the top and bottom rows are fibre sequences of $C$-modules, $\epsilon'\defeq\mathrm{Id}_{C}\otimes_{B}^{\mathbb{L}}\epsilon_{f}$ and $\epsilon''$ is the composition
\begin{displaymath}
    \xymatrix{
    C\otimes_{A}^{\mathbb{L}}\mathrm{cofib}(g)\ar[r]^{\phi} & C\otimes_{B}^{\mathbb{L}}\mathrm{cofib}(g)\ar[r]^{\epsilon_{g}} & \mathbb{L}_{C\big\slash B}
    }
\end{displaymath}
It suffices to show $\mathrm{fib}(\epsilon'')\in\mathbf{C}_{\ge n+a(n)}$
We get a fibre sequence
$$\mathrm{fib}(\phi)\rightarrow\mathrm{fib}(\epsilon'')\rightarrow\mathrm{fib}(\epsilon_{g})$$
Now it suffices to prove that $\mathrm{fib}(\phi)\in\mathbf{C}_{\ge n+a(n)}$, as $\mathrm{fib}(\epsilon_{g})$ is by assumption. But $\mathrm{cofib}(g)$ is $n$-connected, and $\mathrm{cofib}(f)$ is $n$-connected. Thus $\mathrm{fib}(\phi[-n])$ is $n$-connected, so $\mathrm{fib}(\phi)$ is $2n\ge n+a(n)$-connected.
\item
 Also exactly as in the proof of \cite{HA} Theorem 7.4.3.12, part (b). If $f:A\rightarrow B$ is $n$-good and $A\rightarrow A'$ is any map then $A'\otimes_{A}^{\mathbb{L}}B\rightarrow B'$ is $n$-good. 
\item
Again exactly as in part (c) of the proof \cite{HA} Theorem 7.4.3.12, a  morphism $f:A\rightarrow B$ the domain $B\otimes_{A}^{\mathbb{L}}\mathrm{cofib}(f)$ of $\epsilon_{f}$ can be identified with the cofibre of the map $B\rightarrow B\otimes_{A}^{\mathbb{L}}B$ also works as in \cite{HA}. Moreover this latter maps clearly admits a left homotopy inverse. 
\item
As in part (d) of the proof of \cite{HA} Theorem 7.4.3.12, let $M\in\mathbf{C}_{\ge n-1}$ and consider the map $\mathbf{D}(M)\rightarrow\mathbb{I}$ which is adjoint to the zero map $M\rightarrow\mathbb{I}$. Then this is is $n$-good. Indeed we have a fibre sequence
$$\mathbb{I}\otimes_{\mathbf{D}(M)}\mathbb{L}_{\mathbf{D}(M)}\rightarrow\mathbb{L}_{\mathbb{I}}\rightarrow\mathbb{L}_{\mathbb{I}\big\slash\mathbf{D}(M)}$$
This fibre sequence is equivalent to 
$$M\rightarrow 0\rightarrow\mathbb{L}_{\mathbf{D}(M)}$$
so that the codomain of $\epsilon_{f}$ is given by $M[1]$. We then have
$$\mathbb{I}\otimes_{\mathbf{D}(M)}\mathbb{I}\cong\mathbf{D}(M[1])$$
$\mathbb{I}\otimes^{\mathbb{L}}_{\mathbf{D}(M)}\mathrm{cofib(f)}$ can be identified with the cofibre of $\mathbb{I}\rightarrow\mathbf{D}(M[1])\cong\bigoplus_{i=1}^{\infty}\mathbf{D}^{i}(M[1])$. The composition
$$M[1]\cong\mathbf{D}(M[1])\rightarrow\bigoplus_{i>0}\mathbf{D}^{i}(M[1])\rightarrow M[1]\cong\mathbb{L}_{\mathbb{I}\big\slash\mathbf{D}(M)}$$
is then the identity, so that $\mathrm{fib}(\epsilon_{f})\cong\bigoplus_{i\ge 2}\mathbf{D}^{i}(M[1])$.
 Now $M[1]$ is concentrated in degrees $\mathbf{C}_{\ge n}$, so $\mathbf{D}^{k}(M[1])$ is in $\mathbf{C}_{\ge n+a(n)}$ for all $k\ge 2$, as required.
 \item 
 Suppose that $f:A\rightarrow B$ is a map of connective algebras such that $\tau_{\le n+a(n)+1}A\rightarrow\tau_{\le n+a(n)+1}B$ is an equivalence. Then $B\otimes_{A}^{\mathbb{L}}\mathrm{cofib}(f)$ and $\mathbb{L}_{B\big\slash A}$ are both in $\mathbf{C}_{\ge n+a(n)+1}$, and so the fibre of $B\otimes_{A}^{\mathbb{L}}\mathrm{cofib}(f)\rightarrow \mathbb{L}_{B\big\slash A}$ will be in $\mathbf{C}_{\ge n+a(n)}$ as required.
 \end{enumerate}
 The conclusion of the proof also works as in \cite{HA}, using the analog of Lemma 7.4.3.15. So, let $f:A\rightarrow B$ be a map in $\mathbf{Alg_{D}}^{cn}\mathbf{(C)}$ such that $\mathrm{cofib}(f)\in\mathbf{C}_{\ge n}$. By Corollary \ref{cor:algpres} we construct a sequence of objects
 $$A_{n}\rightarrow A_{n+1}\rightarrow A_{n+2}\rightarrow\ldots$$
 in $\mathbf{Alg_{D}}^{cn}(\mathbf{C})_{\big\slash B}$ with
 \begin{enumerate}
     \item 
     $A_{n}=A$
     \item $\mathrm{cofb}(A_{m}\rightarrow B)\in\mathbf{C}_{\ge m}$.
     \item 
     for each $m\ge n$ there is an object $M\in\mathbf{C}_{\ge m-1}$ and a pushout diagram
     \begin{displaymath}
\xymatrix{
\mathbf{D}(M)\ar[d]\ar[r]^{\phi_{m}} & \mathbb{I}\ar[d]\\
A_{m}\ar[r]^{g_{m,m+1}} & A_{m+1}
}         
     \end{displaymath}
     where $g_{j,k}$ denotes the map $A_{j}\rightarrow A_{k}$ in our sequence, and $\phi_{m}$ is adjoint to the zero map $M\rightarrow 0$ in $\mathbf{C}$.
 \end{enumerate}
 Now $A_{n+a(n)+1}\rightarrow B$ is $n$-good. It then suffices to prove that $g_{m,m+1}$ is $n$-good for each $n$, and hence that $\phi_{m}$ is $n$-good for each $m$. But this is what we showed above.
\end{proof}

This result implies some immediate very important corollaries, all of which can be proven exactly as in \cite{HA}.

\begin{cor}[\cite{HA} Corollary 7.4.3.2]
Let $(\mathbf{C}_{\ge0},\mathbf{C}_{\le0})$ be a $t$-structure on a presentably symmetric monoidal stable $(\infty,1)$-category $\mathbf{C}$ which is compatible with the monoidal structure. If $f:A\rightarrow B$ is a map in $\mathbf{Alg}_{\mathbf{D}}(\mathbf{C})$ such that $\mathrm{cofib}(f)$ is $n$-connective for $n\ge0$. Then $\mathbb{L}_{B\big\slash A}$ is $n$-connective. The converse is true if 
\begin{enumerate}
    \item
    $\pi_{0}(A)\rightarrow\pi_{0}(B))$ is an isomorphism.
    \item 
    $a(n)>0$ for $n>0$.
\end{enumerate}
\end{cor}

\begin{proof}
    We have a fibre sequence of $B$-modules.
    $$\mathrm{fib}(\epsilon_{f})\rightarrow B\otimes_{A}^{\mathbb{L}}\mathrm{cofib}(f)\rightarrow\mathbb{L}_{B\big\slash A}$$
    It suffices to prove that $\mathrm{cofib}(f)\otimes_{A}^{\mathbb{L}}B$ is $n$-connective and $\mathrm{fib}(\epsilon_{f})$ is $(n-1)$-connective. Only the latter claim is non-trivial, but this follows because $n+a(n)\ge n\ge n-1$.

    Conversely suppose $\mathrm{cofib}(f)$ is not $n$-connective. Let $n$ be minimal, so that $\mathrm{cofib}(f)$ $\mathrm{cofib}(f)$ is $(n-1)$-connective. We must have $n\ge 2$ by assumption. Thus $\epsilon_{f}$ is $(n-1)+a(n-1)$-connective. Now $(n-1)+a(n-1)\ge n$ since $n-1>0$. Thus
    $$\pi_{n-1}(B\otimes_{A}^{\mathbb{L}}\mathrm{cofib}(f))\rightarrow\pi_{n-1}\mathbb{L}_{B\big\slash A}$$
    is an isomorphism. Since $\mathrm{cofib}(f)$ is $n-1$-connective and $\pi_{0}(A)\cong\pi_{0}(B)$, we have $\pi_{n-1}(\mathrm{cofib}(f))\cong\pi_{n-1}(B\otimes_{A}^{\mathbb{L}}\mathrm{cofib}(f))$. Thus $\pi_{n-1}(\mathrm{cofib}(f))\cong\pi_{n-1}(\mathbb{L}_{B\big\slash A})$ is also an isomorphism, and $\pi_{n-1}(\mathbb{L}_{B\big\slash A})$ is nonzero.
\end{proof}

\begin{cor}[\cite{HA} Corollary 7.4.3.3]
Let $A\in\mathbf{Alg_{D}}^{cn}(\mathbf{C})$. Then $\mathbb{L}_{A}$ is connective.
\end{cor}

\begin{cor}[\cite{HA} Corollary 7.4.3.4]\label{cor:pi0defequiv}
Let $f:A\rightarrow B$ be a map in  $\mathbf{Alg}^{cn}_{\mathbf{D}}(\mathbf{C})$. Suppose that for $n>0$ we have $a(n)>0$. Then $f$ is an equivalence if and only if
\begin{enumerate}
\item
$f$ induces an isomorphism $\pi_{0}(A)\rightarrow\pi_{0}(B)$.
\item
$\mathbb{L}_{B\big\slash A}\cong 0$.
\end{enumerate}
\end{cor}

\begin{cor}[\cite{HA} Corollary 7.4.3.5]\label{cor:pi0L}
Let $f:A\rightarrow B$ be a map in $\mathbf{Alg}^{cn}_{\mathbf{D}}(\mathbf{C})$ such that $\mathrm{cofib}(f)$ is $n$-connective for $n\ge 0$. Then the induced map $\mathbb{L}_{f}:\mathbb{L}_{A}\rightarrow\mathbb{L}_{B}$ has an $n$-connective cofibre. In particular the canonical map $\pi_{0}(\mathbb{L}_{A})\rightarrow\pi_{0}(\mathbb{L}_{\pi_{0}(A)})$ is an isomorphism.
\end{cor}

\begin{proof}
    $\mathbb{L}_{f}$ is given as the composition
    $$\mathbb{L}_{A}\rightarrow B\otimes_{A}^{\mathbb{L}}\mathbb{L}_{A}\rightarrow\mathbb{L}_{B}$$
    Then $\mathrm{cofib}(g)\cong\mathrm{cofib}(f)\otimes_{A}^{\mathbb{L}}\mathbb{L}_{A}$.  Since $\mathbb{L}_{A}$ is connective and $\mathrm{cofib}(f)$ is $n$-connective, $\mathrm{cofib}(g)$ is $n$-connective. We get a fibre sequence
    $$B\otimes_{A}^{\mathbb{L}}\mathrm{cofib}(f)\rightarrow\mathbb{L}_{B\big\slash A}\rightarrow\mathrm{cofib}(\epsilon_{f})$$
    We clearly have that $B\otimes_{A}^{\mathbb{L}}\mathrm{cofib}(f)$ is $n$-connective, and $\mathrm{cofib}(\epsilon_{f})$ is $n+a(n)+1\ge n$-connective.
\end{proof}

%
%

We can now give verifiable conditions under which maps of algebras arise as square-zero extensions. In a stable $(\infty,1)$-category with $t$-structure $(\mathbf{C},\mathbf{C}_{\ge0},\mathbf{C}_{\le0})$ we say that an object is \textit{discrete} if it is contained in $\mathbf{C}^{\heart}$.

\begin{defn}[\cite{lurieDAG4} Definition 3.2.2]
Let $n\ge 0$ be an integer. A map $f:\tilde{A}\rightarrow A$ in $\mathbf{Alg_{D}}$ is said to be an $a(n)$-\textit{ small extension} if the following conditions are satisfied.
\begin{enumerate}
\item
$A$ is connective.
\item
$\pi_{i}\mathrm{Ker}(f)$ vanish unless 
 $n\le i\le n+a(n)$
\item
The map
$$\mathrm{fib}(f)\otimes_{\tilde{A}}\mathrm{fib}(f)\rightarrow\tilde{A}\otimes_{\tilde{A}}\mathrm{fib}(f)\cong\mathrm{fib}(f)$$
is zero in the homotopy category. 
\end{enumerate}
Let $\mathbf{Fun}_{a(n)-sm}(\Delta^{1},\mathbf{Alg_{D}})$ denote the full subcategory of the arrow category $\mathbf{Fun}(\Delta^{1},\mathbf{Alg_{D}})$ consisting of $a(n)$-small extensions. 
\end{defn}

\begin{rem}[c.f. \cite{HA} Remark 3.2.5]\label{rem:nsmallcond}
    Let $f:\tilde{A}\rightarrow A$ such that $A$ is connective and $\pi_{i}\mathrm{Ker}(f)$ vanish unless 
$n\le i\le n+a(n)$. Then $\mathrm{fib}(f)\otimes_{\tilde{A}}^{\mathbb{L}}\mathrm{fib}(f)$ is $(2n)$-connective. Moreover $\mathrm{fib}(f)$ is $n+a(n)\le 2n$ truncated. Thus the multiplication map 
$$\mathrm{fib}(f)\otimes_{\tilde{A}}^{\mathbb{L}}\mathrm{fib}(f)\rightarrow\mathrm{fib}(f)$$
is determined by 
$$\pi_{n}(\mathrm{fib}(f))\otimes_{\pi_{0}(\tilde{A})}^{\mathbb{L}}\pi_{n}(\mathrm{fib}(f))\rightarrow\pi_{2n}\mathrm{fib}(f)$$
and even then this can only possibly be non-trivial if $a(n)=n$.
\end{rem}

\begin{example}[c.f. \cite{HA} Example 7.4.1.21]\label{ex:0small}
Let $f:\tilde{A}\rightarrow A$ be a map between objects in $\mathbf{Alg_{D}}(\mathbf{C})$ with $A\in\mathbf{C}^{\heart}$. Then $f$ is an $a(0)$-small extension if and only if
\begin{enumerate}
\item
$\tilde{A}$ is discrete.
\item
the induced map $\pi_{0}(\tilde{A})\rightarrow\pi_{0}(A)$ is an epimorphism.
\item
if $I=\mathrm{Ker}(\pi_{0}(\tilde{A})\rightarrow\pi_{0}(A))$ then the map $I\otimes I\rightarrow\pi_{0}(\tilde{A})$ is the zero map.
\end{enumerate}
(Note that $0\le a(0)\le 20=0$ so $a(0)=0$)
.\end{example}

The following is clear.

\begin{prop}
A map $f:\tilde{A}\rightarrow A$ in $\mathbf{Alg_{D}}(\mathbf{C})$ $a(n)$-small if and only if $\Theta(f):\Theta(\tilde{A})\rightarrow\Theta(A)$ $a(n)$-small, in $\mathbf{Comm}(\mathbf{C})$. 
\end{prop}

Note that usually $n$-small extensions in $\mathbf{Comm(C)}$ in the sense of \cite{HA} have quite different bounds.

We may identify an object of $\underline{\mathbf{Der}}(\mathbf{Alg_{\mathbf{D}}})$ with a derivation $d:A\rightarrow A\oplus M[-1]$. 

\begin{defn}
A derivation $d:A\rightarrow A\oplus M$ is said to be $a(n)$-\textit{small} if
$M\in\mathbf{C}_{\ge n+1}\cap\mathbf{C}_{\le n+a(n)+1}$.
Let $\underline{\mathbf{Der}}_{a(n)-sm}(\mathbf{Alg_{D}}(\mathbf{C}))$ denote the full subcategory of $\underline{\mathbf{Der}}(\mathbf{Alg_{D}}(\mathbf{C}))$ consisting of $a(n)$-small derivations.
\end{defn}

Again a derivation $d:A\rightarrow A\oplus M$ is $a(n)$-small if and only if the derivation of commutative monoids $d^{\Theta}:\Theta(A)\rightarrow\Theta(A)\oplus M$ is $a(n)$-small. The following can be proven exactly as in \cite{lurieDAG4} Theorem 3.2.7, using our version of \cite{HA} Theorem 7.4.3.12, Theorem \ref{thm:connectivcot}.

\begin{thm}\label{thm:n-smn-der}
Let $0\le a(n)\le 2n$ be an increasing sequence of natural numbers such that 
$$\mathbf{D}^{k}(M)\in\mathbf{C}_{\ge n+a(n)}$$
whenever $M\in\mathbf{C}_{\ge n}$ and $k\ge 2$. Let $0\le\overline{a}(n)\le 2n$ be an increasing sequence of natural numbers such that $\overline{a}(n)\le a(n)$ and $\overline{a}(n)<a(n+1)$.
The functor
$$\Phi:\underline{\mathbf{Der}}(\mathbf{Alg_{D}}(\mathbf{C}))\rightarrow\mathbf{Fun}(N(\Delta^{1}),\mathbf{Alg_{D}}(\mathbf{C}))$$
induces an equivalence
$$\Phi:\underline{\mathbf{Der}}_{\overline{a}(n)-sm}(\mathbf{Alg_{D}}(\mathbf{C}))\rightarrow\mathbf{Fun}_{\overline{a}(n)-sm}(N(\Delta^{1}),\mathbf{Alg_{D}}(\mathbf{C}))$$
\end{thm}

\begin{proof}
First we show that the functor $\Phi$ does indeed send $\underline{\mathbf{Der}}_{a(n)-sm}(\mathbf{Alg_{D}}(\mathbf{C}))$ to $\mathbf{Fun}_{a(n)-sm}(N(\Delta^{1}),\mathbf{Alg_{D}}(\mathbf{C}))$. Indeed let $\eta:A\rightarrow M$ be an object in $\underline{\mathbf{Der}}_{a(n)-sm}(\mathbf{Alg_{D}}(\mathbf{C}))$. Consider the square zero extension $f:A^{\eta}\rightarrow A$. The fibre is $M[-1]$, so clearly $\pi_{i}(M[-1])$ vanishes for $n\le i\le n+\overline{a}(n)$. The multiplication map 
$$A^{\eta}\otimes^{\mathbb{L}} M[-1]\rightarrow M[-1]$$
factors through $A\otimes M[-1]\cong A\otimes^{\mathbb{L}}\mathrm{fib}(f)$, and so is trivial when restricted to $\mathrm{fib}(f)\otimes^{\mathbb{L}}\mathrm{fib}(f)$. That $f$ is an $\overline{a}(n)$-small extension now follows from Remark \ref{rem:nsmallcond}.

    The inverse functor $G$ to $\phi$ may be constructed as in \cite{lurieDAG4} Theorem 3.2.7 as follows. Let $\tilde{A}\rightarrow A$ be an $\overline{a}(n)$-small extension. Consider the composite map
    $$\mathbb{L}_{A}\rightarrow\mathbb{L}_{A\big\slash\tilde{A}}\rightarrow\tau_{\le n+\overline{a}(n)+1}\mathbb{L}_{A\big\slash\tilde{A}}$$
    This classifies a derivation $d:A\rightarrow\tau_{\le n+\overline{a}(n)+1}\mathbb{L}_{A\big\slash\tilde{A}}$. We prove that this derivation is $\overline{a}(n)$-small. This functor is left adjoint to $\Phi$. We denote the unit by $u$. We have a commutative diagram
    \begin{displaymath}
        \xymatrix{
        \mathrm{fib}(f)\ar[d]^{g}\ar[r] & \tilde{A}\ar[d]^{u_{f}}\ar[r]^{f} & A\ar[d]^{\mathrm{Id}}\\
        \tau_{\le n+\overline{a}(n)}(\mathbb{L}_{A\big\slash\tilde{A}}[-1])\ar[r] & A^{\eta}\ar[r] & A
        }
    \end{displaymath}
    We need to prove that $g$ is an equivalence. This proves that $\tilde{A}$ is a square-zero extension by an $\overline{a}(n)$-small derivation. Now $\pi_{i}\mathrm{fib}(f)\cong0$ for $i>n+\overline{a}(n)$, so it suffices to prove that 
    $$\pi_{i}\mathrm{fib}(f)\rightarrow\pi_{i}\tau_{\le n+\overline{a}(n)}(\mathbb{L}_{A\big\slash\tilde{A}}[-1])$$
    is an equivalence for all $i\le n+\overline{a}(n)$. $g$ factors as the composition 
    \begin{displaymath}
        \xymatrix{
\mathrm{fib}(f)\ar[r]^{g'}&\mathrm{ker}(f)\otimes^{\mathbb{L}}_{\tilde{A}}A\ar[r]^{\epsilon_{f}[-1]}& \mathbb{L}_{A\big\slash\tilde{A}}[-1]\ar[r]^{g'''} & \tau_{\le n+\overline{a}(n)}(\mathbb{L}_{A\big\slash\tilde{A}}[-1])
        }
    \end{displaymath}
    We therefore need to show that each of $g'$, $g''=\epsilon_{f}[-1]$, and $g'''$ induce isomorphisms on $\pi_{i}$ for $i\le 
    n+\overline{a}(n)$. For $g'''$ this is trivial. Now for $g'$ we consider the long exact sequence
    \begin{displaymath}
    \xymatrix{
    \pi_{i}(\mathrm{fib}(f)\otimes_{\tilde{A}}^{\mathbb{L}}\mathrm{fib}(f))\ar[r] &\pi_{i}(\mathrm{fib}(f))\ar[r] &\pi_{i}(\mathrm{fib}(f)\otimes_{\tilde{A}}^{\mathbb{L}}A)\ar[r] &  \pi_{i-1}(\mathrm{fib}(f)\otimes_{\tilde{A}}^{\mathbb{L}}\mathrm{fib}(f))
    }
    \end{displaymath}
Since $\mathrm{fib}(f)$ is $n$-connective, $\mathrm{fib}(f)\otimes_{\tilde{A}}^{\mathbb{L}}\mathrm{fib}(f)$ is $2n\ge n+\overline{a}(n)$-connective. Thus for $i<2n$ the map $\pi_{i}(g')$ is an isomorphism. For $i=2n$ we have that $\pi_{i}(g)$ is surjective with kernel generated by the image of 
$$\pi_{n}(\mathrm{fib}(f))\otimes_{\pi_{0}(\tilde{A})}\pi_{n}(\mathrm{fib}(f))\rightarrow\pi_{2n}(\mathrm{fib}(f))$$
which we have assumed to be zero.

    Now the map $f$ has $n$-connective fibre, and so has $n+1$-connective cofibre. Thus $\mathrm{fib}(\epsilon_{f})\in\mathbf{C}_{\ge n+1+a(n+1)}$. So $\mathrm{fib}(\epsilon_{f}[-1])\in \mathbf{C}_{\ge n+a(n+1)}$ and hence is in $\mathbf{C}_{\ge n+\overline{a}(n)+1}$. Thus the map $\epsilon_{f}[-1]$ induces an equivalence on $\pi_{i}$ for $i\le n+\overline{a}(n)$, as required.

    It remains to prove that $G$ is conservative. But this is straightforward, as in \cite{lurieDAG4} Theorem 3.2.7.

\end{proof}

%
%

%
%
%
%
%
The following immediate corollary can be deduced by taking $a(0)=0$ and $a(n)=1$ for $n>0$, and $\overline{a}(n)=0$ for all $n$. It will be crucial for obstruction theory.

\begin{thm}[\cite{HA} Corollary 7.4.1.28]\label{thm:squarezeropostnikov}
Let the $t$-structure be Koszul. Let $A\in\mathbf{Alg}_{\mathbf{D}}^{cn}(\mathbf{C})$. Then for each $n\ge 1$ there is a derivation
$$d^{A}_{n}\in\pi_{0}\mathbf{Der}(\tau_{\le n-1}A,\pi_{n}(A)[n+1])$$
such that $\tau_{\le n}A$ is the square-zero extension of $\tau_{\le n-1}A$ along $d^{A}_{n}$. 
\end{thm}

By Proposition \ref{prop:lifting} and the existence of Postnikov towers, we have the following.

\begin{cor}\label{cor:liftingetalemaps}
Suppose our algebraic context is Koszul. Let $f:A\rightarrow B$ and $g:A\rightarrow C$ be maps of algebras in $\mathbf{Alg_{D}}^{cn}$ with $f:A\rightarrow B$ being formally \'{e}tale. Then
\begin{enumerate}
\item
any map $\pi_{0}(B)\rightarrow\pi_{0}(C)$ lifts, uniquely up to a contractible choice, to a map of $A$-modules $B\rightarrow C$.
\item
if $g:A\rightarrow C$ is also formally \'{e}tale, then any isomorphism $\pi_{0}(B)\rightarrow\pi_{0}(C)$ lifts, uniquely up to a contractible choice, to an isomorphism of $A$-modules $B\rightarrow C$.
\end{enumerate}
\end{cor}

For posterity we include the following definition.

\begin{defn}
\begin{enumerate}
\item
A \textit{Postnikov spectral algebraic pre-context} is a tuple
$$(\mathbf{C},\mathbf{C}_{\ge0},\mathbf{C}_{\le0},\mathbf{C}^{0},\mathbf{D},\theta,\mathbf{G}_{0},\mathbf{S})$$
where
\begin{enumerate}
\item
$(\mathbf{C},\mathbf{D},\theta,\mathbf{G}_{0},\mathbf{S})$ is an $(\infty,1)$-algebra context.
\item
$(\mathbf{C}_{\ge0},\mathbf{C}_{\le0})$ is a compatible $t$-structure  on $(\mathbf{C},\mathbf{D},\theta)$. 
\end{enumerate}
\item
A \textit{connective Postnikov algebraic pre-context} is a tuple
$$(\mathbf{C},\mathbf{C}_{\ge0},\mathbf{C}_{\le0},\mathbf{C}^{0},\mathbf{D},\theta,\mathbf{G}_{0},\mathbf{S})$$
where
\begin{enumerate}
\item
$(\mathbf{C},\mathbf{C}_{\ge0},\mathbf{C}_{\le0},\mathbf{C}^{0},\mathbf{D},\theta)$ is a Postnikov algebraic pre-context.
\item
$(\mathbf{C},\mathbf{D},\theta,\mathbf{G}_{0},\mathbf{S})$ is an $(\infty,1)$-algebra context.
\item $\mathbf{G}_{0}\subseteq\mathbf{C}_{\ge0}$, $\mathbf{S}\subseteq\mathbf{Alg}_{D}(\mathbf{C}_{\ge0})$.
\end{enumerate}
\end{enumerate}
A Postnikov algebraic pre-context/ connective Postnikov algebraic pre-context is said to be a \textit{Koszul algebraic pre-context/ connective Koszul algebraic pre-context} if the $t$-structure on $(\mathbf{C},\mathbf{D},\theta,\mathbf{G}_{0},\mathbf{S})$ is Koszul compatible.
\end{defn}
The point of the definition of a connective Postnikov algebraic context, is simply that we wish to work connectively in some larger stable category, and we want to remember this stable category as part of the data. Note that $\mathbf{C}\cong\mathbf{Stab}(\mathbf{C}_{\ge0})$. When $\mathbf{G}_{0}=\mathbf{C}_{\ge0}$ and $\mathbf{S}=\mathbf{Alg_{D}}$, then we just say that $(\mathbf{C},\mathbf{C}_{\ge0},\mathbf{C}_{\le0},\mathbf{C}^{0},\mathbf{D},\theta)$ is a Postnikov algebraic context (resp. a connective Postnikov algebraic context). 

Usually our Postnikov algebraic contexts will be of the form

$$(\mathbf{C},\mathbf{C}_{\ge0},\mathbf{C}_{\le0},\mathbf{C}^{0},\mathbf{D},\theta,\mathbf{C},\mathbf{Alg_{D}}(\mathbf{C}))$$
or
$$(\mathbf{C},\mathbf{C}_{\ge0},\mathbf{C}_{\le0},\mathbf{C}^{0},\mathbf{D},\theta,\mathbf{C}_{\ge0},\mathbf{Alg_{D}}(\mathbf{C})^{\heart})$$

and our connective Postnikov algebraic contexts will be of the form

$$(\mathbf{C},\mathbf{C}_{\ge0},\mathbf{C}_{\le0},\mathbf{C}^{0},\mathbf{D},\theta,\mathbf{C}_{\ge0},\mathbf{Alg_{D}}(\mathbf{C})^{cn})$$

\subsubsection{Cotangent Complexes in Subcategories}

It will be helpful for us to know when the cotangent complex of a certain class of algebras lies in a certain class of modules. This will be particularly useful when we wnat the cotangent complex to be coherent.

\begin{lem}\label{lem:cotangentsubcat}
Let $\mathbf{H}\subseteq\mathbf{Comm}^{cn}(\mathbf{C})$ be a full subcategory such that whenever $A\in\mathbf{H}$, $\pi_{0}(A)\in\mathbf{H}$. For each $A\in\mathbf{H}\cap\mathrm{Comm}(\mathbf{C}^{\heart})$ fix a full subcategory $\mathrm{F}_{A}$. For $A\in\mathbf{H}$, let $\mathbf{F}_{A}$ denote the full subcategory of ${}_{A}\mathbf{Mod}$ consisting of objects $X$ such that $\pi_{m}(X)\in\mathrm{F}_{\pi_{0}(A)}$ for all $m\in\mathbb{Z}$. Suppose that 
\begin{enumerate}
    \item 
    $A\in\mathbf{H}$ if and only if $\pi_{0}(A)\in\mathbf{H}$ and each $\pi_{n}(A)\in\mathbf{F}_{\pi_{0}(A)}$
    \item 
    whenever $A\rightarrow B$ is a map in $\mathbf{H}$, and $X\in\mathbf{F}_{A}$ is connective, then $B\otimes_{A}^{\mathbb{L}}M$ is in $\mathbf{F}_{B}$.
    \item 
    $\mathrm{F}_{A}$ is thick, and closed under finite limits and colimits, and tensor products.
\end{enumerate}
Let $f:A\rightarrow B$ be a map in $\mathbf{D}$ such that $\mathrm{cofib}(f)$ is $1$-connective. Then $\mathbb{L}_{B\big\slash A}\in\mathbf{F}_{B}$
\end{lem}

\begin{proof}
        Write $A_{0}=A$. Let $M$ be the fibre of 
    $A\rightarrow B$
    Then $M$ is $0$-connective. Note moreover that is is in $\mathbf{F}_{A}$. Consider the pushout
    $$A_{1}\defeq A_{0}\otimes_{\mathbf{D}_{A_{0}}(M)}^{\mathbb{L}} A_{0}\cong\mathbf{D}_{A_{0}}(M[1])$$
The map $A_{1}\rightarrow B$ has $2$-connective cofibre. 
    
    $A_{1}$ is in $\mathbf{H}$. Indeed since $M[1]$ is $1$-connective we have $\pi_{0}(A_{1})\cong A_{0}$.
    
Now $\mathbf{Sym}_{A_{0}}^{m}(M[1])$ is a quotient of the tensor algebra $\mathrm{T}_{A_{0}}^{m}(M[1])$. Moreover each $\pi_{n}(\mathrm{T}_{A_{0}}^{m}(M[1]))\cong 0$ for $m>n$ by the connectivity of $M$, so $\pi_{n}(\mathrm{T}(M[1]))$ is in $\mathrm{F}_{A}$. Again since each finite tensor is in $\mathbf{F}_{A_{0}}$ and this implies that $A_{1}$ is in $\mathbf{H}$. Finally the map $A_{1}\rightarrow A$ now has $2$-connective cofibre. We have $\mathbb{L}_{A_{1}\big\slash A_{0}}\cong A_{1}\otimes_{A_{0}}^{\mathbb{L}}M[1]$ which is in $\mathbf{F}_{A_{1}}$, since $M[1]$ is in $\mathbf{F}_{A_{0}}$. Continuing in this way we produce maps $A_{n}\rightarrow A$ with $n$-connective cofibre such that $\mathbb{L}_{A_{n}\big\slash A_{0}}$ is in $\mathbf{F}_{A_{n}}$. Taking the colimit gives that $\mathbb{L}_{B\big\slash A}$ is in $\mathbf{F}_{B}$. as required.
\end{proof}

\subsubsection{\'{E}tale Maps, and Monomorphisms}

Let $(\mathbf{C},\mathbf{D},\theta,\mathbf{G}_{0},\mathbf{S})$ be a stable $(\infty,1)$-algebra context

\begin{defn}
A map $f:A\rightarrow B$ is said to be 
\begin{enumerate}
\item
\textit{Formally unramified} if $\mathbb{L}_{B\big\slash A}\cong 0$.
\item
\textit{Formally \'{e}tale} if the map $B\otimes^{\mathbb{L}}_{A}\mathbb{L}_{A}\rightarrow\mathbb{L}_{B}$
is an equivalence.
\item
\textit{A homotopy epimorphism} if the map
$$B\otimes^{\mathbb{L}}_{A}B\rightarrow B$$
is an equivalence.
\end{enumerate}
\end{defn}

Let us list some basic properties of the cotangent complex, as in \cite{toen2008homotopical}  Proposition 1.2.1.6.

\begin{prop}
\begin{enumerate}
    \item
    Let 
$$A\rightarrow B\rightarrow C$$
be a sequence of maps in $\mathrm{Alg}_{\mathbf{D}}(\mathbf{C})$. We get a fibre-cofibre sequence in ${}_{C}\mathrm{Mod}$:
$$\mathbb{L}_{B\big\slash A}\otimes^{\mathbb{L}}_{B}C\rightarrow\mathbb{L}_{C\big\slash A}\rightarrow\mathbb{L}_{C\big\slash B}$$
\item Let 
$B'\cong A'\otimes_{A}^{\mathbb{L}}B$. Then the natural morphism
$$\mathbb{L}_{B\big\slash A}\otimes_{B}^{\mathbb{L}}B'\rightarrow\mathbb{L}_{B'\big\slash A'}$$
is an equivalence.
\item  Let 
$B'\cong A'\otimes_{A}^{\mathbb{L}}B$. Then there is a fibre-cofiber sequence
$$\mathbb{L}_{A}\otimes_{A}^{\mathbb{L}}B'\rightarrow\mathbb{L}_{A'}\otimes_{A'}^{\mathbb{L}}B'\oplus\mathbb{L}_{B}\otimes_{B}^{\mathbb{L}}B'\rightarrow\mathbb{L}_{B'}$$
\end{enumerate}

\end{prop}
This immediately implies the following.

\begin{prop}
A map $f:A\rightarrow B$ is formally \'{e}tale if and only if it is formally unramified.
\end{prop}


For $A\rightarrow B$ a map consider the sequence
$$A\rightarrow B\otimes^{\mathbb{L}}_{A}B\rightarrow B$$
where the first map is inclusion in the first factor at the second is multiplication, so that the composite is the identity. We get a fibre-cofibre sequence
$$\mathbb{L}_{B\otimes^{\mathbb{L}}_{A}B\big\slash B}\otimes_{B\otimes^{\mathbb{L}}_{A}B}^{\mathbb{L}}B\rightarrow 0\rightarrow\mathbb{L}_{B\big\slash B\otimes^{\mathbb{L}}_{A}B}$$
On the other hand $B\rightarrow B\otimes^{\mathbb{L}}_{A}B$ is the base change of $A\rightarrow B$. Thus we have
$$\mathbb{L}_{B\otimes^{\mathbb{L}}_{A}B\big\slash B}\otimes_{B\otimes^{\mathbb{L}}_{A}B}^{\mathbb{L}}B\cong\mathbb{L}_{B\big\slash A}$$
and hence
$$\mathbb{L}_{B\big\slash A}\cong\mathbb{L}_{B\big\slash B\otimes^{\mathbb{L}}_{R}B}[-1]$$
This proves the following.
\begin{lem}\label{lem:formallyetiffdiag}
A map $f:A\rightarrow B$ is formally \'{e}tale if and only if $B\otimes^{\mathbb{L}}_{A}B\rightarrow B$ is formally \'{e}tale. In particular homotopy epimorphisms are formally \'{e}tale.
\end{lem}

We finish with one more class of maps.

\begin{defn}[\cite{toen2008homotopical} Definition 1.2.8.1 (1)]
A morphism $f:A\rightarrow B$ in $\mathbf{Alg_{D}}(\mathbf{C})$ is said to be \textit{formally i-smooth} if for any map $f:A\rightarrow C$ in $\mathbf{Alg_{D}}(\mathbf{C})$ with $C$ $\mathbf{G}_{0}$-connective, any $M$ which is $\mathbf{G}_{1}$-connective, and any $d\in\pi_{0}(\mathbf{Der}_{A}(C,M)$ the map
$$\pi_{0}(\mathbf{Map}_{\mathbf{Alg_{D}}({}_{A}\mathbf{Mod})}(B,C\oplus_{d}\Omega M))\rightarrow\pi_{0}(\mathbf{Map}_{\mathbf{Alg_{D}}({}_{A}\mathbf{Mod})}(B,C))$$
is surjective. 
\end{defn}

The following can be proven exactly as in \cite{toen2008homotopical} 1.2.8.3.

\begin{lem}
A map $f:A\rightarrow B$ is formally $i$-smooth if and only if for any morphism $B\rightarrow C$ with $C\in\mathbf{S}$, and any $C$-module $M\in{}_{C}\mathbf{Mod}_{\ge1}$ the natural map 
$$\pi_{0}\mathbf{Map}(\mathbb{L}_{C\big\slash A},M)\rightarrow \pi_{0}\mathbf{Map}(\mathbb{L}_{B\big\slash A},M)$$
is zero. 
\end{lem}

In particular we have the following.

\begin{cor}
Formally \'{e}tale maps are formally $i$-smooth.
\end{cor}

\subsection{Transformations of $(\infty,1)$-Algebraic Contexts}

\begin{defn}\label{defn:transformation}
A \textit{transformation} from an $(\infty,1)$-algebraic context $(\mathbf{C},\mathbf{D},\theta,\mathbf{G}_{0},\mathbf{S})$ to an $(\infty,1)$-algebraic context $(\mathbf{C}',\mathbf{D}',\theta,\mathbf{G}'_{0},\mathbf{S}')$ is a tuple
$$(\mathbf{F},\eta,\delta)$$
where 
\begin{enumerate}
\item
$\mathbf{F}:\mathbf{C}\rightarrow\mathbf{C}'$ is a left adjoint strongly symmetric monoidal functor with right adjoint $\mathbf{R}$.
\item
$\eta:\mathbf{F}\circ\mathbf{D}\rightarrow\mathbf{D}'\circ\mathbf{F}$
is a natural transformation of functors.
\item
considering the diagram
\begin{displaymath}
\xymatrix{
\mathbf{F}\circ\mathbf{Comm}_{\mathbf{C}}(-)\ar[d]^{\mathbf{F}(\theta)}\ar[r]^{\sigma_{\mathbf{F}}} & \mathbf{Comm}_{\mathbf{C}'}(-)\circ\mathbf{F}\ar[d]^{\theta'_{\mathbf{F}}}\\
\mathbf{F}\circ\mathbf{D}\ar[r]^{\eta} & \mathbf{D}'\circ\mathbf{F}
}
\end{displaymath}
$\delta:\theta'_{\mathbf{F}}\circ\sigma_{\mathbf{F}}\Rightarrow\eta\circ\mathbf{F}(\theta)$ is a natural equivalence, where $\sigma_{F}:\mathbf{F}\circ\mathbf{Comm}_{\mathbf{C}}(-)\rightarrow\mathbf{Comm}_{\mathbf{C}'}(-)\circ\mathbf{F}$ is the natural isomorphim arising from the strongly symmetirc monoidal structure on $\mathbf{F}$,
\end{enumerate}
such that
\begin{enumerate}
\item
$\mathbf{F}(\mathbf{G}_{0})\subset\mathbf{G}'_{0}$.
\item
$\mathbf{R}_{\mathbf{D}}(\mathbf{S}')\subseteq\mathbf{S}$, where $\mathbf{R}_{\mathbf{D}}$ is the right adjoint to $\mathbf{F}:\mathbf{Alg_{D}}(\mathbf{C})\rightarrow\mathbf{Alg_{D'}}(\mathbf{C}')$
\end{enumerate}
\end{defn}

\subsubsection{Preservation of Properties of Maps Under Transformations}

Fix a transformation $(\mathbf{F},\eta,\delta)$ from an $(\infty,1)$-algebraic context $(\mathbf{C},\mathbf{D},\theta,\mathbf{G}_{0},\mathbf{S})$ to an $(\infty,1)$-algebraic context $(\mathbf{C}',\mathbf{D}',\theta,\mathbf{G}_{0},\mathbf{S}')$.

\begin{lem}
Let $f:A\rightarrow B$ be a map in $\mathbf{Alg_{D}}(\mathbf{C})$.  If $f$ is a homotopy epimorphism then so is $\mathbf{F}(f)$.
\end{lem}

\begin{proof}
We have 
$$\mathbf{F}(B)\coprod^{\mathbb{L}}_{A}\mathbf{F}(B)\cong\mathbf{F}(B\coprod^{\mathbb{L}}_{A}B)\cong\mathbf{F}(B)$$
\end{proof}

\subsection{Perfect, Dualisable, Compact, and Nuclear Objects}

Throughout the paper we will have cause to consider object which are `small' in various ways. We define perfect and dualisable objects, and  compare them with compact objects. Fix a stable symmetric monoidal $(\infty,1)$-category $\mathbf{C}$ with tensor product $\otimes$ and unit $\mathbb{I}$.

\begin{defn}
An object $M$ of $\mathbf{C}$ is said to be \textit{perfect} if it is a retract of a finite colimit of objects of the form $\mathbb{I}^{\coprod n}$ for $n\in\mathbb{N}_{0}$ where $\mathbb{I}$ is the tensor unit.
\end{defn}

The full subcategory of $\mathbf{C}$ consisting of perfect objects is denoted $\mathbf{Perf}(\mathbf{C})$. 

For $\kappa$ a cardinal denote by $\mathbf{Cpt}_{\kappa}(\mathbf{C})$ the category of $\kappa$-compact objects in $\mathbf{C}$. For $\kappa=\aleph_{0}$ we just write $\mathbf{Cpct}_{\aleph_{0}}(\mathbf{C})=\mathbf{Cpct}(\mathbf{C})$

$$\mathbf{Perf}(\mathbf{C})\subseteq\mathbf{Cpct}\mathbf{C})$$

In our examples of interest, namely complete bornological modules over a Banach ring, this \textit{will not} be an equality.

\begin{defn}
Let $\mathbf{C}$ be a closed symmetric monoidal $(\infty,1)$-category. An object $m\in\mathbf{M}$ is said to be \textit{strongly dualisable} if the map
$$m^{\vee}\otimes^{\mathbb{L}} n\rightarrow\underline{\mathbf{Map}}(m,n)$$
is an equivalence for any $n\in\mathbf{M}$.
\end{defn}

Denote by $\mathbf{Dls}(\mathbf{C})\subset\mathbf{C}$ the full subcategory consisting of strongly dualisable objects. Since $\mathbf{C}$ is stable, $\mathbf{Dls}(\mathbf{C})$ is closed under finite colimits and retracts. Moreover it evidently contains $\mathbb{I}$. Thus $\mathbf{Perf}(\mathbf{C})\subseteq\mathbf{Dls}(\mathbf{C})$. Note moreover that since tensoring commutes with colimits, dualisable objects are $\kappa$-compact if $\mathbb{I}$ is so. In fact by \cite{toen2008homotopical} Corollary 1.2.3.8 in this case they coincide with the perfect objects.

Note that the perfect, dualisable, and $\kappa$-compact objects are stable by base-change.
Precisely if $A\in\mathbf{Comm}(\mathbf{C})$, and $m$ is in $\mathbf{Perf}({}_{A}\mathbf{Mod})$ (resp. $\mathbf{Cpct}_{\kappa}({}_{A}\mathbf{Mod})$/ $\mathbf{Dls}_{\kappa}({}_{A}\mathbf{Mod})$) then $B\otimes^{\mathbb{L}}_{A}m$ is in $\mathbf{Perf}({}_{B}\mathbf{Mod})$ (resp. $\mathbf{Cpct}_{\kappa}({}_{B}\mathbf{Mod})$/ $\mathbf{Dls}_{\kappa}({}_{B}\mathbf{Mod})$). This is also the case for $\mathbf{Dls}(\mathbf{C})$.  Indeed for perfect and $\kappa$-compact objects this is clear. 

As for dualisable objects, let $A,B\in\mathbf{Comm}(\mathbf{C})$, and let $m\in{}_{A}\mathbf{Mod}(\mathbf{C})$ be dualsiable. Let $A\rightarrow B$ be a map, and let $n$ be any object of ${}_{B}\mathbf{Mod}(\mathbf{C})$. We have
\begin{align*} 
\mathbf{Map}_{B}(B\otimes^{\mathbb{L}}_{A} m,n)&\cong\mathbf{Map}_{A}(m,n)\\
&\cong\mathbf{Map}_{A}(m,A)\otimes^{\mathbb{L}} n
\end{align*}
In particular for $n=B$ we have
$$\mathbf{Map}_{B}(B\otimes^{\mathbb{L}}_{A} m,B)\cong\mathbf{Map}_{A}(m,A)\otimes^{\mathbb{L}}_{A} B$$
Hence
\begin{align*}
\mathbf{Map}_{B}(B\otimes_{A}^{\mathbb{L}} m,B)\otimes^{\mathbb{L}}_{B}n &\cong\mathbf{Map}_{A}(m,A)\otimes_{A} n\\
&\cong\mathbf{Map}_{A}(m,n)\\
&\cong\mathbf{Map}_{B}(B\otimes^{\mathbb{L}}_{A}m,n)
\end{align*}
Thus $B\otimes^{\mathbb{L}}_{A}m$ is perfect in ${}_{B}\mathbf{Mod}(\mathbf{C})$.

\subsubsection{Nuclear Objects}

Let $\mathbf{C}$ be a compactly generated closed symmetric monoidal $(\infty,1)$-category. Recall \cite{scholzecomplexlectures}, \cite{ben2020fr} the definition of a nuclear objects and maps in $\mathbf{C}$

\begin{defn}
\begin{enumerate}
\item 
A map $f:X\rightarrow Y$ in $\mathbf{C}$ is said to be \textit{nuclear} or \textit{trace-class} if it is in the image of 
$$\pi_{0}\mathbf{Map}(\mathbb{I},X^{\vee}\otimes^{\mathbb{L}}Y)\rightarrow\pi_{0}\mathbf{Map}(X,Y)$$
    \item 
    An object $X$ of $\mathbf{C}$ is said to be \textit{nuclear} if for any $c\in\mathbf{C}$ compact, the natural map
$$\pi_{0}\mathbf{Map}(\mathbb{I},c^{\vee}\otimes^{\mathbb{L}}X)\rightarrow\pi_{0}\mathbf{Map}(c,X)$$
is an isomorphism.
\item     
An object $X$ of $\mathbf{C}$ is said to be \textit{strongly nuclear} if for any $c\in\mathbf{C}$ compact, the natural map
$$c^{\vee}\otimes^{\mathbb{L}}X\rightarrow\underline{\mathbf{Map}}(c,X)$$
is an equivalence.
\end{enumerate}
\end{defn}

Note that dualisable objects are in particular strongly nuclear.  For $A\in\mathbf{Comm}(\mathbf{C})$, write $\mathbf{M}^{\vee_{A}}\defeq\underline{\mathbf{Map}}_{{}_{A}\mathbf{Mod}(\mathbf{C})}(M,A)$.

In \cite{scholzecomplexlectures} Chapter VIII Clausen and Scholze prove the following.

\begin{thm}[\cite{scholzecomplexlectures}, Theorem 8.6]
    \begin{enumerate}
        \item 
        The category $\mathbf{Nuc}(\mathbf{C})$ is stable, closed under colimits in $\mathbf{C}$, and closed under tensor products in $\mathbf{C}$.
        \item 
        $\mathbf{Nuc}(\mathbf{C})$ is $\aleph_{1}$-compactly presented, with $\aleph_{1}$-compact generators given by the \textit{basic nuclear objects}, i.e. objets which are colimits of countable sequences
        $$C_{0}\rightarrow C_{1}\rightarrow\cdots\rightarrow C_{n}\rightarrow C_{n+1}\rightarrow\cdots$$
        where each $C_{i}$ is compact and each $C_{i}\rightarrow C_{i+1}$ is trace class.
    \end{enumerate}
\end{thm}

The $\aleph_{1}$-compact objects are the so-called \textit{basic nuclear objects} which can be written as sequential colimits
$$\limind_{\mathbb{N}}P_{n}$$
with each $P_{n}$-compact and each $P_{n}\rightarrow P_{n+1}$ nuclear.

\begin{defn}
Let $A,B\in\mathbf{Comm}(\mathbf{C})$. A map $f:A\rightarrow B$ of commutative monoids is said to be \textit{(strongly) nuclear} if $B$ is (strongly) nuclear as a left $A$-module.
\end{defn}

\begin{prop}
    Let $f:A\rightarrow B$ be a map in $\mathbf{Comm}(\mathbf{C})$. If $f:X\rightarrow Y$ is a trace class map in ${}_{A}\mathbf{Mod}$, then $B\otimes_{A}^{\mathbb{L}}X\rightarrow B\otimes_{A}^{\mathbb{L}}Y$ is trace class in ${}_{B}\mathbf{Mod}$.
\end{prop}

\begin{proof}
This follows from the commutative diagram 
\begin{displaymath}
    \xymatrix{
\pi_{0}\mathbf{Map}_{A}(A,X^{\vee_{A}}\otimes^{\mathbb{L}}_{A}Y)\ar[d]\ar[r] & \pi_{0}\mathbf{Map}_{A}(X,Y)\ar[d]\\
\pi_{0}\mathbf{Map}_{B}(B,(B\otimes_{A}^{\mathbb{L}}X)^{\vee_{B}}\otimes^{\mathbb{L}}_{B}B\otimes_{A}^{\mathbb{L}}Y)\ar[r] &\pi_{0}\mathbf{Map}_{B}(B\otimes_{A}^{\mathbb{L}}X,B\otimes_{A}^{\mathbb{L}}Y)
    }
\end{displaymath}
\end{proof}

\begin{cor}\label{cor:nucmodtens}
    Let $f:A\rightarrow B$ be a map in $\mathbf{Comm}(\mathbf{C})$. If $M$ is a nuclear left $A$-module then $B\otimes_{A}^{\mathbb{L}}M$ is a nuclear left $B$-module.
\end{cor}

\begin{lem}\label{lem:nucfacts}
\begin{enumerate}
\item 
Isomorphisms are strongly nuclear.
\item 
If $A\rightarrow B$ is nuclear and $A\rightarrow C$ is any map, then $C\rightarrow B\otimes^{\mathbb{L}}_{A}C$ is nuclear. 
\item 
Let $A\rightarrow B$ be strongly nuclear. If $M$ is nuclear in ${}_{A}\mathbf{Mod}(\mathbf{C})$, then $M$ is nuclear as an object of ${}_{A}\mathbf{Mod}$.
\item
Let $A\rightarrow B$ be strongly nuclear. If $M$ is a $B$-module which is strongly nuclear as an $A$-module, then it is strongly nuclear as a $B$-module.
\item
Let $A\rightarrow B\rightarrow C$ be maps in $\mathbf{Comm}(\mathbf{C})$ such that $A\rightarrow C$ and $A\rightarrow B$ are strongly nuclear. Then $B\rightarrow C$ is strongly nuclear.
\end{enumerate}
\end{lem}

\begin{proof}
    \begin{enumerate}
    \item 
    This is clear.
     \item 
    This follows immediately from Corollary \ref{cor:nucmodtens}.
    \item 
    We may assume that $A=\mathbb{I}$ and that $M$ is basic nuclear. Write $M\cong\limind_{\mathbb{N}}A\otimes^{\mathbb{L}}P_{i}$ where each $P_{i}$ is compact and each $A\otimes^{\mathbb{L}}P_{i}\rightarrow A\otimes^{\mathbb{L}}P_{i+1}$ is trace class as a map of $B$-modules, and $\mathcal{I}$ is $\aleph_{1}$-filtered. It suffices to show that it is trace class as a map in $\mathbf{C}$. Now we have $\underline{\mathbf{Map}}(B\otimes^{\mathbb{L}}P_{i},B)\cong P_{i}^{\vee}\otimes^{\mathbb{L}}B$, since $B$ is strongly nuclear. Since $f_{i}:B\otimes^{\mathbb{L}}P_{i}\rightarrow B\otimes^{\mathbb{L}}P_{i+1}$  is trace class it is in the image of
    \begin{align*}
  \pi_{0}\mathbf{Map}_{\mathbf{C}}(\mathbb{I},(B\otimes^{\mathbb{L}}P_{i})^{\vee}\otimes^{\mathbb{L}}B)&\cong\pi_{0}\mathbf{Map}_{{}_{B}\mathbf{Mod}}(B,P_{i}^{\vee}\otimes^{\mathbb{L}}(B\otimes^{\mathbb{L}}P_{i+1}))  )\\
  &\cong \pi_{0}\mathbf{Map}_{{}_{B}\mathbf{Mod}}(A,(P_{i}^{\vee}\otimes^{\mathbb{L}}B)\otimes_{B}^{\mathbb{L}}(B\otimes^{\mathbb{L}}P_{i+1})\\  &\rightarrow\pi_{0}\mathbf{Map}_{{}_{B}\mathbf{Mod}}(B\otimes^{\mathbb{L}}P_{i},B\otimes^{\mathbb{L}}P_{i+1})\\
    &\cong
\pi_{0}\mathbf{Map}_{\mathbf{C}}(P_{i},B\otimes^{\mathbb{L}}P_{i+1})
    \end{align*}

        \item 
   Let $B\otimes_{A}^{\mathbb{L}}d$ be a compact left $B$-module. We have
   \begin{align*}
   \underline{\mathbf{Map}}_{B}(B\otimes^{\mathbb{L}}d,M)&\cong \underline{\mathbf{Map}}_{A}(A\otimes^{\mathbb{L}}d,M)\\
   &\cong\underline{\mathbf{Map}}_{A}(A\otimes^{\mathbb{L}}d,A)\otimes_{A}^{\mathbb{L}}M\\
   &\cong\underline{\mathbf{Map}}_{A}(A\otimes^{\mathbb{L}}d,A)\otimes_{A}^{\mathbb{L}}B\otimes^{\mathbb{L}}_{B}M\\
   &\cong\underline{\mathbf{Map}}_{A}(A\otimes^{\mathbb{L}} d,B)\otimes^{\mathbb{L}}_{B}M
   \end{align*}
   as required.
        \item 
        This is the same as part (4).
    \end{enumerate}
\end{proof}

Let $A\rightarrow B$ be a map in $\mathbf{Comm}(\mathbf{C})$. By the above we get a functor
$$B\otimes_{A}^{\mathbb{L}}(-):\mathbf{Nuc}_{A}(\mathbf{C})\rightarrow\mathbf{Nuc}_{B}(\mathbf{C})$$
Now $\mathbf{Nuc}_{A}(\mathbf{C})$ is presentable and in fact closed under colimits in ${}_{A}\mathbf{Mod}(\mathbf{C})$. And $B\otimes_{A}^{\mathbb{L}}(-)$ commutes with colimits. Thus it has a right adjoint $|-|_{\mathbf{Nuc}}$. If it so happens that strongly nuclear $A$ and $B$ modules coincide with nuclear modules, and $A\rightarrow B$ is a nuclear map, then this right adjoint $|-|_{\mathbf{Nuc}}$ is just the usual forgetful functor. 

An immediate consequence of Proposition \ref{lem:nucfacts} is the following.

\begin{cor}
    Let $A\rightarrow B$ be strongly nuclear. If any nuclear $A$-module $M$ is strongly nuclear, then any nuclear $B$-module $M$ is strongly nuclear.
\end{cor}

\section{Derived and Spectral Algebraic Pre-contexts}\label{sec:dac}

In this section we define a class of stable $(\infty,1)$-categories with $t$-structures with sufficient structure for generalising definitions and results from algebraic geometry, including those concerning \'{e}tale and smooth maps. These generalise the derived algebraic contexts of \cite{raksit2020hochschild}. 

\subsection{Projectves in Stable $(\infty,1)$-Categories}

We must first recall some basic facts regarding projectives in stable $(\infty,1)$-categories.

\begin{defn}[\cite{HA} 5.5.8.18]\label{defn:stableproj}
\begin{enumerate}
\item
An object $P$ of an $(\infty,1)$-category is said to be \textit{projective} if the functor $\mathbf{Map}(P,-)$ commutes geometric realisations 
\item
Let $(\mathbf{C},\mathbf{C}_{\ge0},\mathbf{C}_{\le0})$ be a stable $(\infty,1)$-category with $t$-structure. An object $P$ of $\mathbf{C}_{\ge0}$ is said to be a \textit{projective in }$(\mathbf{C},\mathbf{C}_{\ge0},\mathbf{C}_{\le0})$  if it is projective when regarded as an object of $\mathbf{C}_{\ge0}$.
\end{enumerate}
\end{defn}

\begin{lem}
Let $(\mathbf{C},\mathbf{C}_{\ge0},\mathbf{C}_{\le0})$ be a stable $(\infty,1)$-category with a left complete $t$-structure and $P$ a projective object of $\mathbf{C}_{\ge0}$. Then $\pi_{0}(P)$ is projective in $\mathbf{C}^{\heart}$
\end{lem}

\begin{proof}
This follows immediately from \cite{lurieDAG1} Proposition 6.16, noting that for $P$ a projective and $A\in\mathbf{C}^{\heart}$
$$\pi_{0}\mathbf{Map}(P,A)\cong\mathrm{Hom}(\pi_{0}(P),A)$$
\end{proof}


\begin{defn}
Let $(\mathbf{C},\mathbf{C}_{\ge0},\mathbf{C}_{\le0})$ be an $(\infty,1)$-category with $t$-structure. A collection of objects $\mathbf{C}^{0}\subset\mathbf{C}_{\ge0}$ is said to be a \textit{projective generating collection} if 
\begin{enumerate}
\item
each $P\in\mathbf{C}^{0}$ is projective
\item
Every object in $\mathbf{C}$ can be expressed as a sifted colimit of direct sums of objects of $\mathbf{C}^{0}$.
\end{enumerate}
$(\mathbf{C},\mathbf{C}_{\ge0},\mathbf{C}_{\le0})$ is said to \textit{have enough projectives} if there exists a projective generating collection. We denote a stable $(\infty,1)$-category with $t$-structure with a prescribed projective generating collection $\mathbf{C}^{0}$ as
$$(\mathbf{C},\mathbf{C}_{\ge0},\mathbf{C}_{\le0},\mathbf{C}^{0})$$
We refer to such a tuple as a \textit{ projectively generated stable} $(\infty,1)$-\textit{category with} $t$-\textit{structure}.
\end{defn}

\begin{lem}[c.f. \cite{HA} Proposition 7.2.2.7]
Let $(\mathbf{C},\mathbf{C}_{\ge0},\mathbf{C}_{\le0},\mathbf{C}^{0})$ be a projectively generated stable $(\infty,1)$-category with a left complete $t$-structure. Then $Q\in\mathbf{C}$ is projective, if and only if $Q$ is a retract of an object of the form $\bigoplus_{i}P_{i}$ for some $P_{i}\in\mathbf{C}^{0}$.
\end{lem}

\begin{proof}
First it is clear that all such modules are projective. Conversely let $Q$ be projective. We may write $Q$ as the geometric realisation of a simplicial diagram $Q_{n}$ with $Q_{n}=\bigoplus_{i_{n}}P_{i_{n}}$. Then $\pi_{0}(\bigoplus_{i_{0}}P_{i_{0}})\rightarrow\pi_{0}(Q)$ is a surjection in $\mathbf{C}^{\heart}$. By \cite{HA} Proposition 7.2.2.6 this map splits, and so $Q$ is a retract of $\bigoplus_{i_{0}}P_{i_{0}}$,
\end{proof}

Let $(\mathbf{Sp},\mathbf{Sp}_{\ge0},\mathbf{Sp}_{\le0})$ denote the $(\infty,1)$-category of spectra equipped with its standard $t$-structure. If $\mathbf{C}$ is a presentable stable $(\infty,1)$-category then it has a $\mathbf{Sp}$-enrichment $\mathbf{Map}^{s}(-,-)$.

\begin{lem}
Let $(\mathbf{C},\mathbf{C}_{\ge0},\mathbf{C}_{\le0})$ be a right complete stable $(\infty,1)$-category with $t$-structure. Let $\mathbf{C}^{0}$ be a collection of projective objects of $\mathbf{C}_{\ge0}$ such that every object in $\mathbf{C}_{\ge0}$ can be written as a colimit of objects in $\mathbf{C}^{0}$. Then an object $X$ of $\mathbf{C}$ is
\begin{enumerate}
\item
in $\mathbf{C}_{\le0}$ if and only if $\mathbf{Map}^{s}(P,N)\in\mathbf{Sp}_{\le0}$ for all $P\in\mathbf{C}^{0}$.
\item
in $\mathbf{C}_{\ge0}$ if and only if $\mathbf{Map}^{s}(P,N)\in\mathbf{Sp}_{\ge0}$ for all $P\in\mathbf{C}^{0}$.
\end{enumerate}
\end{lem}

\begin{proof}
\begin{enumerate}
\item
Since $\mathbf{C}$ is right complete we have $\mathbf{C}\cong\mathbf{Sp}(\mathbf{C}_{\ge0})$. For $X=(X_{0},X_{1},\ldots,X_{n},\ldots)$ an object of $\mathbf{Sp}(\mathbf{C}_{\ge0})$ we have that 
$$\mathbf{Map}^{s}(P,X)\cong(\mathbf{Map}(P,X_{0}),\mathbf{Map}^(P,X_{1}),\ldots,\mathbf{Map}(P,X_{n}),\ldots)$$
Since every object of $\mathbf{C}_{\ge0}$ can be written as a colimit of objects in $\mathbf{C}^{0}$ we have that $X_{0}\cong0$ if and only if $\mathbf{Map}^{s}(P,X_{0})\cong 0$  for all $P\in\mathbf{C}^{0}$. This proves the claim.
\item
Since $P\in\mathbf{C}^{0}$ is projective $\mathbf{Map}(P,-)$ commutes with geometric realisation and hence suspension. Thus $\mathbf{Map}(P,\Sigma^{\infty}(Q))\cong\Sigma^{\infty}\mathbf{Map}(P,Q)$ for any $Q\in\mathbf{C}_{\ge0}$. 
Conversely let $N$ be such that $\mathbf{Map}(P,N)\in\mathbf{Sp}_{\ge0}$ for all $P\in\mathbf{C}^{0}$. Let
$$N'\rightarrow N\rightarrow N''$$
be an exact triangle with $N'\in\mathbf{C}_{\ge0}$ and $N''\in\mathbf{C}_{\le- 1}$. Then
$$\mathbf{Map}^{s}(P,N')\rightarrow\mathbf{Map}^{s}(P,N)\rightarrow\mathbf{Map}^{s}(P,N'')$$
is an exact triangle. Now $\mathbf{Map}^{s}(P,N'),\mathbf{Map}^{s}(P,N)$ are in $\mathbf{Sp}_{\ge0}$, while $\mathbf{Map}(P,N'')\in\mathbf{Sp}_{\le -1}$. We therefore must have \[\mathbf{Map}(P,N'')\cong 0.\] This is true for all $P\in\mathbf{C}^{0}$. Therefore $N''\cong 0$. 
\end{enumerate}
\end{proof}

\subsubsection{Eilenberg-Maclane Categories}

In this subsection we provide conditions on a stable $(\infty,1)$-category such that the Eilenberg-Maclane construction still works.

\begin{defn}
A projectively generated stable $(\infty,1)$-category with $t$-structure
$$(\mathbf{C},\mathbf{C}_{\ge0},\mathbf{C}_{\le0},\mathbf{C}^{0})$$
is said to be \textit{Eilenberg-Maclane} if
\begin{enumerate}
\item
$\mathbf{C}_{\le0}$ is closed under filtered colimits.
\item
Each $P\in\mathbf{C}^{0}$ is tiny in $\mathbf{C}_{\ge0}$, i.e. $\mathbf{Map}(P,-):\mathbf{C}_{\ge0}\rightarrow\mathbf{Sp}$ commutes with filtered colimits.
\item
The $t$-structure is left and right complete.  
\end{enumerate}
\end{defn}

The category $\mathbf{C}^{\heart}$ is a reflective subcategory of $\mathbf{C}_{\ge0}$. Thus $\mathbf{C}^{\heart}$ is generated under ($1$-categorical) sifted colimits by the full subcategory 
$$\pi_{0}(\mathbf{C}^{0})\defeq\{\pi_{0}(P):P\in\mathbf{C}^{0}\}$$ 
Now $\pi_{0}(P)$ is projective in $\mathbf{C}^{\heart}$. Therefore the abelian category $\mathbf{C}^{\heart}$ has enough projectives. If in addition $\mathbf{C}_{\le0}$ is closed under filtered colimits then the functor $\mathbf{C}^{\heart}\rightarrow\mathbf{C}_{\ge0}$ commutes with filtered colimits. In particular if each $P\in\mathbf{C}^{0}$ is tiny then each $\pi_{0}(P)$ is tiny in $\mathbf{C}^{\heart}$. In this case $\mathbf{C}^{\heart}$ is an \textit{elementary} abelian category in the sense of \cite{qacs}. By \cite{kelly2016homotopy} we have proven the following.

\begin{lem}
Let 
$$(\mathbf{C},\mathbf{C}_{\ge0},\mathbf{C}_{\le0},\mathbf{C}^{0})$$ be a projectively generated stable $(\infty,1)$-category with $t$-structure. Then the projective model structure exists on $\mathrm{Ch}_{\ge0}(\mathbf{C}^{\heart})$. If $(\mathbf{C},\mathbf{C}_{\ge0},\mathbf{C}_{\le0},\mathbf{C}^{0})$ is Eilenberg-Maclane then the projective model structure exists on $\mathrm{Ch}(\mathbf{C}^{\heart})$ and the model structures on both $\mathrm{Ch}_{\ge0}(\mathbf{C}^{\heart})$ and $\mathrm{Ch}(\mathbf{C}^{\heart})$ are combinatorial. 
\end{lem}
Let $(\mathbf{C},\mathbf{C}_{\ge0},\mathbf{C}_{\le0},\mathbf{C}^{0})$ be Eilenberg-Maclane. Consider the $(\infty,1)$-categories $\mathbf{Ch}(\mathbf{C}^{\heart})$ and $\mathbf{Ch}_{\ge0}(\mathbf{C}^{\heart})$ presented by the model categories $Ch(\mathbf{C}^{\heart})$ and $Ch_{\ge0}(\mathbf{C}^{\heart})$ respectively. $\mathbf{Ch}_{\ge0}(\mathbf{C}^{\heart})$ is equivalent to the free sifted cocompletion $\mathbf{P}_{\Sigma}(\pi_{0}(\mathbf{C}^{0}))$.  The inclusion
$$\pi_{0}(\mathbf{C}^{0})\rightarrow\mathbf{C}_{\ge0}$$
extends to a functor

$$K:\mathbf{Ch}_{\ge0}(\mathbf{C}^{\heart})\rightarrow\mathbf{C}_{\ge0}$$
Moreover $\mathbf{C}_{\ge0}$ is equivalent to the free sifted cocompletion $\mathbf{P}_{\Sigma}(\mathbf{C}^{0})$. Thus $\pi_{0}:\mathbf{C}^{0}\rightarrow\pi_{0}(\mathbf{C}^{0})$ induces by sifted colimits a functor
$$\tilde{\Gamma}:\mathbf{C}_{\ge0}\rightarrow \mathbf{Ch}_{\ge0}(\mathbf{C}^{\heart})$$
There is an adjunction
$$\adj{\tilde{\Gamma}}{\mathbf{C}_{\ge0}}{\mathbf{Ch}_{\ge0}(\mathbf{C}^{\heart})}{K}$$
which induces by stabilisation of $K$ and the adjoint functor theorem (since $K$ is accessible), an adjunction
$$\adj{\Gamma}{\mathbf{C}}{\mathbf{Ch}(\mathbf{C}^{\heart})}{K}$$
Note that of course `stabilisation' is in general not a functorial construction. We are very much using the fact here that $K$ commutes with limits and is accessible.

\begin{defn}
$K:\mathbf{Ch}(\mathbf{C}^{\heart})\rightarrow\mathbf{C}$ is called the \textit{Eilenberg-Maclane functor of }$\mathbf{C}$.
\end{defn}

\begin{rem}
$$K:\mathbf{Ch}(\mathbf{C}^{\heart})\rightarrow\mathbf{C}$$
 is exact for the given $t$-structure on $\mathbf{C}$ and the standard $t$-structure on $\mathbf{Ch}(\mathbf{C}^{\heart})$. In particular we have $\pi_{i}(K(X))\cong H_{i}(X)$ for any $X\in\mathbf{Ch}(\mathbf{C}^{\heart})$ and any $i\in\mathbb{Z}$. 
\end{rem}

We also have the following explicit description of $\mathbf{C}^{\heart}$

\begin{prop}
Let $(\mathbf{C},\mathbf{C}_{\ge0},\mathbf{C}_{\le0},\mathbf{C}^{0})$ be Eilenberg-Maclane. Then $\mathbf{C}^{\heart}$ is equivalent to the category $\mathrm{Add}(\pi_{0}(\mathbf{C}^{0},\mathrm{Ab})$ of contravariant additive functors from $\pi_{0}(\mathbf{C}^{0})$ to $\mathrm{Ab}$.
\end{prop}

\begin{proof}
Each $\pi_{0}(P)$ is compact and projective in $\mathbf{C}^{\heart}$. Moreover $\pi_{0}(\mathbf{C}^{0})$ clearly generates $\mathbf{C}^{\heart}$ under sifted colimits. This proves that $\mathbf{C}^{\heart}\cong\mathrm{Add}(\pi_{0}(\mathbf{C}^{0})^{op},\mathrm{Ab})$. Note that the functor $\mathbf{C}^{\heart}\rightarrow\mathrm{Add}(\pi_{0}(\mathbf{C}^{0})^{op},\mathrm{Ab})$ is induced by mapping $\pi_{0}(P)$ to $\mathrm{Hom}(-,P)$ for $P\in\mathbf{C}^{0}$. 
\end{proof}
%

%

\section{Spectral Algebraic Pre-contexts}\label{sec:spapc}

The derived algebraic contexts of \cite{raksit2020hochschild} provide a convenient setting for relative geometry. In particular such categories have canonical $t$-structures which interact well with the monoidal structure. This allows one to generalise many of the useful spectral sequences appearing in e.g. \cite{toen2008homotopical}. Most of the $(\infty,1)$-algebra contexts in which we are interested also have the structure of a derived algebraic context. Before recalling the definition, let us introduce a concept which is weaker in two senses - we relax an additivity requirement, and a requirement regarding symmetric products in the heart. 

\subsection{Definitions}

\begin{defn}
A \textit{weak spectral algebraic pre-context} is a tuple $(\mathbf{C},\mathbf{C}_{\ge0},\mathbf{C}_{\le0},\mathbf{C}^{0})$ where
\begin{enumerate}
\item
$\mathbf{C}$ is a stably symmetric monoidal locally presentable $(\infty,1)$-category.
\item
$(\mathbf{C}_{\ge0},\mathbf{C}_{\le0})$ is a complete $t$-structure on $\mathbf{C}$, which is compatible with the monoidal structure. 
\item
$\mathbf{C}^{0}\subset\mathbf{C}_{\ge0}$ is a small full subcategory of compact projective objects such that 
\begin{enumerate}
\item
The monoidal unit $\mathbb{I}$ is in $\mathbf{C}^{0}$
\item
If $P\in\mathbf{C}^{0}$ then $\pi_{*}(P)$ is a projective $\pi_{*}(\mathbb{I})$-module. 
\item
$\mathbf{C}^{0}$ is closed under the formation of finite coproducts.
\item
$\mathbf{C}_{\ge0}$ is generated under sifted colimits by $\mathbf{C}^{0}$. 
\item
If $P,Q\in\mathbf{C}^{0}$ then $P\otimes Q\in\mathbf{C}^{0}$ and $\pi_{0}(P\otimes^{\mathbb{L}} Q)\cong\pi_{0}(P)\otimes\pi_{0}(Q)$.
\end{enumerate}
\end{enumerate}
\end{defn}

In particular if $(\mathbf{C},\mathbf{C}_{\ge0},\mathbf{C}_{\le0},\mathbf{C}^{0})$ is a weak spectral algebraic pre-context then it is Eilenberg-Maclane. $\mathbf{C}^{\heart}$ has an explicit construction as a monoidal elementary abelian category in the sense of \cite{kelly2016homotopy} Definition 2.6.117. 

\begin{prop}
$\mathbf{C}^{\heart}$ is equivalent as a symmetric monoidal category to the category $\mathrm{Add}(\pi_{0}(\mathbf{C}^{0})^{op},\mathrm{Ab})$ of contravariant additive functors from $\pi_{0}(\mathbf{C}^{0})$ to $\mathrm{Ab}$, equipped with the Day convolution symmetric monoidal structure.
\end{prop}

\begin{proof}
It suffices to prove that Day convolution and the monoidal structure on $\mathbf{C}^{\heart}$ agree on $\pi_{0}(\mathbf{C}^{0})$. Now since $\pi_{0}(P)\otimes\pi_{0}(Q)\cong\pi_{0}(P\otimes Q)\in\pi_{0}(\mathbf{C}^{0})$, classical results regarding Day convolution imply that the convolution tensor product of $\mathrm{Hom}(-.\pi_{0}(P))$ and $\mathrm{Hom}(-,\pi_{0}(Q))$ is $\mathrm{Hom}(-,\pi_{0}(P)\otimes\pi_{0}(Q))$, as required.
\end{proof}

The following is also immediate.

\begin{cor}
 If $(\mathbf{C},\mathbf{C}_{\ge0},\mathbf{C}_{\le0},\mathbf{C}^{0})$ is a weak spectral algebraic pre-context the functor $\Gamma$ is strongly monoidal.
\end{cor}

Before continuing we introduce one final definition.

\begin{defn}
A weak spectral algebraic pre-context \[(\mathbf{C},\mathbf{C}_{\ge0},\mathbf{C}_{\le0},\mathbf{C}^{0})\] is said to be a \textit{symmetric weak spectral algebraic pre-context} if for any $P\in\mathbf{C}^{0}$ and any $n\in\mathbb{N}_{0}$ one has $\mathrm{Sym}^{n}_{\mathbf{C}^{\heart}}(\pi_{0}(P))\in\pi_{0}(\mathbf{C}^{0})$. 
\end{defn}

\subsubsection{Strong Morphisms and Strong Modules}

\begin{defn}
Let $A$ be an object of $\mathbf{Ass}(\mathbf{C})$ and $M$ a left $A$-module. $M$ is said to be \textit{strong} if the map
$$\pi_{*}(A)\otimes_{\pi_{0}(A)}\pi_{0}(M)\rightarrow\pi_{*}(M)$$
is an isomorphism of graded modules in $\mathbf{C}^{\heart}$, and \textit{derived strong} if in addition the map
$$\pi_{*}(A)\otimes^{\mathbb{L}}_{\pi_{0}(A)}\pi_{0}(M)\rightarrow\pi_{*}(M)$$
is an equivalence.

A map $f:A\rightarrow B$ in $\mathbf{Ass}(\mathbf{C})$  is said to be \textit{(derived) strong} if $B$ is (derived) strong when regarded as a left $A$-module.
\end{defn}

\begin{rem}\label{rem:flatstrongderivedstrong}
An $A$-module $M$ is derived strong if and only if it is strong and each $\pi_{n}(A)$ is transverse to $\pi_{0}(M)$ over $\pi_{0}(A)$. In particular if $\pi_{0}(M)$ is flat as a $\pi_{0}(A)$-module then $M$ is strong if and only if it is derived strong.
\end{rem}

\begin{defn}
    A weak spectral algebraic pre-context is said to be a \textit{spectral algebraic pre-context} if all projectives in $\mathbf{C}_{\ge0}$ are strong.
 \end{defn}

 This is equivalent to all objects of $\mathbf{C}^{0}$ being strong.

\begin{example}
\begin{enumerate}
\item
Any weak spectral algebraic pre-context with $\mathbf{C}^{0}\subset\mathbf{C}^{\heart}$ is a spectral algebraic pre-context. 
\item
$\mathbf{Sp}$ is a spectral algebraic pre-context. 
\end{enumerate}
\end{example}

\subsubsection{The $\mathbf{Map}$ Spectral Sequences}

In weak spectral algebraic pre-contexts we have a useful spectral sequence for computing homotopy groups of mapping spectra. 

\begin{lem}\label{lem:mapspecprep}
Let $(\mathbf{C},\mathbf{C}_{\ge0},\mathbf{C}_{\le0},\mathbf{C}^{0})$ be a spectral algebraic pre-context. If $N$ is object in $\mathbf{C}_{\ge0}$, and $P$ is projective in $\mathbf{C}$, then we have
$$\mathrm{Hom}^{0}_{\pi_{*}(R)}(\pi_{*}(P),\pi_{*}(N))\cong\pi_{*}(\mathbf{Map}(P,N))$$
\end{lem}

\begin{proof}
This immediately follows from the fact that $\mathrm{Ext}^{i}(P,N)\cong 0$ for $i\neq0$ by \cite{HA} Proposition 7.2.2.6.
\end{proof}

The following can be proven as in \cite{HA} Variant 7.2.1.24.

\begin{lem}
Let $(\mathbf{C},\mathbf{C}_{\ge0},\mathbf{C}_{\le0},\mathbf{C}^{0})$ be a spectral algebraic pre-context. Let $N\in\mathbf{C}$ be such that $\pi_{n}(N)=0$ for $n>>0$. Let $M\in\mathbf{C}_{\ge0}$. Then there is a strongly convergent spectral sequence
$$\pi_{p}\mathbf{Map}_{\pi_{*}(R)}(\pi_{*}(M),\pi_{*}(N))_{q}\Rightarrow\pi_{-p-q}\mathbf{Map}(M,N)$$
\end{lem}

\subsection{Flat (Weak) Spectral Algebraic Pre-contexts}

For the purposes of derived geometry it is convenient to require a flatness assumption on the generators $\mathbf{C}^{0}$.

\begin{defn}
Let \((\mathbf{C},\mathbf{C}_{\geq 0}, \mathbf{C}_{\leq 0},\mathbf{C}^{0})\)  be a weak spectral algebraic pre-context. An object $P$ of $\mathbf{C}_{\geq 0}$ is said to be \textit{homotopy flat} if for any object $N$ of $\mathbf{C}$,  the natural map of graded objects
$$\pi_{*}(M)\otimes^{\mathbb{L}}\pi_{*}(P)\rightarrow\pi_{*}(M\otimes^{\mathbb{L}} P)$$
is an isomorphism.
\end{defn}

\begin{defn}
A (weak) spectral algebraic pre-context 
 \((\mathbf{C},\mathbf{C}_{\geq 0}, \mathbf{C}_{\leq 0},\mathbf{C}^{0})\) 
is said to be \textit{flat} if any object $P$ of $\mathbf{C}^{0}$ is homotopy flat.
\end{defn}

There are several important consequences of a spectral algebraic pre-context being flat. The first is the following spectral sequence. From now on we will assume that we are working in a flat spectral algebraic pre-context.

\begin{lem}\label{lem:flattorspec}
Let $(\mathbf{C},\mathbf{C}_{\ge0},\mathbf{C}_{\le0},\mathbf{C}^{0})$ be a flat spectral algebraic pre-context. Then there is a spectral sequence
$$\mathrm{Tor}_{p}^{\pi_{*}(A)}(\pi_{*}(B),\pi_{*}(C))_{q}\Rightarrow\pi_{p+q}(B\otimes^{\mathbb{L}}_{A}C)$$
\end{lem}

\begin{proof}
This works as in \cite{HA} Proposition 7.2.1.17 once we can show it for modules of the form $B\cong \bigoplus_{i\in\mathcal{I}}A\otimes P_{i}$ for some set $\mathcal{I}$, and some $P_{i}\in\mathbf{C}^{0}$. Then
$$\mathrm{Tor}_{0}^{\pi_{*}(A)}(\pi_{*}(B),\pi_{*}(C))\cong\pi_{*}(B\otimes^{\mathbb{L}}_{A}C)$$
As in \cite{HA} Proposition 7.2.1.17 we may assume that $B\cong A\otimes P$ for some $P\in\mathbf{C}^{0}$. Then the result follows from the homotopy flat assunptions. 
\end{proof}

We will repeatedly make use of this spectral sequence in what follows.





Unless we state otherwise, from now on w fix a spectral algebraic pre-context $(\mathbf{C},\mathbf{C}_{\ge0},\mathbf{C}_{\le0},\mathbf{C}^{0})$.

\subsubsection{Transverslity and Flatness}

\begin{defn}
\begin{enumerate}
\item
Let $M$ be a right $A$-module and $N$ a left $A$-module with $N$ discrete. $N$ is said to be \textit{transverse to }$M$ if $M\otimes^{\mathbb{L}}_{A}N$ is discrete. 
\item
Let $M$ be a right $A$-module. $M$ is said to be \textit{flat} if it is transverse to any discrete left $A$-module $N$.
\end{enumerate}
\end{defn}

\begin{prop}[c.f. \cite{toen2008homotopical} Lemma 2.2.2.2 (2)]\label{transversedisc}
Let $M$ be a connective right $A$-module and $N$ a left $A$-modules with $N$ discrete. If $N$ is transverse to $M$ then the map
$$\pi_{0}(M)\otimes^{\mathbb{L}}_{\pi_{0}(A)}N\rightarrow\pi_{0}(M\otimes^{\mathbb{L}}_{A}N)$$
is an equivalence. 
\end{prop}

\begin{proof}
Let
$$0\rightarrow K\rightarrow P\rightarrow N\rightarrow 0$$
be an exact sequence of $\pi_{0}(A)$-modules. We get a fibre-cofibre sequence
$$M\otimes^{\mathbb{L}}_{A}K\rightarrow M\otimes^{\mathbb{L}}_{A}P\rightarrow M\otimes^{\mathbb{L}}_{A}N$$
The long exact sequence on $\pi_{n}$ and the fact that $M\otimes^{\mathbb{L}}_{A}N$ is discrete then implies that
$$0\rightarrow\pi_{0}(M)\otimes_{\pi_{0}(A)}K\rightarrow\pi_{0}(M)\otimes_{\pi_{0}(A)}P\rightarrow\pi_{0}(M)\otimes_{\pi_{0}(A)}N\rightarrow0$$
is an exact sequence
\end{proof}

\begin{lem}[c.f. \cite{toen2008homotopical} Lemma 2.2.2.2 (2)]\label{lem:connectivetrans}
Let $M$ be a connective right $A$-module and $N$ a connective left $A$-module. Suppose that each $\pi_{n}(N)$ is transverse to $M$. Then there is an equivalence
$$\pi_{i}(M\otimes_{A}^{\mathbb{L}}N)\cong\pi_{0}(M)\otimes^{\mathbb{L}}_{\pi_{0}(A)}\pi_{i}(N)$$
\end{lem}

\begin{proof}
The previous Proposition implies that each $\pi_{n}(N)$ is transverse to $\pi_{0}(M)$ over $\pi_{0}(A)$.

By right exactness of the tensor product we have that
$$\pi_{i}(M\otimes_{A}^{\mathbb{L}}N)\rightarrow \pi_{i}(M\otimes_{A}^{\mathbb{L}}\tau_{\ge n}N)$$
is an equivalence whenever $i\le n$. Thus we may assume that $N$ is bounded. Let $\mathcal{N}$ be the category of connective left $A$-modules such that $\pi_{i}(M\otimes_{A}^{\mathbb{L}}K)\cong\pi_{0}(M)\otimes^{\mathbb{L}}_{\pi_{0}(A)}\pi_{i}(K)$. By induction on the length of $N$ we shall show that $N\in\mathcal{N}$. The previous Proposition proves the claim when $N$ is concentrated in degree $0$, and by shifting, when $N$ is concentrated in any single degree.
Suppose now $N$ has length $n+1$ and consider the fibre-cofibre sequence
$$\pi_{n+1}(N)[n+1]\rightarrow N\rightarrow\tau_{\ge n}(N)$$
By induction we may assume $\tau_{\ge n}(N)\in\mathcal{N}$.
We get a fibre-cofibre sequence
$$(M\otimes^{\mathbb{L}}_{A}\pi_{n+1}(N))[n+1]\rightarrow M\otimes^{\mathbb{L}}_{A}\tau_{\ge n+1}N\rightarrow M\otimes^{\mathbb{L}}_{A}\tau_{\ge n}N$$
An easy argument using the long exact sequence proves the claim.
\end{proof}


\begin{defn}
\begin{enumerate}
\item
A map $f:A\rightarrow B$ in $\mathbf{Ass}(\mathbf{C})$ is said to be \textit{transverse} if every discrete $B$-module is transverse to $B$ as an $A$-module. 
\item
A map $f:A\rightarrow B$ in $\mathbf{Ass}(\mathbf{C})$ is said to be \textit{flat} if every discrete $A$-module is transverse to $B$ as an $A$-module. 
\item
A map $f:A\rightarrow B$ in $\mathbf{Ass}(\mathbf{C})$ is said to be \textit{faithfully flat} if it is flat, and whenever $M$ is an $A$-module such that $B\otimes^{\mathbb{L}}_{A}M\cong 0$, then $M\cong 0$.
\end{enumerate}
\end{defn}

\begin{example}
Let $f:A\rightarrow B$ be a map of algebras. $f$ is said to be \textit{finite} if $B\cong A^{\oplus n}$ as a left $A$-module, for some $n$. Then $A^{\oplus n}\otimes^{\mathbb{L}}_{A}M\cong M^{\oplus n}$. Hence $f$ is clearly faithfully flat. 
\end{example}

\subsubsection{Strongness in Spectral Algebraic Pre-Contexts}

\begin{rem}\label{rem:conntransA}
By Lemma \ref{lem:connectivetrans}, if $A$ and $M$ are connective, then $M$ is a derived strong $A$-module if and only if each $\pi_{n}(A)$ is transverse to $M$.
\end{rem}

\begin{cor}[c.f. \cite{toen2008homotopical} Lemma 2.2.2.2 (2)]\label{cor:flatderstrong}
Let $A$ and $M$ be connective. The following are equivalent.
\begin{enumerate}
\item
$M$ is flat as an $A$-module.
\item
$M$ is strong as an $A$-module and $\pi_{0}(M)$ is flat as a $\pi_{0}(A)$-module.
\item
$M$ is derived strong as an $A$-module and $\pi_{0}(M)$ is flat as a $\pi_{0}(A)$-module.
\item
$M\otimes^{\mathbb{L}}_{A}(-):{}_{A}\mathbf{Mod}(\mathbf{C}_{\ge0})\rightarrow\mathbf{C}_{\ge0}$ commutes with finite limits.
\end{enumerate}
\end{cor}

\begin{proof}
$1)\Rightarrow 2)$. By Remark \ref{rem:conntransA} $M$ is derived strong. By Proposition \ref{transversedisc} $\pi_{0}(M)$ is flat as a $\pi_{0}(A)$-module. 

$2)\Rightarrow 1)$.  Let $N$ be a discrete right $A$-module. We have $\pi_{*}(M\otimes^{\mathbb{L}}_{A}N)\cong\pi_{0}(M)\otimes_{\pi_{0}(A)}\pi_{*}(N)\cong\pi_{0}(M)\otimes_{\pi_{0}(A)}N$ which is clearly discrete. 

The equivalence of $(2)$ and $(3)$ is Remark \ref{rem:flatstrongderivedstrong}.

$3)\Rightarrow 4)$. It suffices to prove that $\tau_{\ge0}(M\otimes^{\mathbb{L}}_{A}B)\cong M\otimes^{\mathbb{L}}_{A}\tau_{\ge0}B$ for each $B\in{}_{A}\mathbf{Mod}(\mathbf{C})$. Now we have $\pi_{*}(M\otimes^{\mathbb{L}}_{A}B)\cong\pi_{0}(M)\otimes_{\pi_{0}(A)}\pi_{*}(B)$. On the other hand we also have $\pi_{*}(M\otimes^{\mathbb{L}}_{A}\tau_{\ge0}(B))\cong\pi_{0}(M)\otimes_{\pi_{0}(A)}\pi_{*}(\tau_{\ge0}B)$. Hence the map
$$\tau_{\ge0}(M\otimes^{\mathbb{L}}_{A}B)\rightarrow M\otimes^{\mathbb{L}}_{A}\tau_{\ge0}B$$
clearly induces an equivalence of homotopy groups.

$4)\Rightarrow 1)$. Let $0\rightarrow N\rightarrow P$ be a left exact sequence in ${}_{A}\mathbf{Mod}(\mathbf{C})^{\heart}$. Then it is a homotopy fibre sequence in ${}_{A}\mathbf{Mod}(\mathbf{C}_{\ge0})$, since the inclusion ${}_{A}\mathbf{Mod}(\mathbf{C})^{\heart}\rightarrow{}_{A}\mathbf{Mod}(\mathbf{C}_{\ge0})$ is a right adjoint. Thus

$$0\rightarrow M\otimes_{A}^{\mathbb{L}}N\rightarrow M\otimes_{A}^{\mathbb{L}}P$$

is a homotopy fibre sequence. Applying $\pi_{0}$ we get a left exact sequence
$$0\rightarrow\pi_{0}(M)\otimes_{\pi_{0}(A)}N\rightarrow\pi_{0}(M)\otimes_{\pi_{0}(A)}P$$
which shows that $\pi_{0}(M)$ is flat as a right $\pi_{0}(A)$-module. Taking $N=0$ also gives that $\pi_{i}(M\otimes_{A}^{\mathbb{L}}P)\cong 0$ for any $i\ge0$. Now Remark \ref{rem:conntransA} implies that $M$ is derived strong as an $A$-module, as required.
\end{proof}

Note that isomorphisms are derived strong maps, and compositions of (derived) strong maps are (derived) strong.

Let $M\in\mathbf{C}_{\ge0}$ be projective and strong, and $N$ be any object of $\mathbf{C}_{\ge0}$. We have
\begin{align*}
\pi_{0}\mathbf{Map}(P,N)&\cong\mathrm{Hom}_{\pi_{*}(\mathbb{I})}(\pi_{*}(P),\pi_{*}(N))_{0}\\
&\cong\mathrm{Hom}_{\pi_{0}(\mathbb{I})}(\pi_{0}(P),\pi_{*}(N))_{0}\\
&\cong\mathrm{Hom}_{\pi_{0}(\mathbb{I})}(\pi_{0}(P),\pi_{0}(N))
\end{align*}

\begin{prop}
Let $\mathbf{C}$ be a flat, weak spectral algebraic pre-context. Let  $A\in\mathbf{Ass}^{cn}(\mathbf{C})$ and $M\in{}_{A}\mathbf{Mod}$ be derived strong. Then the map
$$\pi_{0}(A)\otimes^{\mathbb{L}}_{A}M\rightarrow\pi_{0}(M)$$
is an isomorphism.
\end{prop}

\begin{proof}
Consider the spectral sequence
$$\mathrm{Tor}^{p}_{\pi_{*}(A)}(\pi_{0}(A),\pi_{*}(M))_{q}\Rightarrow\pi_{p+q}(\pi_{0}(A)\otimes^{\mathbb{L}}_{A}M)$$
Now 
$$\pi_{0}(A)\otimes^{\mathbb{L}}_{\pi_{*}(M)}\pi_{*}(M)\cong\pi_{0}(A)\otimes^{\mathbb{L}}_{\pi_{0}(A)}\pi_{0}(M)\cong\pi_{0}(M)$$
Thus the spectral sequences collapses in the first page, and we have $\pi_{0}(M)\cong\pi_{0}(A)\otimes^{\mathbb{L}}_{A}B$.
\end{proof}

\begin{prop}\label{prop:derstronghtpy}
Let $M$ be a derived strong $A$-module and $N$ any $A$-module such that $\pi_{*}(N)$ is transverse to $\pi_{0}(M)$ over $\pi_{0}(A)$. Then there is an isomorphism
$$\pi_{*}(M\otimes^{\mathbb{L}}_{A}N)\cong\pi_{0}(M)\otimes^{\mathbb{L}}_{\pi_{0}(A)}\pi_{*}(N)$$
\end{prop}

\begin{proof}
Consider the spectral sequence
$$\mathrm{Tor}_{\pi_{*}(A)}^{p}(\pi_{*}(M),\pi_{*}(N))_{q}\Rightarrow\pi_{p+q}(M\otimes^{\mathbb{L}}_{A}N)$$
We have \[\pi_{*}(M)\cong\pi_{*}(A)\otimes_{\pi_{0}(A)}\pi_{0}(M),\] so $\pi_{*}(M)\otimes_{\pi_{*}(A)}\pi_{*}(N)\cong\pi_{0}(M)\otimes_{\pi_{0}(A)}\pi_{*}(N)$. Thus the spectral sequence degenerates, and the result follows. 
\end{proof}

\begin{cor}\label{cor:basechangederivestrong}
Let $f:A\rightarrow B$ be a map in $\mathbf{Ass}(\mathbf{C})$ and $M\in{}_{A}\mathbf{Mod}$  a derived strong $A$-module such that $\pi_{*}(B)$ is transverse to $\pi_{0}(M)$ over $\pi_{0}(A)$. Then $B\otimes^{\mathbb{L}}_{A}M$ is a derived strong $B$-module.
\end{cor}

\begin{proof}
Consider the spectral sequence
$$\mathrm{Tor}^{p}_{\pi_{*}(A)}(\pi_{*}(B),\pi_{*}(M))_{q}\Rightarrow\pi_{p+q}(B\otimes^{\mathbb{L}}_{A}M)$$
Now $\pi_{*}(B)\otimes^{\mathbb{L}}_{\pi_{*}(A)}\pi_{*}(M)\cong\pi_{*}(B)\otimes^{\mathbb{L}}_{\pi_{0}(A)}\pi_{0}(M)\cong\pi_{*}(B)\otimes_{\pi_{0}(A)}\pi_{0}(M)$. Thus the spectral sequence degenerates and we hae
$$\pi_{*}(M\otimes_{A}^{\mathbb{L}}B)\cong\pi_{*}(B)\otimes_{\pi_{0}(A)}\pi_{0}(M)\cong\pi_{*}(B)\otimes_{\pi_{0}(B)}^{\mathbb{L}}\pi_{0}(M\otimes_{A}^{\mathbb{L}}B)$$
\end{proof}

\begin{prop}\label{prop:reversederstrong}
    Let $A\in\mathbf{Ass}^{cn}(\mathbf{C})$ and let $M$ be a connective left $A$-module. Suppose that
    \begin{enumerate}
        \item 
        the map $\pi_{0}(A)\otimes^{\mathbb{L}}_{A}M\rightarrow\pi_{0}(M)$ is an equivalence.
        \item 
        $\pi_{n}(A)$ is transverse to $\pi_{0}(M)$ over $\pi_{0}(A)$ for each $n$.
    \end{enumerate}
    Then $M$ is a derived strong $A$-module.
\end{prop}

\begin{proof}
    The proof is similar to \cite{HA} Theorem 7.2.2.15, specifically the $(5)\Rightarrow(3)$ implication of the proof. As in loc. cit. we prove by induction on $n$ that the map
    $$\pi_{n}(A)\otimes_{\pi_{0}(A)}\pi_{0}(M)\rightarrow\pi_{n}(M)$$
    is an isomorphism. By the transversality assumption this will automatically imply that $M$ is derived strong.

    For $n=0$ this is of course trivial. Suppose we have proven the claim for all $n\le m$. Pick a projective resolution by $\pi_{0}(A)$-modules $Q_{\bullet}\rightarrow\pi_{0}(N)$. Again by the transversality assumption we have that 
    $$\pi_{*}(A)\otimes_{\pi_{0}(A)}Q_{\bullet}\rightarrow\pi_{*}(A)\otimes_{\pi_{0}(A)}\pi_{0}(N)$$
    is an equivalence. Moreover 
    $$\pi_{*}(A)\otimes_{\pi_{0}(A)}\pi_{0}(N)\rightarrow \pi_{*}(N)$$
    is an isomorphism
    in graded degrees $\le m$. Extend in higher degrees to a projective resolution of graded $\pi_{*}(A)$-modules
    $$P_{\bullet,*}\rightarrow\pi_{*}(N)$$
    with $P_{\bullet,k}=\pi_{k}(A)\otimes_{\pi_{0}(A)}Q_{\bullet}$ for $k\le m$.

    Now using this resolution, we see that 
    $$\mathrm{Tor}^{\pi_{*}(A)}(\pi_{0}(A),\pi_{*}(N))_{q}$$
    is zero unless $p=q=0$ or $q\ge m$. By analysing the Tor spectral sequence as in the proof of \cite{HA} Theorem 7.2.2.15, it suffices to observe that $\pi_{i+n}(\pi_{0}(A)\otimes^{\mathbb{L}}_{A}M)\cong\pi_{i+n}(\pi_{0}(M))$ vanishes unless $i=n=0$.
\end{proof}

A similar (in fact easier) inductive proof gives the following.

\begin{prop}
  Let $A\in\mathbf{Ass}^{cn}(\mathbf{C})$ and let $M$ be a connective left $A$-module. If $\pi_{0}(A)\otimes_{A}^{\mathbb{L}}M\cong 0$ then $M\cong 0$.
\end{prop}

%

\subsubsection{Projectivity}

\begin{defn}
Let $A\in\mathbf{Ass}^{cn}(\mathbf{C})$. An $A$-module $M$ is said to be 
\begin{enumerate}
\item
$P$-\textit{projective} for $P\in\mathbf{C}_{\ge0}$ a projective object, if $M$ is a retract in $\mathbf{Mod}_{A}$ of $A\otimes P$.
\item
\textit{finite projective} if it is $\mathbb{I}^{\oplus n}$-projective for some $n\in\mathbb{N}_{0}$. 
\end{enumerate}
\end{defn}

\begin{lem}[\cite{toen2008homotopical} Lemma 2.2.2.2]
Let $A\in\mathbf{Ass}^{cn}(\mathbf{C})$, and let $M$ be a connective $A$-module. If $P$ is strong then $M$ is $P$-projective if and only if it is derived strong and $\pi_{0}(M)$ is $\pi_{0}(P)$-projective as a $\pi_{0}(A)$-module.
\end{lem}

\begin{proof}
Let $M$ be a projective $A$-module. Then $M$ is a retract of an object of the form $A\otimes^{\mathbb{L}} P$, where $P$ is projective in $\mathbf{C}_{\ge0}$, and $\pi_{0}(M)$ is a retract of $\pi_{0}(A\otimes^{\mathbb{L}} P)\cong\pi_{0}(A)\otimes \pi_{0}(P)\cong\pi_{0}(A)\otimes^{\mathbb{L}}\pi_{0}(P)$, where we have used that $P$ is homotopy flat. Moreover retracts of strong modules are strong, and free modules on flat objects are strong. Thus $M$ is strong.

Conversely, let $M$ be a strong $A$-module with $\pi_{0}(M)$ a projective $\pi_{0}(A)$-module. Then $\pi_{0}(M)$ is a retract of $\pi_{0}(A)\otimes P\cong\pi_{0}(A\otimes^{\mathbb{L}} P)$ for some projective $P\in\mathbf{C}^{0}$, i.e. there are maps $r:\pi_{0}(A\otimes P)\rightarrow\pi_{0}(M)$ and $i:\pi_{0}(M)\rightarrow\pi_{0}(A\otimes P)$ such that $p=r\circ i$ is the identity. We have
$$\pi_{0}\mathbf{Map}_{\mathbf{Mod}_{A}}(A\otimes^{\mathbb{L}} P,A\otimes^{\mathbb{L}} P)\cong\Hom_{\mathrm{Mod}_{\pi_{0}(A)}}(\pi_{0}(A)\otimes \pi_{0}(P),\pi_{0}(A)\otimes \pi_{0}(P))$$
Thus there is a map $p':A\otimes^{\mathbb{L}} P\rightarrow A\otimes^{\mathbb{L}} P$ such that $p'\circ p'=p'$, and $\pi_{0}(p')=p$. Similarly, there is a map $r':A\otimes^{\mathbb{L}} P\rightarrow M$ such that $\pi_{0}(r')=r$.
We get a split fibre sequence
$$K\rightarrow A\otimes^{\mathbb{L}} P\rightarrow C$$
where $K$ is the fibre of $p'$. $r'$ induces a map $r'':C\rightarrow M$ which is an isomorphism on $\pi_{0}$. As a retract of a free $A$-module on a strong object, $C$ is strong. Since $M$ is strong $r''$ is an equivalence on each $\pi_{n}$, and hence is an equivalence.
\end{proof}


\begin{prop}
Let $A\in\mathbf{Ass}^{cn}(\mathbf{C})$, $M$ be an $A$-module, $P$ a strong projective $A$-module, and $f:Q\rightarrow\pi_{0}(M)$ a map with $Q$ being $\pi_{0}(P)$-projective as a $\pi_{0}(A)$-module. Then there is a $P$-projective $A$-module $\tilde{Q}$  and a map $\tilde{f}:\tilde{Q}\rightarrow M$ such that $\pi_{0}(A)\otimes^{\mathbb{L}}_{A}\tilde{f}\cong f$.
\end{prop}

\begin{proof}
The proof of the previous proposition shows that there is an endomorphism $p':A\otimes^{\mathbb{L}} P\rightarrow A\otimes^{\mathbb{L}}P$ such that $p'\circ p'=p'$ and $\pi_{0}(p')=p$. This again gives a split fibre sequence
$$ K\rightarrow A\otimes^{\mathbb{L}} P\rightarrow C$$
and $C$ works as a choice of $\tilde{Q}$. Now identically to the proof of the previous proposition, we have
$$\pi_{0}\mathbf{Map}_{\mathbf{Mod}_{A}}(A\otimes^{\mathbb{L}} P,M)\cong\Hom_{\mathrm{Mod}_{\pi_{0}(A)}}(\pi_{0}(A)\otimes \pi_{0}(P),\pi_{0}(M))$$
Thus the map $\pi_{0}(A)\otimes \pi_{0}(P)\rightarrow Q\rightarrow \pi_{0}(M)$ gives rise to a map $A\otimes^{\mathbb{L}} P\rightarrow M$. Define $\tilde{f}$ to be the composite
$$\tilde{Q}\rightarrow A\otimes P\rightarrow M$$

\end{proof}

\subsubsection{Spectral (Pre-)Contexts of Modules}

Here we show that modules over algebras internal to spectral pre-contexts give rise to further examples of spectral pre-contexts.

\begin{lem}[(c.f. \cite{lurieDAG9} Lemma 12.11))]\label{lem:qcohfilt}
Let $A\in\mathbf{Ass}^{cn}(\mathbf{C})$ and let $F\in{}_{A}\mathbf{Mod}^{cn}$. Then there exists a convergent filtration
$$0=F(-1)\rightarrow F(0)\rightarrow F(1)\rightarrow\ldots$$
on $F$ such that 
\begin{enumerate}
\item
For each $i\ge 0$ the fibre of the map $F(i-1)\rightarrow F(i)$ is equivalent to an object of the form $A\otimes^{\mathbb{L}}P[i]$ with $P$ projective.
\item
For $0\le i\le m$ the fibre of the map $F(i)\rightarrow F$ is $i$-connective.
\end{enumerate}
\end{lem}

\begin{proof}
The proof is by induction, taking $F(-1)=0$.  Suppose that we have provided a sequence $$0=F(-1)\rightarrow F(0)\rightarrow F(1)\rightarrow \cdots F(m)\rightarrow F$$
satisfying the conditions up to level $m$. Let $F'$ be the fibre of $F(m)\rightarrow F$. There is an exact sequence
$$\pi_{m+1}(F(m))\rightarrow\pi_{m+1}(F)\rightarrow\pi_{m}(F')\rightarrow\pi_{m}(F(m))\rightarrow\pi_{m}(F)$$
Pick an epimorphism $\pi_{0}(A)\otimes\pi_{0}(P)\rightarrow\pi_{m}(F')$. Consider the induced map
$$A\otimes^{\mathbb{L}}P[m]\rightarrow F'$$
which by construction is an epimorphism on $\pi_{m}$. Let $F(m+1)$ denote the cofibre of the composite $A\otimes^{\mathbb{L}}P[m]\rightarrow F'\rightarrow F(m)$. An easy check using long exact sequences shows that $F(m+1)$ satisfies the required conditions. Moreover it is clear that the map
$$\colim F(m)\rightarrow F$$
is an equivalence.
\end{proof}

\begin{lem}
Let $(\mathbf{C},\mathbf{C}_{\ge0},\mathbf{C}_{\le0},\mathbf{C}^{0})$ be a spectral algebraic pre-context, and let $A\in\mathbf{Comm}(\mathbf{C}_{\ge0})$. Then there is a $t$-structure on ${}_{A}\mathbf{Mod}$ whereby $M\in{}_{A}\mathbf{Mod}_{\ge0}$ (resp. $M\in{}_{A}\mathbf{Mod}_{\ge0}$) precisely if the underlying object of $M$ is in $\mathbf{C}_{\ge0}$ (resp. $\mathbf{C}_{\le0}$). Moreover with this $t$-structure
$$({}_{A}\mathbf{Mod},{}_{A}\mathbf{Mod}_{\ge0},{}_{A}\mathbf{Mod}_{\le0},{}_{A}\mathbf{Mod}^{0})$$
is a spectral algebraic pre-context, where ${}_{A}\mathbf{Mod}^{0}$ is the full subcategory spanned by objects of the form $A\otimes^{\mathbb{L}}P$, for $P\in\mathbf{C}^{0}$. It is symmetric if  $(\mathbf{C},\mathbf{C}_{\ge0},\mathbf{C}_{\le0},\mathbf{C}^{0})$ is. Moreover the functor $\pi_{0}$ determines an equivalence between ${}_{A}\mathbf{Mod}^{\heart}$ and ${}_{\pi_{0}(A)}\mathrm{Mod}(\mathbf{C}^{\heart})$.
\end{lem}

\begin{proof}
Using our Lemma \ref{lem:qcohfilt}, the proof of \cite{HA} Proposition 7.1.1.13 generalises to show that the $t$-structure is well-defined and is complete, and that $\pi_{0}$ determines an equivalence between ${}_{A}\mathbf{Mod}^{\heart}$ and ${}_{\pi_{0}}\mathrm{Mod}(\mathbf{C}^{\heart})$. ${}_{A}\mathbf{Mod}^{0}$ generates ${}_{A}\mathbf{Mod}_{\ge0}$ under sifted colimits by \cite{HA} Proposition 4.7.3.14. We have
\begin{equation}
    \begin{split}
\pi_{*}(A\otimes^{\mathbb{L}}P) & \cong\pi_{*}(A)\otimes_{\pi_{*}(R)}\pi_{*}(P)\cong\pi_{*}(A)\otimes_{\pi_{0}(R)}\pi_{0}(P)\\ & \cong\pi_{*}(A)\otimes_{\pi_{0}(A)}\pi_{0}(A)\otimes_{\pi_{0}(R)}\pi_{0}(P) \\ & \cong\pi_{*}(A)\otimes_{\pi_{0}(A)}\pi_{0}(A\otimes^{\mathbb{L}}P)
    \end{split}
\end{equation}
Thus the projective generators are strong. Moreover this also shows that $\pi_{*}(A\otimes^{\mathbb{L}}P)$ is a projective $\pi_{*}(A)$-module. We also have
\begin{equation}
\begin{split}\pi_{0}((A\otimes^{\mathbb{L}}P)\otimes^{\mathbb{L}}_{A}(A\otimes^{\mathbb{L}}Q)) & \cong\pi_{0}(A\otimes^{\mathbb{L}}(P\otimes^{\mathbb{L}}Q))\cong\pi_{0}(A)\otimes\pi_{0}(P)\otimes\pi_{0}(Q) \\ & \cong(\pi_{0}(A)\otimes\pi_{0}(P))\otimes_{\pi_{0}(A)}\otimes(\pi_{0}(A)\otimes\pi_{0}(Q)) \\ & \cong\pi_{0}(A\otimes^{\mathbb{L}} P)\otimes_{\pi_{0}(A)}\otimes(\pi_{0}(A\otimes^{\mathbb{L}} Q))
\end{split}
\end{equation}
A similar computation shows that projectives are homotopy flat.
That ${}_{A}\mathbf{Mod}^{0}$ contains the unit and is closed under coproducts is clear. ${}_{A}\mathbf{Mod}_{\le0}$ is closed under filtered colimits by definition of the $t$-structure, and since the forgetful functor commutes with filtered colimits.

Finally if the ground (pre)-context is symmetric, and $P\in\mathbf{C}^{0}$ then 
$$\mathrm{Sym}_{\pi_{0}(A)}^{n}(\pi_{0}(A\otimes^{\mathbb{L}}P))\cong \pi_{0}(A)\otimes\mathrm{Sym}^{n}(\pi_{0}(P))\cong\pi_{0}(A\otimes^{\mathbb{L}}Q)$$
where $Q\in\mathbf{C}^{0}$ is such that $\pi_{0}(Q)\cong\mathrm{Sym}^{n}(\pi_{0}(P))$. Here $\mathrm{Sym}^{n}_{\pi_{0}(A)}$ denotes the symmetric power in ${}_{\pi_{0}(A)}\mathrm{Mod}$.
\end{proof}

\subsubsection{Eilenberg-Maclane for spectral algebraic pre-contexts}

Let \[(\mathbf{C},\mathbf{C}_{\ge0},\mathbf{C}_{\le0},\mathbf{C}^{0})
\] be a spectral algebraic pre-context.

%
%

\begin{lem}
If projectives are homotopy flat, and $\mathbb{I}$ is projective then there is an exact equivalence of stable categories with $t$-structure.
$$\mathbf{Ch}(\mathbf{C}^{\heart})\cong\Mod_{\pi_{0}(\mathbb{I})}(\mathbf{C})$$
\end{lem}

\begin{proof}
$\mathbf{\Mod}_{\pi_{0}(R)}(\mathbf{C})_{\ge0}$ is freely generated under sifted colimits by the full subcategory $\{\pi_{0}(R)\otimes_{R}P:P\in\mathbf{C}^{0}\}$. Now we have that
$$\pi_{0}(P)\cong\pi_{0}(R)\otimes_{\pi_{*}(R)}\pi_{*}(P)\cong\pi_{*}(\pi_{0}(R)\otimes^{\mathbb{L}}_{R}P)$$
where we have used that $P$ is homotopy flat.
$\pi_{0}(P)$ is compact. Moreover by assumption on the fact that $\pi_{0}(P)$ is projective in $\mathbf{C}^{\heart}$ and the fact that $K_{\ge0}$ commutes with sifted colimits, we have that $\pi_{0}(P)$ is projective in $\Mod_{\pi_{0}(R)}(\mathbf{C}_{\ge0})$. 
\end{proof}

\begin{cor}
If the map $R\rightarrow\pi_{0}(R)$ is an equivalence then the adjunction
$$\adj{\Gamma}{\mathbf{C}}{\mathbf{Ch}(\mathbf{C}^{\heart})}{K}$$
is an equivalence. 
\end{cor}

This in particular happens if $\mathbf{C}^{0}\subset\mathbf{C}^{\heart}$. We will return to this situation shortly, in Section \ref{subsec:deralg}

%
%
%

\subsection{Coherent and Noetherian Algebras}

 Fix a  spectral algebraic pre-context. $(\mathbf{C},\mathbf{C}_{\ge0},\mathbf{C}_{\le0},\mathbf{C}^{0})$. Here we discuss coherence and Noeherian properties of rings in general abelian categories.

\subsubsection{Compactly and Finitely Presented Modules}

\begin{defn}
Let $A\in\mathbf{Ass}(\mathbf{C}^{\heart})$. $M\in{}_{A}\mathrm{Mod}(\mathbf{C}^{\heart})$ is said to be
\begin{enumerate}
\item
\textit{compactly generated} if there is an epimorphism $A\otimes P\rightarrow M$ for some compact projective $P\in\mathbf{C}^{\heart}$.
\item
\textit{compactly presented} if there is an epimorphism $A\otimes P\rightarrow M$ for some compact projective $P\in\mathbf{C}^{\heart}$ such that $\mathrm{Ker}(A\otimes P\rightarrow M)$ is compactly generated.
\item
\textit{finitely generated} if there is an epimorphism $A^{\oplus n}\rightarrow M$ for some $n\ge0$.
\item
\textit{finitely presented} if there is an epimorphism $A^{\oplus n}\rightarrow M$ for some $n\ge0$, and $\mathrm{Ker}(A^{\oplus n}\rightarrow M)$ is finitely generated.
\item
\textit{coherent} if it is finitely generated, and every finitely generated submodule is finitely presented.
\item
\textit{Noetherian} if every submodule is finitely generated.
\end{enumerate}
$A$ is said to be \textit{left coherent (resp. left Noetherian)} if it is coherent (resp. Noetherian) as a left $A$-module over itself.
\end{defn}

Exactly as in the standard result for usual algebras (see 
\cite[\href{https://stacks.math.columbia.edu/tag/0517}{Tag 0517}]{stacks-project}
) we have the following.

\begin{lem}
Let $A\in\mathbf{Ass}(\mathbf{C}^{\heart})$. $M\in{}_{A}\mathrm{Mod}(\mathbf{C}^{\heart})$, and let 
$$0\rightarrow M_{1}\rightarrow M_{2}\rightarrow M_{3}\rightarrow 0$$
be an exact sequence of left $A$-modules. Then
\begin{enumerate}
\item
If $M_{1}$ and $M_{3}$ are compactly generated/ finitely generated/ compactly presented/ finitely presented, then so is $M_{2}$.
\item
if $M_{2}$ is compactly generated/ finitely generated, then so is $M_{3}$.
\item
If $M_{2}$ is compactly presented/ finitely presented and $M_{1}$ is compactly generated/ finitely generated, then $M_{3}$ is compactly presented/ finitely presented.
\item
If $M_{3}$ is compactly presented/ finitely presented, and $M_{2}$ is compactly generated/ finitely generated, then $M_{1}$ is compactly generated/ finitely generated.
\end{enumerate}
\end{lem}

By arguments identical to the standard ones for usual coherent rings (see. e.g. \cite[\href{https://stacks.math.columbia.edu/tag/05CU}{Tag 05CU}]{stacks-project}) we have the following results.

\begin{lem}\label{lem:factscoherent}
Let $A\in\mathbf{Ass}(\mathbf{C}^{\heart})$.
\begin{enumerate}
\item
Any finitely generated submodule of a coherent left $A$-module is coherent.
\item
The kernel and cokernel of a morphism between coherent left $A$-modules is coherent.
\item
The category of coherent left $A$-modules is an extension closed subcategory of ${}_{A}\mathrm{Mod}(\mathbf{C}^{\heart})$.
\item
Any finite sum of coherent left $A$-modules is coherent.
\item
If $A$ is left coherent then a left $A$-module is coherent if and only if it is finitely presented.
\end{enumerate}
\end{lem}

\begin{cor}
Let $A$ be a left coherent ring and $I\subset A$ a finitely generated left ideal. Then $A\big\slash I$ is left coherent.
\end{cor}

\begin{proof}
This follows immediately from the fact that $A\big\slash I$ is coherent as a left $A$-module.
\end{proof}

Let us give a useful source of examples arising from quasi-abelian categories.

\begin{example}
Let $\mathpzc{E}$ be a monoidal elementary quasi-abelian category. Then its left heart $\mathrm{LH}(\mathpzc{E})$ is a monoidal elementary abelian category. Let $A\in\mathrm{Ass}(\mathpzc{E})\subset\mathrm{Ass}(\mathrm{LH}(\mathpzc{E})$. By \cite{qacs} Proposition 1.4.15 ${}_{A}\mathrm{Mod}(\mathpzc{E})\subset{}_{A}\mathrm{Mod}(\mathrm{LH}(\mathpzc{E}))$ is essentially stable by subobjects. Suppose that $A$ is such that if $I\rightarrow A$ is a monomorphism and there exists a strict epimorphism $A^{\oplus n}\rightarrow I$ for some $n\in\mathbb{N}_{0}$ then 
\begin{enumerate}
\item
$I\rightarrow A$ is a strict monomorphism. 
\item
there is a strict exact sequence
$$R^{m}\rightarrow R^{n}\rightarrow I\rightarrow 0$$
for some $m$.

Specifically in this context, we will say that an algebra $A\in\mathrm{Ass}(\mathpzc{E})$ is \textit{strongly left coherent} if it is left coherent, and whenever $I\rightarrow A$ is a monomorphism in ${}_{A}\mathrm{Mod}(\mathpzc{E})$ such that there exists an admissible epimorphism $A^{\oplus n}\rightarrow I$ for some $n$, then $I\rightarrow A$ is an admissible monomorphism.
\end{enumerate}
Then $A$ is left coherent as an object of $\mathrm{Ass}(\mathrm{LH}(\mathpzc{E}))$. 

If whenever $I\rightarrow A$ is a monomorphism \textit{then} it is a strict monomorphism and there is an admissible epimorphism $A^{\oplus n}\rightarrow I$, then $A$ is Noetherian. 

The algebras of interest in analytic geometry will in typically be coherent - this will essentially amount to them having the property that ideals are finitely generated precisely if they are closed, and the fact that the underlying algebra is coherent in the usual sense.
\end{example}

\subsubsection{Presentations of Modules}

An important application of the theory of coherent algebras is to colimits of algebras. This will be particularly useful when discussing algebras appearing in analytic geometry later.

\begin{lem}
Let $\mathrm{C}$ be a monoidal elementary abelian category, $\mathcal{I}$ a filtered category, and $F:\mathcal{I}\rightarrow\mathrm{Ass(C)}$ a diagram. Write $R_{i}=F(i)$ and $R\defeq\limind_{\mathcal{I}}R_{i}$.
\begin{enumerate}
\item
Let $M$ be an $R$-module. Then

\begin{enumerate}
\item
for each $j\in\mathcal{I}$ there is an $R_{j}$-module $M_{j}$
\item
for each $k\ge j$ there is  a map of $R_{j}$-modules $f_{j,k}:M_{j}\rightarrow M_{k}$
\end{enumerate}
such that 
\begin{enumerate}
\item
for each $l\ge k\ge j$ we have $f_{j,l}=f_{k,l}\circ f_{j,k}$.
\item
there is an isomorphism of $R$-modules $\limind_{\mathcal{I}}M_{j}\cong M$.
\end{enumerate}
Moreover if $M$ is finitely generated then each $M_{j}$ can be chosen to be finitely generated.
\item
Let $M\rightarrow N$ be a map of left $R$-modules. Suppose $M$ has a presentation $M\cong\limind_{\mathcal{I}}M_{i}$ with each $M_{i}$ a compact $R_{i}$-module. Then there is a cofinal functor $F:\mathcal{K}\rightarrow\mathcal{I}$ such that $M\cong\colim M_{k}$, $N\cong\colim N_{k}$, each $M_{k}$ and $N_{k}$ are $R_{F(k)}$-modules, and $M\rightarrow N$ is induced by maps of $R_{F(k)}$-modules $M_{k}\rightarrow N_{k}$. If each $R_{i}$ is coherent $M$ and $N$ are finitely generated then each $M_{k}$ and each $N_{k}$ can be chosen to be finitely generated.
\item
Let 
$$0\rightarrow L\rightarrow M\rightarrow N\rightarrow 0$$
be an exact sequence of $R$-modules with $N$ compact and $M$ compactly generated. There is a cofinal functor $F:\mathcal{K}\rightarrow\mathcal{I}$ and a $\mathcal{K}$-indexed diagram of exact sequences
$$0\rightarrow L_{k}\rightarrow M_{k}\rightarrow N_{k}\rightarrow0$$
of $R_{k}$-modules whose colimit is the original exact sequence.  If each $R_{k}$ is coherent, and $M$ and $N$ are finitely generated, then each $L_{k},M_{k}$ and each $N_{k}$ can be chosen to be finitely generated.
\end{enumerate}
\end{lem}

\begin{proof}
\begin{enumerate}
\item
Let $R\otimes P\rightarrow M$ be an epimorphism with $P$ projective. Define $M_{i}$ to be the image of $R_{i}\otimes P\rightarrow M$ This clearly suffices. 
\item

Let $M\cong\colim M_{i}$, $N\cong \colim N_{i}$ be presentations as in the first part, with $M_{i}$ compact. Let $f:M\rightarrow N$ be a map of $R$-modules. Now since each $M_{i}$ is compact as an $R_{i}$-module, for each $i\in\mathcal{I}$ the map of $R_{i}$-modules $M_{i}\rightarrow N$ factors through some $N_{j(i)}$. 

Let $\mathcal{I}^{f}$ denote the category defined as follows. An object is a pair $(i,f_{i'})$ where $i\le i'\in\mathcal{I}$, and $f_{i'}:M_{i}\rightarrow N_{i'}$ is a map of $R_{i}$-modules such that the composition $M_{i}\rightarrow N_{i'}\rightarrow N$  coincides with the composition $M_{i}\rightarrow M\rightarrow N$. Say that $(i,f_{i'})\le (j,f_{j'})$ if $i\le j$, $i'\le j'$, and there is a commutative diagram
\begin{displaymath}
\xymatrix{
M_{i}\ar[d]^{f_{i'}}\ar[r] & M_{j}\ar[d]^{f_{j'}}\\
N_{i'}\ar[r] & M_{j'}
}
\end{displaymath}
Define
$$M_{(i,f_{i'})}\defeq M_{i}$$
and 
$$N_{(k,f_{k'})}\defeq N_{k'}$$
where here $N_{k'}$ is regarded as an $R_{k}$-module.

The functor $\mathcal{I}^{f}\rightarrow\mathcal{I}$, $(i,f_{i'})\mapsto i'$ is clearly cofinal. Thus $N\cong\limind_{\mathcal{I}^{f}}N_{(k,f_{k'})}$ and $M\cong\limind_{\mathcal{I}^{f}} R_{k'}\otimes_{R_{k}}M_{k}$. 
\item
Let
$$0\rightarrow L\rightarrow M\rightarrow N\rightarrow 0$$
be an exact sequence of left $R$-modules. We may assume that $L\cong\colim \overline{L}_{i}$, $M\cong\colim  M_{i}$ with each $\overline{L}_{i}$ and each $M_{i}$ being an $R_{i}$-module, and that $L\rightarrow M$ is induced by maps of $R_{i}$ modules $\overline{L}_{i}\rightarrow M_{i}$. Define $N_{i}\defeq\coker (\overline{L}_{i}\rightarrow M_{i})$, and $L_{i}\defeq\mathrm{Ker}(M_{i}\rightarrow N_{i})$. This proves the claim since in an elementary, and hence Grothendieck, abelian category, filtered colimits commute with kernels.
\end{enumerate}
\end{proof}

\begin{lem}[c.f. \cite{364690}]
Let $\mathrm{C}$ be a monoidal elementary abelian category, $\mathcal{I}$ a filtered category, and $F:\mathcal{I}\rightarrow\mathrm{Ass(C)}$ a diagram such that each $R_{i}\defeq F(i)$ is left coherent. If each map $R_{i}\rightarrow R$ is transverse to finitely presented $R_{i}$-modules. Then $R$ is coherent.
\end{lem}

\begin{proof}
Let $I\subseteq R$ be a submodule which is finitely generated. In particular there is an epimorphism $R^{\oplus n}\rightarrow I$. Consider the composite $R^{\oplus n}\rightarrow I \rightarrow R$. Note that $\Hom(R^{\oplus n},R)\cong\Hom(R,R)^{\oplus n}\cong\Hom(\mathbb{I},R^{\oplus n})$ where $\mathbb{I}$ is the unit of $\mathrm{C}$. Since $\mathbb{I}$ is tiny there is some $R_{j}$ such that the map in $\Hom(\mathbb{I},R^{\oplus n})$ corresponding to $R^{\oplus n}$ factors through some $R_{j}^{\oplus n}$. Consider the corresponding map $R_{j}^{\oplus n}\rightarrow R_{j}$. Denote by $I_{j}\subseteq R_{j}$ the image of this map. This is a finitely generated left sub-module of $R_{j}$. Thus there is an exact sequence
$$R_{j}^{\oplus m}\rightarrow R_{j}^{\oplus n}\rightarrow I_{j}\rightarrow A\big\slash I_{j}\rightarrow0$$
We have $I\cong R\otimes_{R_{j}}I_{j}$, so there is an exact sequence
$$R^{\oplus m}\rightarrow R^{\oplus n}\rightarrow I\rightarrow A\big\slash I\rightarrow 0$$
Thus $I$ is finitely presented, as required.
\end{proof}

\subsubsection{Derived Coherent Algebras}

Let us now consider derived analogues of coherent rings.

\begin{defn}
Let $A\in\mathbf{Ass}(\mathbf{C})$. $A$ is said to be \textit{left coherent} if $\pi_{0}(A)$ is left coherent and each $\pi_{*}(A)$ is finitely generated as a left $\pi_{0}(A)$-module.
\end{defn}

\begin{prop}\label{prop:leftcohtens}
Let $A,B\in\mathbf{Comm}^{cn}(\mathbf{C})$ be left coherent with either $A$ or $B$ homotopy flat. Then $A\otimes^{\mathbb{L}}B$ is left coherent if and only if $\pi_{0}(A)\otimes\pi_{0}(B)$ is left coherent.
\end{prop}

\begin{proof}
We have
$$\pi_{*}(A\otimes^{\mathbb{L}}B)\cong\pi_{*}(A)\otimes\pi_{*}(B)$$
and the result is immediate.
\end{proof}

\begin{lem}[\cite{lurieDAG9} Lemma 12.11]\label{lem:cohfilt}
Let $A\in\mathbf{Ass}^{cn}(\mathbf{C})$ be such that $\pi_{0}(A)$ is left coherent, and each $\pi_{m}(A)$ is a finitely generated $\pi_{0}(A)$-module. Let $F\in{}_{A}\mathbf{Mod}^{cn}$ be such that each $\pi_{m}(F)$ is finitely generated as a $\pi_{0}(A)$-module. Then there exists a convergent filtration
$$0=F(-1)\rightarrow F(0)\rightarrow F(1)\rightarrow\ldots$$
on $F$ such that 
\begin{enumerate}
\item
For each $i\ge 0$ the fibre of the map $F(i-1)\rightarrow F(i)$ is equivalent to a direct sum of finitely many copies of $A[i]$.
\item
For $0\le i\le m$ the fibre of the map $F(i)\rightarrow F$ is $i$-connective.
\item
For $-1\le i\le m$ the homotopy groups $\pi_{j} F(i)$ are coherent $\pi_{0}(A)$-modules.
\end{enumerate}
\end{lem}

\begin{proof}
The proof is again by induction, taking $F(-1)=0$. In fact the proof is essentially identical to the previous one - we just have to note that at each stage, because of the coherence conditions, we may take $P=R^{\oplus n}$ for some $n$. Suppose that we have provided a sequence $$0=F(-1)\rightarrow F(0)\rightarrow F(1)\rightarrow\ldots F(m)\rightarrow F$$
satisfying the conditions up to level $m$. Let $F'$ be the fibre of $F(m)\rightarrow F$. There is an exact sequence
$$\pi_{m+1}(F(m))\rightarrow\pi_{m+1}(F)\rightarrow\pi_{m}(F')\rightarrow\pi_{m}(F(M))\rightarrow\pi_{m}(F)$$
Thus $\pi_{m}(F')$ is finitely generated as a $\pi_{0}(A)$-module. Consider the induced map
$$A^{\oplus n}[m]\rightarrow F'$$
which by construction is an epimorphism on $\pi_{m}$. Let $F(m+1)$ denote the cofibre of the composite $A^{n}[m]\rightarrow F'\rightarrow F(m)$. An easy check using long exact sequences shows that $F(m+1)$ satisfies the required conditions. Moreover it is clear that the map
$$\colim F(m)\rightarrow F$$
is an equivalence.
\end{proof}

\begin{cor}\label{cor:sthngleftcoh}
Let $f:A\rightarrow B$ be a map in $\mathbf{Alg_{D}}^{cn}(\mathbf{C})$ with $\pi_{0}(A)$ and $\pi_{0}(B)$ both left coherent, and each $\pi_{m}(A)$ and $\pi_{m}(B)$ being finitely generated modules over $\pi_{0}(A)$ and $\pi_{0}(B)$ respectively. Let $F$ be a connective $A$-module such that each $\pi_{m}(F)$ is finitely generated as a $\pi_{0}(A)$-module. Then each $\pi_{k}(B\otimes^{\mathbb{L}}_{A}F)$ is finitely generated as a $\pi_{0}(B)$-module.
\end{cor}

\begin{proof}
Let 
$$0=F(-1)\rightarrow F(0)\rightarrow F(1)\rightarrow\cdots$$
be a filtration of $F$ as in Lemma \ref{lem:cohfilt}. Note that the fibre of 
$$B\otimes^{\mathbb{L}}_{A}F(k+1)\rightarrow B\otimes^{\mathbb{L}}_{A}F$$
is $k+1$-connective. Thus it suffices to prove that each $\pi_{k}(B\otimes^{\mathbb{L}}_{A}F(k+1))$ is finitely generated as a $\pi_{0}(B)$-module for $-1\le i\le k+1$. The case $i=-1$ is obvious, and the inductive step follows from the long-exact sequence associated to the existence of a fibre sequence
$$B^{n}[i]\rightarrow B\otimes^{\mathbb{L}}_{A}F(i)\rightarrow B\otimes^{\mathbb{L}}_{A}F(i+1)$$
\end{proof}

For $A\in\mathbf{DAlg}^{cn}(\mathbf{C})$ denote by ${}_{A}\mathbf{Mod}_{+}^{h-fg}(\mathbf{C})$ the class of connective $A$-modules $M$ such that $\pi_{m}(M)$ is a finitely generated $\pi_{0}(A)$-module for each $m$. An argument using long exact sequences gives the following.

\begin{prop}\label{prop:cohthick}
    Let $A\in\mathbf{DAlg}^{cn}(\mathbf{C})$ be such that $\pi_{0}(A)$ is left coherent. The class ${}_{A}\mathbf{Mod}_{+}^{h-fg}(\mathbf{C})$ is thick, and closed under finite limits and colimits.
\end{prop}

\begin{cor}
Let $A,B,C\in\mathbf{Comm}^{cn}(\mathbf{C})$ be left coherent. Suppose further that
\begin{enumerate}
\item
$A$ and one of either $B$ or $C$ is homotopy flat.
 \item
  $\pi_{0}(A)\otimes \pi_{0}(A)$, $\pi_{0}(B)\otimes \pi_{0}(C)$, and $\pi_{0}(B)\otimes_{\pi_{0}(A)}\pi_{0}(C)$ are all left coherent.
 \end{enumerate}
 Then $B\otimes_{A}^{\mathbb{L}}C$ is left coherent.
\end{cor}

\begin{proof}
We have $B\otimes^{\mathbb{L}}_{A}C\cong A\otimes^{\mathbb{L}}_{A\otimes^{\mathbb{L}}A}B\otimes^{\mathbb{L}}C$. Now $A\otimes^{\mathbb{L}}A$ and $B\otimes^{\mathbb{L}}C$ are left coherent by Proposition \ref{prop:leftcohtens}. $\pi_{*}(A)$ is a retract of $\pi_{*}(A\otimes^{\mathbb{L}}A)\cong\pi_{*}(A)\otimes\pi_{*}(A)$ and is therefore finitely generated as a $\pi_{0}(A)\otimes\pi_{0}(A)$-module. Now the claim follows from Corollary \ref{cor:sthngleftcoh}.
\end{proof}

%
%
%


\section{Spectral Algebraic Contexts}\label{sec:spac}

Spectral algebraic contexts combine weak spectral algebraic contexts and Postnikov algebraic contexts.

\begin{defn}
\begin{enumerate}
\item
A \textit{spectral algebraic context} is a tuple
$$(\mathbf{C},\mathbf{C}_{\ge0},\mathbf{C}_{\le0},\mathbf{C}^{0},\mathbf{D},\theta,\mathbf{G}_{0},\mathbf{S})$$
where
\begin{enumerate}
\item
$(\mathbf{C},\mathbf{C}_{\ge0},\mathbf{C}_{\le0},\mathbf{C}^{0})$ is a spectral algebraic pre-context.
\item
$(\mathbf{C},\mathbf{C}_{\ge0},\mathbf{C}_{\le0},\mathbf{D},\theta,\mathbf{G}_{0},\mathbf{S})$ is a Postnikov algebraic context.
\end{enumerate}
\item
A \textit{connective spectral algebraic context} is a tuple
$$(\mathbf{C},\mathbf{C}_{\ge0},\mathbf{C}_{\le0},\mathbf{C}^{0},\mathbf{D},\theta,\mathbf{G}_{0},\mathbf{S})$$
where
\begin{enumerate}
\item
$(\mathbf{C},\mathbf{C}_{\ge0},\mathbf{C}_{\le0},\mathbf{C}^{0},\mathbf{D},\theta)$ is a spectral algebraic pre-context.
\item
$(\mathbf{C},\mathbf{C}_{\ge0},\mathbf{C}_{\le0},\mathbf{D},\theta,\mathbf{G}_{0},\mathbf{S})$ is a connective Postnikov algebraic context.
\end{enumerate}
\end{enumerate}
A (connective) spectral algebraic context is said to be a \textit{Koszul (connective) spectral algebraic context} if the underlying (connective) Postnikov algebraic context is Koszul.
\end{defn}

\begin{defn}
A \textit{transformation} from a spectral algebraic context 
$$(\mathbf{C},\mathbf{C}_{\ge0},\mathbf{C}_{\le0},\mathbf{C}^{0},\mathbf{D},\theta,\mathbf{G}_{0},\mathbf{S})$$
 to a spectral algebraic context 
 $$(\mathbf{C}',\mathbf{C}'_{\ge0},\mathbf{C}'_{\le0},(\mathbf{C}')^{0},\mathbf{D}',\theta',\mathbf{G}'_{0},\mathbf{S}')$$
 is a transformation $(\mathbf{F},\eta,\delta)$ from $(\mathbf{C},\mathbf{D},\theta,\mathbf{G}_{0},\mathbf{S})$ to $(\mathbf{C}',\mathbf{D'},\theta',\mathbf{G}'_{0},\mathbf{S}')$, in the sense of Definition \ref{defn:transformation}, such that $\mathbf{F}(\mathbf{C}^{0})\subset(\mathbf{C}')^{0}$. 
\end{defn}

\subsection{(Weak) Derived Algebraic Contexts and Derived Commutative Rings}\label{subsec:deralg}

At this point it is convenient to recall the definition of a \textit{derived algebraic context} from \cite{raksit2020hochschild}. We will give a definition in our language of spectral algebraic pre-contexts which is evidently equivalent.

\begin{defn}
A (weak) spectral algebraic pre-context $(\mathbf{C},\mathbf{C}_{\ge0},\mathbf{C}_{\le0},\mathbf{C}^{0})$ is said to be an $E_{\infty}$-\textit{derived algebraic context} if
\begin{enumerate}
\item
$\mathbf{C}^{0}\subset\mathbf{C}^{\heart}$
\item
Each $P\in\mathbf{C}^{0}$ is projective in $\mathbf{C}_{\ge0}$
\end{enumerate}
If in addition it is symmetric as a weak spectral algebraic pre-context, then it is said to be a \textit{ derived algebraic context}.
\end{defn}

\begin{rem}
\begin{enumerate}
\item
$E_{\infty}$-derived algebraic contexts are in particular spectral algebraic pre-contexts.
\item
Let $(\mathbf{C},\mathbf{C}_{\ge0},\mathbf{C}_{\ge0},\mathbf{C}^{0})$ be an $E_{\infty}$-derived algebraic context. Any projective $P$ in $\mathbf{C}^{\heart}$ is projective in the sense of Definition \ref{defn:stableproj}. This follows immediately since $R\cong\pi_{0}(R)$ is the unit of the monoidal structure on $\mathbf{C}^{\heart}$. 
\end{enumerate}
\end{rem}

In \cite{raksit2020hochschild}, following \cite{brantner2019deformation}, Raksit constructs categories of both non-connective and connective \textit{derived commutative rings}, which codify simplicial commutative rings in $(\infty,1)$-categorical language as follows. Consider the symmetric algebra functor
$$\mathrm{Sym}:\mathbf{C}^{0}\rightarrow\mathbf{C}^{\heart}$$
This extends by sifted colimits to a (monadic) functor
$$\mathbf{LSym}:\mathbf{C}_{\ge0}\rightarrow\mathbf{C}_{\ge0}$$
There are similar constructions for derived exterior powers $\mathbf{L}\Lambda(-)$ and derived divided powers $\mathbf{L}\Gamma(-)$, which are particularly well-behaved in the case that the $1$-categorical constructions $\Lambda^{n}(P)$ and $\Gamma^{n}(P)$ are also projective for a given projective $P$.

\begin{prop}
 Let  $(\mathbf{C},\mathbf{C}_{\ge0},\mathbf{C}_{\le0},\mathbf{C}^{0})$ be a derived algebraic context. 
    If $M$ is in $\mathbf{C}_{\ge m}$ then $\mathbf{LSym}^{n}(M)$, $\mathbf{L}\Lambda^{n}(M)\in\mathbf{C}_{\ge m}$, and $\mathbf{L}\Gamma^{n}(M)$ are in $\mathbf{C}_{\ge m}$
\end{prop}

\begin{proof}
    The proof of this works as in \cite{lurie2018spectral} Proposition 25.2.4.3. We prove the $\mathbf{LSym}^{n}$ cases, the others being identical. The proof is by induction on $m$. For $m=0$ it is clear by definition. Suppose that $M$ is $m$-connective for $m>0$ and let $N_{\bullet}$ be the \v{C}ech nerve of the map $0\rightarrow M$, so that we have 
    $$\mathbf{LSym}^{n}(M)\cong|\mathbf{LSym}^{n}(N_{\bullet})|$$
    Now $N_{k}\cong(\Sigma^{-1}M)^{k}$ so is $m-1$-connective. Thus $\mathbf{LSym}^{n}(M)$ is $m-1$-connective. Moreover $\pi_{m-1}\mathbf{LSym}^{n}(M)$ is a quotient of $\pi_{m-1}(\mathbf{LSym}^{n}(N_{0}))\cong 0$.
\end{proof}

It is shown in \cite{raksit2020hochschild} Section 4.2 that $\mathbf{LSym}$ extends to a monad on all of $\mathbf{C}$, called the \textit{derived symmetric algebra monad on }$\mathbf{C}$. A module over this monad is called a \textit{derived commutative algebra object of }$\mathbf{C}$. The category of these modules is denoted $\mathbf{DAlg}(\mathbf{C})$. We also write
$$\mathbf{DAlg}^{cn}(\mathbf{C})\defeq\mathbf{DAlg}(\mathbf{C})\times_{\mathbf{C}}\mathbf{C}_{\ge0}$$
for the category of \textit{connective derived commutative algebras}. There is a functor

$$\Theta:\mathbf{DAlg}(\mathbf{C})\rightarrow\mathrm{Comm}(\mathbf{C})$$
which by \cite{raksit2020hochschild} Proposition 4.2.27 commutes with all small limits and colimits.

There is an important distinction between (connective) derived commutative rings and commutative rings internal to the symmetric monoidal $(\infty,1)$-category ($\mathbf{C}_{\ge0}$) $\mathbf{C}$. Briefly, $\mathbf{C}_{\ge0}$ is presented by the model category $\mathrm{s}\mathbf{C}^{\heart}$ of simplicial objects in $\mathbf{C}^{\heart}$, and $\mathbf{C}$ by the category $\mathrm{Sp}_{\Sigma}(\mathrm{s}\mathbf{C}^{\heart})$ of symmetric spectra internal to this category (as described in \cite{MR1860878} Section 2). Consider the commutative operad $\mathpzc{Comm}$, and a $\Sigma$-cofibrant replacement $E_{\infty}$ of 
$\mathpzc{Comm}$. As explained in the next chapter, these are admissible operads, i.e. the transferred model structures exist on both $\mathrm{Alg}_{\mathpzc{Comm}}(\mathrm{Sp}_{\Sigma}(\mathrm{s}\mathbf{C}^{\heart}))$ and $\mathrm{Alg}_{E_{\infty}}(\mathrm{Sp}_{\Sigma}(\mathrm{s}\mathbf{C}^{\heart}))$.  The Quillen adjunctions
$$\adj{\mathrm{Free}_{\mathpzc{Comm}}}{\mathrm{Sp}_{\Sigma}(\mathrm{s}\mathbf{C}^{\heart})}{\mathrm{Alg}_{\mathpzc{Comm}}(\mathrm{Sp}_{\Sigma}(\mathrm{s}\mathbf{C}^{\heart}))}{|-|_{\mathpzc{Comm}}}$$
$$\adj{\mathrm{Free}_{E_{\infty}}}{\mathrm{Sp}_{\Sigma}(\mathrm{s}\mathbf{C}^{\heart})}{\mathrm{Alg}_{E_{\infty}}(\mathrm{Sp}_{\Sigma}(\mathrm{s}\mathbf{C}^{\heart}))}{|-|_{E_{\infty}}}$$
present monadic adjunctions of $(\infty,1)$-categories. The former gives the adjunction
$$\adj{\mathbf{LSym}}{\mathbf{C}}{\mathbf{DAlg(C)}}{|-|_{\mathbf{LSym}}}$$
and the latter the adjunction
$$\adj{\mathbf{Sym}}{\mathbf{C}}{\mathbf{Comm(C)}}{|-|_{\mathbf{Sym}}}$$ 
There is an implicit map of operads 
$$\theta:E_{\infty}\rightarrow\mathpzc{Comm}$$
which gives the map
$$\mathbf{Sym}\rightarrow\mathbf{LSym}$$
If $\mathpzc{C}^{\heart}$ is enriched over $\mathbb{Q}$, then $\mathpzc{Comm}$ is already $\Sigma$-cofibrant. In this case the map
$$\mathbf{Sym}\rightarrow\mathbf{LSym}$$
is an equivalence of monads, and so the functor $\Theta:\mathbf{DAlg}(\mathbf{C})\rightarrow\mathbf{Comm}(\mathbf{C})$ is an equivalence. 


\begin{rem}\label{rem:modulesspectra}
Let $A\in\mathbf{DAlg}^{cn}(\mathbf{C})$. As pointed out in the introduction of \cite{porta2017representability},
it is generally not the case that one has an equivalence
$${}_{A}\mathbf{Mod(C)}\cong\mathbf{Stab}(\mathbf{DAlg}^{cn}(\mathbf{C})_{\big\slash A})$$
However exactly as in loc. cit. Appendix 8.1, it is possible to prove the following equivalence
$${}_{A}\mathbf{Mod(C)}\cong\mathbf{Stab}(\mathbf{Ab}(\mathbf{DAlg}^{cn}(\mathbf{C})_{\big\slash A}))$$
The only difference is in the proof of \cite{porta2017representability}
Proposition 8.1, for $A\in\mathbf{DAlg}^{cn}(\mathbf{C})$ one is required to use the pregeometry $\mathcal{T}^{\mathbf{C}}_{A}$ given by the full subcategory of $\mathbf{DAlg}^{cn}(\mathbf{C})^{op}$ on objects of the form
$$\mathrm{LSym}(A\otimes Q)$$
for $Q\in\mathcal{Q}$. 
\end{rem}

Now we have a commutative diagram

\begin{displaymath}
\xymatrix{
{}_{A}\mathbf{Mod(C)}\ar[dd]^{\mathrm{Id}}\ar[r] &\mathbf{Stab}(\mathbf{Ab}(\mathbf{DAlg}^{cn}(\mathbf{C})_{\big\slash A}))
\ar[d]\ar[r]&\mathbf{Ab}(\mathbf{DAlg}^{cn}(\mathbf{C})_{\big\slash A})\ar[d]\\
&\mathbf{Stab}(\mathbf{DAlg}^{cn}(\mathbf{C})_{\big\slash A})\ar[r]\ar[d]&\mathbf{DAlg}^{cn}(\mathbf{C})_{\big\slash A})\ar[d]
\\
{}_{\Theta(A)}\mathbf{Mod(C)}\ar[r] & \mathbf{Stab}(\mathbf{Comm}^{cn}(\mathbf{C})_{\big\slash \Theta(A)})\ar[r] & \mathbf{Comm}^{cn}(\mathbf{C})_{\big\slash \Theta(A)}
}
\end{displaymath}

where the first top and bottom horizontal maps are equivalences. It follows that $\theta_{A}:\Theta(\mathrm{sqz}_{A}^{\mathbf{LSym}})\rightarrow\mathrm{sqz}_{\Theta(A)}$ is an equivalence.

\begin{defn}\label{defn:KoszulDAC}
    A derived algebraic context  $(\mathbf{C},\mathbf{C}_{\ge0},\mathbf{C}_{\le0},\mathbf{C}^{0})$ is said to be \textit{Koszul} if whenever $P\in\mathbf{C}^{\heart}$ is projectve
    \begin{enumerate}
        \item
         for each $n$, the $n$th exterior power $\Lambda^{n}(P)$ is projective.
         \item 
         The Koszul complex
         $$\Lambda(P)\otimes\mathrm{Sym}(P)\rightarrow\mathbb{I}$$
         ia a graded resolution.
    \end{enumerate}
\end{defn}

It suffices to check this condition for $P\in\mathbf{C}^{0}$. 

\begin{lem}[Illusie, \cite{lurie2018spectral} Proposition 25.2.4.2]
    Let  $(\mathbf{C},\mathbf{C}_{\ge0},\mathbf{C}_{\le0},\mathbf{C}^{0})$ be a Koszul derived algebraic context. Then for any $A\in\mathbf{DAlg}^{cn}(\mathbf{C})$, for every connective $A$-module $M$, and for every $n\ge 0$ we have canonical equivalences
    $$\mathbf{LSym}^{n}(\Sigma M)\cong\Sigma^{n}\Lambda^{n}_{A}(M)$$
\end{lem}

\begin{proof}
    First assume that $A=R$ is discrete. Consider a short exact sequence
    \begin{displaymath}
    \xymatrix{
        0\ar[r] & M'\ar[r]^{\rho} & M\ar[r] & M''\ar[r] &0
        }
    \end{displaymath}
    with $M',M,M''$ being discrete projective $R$-modules. Consider the Koszul resolution 
    $$\Lambda M'\otimes\mathrm{Sym}(M')\rightarrow\mathbb{I}$$
    Since this is a complex of projective modules we may tensor with $\mathrm{Sym}(M'')$ to get a resolution
    $$\Lambda M'\otimes\mathrm{Sym}(M')\otimes\mathrm{Sym}(M'')\rightarrow\mathrm{Sym}(M'')$$
    But $\mathrm{Sym}(M')\otimes\mathrm{Sym}(M'')\cong\mathrm{Sym}(M)$ so we in fact get a resolution
        $$\Lambda M'\otimes\mathrm{Sym}(M)\rightarrow\mathrm{Sym}(M'')$$
        This is a graded resolution, and looking at the $n$th graded piece we get a resolution
        $$0\rightarrow\Lambda^{n}M'\otimes_{R}\mathrm{Sym}^{0}(M)\ldots\Lambda^{0}M'\otimes\mathrm{Sym}^{n}(M)\rightarrow\mathrm{Sym}^{n}(M'')\rightarrow0$$
        Working now as in the proof of \cite{lurie2018spectral} Proposition 25.2.4.2 gives the result.
\end{proof}

\begin{defn}
    A Koszul derived algebraic context is said to be \textit{Illusie}
 if in addition for any connective $A$-module $M$ there are canonical equivalences
 $$\Lambda_{A}^{n}(\Sigma M)\cong\Sigma^{n}\Gamma_{A}^{n}(M)$$

 \end{defn}

The property of being Illusie usually boils down to the existence of for each $n>0$ an exact sequence
$$\Gamma^{n-\bullet}(M')\otimes\Lambda^{\bullet}$$

\begin{cor}
    Let  $(\mathbf{C},\mathbf{C}_{\ge0},\mathbf{C}_{\le0},\mathbf{C}^{0})$ be a Koszul derived algebraic context. Let $M\in\mathbf{C}_{\ge m}$ for $m>0$. Then for $k\ge 1$ we have $\mathbf{LSym}^{k}(M)\in\mathbf{C}_{\ge m+k-1}$. If the derived algebraic context is Illusie and $m>1$ then we have $\mathbf{LSym}^{k}(M)\in \mathbf{C}_{\ge m+2k-2}$.
\end{cor}

\begin{proof}
    We have $\Sigma^{-1}(M)$ is $n-1$-connective. Thus we have
    $$\mathbf{LSym}^{m}_{A}(M)\cong\Sigma^{m}\Lambda^{m}_{A}(\Sigma^{-1}M)$$
    If the context is Illusie and $m\ge 2$, then we further have 
    $$\mathbf{LSym}^{n}_{A}(M)\cong\Sigma{2m}\Lambda^{m}\Gamma^{n}(\Sigma^{-2}M)$$
\end{proof}

In particular we have the following.

\begin{lem}
Let $(\mathbf{C},\mathbf{C}_{\ge0},\mathbf{C}_{\le0},\mathbf{C}^{0})$ be a derived algebraic context. Then 
$$(\mathbf{C},\mathbf{C}_{\ge0}\mathbf{C}_{\le0},\mathbf{C},\mathbf{LSym},\theta)$$
is spectral algebraic context. It is a Koszul spectral algebraic context if $\mathbf{C}^{\heart}$ is Koszul.
\end{lem}

\subsubsection{Derived Coherent Algebras in Derived Algebraic Contexts}

Let $(\mathbf{C},\mathbf{C}_{\ge0},\mathbf{C}_{\le0})$ be a derived algebraic context. We get useful presentations of coherent algebras.

\begin{prop}\label{prop:altLcoh}
Let $A\in\mathbf{DAlg}^{\heart}(\mathbf{C})$ be coherent as an algebra in $\mathbf{C}^{\heart}$. Let $B\in\mathbf{DAlg}^{cn}(\underline{\mathbf{C}})$. Suppose that there is an epimorphism $A\rightarrow \pi_{0}(B)$ with finitely generated ideal, and that each $\pi_{n}(B)$ is finitely presented as an $A$-module. Then there is a presentation
$$B\cong\colim_{m\in\mathbb{N}} X_{m}$$ 
where
\begin{enumerate}
\item 
$X_{0}\cong A$
\item 
for each $n+1$
there is a pushout diagram
$$X_{m+1}\cong X_{m}\otimes^{\mathbb{L}}_{\Sigma^{n_{m}}\mathrm{Sym}_{A}(A)}A$$
or 
$$X_{m+1}\cong X_{m}\otimes^{\mathbb{L}}\Sigma^{n_{m}}\mathrm{Sym}_{A}(A)$$
\end{enumerate}
where for each $n\in\mathbb{N}$ there are only finitely many $m$ such that $n_{m}=n$. Moreover each $\pi_{n}(\mathbb{L}_{B\big\slash A})$ is also finitely presented as an $A$-module.
\end{prop}

\begin{proof}
    The presentation works as in Lemma \ref{lem:Tconstructsequence}. The finiteness condition is a consequence of the fact that to get each homotopy group $\pi_{n}(B)$ we only need to add finitely many generators and kill finitely many relations. Here we are using that $A$ is coherent. 

    For coherence of the cotangent complex we verify the conditions of Lemma \ref{lem:cotangentsubcat} for the class $\mathbf{H}$ of coherent algebras and the class $\mathbf{F}_{A}$ of objects $M$ such that $\pi_{m}(M)$ is finitely presented as a $\pi_{0}(A)$-module. The first condition is clear. Since $\pi_{0}(A)$ is coherent the category of finitely presented $\pi_{0}(A)$-modules is closed under finite limits and colimits and under extensions. Let $M$ be an $A$-module with finitely presented homotopy groups, and let $A\rightarrow B$ be a map with $B$ a derived coherent algebra. Then $B\otimes_{A}^{\mathbb{L}}M$ is such that each $\pi_{m}(B\otimes_{A}^{\mathbb{L}}M)$ is finitely generated as a $\pi_{0}(B)$-module by Corollary \ref{cor:sthngleftcoh}.
    
\end{proof}

\subsubsection{Graded Contexts and The Dirac Category}
Let 
$$\underline{\mathbf{C}}\defeq(\mathbf{C},\mathbf{C}_{\ge0},\mathbf{C}_{\le0},\mathbf{C}^{0})$$
be a spectral algebraic pre-context. Regard $\mathbb{Z}$ as a discrete category, and define
$$\mathbf{Gr}(\mathbf{C})\defeq\mathbf{Fun}(N(\mathbb{Z}),\mathbf{C})$$
to be the category of graded objects in $\mathbf{C}$. This is also a stable $(\infty,1)$-category. As in \cite{HessPstra} Section 2.1 It may be equipped with the \textit{Beilinson }$t$-\textit{structure}, in which $(X_{k})\in\mathbf{Gr}(\mathbf{C})_{\ge0}$ precisely if each $X_{k}\in\mathbf{C}_{\ge -k}$. The category $\mathbf{Gr}(\mathbf{C})_{\ge0}$ is further projectively generated by $\mathbf{Gr}(\mathbf{C})^{0}\defeq\{P(k):P\in\mathbf{C}^{0},k<0\}$. Here $P(k)$ denotes the graded object which is $P$ in degree $k$ and $0$ in all other degrees. It is straightforward to show that 
$$\mathbf{Gr}(\underline{\mathbf{C}})\defeq(\mathbf{Gr}(\mathbf{C}),\mathbf{Gr}(\mathbf{C})_{\ge0},\mathbf{Gr}(\mathbf{C})_{\le0},\mathbf{Gr}(\mathbf{C})^{0})$$
is itself a spectral algebraic pre-context. We wrote

$$\mathrm{Dir}(\underline{\mathbf{C}})\defeq\mathbf{Gr}(\underline{\mathbf{C}})^{\heart}$$
for the heart of $\mathbf{Gr}(\mathbf{C})$ equipped with the symmetric monoidal structure inherited from the one on $\mathbf{Gr}(\mathbf{C})$. In general this will \textit{not be symmetric spectral algebraic pre-context} as explained in \cite{HessPstra} on page 7 for the case of spectra. However it will of course be an $E_{\infty}$-derived algebraic context.

\subsection{Embedding Algebraic Modules}\label{subsec:embeddingalgebra}

Let $\underline{\mathbf{C}}\defeq(\mathbf{C},\mathbf{C}_{\ge0},\mathbf{C}_{\le0},\mathbf{C}^{0})$ be a derived algebraic context. Let $R=\mathbf{Map}(\mathbb{I},\mathbb{I})\cong\mathrm{Hom}(\mathbb{I},\mathbb{I})$. This is a unital commutative ring. Consider the derived algebraic context
$$\mathbf{C}_{R}\defeq(\mathbf{Ch}({}_{R}\mathrm{Mod}),\mathbf{Ch}_{\ge0}({}_{R}\mathrm{Mod}),\mathbf{Ch}_{\le0}({}_{R}\mathrm{Mod)},{}_{R}\mathrm{Mod}^{ffg})$$
Where ${}_{R}\mathrm{Mod}^{ffg}$ is the category of free and finitely generated $R$-modules. By \cite{kelly2021analytic} 3.12 there is a morphism of derived algebraic contexts
$$\mathbb{L}i_{R}:\mathbf{C}_{R}\rightarrow\underline{\mathbf{C}}$$
whose restriction to 
$$\mathbf{Ch}({}_{R}\mathrm{Mod})\cong\mathbf{P}_{\Sigma}(\{R^{\oplus n}\}_{n\in\mathbb{N}_{0}})$$
just sends $R^{\oplus n}$ to $\mathbb{I}^{\oplus n}$ and extends by sifted colimits. Moreover this is fully faithful. In plain terms, any complex $X_{\bullet}\in\mathrm{Ch}({}_{R}\mathrm{Mod})$ may be represented as a cofibrant complex with $X_{n}\cong\bigoplus_{i_{n}\in\mathcal{I}_{n}}R^{i_{n}}$. $\mathbb{L}i_{R}$ sends this complex to $\bigoplus_{i_{n}\in\mathcal{I}_{n}}\mathbb{I}^{i_{n}}$. This is in fact a left Quillen functor.

Now let $B\in\mathbf{DAlg}^{cn}(\mathbf{C})$. Consider the category ${}_{B}\mathbf{Mod}(\mathbf{C})$. Write $|B|_{alg}\defeq\mathbf{Map}_{\mathbf{DAlg}^{cn}(\mathbf{C})}(\mathrm{Sym}(\mathbb{I}),B)\in\mathbf{DAlg}^{cn}(\mathbf{C}_{R})$. Now $\mathbb{L}i_{R}$ is strongly monoidal. Thus we get an induced functor
$$\mathbb{L}\tilde{i}_{B}:{}_{|B|_{alg}}\mathbf{Mod}(\mathbf{C}_{R})\rightarrow{}_{\mathbb{L}i_{R}(|B|_{alg})}\mathbf{Mod}(\mathbf{C})$$
There is a natural map of algebras $\mathbb{L}i_{R}(|B|_{alg})\rightarrow B$, so by base change we get a functor
$$\mathbb{L}i_{B}:{}_{|B|_{alg}}\mathbf{Mod}(\mathbf{C}_{R})\rightarrow{}_{B}\mathbf{Mod}(\mathbf{C})$$
This functor is also strongly monoidal and fully faithful. We write
$${}_{B}\mathbf{Mod}^{alg}(\mathbf{C})\defeq{}_{|B|_{alg}}\mathbf{Mod}(\mathbf{C}_{R})$$

We would also like to understand when this functor is \textit{exact}. Note that by cosntruction it is always right exact.


\begin{thm}\label{thm:algembedexact}
    Suppose that $B$ is left coherent. Then $\mathbb{L}i_{B}$ is exact.
\end{thm}

\begin{proof}
It is enough to prove exactness for finitely presented $\pi_{0}(|B|_{alg})$-modules.
Let $M$ be a finitely presented $\pi_{0}(|B|_{alg})$-module. There is an exact sequence
$$|\pi_{0}(B)|_{alg}^{m}\rightarrow|\pi_{0}(B)|_{alg}^{n}\rightarrow M\rightarrow 0$$
Since $i_{\pi_{0}(B)}$ is right-exact we get an exact sequence
$$\pi_{0}(B)^{m}\rightarrow\pi_{0}(B)^{n}\rightarrow i_{\pi_{0}(B)}(M)\rightarrow 0$$
Thus $i_{\pi_{0}(B)}(M)$ is finitely presented as a $\pi_{0}(B)$-module. We may therefore construct a resolution by the totalisation of a double complex $i_{\pi_{0}(B)}(M)\cong|B^{i_{n}}|$ with $i_{n}$ finite. Since $\mathbf{Map}(\mathbb{I},-)$ commutes with filtered colimits and $\mathbf{Map}(\mathbb{I},i_{\pi_{0}(B)}(M))\cong M$, we have that 
$$M\cong||B|_{alg}^{i_{n}}|$$
Then $\mathbb{L}i_{B}(M)\cong|B^{i_{n}}|\cong i_{\pi_{0}(B)}(M)$. In particular $\mathbb{L}i_{B}(M)$ is discrete for $M$ discrete, and this suffices to prove that $\mathbb{L}i_{B}$ is $t$-exact.
\end{proof}

Later we will see that there are other important situations under which this functor is exact, in particular in functional analytic settings when finitely presented modules are Fr\'{e}chet spaces, the open mapping theorem implies that this functor is exact.


\subsection{Derived Quotients}

When we discuss localisation later, we will need some results about derived quotients. Let $(\mathbf{C},\mathbf{C}_{\ge0},\mathbf{C}_{\le0},\mathbf{C}^{0})$ be a flat spectral algebraic pre-context and $\mathbf{D}:\mathbf{C}\rightarrow\mathbf{C}$ a monad together with a map $\theta:\mathbf{Comm}(-)\rightarrow\mathbf{D}$ such that $(\mathbf{C},\mathbf{D},\theta)$ is a stable $(\infty,1)$-algebraic context, and $(\mathbf{C}_{\ge0},\mathbf{C}_{\le0})$ is a compatible $t$-structure. Let $A\in\mathbf{Alg_{D}}^{cn}(\mathbf{C})$, and write $\mathbb{I}=\mathbf{D}(0)$.

For a collection of maps
$$a_{1},\ldots,a_{n}:\mathbb{I}\rightarrow A$$
denote by 
$$A\big\slash\big\slash(a_{1},\ldots,a_{n})\defeq A\coprod_{\mathbf{D}(\mathbb{I}^{\oplus n})}\mathbb{I}$$
where the map $\mathbf{D}(\mathbb{I}^{\oplus n})\rightarrow \mathbb{I}$ is the canonical augmentation, and the map $R^{\oplus n}$ is induced from the maps 
$$a_{i}:\mathbb{I}\rightarrow A$$
for $1\le i\le n$. Note that the underlying object of $A\big\slash\big\slash(a_{1},\ldots,a_{n})$ is $ A\otimes^{\mathbb{L}}_{\mathrm{Sym}(\mathbb{I}^{\oplus n})}R$.

%

\section{Classes of Maps}

Spectral algebraic contexts are sufficiently rich to define smooth maps between algebras. Here we provide the definition and illustrate some of the basic properties of such maps.

\subsection{Formally Smooth and \'{E}tale Maps}

\begin{defn}[\cite{toen2008homotopical}]
Let $P\in\mathbf{C}_{\ge0}$ be a projective. A map $f:A\rightarrow B$ in $\mathbf{Alg_{D}}(\mathbf{C})$ is said to be \textit{formally }$P$-\textit{smooth} if
\begin{enumerate}
\item
$\mathbb{L}_{B\big\slash A}$ is $P$-projective.
\item
The morphism $\mathbb{L}_{A}\otimes^{\mathbb{L}}_{A}B\rightarrow\mathbb{L}_{B}$ has a retraction in ${}_{B}\mathbf{Mod}^{cn}$. 
\end{enumerate}
$f$ is said to be \textit{formally \'{e}tale} if it is formally $0$-smooth. 
$f$ is said to be \textit{formally smooth} if it is formally $P$-smooth for some $P$. 
\end{defn}

Later $\mathbf{P}$ will usually be a class of \'{e}tale maps determined by a Lawvere theory.


The following can be proven exactly as in \cite{toen2008homotopical} 1.2.8.3.

\begin{lem}
Let $\mathbf{C}$ be a spectral algebraic pre-context. A map $f:A\rightarrow B$ is formally $i$-smooth if and only if for any morphism $B\rightarrow R$ and any $R$-module $M\in{}_{R}\mathbf{Mod}_{\ge1}$ the natural map 
$$\pi_{0}(\mathbf{Map}(\mathbb{L}_{R\big\slash A}),M) \rightarrow \pi_{0}\mathbf{Map}(\mathbb{L}_{B\big\slash A},M)$$
is zero. 
\end{lem}

\begin{lem}
Let $\mathbf{C}$ be a spectral algebraic pre-context. If a map $f:A\rightarrow B$ in $\mathbf{Alg_{D}}$ is formally smooth then for every square-zero extension $S\rightarrow S'$ the map
$$\pi_{0}(\mathbf{Der}_{A}(B,S))\rightarrow \pi_{0}(\mathbf{Der}_{A}(B,S'))$$
is an epimorphism.
\end{lem}

\begin{proof}
Since $\mathbb{L}_{B\big\slash A}$ is projective we have
$$\pi_{0}\mathbf{Map}(\mathbb{L}_{B\big\slash A},M)\cong\mathrm{Hom}(\pi_{0}(\mathbb{L}_{B\big\slash A}),\pi_{0}(M))\cong 0$$
\end{proof}

\begin{cor}
Suppose the spectral algebraic context is connective. Then any formally smooth map is formally $i$-smooth.
\end{cor}

\begin{lem}
\begin{enumerate}
\item
Isomorphisms are formally $0$-smooth (i.e. \'{e}tale).
\item
Let $f:A\rightarrow B$ be formally $P$-smooth and $g:B\rightarrow C$ be formally $Q$-smooth. Then $g\circ f$ is formally $P\oplus Q$-smooth.
\item
Let $f:A\rightarrow B$ be formally $P$-smooth and $A\rightarrow C$ any map. Then $C\rightarrow C\otimes^{\mathbb{L}}_{A}B$ is formally $P$-smooth.
\end{enumerate}
\end{lem}

\begin{proof}
\begin{enumerate}
\item
This is clear.
\item
We have a cofibre sequence
$$\mathbb{L}_{B\big\slash A}\otimes_{B}C\rightarrow\mathbb{L}_{C\big\slash A}\rightarrow\mathbb{L}_{C\big\slash B}$$
Thus we have
$$\mathbb{L}_{C\big\slash A}\cong(\mathbb{L}_{B\big\slash A}\otimes_{B}C)\oplus(\mathbb{L}_{C\big\slash B})$$
The result now follows 
\item
We have $\mathbb{L}_{C\otimes^{\mathbb{L}}_{A}B\big\slash C}\cong (C\otimes^{\mathbb{L}}_{A}B)\otimes^{\mathbb{L}} _{B}\mathbb{L}_{B\big\slash A}$ is $P$-projective as a $C\otimes^{\mathbb{L}}_{A}B$-module. Moreover as a pushout of a split map, 
$$\mathbb{L}_{C}\otimes^{\mathbb{L}}_{C}(C\otimes_{A}^{\mathbb{L}}B)\rightarrow\mathbb{L}_{C\otimes^{\mathbb{L}}_{A}B}$$
is split.
\end{enumerate}
\end{proof}

\begin{defn}
Let $\mathbf{P}$ be a class of formally \'{e}tale maps A map $f:A\rightarrow B$ is said to be \textit{standard }$(P,\mathbf{P})$-\textit{smooth} if there is a factorisation
\begin{displaymath}
\xymatrix{
A\ar[r]^{g} & A\otimes^{\mathbb{L}}\mathrm{Sym}(P)\ar[r]^{h} & B
}
\end{displaymath}
where $h\in\mathbf{P}$.
\end{defn}

Note that such maps are formally $P$-smooth.

\begin{lem}
\begin{enumerate}
\item
If $\mathbf{P}$ contains isomorphisms then isomorphisms are standard $(0,\mathbf{P})$-smooth.
\item
Suppose $\mathbf{P}$ is stable by composition. If $f:A\rightarrow B$ is standard $(P,\mathbf{P})$-smooth and $g:B\rightarrow C$ is standard $(Q,\mathbf{P})$-smooth-smooth then $g\circ f$ is standard $(P\oplus Q,\mathbf{P})$-smooth.
\item
Suppose $\mathbf{P}$ is stable by pushout.  If $f:A\rightarrow B$ is standard $(P,\mathbf{P})$-smooth and $A\rightarrow C$ is any map, then
$$C\rightarrow C\otimes^{\mathbb{L}}_{A}B$$
is standard $(P,\mathbf{P})$-smooth.
\end{enumerate}
\end{lem}

\begin{proof}
The only claim which is not completely trivial is the second one. We can write $F=A\rightarrow A\otimes^{\mathbb{L}} C(P)\rightarrow B$ and $g=B\rightarrow B\otimes^{\mathbb{L}} \mathbf{D}(Q)\rightarrow D$ where the maps $A\otimes^{\mathbb{L}} \mathbf{D}(P)\rightarrow B$ and $B\otimes^{\mathbb{L}} \mathbf{D}(Q)\rightarrow D$ are formally \'{e}tale. Then we can factor $g\circ f$ as 
$$A\rightarrow A\otimes^{\mathbb{L}} \mathbf{D}(P)\otimes^{\mathbb{L}} \mathbf{D}(Q)\rightarrow B$$
where it is straightforward to check that $A\otimes \mathbf{D}(P)\otimes^{\mathbb{L}} \mathbf{D}(Q)\rightarrow B$ is formally \'{e}tale.
\end{proof}

Let $\mathbf{Q}$ be a class of projectives. Say that a map $f:A\rightarrow B$ is \textit{formally }$(\mathbf{Q},\mathbf{P})$-\textit{smooth} if it is formally $(Q,\mathbf{P})$-smooth for some $Q\in\mathbf{Q}$. Then the class of formally $(\mathbf{Q},\mathbf{P})$-smooth maps is stable under equivalence, composition, and pushout.

\subsubsection{Transformations of Spectral Aglebraic Contexts and Smooth Maps}

\begin{lem}
Let $(\mathbf{F},\eta,\delta)$ be a transformation from
$$(\mathbf{C},\mathbf{C}_{\ge0},\mathbf{C}_{\le0},\mathbf{C}^{0},\mathbf{D},\theta,\mathbf{G}_{0},\mathbf{S})$$
 to 
 $$(\mathbf{C}',\mathbf{C}'_{\ge0},\mathbf{C}'_{\le0},\mathbf{C}^{`0},\mathbf{D}',\theta',\mathbf{G}'_{0},\mathbf{S}')$$
Let $f:A\rightarrow B$ be a map in $\mathbf{Alg_{D}}^{cn}(\mathbf{C})$. There is a natural equivalence
$$\mathbb{L}\mathbf{F}(\mathbb{L}_{B\big\slash A})\rightarrow\mathbb{L}_{\mathbf{F}(B)\big\slash\mathbf{F}(A)}$$
\end{lem}

\begin{proof}
The existence of a natural transformation $\mathbf{F}(\mathbb{L}_{B\big\slash A})\rightarrow\mathbb{L}_{\mathbf{F}(B)\big\slash\mathbf{F}(A)}$ follows easily from universal properties. To prove it is an equivalence, without loss of generality we may assume that wer are in the absolute case and $A$ is the monoidal unit. Pick a simplicial resolution $|\mathbf{D}(P_{\bullet})|\cong B$ with $P_{\bullet}$ being projective. Then $\mathbf{F}(B)\cong |\mathbf{D}'(\mathbf{F}(P_{\bullet}))|$. $\mathbf{F}(P_{\bullet})$ is projective. So we get
\begin{align*}
\mathbf{F}(\mathbb{L}_{B})&\cong\mathbf{F}(|B\otimes\mathbf{D}(P_{\bullet})\otimes P_{\bullet}|)\\
&\cong|\mathbf{F}(B)\otimes\mathbf{F}(\mathbf{D}(P_{\bullet}))\otimes\mathbf{F}(P_{\bullet})\\
&\cong |\mathbf{F}(B)\otimes\mathbf{D}'(\mathbf{F}(P_{\bullet}))\otimes\mathbf{F}(P_{\bullet})\\
&\cong\mathbb{L}_{\mathbf{F}(B)}
\end{align*}
\end{proof}

\begin{cor}
Let $(\mathbf{F},\eta,\delta)$ be a transformation from
$$(\mathbf{C},\mathbf{C}_{\ge0},\mathbf{C}_{\le0},\mathbf{C}^{0},\mathbf{D},\theta,\mathbf{G}_{0},\mathbf{S})$$
 to 
 $$(\mathbf{C}',\mathbf{C}'_{\ge0},\mathbf{C}'_{\le0},\mathbf{C}^{`0},\mathbf{D}',\theta',\mathbf{G}'_{0},\mathbf{S}')$$
 If a map $f:A\rightarrow B$ in $\mathbf{Alg_{D}}^{cn}$ is formally $P$-smooth, then $\mathbf{F}(f)$ is formally $\mathbf{F}(P)$-smooth.
\end{cor}

\subsubsection{Homotopy Epimorphisms}

Recall that a map $f:A\rightarrow B$ in $\mathbf{Alg_{D}}(\mathbf{C})$ is said to be a homotopy epimorphism if the map 
$$B\otimes_{A}^{\mathbb{L}}B\rightarrow B$$
is an equivalence.

They are in particular formally \'{e}tale. 

\begin{prop}\label{prop:projhtpy}
Let $A$ and $B$ be objects of  $\mathbf{Alg_{D}}(\mathbf{C})$. The projection $A\times B\rightarrow A$ is a derived strong homotopy epimorphism.
\end{prop}

\begin{proof}
Clearly $\pi_{0}(A\times B)\cong \pi_{0}(A)\times\pi_{0}(B)\rightarrow\pi_{0}(A)$ is flat and an epimorphism. It remains to prove that it is strong. But we have $\pi_{n}(A\times B)\cong \pi_{n}(A)\oplus\pi_{n}(B)$, and 
$$\pi_{0}(A)\otimes_{\pi_{0}(A)\times\pi_{0}(B)}(pi_{n}(A)\oplus\pi_{n}(B))\cong\pi_{n}(A)$$
as required. 
\end{proof}

\section{Discrete Maps}

Let $(\mathbf{C},\mathbf{C}_{\ge0},\mathbf{C}_{\le0},\mathbf{C}^{0})$ be a flat spectral algebraic pre-context. In this section we define various classes of maps in $\mathbf{C}^{\heart}$, and compare them to truncations of the class of maps defined above.

\subsubsection{Algebraic Zariski Localisations}

We begin with a generalisation of Zariski open immersions. Throughout this part we will work in the heart $\mathbf{C}^{\heart}$, and $\mathrm{Sym}$ will denote the underived symmetric algebra functor. We shall denote the unit of the monoidal structure on $\mathbf{C}^{\heart}$ by $\pi_{0}(R)$.

\begin{defn}
Let $A\in\mathrm{Comm}(\mathbf{C}^{\heart})$ and let $f_{0},\ldots,f_{n}:\pi_{0}(R)\rightarrow A$ be maps. The \textit{algebraic Zariski localisation} of $A$ at $f_{0},\ldots,f_{n}$ is 
$$A\otimes^{\mathbb{L}}\mathrm{Sym}(\pi_{0}(R)^{\oplus n})\big\slash\big\slash(1-f_{0},\ldots,1-f_{n})$$ 
\end{defn}

\begin{lem}
Algebraic Zariski localisations are flat epimorphisms.
\end{lem}

\begin{proof}
We clearly may assume that $n=0$. Let $A\in\mathrm{Comm}(\mathbf{C}^{\heart})$ and $g:\pi_{0}(\mathbb{I})\rightarrow A$ a map. Consider the complex
$$(1-gx)\times:A\otimes\mathrm{Sym}(\pi_{0}(R))\rightarrow A\otimes\mathrm{Sym}(\pi_{0}(\mathbb{I}))$$
We claim that this map is a monomorphism. Indeed we may write 
$$A\otimes\mathrm{Sym}(\pi_{0}(\mathbb{I})\cong\bigoplus_{i=0}^{\infty}A$$
The map $xg\times$ shfits and multiplies by $g$. An easy inductive proof then shows that this map is injective. This gives a resolution of $A_{f}$ by flat $A$-modules. Thus for any $A$-module $M$ we have that $M\otimes^{\mathbb{L}}_{A}A_{g}$
is computed by the complex
$$[M\otimes\mathrm{Sym}(\pi_{0}(\mathbb{I}))\rightarrow M\otimes\mathrm{Sym}(\pi_{0}(\mathbb{I}))]$$
and an identical proof shows that the map $M\otimes\mathrm{Sym}(\pi_{0}(\mathbb{I}))\rightarrow M\otimes\mathrm{Sym}(\pi_{0}(\mathbb{I}))$ is injective. This suffices to complete the proof.

%
\end{proof}

We write $A[g^{-1}]$ for the localisation of $A$ at a map $g:\pi_{0}(\mathbb{I})\rightarrow A$

\begin{example}
Let $e\in\mathrm{Hom}_{{}_{A}\mathrm{Mod}}(A,A)$ be such that $e^{2}=e$. Then $A\otimes\mathrm{Sym}(\pi_{0}(\mathbb{I}))\big\slash\big\slash(1-xe)\cong A\big\slash\mathrm{Im}(1-e)$. 
\end{example}

\begin{cor}\label{cor:idempotentloc}
Let $\mathpzc{E}$ be a monoidal elementary abelian category, $A\in\mathrm{Comm}(\mathpzc{E})$, and $I\subset A$ a finitely generated ideal such that the map
$$I\otimes I\rightarrow I$$
is an epimorphism. Then there is a map $e:A\rightarrow A$ such that $A\big\slash I\cong Aotimes^{\mathbb{L}}\mathrm{Sym}(\mathbb{I})\big\slash\big\slash(1-xe)$. 
\end{cor}

\begin{proof}
 We have an epimorphism
$$A^{\oplus n}\rightarrow I$$
for some $n$ such that the composite map
$$A^{\oplus n^{2}}\rightarrow I\otimes I\rightarrow I$$
is also an epimorphism. Write $R=\mathrm{Hom}_{A}(A,A)$ and $I_{R}=\mathrm{Hom}_{A}(A,I)$. Then $I_{R}$ is an ideal of $R$. Moreover the composite map
$$R^{\oplus n^{2}}\rightarrow I_{R}\otimes I_{R}\rightarrow I_{R}$$
is still an epimorphism. Thus $I_{R}\otimes I_{R}\rightarrow I_{R}$ is an epimorphism.

The (standard) proof for modules over a commutative ring (e.g. \cite{42735}) implies that there is $e\in\mathrm{Hom}_{A}(A,I)\subset\mathrm{Hom}_{A}(A,A)$ such that the composition
\begin{displaymath}
\xymatrix{
I\ar[r] & A\ar[r]^{1-e} & I
}
\end{displaymath}
is zero. Moreover $e^{2}=e$ (regarded as a map $A\rightarrow A$) and $e:A\rightarrow I$ is an epimorphism. Thus we have $A\otimes^{\mathbb{L}}\mathrm{Sym}(\mathbb{I})\big\slash I\cong A\big\slash e\cong A\big\slash\big\slash (1-(1-e)x)$
\end{proof}

\subsection{Discrete Smooth and \'{E}tale Maps}

We are now going to define discrete smooth and \'{e}tale maps. We first need to introduce Kaehler differentials. Through this section we fix a monoidal elementary abelian category $\mathbf{C}^{\heart}$ with set of compact projective generators $\mathbf{C}^{0}$, and work in the spectral algebraic pre-context
$$(\mathbf{Ch}(\mathbf{C}^{\heart}),\mathbf{Ch}_{\ge0}(\mathbf{C}^{\heart}),mathbf{Ch}_{\le0}(\mathbf{C}^{\heart}),\mathbf{C}^{0})$$
In particular our algebras will be $E_{\infty}$-algebras.

\subsubsection{The Module of Kaehler Differentials}

\begin{defn}
 Let $f:A\rightarrow B$ be a map in $\mathrm{Comm}(\mathbf{C}^{\heart})$. Define
$$\Omega_{B\big\slash A}\defeq\pi_{0}(\mathbb{L}_{B\big\slash A})$$
to be the \textit{module of relative Kaehler differentials}.
\end{defn}

By Corollary \ref{cor:pi0L} we have the following result.

\begin{cor}
Let $f:A\rightarrow B$ be a map in $\mathbf{Alg_{D}}^{cn}(\mathbf{C})$. Then the natural map
$$\pi_{0}(\mathbb{L}_{B\big\slash A})\rightarrow\Omega_{\pi_{0}(B)\big\slash\pi_{0}(A)}$$
is an equivalence.
\end{cor}

\begin{lem}\label{lem:Kaehlerprop}
Let $f:A\rightarrow B$ be a map in $\mathbf{Alg_{D}}^{\heart}(\mathbf{C})$. The module $\Omega_{B\big\slash A}$ corepresents the functor
$$\mathrm{Der}_{A}(B,-):{}_{A}\mathrm{Mod}\rightarrow\mathrm{Set},\;\; M\mapsto\Hom_{{}_{A\big\backslash}\mathrm{Comm}(\mathbf{C}^{\heart})_{\big\slash B}}(B,B\oplus M)$$
%
\end{lem}

Let $P\in\mathbf{C}^{\heart}$ be projective. Consider the free $E_{\infty}$-algebra $E_{\infty}(P)$, and the free commutative algebra $\mathrm{Sym}(P)$. There is a natural map
$$E_{\infty}(P)\rightarrow\mathrm{Sym}(P)$$

\begin{lem}
The map
$$\pi_{0}(\mathbb{L}_{E_{\infty}(P)})\rightarrow\Omega_{\mathrm{Sym}(P)}$$
is an isomorphism.
\end{lem}

\begin{proof}
We have $\mathbb{L}_{E_{\infty}(P)}\cong E_{\infty}(P)\otimes P$ and $\Omega_{\mathrm{Sym}(P)}\cong\mathrm{Sym}(P)\otimes P$. The claim follows easily from the fact that $\pi_{0}(E_{\infty}(P))\rightarrow\mathrm{Sym}(P)$ is an isomorphism.
\end{proof}


The universal property of $\Omega_{B\big\slash A}$ immediately implies some very useful basic properties. 

\begin{lem}
\begin{enumerate}
\item
Let $A\rightarrow B\rightarrow C$ be a sequence of maps in $\mathrm{Comm}(\mathbf{C}^{\heart})$. Then we have an exact sequence
$$\Omega_{B\big\slash A}\otimes_{B}C\rightarrow\Omega_{C\big\slash A}\rightarrow\Omega_{C\big\slash B}\rightarrow0$$
\item
Let $A\rightarrow B$ and $A\rightarrow C$ be maps in $\mathrm{Comm}(\mathbf{C}^{\heart})$. We have an isomorphism
$$\Omega_{C\otimes_{A}B\big\slash C}\cong\Omega_{B\big\slash A}\otimes_{B}(C\otimes_{A}B)\cong\Omega_{B\big\slash A}\otimes_{A}C$$
\item
Let $A\rightarrow B$ and $A\rightarrow C$ be maps in $\mathrm{Comm}(\mathbf{C}^{\heart})$. The diagram
\begin{displaymath}
\xymatrix{
\Omega_{A}\otimes_{A}(C\otimes_{A}B)\ar[d]\ar[r] &\Omega_{B}\otimes_{A}C\ar[d]\\
\Omega_{C}\otimes_{A}B\ar[r] & \Omega_{C\otimes_{A}B}
}
\end{displaymath}
is a pushout in ${}_{C\otimes_{A}B}\mathrm{Mod}^{\heart}$.
\end{enumerate}
\end{lem}


\begin{cor}
Let $(\mathbf{C},\mathbf{C}_{\ge0},\mathbf{C}_{\le0},\mathbf{D},\theta)$ be a Koszul algebraic context. 
Let $f:A\rightarrow B$ be a map in $\mathbf{Alg_{D}}(\mathbf{C})$ which is an epimorphism of modules, and such that $A\rightarrow B$ has fibre $\overline{I}$. Then 
$$\tau_{\le1}\mathbb{L}_{B\big\slash A}\cong (\pi_{0}(\overline{I}\big\slash \overline{I}^{2})\rightarrow 0)$$
with $\pi_{0}(\overline{I}\big\slash \overline{I}^{2})$ in degree $1$. In particular if $\pi_{1}(B)\cong 0$ then 
$$\tau_{\le1}\mathbb{L}_{B\big\slash A}\cong (I\big\slash I^{2}\rightarrow 0)$$
where $I$ is the kernel of $\pi_{0}(A)\rightarrow\pi_{0}(B)$.
\end{cor}

\begin{proof}
    The map $A\rightarrow B$ has $1$-connective cofibre. Therefore the map
    $$B\otimes_{A}^{\mathbb{L}}\mathrm{cofib}(f)\rightarrow\mathbb{L}_{B\big\slash A}$$
    has a $2$-connective fibre. In particular it induces an equivalence in degrees $0$ and $1$. 
    In degree $0$ $\pi_{0}(B\otimes_{A}^{\mathbb{L}}\mathrm{cofib}(f))\cong\pi_{0}(B)\otimes_{\pi_{0}(A)}0\cong0$. Now we have a cofibre sequence
    $$A\rightarrow B\rightarrow\overline{I}[1]$$
  Tensoring with $B$ we get cofibre sequence
    $$B\rightarrow B\otimes_{A}^{\mathbb{L}}B\rightarrow B\otimes_{A}^{\mathbb{L}}\overline{I}[1]$$
    Equivalently we have a fibre sequence
$$B\otimes_{A}^{\mathbb{L}}\overline{I}\rightarrow B\rightarrow B\otimes_{A}^{\mathbb{L}}B$$
Thus we get a long exact sequence
$$\pi_{1}(B)\rightarrow\pi_{1}(B\otimes_{A}^{\mathbb{L}}B)\rightarrow\pi_{0}(B\otimes_{A}^{\mathbb{L}}\overline{I})\rightarrow\pi_{0}(B)\rightarrow\pi_{0}(B\otimes_{A}^{\mathbb{L}}B)$$
But since $\pi_{0}(B)\rightarrow\pi_{0}(B\otimes_{A}^{\mathbb{L}}B)$ has a section it is injective, so we have
$$\pi_{0}(B\otimes_{A}^{\mathbb{L}}\overline{I})\cong\mathrm{coker}(\pi_{1}(B)\rightarrow\pi_{1}(B\otimes_{A}^{\mathbb{L}}B))\cong\pi_{1}(\mathbb{L}_{B\big\slash A})$$
Tensoring with $\overline{I}\otimes_{A}(-)$ we also get a fibre sequence
$$\overline{I}\otimes_{A}^{\mathbb{L}}\overline{I}\rightarrow\overline{I}\rightarrow\overline{I}\otimes_{A}^{\mathbb{L}}B$$
Thus $\pi_{0}(\overline{I}\otimes_{A}B)\cong\pi_{0}(\overline{I})\big\slash\pi_{0}(\overline{I}^{2})$ as required.
\end{proof}

We can actually compute this more explicilty.
Let $\mathfrak{P}$ denote either a cofibrant replacement of the commutative operad, or just the commutative operad in the case that the latter is admissible. We are going to compute the low degree terms of $\mathbb{L}_{B\big\slash A}$ for $A\rightarrow B$ a map in $\mathrm{Comm}(\mathbf{C}^{\heart})$. Note that by Lemma \ref{lem:Kaehlerprop}, $\Omega_{B\big\slash A}$ is independent of the choice of spectral algebraic pre-context of which it is the heart. Thus we will work in the $E_{\infty}$-algebraic context 
$$(\mathbf{Ch}(\mathbf{C}^{\heart}),\mathbf{Ch}_{\ge0}(\mathbf{C}^{\heart}),\mathbf{Ch}_{\le0}(\mathbf{C}^{\heart}),\mathcal{P})$$
where $\mathcal{P}$ is a set of tiny projective generators of $\mathbf{C}^{\heart}$.

\begin{lem}
Let $\mathfrak{P}\rightarrow\mathpzc{Comm}$ be a weak equivalence of operads with $\mathfrak{P}$ admissible, i.e. the transferred model structure exists on algebras. Let $f:A\rightarrow B$ be a map in $\mathrm{Alg}_{\mathfrak{P}}(\mathrm{s}\mathbf{C}^{\heart})$. Let $Q_{\bullet}\rightarrow B$ be a simplicial resolution in $\mathrm{Alg}_{\mathfrak{P}}({}_{A}\mathrm{Mod}(\mathrm{s}\mathbf{C}^{\heart}))$ where each $Q_{\bullet}$ is free on a cofibrant $A$-module. Then there is an equivalence
$$\mathbb{L}_{B\big\slash A}\cong |\mathbb{L}_{Q_{\bullet}\big\slash A}\otimes_{Q_{\bullet}}B|$$
In particular if $\mathfrak{P}$ is an admissible operad concentrated in degree $0$, then there is an equivalence
$$\mathbb{L}_{B\big\slash A}\cong |\Omega_{Q_{\bullet}\big\slash A}\otimes_{Q_{\bullet}}B|$$
\end{lem}

This allows us to compute the low degree terms of $\mathbb{L}_{B\big\slash A}$ for $A\rightarrow B$ a map in $\mathrm{Comm}(\mathbf{C}^{\heart})$. For general maps $A\rightarrow B$ in $\mathrm{Alg}_{\mathfrak{P}}(\mathrm{s}\mathbf{C}^{\heart})$ we construct a $2$-truncated simplicial resolution of $B$ as an $A$-algebra similarly to \cite[\href{https://stacks.math.columbia.edu/tag/08PV}{Tag 08PV}]{stacks-project} Example 91.5.9. 
Choose a surjection $\phi_{0}:\mathfrak{P}_{A}(P_{0})\rightarrow B$ with $P_{0}$ projective. Set $Q_{0}\defeq\mathfrak{P}_{A}(P_{0})$. Let $I=Ker(\phi_{0})$. Pick a surjection $P_{1}\rightarrow I$, with $P_{1}$ cofibrant, and set $Q_{1}\defeq\mathfrak{P}_{A}(P_{0}\oplus P_{1})$. The face maps $d_{0},d_{1}:Q_{1}\rightarrow Q_{0}$ are the unique $A$-algebra maps defined as follows. The restriction of both $d_{0}$ and $d_{1}$ to $P_{0}$ is the identity. The restriction of $d_{0}$ to $P_{1}$ is $0$, and the restriction of $d_{1}$ to $P_{1}$ is the composition $P_{1}\rightarrow I\rightarrow Q_{0}$. The degeneracy map $s_{0}:Q_{0}\rightarrow Q_{1}$ is induced by the inclusion $P_{0}\rightarrow P_{0}\oplus P_{1}$. There is an exact sequence of double complexes
\begin{displaymath}
\xymatrix{
Q_{1}\ar[r]^{d_{0}-d_{1}} & Q_{0}\ar[r] & B\ar[r] & 0
}
\end{displaymath}
so that the induced map
$$(d_{0},d_{1}):Q_{1}\rightarrow Q_{0}\times_{B}Q_{0}$$
is an epimorphism. We have constructed a $1$-skeletal simplicial $A$-algebra $Q_{\le 1}$ such that the induced map
$$\mathrm{cosk}_{1}(Q_{\le 1})\rightarrow B$$
is an equivalence.

Now we construct an epimorphism $Q_{2}\rightarrow\mathrm{cosk}_{1}(Q_{\le 1})_{2}$. Note that for a simplicial diagram in an abelian category, $\mathrm{cosk}_{1}(Q_{\le 1})_{2}$ can be written as a sum
$$Q_{1}\oplus\mathrm{Ker}(d_{0})\oplus\mathrm{Ker}(d_{1})\cap\mathrm{Ker}(d_{0})$$
Note that in our case $\mathrm{Ker}(d_{0})$ is the ideal of $\mathfrak{P}_{A}(P_{0}\oplus P_{1})$ generated by $P_{1}$. By \cite{kellythesis} Proposition B.3.5 this is given by the image of the map
$$\mathfrak{P}_{A}(P_{0}\oplus P_{1};P_{1})\rightarrow\mathfrak{P}_{A}(P_{0}\oplus P_{1})$$
Thus the image of the map $\mathfrak{P}_{A}((P_{0}\oplus P_{1};P_{1}))\rightarrow Q_{1}$ contains $\mathrm{Ker}(d_{0})$.

Finally let $P_{2}\rightarrow\mathrm{Ker}(d_{0})\cap\mathrm{Ker}(d_{1})$ be an epimorphism with $P_{2}$ cofibrant. Define
$$Q_{2}\defeq\mathfrak{P}_{A}(P_{0}\oplus P_{1}\oplus \mathfrak{P}_{A}(P_{0}\oplus P_{1};P_{1})\oplus P_{2})$$
We have constructed an epimorphism

$$Q_{2}\rightarrow\mathrm{cosk}_{1}(Q_{\le 1})_{2}$$

Now let $f:A\rightarrow B$ be a map in $\mathrm{Alg}_{\mathfrak{P}}^{cn}(\mathbf{C})$ which is an epimorphism of modules, and such that $A\rightarrow B$ has fibre $\overline{I}$. 

Observe that we have an exact sequence
$$0\rightarrow \overline{I}\rightarrow A\rightarrow B\rightarrow 0$$
so after tensoring with $\overline{I}\otimes_{A}(-)$ we get an exact sequence
$$\overline{I}\otimes_{A}\overline{I}\rightarrow \overline{I}\rightarrow \overline{I}\otimes_{A}B\rightarrow 0$$
Thus $\overline{I}\big\slash \overline{I}^{2}\cong \overline{I}\otimes_{A}B$. 

We can take $P_{0}=0$, and $P_{1}\rightarrow I$ an epimorphism. In this case we may actually take
$$Q_{2}\defeq\mathfrak{P}_{A}(P_{1}\oplus P_{1}\oplus P_{2})$$
Now the simplicial object $\mathbb{L}_{Q_{\bullet}\big\slash A}$ is given in low degrees by
$$\mathfrak{P}_{A}( P_{1}\oplus P_{1}\oplus P_{2})\otimes^{\mathbb{L}}_{A}(P_{1}\oplus P_{1}\oplus P_{2})\rightarrow\mathfrak{P}_{A}(P_{1})\otimes^{\mathbb{L}}_{A}P_{1}\rightarrow 0$$
Pulling back to $B$ gives 

$$B\otimes^{\mathbb{L}}_{A}(P_{1}\oplus P_{1}\oplus P_{2})\rightarrow B\otimes^{\mathbb{L}}_{A}P_{1}\rightarrow 0$$
But this is just 
$$B\otimes_{A}(P_{1}\oplus P_{1}\oplus P_{2})\rightarrow B\otimes_{A}P_{1}\rightarrow 0$$
We consider the normalised Moore complex associated to this simplicial complex which (again in the same low degrees) is given by
$$B\otimes_{A}P_{2}\rightarrow B\otimes_{A}P_{1}\rightarrow0$$
The spectral sequence associated to this double complex implies that $\pi_{0}(\mathbb{L}_{B\big\slash A})\cong 0$, and $\pi_{1}(\mathbb{L}_{B\big\slash A})$ is the cokernel of the map
$$\pi_{0}(B\otimes_{A} P_{2})\rightarrow \pi_{0}(B\otimes_{A} P_{1})$$
Consider the map $B\otimes_{A}P_{1}\rightarrow B\otimes_{A}I$. Note that by construction the map $B\otimes_{A}P_{2}\rightarrow\mathrm{Ker}(B\otimes_{A}P_{1}\rightarrow B\otimes_{A}\overline{I})$ is an epimorphism. Thus we have an exact sequence of complexes
$$B\otimes_{A}P_{2}\rightarrow B\otimes_{A}P_{1}\rightarrow B\otimes_{A}\overline{I}\rightarrow 0$$
In particular we see that $\mathrm{coker}(B\otimes_{A}P_{2}\rightarrow B\otimes_{A}P_{1})\cong B\otimes_{A}\overline{I}\cong \overline{I}\big\slash \overline{I}^{2}$ which proves the first claim. The second claim follows from the long exact sequence associated to the exact sequence of complexes.
$$0\rightarrow\overline{I}\rightarrow A\rightarrow B\rightarrow 0$$

Let $f:A\rightarrow B$ be a map in $\mathrm{Alg}_{\mathfrak{P}}^{cn}(\mathbf{C})$. Let $P\rightarrow B$ be an epimorphism with $P$ a projective $A$-module. In particular there is an induced epimorphism of $A$-algebras, $Q\defeq\mathfrak{P}_{A}(P)\rightarrow B$. We get a fibre sequence
$$B\otimes^{\mathbb{L}}_{Q}\mathbb{L}_{Q\big\slash A}\rightarrow\mathbb{L}_{B\big\slash A}\rightarrow\mathbb{L}_{B\big\slash Q}$$
In particular we have an equivalence
$$\mathbb{L}_{B\big\slash A}\cong\mathrm{cone}(\mathbb{L}_{B\big\slash Q}[-1]\rightarrow B\otimes^{\mathbb{L}}_{Q}\mathbb{L}_{Q\big\slash A}\cong B\otimes^{\mathbb{L}}_{A}P)$$
We therefore get a map 
$$\tau_{\le1}(\mathbb{L}_{B\big\slash A})\rightarrow (\pi_{0}(\overline{I})\big\slash\pi_{0}(\overline{I}^{2})\rightarrow \pi_{0}(B)\otimes_{\pi_{0}(A)}\pi_{0}(P))$$
If $\pi_{1}(B)\cong 0$ then in fact we get an equivalence
$$\tau_{\le1}(\mathbb{L}_{B\big\slash A})\rightarrow (I\big\slash I^{2}\rightarrow \pi_{0}(B)\otimes_{\pi_{0}(A)}\pi_{0}(P))$$
where $I$ is the kernel of the map $\pi_{0}(A)\rightarrow\pi_{0}(B)$. 

The important corollaries below are straightforward adaptations of the classical ones, found for example in \cite{stacks-project} Section 10.131.

\begin{cor}\label{cor:naivecot}
Let $f:A\rightarrow B$ be a map in $\mathrm{Comm}(\mathbf{C}^{\heart})$, and $g:Q=\mathrm{Sym}(P)\rightarrow B$ a surjective map with kernel $I$, where $P$ is projective. The map 
$$\tau_{\le1}\mathbb{L}_{B\big\slash A}\rightarrow(I\big\slash I^{2}\rightarrow \Omega_{Q\big\slash A}\otimes_{Q}B)$$
is an equivalence.	
\end{cor}

\begin{cor}
Let $f:A\rightarrow B$ be a map in $\mathrm{Comm}(\mathbf{C}^{\heart})$, and $g:\mathrm{Sym}(P)\rightarrow B$ a surjective map with kernel $I$, where $P$ is projective. Then there is a short exact sequence of $B$-modules
$$I\big\slash I^{2}\rightarrow\Omega_{Q\big\slash A}\otimes_{Q}B\rightarrow\Omega_{B\big\slash A}\rightarrow 0$$
\end{cor}

\begin{proof}
This follows immediately by considering the long exact sequence associated to the fibre sequence
$$B\otimes_{Q}\mathbb{L}_{Q\big\slash A}\rightarrow\mathbb{L}_{B\big\slash A}\rightarrow\mathbb{L}_{B\big\slash Q}$$
\end{proof}

\begin{cor}\label{cor:squarezerosio}
Let $f:A\rightarrow B$ be a map in $\mathrm{Comm}(\mathbf{C}^{\heart})$ and $I\subset B$ an ideal. Write $B'=B\big\slash I^{n+1}$. The map $\Omega_{B\big\slash A}\rightarrow\Omega_{B'\big\slash A}$ induces an isomorphism
$$\Omega_{B\big\slash A}\otimes_{B}B\big\slash I^{n}\rightarrow \Omega_{B'\big\slash A}\otimes_{B'}B\big\slash I^{n}$$
\end{cor}

\begin{proof}
Consider the exact sequence
$$I^{n+1}\big\slash I^{2n+2}\rightarrow\Omega_{B\big\slash A}\otimes_{B}B'\rightarrow\Omega_{B'\big\slash A}$$
Tensoring over $B'$ with $B\big\slash I^{n}$ gives the result.
\end{proof}

\begin{lem}\label{lem:squarezerosplit}
Let $f:A\rightarrow B$ be a map in $\mathrm{Comm}(\mathbf{C}^{\heart})$. Let $f:B\rightarrow B'$ be a map of $A$-algebras with kernel $I$ which has a right inverse $B'\rightarrow B$. Then the sequence
$$I\big\slash I^{2}\rightarrow\Omega_{B\big\slash A}\otimes_{B}B'\rightarrow\Omega_{B'\big\slash A}\rightarrow 0$$
is in fact a short split exact sequence.
\end{lem}

\begin{proof}
It suffices to show that the map $d:I\big\slash I^{2}\rightarrow\Omega_{B\big\slash A}\otimes_{B}B'$ is a split injection. Let $\beta$ be a splitting of $f$, and consider the map
$$\tilde{D}:B\rightarrow I$$
defined by $\tilde{D}\defeq Id-\beta\circ f$, and the map 
$$D:I\big\slash I^{2}$$
induced by composing with the projection $I\rightarrow I\big\slash I^{2}$. $D$ is a derivation, and so induces a map 
$$\tau:\Omega_{B\big\slash A}\otimes_{B}B'\rightarrow I\big\slash I^{2}$$ which is a left inverse to
$$I\big\slash I^{2}\rightarrow \Omega_{B\big\slash A}\otimes_{B}B'$$
\end{proof}

\begin{cor}[c.f. \cite{stacks-project} 10.131.3]
Let $f:R\rightarrow S$ be a map in $\mathrm{Comm}(\mathbf{C}^{\heart})$ and let $I=\mathrm{Ker}(S\otimes_{R}S\rightarrow S)$. Then $\Omega_{R\big\slash S}\cong I\big\slash I^{2}$.
\end{cor}

\begin{proof}
Applying Lemma \ref{lem:squarezerosplit} to the split surjection map $S\otimes_{R}S\rightarrow S$ and using the fact that $\Omega_{S\big\slash S}\cong 0$ to give an isomorphism
$$I\big\slash I^{2}\cong\Omega_{S\otimes_{R}S\big\slash S}\otimes_{S\otimes_{R}S}S$$
On the other hand we have
$$\Omega_{S\otimes_{R}S\big\slash S}\cong\Omega_{S\big\slash R}\otimes_{S}(S\otimes_{R}S)$$
and hence
$$I\big\slash I^{2}\cong\Omega_{S\otimes_{R}S\big\slash S}\otimes_{S\otimes_{R}S}S\cong\Omega_{S\big\slash R}$$
\end{proof}

\subsubsection{Discrete Unramified, Smooth and \'{E}tale Maps}

We can now define discrete unramified, smooth and \'{e}tale maps.

\begin{defn}
A map $f:A\rightarrow B$ in $\mathrm{Comm}(\mathbf{C}^{\heart})$ is said to be
\begin{enumerate}
\item
 \textit{discrete formally unramified} if $\Omega_{B\big\slash A}\cong 0$.
 \item
 \textit{discrete formally  \'{e}tale} if the map $B\otimes_{A}\Omega_{A}\rightarrow\Omega_{B}$ is an isomorphism.
 \item
 \textit{discrete formally }$P$-\textit{smooth} for $P\in\mathbf{C}^{\heart}$ projective if $\Omega_{B\big\slash A}$ is a retract of $B\otimes P$, and the map
 $$\tau_{\le 1}(\mathbb{L}_{B\big\slash A})\rightarrow \Omega_{B\big\slash A}$$
 is an equivalence. 
 \end{enumerate}
\end{defn}

\begin{rem}
Suppose that $A\rightarrow B$ is discrete formally $P$-smooth. By the long exact sequence associated to the fibre-cofibre sequence
$$B\otimes^{\mathbb{L}}_{A}\mathbb{L}_{A}\rightarrow\mathbb{L}_{B}\rightarrow\mathbb{L}_{B\big\slash A}$$
we have that
$$0\rightarrow B\otimes_{A}\Omega_{A}\rightarrow\Omega_{B}\rightarrow\Omega_{B\big\slash A}\rightarrow 0$$
is short exact.
\end{rem}

The following is immediate from the definitions.

\begin{prop}
Formally \'{e}tale maps are discrete formally unramified and discrete formally $0$-smooth.
\end{prop}

\begin{lem}\label{lem:babyNakayama}
Let $f:R\rightarrow S$ be a discrete formally unramified map such that $I\defeq\mathrm{Ker}(S\otimes_{R}S\rightarrow S)$ is finitely generated over $S$. Then the map $S\otimes_{R}S\rightarrow S$ is a Zariski localisation.
\end{lem}

\begin{proof}
We have $\Omega_{S\big\slash R}\cong I\big\slash I^{2}\cong 0$. Thus $I^{2}\cong I$, and we use Corollary \ref{cor:idempotentloc}.
\end{proof}

\begin{lem}
A map $f:A\rightarrow B$ in $\mathrm{Comm}(\mathbf{C}^{\heart})$ is $P$-smooth if and only if the map $\mathbb{L}_{B\big\slash A}\rightarrow \Omega_{B\big\slash A}$ is an equivalence and $\Omega_{B\big\slash A}$ is $P$-projective.
\end{lem}

\begin{proof}
One direction is clear.

Suppose that $f:A\rightarrow B$ is $P$-smooth. Then $\mathbb{L}_{B\big\slash A}$ is a retract of $B\otimes P$. But $B$ is discrete and $P$ is both flat and strong, so $B\otimes P\cong\pi_{0}(B\otimes P)$. Thus $\mathbb{L}_{B\big\slash A}$ is discrete. 
\end{proof}

\begin{lem}[c.f. \cite{stacks-project} Proposition 10.138.8]\label{lem:discsmoothequiv}
Let $f:A\rightarrow B$ be a map in $\mathrm{Comm}(\mathbf{C}^{\heart})$. The following are equivalent.
\begin{enumerate}
\item
$f$ is discrete formally smooth.
\item
$f$ has the left lifting property with respect to maps of the form $T\rightarrow T\big\slash I$, where $I\subset T$ is a square-zero ideal.
\item
Let $Q=\mathrm{Sym}_{A}(P)$ with $P$ projective, and $Q\rightarrow B$ be an epimorphism with kernel $J$. Then there is an $A$-algebra map $\sigma:S\rightarrow Q\big\slash J^{2}$ which is right-inverse to the surjection $P\big\slash J^{2}\rightarrow S$.
\item
 Let $\phi:Q\rightarrow B$ be a map of $A$-algebras with $Q$ being projective and the map being an epimorphism of $B$-modules. Then the sequence 
$$0\rightarrow I\big\slash I^{2}\rightarrow\Omega_{Q\big\slash A}\otimes_{Q}B\rightarrow\Omega_{B\big\slash A}\rightarrow 0$$
is split exact.
\end{enumerate}
\end{lem}

\begin{proof}
We clearly have $1)\Leftrightarrow 4)$ follows from Corollary \ref{cor:naivecot}.

$2)\Rightarrow 3)$ is clear. Let us prove $3)\Rightarrow 2)$. Suppose we have a commutative diagram
\begin{displaymath}
\xymatrix{
A\ar[d]\ar[r] & T\ar[d]\\
B\ar[r] & T\big\slash I
}
\end{displaymath}
with $I$ a square-zero ideal. Since $Q$ is cofibrant, there is a homomorphism of $R$-algebras $\psi:Q\rightarrow T$ lifting the map $Q\rightarrow B\rightarrow T\big\slash I$. Since $\psi|_{J}=0$ and $I^{2}=0$, we have $\psi|_{J^{2}}=0$. Thus $\psi$ factors through a map $Q\big\slash J^{2}\rightarrow A$, and the composite map $B\rightarrow Q\big\slash J^{2}\rightarrow A$ gives the lift.

Now suppose 3) holds. We have
$$\Omega_{Q\big\slash A}\otimes_{Q}B\cong\Omega_{Q\big\slash I^{2}\big\slash A}\otimes_{Q\big\slash I^{2}}B$$
The map $Q\big\slash I^{2}\rightarrow B$ has a section. By Lemma \ref{lem:squarezerosplit} and Corollary \ref{cor:squarezerosio} 4) holds.

Now suppose 4) holds and let $\sigma:\Omega_{B\big\slash A}\rightarrow\Omega_{Q\big\slash A}\otimes_{Q}B$ be a splitting. Consider the map $d-\sigma\circ d\circ\phi:Q\rightarrow\Omega_{Q\big\slash A}\otimes_{Q}B$. Now the composite
\begin{displaymath}
\xymatrix{
Q\ar[rr]^{d-\sigma\circ d\circ\phi}& & \Omega_{Q\big\slash A}\otimes_{Q}B\ar[r] & \Omega_{B\big\slash A}
}
\end{displaymath}
is $0$. Thus the map $d-\sigma\circ d\circ\phi$ factors through $K\defeq\mathrm{Ker}( \Omega_{Q\big\slash A}\otimes_{Q}B\rightarrow\Omega_{B\big\slash A})$. Since $J\rightarrow K$ is an epimorphism and $Q$ is projective as an $R$-module we get a commutative diagram of $R$-modules
\begin{displaymath}
\xymatrix{
 & Q\ar[dl]_{\psi}\ar[d]^{d-\sigma\circ d\circ\phi}\\
 J\ar[r]^{d} & K
}
\end{displaymath}
Define $\overline{s}$ to be the composite of $Id-\psi:Q\rightarrow Q$ with the projection $\pi:Q\rightarrow Q\big\slash J^{2}$. The restriction of $\pi\circ\psi$ to $J$ is $\pi$, so the restriction of $\overline{s}$ to $J$ is zero. Therefore there is an induced map $s:B\rightarrow Q\big\slash J^{2}$. A straightforward computation shows that this is in fact a map of algebras.
\end{proof}

\begin{cor}\label{cor:ismoothdiscmsooth}
Let $\mathbf{C}$ be a connective spectral algebraic context. Let $f:A\rightarrow B$ be a map in $\mathbf{Comm}^{cn}(\mathbf{C})$ which is formally $i$-smooth. Then $\pi_{0}(f):\pi_{0}(A)\rightarrow\pi_{0}(B)$ is discrete formally smooth. 
\end{cor}

\begin{proof}
This follows immediately from Lemma \ref{lem:discsmoothequiv} (2), Example \ref{ex:0small}, Theorem \ref{thm:n-smn-der}, and the definition of formal $i$-smoothness.
\end{proof}

\begin{lem}
Let $S\in\mathrm{Comm}(\mathbf{C}^{\heart})$ and $A\in\mathrm{Comm}({}_{S}\mathrm{Mod})$ be such that $\mathbb{L}_{A\big\slash S}\cong A^{\oplus n}$. Let $I\subset A$ be an ideal. The map $S\rightarrow A\big\slash I$ is discrete formally \'{e}tale if and only if the map 
$$I\big\slash I^{2}\rightarrow (A\big\slash I)^{\oplus n}$$
is an isomorphism. 
\end{lem}

\begin{proof}
The exact sequence
$$I\big\slash I^{2}\rightarrow\Omega_{A\big\slash S}\otimes_{A}A\big\slash I\rightarrow\Omega_{(A\big\slash I)\big\slash A}\rightarrow 0$$
corresponds precisely to
$$I\big\slash I^{2}\rightarrow (A\big\slash I)^{\oplus n}\rightarrow\Omega_{(A\big\slash I)\big\slash R}$$
It is smooth if and only if $I\big\slash I^{2}\rightarrow (A\big\slash I)^{\oplus n}$ is split, and further it is \'{e}tale if and only if $\Omega_{(A\big\slash I)\big\slash R}\cong 0$, i.e. $I\big\slash I^{2}\rightarrow (A\big\slash I)^{\oplus n}$ is an isomorphism. 
\end{proof}

\begin{example}
Let $\mathbf{C}=\mathbf{Ch}({}_{R}\mathrm{Mod})$ for some unital commutative ring $R$. Let $S\in{}_{R}\mathrm{Mod}$, and let $f_{1},\ldots,f_{n}\in S[x_{1},\ldots,x_{n}]$. Then $S\rightarrow S[x_{1},\ldots,x_{n}]\big\slash(f_{1},\ldots,f_{n})$ is formally \'{e}tale if and only if the determinant of the Jacobian $\Bigr(\frac{\partial f_{i}}{\partial x_{j}}\Bigr)_{i,j}$ is a unit in $S[x_{1},\ldots,x_{n}]\big\slash(f_{1},\ldots,f_{n})$. 
\end{example}

\subsubsection{Comparison Between Derived and Discrete Smooth and \'{E}tale Maps}

Let $\mathbf{C}$ be a connective spectral algebraic context.

\begin{prop}[c.f. \cite{toen2008homotopical} 2.2.2.6]\label{prop:discretestrongsmooth}
\begin{enumerate}
\item
Let $f:A\rightarrow B$ be a derived strong map such that $\pi_{0}(A)\rightarrow\pi_{0}(B)$ is formally $P$-smooth. Then $f$ is formally $P$-smooth. 
\item
If $f$ is $P$-smooth then $\pi_{0}(A)\rightarrow\pi_{0}(B)$ is discrete $P$-smooth. 
\item
If the context is Koszul, $f$ and $\pi_{0}(f)$ are formally smooth, and $\pi_{n}(A)$ is transverse to $\pi_{0}(B)$ over $\pi_{0}(A)$ for each $n$, then $f$ is derived strong. 
\end{enumerate}
\end{prop}

\begin{proof}
\begin{enumerate}
\item
Since $f$ is strong we have \[\mathbb{L}_{B\big\slash A}\otimes^{\mathbb{L}}_{B}\pi_{0}(B)\cong\mathbb{L}_{\pi_{0}(B)\big\slash\pi_{0}(A)}.\] Since $\pi_{0}(A)\rightarrow\pi_{0}(B)$ is formally $\pi_{0}(P)$-smooth $\mathbb{L}_{\pi_{0}(B)\big\slash\pi_{0}(A)}$ is a retract of $\pi_{0}(B)\otimes P$. In particular it is $P$-projective. Thus it lifts to a morphism
$$Q\rightarrow\mathbb{L}_{B\big\slash A}$$
with $Q$ a $P$-projective $B$-module, such that \[P\otimes^{\mathbb{L}}_{B}\pi_{0}(B)\cong\mathbb{L}_{\pi_{0}(B)\big\slash\pi_{0}(A)}.\] Let $K$ be the cofiber of $Q\rightarrow\mathbb{L}_{B\big\slash A}$. Then \[K\otimes_{B}^{\mathbb{L}}\pi_{0}(B)\cong 0.\] Hence $\pi_{*}(K)\cong 0$ and $Q\rightarrow\mathbb{L}_{B\big\slash A}$ is an isomorphism. In particular $\mathbb{L}_{B\big\slash A}$ is $P$-projective. 

Now consider the cofibre sequence in $\mathbf{Mod}_{B}$
$$\mathbb{L}_{A}\otimes_{A}^{\mathbb{L}}B\rightarrow\mathbb{L}_{B}\rightarrow\mathbb{L}_{B\big\slash A}$$
and the shifted cofibre sequence
$$\mathbb{L}_{B}\rightarrow\mathbb{L}_{B\big\slash A}\rightarrow\mathbb{L}_{A}\otimes_{A}^{\mathbb{L}}B[-1]$$
Since $\mathbb{L}_{B\big\slash A}$ is a retract of $B\otimes P$, $\mathbf{Map}_{\mathbf{Mod}_{B}}(\mathbb{L}_{B\big\slash A},\mathbb{L}_{A}\otimes_{A}^{\mathbb{L}}B[-1])$ is a retract of \[\mathbf{Map}_{\mathbf{C}}(P,\mathbb{L}_{A}\otimes_{A}^{\mathbb{L}}B[-1])\cong\mathrm{Hom}_{\mathbf{C}^{\heart}}(\pi_{0}(P),\pi_{0}(\mathbb{L}_{A}\otimes_{A}^{\mathbb{L}}B[-1]))\cong 0.\] Thus the map $\mathbb{L}_{B\big\slash A}\rightarrow\mathbb{L}_{A}\otimes_{A}^{\mathbb{L}}B[-1]$ is trivial and the sequence
$$\mathbb{L}_{A}\otimes_{A}^{\mathbb{L}}B\rightarrow\mathbb{L}_{B}\rightarrow\mathbb{L}_{B\big\slash A}$$
splits.
\item
This follows immediately from Corollary \ref{cor:ismoothdiscmsooth}. 
\item
For the converse, let $f:A\rightarrow B$ be a formally $P$-smooth morphism. Since $\mathbb{L}_{B\big\slash A}$ is $P$-projective it is derived strong as a $B$-module. Thus 
$$\pi_{*}\mathbb{L}_{B\big\slash A}\cong\pi_{0}(\mathbb{L}_{B\big\slash A})\otimes^{\mathbb{L}}_{\pi_{0}(B)}\pi_{*}(B)$$
and 
$$\mathbb{L}_{B\big\slash A}\otimes_{B}\pi_{0}(B)\cong\pi_{0}(\mathbb{L}_{B\big\slash A})$$
By Lemma \ref{cor:ismoothdiscmsooth} we know that $\pi_{0}(A)\rightarrow\pi_{0}(B)$ is formally $P$-smooth. Now consider the pushout
\begin{displaymath}
\xymatrix{
A\ar[d]\ar[r] & B\ar[d]\\
\pi_{0}(A)\ar[r] & C
}
\end{displaymath}
By base change $\pi_{0}(A)\rightarrow C$ is formally $P$-smooth. We claim that $C\rightarrow\pi_{0}(B)$ is an isomorphism. Consider the sequence of maps
$$\pi_{0}(A)\rightarrow C\rightarrow\pi_{0}(B)$$
We get a fibre-cofibre sequence
$$\mathbb{L}_{C\big\slash\pi_{0}(A)}\otimes_{C}^{\mathbb{L}}\pi_{0}(B)\rightarrow\mathbb{L}_{\pi_{0}(B)\big\slash\pi_{0}(A)}\rightarrow\mathbb{L}_{\pi_{0}(B)\big\slash C}$$
Now since $\mathbb{L}_{C\big\slash\pi_{0}(A)}$ is projective it is strong, and we have
$$\mathbb{L}_{C\big\slash\pi_{0}(A)}\otimes_{C}^{\mathbb{L}}\pi_{0}(B)\cong\pi_{0}(\mathbb{L}_{C\big\slash\pi_{0}(A)})\cong\Omega_{\pi_{0}(B)\big\slash\pi_{0}(A)}\cong\mathbb{L}_{\pi_{0}(B)\big\slash\pi_{0}(A)}$$
Thus $\mathbb{L}_{\pi_{0}(B)\big\slash C}\cong 0$. The result now follows from Corollary \ref{cor:pi0defequiv} and Proposition \ref{prop:reversederstrong}. 
\end{enumerate}
\end{proof}

%
%
%

\subsubsection{Quasi-Finitely Presented \'{E}tale Maps}

In this short section we generalise the classical result that if $f:A\rightarrow B$ is a map of unital commutative rings such that $B$ is finitely presented as an $A$-algebra, then $f$ is \'{e}tale if and only if it is \'{e}tale in the derived sense.

\begin{defn}
A map $f:A\rightarrow B$ in $\mathrm{Comm}(\mathbf{C}^{\heart})$ is said to be \textit{quasi-finitely presented} if the kernel of the map $B\otimes_{A}B\rightarrow B$ is finitely generated. 
\end{defn}

\begin{lem}\label{lem:flattransformet}
Let $A\rightarrow B$ be a map in $\mathrm{Comm}(\mathbf{C}^{\heart})$ such that
\begin{enumerate}
\item
The map $B$ is transverse to itself over $A$.
\item
The map $B\otimes^{\mathbb{L}}_{B\otimes_{A}B}B\rightarrow B$ is an equivalence
\end{enumerate}
Then $A\rightarrow B$ is formally \'{e}tale.
\end{lem}

\begin{proof}
In this case we have $B\otimes_{A}^{\mathbb{L}}B\cong B\otimes_{A}B\rightarrow B$ is a homotopy epimorphism. Now the result follows from Lemma \ref{lem:formallyetiffdiag}.
\end{proof}

\begin{lem}\label{lem:qfpI}
Let $f:A\rightarrow B$ be a quasi-finitely presented map in $\mathrm{Comm}(\mathbf{C}^{\heart})$. 
\begin{enumerate}
\item
If $f:A\rightarrow B$ is discrete formally unramified then $B\otimes_{A}B\rightarrow B$ is an algebraic Zariski localisation.
\item
If $B$ is transverse to itself over $A$ then $f$ is formally \'{e}tale if and only if it is discrete formally unramified. 
\end{enumerate}
\end{lem}

\begin{proof}
\begin{enumerate}
\item
Let $I=\mathrm{Ker}(B\otimes_{A}B\rightarrow B)$. Then $0\cong\Omega_{B\big\slash A}\cong I\big\slash I^{2}$. The claim now follows from Corollary \ref{cor:idempotentloc}.
\item
This follows immediately from the first part and Lemma \ref{lem:flattransformet}.
\end{enumerate}
\end{proof}

\begin{lem}
Let $(\mathbf{C},\mathbf{C}_{\ge0},\mathbf{C}_{\le0},\mathbf{D},\theta)$ be a Koszul spectral algebraic context. Let $A\rightarrow B$ be a formally \'{e}tale map such that $\mathrm{Ker}(\pi_{0}(B\otimes^{\mathbb{L}}_{A}B)\rightarrow\pi_{0}(B))$ is finitely generated. Then $B\otimes_{A}^{\mathbb{L}}B\rightarrow B$ is a derived strong Zariski localisation.
\end{lem}

\begin{proof}
    The map $B\otimes_{A}^{\mathbb{L}}B\rightarrow B$ is the base change of a formally \'{e}tale map and is therefore formally \'{e}tale. By Lemma \ref{lem:qfpI}. $\pi_{0}(B\otimes^{\mathbb{L}}_{A}B)\rightarrow\pi_{0}(B)$ is an algebraic Zariski localisation. In particular it is formally \'{e}tale. Hence $A\rightarrow B$ is derived strong, and the result follows. 
\end{proof}

\subsubsection{Epimorphisms and Homotopy Epimorphisms}

We conclude this subsection with a discussion of homotopy epimorphisms and discrete epimorphisms. Recall that a map $A\rightarrow B$ in $\mathrm{Comm}(\mathbf{C}^{\heart})$ is an epimorphism precisely if the map
$$B\otimes_{A}B\rightarrow B$$
is an isomorphism.


\begin{prop}\label{prop:htypepistrong}
\begin{enumerate}
\item
Let $f:A\rightarrow B$ be a homotopy epimorphism in $\mathbf{Comm}^{cn}(\mathbf{C})$. Then $\pi_{0}(A)\rightarrow\pi_{0}(B)$ is an epimorphism.
\item
Let $f:A\rightarrow B$ be a derived strong map such that $\pi_{0}(A)\rightarrow\pi_{0}(B)$ is an epimorphism. Then $f$ is a homotopy epimorphism.
\item 
Let $f:A\rightarrow B$ be a homotopy epimorphism such that $\pi_{0}(f)$ is a homotopy epimorphism, and each $\pi_{n}(A)$ is transverse to $\pi_{0}(B)$ over $\pi_{0}(A)$. Then $f$ is a derived strong homotopy epimorphism. 
\end{enumerate}
\end{prop}

\begin{proof}
\begin{enumerate}
\item
We have $\pi_{0}(B)\cong\pi_{0}(B\otimes_{A}^{\mathbb{L}}B)\cong\pi_{0}(B)\otimes_{\pi_{0}(A)}\pi_{0}(B)$.
\item
First note that we have
\begin{equation*}
\begin{split}
\pi_{*}(B)\otimes^{\mathbb{L}}_{\pi_{*}(A)}\pi_{*}(B) & \cong\pi_{*}(A)\otimes_{\pi_{0}(A)}^{\mathbb{L}}\pi_{0}(B)\otimes^{\mathbb{L}}_{\pi_{0}(A)}\pi_{0}(B) \\ & \cong \pi_{*}(A)\otimes^{\mathbb{L}}_{\pi_{0}(A)}\pi_{0}(B)\cong\pi_{*}(B)
\end{split}
\end{equation*}
Consider the spectral sequence
$$\mathrm{Tor}_{\pi_{*}(A)}^{p}(\pi_{*}(B),\pi_{*}(B))_{q}\Rightarrow\pi_{p+q}(B\otimes^{\mathbb{L}}_{A}B)$$
The left-hand side degenerates on the first page to give 
$$\pi_{*}(B)\cong\pi_{*}(B\otimes^{\mathbb{L}}_{A}B)$$
and hence $B\otimes^{\mathbb{L}}_{A}B\cong B$, as required. 
\item
This follows immediately from the assumptions and Proposition \ref{prop:discretestrongsmooth}.
\end{enumerate}
\end{proof}

\begin{cor}
    Let $A\rightarrow B$ be a map in $\mathbf{Alg_{D}}^{cn}(\mathbf{C})$ such that $\pi_{0}(A)\rightarrow\pi_{0}(B)$ is discrete formally unramified and $\mathrm{Ker}(\pi_{0}(B)\otimes_{\pi_{0}(A)}\pi_{0}(B)\rightarrow\pi_{0}(B))$ is finitely generated over $\pi_{0}(A)$. Then the map $B\otimes_{A}^{\mathbb{L}}B\rightarrow B$ is a derived strong Zariski localisation.
\end{cor}


\chapter{Model Category Presentations}\label{MCP}

\section{HA Contexts}

It will be convenient to have presentations of $(\infty,1)$-contexts by model categories. These are furnished by To\"{e}n and Vezzosi's notion of an HA context.

\subsection{Basic Definitions and Properties}
We begin by recalling the definition. Let $\mathcal{C}$ be a $\mathbb{U}$-combinatorial symmetric monoidal model category. For the remainder of this subsection fix a subcategory $\mathcal{C}_{0}$ satisfying the following conditions. 

\begin{ass}\label{ass:c0}
\begin{enumerate}
\item
$\id \in\mathcal{C}_{0}$.
\item
$\mathcal{C}_{0}\subset\mathcal{C}$ is stable by equivalences and $\mathbb{U}$-small homotopy colimits.
\item
$\mathrm{Ho}(\mathcal{C}_{0})\subset\mathrm{Ho}(\mathcal{C})$ is closed under the tensor product $\otimes^{\mathbb{L}}$. 
\end{enumerate}
\end{ass}

For $n\ge 0$ write $\mathcal{C}_{n}$ for the full subcategory of $\mathcal{C}$ consisting of objects which are the $n$-fold suspension of some object of $\mathcal{C}$. Denote by $\mathrm{Comm}(\mathcal{C})_{n}$ the full subcategory of $\mathrm{Comm}(\mathcal{C})$ consisting of commutative monoids whose underlying object lies in $\mathcal{C}_{n}$. For $A\in\mathrm{Comm}(\mathcal{C})$ and $n\ge0$ denote by ${}_{A}\Mod_{n}$ the full subcategory of ${}_{A}\Mod$ consisting of $A$-modules whose underlying object is in $\mathcal{C}_{n}$. Finally, denote by $A-\Comm(\mathcal{C})_{n}$ the full subcategory of $A-\Comm(\mathcal{C})$ consisting of commutative $A$-algebras whose underlying object is in $\mathcal{C}_{n}$.

Equip the category ${}_{A}\mathrm{Mod}$ with the transferred model structure, and the category $\mathrm{Fun}({}_{A}\mathrm{Mod}_{0},\mathrm{sSet})$ with the left Bousfield localisation of the projective model structure on $\mathrm{Fun}({}_{A}\mathrm{Mod}_{0},\mathrm{sSet})$ along equivalences in ${}_{A}\mathrm{Mod}_{0}$. 

Consider the derived functor
$$\mathbb{R}\mathrm{h}_{0}^{-}:{}_{A}\mathrm{Mod}^{op}\rightarrow\mathrm{Fun}({}_{A}\Mod_{0},\mathrm{sSet})$$
of the restricted Yoneda embedding defined by the composition
$$\mathrm{h}_{A,0}^{-}:{}_{A}\mathrm{Mod}^{op}\rightarrow\mathrm{Fun}({}_{A}\Mod,\mathrm{sSet})\rightarrow\mathrm{Fun}({}_{A}\Mod_{0},\mathrm{sSet})$$
where the first functor is the Yoneda embedding, and the second the obvious restriction functor. 

\begin{defn}
$A\in\mathrm{Comm}(\mathcal{C})$ is said to be \textit{good with respect to} $\mathpzc{C}_{0}$ if the functor 
$$\mathbb{R}\mathrm{h}_{A,0}^{-}:\mathrm{Ho}({}_{A}\mathrm{Mod})\rightarrow\mathrm{Ho}(\mathrm{Fun}({}_{A}\mathrm{Mod}_{0},\mathrm{sSet}))$$ 
is fully faithful. 
\end{defn}

\begin{defn}[\cite{toen2008homotopical} Definition 1.1.0.1.1]\label{defn:HA}
A \textit{Homotopical Algebra context} (HA context)  is a triple $(\mathcal{C},\mathcal{C}_{0},\mathcal{A})$ where
$\mathcal{C}$ is a $\mathbb{U}$-combinatorial symmetric monoidal model category, and $\mathcal{C}_{0}\subset \mathcal{C}$ and $\mathcal{A}\subset \Comm(\mathcal{C})$ are full subcategories which are stable by equivalences. \begin{enumerate}
\item
The model category $\mathcal{C}$ is proper, pointed and for any two objects $X$ and $Y$ in $\mathcal{C}$ the natural morphisms
$$QX\coprod QY\rightarrow X\coprod Y\rightarrow RX\times RY$$
are equivalences.\\
\item
$\Ho(\mathcal{C})$ is an additive category.
\item
With the transferred model structure and monoidal structure $\ootimes_{A}$, the category ${}_{A}\Mod$ is a $\mathbb{U}$-combinatorial, proper, symmetric monoidal model category.
\item
For any cofibrant object $M\in{}_{A}\Mod$ the functor
$$-\ootimes_{A}M:{}_{A}\Mod\rightarrow{}_{A}\Mod$$
preserves equivalences.
\item
With the transferred model structures $A-\Comm(\mathcal{C})$ and $A-\Comm_{nu}(\mathcal{C})$ are $\mathbb{U}$-combinatorial proper model categories.
\item
If $B$ is cofibrant in $A-\Comm(\mathcal{C})$ then the functor
$$B\ootimes_{A}-:{}_{A}\Mod\rightarrow\Mod(B)$$
preserves equivalences.
\item
$\mathcal{C}_{0}$ satisfies Assumptions \ref{ass:c0}
\item
Every $A\in\mathcal{A}$ is $\mathcal{C}_{0}$-good.
\end{enumerate}
\end{defn}

The first two conditions will for us be guaranteed by the following result. 

\begin{prop}
Let $\mathpzc{M}$ be a model category which is additive as a category. Then $\mathrm{Ho}(\mathpzc{M})$ is additive, and for any objects $X$ and $Y$ of we have that the maps
$$QX\oplus QY\rightarrow X\oplus Y\rightarrow RX\oplus RY$$
are all equivalences.
\end{prop}

\begin{proof}
Let $X$ and $Y$ be fibrant-cofibrant objects. Left/ right homotopy is clearly a congruence relation on $\mathrm{Hom}(X,Y)$. Thus $\mathrm{Hom}(X,Y)\big\slash\sim$ where $\sim$ denotes the (left or right) homotopy equivalence relation is naturally an abelian group. It is easily seen that this makes $\mathrm{Ho}(\mathpzc{M})$ a pre-additive category. In any model category products of acyclic fibrations are acyclic fibrations, and coproducts of acyclic cofibrations are acyclic cofibrations. Thus in an additive category, finite biproducts of weak equivalences are weak equivalences. This also means they are biproducts in $\mathrm{Ho}(\mathpzc{M})$. Moreover for $X$ and $Y$ any objects of $\mathpzc{M}$ we clearly have that
$$QX\oplus QY\rightarrow X\oplus Y\rightarrow RX\oplus RY$$
are all equivalences. 
\end{proof}

\subsubsection{The Cotangent Complex in HA Contexts}

By \cite{toen2008homotopical} Section 1.2 a homotopical algebra context has sufficient structure to define the relative cotangent complex of a map $f:A\rightarrow B$ in $\mathrm{Alg}_{\mathpzc{Comm}}(\mathpzc{M})$. Let us briefly recall the discussion here. For a commutative monoid $B$ write $\mathrm{Alg}_{\mathpzc{Comm}}^{aug}({}_{B}\mathrm{Mod})\defeq \mathrm{Alg}_{\mathpzc{Comm}}({}_{B}\mathrm{Mod})\big\slash B$ for the category of augmented commutative $B$-algebras. There is an adjunction
$$\adj{K}{\mathrm{Alg}_{\mathpzc{Comm}^{nu}}({}_{B}\mathrm{Mod})}{\mathrm{Alg}_{\mathpzc{Comm}}^{aug}({}_{B}\mathrm{Mod})}{I}$$
where $K$ is the trivial extension functor and $I$ sends an algebra $C$ to the kernel of the map $C\rightarrow B$. This is both an equivalence of categories and a Quillen equivalence of model categories. There is also a Quillen adjunction
$$\adj{Q}{\mathrm{Alg}_{\mathpzc{Comm}^{nu}}({}_{B}\mathrm{Mod})}{{}_{B}\mathrm{Mod}}{Z}$$
where $Q(C)$ is defined by the pushout
\begin{displaymath}
\xymatrix{
C\otimes_{B}C\ar[d]\ar[r] & C\ar[d]\\
\bullet\ar[r] & Q(C)
}
\end{displaymath}
and $Z$ just equips a module $M$ with the trivial non-unital commutative monoid structure. Now given a map $f:A\rightarrow B$ in $\mathrm{Comm}(\mathpzc{E})$ we define the \textbf{relative cotangent complex} by $$\mathbb{L}_{B\big\slash A}\defeq\mathbb{L}Q\mathbb{R}I(B\otimes_{A}^{\mathbb{L}}B)$$
It is shown in \cite{toen2008homotopical} that $\mathbb{L}_{B\big\slash A}$ corepresents the functor of $(\infty,1)$-categories 
$${}_{B}\textbf{Mod}\rightarrow\textbf{sSet},\;M\mapsto \Map_{\mathrm{Comm}({}_{A}\mathrm{Mod})\big\slash B}(B,B\ltimes M)$$
where $B\ltimes M$ is the square-zero extension of $B$ by $M$. In this context it is define as the algebra
$$A\coprod M$$
We have
$$(A\coprod M)\otimes(A\coprod M)\cong A\otimes A\coprod A\otimes M\coprod M\otimes A\coprod M\otimes M$$
The restriction of the multiplication to $A\otimes A$ is given by multiplication on $A$, to $A\otimes M$ and $M\otimes A$ is given by the module structure on $M$, and to $M\otimes M$ is given by the zero map. Now let $C$ be any $A$-algebra and consider the category $\mathrm{Comm}({}_{A}\mathrm{Mod})\big\slash C$. There is a functor
$$\mathrm{Comm}({}_{A}\mathrm{Mod})\big\slash C\rightarrow{}_{C}\textbf{Mod},\; B\mapsto\mathbb{L}_{B\big\slash A}\otimes^{\mathbb{L}}_{B}C$$
It is left adjoint to the functor sending a $C$-module $M$ to the square zero extension $C\ltimes M$. When $C=A=k$ so that $\mathrm{Comm}({}_{A}\mathrm{Mod})\big\slash k=\mathrm{Comm}^{aug}$ we denote this functor by $\mathbb{L}_{0}$.   We will make use of the following facts which constitute Proposition 1.2.1.6 in \cite{toen2008homotopical}. 

\begin{rem}
The (relative) cotangent complex here does not necessarily live in a stable category. Therefore there is in principle a distinction between formally unramified and formally \'{e}tale (the latter implies the former but the converse is not generally true).
\end{rem}

\subsection{Some Homotopical Algebra}

A large class of HA contexts is furnished by the following result.

\begin{thm}\label{thm:HAcontexts}
Let $\mathpzc{C}$ be a $\mathbb{U}$-combinatorial proper symmetric $h$-monoidal model category satisfying the monoid axiom and strong commutative monoid axiom. Suppose further that the domains of the generating cofibrations in $\mathpzc{C}$ are cofibrant and that cofibrant objects are $K$-flat. Then  conditions 3)-6) of Definition \ref{defn:HA} are satisfied.
\end{thm}

We of course need to define some terms in this statement.

\subsubsection{Monoidal Model Categories and Monoid Axioms}

Suppose that $f:X\rightarrow Y$ and $g:X'\rightarrow Y'$ are maps in $\mathpzc{C}$. We define the \textit{pushout-product} map to be the unique map $X'\otimes Y\coprod_{X\otimes Y}X\otimes Y'\rightarrow X'\otimes Y'$ determined by the maps $X'\otimes Y\rightarrow X'\otimes Y'$ and $X\otimes Y'\rightarrow X'\otimes Y'$. We write this map as $f\Box g$. Note that there is a natural action of $\Sigma_{n}$ on $f^{\Box n}$, where $f$ denotes the iterated pushout-product of $f$ with itself $n$-times. 

\begin{defn}
Let $\mathpzc{C}$ be a model category which is equipped with a symmetric monoidal structure. Let $\mathcal{S}$ be a class of maps in $\mathpzc{C}$.
\begin{enumerate}
\item
$\mathcal{S}$ is said to \textit{satisfy the pushout-product axiom} if it is closed under arbitrary pushout-products.
\item
$\mathcal{S}$ is said to \textit{satisfy the weak pushout-product axiom} if whenever $s_{1}\Box\cdots\Box s_{n}$ is an iterated pushout-product of maps in $\mathcal{S}$, and one of the $s_{i}$ is a weak equivalence, then $s_{1}\Box\cdots\Box s_{n}$ is a weak equivalence. 
\item
$\mathcal{S}$ is said to satisfy the \textit{strong commutative monoid axiom} if $s^{\Box n}\big\slash \Sigma_{n}\in\mathcal{S}$ for any $s\in\mathcal{S}$ and any $n\in\mathbb{N}_{0}$.
\item
$\mathpzc{C}$ is said to be a \textit{weak monoidal model category} if cofibrations satisfy the weak pushout-product axiom.
\item
$\mathpzc{C}$ is said to be a \textit{monoidal model category} if cofibrations satisfy the pushout-product axiom and the weak pushout-product axiom.
\item
$\mathpzc{C}$ is said to satisfy the \textit{strong commutative monoid axiom} if both the class of cofibrations and the class of acyclic cofibrations satisfy the strong commutative monoid axiom.
\end{enumerate}
\end{defn}

\begin{defn}
A monoidal model category $(\mathcal{V},\otimes,k)$ is said to satisfy the \textit{monoid axiom} if every morphism which is obtained as a transfinite composition of pushouts of tensor products of acyclic cofibrations with any object is a weak equivalence.
\end{defn} 

If $\mathpzc{C}$ is a monoidal model category which satisfies the monoid axiom, then by \cite{schwede} Theorem 4.1/ Remark 4.2 for any associative monoid $A$ in $\mathpzc{C}$ the transferred model structure exists on ${}_{A}\mathrm{Mod}$. Moreover (also by \cite{schwede} Theorem 4.1), the transferred model structure also exists on the category $\mathrm{Alg}_{\mathpzc{Ass}}(\mathpzc{C})$ of unital associative monoids internal to $\mathpzc{C}$. By \cite{white2017model}  Theorem 3.2 if further $\mathpzc{C}$ satisfies the strong commutative monoid axiom then the transferred model structure exists on the category $\mathrm{Alg}_{\mathpzc{Comm}}(\mathpzc{C})$ of unital commutative monoids, and by \cite{white2017model}, the category $\mathrm{Alg}_{\mathpzc{Comm}^{nu}}(\mathpzc{C})$ of non-unital commutative associative monoids is also equipped with the transferred model structure.

\subsubsection{Properness}

Most of material on properness here is from \cite{kelly2016homotopy} Appendix A, which in turn is largely based on the sources cited below.

\begin{defn}[\cite{white2014monoidal} Definition 8.1.]\label{defn:hcof}
\begin{enumerate}
\item
A map $f:X\rightarrow Y$ in a model category $\mathpzc{C}$
is said to be a $h$-\textit{cofibration} if whenever 
\begin{displaymath}
\xymatrix{
X\ar[d]^{f}\ar[r]^{g} & A\ar[d]\ar[r]^{w}& B\ar[d]\\
Y\ar[r]^{g'} & A'\ar[r]^{w'} & B'
}
\end{displaymath}
is a commutative diagram in which both squares are pushouts, and $w$ a weak equivalence, then $w'$ is a weak equivalence. 
\item
A model category is said to be \textit{left proper} if cofibrations are $h$-cofibrations.
\end{enumerate}
One defines $h$-\textit{fibrations} of maps and \textit{right properness} of model categories dually.
\end{defn}
\begin{rem}
$h$-cofibrations were introduced by Hopkins under the name of \textit{flat maps}. $h$-fibrations are referred to as \textit{sharp maps} in \cite{Rezk1998FibrationsAH} Section 2.
\end{rem}
\begin{prop}[\cite{batanin2013homotopy} Proposition 1.6]
A map $f:X\rightarrow Y$ in a left proper model category $\mathpzc{C}$ is a $h$-cofibration if and only if any pushout diagram
\begin{displaymath}
\xymatrix{
X\ar[d]^{f}\ar[r] & A\ar[d]\\
Y\ar[r] & P
}
\end{displaymath}
is a homotopy pushout.
\end{prop}
This motivates the following definition.
\begin{defn}\label{defn:lprop}
Let $\mathpzc{C}$ be a model category and $\mathcal{P}\subset\mathpzc{C}$ a subcategory of $\mathpzc{C}$. A map $f:X\rightarrow Y$ in $\mathpzc{C}$ is said to be \textit{left proper relative to }$\mathcal{P}$ if any pushout diagram
\begin{displaymath}
\xymatrix{
X\ar[d]^{f}\ar[r] & A\ar[d]\\
Y\ar[r] & P
}
\end{displaymath}
with $A\in\mathcal{P}$ is a homotopy pushout. If $\mathcal{P}=\mathpzc{C}$, then $f$ is said to be \textit{left proper}. One defines relative right properness dually.
\end{defn}
In particular, in a proper model category $h$-cofibrations are left proper. Let us make some straightforward observations. The next proposition is Definition 1.1 in \cite{batanin2013homotopy}. The equivalence to Definition \ref{defn:hcof} is proved in loc. cit immediately thereafter.
\begin{prop}\label{prop:pushoutstableftprop}
Let $\mathcal{S}$ be a pushout-stable class of maps such that any pushout of a weak equivalence along a map in $\mathcal{S}$ is still a weak equivalence. Then any map in $\mathcal{S}$ is a $h$-cofibration.
\end{prop}
The following is just a consequence of the $2$-out-of-$3$ property of weak equivalences.
\begin{prop}\label{prop:pushoutstabweq}
Let $f$ be a weak equivalence such that any pushout of $f$ is a weak equivalence. Then $f$ is left proper, and a $h$-cofibration.
\end{prop}

Suppose now that $\mathpzc{C}$ is a monoidal model category, with tensor product $\otimes$.

\begin{defn}
A map $f:X\rightarrow Y$ in a monoidal model category is said to be a \textit{monoidal (trivial)} $h$-\textit{cofibration} if for any object $Z$ of $\mathpzc{C}$, $Z\otimes f$ is a (trivial) $h$-cofibration.
\end{defn}


\begin{defn}[\cite{white2014monoidal} Definition 4.15]
A monoidal model category is said to be $h$-\textit{monoidal} if for any (trivial) cofibration $f$ and any object $Z$, the map $f\otimes Z$ is a (trivial) $h$-cofibration. 
\end{defn}

\subsubsection{$K$-Flat Objects}

\begin{defn}
    Let $\mathpzc{C}$ be a monoidal model category. An object $X$ of $\mathpzc{C}$ is said to be $K$-\textit{flat} if for any $Y$ in $\mathpzc{C}$ the natural map
    $$X\otimes^{\mathbb{L}}Y\rightarrow X\otimes Y$$
    is an equivalence.
\end{defn}

\begin{lem}\label{Lem:pushouttransK-flat}
\begin{enumerate}
\item
Let
\begin{displaymath}
\xymatrix{
A\ar[d]^{f}\ar[r] & B\ar[d]\\
C\ar[r] & D
}
\end{displaymath}
be a pushout diagram where $f$ is a monoidal $h$-cofibration, and $A,B$, and $C$ are $K$-flat. Then $D$ is $K$-flat.\\
\item
Let $\mathcal{M}$ be a class of maps in $\mathpzc{C}$ such that for any ordinal $\lambda$ ,any $\lambda$-indexed transfinite composition of maps in $\mathcal{M}$ presents the homotopy colimit. Suppose further that if $f\in\mathcal{M}$ and $Q$ is any object of $\mathpzc{C}$, then $Q\otimes f$ is in $\mathcal{M}$. Let 
\begin{displaymath}
\xymatrix{
X_{0}\ar[r] & X_{1}\ar[r] & \cdots\ar[r] & X_{\beta}\ar[r] & \cdots
}
\end{displaymath}
be a $\lambda$-indexed transfinite sequence where each $X_{\beta}$ is $K$-flat and each map in the diagram is in $\mathcal{M}$. Then $\limind_{\beta<\lambda} X_{\beta}$ is $K$-flat. 
\end{enumerate}
\end{lem}

\begin{proof}
\begin{enumerate}
\item
By \cite{white2017model} the diagram is in fact a homotopy pushout. Let $Q$ be an object of $\mathpzc{C}$. Consider the following commutative diagram
\begin{displaymath}
\xymatrix{
& A\otimes Q\ar[dd]\ar[rr] & &B\otimes Q\ar[dd]\\
A\otimes^{\mathbb{L}}Q\ar[dd]\ar[ur]\ar[rr] & & B\otimes^{\mathbb{L}}Q\ar[dd]\ar[ur]\\
& C\otimes Q\ar[rr] && D\otimes Q\\
C\otimes^{\mathbb{L}}Q\ar[ur]\ar[rr] &  & D\otimes^{\mathbb{L}}Q \ar[ur]
}
\end{displaymath}
Now since $f$ is a monoidal $h$-cofibration the near and far faces of the cube are both homotopy pushouts. Now the maps connecting the top left, bottom left, and top right vertices of the near and far faces are weak equivalences. Therefore the map $D\otimes^{\mathbb{L}}Q\rightarrow D\otimes Q$ is an equivalence. Since $Q$ was arbitrary this implies that $D$ is $K$-flat. 
\item
Consider a diagram
\begin{displaymath}
\xymatrix{
X_{0}\ar[r] & X_{1}\ar[r] & \cdots\ar[r] & X_{\beta}\ar[r] & \cdots
}
\end{displaymath}
where each $X_{\beta}$ is $K$-flat and each map in the diagram is in $\mathcal{M}$. 
We have
\begin{align*}
Q\otimes^{\mathbb{L}}\colim X_{\beta}&\cong Q\otimes^{\mathbb{L}}\hocolim X_{\beta}\\
&\cong\hocolim Q\otimes^{\mathbb{L}}X_{\beta}\\
&\cong\hocolim Q\otimes X_{\beta}\\
&\cong\colim Q\otimes X_{\beta}\\
&\cong Q\otimes \colim X_{\beta}
\end{align*}
\end{enumerate}
\end{proof}

\begin{cor}\label{cor:K-flatcofib}
Let $\mathpzc{C}$ be a cofibrantly generated $h$-monoidal model category in which the domains and codomains of generating cofibrations are $K$-flat. Let $f:X\rightarrow Y$ be a cofibration such that $X$ is $K$-flat. Then $Y$ is $K$-flat.
\end{cor}

\subsubsection{Proof of the Theorem}

We are now ready to proof Theorem  \ref{thm:HAcontexts}.

\begin{proof}[Proof of Theorem \ref{thm:HAcontexts}]
First let us show that for any $A\in\Comm(\mathpzc{C})$, ${}_{A}\mathrm{Mod}$ is also a $\mathbb{U}$-combinatorial proper symmetric $h$-monoidal model category satisfying the monoid axiom and strong commutative monoid axiom, that the domains of the generating cofibrations are cofibrant, and cofibrant objects are $K$-flat.

Let $I$ denote the class of generating cofibrations in $\mathpzc{C}$. Then $A\otimes I\defeq\{A\otimes f:f\in I\}$ is a set of generating cofibrations in $\Comm({}_{A}\mathrm{Mod}(\mathpzc{C}))$. Since free modules on cofibrant objects are cofibrant, the domains of maps in $A\otimes I$ are cofibrant.

We first show that ${}_{A}\mathrm{Mod}$ is $h$-monoidal. Let $\tilde{f}$ be a generating (acyclic) cofibration in ${}_{A}\mathrm{Mod}$. We may assume it is of the form $A\otimes f$ for $f:X\rightarrow Y$ a generating (acyclic) cofibration in $\mathpzc{C}$. Let $Q$ be an object of ${}_{A}\mathrm{Mod}$. Then $Q\otimes_{A} f\cong Q\otimes\tilde{f}$. This is an (acyclic) $h$-cofibration in $\mathpzc{C}$, and hence an (acyclic) $h$-cofibration in ${}_{A}\mathrm{Mod}$.

 Let us prove that cofibrant objects in ${}_{A}\mathrm{Mod}$ are $K$-flat. By Corollary \ref{cor:K-flatcofib} it suffices to observe that the domains and codomains of generating cofibrations are of the form $A\otimes X$ for $X$ a cofibrant, and hence $K$-flat object of $\mathpzc{C}$. Therefore $A\otimes X$ is $K$-flat in ${}_{A}\mathrm{Mod}$.

The fact that $\Comm({}_{A}\mathrm{Mod})$ is a $\mathbb{U}$-combinatorial proper model categories now follows immediately from Theorem 4.17 in \cite{white2014monoidal}. The fact that $\Comm^{nu}({}_{A}\mathrm{Mod})$ is, is \cite{kelly2016homotopy} Lemma 6.3.31.

Since cofibrant objects are $K$-flat, the functor
$$-\otimes_{A}M:{}_{A}\mathrm{Mod}\rightarrow{}_{A}\mathrm{Mod}$$
sends equivalences to equivalences whenever $M$ is cofibrant. Now by Proposition 3.5 in \cite{white2014monoidal} if $B$ is cofibrant in $\Comm({}_{A}\mathrm{Mod})$ then it is cofibrant as an $A$-module. Therefore the functor 
$$B\otimes_{A}(-):{}_{A}\mathrm{Mod}\rightarrow{}_{B}\mathrm{Mod}$$
sends weak equivalences to weak equivalences. 

In an earlier version of this we neglected to include the proof that for $A\in\mathrm{Comm}(\mathpzc{C})$ then ${}_{A}\mathrm{Mod}$ satisfies the (strong) commutative argument. This gap was resolved by Savage in the proof of a related result - we give their proof here. Let $f\in{}_{A}\mathrm{Mod}$ be a cofibration or an acyclic cofibration. Without loss of generality we may assume that it is generating, so of the form $A\otimes\tilde{f}$ for some $\tilde{f}$ a cofibration or acyclic cofibration in $\mathpzc{C}$. But then $f^{\boxtimes n}\big\slash\Sigma_{n}\cong A\otimes(\tilde{f}^{\boxtimes n}\big\slash\Sigma_{n})$. By the (strong) commutative monoid axiom in $\mathpzc{C}$ $\tilde{f}^{\boxtimes n}\big\slash\Sigma_{n}$ is a cofibration (resp. an acyclic cofibration). Then $f$ clearly is as well. 
\end{proof}


\subsection{Homotopical Algebra in Exact Categories}\label{sec:haec}

In this section we can explain how one can construct many examples of HA contexts using (left) exact structures on additive categories. 
\subsubsection{Relative Homological Algebra}

\begin{defn}[\cite{christensen2002quillen} Section 1]
Let $\mathcal{C}$ be an additive category and $\mathcal{P}$ a full subcategory of $\mathcal{C}$.
\begin{enumerate}
\item
A map $f:X\rightarrow Y$ in $\mathcal{C}$ is said to be a $\mathcal{P}$-\textit{epimorphism} if for any $P\in\mathcal{P}$ the map of abelian groups $\mathrm{Hom}(P,f):\mathrm{Hom}(P,X)\rightarrow\mathrm{Hom}(P,Y)$ is an epimorphism. The class of all $\mathcal{P}$-epimorphisms is denoted $\mathcal{E}_{\mathcal{P}}$.
 \item
$(\mathcal{P},\mathcal{E}_{\mathcal{P}})$ is said to be a \textit{projective class} if whenever $\mathrm{Hom}(Q,f)$ is an epimorphism for any $f\in\mathcal{E}_{\mathcal{P}}$ then $Q\in\mathcal{P}$.
\item
A projective class $(\mathcal{P},\mathcal{E}_{\mathcal{P}})$ is said to \textit{have enough projectives} if for any object $C\in\mathcal{C}$ there is a map $P\rightarrow C$ in $\mathcal{E}_{\mathcal{P}}$ with $P$ in $\mathcal{P}$
 \end{enumerate}
\end{defn}

\begin{rem}
A projective class $(\mathcal{P},\mathcal{E}_{\mathcal{P}})$ on $\mathcal{C}$ determines a \textit{left exact category} structure on $\mathcal{C}$ in the sense of \cite{bazzoni2013one}, with $\mathcal{E}_{\mathcal{P}}$ being the class of admissible epimorphisms and $\mathcal{P}$ being the class of projectives. 
\end{rem}

Let us introduce some useful notation and terminology. For $\mathcal{G}$ a class of objects in $\mathpzc{E}$ denote by 

\begin{enumerate}
    \item 
    $\bigoplus\mathcal{G}$ the class of objects which are coproducts of objects in $\mathcal{G}$
    \item 
    $C(\mathcal{G})$ the class of objects which are summands of coproducts of objects of $\mathcal{G}$
\end{enumerate}

\begin{defn}
\begin{enumerate}
\item
A $\mathcal{P}$-epimorphism will also be called an \textit{admissible epimorphism}
\item
Let $(\mathcal{P},\mathcal{E}_{\mathcal{P}})$  be a projective class on $\mathcal{C}$. A map $f:X\rightarrow Y$ in $\mathcal{C}$ is said to be a \textit{weakly admissible }$\mathcal{P}$-\textit{monomorphism} if
\begin{enumerate}
\item
 for any $P\in\mathcal{P}$ the map \[\mathrm{Hom}(P,f):\mathrm{Hom}(P,X)\rightarrow\mathrm{Hom}(P,Y)\] is a monomorphism of abelian groups.
 \item
the functor $\mathrm{Hom}(P,-)$ commutes with pushouts along $f$.
 \end{enumerate}
 \end{enumerate}
\end{defn}

Note that if the projective class determines a full Quillen exact structure then any admissible monomorphism for the exact structure is weakly admissible. 
\begin{defn}
Let $(\mathcal{P},\mathcal{E}_{\mathcal{P}})$  be a projective class on $\mathcal{C}$, $\kappa$ a cardinal, and $\mathcal{S}$ a class of morphisms in $\mathcal{C}$
\begin{enumerate}
\item
A subcategory $\mathcal{G}$ is said to \textit{generate }$\mathcal{C}$ if for each $X\in\mathcal{C}$ there is $G\in\mathcal{G}$ and a $\mathcal{P}$-epimorphism $G\rightarrow X$.
\item
An object $X$ of $\mathcal{C}$ is said to be $(\kappa,\mathcal{S})$-small if for any cardinal $\lambda$ with $\mathrm{cofin}(\lambda)\ge\kappa$ and any $\lambda$-indexed transfinite sequence $F:\lambda\rightarrow\mathcal{C}$ with $F_{i}\rightarrow F_{i+1}$ in $\mathcal{S}$ for $i\le i+1<\lambda$, the map
$$\limind_{\beta\in\lambda}\mathrm{Hom}(X,F_{\beta})\rightarrow\mathrm{Hom}(X,\limind_{\beta\in\lambda}F_{\beta})$$
is an isomorphism.
\item 
An object $E$ of $\mathpzc{E}$ is said to be \textit{tiny} if $\mathrm{Hom}(P,-):\mathpzc{E}\rightarrow\mathrm{Ab}$ commutes with filtered colimits.
\item
$\mathcal{S}$ is said to be $(\kappa,\mathcal{S})^{small}$-elementary if it has a small generating subcategory consisting of projective $(\kappa,\mathcal{S})$-small objects.
\end{enumerate}
\end{defn}

Generally we will be interested in the class $\mathbf{SplitMon}$ of split monomorphisms.

\begin{prop}
A set $\mathcal{Z}$ of objects of $\mathcal{C}$ is a generating subcategory if and only if for all objects $E$ of $\mathcal{C}$ there exists an object $C$ of $C(\mathcal{Z})$ together with an admissible epimorphism $C\rightarrow E$. 
\end{prop}

\begin{proof}
Suppose that $\mathcal{Z}$ be a generating subcategory and let $E$ be an object of $\mathcal{C}$. Then there exists a set $\{Z_{i}:i\in\mathcal{I}_{E}\}$  of objects of $\mathcal{Z}$ together with an admissible epimorphism
$\bigoplus_{i\in\mathcal{I}_{E}}Z_{i}\rightarrow E$. But $\bigoplus_{i\in\mathcal{I}_{E}}Z_{i}\in C(\mathcal{Z})$.

Conversely let $E$ be an object of $\mathcal{E}$ and suppose there is $C$ in $C(\mathcal{Z})$ and an admissible epimorphism $C\rightarrow E$. But $C$ is a retract of an object of the form $\bigoplus_{i\in\mathcal{I}}Z_{i}$ where $\{Z_{i}:i\in\mathcal{I}\}$ is a subset of $\mathcal{Z}$. Thus there is an admissible epimorphism $\bigoplus_{i\in\mathcal{I}}Z_{i}\rightarrow C$, and by composition an admissible epimorphism $\bigoplus_{i\in\mathcal{I}}Z_{i}\rightarrow E$. 
\end{proof}
\begin{lem}\label{lem:gencproj}
Let $\mathcal{Q}$ be a set of objects in $\mathcal{C}$. Then $\mathcal{Q}$ is a generating subcategory consisting of projective objects if and only if $\mathcal{C}$ has enough projectives and $C(\mathcal{Q})=\mathcal{P}$, where $\mathcal{P}$ is the class of projectives. 
\end{lem}

\begin{proof}
Suppose that $\mathcal{C}$ has enough projectives and $C(\mathcal{Q})=\mathcal{P}$. Since $\mathcal{Q}\subset C(\mathcal{Q})$, $\mathcal{Q}$ consists of projective objects. Since $\mathcal{C}$ has enough projectives, for every object $E$ of $\mathcal{C}$ there is a $P\in\mathcal{P}=C(\mathcal{Q})$ together with an admissible epimorphism $P\rightarrow E$. Thus $\mathcal{Q}$ is a generating subcategory.

Conversely suppose that $\mathcal{Q}$ is a generating subcategory consisting of projective objects. Since $\mathcal{P}$ is closed under taking coproducts, $\mathcal{C}$ has enough projectives. Moreover $\mathcal{P}$ is also closed under retracts, so $C(\mathcal{Q})\subset\mathcal{P}$. Now let $P\in\mathcal{P}$. There is a collection $\{Q_{i}:i\in\mathcal{Q}_{P}\}$ and an admissible epimorphism $\bigoplus_{i\in\mathcal{Q}_{P}}Q_{i}\rightarrow P$. Since $P$ is projective this epimorphism is split, so $P$ is a summand of $\bigoplus_{i\in\mathcal{Q}_{P}}Q_{i}$. In particular it is in $C(\mathcal{Q})$. 
\end{proof}

\subsubsection{The Model Structure}

Let us now explain how to equip the category $\mathrm{s}\mathcal{C}$ with the projective model structure, and enumerate its properties.

\begin{lem}\label{lem:sabprop}
Suppose that $\mathrm{sAb}$ is equipped with a model structure such that the weak equivalences are the weak homotopy equivalences of simplicial abelian groups.  Then level-wise monomorphisms are $h$-cofibrations and maps which are epimorphisms in each strictly positive level are $h$-fibrations.
\end{lem}

\begin{proof}
Using the exactness of the Dold-Kan correspondence we may pass to $\mathrm{Ch}_{\ge0}(\mathrm{Ab})$. Here it suffices to prove that for the the pushout of a weak equivalence along a monomorphism is a weak equivalence and the pullback of a weak equivalence along a map which is an epimorphism in each strictly positive degree is a weak equivalence. This is a straightforward computation.
\end{proof}

\begin{thm}
Let $\mathcal{C}$ be a complete and cocomplete additive category, and $(\mathcal{P},\mathcal{E}_{\mathcal{P}})$ a projective class on $\mathcal{C}$ with enough projectives. Consider the category $\mathrm{s}\mathcal{C}$ of simplicial objects in $\mathcal{C}$. There exists a model structure on $\mathrm{s}\mathcal{C}$ in which
\begin{enumerate}
\item
A map $f:X\rightarrow Y$ is a weak equivalence precisely if for each $P\in\mathcal{P}$ the map $\mathrm{Hom}(P,f)$ is a weak equivalence of simplicial abelian groups.
\item
A map $f:X\rightarrow Y$ is a fibration precisely if for each $P\in\mathcal{P}$ the map $\mathrm{Hom}(P,f)$ is a fibration of simplicial abelian groups.
\end{enumerate}
The model category structure is cofibrantly generated if $\mathcal{C}$ is $(\aleph_{0},\mathbf{SplitMon})^{\mathrm{small}}$, where $\aleph_{0}$ is the first uncountable cardinal. Moreover in this model structure degree-wise weakly admissible monomorphisms are $h$-cofibrations and maps which are $\mathcal{P}$-epimorphisms in each strictly positive level are $h$-fibrations. In particular $\mathcal{C}$ is both left and right proper. 
\end{thm}

\begin{proof}
The existence of the model structure and the claim about cofibrant generation is precisely Theorem 6.3 in \cite{christensen2002quillen}. Let $f:X\rightarrow Y$ be a $\mathcal{P}$-epimorphism in each strictly positive level. Then for each $P\in\mathcal{P}$ $\mathrm{Hom}(P,f)$ is an epimorphism in each strictly positive level. Since the functor $\mathrm{Hom}(P,-)$ commutes with limits the claim about $h$-fibrations follows from Lemma \ref{lem:sabprop}. The claim about $h$-cofibrations is similar. Right properness follows since fibrations are in particular $\mathcal{P}$-epimorphsims in each degree. Left properness follows because the proof of Theorem 6.3 shows that cofibrations are, in each degree, split monomorphisms. Thus a pushout along a cofibration $f$ is also a split monomorphism and the cokernel coincides with the cokernel of $f$. Since $\mathcal{C}$ is additive $\mathrm{Hom}(P,-)$ clearly commutes with such pushouts. Moreover $\mathrm{Hom}(P,f)$ is a monomorphism of simplicial abelian groups. Hence cofibrations are in particular degree-wise weakly admissible monomorphisms.
\end{proof}

\begin{rem}\label{rem:cofgens}
Let $\mathcal{G}$ be a generating set of projetives. By \cite{christensen2002quillen} Theorem 6.3, in the case that objects of $\mathcal{G}$ are $(\aleph_{0},\mathbf{SplitMon})$-small, generating fibrations and cofibrations can be described explicitly. The generating cofibrations and generating acyclic cofibrations are given by the following sets respectively
  \[I = 
  \{ P \otimes \partial \Delta[n] 
  \longrightarrow 
  P \otimes \Delta[n] \ \ | \ \  P \in \mathcal{G} , n \geq 0 \} 
  \]
  and
  \[J = \{ P \otimes V[n,k] \longrightarrow P\otimes \Delta[n] \ \  | \ \  P \in \mathcal{G} , n \geq k \geq 0, n>0 \},
  \]
  where $\partial \Delta[n]$ is the boundary of $\Delta[n]$ and $V[n,k]$ is the sub complex of $\Delta[n]$ resulting from removing the $n$-cell and $k$th face.
\end{rem}

Let $\mathcal{C}$ be a complete and cocomplete  additive category. Then $\mathcal{C}$ is tensored over abelian groups. For an object $C$ of $\mathcal{C}$ and a free abelian group of finite rank $\mathbb{Z}^{n}$ we define 
$$\mathbb{Z}^{n}\otimes C\defeq C^{\oplus n}$$
The category of abelian groups is freely generated under sifted colimits by free abelian groups of finite rank, so this definition extends essentially uniquely to a colimit preserving functor
$$\otimes:\mathrm{Ab}\times\mathcal{C}\rightarrow\mathcal{C}$$
There is also a bifunctor
$$(-)^{(-)}:\mathcal{C}\times\mathrm{Ab}^{op}\rightarrow\mathcal{C}$$
defined on $\mathbb{Z}^{n}$ by $(C)^{\mathbb{Z}^{\oplus n}}\defeq C^{\oplus n}$, and extending to all abelian groups by \textit{limits}. Then
$$\otimes:\mathrm{Ab}\times\mathcal{C}\rightarrow\mathcal{C}$$
$$\Hom(-,-):\mathcal{C}^{op}\times\mathcal{C}\rightarrow\mathrm{Ab}$$
$$(-)^{(-)}:\mathcal{C}\times\mathrm{Ab}^{op}\rightarrow\mathcal{C}$$
is a two-variable adjunction. We can extend this to a two-variable adjunction on the level of simplicial objects as follows.
$$\otimes:\mathrm{s}\mathrm{Ab}\times\mathrm{s}\mathcal{C}\rightarrow\mathrm{s}\mathcal{C}$$
is defined by
$$(X_{\bullet}\otimes C_{\bullet})_{n}\defeq X_{n}\otimes C_{n}$$
Define
$$\mathrm{Hom}(\mathbb{Z}[\Delta^{n}],C_{\bullet})\defeq C_{n}$$
and extend to all simplicial abelian groups by colimits. Then for $X_{\bullet}$ a simplicial abelian group define a simplicial object $C_{\bullet}^{X_{\bullet}}$ by
$$(C_{\bullet}^{X_{\bullet}})_{n}\defeq\mathrm{Hom}(X\otimes\Delta_{n},C_{\bullet})$$

For $C_{\bullet},D_{\bullet}\in\mathrm{s}\mathpzc{C}$, define 
$$\underline{\Hom}_{\mathpzc{C}}(D_{\bullet},C_{\bullet})$$ 
to be the simplicial abelian group given by
$$(\underline{\Hom}_{\mathpzc{C}}(D_{\bullet},C_{\bullet})_{n}\defeq\Hom_{\mathpzc{C}}(D_{\bullet}\otimes\mathbb{Z}[\Delta^{n}],C_{\bullet})$$
Then

$$\otimes:\mathrm{s}\mathrm{Ab}\times\mathrm{s}\mathcal{C}\rightarrow\mathrm{s}\mathcal{C}$$
$$\underline{\Hom}(-,-):(\mathrm{s}\mathcal{C})^{op}\times\mathrm{s}\mathcal{C}\rightarrow\mathrm{s}\mathrm{Ab}$$
$$(-)^{(-)}:\mathrm{s}\mathcal{C}\times(\mathrm{s}\mathrm{Ab})^{op}\rightarrow\mathrm{s}\mathcal{C}$$
is a two-variable adjunction.

\begin{obs}
Let $\mathcal{P}$ be a class of objects of $\mathcal{C}$ such that  for any $P\in\mathcal{P}$ $\Hom(P,-):\mathcal{C}\rightarrow\mathrm{Ab}$ commutes with filtered colimits indexed by some filtered category $\mathcal{I}$. Then if $P_{\bullet}$ is an $n$-coskeletal (for some $n$) object of $\mathrm{s}\mathcal{C}$ such that for each $P_{n}\in\mathcal{P}$, the functor
$$\underline{\mathrm{Hom}}(P_{\bullet},-):\mathrm{s}\mathcal{C}\rightarrow\mathrm{sAb}$$
commutes with filtered colimits in $\mathrm{s}\mathcal{C}$ indexed by $\mathcal{I}$. 
\end{obs}

\begin{lem}\label{lem:tensoredoversimp}
Let $\mathcal{C}$ be a complete and cocomplete additive category, and $(\mathcal{P},\mathcal{E}_{\mathcal{P}})$ a projective class on $\mathcal{C}$ with enough projectives and consider the induced model structure on $\mathrm{s}\mathcal{C}$. 
\begin{enumerate}
\item
$$\otimes:\mathrm{s}\mathrm{Ab}\times\mathrm{s}\mathcal{C}\rightarrow\mathrm{s}\mathcal{C}$$
$$\underline{\Hom}(-,-):(\mathrm{s}\mathcal{C})^{op}\times\mathrm{s}\mathcal{C}\rightarrow\mathrm{s}\mathrm{Ab}$$
$$(-)^{(-)}:\mathrm{s}\mathcal{C}\times(\mathrm{s}\mathrm{Ab})^{op}\rightarrow\mathrm{s}\mathcal{C}$$
is a two-variable Quillen adjunction. In particular the model structure on $\mathrm{s}\mathcal{C}$ is simplicial.
\item
Let $K$ be a finite simplicial set and $f$ a weak equivalence in $\mathrm{s}\mathcal{C}$. Then $K\otimes f$ is a weak equivalence. 
\item
If objects of $\mathcal{P}$ are $(\aleph_{0},\mathbf{SplitMon})$-tiny then $K\otimes f$ is a weak equivalence for any simplicial set $K$.
\end{enumerate}
\end{lem}

\begin{proof}
\begin{enumerate}
\item
It suffices to prove that $\otimes:\mathrm{s}\mathrm{Ab}\times\mathrm{s}\mathcal{C}\rightarrow\mathrm{s}\mathcal{C}$ is a left Quillen bifunctor.  Let $f:X\rightarrow Y$ be a cofibration in $\mathrm{s}\mathrm{Ab}$ and $g:C\rightarrow D$ a cofibration in $\mathrm{s}\mathcal{C}$.  We may assume that $g=P\otimes g'$ for some $g'$ a cofibration in $\mathrm{sAb}$. Then
$$f\otimes g=P\otimes (f\otimes_{\mathbb{Z}} g')$$
where $\otimes_{\mathbb{Z}} $ is the tensor product of simplicial abelian groups. If $f$ and $g'$ are cofibrations then $f\otimes_{\mathbb{Z}} g'$ is a cofibration, and is an acyclic cofibration if either $f$ or $g'$ is. Thus it suffices to show that $P\otimes:\mathrm{sAb}\rightarrow\mathrm{s}\mathcal{C}$ sends cofibrations to cofibrations, and acyclic cofibrations to ayclic cofibrations. We may test this on generating (acyclic) cofibrations, and then by Remark \ref{rem:cofgens} this is clear.
\item
Let $P\in\mathcal{P}$ and $K$ a finite simplicial set. Then $\Hom(P,K\otimes f)\cong K\otimes\Hom(P,f)$. This is a weak equivalence of simplicial abelian groups. 
\item
If objects of $\mathcal{P}$ are $\mathbf{SplitMon}$-tiny, then since any simplicial set is obtained as a transfinite composition of pushouts of maps between finite simplicial sets, one has $\Hom(P,K\otimes f)\cong K\otimes\Hom(P,f)$ for any simplicial set $K$. Again this is a weak equivalence of simplicial abelian groups. 
\end{enumerate}
\end{proof}

\begin{lem}\label{lem:transfinitecompequiv}
Let $\mathcal{C}$ be a complete and cocomplete additive category, and $(\mathcal{P},\mathcal{E}_{\mathcal{P}})$ a projective class on $\mathcal{C}$ with enough projectives. Consider the induced model structure on $\mathrm{s}\mathcal{C}$ and let $\mathcal{Q}$ be a generating set of $(\aleph_{0},\mathbf{SplitMon})^{\mathrm{tiny}}$ projectives. Then transfinite compositions of split monomorphisms are homotopy colimits.
\end{lem}

\begin{proof}
Note first that the transition maps in cofibrant filtered diagrams are cofibrations and in particular are level-wise split monomorphisms. Thus it suffices to prove that if
\begin{displaymath}
\xymatrix{
X_{0}\ar[r]\ar[d]^{f_{0}} & X_{1}\ar[d]^{f_{1}}\ar[r] & \cdots\ar[r] & X_{\beta}\ar[d]^{f_{\beta}}\ar[r] & \cdots\\
Y_{0}\ar[r] & Y_{1}\ar[r] &\cdots\ar[r] & Y_{\beta}\ar[r] & \cdots
}
\end{displaymath}
is a commutative diagram with each map $X_{\beta}\rightarrow X_{\beta+1}$ and $Y_{\beta}\rightarrow Y_{\beta+1}$ being a level-wise split monomorphism, and each $f_{\beta}$ an equivalence, then the induced map $\colim X_{\beta}\rightarrow\colim Y_{\beta}$ is an equivalence. 
be a transfinite sequence with each $X_{\beta}\rightarrow X_{\beta+1}$ being a level-wise split monomorphism. By the fact that $P$ is $(\aleph_{0},\mathbf{SplitMon})^{\mathrm{tiny}}$ the map
$$\colim\Hom(P,X_{\beta})\rightarrow\Hom(P,\colim X_{\beta})$$
is an isomorphism. The claim now follows from the corresponding claim for $\mathrm{sAb}$ which can be proven, for example, by the Dold-Kan correspondence and the fact that filtered colimits in $\mathrm{Ab}$ are exact.
\end{proof}

\begin{prop}\label{prop:PSigmaE}
Let $\mathcal{C}$ be a complete and cocomplete locally presentable additive category, and $(\mathcal{P},\mathcal{E}_{\mathcal{P}})$ a projective class on $\mathcal{C}$ with enough projectives. Let $\mathcal{Q}$ be a generating set of $(\aleph_{0},\mathbf{SplitMon})^{\mathrm{tiny}}$ projectives. Then $\mathrm{L}^{H}(\mathrm{s}\mathcal{C})$ is equivalent to $\mathcal{P}_{\Sigma}(N(\mathcal{Q}))$.
\end{prop}

\begin{proof}
The same proof as in \cite{kelly2021analytic} Proposition 4.27 shows that any object of $\mathcal{P}$ is projective as an object of $\mathrm{L}^{H}(\mathrm{s}\mathcal{C})$, and that any object of $\mathrm{s}\mathcal{C}$ may be written as a homotopy sifted colimit of objects of $\mathcal{Q}$. Moreover since cofibrant filtered diagrams in $\mathrm{s}\mathcal{C}$ have cofibrant, and hence degree-wise split, maps between objects, any object of $\mathcal{Q}$ is also tiny in $\mathrm{L}^{H}(\mathrm{s}\mathcal{C})$. This completes the proof. 
\end{proof}

\subsubsection{Monoidal Structures}

Here we fix a complete and cocomplete additive category $\mathpzc{C}$ equipped with a closed symmetric monoidal structure $(k,\ootimes,\underline{\mathrm{Hom}})$, and a projective class $(\mathcal{P},\mathcal{E}_{\mathcal{P}})$ which has enough projectives.

\begin{defn}
    An object $Q$ of $\mathpzc{E}$ is said to be \textit{flat} if $Q\otimes(-)$ sends exact sequences to exact sequences, and \textit{strongly flat} if it is flat, and $Q\otimes(-)$ commutes with finite limits.
\end{defn}

\begin{defn}
$(\mathcal{C},\mathcal{P},\mathcal{E}_{\mathcal{P}})$ is said to be \textit{quasi-projectively monoidal} if 
\begin{enumerate}
\item
$k\in\mathcal{P}$.
\item
If $P,P'\in\mathcal{P}$ then $P\ootimes P'$ is in $\mathcal{P}$.
\end{enumerate}
$(\mathpzc{C},\mathcal{P},\mathcal{E}_{\mathcal{P}})$ is said to be \textit{projective monoidal} if in addition projectives are flat, and \textit{strongly projectively monoidal} if projectives are strongly flat. 
\end{defn}

\begin{prop}
Let $\mathcal{Q}$ be a projective generating subcategory for $\mathpzc{C}$. Then $\mathpzc{C}$ is 
\begin{enumerate}
\item
quasi-projectively monoidal if and only if $k\in C(\mathcal{Q})$ and for any $Q,Q'\in\mathcal{Q}$, $Q\ootimes Q'\in C(\mathcal{Q})$.
\item
projectively monoidal if and only if it is quasi-projectively monoidal and for each $Q\in\mathcal{Q}$, $Q$ is flat.
\item
strongly projectively monoidal if and only if it is quasi-projectively monoidal, coproducts are strongly exact, and for each $Q\in\mathcal{Q}$, $Q$ is strongly flat. 
\end{enumerate}
\end{prop}

\begin{proof}
\begin{enumerate}
\item
Since $\mathcal{Q}$ is a projective generating subcategory, by Lemma \ref{lem:gencproj} we have that $C(\mathcal{Q})=\mathcal{P}$. Thus if $\mathcal{C}$ is quasi-projectively monoidal then $k\in C(\mathcal{Q})$ and $Q\ootimes Q'\in C(\mathcal{Q})$ whenever $Q,Q'\in\mathcal{Q}\subset\mathcal{P}$. 

Conversely suppose $k\in C(\mathcal{Q})$ and $Q\ootimes Q'\in C(\mathcal{Q})$ whenever $Q,Q'\in\mathcal{Q}$. Let $P$ and $P'$ be projectives. Then $P$ is a summand of an object of the form $\bigoplus_{i\in\mathcal{I}}Q_{i}$ and $P'$ is a summand of an object of the form $\bigoplus_{j\in\mathcal{J}}Q'_{j}$ with each $Q_{i}$ and $Q'_{j}\in\mathcal{Q}$. Thus $P\ootimes P'$ is a summand of $\bigoplus_{(i,j)\in\mathcal{I}\times\mathcal{J}}Q_{i}\ootimes Q'_{j}$. Now $Q_{i}\ootimes Q'_{j}$ is in $\mathcal{P}$ and therefore $\bigoplus_{(i,j)\in\mathcal{I}\times\mathcal{J}}Q_{i}\ootimes Q'_{j}\in\mathcal{P}$. Hence $P\ootimes P'\in\mathcal{P}$.
\item
This is immediate from the first part, and that fact that the class of flat objects is closed under coproducts and summands.
\item
This is immediate from the first part, and that fact that the class of strongly flat objects is closed under coproducts and summands.
\end{enumerate}
\end{proof}

The following definitions and results can be found in \cite{kelly2021analytic} Section 4.
\begin{defn}
$(\mathcal{C},\mathcal{P},\mathcal{E}_{\mathcal{P}})$ is said to have \textit{symmetric projectives} if for any projective $P\in\mathcal{E}$ and any $n\in\mathbb{N}$, $\mathsf{Sym}_{\mathcal{E}}^{n}(P)$ is projective.
\end{defn}

\begin{rem}
If $(\mathcal{C},\mathcal{P},\mathcal{E}_{\mathcal{P}})$ is quasi-projectively monoidal and enriched over $\mathbb{Q}$ then it has symmetric projectives. Indeed in this case $\mathsf{Sym}^{n}(P)$ is a summand of $P^{\otimes n}$, which is projective. 
\end{rem}

\begin{defn}
Call a projective generating set \(\mathcal{P}\) in $(\mathcal{C},\mathcal{P},\mathcal{E}_{\mathcal{P}})$ \textit{symmetrically closed} if the full subcategory on $\mathcal{P}$ is closed under both tensor and symmetric powers, with respect to the tensor product of \(\mathcal{C}\).
\end{defn}


\begin{prop}\label{prop:projclosed}
Let $(\mathcal{C},\mathcal{P},\mathcal{E}_{\mathcal{P}})$ be quasi-projectively monoidal elementary with symmetric projectives. Suppose further that the tensor product of two tiny objects is tiny. Then there exists a set $\mathcal{P}$ of tiny projective generators of $\mathcal{E}$ such that 
\begin{enumerate}
\item
for any $P\in\mathcal{P}$ and any $n\in\mathbb{N}$, $\mathsf{Sym}_{\mathcal{E}}^{n}(P)\in\mathcal{P}$
\item
for any $P,Q\in\mathcal{P}$, $P\otimes Q\in\mathcal{P}$
\item
$\mathcal{P}$ is closed under finite coproducts.
\end{enumerate}
\end{prop}

\begin{thm}
Let $(\mathcal{C},\mathcal{P},\mathcal{E}_{\mathcal{P}})$ be quasi-projectively monoidal. 
\begin{enumerate}
\item
The projective model structure on $\mathrm{s}\mathcal{C}$ is monoidal and satisfies the monoid axiom. 
\item
 If $(\mathcal{C},\mathcal{P},\mathcal{E}_{\mathcal{P}})$ has symmetric projectives then the strong commutative monoid axiom is satisfied.
 \item
If $(\mathcal{C},\mathcal{P},\mathcal{E}_{\mathcal{P}})$  is projectively monoidal and has a generating set consisting of $(\aleph_{0},\mathbf{SplitMon})^{\mathrm{tiny}}$ projectives, then cofibrant objects are $K$-flat.
 \end{enumerate}
\end{thm}

\begin{proof}
First let us check that it is monoidal. The unit axiom is clear. For the pushout-product axiom let $f:X\rightarrow Y$ and $f':X'\rightarrow Y'$ be cofibrations. Without loss of generality we may assume that $f$ and $f'$ are generating cofibrations. In particular they are of the form $f=Id_{Z}\ootimes g$, $f'=Id_{Z}\ootimes g'$ where $Z,Z\in\mathcal{Z}$ and $g:W\rightarrow U,g':W'\rightarrow U'$ are generating cofibrations of simplicial sets. Then the map
$$X\ootimes Y'\coprod_{X\ootimes X'}X\ootimes Y\rightarrow X'\ootimes Y'$$
is isomorphic to the map
$$Z\ootimes Z'\ootimes ((W\times U')\coprod_{W\times W'}(W'\times U))\rightarrow (Z\ootimes Z')\ootimes (U'\times U)$$
The pushout-product axiom is known for simplicial sets so \[(W\times U')\coprod_{W\times W'}(W'\times U)\rightarrow U'\times U\] is a cofibration. Since $Z\ootimes Z'\in C(\mathcal{Z})$, \[Z\ootimes Z'\ootimes((W\times U')\coprod_{W\times W'}(W'\times U))\rightarrow (Z\ootimes Z')\ootimes (U'\times U)\] is a cofibration in $\textrm{s}\mathpzc{C}$, which is an acyclic cofibration whenever either $g$ or $g'$ is.

Now let us prove the monoid axiom. First observe that if $f$ is a generating acyclic cofibration in $\mathpzc{sC}$ then it is a homotopy equivalence, so if $X$ is any object in $\mathpzc{sC}$ then $Id_{X}\ootimes f$ is a degree-wise split homotopy equivalence. The pushout of a degree-wise split homotopy equivalence is itself a degreewise split homotopy equivalence. By Lemma \ref{lem:transfinitecompequiv}, any transfinite composition of degree-wise split weak equivalences is a degree-wise split weak equivalence. 

Finally let us prove that if $\mathcal{C}$ is projectively monoidal then cofibrant objects are $K$-flat. Let $X$ be a trivial object of $\mathrm{s}\mathcal{C}$, and $C$ a cofibrant object. $C$ may be written as a transfinite composition of pushouts along maps of the form
$$P\otimes\partial\Delta^{n}\rightarrow P\otimes\Delta^{n}$$
where $P$ is projective and $n\in\mathbb{N}$.
It suffices to show that $P\otimes K$ is $K$-flat for any (finite) simplicial set $K$. Let $f:X\rightarrow Y$ be a weak equivalence. By Lemma \ref{lem:tensoredoversimp} $K\otimes f$ is a weak equivalence. Thus it suffices to show that $P\otimes f$ is a weak equivalence. This follows immediately from the Dold-Kan correspondence, and the fact that $P$ is flat. 
\end{proof}

In particular, this implies the following.

\begin{thm}\label{thm:HAcontextsimplicial}
Let $(\mathcal{P},\mathcal{E}_{\mathcal{P}})$ be a projective class on a locally presentably closed symmetric monoidal additive category $\mathcal{C}$. Suppose that $(\mathcal{C},\mathcal{P},\mathcal{E}_{\mathcal{P}})$ is projectively monoidal, has symmetric projectives, and is $\mathbf{SplitMon}$-elementary. Then $\mathrm{s}\mathcal{C}$ is a HA context.
\end{thm}

It is also a derived algebraic context.

\begin{thm}[\cite{kelly2021analytic} Theorem 4.30]
    Let $\mathpzc{E}$ be a quasi-projectively monoidal elementary exact category. Then $$(\mathbf{Stab}(\mathbf{s}\mathpzc{E}),\mathbf{Stab}_{\ge0}(\mathbf{s}\mathpzc{E}),\mathbf{Stab}_{\le0}(\mathbf{S}\mathpzc{E},\mathcal{P})$$
    is a derived algebraic context, where $\mathcal{P}$ is a generating set of tiny projective generators closed under finite direct sums, tensor products, and symmetric powers.
\end{thm}

\begin{rem}
Using Remark \ref{rem:cofgens} it is easy to write down generating sets of cofibrations and of acyclic cofibrations in $\mathrm{Comm(sC)}$. Indeed for $K$ a simplicial set let $\mathbb{Z}[K]$ denote the simplicial unital commutative ring which in degree $n$ is the polynomial algebra on $K_{n}$. Let $\mathcal{Q}$ be a generating set of projectives. Then the generating sets are

  \[\mathrm{Sym}(I) = 
  \{ \mathrm{Sym}(P) \otimes \mathbb{Z}[\partial \Delta[n]] 
  \longrightarrow 
  \mathrm{Sym}(P) \otimes \mathbb{Z}[\Delta[n]] \ \ | \ \  P \in \mathcal{Q} , n \geq 0 \} 
  \]
  and
  \[\mathrm{Sym}(J) = \{ \mathrm{Sym}(P) \otimes \mathbb{Z}[V[n,k]] \longrightarrow \mathrm{Sym}(P)\otimes \mathbb{Z}[\Delta[n]] \ \  | \ \  P \in \mathcal{Q} , n \geq k \geq 0, n>0 \}.
  \]
  \end{rem}

\subsubsection{Perfect, Compact, and Nuclear Objects}

Consider a monoidal elementary exact category of the form $\mathpzc{E}$ with tensor product $\otimes$ and monoidal unit $\mathbb{I}$. Let $\mathcal{Q}$ be a set of compact projective generators for $\mathpzc{E}$. Consider the symmetric monoidal $(\infty,1)$-category $\mathbf{Ch}(\mathpzc{E})$.  By definition we have that $\mathbf{Perf}(\mathbf{Ch}(\mathrm{Ind}(\mathpzc))$ consists of complexes equivalent to bounded complexes $X_{\bullet}$ such that each $X_{n}$ is a retract of $\mathbb{I}^{\oplus m_{n}}$ for some $m_{n}$. 

By an identical argument to \cite{stacks-project} Proposition 15.78.3 we have the following. 

\begin{lem}
$\mathbf{Cpct}_{\aleph_{0}}(\mathbf{Ch}(\mathpzc{E}))$ consists of complexes which are equivalent to bounded complexes of compact projectives.
\end{lem}

It follows that $\mathbf{Perf}(\mathpzc{E})=\mathbf{Cpct}_{\aleph_{0}}(\mathbf{Ch}(\mathpzc{E})$ precisely if $\mathbb{I}$ is a (compact) projective generator for $\mathpzc{E}$. This would imply that there is an exact symmetric monoidal equivalence ${}_{R}\mathrm{Mod}(\mathrm{Ab})\cong\mathpzc{E}^{\heart}$, where $R\defeq\mathrm{End}(\mathbb{I})$ and $\mathpzc{E}^{\heart}$ is the heart of the projective $t$-structure \cite{kelly2021analytic}. 

\begin{lem}
$\mathbf{Dls}(\mathbf{Ch}(\mathpzc{E})$ consists of those complexes which are equivalent to bounded complexes of compact projectives $K_{\bullet}$, such that each $K_{n}$ is dualisable as an object of $\mathpzc{E}$, and $K_{n}^{\vee}$ is flat.
\end{lem}

\begin{proof}
Let $K_{\bullet}\in\mathbf{Dls}(\mathbf{Ch}(\mathpzc{E})$. In particular it is compact so we may assume it is a bounded complex of compact projectives. By shifting we may assume that $K_{\bullet}$ is concentrated in the interval $[b,a]$. The proof is by induction on $n=b-a$. Let $n=0$, so that $K_{\bullet}=K_{n}[n]$. Without loss of generality we may assume that $n=0$ and we write $K=K_{0}$. Let $X_{\bullet}$ be an aribtrary complex. We have
\begin{align*}
\underline{\mathrm{Hom}}(K,X_{\bullet})&\cong\mathbb{R}\underline{\mathrm{Hom}}(K ,X_{\bullet})\\
&\cong \mathbb{R}\underline{\mathrm{Hom}}(K,\mathbb{I})\otimes^{\mathbb{L}}X_{\bullet}\\
&\cong K^{\vee}\otimes^{\mathbb{L}}X_{\bullet}
\end{align*}
Now for $X_{\bullet}$ concentrated in degree $0$ we find that $K$ is dualisable as an object of $\mathpzc{E}$, and $K^{\vee}$ is flat. Conversely suppose that $K$ is dualisable in $\mathpzc{E}$ and $K^{\vee}$ is flat. Then a similar chain of isomorphisms to the above gives that $K$ is in $\mathbf{Dls}(\mathbf{Ch}(\mathpzc{E}))$. 

Suppose we have proven the claim whenever $K'_{\bullet}$ is concentrated in $[b-1,a]$. Let $K_{\bullet}$ be concentrated in $[b,a]$. We have a fibre-cofibre sequence
$$K'_{\bullet}\rightarrow K_{\bullet}\rightarrow K_{b}[b]$$
By the inductive hypotheses the claim is true for $K'_{\bullet}$, and by the base case it is true for $K_{b}[b]$. Since the validity of the claim is stable under extensions, it is also true for $K_{\bullet}$. 
\end{proof}

\subsubsection{Nuclear Maps and Nuclear Objects}

Here we discuss nuclearity in closed symmetric monoidal elemenatry exact categories $\mathpzc{E}$.

\begin{defn}[\cite{ben2020fr} Definition 4.11]
An object $W$ of $\mathrm{Ind}(\mathpzc{E})$ is said to be 
\begin{enumerate}
\item
\textit{nuclear} if for any compact object $V\in\mathpzc{E}$ the map
$$\mathrm{Hom}(\mathbb{I},V^{\vee}\otimes W)\rightarrow\mathrm{Hom}(V,W)$$
is a bijection.
\item
\textit{projectively nuclear} if for any compact projective object $V\in\mathpzc{E}$ the map
$$\mathrm{Hom}(\mathbb{I},V^{\vee}\otimes W)\rightarrow\mathrm{Hom}(V,W)$$
is a bijection.
\item
\textit{strongly nuclear} if for any compact object $V\in\mathpzc{E}$ the map
$$V^{\vee}\otimes W\rightarrow\underline{\mathrm{Hom}}(V,W)$$
is an isomorphism. 
\item
\textit{strongly projectively nuclear} if for any compact projective object $V\in\mathpzc{E}$ the map
$$V^{\vee}\otimes W\rightarrow\underline{\mathrm{Hom}}(V,W)$$
is an isomorphism. 
\end{enumerate}
\end{defn}


Nuclear objects arise in the following way.

\begin{defn}[\cite{ben2020fr} Definition 4.2]
A map $f:E\rightarrow F$ in $\mathpzc{E}$ is said to be \textit{nuclear} if it is in the image of
$$\mathrm{Hom}(\mathbb{I},E^{\vee}\otimes F)\rightarrow\mathrm{Hom}(E,F)$$
\end{defn}

\begin{lem}[\cite{ben2020fr} Lemma 4.16]\label{lem:nucdiscind}
An object $V$ of $\mathrm{Ind}(\mathpzc{C})$ is strongly nuclear if and only if it has a presentation $V\cong``\limind_{\mathcal{I}}"V_{i}$ such that each structure map $V_{i}\rightarrow V_{j}$ is a nuclear map in $\mathpzc{C}$.
\end{lem}


\subsubsection{$(\infty,1)$-Categories of Nuclear and Strong Nuclear Objects}

Let 
$$\mathpzc{C}=(\mathbf{Ch}(\mathpzc{E}),\mathbf{Ch}_{\ge0}(\mathpzc{E}),\mathbf{Ch}_{\le0}(\mathpzc{E}))$$
be a derived algebraic contexts. Let us analyse nuclear objects in such categories. The following is clear.

\begin{prop}\label{prop:basicnuclearderived}
    Let $X$ be flat as an object of $\mathpzc{E}$. Then $X$ is strongly projectively nuclear as an object of $\mathbf{Ch}(\mathpzc{E})$ if and only if it is strongly nuclear as an object of $\mathpzc{E}$.
\end{prop}

\begin{prop}
A basic nuclear object $X$ in $\mathbf{C}$ has a presentation
$$X_{\bullet}\cong\limind_{n}X^{n}_{\bullet}$$
where
\begin{enumerate}
    \item 
    Each $X^{n}_{\bullet}$ is a bounded complex of compact projectives.
    \item 
    Each $X^{n}_{m}\rightarrow X^{n+1}_{m}$ is trace class.
\end{enumerate}
In particular each $X_{m}$ is itself a (discrete) basic nuclear object. 
\end{prop}

\begin{proof}
    The proof of \cite{scholzecomplexlectures} Lemma 8.7 goes through identically, noting the important facts that $\mathbf{Map}(P,\mathbb{I})$ is concentrated in degree $0$ for each projective $P$, and that projectives are flat.
\end{proof}

\begin{cor}
Let $\mathbf{C}$ be a derived algebraic context, and let $A\in\mathbf{Ass}(\mathbf{C})$ be an algebra whose underlying object is nuclear. Suppose that discrete basic nuclear objects are strongly projectively nuclear in $\mathbf{C}^{\heart}$. Then any nuclear object of ${}_{A}\mathbf{Mod}(\mathbf{C})$ is strongly nuclear.
\end{cor}

\begin{proof}
First assume that $A=\mathbb{I}$, let $X$ be a nuclear object, and let $c$ be a compact object. We may assume that $c$ is a compact projective concentrated in degree $0$. Since the functors $\underline{\mathbf{Map}}(c,-)$ and $c^{\vee}\otimes^{\mathbb{L}}(-)$ both commute with filtered colimits, we may assume that $X$ is basic nuclear. Since each $X_{n}$ is then basic nuclear, we may also assume that $X$ is discrete. Since discrete basic nuclear objects are flat, and are strongly projetively nuclear by assumption, the claim follows from Proposition \ref{prop:basicnuclearderived}.

Now let $M$ be a nuclear object of ${}_{A}\mathbf{Mod}(\mathbf{C})$. Then $M$ is a nuclear object of $\mathbf{C}$ by Lemma \ref{lem:nucfacts}. We then have
\begin{align*}
    \underline{\mathbf{Map}}_{{}_{A}\mathbf{Mod}}(A\otimes^{\mathbb{L}}c,M)&\cong  \underline{\mathbf{Map}}_{\mathbf{C}}(c,M)\\
    &\cong c^{\vee}\otimes^{\mathbb{L}}M\\
    &\cong c^{\vee}\otimes^{\mathbb{L}}A\otimes_{A}^{\mathbb{L}}M\\
    &\cong \underline{\mathbf{Map}}_{{}_{A}\mathbf{Mod}}(A\otimes^{\mathbb{L}}c,A)\otimes_{A}^{\mathbb{L}}M
\end{align*}
as required.
\end{proof}

By Lemma \ref{lem:nucdiscind} we get the following.

\begin{cor}\label{cor:strongnucnuc}
    Let $\mathpzc{E}$ be a closed monoidal quasi-abelian category with enough projectives, such that tensor products and symmetric powers of projectives are projectives. Consider the derived algebraic context
    $$(\mathbf{Ch}(\mathrm{Ind}(\mathpzc{E})),\mathbf{Ch}_{\ge0}(\mathrm{Ind}(\mathpzc{E})),\mathbf{Ch}_{\le0}(\mathrm{Ind}(\mathpzc{E})),\mathcal{P}^{0})$$
    Let $A$ be a strongly nuclear object of $\mathbf{Ass}(\mathbf{Ch}(\mathrm{Ind}(\mathpzc{E}))$. Then an $A$-module $M$ is strongly nuclear if and only if it is nuclear.
\end{cor}

\subsection{Commuting Tensor Products and Limits}

For many descent results later, we will need some tensor products to commute with some infinite products.

\begin{defn}\label{defn:formallykfilt}
Let $(\mathbf{C},\otimes)$ be a symmetric monoidal $(\infty,1)$-category, and let $\kappa$ be a cardinal. Let $\mathbf{D}\subset\mathbf{C}$ be a subcategory. 
\begin{enumerate}
\item
An object $c$ of $\mathbf{C}$ is said to be \textit{formally }$\kappa$-\textit{filtered relative to }$\mathbf{D}$ if for any $\kappa$-small family $\{d_{i}\}_{i\in\mathcal{I}}$ of objects of $\mathbf{D}$ whose product exists, the map
$$c\otimes^{\mathbb{L}}\prod_{i\in\mathcal{I}}d_{i}\rightarrow\prod_{i\in\mathcal{I}}c\otimes^{\mathbb{L}} d_{i}$$
is an equivalence. 
The class of formally $\kappa$-filtered objects is denoted $\mathbf{C}^{\kappa}$
\item
An object $c$ of $\mathbf{C}$ is said to be \textit{strongly formally }$\kappa$-\textit{filtered relative to }$\mathbf{D}$ if for any $\kappa$-small diagram $d:\mathcal{I}\rightarrow\mathrm{D}$ whose limit exists, the map
$$c\otimes^{\mathbb{L}}\lim_{i\in\mathcal{I}}d_{i}\rightarrow\lim_{i\in\mathcal{I}}c\otimes^{\mathbb{L}} d_{i}$$
is an equivalence. 
\end{enumerate}
An object is said to be (strongly) formally $\kappa$-filtered if it is (strongly) formally $\kappa$-filtered relative to $\mathrm{Ob}(\mathbf{C})$
\end{defn}

\begin{rem}
If $\mathbf{C}$ is stable then any formally $\kappa$-filtered object (relative to any subcategory) is strongly formally $\kappa$-filtered (relative to that same subcategory).
\end{rem}

The obvious corresponding notion for $1$-categories is also useful.

\begin{defn}\label{defn:1formallykappafiltered}
Let $\mathpzc{C}$ be a $1$-category  and $\mathpzc{D}\subset\mathrm{C}$ a subcategory. An object $c$ of $\mathpzc{C}$ is said to be\textit{formally }$\kappa$-\textit{filtered relative to }$\mathpzc{D}$ if for any $\kappa$-small family $\{d_{i}\}_{i\in\mathcal{I}}$ of objects of $\mathpzc{D}$ whose product exists, the map
$$c\otimes\prod_{i\in\mathcal{I}}d_{i}\rightarrow\prod_{i\in\mathcal{I}}c\otimes d_{i}$$
is an isomorphism. 
\end{defn}


\begin{lem}
Let $\mathbf{C}$ be an $E_{\infty}$-algebraic context with exact countable products. Let $F$ be a homologically bounded object such that each $\pi_{n}(F)$ is formally $\kappa$-filtered relative to $\mathbf{C}^{\heart}$. Suppose further each  $\pi_{n}(F)$ has a resolution by formally $\kappa$-filtered flat objects. Let $\mathbf{C}^{\otimes-F}$ denote the full subcategory consisting of objects $X$ such that $\pi_{m}(X)\otimes^{\mathbb{L}}\pi_{n}(F)\cong \pi_{m}(X)\otimes\pi_{n}(F)$ for all $n,m\in\mathbb{Z}$. Then $F$ is formally $\kappa$-filtered relative to  $\mathbf{C}^{\otimes-F}$.

\end{lem}

\begin{proof}
Since we are working in an $E_{\infty}$-algebraic context, it is equivalent to $\mathbf{Ch}(\mathbf{C}^{\heart})$, so we can just work with chain complexes. By an easy induction we may assume that $F$ is in fact concentrated in degree $0$. 

We first claim that $F$ is transverse to each $\pi_{*}(\prod_{i} X^{i}_{\bullet})\cong\prod_{i}\pi_{*}(X^{i}_{\bullet})$. For this we may assume also that each $X^{i}$ is concentrated in degree $0$. Pick a resolution $P_{\bullet}\rightarrow F$ by flat formally $\kappa$-filtered objects. We then have
\begin{align*}
F\otimes^{\mathbb{L}}\prod_{i}X^{i}&\cong P_{\bullet}\otimes\prod_{i} X^{i}\\
&\cong\prod_{i}P_{\bullet}\otimes X^{i}\\
&\cong \prod_{i}F\otimes^{\mathbb{L}}X^{i}\\
&\cong\prod_{i}F\otimes X^{i}
\end{align*}
Using the Tor spectral sequence we have 
$$\pi_{p}(F\otimes^{\mathbb{L}}\prod_{i} X^{i}_{\bullet})\cong F\otimes\prod_{i}\pi_{p}(X^{i}_{\bullet})\cong\prod_{i}(F\otimes\pi_{p}(X^{i}_{\bullet}))$$
But the tor spectral sequence also implies $F\otimes\pi_{p}(X^{i}_{\bullet})\cong\pi_{p}(F\otimes^{\mathbb{L}}X^{i}_{\bullet})$. Thus we have
$$\pi_{*}(F\otimes^{\mathbb{L}}\prod_{i} X^{i}_{\bullet})\cong\pi_{*}(\prod_{i}(F\otimes^{\mathbb{L}}X^{i}_{\bullet})$$

\end{proof}

\begin{cor}
Let $F^{\bullet}$ be a bounded complex of formally $\kappa$-filtered objects in $\mathpzc{A}$. Suppose that each $F^{n}$ has a resolution by a bounded complex of flat objects. Then $F^{\bullet}$ is formally $\kappa$-filtered. 
\end{cor}


\begin{cor}
Let $X^{\bullet}$ be a complex which has a bounded resolution by formally $\aleph_{1}$-filtered flat objects, and let $G:\Delta\rightarrow\mathrm{Ch}(\mathcal{A})$ be a cosimplicial object. Then the map
$$X^{\bullet}\otimes^{\mathbb{L}}\mathbf{lim}\;G\rightarrow\mathbf{lim}X^{\bullet}\otimes^{\mathbb{L}}G$$
is an equivalence.
\end{cor}


\subsubsection{$\kappa$-Filtered Objects}

Let $\mathpzc{E}$ be a finitely complete and cocomplete  closed symmetric monoidal exact category with enough projectives. Consider the monoidal elementary exact category $\mathrm{Ind}(\mathpzc{E})$. 

\begin{defn}\label{defn:metrisable}
An object $F$ of $\mathrm{Ind}(\mathpzc{E})$  is said to be $\kappa$-metrisable if it can be written as 
$$M\cong\mathrm{colim}_{\Lambda}( E^{\lambda})$$
where $\Lambda$ is $\kappa$-filtered and $E^{\lambda}\in\mathpzc{E}$. 
\end{defn}

The following is a formal consequence of commutativity of $\kappa$-filtered colimits with $\kappa$-small limits.

\begin{prop}
If objects of $\mathpzc{E}$ are formally $\kappa$-filtered in $\mathrm{Ind}(\mathpzc{E})$ then $\kappa$-metrisable objects are formally $\kappa$-filtered in $\mathrm{Ind}(\mathpzc{E})$.
\end{prop}

\begin{defn}
Let $A\in\mathrm{Ch}(\mathrm{Ind}(\mathpzc{E}))$. An object $M$ of ${}_{A}\Mod$ is said to be 
\begin{enumerate}
\item
\textit{weakly }$\kappa$-\textit{metrisable} if it can be written as 
$$M\cong\mathrm{colim}_{\Lambda}(A\otimes P^{\lambda})$$
where $\Lambda$ is $\kappa$-filtered and $P^{\lambda}$ is a  complex of flat objects in $\mathpzc{E}$.
\item
$\kappa$-\textit{metrisable}  if it can be written as 
$$M\cong\mathrm{colim}_{\Lambda}(A\otimes P^{\lambda})$$
where $\Lambda$ is $\kappa$-filtered and $P^{\lambda}$ is a bounded complex of flat objects in $\mathpzc{E}$.
\end{enumerate}
\end{defn}

\begin{lem}
Suppose objects of $\mathpzc{E}$ are formally $\kappa$-filtered in $\mathrm{Ind}(\mathpzc{E})$.
\begin{enumerate}
\item
$\kappa$-metrisable objects of ${}_{A}\Mod$ are formally $\kappa$-filtered in both ${}_{A}\Mod$ and the $(\infty,1)$-category ${}_{A}\mathbf{Mod}$. 
\item
If $M\in{}_{A}\Mod$ is weakly $\kappa$-metrisable, and $\{F_{n}\}$ is a $\kappa$-small collection of $A$-modules all of which are homologically concentrated in a fixed interval $[a,b]$, then the map
$$M\otimes^{\mathbb{L}}_{A}\prod_{n}F_{n}\rightarrow \prod_{n}(M\otimes^{\mathbb{L}}_{A}F_{n})$$
is an equivalence. 
\end{enumerate}
\end{lem}

\begin{proof}
The first claim works exactly as in the $1$-category case. The second is similar - the boundedness assumptions on $\{F_{n}\}$ means we can still commute tensor products with infinite products.
\end{proof}

\begin{prop}
Let $f:A\rightarrow B$ be a map in $\Comm(\mathrm{Ch}(\mathpzc{E}))$. If $M$ is a (weakly) $\kappa$-metrisable $A$-module then $B\otimes^{\mathbb{L}}_{A}M\cong B\otimes_{A}M$ is (weakly) $\kappa$-metrisable. If $N$ is (weakly) $\kappa$-metrisable as a $B$-module, and $B$ is (weakly) $\kappa$-metrisable as an $A$-module, then $N$ is (weakly) $\kappa$-metrisable as an $A$-module.
\end{prop}

\begin{proof}
Write 
$$M=\mathrm{colim}_{\Lambda}(A\otimes P^{\lambda})$$
 Note that $M$ is $K$-flat as an $A$-module. Then
$$B\otimes^{\mathbb{L}}_{A}M\cong B\otimes_{A}M\cong\mathrm{colim}_{\Lambda}(B\otimes P^{\lambda})$$
is (weakly) $\kappa$-metrisable. 

Conversely let $N$ be (weakly) $\kappa$-metrisable as a $B$-module, and $B$  (weakly) $\kappa$-filtered as an $A$-module. Write 
$$B=\mathrm{colim}_{\Lambda}A\otimes P^{\lambda}$$
$$N=\mathrm{colim}_{\Lambda'}B\otimes Q^{\lambda'}$$
Then 
$$N\cong\mathrm{colim}_{\Lambda\times\Lambda'}A\otimes P^{\lambda}\otimes Q^{\lambda'}$$
is (weakly) $\kappa$-metrisable as an $A$-module.
\end{proof}

\begin{lem}
Suppose that $A$ is formally $\kappa$-filtered as an object of $\mathbf{Ch}(\mathrm{Ind}(\mathpzc{E}))$. If $\kappa$-small products of projective objects of $\mathpzc{E}$ are formally $\kappa$-filtered in $\mathrm{Ind}(\mathpzc{E})$, then $\kappa$-small products of weakly $\kappa$-metrisable objects of ${}_{A}\mathbf{Mod}$ are weakly $\kappa$-metrisable in ${}_{A}\mathbf{Mod}$.
\end{lem}

\begin{proof}
Let $\{\limind_{\Lambda_{n}}A\otimes P_{\lambda_{n}}\}$ be a $\kappa$-small collection of weakly $\kappa$-metrisable $A$-modules. Without loss of generality we may assume that there is a $\kappa$-filtered $\Lambda$ with $\Lambda=\Lambda_{n}$ for all $n$. We then write $P_{\lambda_{n}}=P^{n}_{\lambda}$. Then
$$\prod_{n}\limind_{\Lambda}A\otimes P^{n}_{\lambda}\cong\limind_{\Lambda}\prod_{n}A\otimes P^{n}_{\lambda}\cong\limind_{\Lambda}A\otimes \prod_{n}P^{n}_{\lambda}$$
The result now follows since each $\prod_{n}P^{n}_{\lambda}$ is weakly $\kappa$-metrisable.
\end{proof}

\begin{rem}
    The terminology `metrisable' here comes from Meyer's work on bornological analysis \cite{meyer2004bornological}.
\end{rem}

\subsection{Monoidal Exact Categories and Derived Algebraic Contexts}
Let $(\mathcal{C},\mathcal{P},\mathcal{E}_{\mathcal{P}})$ be a complete and cocomplete monoidal additive category equipped with a projective class, which is $\mathbf{SplitMon}$-elementary, monoidal elementary, and has symmetric projectives. 
Consider the derived algebraic context 
$$\mathbf{C}(\mathcal{Q})$$
As explained in Proposition \ref{prop:PSigmaE}, $\mathbf{C}_{\le0}(\mathcal{Q})$ is equivalent to $\mathrm{L}^{H}(\mathrm{s}\mathcal{C})$.

\begin{prop}
$\mathbf{DAlg}^{cn}(\mathbf{C})$ is equivalent to $\mathrm{L}^{H}(\mathrm{Comm}(\mathrm{s}\mathcal{C}))$.
\end{prop}

\begin{proof}
They are both equivalent to the $(\infty,1)$-categorical free sifted cocompletion of the subcategory $\mathrm{Sym}(\mathcal{Q})$ of $\mathrm{Comm}(\mathcal{Q})$ consisting of free commutative monoids on objects of $\mathcal{Q}$. 
\end{proof}

Similarly, the category $\mathbf{Comm}^{cn}(\mathbf{C})$ has the following presentation. Let $E_{\infty}$ be a cofibrant resolution of the operad $\mathpzc{Comm}$ in $\mathrm{s}\mathcal{C}$. The transferred model structured exists on the category $\mathrm{Alg}_{E_{\infty}}(\mathrm{s}\mathcal{C})$ and there is an equivalence $\mathbf{Comm}^{cn}(\mathbf{C})\cong\mathrm{L}^{H}(\mathrm{Alg}_{E_{\infty}}(\mathrm{s}\mathcal{C}))$. The map of operads $\phi:E_{\infty}\rightarrow\mathpzc{Comm}$ determines a Quillen adjunction
$$\adj{\phi_{!}}{\mathrm{Alg}_{E_{\infty}}(\mathrm{s}\mathcal{C})}{\mathrm{Alg}_{\mathpzc{Comm}}(\mathrm{s}\mathcal{C})}{\phi^{*}}$$

\begin{prop}
If $\mathcal{C}$ is enriched over $\mathbb{Q}$ then this adjunction is a Quillen equivalence.
\end{prop}

\begin{proof}
This follows from \cite{Axiomatic homotopy theory for operads} Theorem 4.1 and the fact that over $\mathbb{Q}$ the commutative operad is $\Sigma$-cofibrant.
\end{proof}

\begin{rem}
If $\mathcal{C}$ is enriched over $\mathbb{Q}$ then $\mathpzc{C}$ can also be presented by the projective model structure on $\mathrm{Ch}(\mathpzc{E})$ (see. e.g. \cite{kelly2016homotopy}). Again $\mathrm{Comm}(-)$ is an admissible. We have $\mathbf{DAlg}(\mathbf{C})\cong\mathrm{L^{H}}(\mathrm{Alg}_{\mathbf{Comm}}(\mathbf{C}))$ and $\mathbf{Comm}(\mathbf{C})\cong\mathrm{L^{H}}(\mathrm{Alg}_{E_{\infty}}(\mathbf{C}))$, and the map
$$\Theta: \mathbf{DAlg}(\mathbf{C})\rightarrow\mathbf{Comm}(\mathbf{C})$$
is an equivalence.
\end{rem}

\begin{prop}
If $(\mathcal{C},\mathcal{P},\mathcal{E}_{\mathcal{P}})$ is strongly quasi-rojectively monoidal elementary then $\mathbf{C}(\mathcal{Q})$ is a flat derived algebraic context.
\end{prop}

\begin{proof}
Since each $Q\in\mathcal{Q}$ is strongly flat, $Q\otimes(-)$ clearly commutes with the truncation functors for the left $t$-structure on $\mathcal{E}$. It follows immediately that it also commutes with the homology functors for the $t$-structure. 
\end{proof}

\subsubsection{Derived Quotients and The Koszul Complex}
Here we consider derived quotients in the category $\mathbf{DAlg}^{cn}(\mathbf{C}^{\heart})$ where $(\mathbf{C},\mathbf{C}^{\le0},\mathbf{C}^{\ge0},\mathbf{C}^{0})$ is a derived algebraic context.
We followd \cite{MC} when writing this subsection.
Let $A$ be an augmented object of $\mathrm{sComm}(\mathbf{C}^{\heart})$. 
Let $\mathrm{Sym}(R)\rightarrow A$ be a homotopy epimorphism where $A$ is flat over $R$. Consider the `multiplication by $x$' map
$$\mathrm{Sym}(R)\rightarrow\mathrm{Sym}(R)$$
induced by the canonical map
$$R\rightarrow\mathrm{Sym}(R)$$
Tensoring with $A$ gives a `multiplication by $x$' map
$$A\rightarrow A$$
Consider the Bar complex.
$$B^{A}_{n}=A\otimes A^{\otimes n}\otimes_{R}R$$
where the face maps are defined by the $A$ module structures, and the degeneracy maps are defined by inserting units. There is a map of simplicial commutative monoids.
$$B^{A}\rightarrow R$$

\begin{lem}
\begin{enumerate}
\item
$B^{A}\rightarrow R$ is an equivalence. 
\item
Let $S$ be a simplicial commutative algebra, and $s:R\rightarrow S$ a map. If the induced map $\mathrm{Sym}(R)\rightarrow S$ factors through $\mathrm{Sym}(R)\rightarrow A$, then
$$S\big\slash\big\slash(s)\cong |S\otimes_{A}B^{A}_{\bullet}|$$
\end{enumerate}
\end{lem}

\begin{proof}
\begin{enumerate}
\item
Under the Dold-Kan correspondence the complex becomes the two-term complex
$$A\rightarrow A$$
where the map is given by the multiplication by $x$ map. Now we clearly have that 
$$\mathrm{Sym}(R)\rightarrow\mathrm{Sym}(R)$$
is quasi-isomorphic to $R$. Thus 
$$R\cong A\otimes_{\mathrm{Sym}(R)}^{\mathbb{L}}R$$
is computed by the complex
$$A\rightarrow A$$
where we have used that $\mathrm{Sym}(R)\rightarrow A$ is a homotopy epimorphism in the equivalence above. Applying Dold-Kan again proves the claim.
\item
Since $\mathrm{Sym}(R)\rightarrow A$ is a homotopy epimorphism we have
$$S\hat\otimes^{\mathbb{L}}_{\mathrm{Sym}(R)}R\cong S\otimes^{\mathbb{L}}_{A}R\cong |S\otimes_{A} B^{A}|$$
\end{enumerate}
\end{proof}

\chapter{Lawvere Theories}\label{LT}
Existing approaches to realising analytic and smooth geometry as some form of generalised `algebraic' geometry work by considering affines as algebras relative to certain Lawvere theories (for example, the entire functional calculus in the former case \cite{MR4036665}, and $C^{\infty}$-rings in the latter case \cite{spivak2010derived}). In this section we will give general results establishing when categories of such algebras can be realised, at least in a derived sense, as full subcategories of commutative monoids internal to a monoidal model category. In the context of bornological geometry, these Lawvere theories serve a related purpose. They allow us to define generalised concepts of finite presentation in more general geometries.

\section{Multi-Sorted Lawvere Theories}\label{sec:multi-sort}

Particularly in the setting of (dagger) affinoid geometry, in which we need to take varying radii into account, we will need to work with \textit{multi-sorted} Lawvere theories. Let $\Lambda$ be a set. Consider the overcategories $\mathrm{Set}_{\big\slash\Lambda}$ and $\mathrm{FinSet}_{\big\slash\Lambda}$ where $\mathrm{FinSet}$ is the category of finite sets.  $\mathrm{FinSet}_{\big\slash\Lambda}$ is the free category with finite coproducts generated by $\Lambda$ \cite{nlab:multisorted Lawvere theories}. Then $(\mathrm{FinSet}_{\big\slash\Lambda})^{op}$ is the free category with finite products generated by $\Lambda$.

\begin{defn}
A $\Lambda$-\textit{sorted Lawvere theory} is a category $\mathrm{T}$ with finite products and a finite-product preserving functor
$$k:(\mathrm{FinSet}_{\big\slash\Lambda})^{op}\rightarrow \mathrm{T}$$
which is essentially surjective on objects. When $\Lambda=\{\bullet\}$ we say hat $\mathrm{T}$ is a \textit{one-sorted Lawvere theory}.
\end{defn}

\begin{rem}
    A one-sorted Lawvere theory can be equivalently described as a category with finite products, generated by a single element $x$ under finite products. Let $\mathrm{T}$ be a one-sorted Lawvere theory. We usually write $T(1)$ for a generator, and $T(n)=T(1)^{\times n}$ for $n\in\mathbb{N}_{0}$. 
\end{rem}

Note that, up to passing to the skeleton of $\mathrm{FinSet}$, an object of $\mathrm{FinSet}_{\big\slash\Lambda}$ is fully determined by the data of a map $\phi:\Lambda\rightarrow\mathbb{N}_{0}$ such that $f(\lambda)=0$ for all but finitely many $\lambda\in\Lambda$. If $f:K\rightarrow\Lambda\in \mathrm{FinSet}_{\big\slash\Lambda}$ then the corresponding function $\phi(f)$ is defined by
$$\phi(f)(\lambda)\defeq |f^{-1}(\lambda)|$$
For $\lambda\in\Lambda$, we will denote $k(\{\bullet\}\rightarrow\Lambda,\bullet\mapsto\lambda)$ by $(\lambda)$. Denote by $\underline{n}$ the set $\{1,2,\ldots,n\}$. If $\lambda_{1},\ldots,\lambda_{n}$ is a collection of elements of $\Lambda$, we write $\underline{\lambda}=(\lambda_{1},\ldots,\lambda_{n}):\underline{n}\rightarrow\Lambda$ for the map that sends $i$ to $\lambda_{i}$. 
\newline
\\
If $\Lambda\rightarrow\Lambda'$ is a map of sets, and $\mathrm{T}$ is a $\Lambda'$-sorted Lawvere theory, then by composition
$$(\mathrm{FinSet}_{\big\slash\Lambda})^{op}\rightarrow (\mathrm{FinSet}_{\big\slash\Lambda'})^{op}\rightarrow \mathrm{T}$$
we may regard $\mathrm{T}$ as a $\Lambda$-sorted Lawvere theory. In particular any (one-sorted) Lawvere theory $\mathrm{T}$ may be regarded as a $\Lambda$-sorted Lawvere theory, which we write as $\mathrm{T}^{\Lambda}$. On the other hand, if $\mathrm{T}$ is a $\Lambda$-sorted Lawvere theory, then each element $\lambda\in\Lambda$ gives a one-sorted Lawvere theory $\mathrm{T}^{\lambda}$. 

\begin{defn}
Let $\mathrm{T}$ be a $\Lambda$-sorted Lawvere theory and $\mathrm{C}$ a locally presentable category such that the tensoring over $\mathrm{Set}$ commutes with products in $\mathrm{Set}$. The category $\mathrm{Alg}_{\mathrm{T}}(\mathrm{C})$ of $\mathrm{T}$-\textit{algebras in} $\mathrm{C}$, is the category $\mathrm{Fun}^{\times}(\mathrm{T},\mathrm{C})$ of finite product preserving functors $\mathrm{T}\rightarrow\mathrm{C}$.
\end{defn} 

Note that such functors are completely determined by what they do to elements of the form $\underline{\lambda}=(\lambda_{1},\ldots,\lambda_{n})\in\Lambda^{n}$ for $n\in\mathbb{N}$. 

\begin{example}
Let $\mathpzc{P}$ be an operad in a presentably symmetric monoidal category $\mathrm{M}$ with unit $\mathbb{I}$. Define $\mathrm{T}_{\mathpzc{P}}(n)\defeq\mathpzc{P}(\mathbb{I}^{\oplus n})$ and \[\mathrm{Hom}_{\mathrm{T}}(m,n)\defeq\mathrm{Hom}_{\mathrm{Alg}_{\mathpzc{P}}(\mathrm{M})}(\mathpzc{P}(n),\mathpzc{P}(m)).\] One can do something similar for coloured operads. Then $\mathrm{T}$ will be a multi-sorted theory. If $\mathpzc{P}=\mathrm{Comm}$ is the commutative operad in $\mathrm{Set}$, and $\mathrm{C}$ is a symmetric monoidal $(\infty,1)$-category then $\mathrm{Alg}_{\mathrm{T}_{\mathrm{Comm}}}(\mathrm{C})\cong\mathrm{Comm(C)}$. Coloured operads give examples of multi-sorted Lawvere theories.
\end{example}

\begin{example}
    A particularly important $1$-sorted Lawvere theory for us is the one controlling unital commutative algebras over a unital commutative ring $R$. This is defined by
    $$\mathrm{Hom}_{T}(T(m),T(n))\defeq\mathrm{Hom}_{\mathrm{Comm}({}_{R}\mathrm{Mod})}(\mathrm{Sym}(R^{\oplus n}),\mathrm{Sym}(R^{\oplus m}))$$
    We denote this lawvere theory by $\mathbb{A}_{R}$. An algebra over this Lawvere theory in $\mathrm{Set}$ is precisely a unital commutative $R$-algebra.
\end{example}

\subsection{Basic Properties}

For each object $(f:S\rightarrow\Lambda)$ of $(\mathrm{FinSet}\big\slash_{\Lambda})^{op}$ and $c\in\mathrm{C}$ one has the representable $\mathrm{T}$-algebra
$$\mathrm{Free}_{\mathrm{T},c}(f)\defeq\mathrm{Hom}_{T}(k(f,-)\otimes c$$

When $\mathrm{C}=\mathrm{Set}$, and $c=\{\bullet\}$ we just write 
$$\mathrm{Free}_{\mathrm{T}}(f)\defeq\mathrm{Free}_{\mathrm{T},\{\bullet\}}(f)$$
Let $\mathrm{Alg}_{\mathrm{T}}^{ffg}$ denote the full subcategory of $\mathrm{Alg}_{\mathrm{T}}$ consisting of algebras of the form $\mathrm{Free}_{T}(f)$, where $f\in\mathrm{FinSet}\big\slash_{\Lambda}$.

Denoting by $\mathcal{P}_{\Sigma}(\mathrm{C})$ the ($1-$categorical) free sifted cocompletion of a category, general nonsense using the Yoneda embedding gives  the following. 

\begin{prop}
 There is an equivalence of categories
 $$\mathrm{Alg}_{\mathrm{T}}\cong\mathcal{P}_{\Sigma}(\mathrm{Alg}_{\mathrm{T}}^{ffg})$$
 \end{prop}

Since $\mathrm{Set}_{\big\slash\Lambda}$ is the free completion of the discrete category $\Lambda$, the functor
$$\mathrm{Free}_{T}:\mathrm{FinSet}_{\big\slash\Lambda}\rightarrow\mathrm{Alg}_{T}$$
actually extends to a sifted colimit preserving functor
$$\mathrm{Free}_{T}:\mathrm{Set}_{\big\slash\Lambda}\rightarrow\mathrm{Alg}_{\mathrm{T}}$$
This functor is left adjoint to the forgetful functor
$\mathrm{Alg}_{\mathrm{T}}\rightarrow\mathrm{Fun}^{\times}(\mathrm{FinSet}_{\big\slash\Lambda},\mathrm{Set})\cong\mathrm{Set}_{\big\slash\Lambda}$
We denote by $\mathrm{Alg}_{T}^{f}\subset\mathrm{Alg}_{T}$ the essential image of the functor
$$\mathrm{Free}_{T}:\mathrm{Set}_{\big\slash\Lambda}\rightarrow\mathrm{Alg}_{\mathrm{T}}$$
Note that $\mathrm{Alg}^{f}_{T}$ is generated by $\mathrm{Alg}^{ffg}_{T}$ under \textit{filtered} colimits.

\subsubsection{Abelian Group Objects and Beck Modules}

Let $\mathrm{C}$ be a combinatorial category, and consider the category $\mathrm{Ab}(\mathrm{C})$ of abelian group objects in $\mathrm{C}$. This can equivalently be described as the category
$$\mathrm{Fun}^{\times}(\mathrm{T}_{\mathrm{Ab}},\mathrm{C})$$
of finite product reserving functors $\mathrm{T}_{\mathrm{Ab}}\rightarrow\mathrm{C}$, where $\mathrm{T}_{\mathrm{Ab}}$ is the Lawvere theory defining abelian groups. For $A\in\mathrm{C}$ the category of \textit{Beck modules over }$A$, introduced in \cite{beck1967triples}, is
$${}_{A}\mathrm{Mod^{Beck}}\defeq\mathrm{Ab}(\mathrm{C}_{\big\slash A})$$
 For $\mathrm{T}$ a multisorted Lawvere theory and $A\in\mathrm{Alg_{T}}$ we define
$${}_{A}\mathrm{Mod}\defeq{}_{A}\mathrm{Mod^{Beck}}(\mathrm{Alg_{T}})$$

\subsection{Simplicial $\mathrm{T}$-Algebras}

 Let $\mathrm{T}$ be a $\Lambda$-sorted Lawvere theory. Consider the category $\mathrm{sAlg}_{\mathrm{T}}$ of \textit{simplicial }$\mathrm{T}$-algebras, i.e. product-preserving functors $\mathrm{T}\rightarrow\mathrm{sSet}$. 
 
 \subsubsection{Evaluation Functors}

Let $\mathrm{T}$ be a $\Lambda$-sorted Lawvere theory. For $\lambda\in\Lambda$ we define the \textit{evaluation functor}
$$U_{\lambda}:\mathrm{Alg_{T}}\rightarrow\mathrm{Set}, A\mapsto \mathrm{Hom}_{\mathrm{Alg_{T}}}(\mathrm{Free}_{\mathrm{T}}(\lambda),A)$$
This functor has a left adjoint, which we denote by $\mathrm{F}_{\lambda}$.
 
 \subsubsection{The Model Structure}
 By \cite{MR2263055} Theorem 4.7 there is a model category structure on $\mathrm{sAlg}_{\mathrm{T}}$ in which a map 
$$f:F\rightarrow G$$
is a weak equivalence (resp. a fibration) precisely if $U_{\lambda}(F)\rightarrow U_{\lambda}(G)$ is a weak equivalence (resp. fibration) for every $\lambda\in\Lambda$. This is naturally a simplicial model category. A generating set of cofibrations of $\mathrm{sAlg}_{\mathrm{T}}$ can be defined as follows. 
For any $G\in\mathrm{Alg_{T}}$ define for each $\lambda\in\Lambda$, as in \cite{MR2263055}, the functor
$$U_{\lambda}:\mathrm{sAlg_{T}}\rightarrow\mathrm{sSet},\;\; A\mapsto\mathrm{Map}(\mathrm{Free_{T}}(\alpha),A)$$
has a simplicially enriched left adjoint $F_{\lambda}$. A generating set of cofibrations is given by

$$\{F_{\lambda}(\partial\Delta^{n})\rightarrow F_{\lambda}(\Delta^{n})\}_{n\in\mathbb{N}_{\ge0},\lambda\in\Lambda}$$
and a generating set of acyclic cofibrations is given by 
$$\{\F_{\lambda}(\Lambda_{i})\rightarrow F_{\lambda}(\Delta^{n})\}_{n\in\mathbb{N}_{\ge0},\lambda\in\Lambda,0\le i\le n}$$
The following is essentially the major result of  \cite{MR2263055} (Theorem 1.2), and is a particular case of Lemma \ref{lem:rigidification}.

\begin{thm}
There is an equivalence of $(\infty,1)$-categories
$$\mathrm{L}^{H}(\mathrm{sAlg}_{\mathrm{T}})\cong\mathbf{P}_{\Sigma}(\mathrm{N}(\mathrm{Alg}^{ffg}_{\mathrm{T}}))$$
\end{thm}

In particular this implies that $\mathrm{sAlg}_{\mathrm{T}}$ presents the $(\infty,1)$-category of $(\infty,1)$-algebras over $\mathrm{T}$, i.e. finite product-preserving $(\infty,1)$-functors
$$N(\mathrm{T})\rightarrow\mathbf{sSet}$$

\subsubsection{Simplicial Beck Modules}

Let $\mathrm{T}$ be a $\Lambda$-sorted Lawvere theory and $A$ a simplicial $\mathrm{T}$-algebra. Also by Lemma \ref{lem:rigidification}, there is a model category structure on 
$$\mathrm{Ab}(\mathrm{sAlg_{T\big\slash A}})\cong\mathrm{Fun}^{\times}(\mathrm{T_{Ab}},\mathrm{sAlg_{T}}_{\big\slash A})$$
whereby a map $f:(\pi:B\rightarrow A)\rightarrow (\psi:C\rightarrow A)$ is a fibration (resp. an acyclic fibration) precisely if the underlying map $B\rightarrow C$ of simplicial $\mathrm{T}$-algebras is a fibration (resp. an acyclic fibration). Geometric realisations of simplicial $\mathrm{T}$-algebras commute with products. Thus by \ref{lem:rigidification}, we see that this model structure presents the $(\infty,1)$-category.
$${}_{A}\mathbf{Mod}^{\mathrm{Beck}}\defeq\mathbf{Ab}(\mathbf{sAlg}_{T\big\slash A})$$

\subsection{Filtered Lawvere Theories}
Let us introduce a convenient variant of $\Lambda$-sorted Lawvere theories, which we will need for discussing the theory defining affinoids later.

\begin{defn}
Let $\underline{\Lambda}$ be a small category with set of objects $\Lambda$. A $\underline{\Lambda}$-\textit{sorted Lawvere theory} is a $\Lambda$-sorted Lawvere theory $\mathrm{T}$, together with an extension of the assignment
$$\Lambda\rightarrow\mathrm{T},\lambda\mapsto (\lambda)$$
to a a functor
$$\underline{\Lambda}\rightarrow\mathrm{T}$$
\end{defn}

\begin{defn}
If $T$ is a $\underline{\Lambda}$-sorted Lawvere theory then we define the \textit{ind-associated}  $1$-\textit{sorted Lawvere theory} by 
$$\mathrm{Hom}_{\mathrm{T}^{Ind}}(m,n)\defeq\limpro_{(\underline{\Lambda}^{m})^{op}}\limind_{\underline{\Lambda}^{n}}\mathrm{Hom}_{\mathrm{T}}((\lambda_{1},\ldots,\lambda_{m}),(\gamma_{1},\ldots,\gamma_{n}))$$ 
\end{defn}

There is a left-adjoint functor

$$C:\mathrm{Alg_{T}}\rightarrow\mathrm{Alg_{T^{ind}}}$$
defined on free algebras by 
$$C(\mathrm{Free_{T}}(\lambda_{1},\ldots,\lambda_{m}))(n)\defeq\limind_{\underline{\Lambda}^{n}}\mathrm{Hom_{T}}((\lambda_{1},\ldots,\lambda_{m}),(\gamma_{1},\ldots,\gamma_{n}))$$

\begin{defn}
Let $R$ be a unital commutative ring and $\underline{\Lambda}$ a filtered category. A $\underline{\Lambda}$-\textit{sorted} $R$-\textit{polynomial Lawvere theory} is a $\underline{\Lambda}$-sorted Lawvere theory $\mathrm{T}$ such that
\begin{enumerate}
\item
For each $f:S\rightarrow\Lambda$, $f':S\rightarrow\Lambda$ with $|S|=n$ and $|S|=n'$, 
$$\mathrm{Hom}_{\mathrm{T}}(f,f')\subset\mathrm{Hom}_{\mathrm{Alg}_{R}}(R[y_{1},\ldots,y_{n'}],R[x_{1},\ldots,x_{n}])$$
\item
The induced functor
$$(\mathrm{T})^{ind}\rightarrow\mathbb{A}_{R}$$
is an equivalence.
\item
For any $\phi\in\mathrm{Hom}_{\mathrm{Alg}_{R}}(R[y_{1},\ldots,y_{n'}],R[x_{1},\ldots,x_{n}])$ there is some $f:S\rightarrow\Lambda$, $f':S\rightarrow\Lambda$ with $|S|=n$ and $|S|=n'$ such that $f\in\mathrm{Hom}_{\mathbb{A}_{R}^{\Lambda}}(f,f')$. 
\end{enumerate}
\end{defn}

We will be particularly interested in the case that $\underline{\Lambda}$ has the structure of a partially ordered commutative monoid. In this case we have in some sense a canonical $\underline{\Lambda}$-sorted $R$-polynomial Lawvere theory for any unital commutative ring $R$. Throughout we will ue multiplicative notation for mnoids.

\begin{example}
Let $\underline{\Lambda}$ be a directed commutative monoid, which we regard as a symmetric monoidal filtered category. We consider the abelian category 
$$\mathrm{Fun}(\underline{\Lambda},{}_{R}\mathrm{Mod})$$
of functors from $\underline{\Lambda}$ to ${}_{R}\mathrm{Mod}$. This is a closed symmetric monoidal abelian category with respect to Day convolution. Moreover it is in fact a monoidal elementary abelian category. A generating set of projectives is given by $\{F_{i}R\}_{i\in\underline{\Lambda}}$ where $(F_{i}R)(\lambda)=R$ if $\lambda\ge i$ and $F_{i}(R)(\lambda)=0$ otherwise. For $\underline{\lambda}=(\lambda_{1},\ldots,\lambda_{n})\in\underline{\Lambda}^{n}$ write 
$$\mathrm{Sym}(F_{\underline{\lambda}}R)\defeq\mathrm{Sym}(F_{\lambda_{1}}R\oplus\ldots\oplus F_{\lambda_{m}}R)$$
 Define a $\underline{\Lambda}^{op}$-filtered Lawvere theory by
$$\mathrm{Hom}(\underline{\lambda},\underline{\gamma})\defeq\mathrm{Hom}_{\mathrm{Comm}(\mathrm{Fun}(\underline{\Lambda},{}_{R}\mathrm{Mod}))}(\mathrm{Sym}(F_{\underline{\gamma}}R),\mathrm{Sym}(F_{\underline{\lambda}}R))$$
This is a $\underline{\Lambda}^{op}$-sorted Lawvere theory which we denote by $\mathbb{A}_{R}^{\underline{\Lambda}}$. Since $\underline{\Lambda}$ is directed we have 
$$(\mathbb{A}_{R}^{\underline{\Lambda}})^{ind}\cong\mathbb{A}_{R}$$
More generally for a $\underline{\Lambda}$-filtered $R$-algebra $S$ we define a Lawvere theory $\mathbb{A}^{\underline{\Lambda}}_{S}$ by 
$$\mathrm{Hom}(\underline{\lambda},\underline{\gamma})\defeq\mathrm{Hom}_{\mathrm{Comm}({}_{S}\mathrm{Mod}(\mathrm{Fun}(\underline{\Lambda},{}_{R}\mathrm{Mod})))}(S\otimes\mathrm{Sym}(F_{\underline{\gamma}}R),S\otimes\mathrm{Sym}(F_{\underline{\lambda}}R))$$
Note that since both categories are the free cocompletion of the category $\mathbb{A}_{R}^{\underline{\Lambda}}$, we have an equivalence.
$$\mathrm{Alg}_{\mathbb{A}_{S}^{\underline{\Lambda}}}\cong{}_{S}\mathrm{Mod}(\mathrm{Fun}(\underline{\Lambda},{}_{R}\mathrm{Mod})))$$
Moreover we have 
$$(\mathbb{A}_{S}^{\underline{\Lambda}})^{ind}\cong\mathbb{A}_{\limind_{\underline{\Lambda}}S_{\lambda}}$$
\end{example}

\subsubsection{Open Filtrations}

\begin{defn}
A $\underline{\Lambda}$-sorted Lawvere theory $\mathrm{T}$ is said to \textit{have an open filtration} if for each $\lambda>\gamma$ the map 
$$\mathrm{Free}_{\mathrm{T}}(\lambda)\rightarrow\mathrm{Free}_{\mathrm{T}}(\gamma)$$
is a homotopy epimorphism in $\mathbf{sAlg}_{\mathrm{T}}$. 
\end{defn}

Let $A\rightarrow B$ be a homotopy epimorphism in $\mathbf{sAlg}_{\mathrm{T}}$. We think of $\mathbf{sAlg}_{\mathrm{T}}^{op}$ as a category of affines. Writing $\mathrm{Spec}(A)$ for the object $A$ considered in the opposite category , we have that $\mathrm{Spec}(B)\rightarrow\mathrm{Spec}(A)$ is a homotopy monomorphism. We think of such maps as open embeddings. We will make this precise in the context of bornological algebras later.

Let $\mathrm{T}$ be a $\underline{\Lambda}$-sorted Lawvere and let $A\in\mathbf{sAlg}_{\mathrm{T}}$. Define 
$$|A|_{\lambda}\defeq\mathbf{Map}_{\mathbf{sAlg}_{\mathrm{T}}}(\mathrm{Free_{T}}(\lambda),A)$$
and
$$|A|\defeq\limind_{\lambda} |A|_{\lambda}$$
if $\mathrm{T}$ has an open filtration then the map 
$$|A|_{\lambda}\rightarrow |A|$$
is a homotopy monomorphism


\section{Fermat Theories}\label{subsec:fermat}

In this section we introduce filtered Fermat theories, following the original definition of one-sorted Fermat theories from \cite{MR0744457}. We must restrict to the case of theories sorted by partially ordered \textit{abelian groups}. We will use multiplicative notation for such groups. In particular for two elements $a,b$ we will write $\frac{a}{b}\defeq ab^{-1}$.

\subsection{Definitions}

\begin{defn}
Let $\Gamma$ be a partially ordered abelian group, $R$ a unital commutative ring, $\mathrm{T}$ a $\Gamma$-sorted Lawvere theory, and $\mathbb{A}^{\mathrm{T}}$ a $\Gamma$-sorted $R$-polynomial Lawvere theory. A map of $\Gamma$-sorted Lawvere theories
$$\mathbb{A}^{\mathrm{T}}\rightarrow\mathrm{T}$$
 is said to \textit{define a} $\Gamma$-\textit{sorted} \textit{Fermat theory} if for any $\delta$, any $(\lambda,\gamma_{1},\ldots,\gamma_{n})$, and any $f(x,z_{1},\ldots,z_{n})\in \mathrm{Free_{T}}(\lambda,\gamma_{1},\ldots,\gamma_{n})(\delta)$, there is a unique $g(x,y,z_{1},\ldots,z_{n})\in \mathrm{Free_{T}}(\lambda,\lambda,\gamma_{1},\ldots,\gamma_{n})(\frac{\delta}{\lambda})$
such that
$$f(x,z_{1},\ldots,z_{n})-f(y,z_{1},\ldots,z_{n})=(x-y)g(x,y,z_{1},\ldots,z_{n})$$
If $S$ is a $\Gamma$-filtered ring a \textit{strict }$\Gamma$-\textit{sorted }$S$-\textit{Fermat theory} is a $\Gamma$-sorted Fermat theory of the form
$$\mathbb{A}^{\underline{\Gamma}}_{S}\rightarrow\mathrm{T}$$
If such a map is implicit, we will say that $\mathrm{T}$ is a strict $\Gamma$-sorted $S$ Fermat theory.
\end{defn}

Consider the diagonal map $\lambda\rightarrow\lambda\times\lambda$ which induces a map 
$$\mathrm{Free_{T}}(\lambda,\lambda,\underline{\lambda})\rightarrow\mathrm{Free_{T}}(\underline{\lambda})$$
We denote by $\partial_{1}f$ the image of $g(x,y,z_{1},\ldots,z_{n})$ under this map. By the inherent symmetry we can define $\partial_{i}f$ for $1\le i\le n+1$.

\begin{example}
$\mathbb{A}_{S}^{\underline{\Gamma}}$ is a strict $\Gamma$-sorted Fermat theory. 
\end{example}

\begin{lem}
Let $\Gamma$ be a partially ordered abelian group, $R$ a unital commutative rirg, $\mathrm{T}$ an $\Gamma$-sorted Lawvere theory, and $\mathbb{A}^{\mathrm{T}}$ a $\Gamma$-sorted $R$-polynomial Lawvere theory. If
$$\mathbb{A}^{\mathrm{T}}\rightarrow\mathrm{T}$$
is a $\Gamma$-sorted Fermat theory then $\mathrm{T}^{ind}$ is a Fermat theory.
\end{lem}

\begin{proof}
Let $f(x,z_{1},\ldots,z_{n})\in\mathrm{T}^{ind}(n+1)$, and let $\mathrm{F}(x,z_{1},\ldots,z_{n})\in\mathrm{T}(\lambda_{1}\ldots\lambda_{n})$ be a representative for $f(x,z_{1},\ldots,z_{n})$ in the colimit. There is a $\overline{g}(x,y,z_{1},\ldots,z_{n})\in T(\lambda,\lambda,\gamma_{1},\ldots,\gamma_{n})(\frac{\delta}{\lambda})$
such that
$$\mathrm{F}(x,z_{1},\ldots,z_{n})-\mathrm{F}(y,z_{1},\ldots,z_{n})=(x-y)\overline{g}(x,y,z_{1},\ldots,z_{n})$$
Let $g(x,y,z_{1},\ldots,z_{n})\in\mathrm{T}^{ind}(n+2)$ be the equivalence class of $\overline{g}$. Then
$$f(x,z_{1},\ldots,z_{n})-f(y,z_{1},\ldots,z_{n})=(x-y)g(x,y,z_{1},\ldots,z_{n})$$
\end{proof}

Note that whenever $\mathrm{T}$ is a strict $\Gamma$-filtered Lawvere theory
there is a forgetful functor
$$|-|_{\Gamma-alg}:\mathrm{Alg_{T}}\rightarrow\mathrm{Comm}({}_{\mathrm{S}}\mathrm{Mod}(\mathrm{Filt}_{\Gamma}(\mathrm{Ab})))$$
to filtered $S$-algebras. If $A$ is a $\mathrm{T}$-algebra we denote the underlying filtered ring by $|A|_{\Gamma-alg}$.

\subsection{Beck Modules and Derivations}

We know define Fermat-theoretic derivations and the cotangent complex, following closely the $1$-sorted exposition in \cite{MR4036665}. Fix a strict $\Gamma$-sorted $S$-Fermat theory $\mathrm{T}$.

\subsubsection{Comparison with Beck Modules}

Let $\mathrm{T}$ be a strict $\Gamma$-sorted $S$-Fermat theory. Let $A\in\mathrm{Alg_{T}}(\mathrm{Set})$ and $M$ a $\Gamma$-filtered $|A|_{\Gamma-alg}$-module in $\mathrm{Set}$. As in \cite{MR4036665}, define a $\mathrm{T}$-algebra structure on the $S$-algebra $A\oplus M$ by
$$\Phi_{f}^{A\oplus M}(a_{1}+m_{1},\ldots,a_{n}+m_{n})\defeq\Phi_{f}^{A}(a_{1},\ldots,a_{n})+\sum_{i=1}^{n}\Phi_{\partial_{i}f}(a_{1},\ldots,a_{n})m_{i}$$
This is an abelian group object of $(\mathrm{Alg}_{\mathrm{T}})_{\big\slash A}$ by defining
$$(a,m)(a,m')\defeq(a,m+m')$$
$$0\defeq(a,0)$$
$$-(a,m)\defeq(a,-m)$$

This extends to a functor.
$$A\oplus(-):{}_{|A|_{\Gamma-alg}}\mathrm{Mod}\rightarrow(\mathrm{Ab}(\mathrm{Alg}_{\mathrm{T}_\big\slash A})$$
$$M\mapsto A\oplus M$$

%

\begin{lem}\label{lem:1beck}
Let $\mathrm{T}$ be a strict $\Gamma$-sorted $S$-Fermat theory. The functor
$$\mathrm{sqz}_{A}\defeq A\oplus(-):{}_{|A|_{\Gamma-alg}}\mathrm{Mod}\rightarrow{}_{A}\mathrm{Mod^{Beck}}$$
is an equivalence of categories.
\end{lem}

\begin{proof}
Let $(\pi:B\rightarrow A)\in{}_{A}\mathrm{Mod^{Beck}}$. Define $\rho(\pi)\defeq\mathrm{Ker}(\pi)$. This is a well-defined functor and $\rho\circ\mathrm{sqz}_{A}$ is isomorphic to the identity. Now let $(\pi:B\rightarrow A)\in{}_{A}\mathrm{Mod^{Beck}}$. The $0$ map for the abelian group structure $(A\rightarrow A)\rightarrow (B\rightarrow A)$ gives a section of $q$. Thus, as filtered $S$-modules, we may write $B\cong A\oplus\mathrm{Ker}(\pi)$. A straightforward computation as in e.g. \cite{nlab:beck_module} shows that this is in fact an isomorphism in ${}_{A}\mathrm{Mod^{Beck}}$. 
\end{proof}

Let $S$ be a filtered ring and $T$ a strict $\Gamma$-sorted $S$-Fermat theory. The functor 
$$|-|_{\Gamma-alg}:\mathrm{sAlg}_{\mathrm{T}}\rightarrow\mathrm{sComm}({}_{\mathrm{S}}\mathrm{Mod}(\mathrm{Filt}(\mathrm{Ab})))$$
is right Quillen. Now as in Lemma \ref{lem:1beck} we also get an equivalence of categories of simplicial objects:

$$\adj{\mathrm{sqz}_{A}}{{}_{|A|_{\Gamma-alg}}\mathrm{sMod}}{{}_{A}\mathrm{sMod^{Beck}}}{\rho}$$

${}_{A}\mathrm{sMod^{Beck}}$ may be equipped with a model structure whereby a map $f:(\pi:B\rightarrow A)\rightarrow (\psi:C\rightarrow A)$ is a fibration (resp. an acyclic fibration) precisely if the map $B\rightarrow C$ is a fibration (resp. an acyclic fibration) of objects of $\mathrm{sAlg_{T}}$. As explained in the proof of Lemma \ref{lem:1beck}, $\pi:B\rightarrow A$ and $\psi:C\rightarrow A$ have sections, so they are fibrations. Thus the map $\mathrm{Ker}(\pi)\rightarrow\mathrm{Ker}(\psi)$ is a fibration (resp. an acyclic fibration) of simplicial objects whenever $f$ is. Hence $\rho$ is right Quillen.  In fact since all objects are fibrant, $\rho$ preserves weak equivalences between all objects. Now the functor $\mathrm{sqz}_{A}$ also preserves all equivalences. We have proven the following.

\begin{lem}\label{lem:modelbeck}
Let $S$ be a filtered unital commutative ring $\mathrm{T}$ be a strict $\Gamma$-sorted $S$-Fermat theory, and $A$ a simplicial $\mathrm{T}$-algebra. The adjunction
$$\adj{\mathrm{sqz}_{A}}{{}_{|A|_{\Gamma-alg}}\mathrm{sMod}}{{}_{A}\mathrm{sMod^{Beck}}}{\rho}$$
is a Quillen equivalence. In particular there is an equivalence of $(\infty,1)$-categories
$${}_{|A|_{\Gamma-alg}}\mathbf{Mod}\cong{}_{A}\mathbf{Mod}^{\mathrm{Beck}}$$
\end{lem}

 \subsubsection{Derivations}

\begin{defn}
Let $A$ be a $\mathrm{T}$-algebra and $M$ a $\Gamma$-filtered $A_{\Gamma-alg}$-module. A $\mathrm{T}$-\textit{derivation} is a map $\delta:A\rightarrow M$  such that 
$$A\rightarrow A\oplus M,\;\; a\mapsto (a,\delta a)$$
is a homomorphism of $\mathrm{T}$-algebras.
\end{defn}

Let $\mathbb{T}_{\mathrm{Alg_{T}}}$ denote the category of pairs $(A,M)$ where $A\in\mathrm{Alg_{T}}$ and $M\in{}_{A}\Mod$. There is a functor
$$Z:\mathbb{T}_{\mathrm{Alg_{T}}}\rightarrow\mathrm{Alg_{T}},\;\; (A,M)\mapsto A\oplus M$$
This functor has a left adjoint $A\mapsto (A,\Omega^{1}_{A})$. 

\begin{example}
Consider a sequence $(\lambda_{1},\ldots,\lambda_{n})\in\Lambda^{n}$ and the free $\mathrm{T}$-algebra $\mathrm{Free_{T}}(\lambda_{1},\ldots,\lambda_{n})$. For each $\lambda_{i}$ let $S_{\lambda_{i}}$ denote the filtered $S$-module generated freely by an element $dx_{i}$ of degree $\lambda_{i}$. Then
$$\Omega^{1}_{\mathrm{Free_{T}}(\lambda_{1},\ldots,\lambda_{n})}\cong \bigoplus_{i=1}^{n}\mathrm{Free_{T}}(\lambda_{1},\ldots,\lambda_{n})\otimes S_{\lambda_{i}}$$
The map 
$$\mathrm{Free_{T}}(\lambda_{1},\ldots,\lambda_{n})\rightarrow\Omega^{1}_{\mathrm{Free_{T}}(\lambda_{1},\ldots,\lambda_{n})}$$
is given by
$$f\mapsto\sum_{i=1}^{n}\frac{\partial f}{\partial x_{i}}dx_{i}$$
\end{example}

\subsection{Fermat Theories and Chain Complexes}

Let $\mathrm{T}$ be a strict $\Gamma$-sorted $S$-Fermat theory. Similarly to \cite{MR4036665} Section 2.2 we define differential-graded $T$-algebras. Throughout this section we will assume that everything is defined over $\mathbb{Q}$.

\begin{defn}
A $\mathrm{T}$-\textit{differential graded algebra} is an object $A$ of $\mathrm{Comm}(Ch_{\ge0}({}_{S}\mathrm{Mod}))$ such that $A_{0}$ is equipped with a $\mathrm{T}$-algebra structure extending its $S$-algebra structure. 
\end{defn}

There is an obvious notion of morphism of $\mathrm{T}$-differential graded algebra, so that $\mathrm{T}$-differential graded algebras arrange into a category $\mathrm{Alg_{T}}(Ch_{\ge0}({}_{S}\mathrm{Mod}))$. The category  $\mathrm{Alg_{T}}(Ch_{\ge0}({}_{S}\mathrm{Mod}))$ has all small limits and colimits. The forgetful functor
$$|-|\mathrm{Alg_{T}}:(Ch_{\ge0}({}_{S}\mathrm{Mod}))\rightarrow Ch_{\ge0}({}_{S}\mathrm{Mod})$$
commutes with limits and filtered colimits. In particular this implies that $|-|$ has a left adjoint, which we denote by $\mathrm{Sym_{T}}$.

Following \cite{MR4036665}, pushouts may be computed as follows. Let $B\rightarrow A$ and $B\rightarrow C$ be maps in $\mathrm{Alg_{T}}(Ch_{\ge0}({}_{S}\mathrm{Mod}))$. Let $\mathrm{Free_{\mathrm{T}}}(\lambda_{1},\ldots,\lambda_{m})\rightarrow A_{0}$ and $\mathrm{Free_{\mathrm{T}}}(\gamma_{1},\ldots,\gamma_{n})\rightarrow B_{0}$ be epimorphisms. Then we have
$$A\coprod_{B}C\cong\mathrm{Free_{\mathrm{T}}}(\lambda_{1},\ldots,\lambda_{m},\gamma_{1},\ldots,\gamma_{n})\otimes_{\mathrm{T}(\lambda_{1},\ldots,\lambda_{m})\otimes\mathrm{T}(\gamma_{1},\ldots,\gamma_{n})}A\otimes_{B}C$$

%

\subsubsection{The de Rham Algebra}
Let $\mathrm{T}$ be a $\Gamma$-sorted $S$-Fermat theory. 

\begin{defn}
Let $\underline{\lambda}=(\lambda_{1},\ldots,\lambda_{n})\in\Gamma^{n}$. Just as in the $1$-sorted case in \cite{carchedi2012homological} we define the de Rham algebra $\Omega^{\bullet}_{\underline{\lambda}}$ to be the differential graded $\mathrm{T}$-algebra constructed as follows.
$$\Omega^{0}_{\underline{\lambda}}\defeq\mathrm{Free_{\mathrm{T}}}(\underline{\lambda})$$
$$\Omega^{n}_{\underline{\lambda}}\defeq\bigwedge^{n}\Omega_{T(\underline{\lambda})}$$
This has an evident differential-graded $\mathrm{T}$-algebra structure. 
\end{defn}

For $\underline{a}=(a_{1},\ldots,a_{n})\in\prod_{i=1}^{n}\mathrm{Free_{\mathrm{T}}}(\emptyset)_{\lambda_{i}}$ we define the evaluation map 
$$p_{\underline{a}}:\Omega^{\bullet}_{\underline{\lambda}}\rightarrow\mathrm{Free_{\mathrm{T}}}(\emptyset)_{\lambda_{i}}$$
to be the composite
\begin{displaymath}
\xymatrix{
\Omega^{\bullet}_{\underline{\lambda}}\ar[r] &\mathrm{Free_{\mathrm{T}}}(\underline{\lambda})\ar[r]^{ev_{\underline{a}}} & \mathrm{Free_{\mathrm{T}}}(\emptyset)
}
\end{displaymath}
where the first map is the projection and the second the evaluation map. In particular we get maps
$$p_{\overline{a}}:A\coprod\Omega^{\bullet}_{\underline{\lambda}}\rightarrow A$$

\begin{defn}
$\mathrm{T}$ is said to \textit{satisfy the homotopy invariance property} if the map
$$A\rightarrow A\coprod\Omega_{(1)}^{\bullet}$$
is a quasi-isomorphism, where $1$ is the unit of the monoid $\Gamma$.
\end{defn}

Let
$$p:A\coprod\Omega^{n}_{(1)}\rightarrow A\times A$$
be given by $(p_{0},p_{1})$, where $0$ and $1$ are the additive and multiplicative identities of the ring $\mathrm{Free_{T}}(1)$ respectively. Let $t\in F_{1}\mathrm{Free_{T}}(1)$ be the element induced by the canonical map  $(1)\rightarrow\mathrm{Free_{T}}(1)$. Consider also the element $1-t\in F_{1}(A\coprod\mathrm{Free_{T}}(1))$. These give, via multiplication a map
$$i:(\times(1-t),\times t):A\times A\rightarrow A\coprod\Omega^{n}_{1}$$
which satisfies $q\circ i=\mathrm{Id}$. Thus if $\mathrm{T}$ satisfies the homotopy invariance property, we get a factorisation
$$A\rightarrow A\coprod\Omega^{\bullet}_{(1)}\rightarrow A\times A$$
with the first map being a quasi-isomorphism and the second map being a fibration of the underlying complexes. Exactly as in \cite{carchedi2012homological} Section 6, this proves the following.

\begin{thm}
The transferred model structure exists along the adjunction
$$\adj{\mathrm{Sym_{T}}}{\mathrm{Ch}_{\ge0}({}_{S}\mathrm{Mod})}{\mathrm{Alg_{T}}(\mathrm{Ch}({}_{S}\mathrm{Mod}))}{|-|}$$
\end{thm}

Moreover, exactly as in \cite{nuiten2018lie} Corollary 2.2.10 we get the following.

\begin{thm}
If $\mathrm{T}$ satisfies the homotopy invariant property then the Dold-Kan adjunction induces a Quillen equivalence between $\mathrm{Alg_{T}}(\mathrm{Ch}({}_{S}\mathrm{Mod}))$ and $\mathrm{sAlg_{T}}$.
\end{thm}

\subsubsection{Integration}

One situation in which the homotopy invariance property is satisfied is when $\mathrm{T}$ admits integration.

\begin{defn}
Let $\mathrm{T}$ be a $\Gamma$-sorted $S$-Fermat theory. $\mathrm{T}$ is said to \textit{admit integration} if 
\begin{enumerate}
\item
for any $(\gamma_{1},\ldots,\gamma_{n})$, any, and any $\delta$, the map
$$\partial_{1}:T(1,\gamma_{1},\ldots,\gamma_{n})(\delta)\rightarrow T(1,\gamma_{1},\ldots,\gamma_{n})(\delta)$$
is a surjection.
\item
If $F,G\in T(1,\gamma_{1},\ldots,\gamma_{n})(\frac{\delta}{\lambda})$ satisfy $\partial_{1}F=\partial_{1}G$ then $F-G\in T(\gamma_{1},\ldots,\gamma_{n})(\delta)$
\end{enumerate}
\end{defn}

The following can be proven exactly as in \cite{carchedi2012homological} Corollary 5.20.

\begin{lem}
If $\mathrm{T}$ admits integration then it satisfies the homotopy invariance property.
\end{lem}

\subsubsection{Presentations and Coherent Algebras Over Fermat Theories}

Let $\mathrm{T}$ be a strict $\Gamma$-sorted $S$-Fermat theory.

\begin{lem}[c.f. \cite{nuiten2018lie} Lemma 2.2.7]\label{lem:algpres}
 Let $A\rightarrow B$ be a map in $\mathrm{Alg_{T}}(\mathrm{Ch_{\ge0}}({}_{S}\mathrm{Mod})$ Then there is a sequence
$$A=A^{-1}_{\bullet}\rightarrow A^{0}_{\bullet}\rightarrow A_{\bullet}^{1}\rightarrow\cdots\rightarrow A_{\bullet}^{n}\rightarrow\cdots\rightarrow B$$
such that
\begin{enumerate}
\item 
each $A^{(n)}\rightarrow B$ induces an isomorphism on homology in degree $<n$ and a surjection on homology in degree $n$. 
\item 
$A^{(n+1)}$ is a pushout of $A^{(n)}$ along a map of the form $\mathrm{Free_{T}}(\bigoplus_{\alpha}S^{n}(S_{\lambda_{\alpha}}))\rightarrow\mathrm{Free_{T}}(\bigoplus_{\alpha}D^{n+1}(S_{\lambda_{\alpha}}))$
\end{enumerate}
If $A=\mathrm{Free_{T}}(\emptyset)$ is the initial algebra and $B$ sis such that $\pi_{0}(B)$ is coherent as an algebra, and each $\pi_{n}(A)$ is finitely generated, then at each stage we can in fact take a sum over finitely many $\alpha$.
\end{lem}

Recall that $S_{\lambda}$ denotes $F_{\lambda}\mathbb{Z}\otimes S$.

\begin{proof}
     As in \cite{nuiten2018lie} Lemma 2.2.7. we may construct a factorisation
$$A\rightarrow\tilde{B}\rightarrow B$$
with $A\rightarrow\tilde{B}$ being a cofibration, and $\tilde{B}\rightarrow B$ being an equivalence.  In fact by transfinite induction we will construct a sequence of cofibrations
$$A=A^{(-1)}\rightarrow A^{(0)}\rightarrow\ldots\rightarrow B$$
where each map $A^{(n)}\rightarrow B$ induces an isomorphism on homology in degree $<n$ and a surjection on homology in degree $n$. Given $A^{(n)}\rightarrow B$ define $A^{(n+1)}$ to be the pushout
\begin{displaymath}
\xymatrix{
\mathrm{Free_{T}}(\bigoplus_{\alpha}S^{n}(S_{\lambda_{\alpha}}))\ar[r]\ar[d] & A^{(n)}\ar[d]\\
\mathrm{Free_{T}}(\bigoplus_{\alpha}D^{n-1}(S_{\lambda_{\alpha}}))\ar[r] & A^{(n+1)}\ar[r] & B
}
\end{displaymath}
where the sum is over all such possible commutative diagrams
\begin{displaymath}
\xymatrix{
\mathrm{Free_{T}}(S^{n}(S_{\lambda_{\alpha}}))\ar[r]\ar[d] & A^{(n)}\ar[d]\\
\mathrm{Free_{T}}(D^{n-1}(S_{\lambda_{\alpha}}))\ar[r]  & B
}
\end{displaymath}

\end{proof}

\begin{rem}
We suspect that this works simplicially as well, though we haven't worked out the details.
\end{rem}

\subsection{The Cotangent Complex}

Following the $1$-sorted exposition in \cite{MR4036665} Section 3.2 we construct the cotangent complex of  a map of $\mathrm{T}$-algebras.
Let $\mathrm{T}$ be a strict $\Gamma$-sorted $S$-Fermat theory.. The category $\mathrm{s}\mathbb{T}_{\mathrm{Alg_{T}}}$ also has a cofibrantly generated model structure whereby a map $(A,M)\rightarrow (A',M')$ is a fibration (resp. weak equivalence) precisely if $A\rightarrow A'$ is a fibration (resp. weak equivalence) in $\mathrm{sAlg_{T}}$, and $M\rightarrow M'$ is a fibration (resp. weak equivalence) of simplicial abelian groups. The following is then clear.

\begin{lem}
The adjunction
$$\adj{\Omega^{1}}{\mathrm{sAlg_{T}}}{\mathbb{T}_{\mathrm{sAlg_{T}}}}{Z}$$
 is Quillen.
\end{lem}

We denote by $A\mapsto (A,\mathbb{L}_{A})$ the left derived functor of $\Omega^{1}$. 

\begin{rem}
If $\mathrm{T}$ satisfies the homotopy invariance property then we can consider $\mathrm{Alg_{T}}(\mathrm{Ch}({}_{S}\mathrm{Mod}))$ rather than $\mathrm{sAlg_{T}}$. We then define $\mathbb{T}_{\mathrm{Alg_{T}}(\mathrm{Ch}({}_{S}\mathrm{Mod}))}$ to be the category of pairs $(A,M)$ where $A\in\mathrm{Alg_{T}}(\mathrm{Ch}({}_{S}\mathrm{Mod}))$ and $M$ is an $A_{alg}$-module. This is Quillen equivalent to $\mathrm{s}\mathbb{T}_{\mathrm{Alg_{T}}}$. Moreover, again get a Quillen adjunction
$$\adj{\Omega^{1}}{\mathrm{Alg_{T}}(\mathrm{Ch}({}_{S}\mathrm{Mod}))}{\mathbb{T}_{\mathrm{Alg_{T}}}(\mathrm{Ch}({}_{S}\mathrm{Mod}))}{Z}$$
and we define $A\mapsto (A,\mathbb{L}_{A})$ to be the left derived functor of $\Omega_{1}$.
At the level of $(\infty,1)$-categories both the simplicial and chain complex construction are equivalent. 
\end{rem}





\subsubsection{Comparison With the Algebraic Cotangent Complex}

Let $\mathrm{U}\rightarrow\mathrm{T}$ be a map of strict $\Gamma$-sorted $S$-Fermat theories, and let $A\rightarrow B$ be a map in $\mathrm{sAlg_{T}}$. Write $A_{\mathrm{U}}$ for $A$ regarded as a $\mathrm{U}$-algebra, and similarly for $B_{\mathrm{U}}$. There is a natural map in ${}_{A_{U}}\mathbf{Mod}\cong{}_{A}\mathbf{Mod}$
$$\mathbb{L}_{B_{U}\big\slash A_{U}}\rightarrow\mathbb{L}_{B\big\slash A}$$

In particular for $\mathrm{U}=\mathbb{A}^{\Gamma}_{S}$ we write
$$A_{alg}\defeq A_{\mathbb{A}^{\Gamma}_{S}}$$
We are now going to establish circumstances under which the map
$$\mathbb{L}_{B_{alg}\big\slash A_{alg}}\rightarrow\mathbb{L}_{B\big\slash A}$$
is an equivalence.

Suppose that $\mathrm{T}$ is a $\Gamma$-sorted $S$-Fermat theory satisfying the homotopy invariance property.

\begin{lem}[\cite{MR4036665} Lemma 3.20]
Suppose that $\mathrm{T}$ is a $\Gamma$-sorted $S$-Fermat theory satisfying the homotopy invariance property. and $A\rightarrow B$ is a map in $\mathrm{Alg_{T}}(\mathrm{Ch}_{\ge0}({}_{S}\mathrm{Mod}))$ such that $H_{0}(A)\rightarrow H_{0}(B)$ is an epimorphism, then the natural map
$$\mathbb{L}_{B_{alg}\big\slash A_{alg}}\rightarrow\mathbb{L}^{alg}_{B\big\slash A}$$
is an equivalence.
\end{lem}

\begin{proof}
We follow a similar proof to \cite{nuiten2018lie} Lemma 2.2.7.
Let $A\rightarrow B$ and $A\rightarrow C$ be maps between cofibrant objects in $\mathrm{sAlg_{T}}$, and suppose that $A\rightarrow B$ induces a surjection on $H_{0}$.
In the situation of Lemma \ref{lem:algpres}, when $H_{0}(A)\rightarrow H_{0}(B)$ is surjective we can take $A^{(0)}=A$.
We then clearly have
$$\mathbb{L}_{A^{(0)}_{alg}\big\slash A_{alg}}\cong\mathbb{L}^{alg}_{A^{(0)}\big\slash A}$$
Suppose we have shown that 
$$\mathbb{L}_{A_{alg}^{(m+1)}\big\slash A^{(m)}_{alg}}\rightarrow \mathbb{L}^{alg}_{A^{(m+1)}\big\slash A^{(m)}}$$
for some $n\ge 0$ and all $m<n$.

 Now we have 
$$\mathbb{L}_{\mathrm{Free_{T}}(\bigoplus_{\alpha}S^{n}(S))}\cong \mathrm{Free_{T}}(\bigoplus_{\alpha}S^{n}(S)) \otimes^{\mathbb{L}}S^{n}(S)\cong\mathrm{Sym}(\bigoplus_{\alpha}S^{n}(S))\otimes^{\mathbb{L}}S^{n}(S)$$
since $n<0$. Similarly we have 
$$\mathbb{L}_{\mathrm{Free_{T}}(\bigoplus_{\alpha}S^{n}(S))}\cong \mathrm{Free_{T}}(\bigoplus_{\alpha}D^{n-1}(S)) \otimes^{\mathbb{L}}D^{n-1}(S)\cong\mathrm{Sym}(\bigoplus_{\alpha}D^{n-1}(S))\otimes^{\mathbb{L}}D^{n-1}(S)$$
This clearly implies that the map

$$\mathbb{L}_{\mathrm{Free_{T}}(\bigoplus_{\alpha}D^{n-1}(S))_{alg}\big\slash\mathrm{Free_{T}}(\bigoplus_{\alpha}S^{n}(S)_{alg})}\rightarrow\mathbb{L}^{alg}_{\mathrm{Free_{T}}(\bigoplus_{\alpha}D^{n-1}(S))\big\slash \mathrm{Free_{T}}(\bigoplus_{\alpha}S^{n}(S))}$$
is an equivalence and therefore using the pushout diagram, that

$$\mathbb{L}_{A_{alg}^{(n+1)}\big\slash A^{(n)}_{alg}}\rightarrow \mathbb{L}^{alg}_{A^{(n+1)}\big\slash A^{(n)}}$$
is an equivalence. An induction using the fibre-cofibre sequence
$$\mathbb{L}_{A^{(n)}\big\slash A}\otimes^{\mathbb{L}}_{A^{(n)}}A^{(n+1)}\rightarrow\mathbb{L}_{A^{(n+1)}\big\slash A}\rightarrow\mathbb{L}_{A^{(n+1)}\big\slash A^{(n)}}$$ then implies that 
$$\mathbb{L}_{A_{alg}^{(n)}\big\slash A_{alg}}\rightarrow \mathbb{L}^{alg}_{A^{(n)}\big\slash A}$$
is an equivalence for each $n$. Taking the (homotopy) colimit, this suffices to establish the result.
\end{proof}





%

\section{Lawvere Theories and Homotopical Commutative Algebra}\label{sec:LAWVFunc}

Let $(\mathrm{C},\otimes,\mathbb{I},\underline{\mathrm{Hom}})$ be a closed symmetric monoidal category which is complete and cocomplete. Consider the Lawvere theory $\mathbb{A}_{\mathrm{C}}$ defined by
$$\mathrm{Hom}_{\mathbb{A}_{\mathrm{C}}}(m,n)\defeq\mathrm{Hom}_{\mathrm{Comm(\mathrm{C})}}(\mathrm{Sym}(\mathrm{Sym}(\mathbb{I}^{\oplus n}),\mathrm{Sym}(\mathbb{I}^{\oplus n})$$

Let $\underline{\Lambda}$ be a directed set. We consider the constant $\underline{\Lambda}$-sorted Lawvere theory $\mathbb{A}^{\Lambda}_{\mathrm{C}}$ .

There is an obvious functor $P:\mathrm{Alg}_{\mathbb{A}_{C}}\rightarrow\mathrm{Comm}(\mathrm{C})$. It sends $\mathrm{Free}_{\mathbb{A}_{C}}(\underline{\lambda})$ to $\mathrm{Sym}(\mathbb{I}^{\oplus |\lambda|})$ and then extends by sifted colimits. The following is clear.

\begin{prop}\label{prop:tinyTff}
If $k$ is tiny (i.e. mapping out of it commutes with all filtered colimits) then the functor
$$\mathrm{Alg}_{\mathbb{A}_{C}}^{f}\rightarrow\mathrm{Comm(C)}$$
is fully faithful.
\end{prop}

%
%

\begin{defn}
\begin{enumerate}
\item
A $\underline{\Lambda}$-sorted Lawvere theory $\mathrm{T}$ is said to be of \textit{weak} $\mathrm{C}$-\textit{polynomial type} if 
\begin{enumerate}
\item
There is a coproduct-preserving functor \[\mathrm{F}:\mathrm{Alg_{T}}^{ffg}\rightarrow\mathrm{Comm(C))}.\]
\item
There is a natural transformation of functors $\Lambda^{op}\rightarrow\mathrm{Comm(C)}$
$$P\circ\mathrm{Free}_{\mathbb{A}^{\Lambda}_{C}}\rightarrow \mathrm{F}(\mathrm{Free}_{T})$$
\end{enumerate}
\item
A $\underline{\Lambda}$-sorted Lawvere theory of weak $\mathrm{C}$-polyonomial type is said to be of $\mathrm{C}$-\textit{polyonomial type} if in addition
\begin{enumerate}
\item
$\mathrm{F}:\mathrm{Alg_{T}}^{ffg}\rightarrow\mathrm{Comm(Fun(\Lambda,C))}$ is fully faithful.
\item
For each $f:S\rightarrow\Lambda$ in $\mathrm{FinSet}_{\big\slash\Lambda}$ the map
$$P\circ\mathrm{Free}_{\mathbb{A}^{\Lambda}_{C}}(f)\rightarrow \mathrm{F}\circ\mathrm{Free}_{T}(f)$$
is an epimorphism.
\end{enumerate}
\end{enumerate}
\end{defn}

The idea we are trying to capture here is that the Lawvere theory $\mathrm{T}$ is somehow generated by algebras $A_{\lambda}$ internal to $\mathrm{C}$ correcsonding to discs of a certain radius $\lambda$, and the map $\mathbb{A}^{\Lambda}_{\mathrm{C}}\rightarrow\mathrm{T}$ determine maps $\mathrm{Sym}(\mathbb{I})\rightarrow A_{\lambda}$. Later in the bornological setup, $A_{\lambda}$ will be certain completions of the bornological polynomial algebra.

\begin{defn}
Let $\mathrm{T}$ be a $\underline{\Lambda}$-sorted Lawvere theory of $\mathrm{C}$-polynomial type. $\mathrm{T}$ is said to be \textit{concretely of} $\mathrm{C}$-\textit{polynomial type} if 
for any $A\in\mathrm{Alg_{T}}$ the map
$$\limind_{\lambda}\mathrm{Hom}(\mathrm{Free}_{\mathrm{T}}(\lambda),A)\rightarrow\mathrm{Hom}(\mathrm{F}(\mathrm{Free_{T}}(\lambda)),\mathrm{F}(A))\rightarrow\mathrm{Hom}(\mathrm{Sym}(\mathbb{I}),\mathrm{F}(A))$$
is an isomorphism.
\end{defn}


\begin{prop}\label{prop:freeff}
Let $\mathrm{T}$ be a Lawvere theory concretely of $\mathrm{C}$-polynomial type, with $\mathrm{C}$ additive. Suppose that $\mathbb{I}$ is tiny in $\mathrm{C}$. Then the functor
$$\mathrm{Alg}_{\mathrm{T}}^{f}\rightarrow\mathrm{Comm(\mathrm{C})}$$
is fully faithful.
\end{prop}

\begin{proof}
This proof is a generalisation of the proof of Proposition 25 in \cite{borisov2017quasi}, and follows the same strategy. We need to check that for $S\in\mathrm{Set}_{\big\slash\Lambda}$, the map
\begin{align*}
& \mathrm{Hom}_{\mathrm{Comm}(\mathrm{C})}(\mathrm{F}\circ\mathrm{Free}_{T}(K),\limind_{S'\subset S, S'\mathrm{ finite}}\mathrm{F}\circ\mathrm{Free}_{T}(S'))\\
&\rightarrow\limind_{S'\subset S, S'\mathrm{ finite}}\mathrm{Hom}_{\mathrm{Comm}(\mathrm{C})}(\mathrm{F}\circ\mathrm{Free}_{T}(K),\mathrm{F}\circ\mathrm{Free}_{T}(S))
\end{align*}
is an isomorphism for any finite $K$. Now the map \[\mathrm{Sym}(\mathbb{I}^{\oplus |K|})\rightarrow \mathrm{F}\circ\mathrm{Free}_{T}(K)\rightarrow\limind_{S'\subset S, S\mathrm{ finite}}\mathrm{F}\circ\mathrm{Free}_{T}(S')\] does factor through some $\mathrm{F}\circ\mathrm{Free}_{T}(S')$ since $\mathbb{I}^{\oplus |K|}$ is tiny. For $S'\subset S$ with  the map $\mathrm{F}\circ\mathrm{Free}_{T}(S')\rightarrow \mathrm{F}\circ\mathrm{Free}_{T}(S)$ is a split monomorphism. Indeed a left inverse is induced by the restriction map. In particular it is a regular monomorphism. Thus the regular image of the map $\mathrm{Sym}(\mathbb{I}^{\oplus |K|})\rightarrow \mathrm{F}\circ\mathrm{Free}_{T}(S)$ factors through $\mathrm{F}\circ\mathrm{Free}_{T}(S')$. But $\mathrm{Sym}(\mathbb{I}^{\oplus |K|})\rightarrow \mathrm{F}\circ\mathrm{Free}_{T}(K)$ is an epimorphism, so its regular image is just $\mathrm{F}\circ\mathrm{Free}_{T}(K)$. Hence $\mathrm{F}\circ\mathrm{Free}_{T}(K)\rightarrow\limind_{S'\subset S, S\mathrm{ finite}} \mathrm{F}\circ\mathrm{Free}_{T}(S')$ factors through some $\mathrm{F}\circ\mathrm{Free}_{T}(S')$ as required.
\end{proof}

\begin{lem}\label{lem:1sortedff}
    Let $\mathrm{T}$ be a one-sorted Lawvere theory concretely of $\mathrm{C}$-polynomial type. If $\mathbb{I}$ is tiny and projective in $\mathrm{C}$, i.e. $\mathrm{Hom}(\mathbb{I},-)$ commutes with sifted colimits. Then the functor
    $$\mathrm{Alg}_{\mathrm{T}}\rightarrow\mathrm{Comm(\mathrm{C})}$$
    is fully faithful.
\end{lem}

\begin{proof}
We need to show that the map
$$\mathrm{Hom}(A,B)\rightarrow\mathrm{Hom}(\mathrm{F}(A),\mathrm{F}(B))$$
is an isomorphism. Since any object of $\mathrm{Alg}_{\mathrm{T}}$ may be written as a colimit of $\mathrm{Free_{T}}(1)$ we may assume that $A\cong\mathrm{Free_{T}}(1)$. We may also write $B\cong\limind_{\mathcal{I}}\mathrm{Free_{T}}(n_{i})$ as a sifted colimit of free algebras. Now we have that the composite
\begin{align*}
    \limind_{\mathcal{I}}\mathrm{Hom}(\mathrm{Free_{T}}(1),\mathrm{Free_{T}}(n_{i}))&\cong\mathrm{Hom}(\mathrm{Free_{T}}(1),B)\\
    &\rightarrow\mathrm{Hom}(\mathrm{F}\circ\mathrm{Free_{T}}(1),F(B))\\
    &\cong\mathrm{Hom}(\mathrm{F}\circ\mathrm{Free_{T}}(1),\limind_{\mathcal{I}} \mathrm{F}\circ\mathrm{Free_{T}}(n_{i}))\\
    &\rightarrow\mathrm{Hom}(\mathrm{Sym}(\mathbb{I}),\limind_{\mathcal{I}} \mathrm{F}\circ\mathrm{Free_{T}}(n_{i}))\\
    &\cong\limind_{\mathcal{I}}\mathrm{Hom}(\mathrm{Sym}(\mathbb{I}),\mathrm{F}\circ\mathrm{Free_{T}}(n_{i}))
\end{align*}
is an isomorphism. Thus for any $B$ the composite
$$\mathrm{Hom}(\mathrm{Free_{T}}(1),B)\rightarrow\mathrm{Hom}(\mathrm{F}\circ\mathrm{Free_{T}}(1),F(B))\rightarrow\mathrm{Hom}(\mathrm{Sym}(\mathbb{I}),F(B))$$
is an isomorphism. In particular $\mathrm{Hom}(\mathrm{F}\circ\mathrm{Free_{T}}(1),F(B))\rightarrow\mathrm{Hom}(\mathrm{Sym}(\mathbb{I}),F(B))$ is an epimorphism. But we know that this is a monomorphism since $\mathrm{Sym}(\mathbb{I})\rightarrow \mathrm{F}\circ\mathrm{Free_{T}}(1)$ is an epimorphism. Thus  $\mathrm{Hom}(\mathrm{F}\circ\mathrm{Free_{T}}(1),F(B))\rightarrow\mathrm{Hom}(\mathrm{Sym}(\mathbb{I}),F(B))$ is an isomorphism. Since the composite is an isomorphism, $\mathrm{Hom}(\mathrm{Free_{T}}(1),B)\rightarrow\mathrm{Hom}(\mathrm{F}(\mathrm{Free_{T}}(1)),F(B))$ is also an isomorphism.  
\end{proof}

%

\begin{defn}
    Let $\mathpzc{C}$ be a combinatorial simplicial model category with simplicial mapping space functor $\underline{\mathrm{Hom}}$, and initial object $\mathbb{I}$. A fibrant object $X$ of $\mathpzc{C}$ is said to be \textit{weakly discrete} if $\underline{\mathrm{Hom}}(\mathbb{I},X)$ is discrete as a simplicial set, i.e. it is a disjoint union of contractible simplicial sets.
\end{defn}

\begin{defn}\label{defn:indhtpypoly}
Let $(\mathpzc{C},\otimes,\mathrm{Id},\underline{\mathrm{Hom}})$ be a combinatorial simplicial monoidal model category satisfying the monoid axiom and the commutative monoid axiom.


A $\Lambda$-sorted Lawvere theory of weak $\mathpzc{C}$-polyonomial type is said to be of \textit{homotopy} $\mathpzc{C}$-\textit{polyonomial type} if in addition
\begin{enumerate}
\item
$F:\mathrm{Alg}^{ffg}_{T}\rightarrow\mathrm{Comm(\mathpzc{C})}$ is fully faithful.
\item
For each $f:S\rightarrow\Lambda$ in $\mathrm{FinSet}_{\big\slash\Lambda}$ the map
$$P\circ\mathrm{Free}_{\mathbb{A}^{\Lambda}_{\mathpzc{C}}}(f)\rightarrow \mathrm{F}\circ\mathrm{Free}_{T}(f)$$
is a homotopy epimorphism.
\item
For each $f:S\rightarrow\Lambda$ in $\mathrm{FinSet}_{\big\slash\Lambda}$, $\mathrm{F}\circ\mathrm{Free}_{T}(f)$ is weakly discrete and $K$-flat as an object of $\mathpzc{C}$.
\end{enumerate}

\end{defn}

\begin{lem}[following \cite{ben2021analytification}]\label{lem:htpyepiT}
Let $\mathpzc{C}$ be a combinatorial monoidal model category which satisfies the monoid axiom and the commutative monoid axiom, and such that the transferred model structure on $\mathrm{Comm}(\mathpzc{C})$ is left proper. Let $\mathrm{T}$ be a $\Lambda$-sorted Lawvere theory of homotopy $\mathpzc{C}$-polynomial type. Suppose further that
\begin{enumerate}
    \item 
    filtered colimits in $\mathpzc{C}$ present homotopy filtered colimits.
    \item 
    geometric realisations of simplicial objects in $\mathrm{s}\mathpzc{C}$ present homotopy geometric realisations.
\end{enumerate}
Then
\begin{enumerate}
\item
Let $S$ be any object of $\mathrm{Set}_{\big\slash\Lambda}$. The map $\mathrm{S}(S)\rightarrow |F(\mathrm{T}_{S})|$ is a homotopy epimorphism.
\item
Let $K$ be any object of $\mathrm{sSet}_{\big\slash\Lambda}$. The map $\mathrm{S}(K)\rightarrow |F(\mathrm{T}_{K})|$ is a homotopy epimorphism.
\end{enumerate}
\end{lem}

\begin{proof}
\begin{enumerate}
\item
Homotopy epimorphisms are stable by homotopy colimits. $F(\mathrm{T}_{S})\cong\limind_{S'\subset S, S\mathrm{ finite}}F(\mathrm{T}_{S'})$ is a filtered colimit. Moreover it is a homotopy filtered colimit by assumption. \[\mathrm{S}(S)\cong\limind_{S'\subset S, S\mathrm{ finite}}\mathrm{S}(S')\] is also a homotopy filtered colimit. This proves the claim. 
\item
 It suffices to observe that the map $\mathrm{S}(K)\rightarrow F(\mathrm{T}_{K})$ is the homotopy colimit of the simplicial diagram of maps $\mathrm{S}(K_{n})\rightarrow F(\mathrm{T}_{K_{n}})$, each of which is a homotopy epimorphism by part 1). 
\end{enumerate}
\end{proof}

Now let $\mathrm{T}$ be a $\Gamma$-sorted Lawvere theory of homotopy $\mathpzc{C}$-polynomial type, with corresponding functor
$$F:\mathrm{Alg_{T}^{ffg}}\rightarrow\mathrm{Comm}(\mathpzc{C})$$
Now we have a homotopy epimorphism
$$\mathrm{Sym}(\mathbb{I})\rightarrow\mathrm{F}(\mathrm{Free_{T}}(\lambda))$$
Thus for any fibrant $A\in\mathrm{Comm}(\mathpzc{C})$ the map
$$\mathbf{Map}(\mathrm{F}(\mathrm{Free_{T}}(\lambda)),A)\rightarrow\mathbf{Map}(\mathrm{Sym}(\mathbb{I}),A)\cong\underline{\mathrm{Hom}}(\mathbb{I},A)$$
is a homotopy monomorphism. In particular if $A$ is weakly discrete then $\mathbf{Map}(\mathrm{F}(\mathrm{Free_{T}}(\lambda)),A)$ is a discrete simplicial set. Hence 
$$\mathbf{Map}(\mathrm{F}(\mathrm{Free_{T}}(\lambda)),\mathrm{F}(\mathrm{Free_{T}}(\gamma)))\rightarrow\pi_{0}(\mathbf{Map}(\mathrm{F}(\mathrm{Free_{T}}(\lambda)),\mathrm{F}(\mathrm{Free_{T}}(\gamma)))$$
is an equivalence for any $\lambda,\gamma$. Thus we get a well-defined map in the $(\infty,1)$-category $\mathbf{sSet}$, natural in $\mathrm{F}(\mathrm{Free_{T}}(\lambda))$ and $\mathrm{F}(\mathrm{Free_{T}}(\gamma))$,
$$\mathrm{Hom}_{\mathpzc{C}}(\mathrm{F}(\mathrm{Free_{T}}(\lambda)),\mathrm{F}(\mathrm{Free_{T}}(\gamma))\rightarrow\mathbf{Map}(\mathrm{F}(\mathrm{Free_{T}}(\lambda)),\mathrm{F}(\mathrm{Free_{T}}(\gamma)))$$
This gives a functor
$$\mathbf{F}:\mathrm{N}(\mathrm{Alg_{T}}^{ffg})\rightarrow\mathrm{L^{H}}(\mathrm{Comm}(\mathpzc{C}))$$
which extends by sifted colimits to a functor
$$\mathbf{F}:\mathbf{sAlg}_{\mathrm{T}}\rightarrow\mathrm{L^{H}}(\mathrm{Comm}(\mathpzc{C}))$$

\begin{lem}
Let $\mathrm{T}$ be of homotopy $\mathpzc{C}$-polynomial type (resp. of homotopy $\mathpzc{C}$-polynomial type). The functor $ \mathbf{sAlg}_{T}\rightarrow\mathrm{L^{H}(\mathrm{Comm}(Fun(\Lambda,\mathpzc{C}))}$ (resp. the functor $ \mathbf{sAlg}_{T}\rightarrow\mathrm{L^{H}(\mathrm{Comm}(\mathpzc{C})}$) commutes with colimits. In particular it is a left adjoint functor.
\end{lem}

\begin{proof}
It suffices to prove that it commutes coproducts. Let $C,D\in\mathbf{sAlg}_{T}$. Without loss of generality we may assume that $C$ and $D$ are cofibrant. Then
$$F(C)\otimes^{\mathbb{L}}F(D)\cong F(C)\otimes F(D)\cong F(C\coprod D)$$
as required. That $\mathbf{F}$ is a left adjoint follows from the fact that both categories are locally presentable $(\infty,1)$-categories.
\end{proof}

We will denote the right adjoint functor by $\mathbf{R}_{\mathrm{T}}$.

\subsection{Lawvere Theories Generated by Algebras}\label{sec:Lawveregenalg}

Here we give some conditions on a collection of algebras in $\mathrm{Comm}(\mathrm{C})$ such that the Lawvere theory they generate is of homotopy $\mathrm{C}$-polynomial type

Let $\mathrm{C}$ be a category and $\underline{\Lambda}$ a partially ordered set. Let  $A:\underline{\Lambda}^{op}\rightarrow\mathrm{Comm(\mathrm{C})},\lambda\mapsto A_{(\lambda)}$ be a functor. Define a $\underline{\Lambda}$-sorted Lawvere theory $\mathrm{T}^{A}$ as follows. Let $f:K\rightarrow\Lambda$ be an object of $\mathrm{FinSet}_{\big\slash\Lambda}$ Define
$$A_{(f)}\defeq \bigotimes_{\lambda\in\Lambda}A^{\otimes|f^{-1}(\lambda)|}_{(\lambda)}$$
Then define
$$\Hom_{\mathrm{T}^{A}}(f:K\rightarrow\Lambda,g:L\rightarrow\Lambda)\defeq\Hom_{\mathrm{Comm(C)}}(A_{(g)},A_{(f)})$$
This is a $\underline{\Lambda}$-sorted Lawvere theory. Moreover there is clearly a fully faithful coproduct-preserving functor $F^{A}:\mathrm{Alg_{T^{A}}^{ffg}}\rightarrow\mathrm{Comm(C)}$.

Suppose that for each $\lambda\in\Lambda$ there is a map $\tilde{p}_{\lambda}:\mathbb{I}\rightarrow A_{(\lambda)}$ in $\mathrm{C}$  such that for any $\lambda'\le\lambda$ the diagram 
\begin{displaymath}
\xymatrix{
\mathbb{I}\ar[dr]^{\tilde{p}_{\lambda'}}\ar[r]^{\tilde{p}_{\lambda}} & A_{(\lambda)}\ar[d]\\
& A_{(\lambda')}
}
\end{displaymath}
commutes. Each $\tilde{p}_{\lambda}$ extends to a unique map of commutative monoids $p_{\lambda}:\mathrm{Sym}(\mathbb{I})\rightarrow A_{(\lambda)}$ such that
\begin{displaymath}
\xymatrix{
\mathrm{Sym}(\mathbb{I})\ar[dr]^{p_{\lambda'}}\ar[r]^{p_{\lambda}} & A_{(\lambda)}\ar[d]\\
& A_{(\lambda')}
}
\end{displaymath}
By tensoring, for each $f:K\rightarrow\Lambda$ in $\mathrm{FinSet}_{\big\slash\Lambda}$ we get a map
$$p_{f}:\mathrm{Sym}(\mathbb{I}^{\oplus |K|})\rightarrow A_{(f)}.$$
Moreover if $\phi:(f:K\rightarrow\Lambda)\rightarrow (g:L\rightarrow\Lambda)$ is a map in $\mathrm{FinSet}_{\big\slash\Lambda}$ then the diagram 
\begin{displaymath}
\xymatrix{
\mathrm{Sym}(\mathbb{I}^{\oplus |K|})\ar[dr]^{p_{g}}\ar[r]^{p_{f}} & A_{(f)}\ar[d]\\
& A_{(g)}
}
\end{displaymath}
commutes.

Note that by construction the maps $p_{\lambda}$, induce a natural transformation
$$P\circ\mathrm{Free}_{\mathbb{A}^{\Lambda}_{\mathrm{C}}}\rightarrow\mathrm{F}\circ\mathrm{Free_{T}}$$

Now let $(\mathpzc{C},\otimes,\mathrm{Id},\underline{\mathrm{Hom}})$ be a combinatorial closed monoidal model category satisfying the monoid axiom and the commutative monoid axiom.
\begin{lem}
If
\begin{enumerate}
\item
 for each $\lambda\in\Lambda$ $p_{\lambda}:\mathrm{Sym}(\mathbb{I})\rightarrow A_{(\lambda)}$ is a homotopy epimorphism
 \item
 for each $\lambda\in\Lambda$ $A_{(\lambda)}$ is $K$-flat as an object of $\mathpzc{C}$.
 \item 
 for each $\underline{\lambda}\in\Lambda^{n}$ $A_{\underline{\lambda}}$ is weakly discrete as an object of $\mathpzc{C}$

 \end{enumerate}
 then $T^{A}$ is of homotopy $\mathpzc{C}$-polynomial type.
\end{lem}

\begin{proof}
It suffices to observe that by the transversality assumptions, the fact that $p_{\lambda}:\mathrm{Sym}(\mathbb{I})\rightarrow A_{(\lambda)}$ is a homotopy epimorphism implies that for any $f:K\rightarrow\Lambda$ in $\mathrm{FinSet}_{\big\slash\Lambda}$ the map 
$$\mathrm{Sym}(\mathbb{I}^{\oplus |K|})\rightarrow A_{f}$$
is a homotopy epimorphism. 
\end{proof}

\begin{example}
    For example, if $\mathrm{C}$ is locally presentable closed monoidal additive category, and $\mathcal{Z}$ a set of tiny objects such that the model structure induced by $\mathcal{Z}$ makes $\mathrm{sC}$ a combinatorial monoidal model category satisfying the monoid axiom and the commutative monoid axiom, then any Lawvere theory determined by $K$-flat algebras $A_{(\lambda)}$ in $\mathrm{C}$ equipped with compatible homotopy epimorphisms $p_{\lambda}:\mathrm{Sym}(\mathbb{I})\rightarrow A_{(\lambda)}$, then $\mathrm{T^{A}}$ is of homotopy $\mathrm{sC}$-polynomial type. 
\end{example}

Can we make this a Fermat theory in a natural way? We define an $\underline{\Lambda}$-filtered Lawvere-theory $\mathbb{A}^{A}_{\mathrm{C}}$ as follows. For $\underline{\lambda}=(\lambda_{1},\ldots,\lambda_{n})$ and $\underline{\gamma}=(\gamma_{1},\ldots,\gamma_{m})$ we define $\mathrm{Hom}(\underline{\lambda},\underline{\gamma})$ to be the fibre product
\begin{align*}
    \xymatrix{
    \mathrm{Hom}(\underline{\lambda},\underline{\gamma})\ar[d]\ar[r] & \mathrm{Hom}_{\mathrm{Comm}(\mathrm{C})}(A_{\gamma_{1}}\otimes\ldots\otimes A_{\gamma_{m}},A_{\lambda_{1}}\otimes\ldots A_{\lambda_{n}})\ar[d]\\
     \mathrm{Hom}_{\mathrm{Comm}(\mathrm{C})}(\mathrm{Sym}(\mathbb{I}^{\oplus m}),\mathrm{Sym}(\mathbb{I}^{\oplus n}))\ar[r] &  \mathrm{Hom}_{\mathrm{Comm}(\mathrm{C})}(\mathrm{Sym}(\mathbb{I}^{\oplus m}),A_{\lambda_{1}}\otimes\ldots A_{\lambda_{n}})
    }
\end{align*}
Note that $ \mathrm{Hom}(\underline{\lambda},\underline{\gamma})\rightarrow   \mathrm{Hom}_{\mathrm{Comm}(\mathrm{C})}(\mathrm{Sym}(\mathbb{I}^{\oplus m}),\mathrm{Sym}(\mathbb{I}^{\oplus n}))$ is a monomorphism. We shall assume that for each $\underline{\lambda}\in\Lambda^{n}$ the natural map 
$$\colim_{\gamma\in\underline{\Lambda}} \mathrm{Hom}_{\mathrm{Comm}(\mathrm{C})}(A_{\gamma},A_{\lambda_{1}}\otimes\ldots A_{\lambda_{n}})\rightarrow \mathrm{Hom}_{\mathrm{Comm}(\mathrm{C})}(\mathrm{Sym}(\mathbb{I}),A_{\lambda_{1}}\otimes\ldots A_{\lambda_{n}})$$
is an isomorphism. In particular this means that for fixed $\underline{\lambda}$ the map 
$$\colim_{\gamma\in\Lambda^{m}}\mathrm{Hom}(\underline{\lambda},\gamma)\rightarrow\mathrm{Hom}_{\mathrm{Comm(C)}}(\mathrm{Sym}(\mathbb{I}),\mathrm{Sym}(\mathbb{I}^{\oplus n}))$$
is an isomorphism, and thus $\mathbb{A}^{A}_{\mathrm{C}}$ is an $\underline{\Lambda}$-filtered polynomial Lawvere theory. There is an obvious functor
$$P_{A}:\mathrm{Alg}_{\mathbb{A}^{A}_{\mathrm{C}}}^{ffg}\rightarrow\mathrm{Comm(C)}$$ 
and by construction a natural transformation
$$P_{A}\circ\mathrm{Free}_{\mathbb{A}_{\mathrm{C}}^{\Lambda}}\rightarrow  F^{A}\circ\mathrm{Free}_{T^{A}}$$
Since coproducts of epimorphisms are epimorphisms, $T^{A}$ will be of $\mathrm{C}$-polynomial type precisely if each $p_{\lambda}$ is an epimorphism in $\mathrm{Comm(C)}$.

\begin{defn}
    We say a $\Gamma$-sorted Lawvere theory $\mathrm{T}^{A}$ generated by algebras $A_{(\gamma)}$ in $\mathrm{C}$ is \textit{Fermat} if the map
    $$\mathbb{A}^{A}_{\mathrm{C}}\rightarrow\mathrm{T}^{A}$$
    realises it as a $\Gamma$-sorted Fermat theory. 
\end{defn}

\subsection{Embedding Lawvere Theories in Categories of Algebras}

Let $\mathrm{T}$ be a $\Gamma$-sorted Lawvere theory of ind homotopy $\mathpzc{C}$-polynomial type. We would like to understand when the functor 
$$\mathbf{F}:\mathbf{sAlg}_{\mathrm{T}}\rightarrow\mathrm{L^{H}}(\mathrm{Comm}(\mathpzc{C}))$$
is fully faithful.

\begin{defn}
Let $\mathrm{T}$ be a $\Gamma$-sorted Lawvere theory of homotopy $\mathpzc{C}$-polynomial type. $\mathrm{T}$ is said to be \textit{concretely of homotopy} $\mathpzc{C}$-\textit{polynomial type} if the map
$$\limind_{\lambda}\mathbf{Map}(\mathrm{Free}_{T}(\lambda),\mathrm{Free}_{T}(\gamma))\rightarrow\mathbf{Map}(\mathbf{F}(\mathrm{Free}_{T}(\lambda)),\mathbf{F}(\mathrm{Free}_{T}(\gamma)))\rightarrow\mathbf{Map}(\mathrm{Sym}(\mathbb{I}),\mathbf{F}(\mathrm{Free}_{T}(\gamma)))$$
is an equivalence for any $\gamma\in\Lambda$.
\end{defn}

\begin{lem}\label{lem:Temb}
    Let $\mathrm{T}$ be a $1$-sorted Lawvere theory which is concretely of homotopy $\mathpzc{C}$-polynomial type. Suppose that $\mathbf{Map}(\mathbb{I},-)$ commutes with sifted colimits in $\mathrm{L^{H}}(\mathpzc{C})$. Then the functor 
    $$\mathbf{F}:\mathbf{sAlg}_{\mathrm{T}}\rightarrow\mathrm{L^{H}}(\mathrm{Comm}(\mathpzc{C}))$$
    is fully faithful.
\end{lem}

\begin{proof}
We need to show that for any $A,B\in\mathbf{sAlg}_{\mathrm{T}}$ the map
$$\mathbf{Map}(A,B)\rightarrow\mathbf{Map}(\mathbf{F}(A),\mathbf{F}(B))$$
is an equivalence.
  Since any object $A$ of $\mathbf{sAlg}_{\mathrm{T}}$ is a colimit of objects of the form $\mathrm{Free_{T}}(1)$, we may assume that $A\cong\mathrm{Free_{T}}(1)$. Now let $B\in\mathbf{sAlg}_{\mathrm{T}}$. We may write $B$ as a sifted colimit $B\cong\colim_{\mathcal{I}}\mathrm{Free_{T}}(n_{i})$. Now consider the composite
\begin{align*}
\colim_{\mathcal{I}}\mathbf{Map}(\mathrm{Free_{T}}(1),\mathrm{Free_{T}}(n_{i}))&\cong\mathbf{Map}(\mathrm{Free_{T}}(1),\colim_{\mathcal{I}}\mathrm{Free_{T}}(n_{i}))\\
  &\rightarrow\mathbf{Map}(\mathbf{F}(\mathrm{Free_{T}}(1)),\colim_{\mathcal{I}}\mathbf{F}(\mathrm{Free_{T}}(n_{i})))\\
  &\rightarrow\mathbf{Map}(\mathrm{Sym}(\mathbb{I}),\colim_{\mathcal{I}}\mathbf{F}(\mathrm{Free_{T}}(n_{i})))\\
  &\cong\colim_{\mathcal{I}}\mathbf{Map}(\mathrm{Sym}(\mathbb{I}),\mathbf{F}(\mathrm{Free_{T}}(n_{i})))
\end{align*}
Now this composite is an equivalence. Thus
$$\mathbf{Map}(\mathrm{Free_{T}}(1),B)\rightarrow\mathbf{Map}(\mathbf{F}(\mathrm{Free_{T}}(1)),\mathbf{F}(B))\rightarrow\mathbf{Map}(\mathrm{Sym}(\mathbb{I}),\mathbf{F}(B))$$
is an equivalence for any $B$. Now also for any $B$ the map $\mathbf{Map}(\mathbf{F}(\mathrm{Free_{T}}(1)),\mathbf{F}(B))\rightarrow\mathbf{Map}(\mathrm{Sym}(\mathbb{I}),\mathbf{F}(B))$ is a homotopy monomorphism. Thus it is a monomorphism on $\pi_{0}$ and an equivalence on $\pi_{n}$ for any $n$ (and any basepoint). Since the composite
$$\pi_{0}\mathbf{Map}(\mathrm{Free_{T}}(1),B)\rightarrow\pi_{0}\mathbf{Map}(\mathbf{F}(\mathrm{Free_{T}}(1)),\mathbf{F}(B))\rightarrow\pi_{0}\mathbf{Map}(\mathrm{Sym}(\mathbb{I}),\mathbf{F}(B))$$
is an isomorphism, the right-hand map is also an epimorphism. Thus \[\mathbf{Map}(\mathbf{F}(\mathrm{Free_{T}}(1)),\mathbf{F}(B))\rightarrow\mathbf{Map}(\mathrm{Sym}(\mathbb{I}),\mathbf{F}(B))\] is also an epimorphism on $\pi_{0}$ and hence is an equivalence. Thus the map \[\mathbf{Map}(\mathrm{Free_{T}}(1),B)\rightarrow\mathbf{Map}(\mathbf{F}(\mathrm{Free_{T}}(1)),\mathbf{F}(B))\] is also an equivalence.
\end{proof}


\section{Lawvere Theories and Derived Algebraic Contexts}

Let us fix a derived algebraic context
$$\underline{\mathbf{C}}\defeq(\mathbf{C},\mathbf{C}_{\ge0},\mathbf{C}_{\le0},\mathbf{LSym},\theta)$$
The heart $\mathbf{C}^{\heart}$ is a monoidal elementary abelian category with symmetric projectives, and the model category $\mathrm{s}\mathbf{C}^{\heart}$ presents $\mathbf{C}_{\ge0}$. In this context we will say that a Lawvere theory is of homotopy $\mathbf{C}^{\heart}$-polynomial type if it is of homotopy $\mathrm{s}\mathbf{C}^{\heart}$-polynomial type. 

\subsection{Homotopy Polynomial Type in Derived Algebraic Contexts}

Let us give a simple formulation of what it means to be of homotopy polynomial type. Let $\mathrm{T}$ be a $\Lambda$-sorted Lawvere theory of weak $\mathrm{C}$-polynomial type, with corresponding functor $\mathrm{F}:\mathrm{Alg_{T}^{ffg}}\rightarrow\mathrm{Comm(C)}$. To ease notational burden we will write $\mathrm{F}^{\mathrm{C}}(\underline{\lambda})\defeq\mathrm{F}(\mathrm{Free_{T}}(\underline{\lambda})$. Suppose that each $\mathrm{F}^{\mathrm{C}}(\underline{\lambda})$ is $K$-flat and that $\limind_{\Lambda}\circ\mathrm{F}$ is fully faithful. Since tensor products of homotopy epimorphisms are homotopy epimorphisms, to show that $\mathrm{Sym}(\mathbb{I}^{\oplus n})\rightarrow\mathrm{F^{C}}(\underline{\lambda})$ is a homotopy epimorphism for each $\underline{\lambda}\in\Lambda^{n}$, it suffices to show that each
$$\mathrm{Sym}(\mathbb{I})\rightarrow\mathrm{F^{C}}(\lambda)$$
is a homotopy epimorphism for each $\lambda\in\Lambda$. 

Now we have an exact sequence
$$0\rightarrow\mathrm{Sym}(\mathbb{I})\otimes\mathrm{Sym}(\mathbb{I})\rightarrow\mathrm{Sym}(\mathbb{I})\otimes\mathrm{Sym}(\mathbb{I})\rightarrow \mathrm{Sym}(\mathbb{I})\rightarrow0$$
where the first map is the `multiplication by $x-y$' map. Thus the $2$-term complex
$$\mathrm{Sym}(\mathbb{I})\otimes\mathrm{Sym}(\mathbb{I})\rightarrow\mathrm{Sym}(\mathbb{I})\otimes\mathrm{Sym}(\mathbb{I}))$$
is a resolution. Again since $\mathrm{F}^{C}(\lambda)$ is flat, $\mathrm{F}^{C}(\lambda)\otimes^{\mathbb{L}}_{\mathrm{Sym}(\mathbb{I})}\mathrm{F}^{C}(\lambda)$ may be computed as the $2$-term complex
$$\mathrm{F^{C}}(\lambda,\lambda)\rightarrow\mathrm{F^{C}}(\lambda,\lambda)$$
where again the map is multiplication by $x-y$. Thus $\mathrm{T}$ being of homotopy $\mathpzc{C}$-polynomial type boils down to showing that this complex is quasi-isomorphic to $\mathrm{F^{C}}(\lambda)$ for each $\lambda\in\Lambda$.

\subsection{Concreteness in Derived Algebraic Contexts}

\begin{lem}\label{lem:stricT}
Let $(\mathrm{C},\otimes,\mathbb{I},\underline{\mathrm{Hom}})$ be a locally presentable closed monoidal additive category with $\mathbb{I}$ tiny. Let $\mathcal{Z}$ be a set of tiny objects containing the monoidal unit such that the $\mathcal{Z}$-model structure on $\mathrm{sC}$ is a combinatorial monoidal model structure which satisfies the monoid axiom and the commutative monoid axiom.
Let $\mathrm{T}$ be a $\Gamma$-sorted Lawvere theory of homotopy $\mathpzc{C}$-polynomial type. Suppose that 
$$\mathrm{F}:\mathrm{sAlg}_{\mathrm{T}}\rightarrow\mathrm{sComm(C)}$$
preserves weak equivalences between fibrant-cofibrant objects. If $\mathrm{T}$ is concreteley of $\mathrm{C}$-polynomial type then for any fibrant-cofibrant $\mathrm{T}$-algebra $Y$ the map
$$\limind_{\Gamma}\underline{\mathrm{Hom}}(\mathrm{Free}_{T}(\gamma),Y)\rightarrow\underline{\mathrm{Hom}}(\mathrm{Sym}(\mathbb{I},\mathrm{F}(Y))$$
is an equivalence. In particular $\mathrm{T}$ is concretely of homotopy $\mathpzc{C}$-polynomial type.
\end{lem}

\begin{proof}
Let $Y$ be a fibrant-cofibrant simplicial $\mathrm{T}$-algebra. We need to show that the map
$$\underline{\mathrm{Hom}}(\mathrm{T}(1),Y)\rightarrow\underline{\mathrm{Hom}}_{\mathrm{sComm(C)}}(\mathrm{Sym}(\mathbb{I}),\mathbf{F}(Y))$$
is an equivalence. Since $Y$ is cofibrant it is degree-wise free. By the proof of Proposition \ref{prop:tinyTff} we may assume the each $Y_{n}$ is free of finite type. Then the claim follows by assumption.
\end{proof}

\subsection{The Induced Map On $\pi_{0}$}

Let $\mathrm{T}$ be a $\Lambda$-sorted Lawvere theory of homotopy $\mathbf{C}^{\heart}$ polynomial type.


There is a functor
$$\mathrm{F_{T}}:\mathrm{Alg_{T}}\cong\mathcal{P}_{\Sigma}(\{\mathrm{Free_{T}}(\underline{\lambda}):\lambda\in\Lambda^{n},n\in\mathbb{N}\})\rightarrow\mathrm{Comm}(\mathbf{C}^{\heart})$$
which extends by $1$-categorical colimits the functor sending $\mathrm{Free_{T}}(\underline{\lambda})$ to $\mathbf{F}(\mathrm{Free_{T}}(\underline{\lambda}))$.

\begin{prop}
Let $B\in\mathbf{sAlg}_{\mathrm{T}}$. There is a natural equivalence 
$$\pi_{0}(\mathbf{F}(B))\cong\mathrm{F}(\pi_{0}(B))$$
\end{prop}

\begin{proof}
We may assume that $B$ is cofibrant, so that it has a presentation which degree-wise is $\mathrm{T}(g_{n})$ for some function $g:X\rightarrow\Lambda$, and all face maps send generators to generators. Then $\mathbf{F}(B)\cong |\mathrm{F}(\mathrm{T}(g_{n}))|$ so that $\pi_{0}(\mathbf{F}(B))$ is given by the (reflexive) coequaliser of the maps $d_{0},d_{1}:\mathrm{F}(\mathrm{T}(g_{1}))\rightarrow \mathrm{F}(\mathrm{T}(g_{0}))$. But this is precisely $\mathrm{F}(\pi_{0}(B))$.
\end{proof}

\subsection{Functors Between Categories of Modules}

Let $\mathrm{T}$ be a $\Gamma$-sorted Fermat theory of homotopy $\mathbf{C}^{\heart}$-polynomial type. Let $B\in\mathbf{DAlg}^{cn}(\mathbf{C})$, and consider the category
$${}_{B}\mathbf{Mod}(\mathbf{C}_{\ge0})\cong\mathbf{Ab}(\mathbf{DAlg}^{cn}(\mathbf{C})_{/B})$$
The functor 
$$\mathbf{R}_{\mathrm{T}}:\mathbf{DAlg}^{cn}(\mathbf{C})\rightarrow\mathbf{sAlg}_{\mathrm{T}}$$
induces a right adjoint functor
$$\mathbf{R}_{B}:{}_{B}\mathbf{Mod}(\mathbf{C}_{\ge0})\cong\mathbf{Ab}(\mathbf{DAlg}^{cn}(\mathbf{C})_{/B})\rightarrow\mathbf{Ab}(\mathbf{sAlg}_{T_{/\mathbf{R}_{\mathrm{T}}(B)}})\cong{}_{\mathbf{R}_{\mathrm{T}}(B)}\mathbf{Mod}^{\mathrm{Beck}}$$
Now as $\mathrm{T}$ is a Fermat theory we have 
$${}_{\mathbf{R}_{\mathrm{T}}(B)}\mathbf{Mod}^{\mathrm{Beck}}\cong{}_{|\mathbf{R}_{\mathrm{T}}(B)|_{\Gamma-alg}}\mathbf{Mod}(\mathbf{Ch}_{\ge0}(\mathrm{Filt}(\mathrm{Ab})))$$ 

 Let $C\in\mathbf{DAlg}^{cn}(\mathbf{C})$. Then $|\mathbf{R}_{\mathrm{T}}(C)|_{\Gamma-alg}$ is the filtered algebra with underlying algebra $|C|_{alg}\defeq\mathbf{Map}(\mathrm{Sym}(\mathbb{I}),C)$ and filtration given by
$$C_{\lambda}\defeq\mathbf{Map}(\mathrm{Free_{T}}(\lambda),C)$$
for $\lambda\in\Gamma$. Now let $(\pi:C\rightarrow B)\in\mathbf{Ab}(\mathbf{DAlg}^{cn}(\mathbf{C})_{\big\slash B})$. Then $(|\mathbf{R}_{\mathrm{T}}(\pi)|_{\Gamma-alg}:(|\mathbf{R}_{\mathrm{T}}(C)|_{\Gamma-alg}\rightarrow|\mathbf{R}_{\mathrm{T}}(B)|_{\Gamma-alg})$ is in $\mathbf{Ab}(\mathbf{DAlg}^{cn}(\mathbf{Filt}(\mathbf{Ab}))_{\big\slash \overline{B}})$. $\mathrm{Ker}(\overline{\pi})$ is a filtered $|\mathbf{R}_{\mathrm{T}}(B)|_{\Gamma-alg}$-module. Thus the functor $\mathbf{R}_{B}$ equips $\overline{\mathrm{Ker}(\pi)}\defeq\mathbf{Map}(\mathbb{I},\mathrm{Ker}(\pi))$ with the induced subfiltration of $\overline{C}$.

Consider also the functor
$$|-|_{Alg}:{}_{B}\mathbf{Mod}(\mathbf{C}_{\ge0}))\rightarrow{}_{|B|_{alg}}\mathbf{Mod}(\mathbf{Ch}_{\ge0}(\mathrm{Ab}))$$
sending $M$ to $|M|_{alg}\defeq\mathbf{Map}_{{}_{B}\mathbf{Mod}(\mathbf{C}_{\ge0}}(B,M)$. This is also a right-adjoint functor with left-adjoint by $\mathbb{L}i_{B}$ from Subsection \ref{subsec:embeddingalgebra}.
Note that if $\mathrm{T}$ is one-sorted then $|-|_{alg}\cong\mathbf{R}_{B}$, and all these constructions fit together well.

\section{$\mathrm{T}$-Derived Quotients and Localisations}\label{sec:Tderloc}

In this section we fix a \textit{Koszul} derived algebraic context
$$(\mathbf{C},\mathbf{C}_{\ge0},\mathbf{C}_{\le0},\mathbf{LSym},\theta)$$
(Definition \ref{defn:KoszulDAC}).
Let $\mathrm{T}$ be a $\Gamma$-sorted Lawvere theory of homotopy $C$-polynomial type. We are going to consider some important examples of homotopy epimorphisms between objects in $\mathrm{Comm(C)}$ determined by $\mathrm{T}$.  For $\lambda\in\Gamma$ write $A_{(\lambda)}\defeq \mathrm{F}(\mathrm{Free}_{T}(\lambda))$. We shall assume that each $A_{(\lambda)}$ is flat over $R$. 

\begin{rem}
If $F(\mathrm{Free_{T}}(g))\rightarrow F(\mathrm{Free_{T}}(h))$ is a homotopy epimorphism (or indeed formally \'{e}tale) with $g$ and $h$ finite, then $|g|=|h|$. Indeed one has 
$$\mathbb{L}_{F(\mathrm{Free_{T}}(k))}\cong F(\mathrm{Free_{T}}(k))^{|k|}$$
Thus we get an isomorphism
$$F(\mathrm{Free_{T}}(h))^{|g|}\cong F(\mathrm{Free_{T}}(h))^{|h|}$$
Since $F(\mathrm{Free_{T}}(h))$ is commutative this implies that $|g|=|h|$.
\end{rem}

\begin{defn}
Let $A\in\mathbf{DAlg}^{cn}(\mathbf{C})$
\begin{enumerate}
    \item 
    A collection of maps
$$(a_{1},\ldots,a_{n}):\mathbb{I}^{\oplus n}\rightarrow |A|$$
in $\mathbf{C}$ is said to be $\mathrm{T}$-\textit{extendable} if there is $\underline{\lambda}\in\Gamma^{n}$ such that 
$$\mathrm{Sym}(\mathbb{I}^{\oplus n})\rightarrow A$$ 
factors through the map $\mathrm{Sym}(\mathbb{I}^{\oplus n})\rightarrow \mathrm{F}(\mathrm{Free}_{T}(\underline{\lambda}))$.
    \item 
   $A$ is said to be $\mathrm{T}$-\textrm{ideal complete} if any collection of maps
    $$(a_{1},\ldots,a_{n}):\mathbb{I}^{\oplus n}\rightarrow |A|$$
    is $\mathrm{T}$-extendable.
\end{enumerate}
\end{defn}

\begin{prop}\label{prop:Textend}
Let $\mathrm{T}$ be a $\Gamma$-sorted Fermat theory which is concretely of $\mathbf{C}^{\heart}$-polynomial type and of $\mathbf{C}^{\heart}$-homotopy polynomial type. Then any algebra $A$ such that $\pi_{0}(A)\cong\mathrm{F}(X)$ for some $X\in\mathrm{Alg_{T}}$ is $\mathrm{T}$-ideal complete.
\end{prop}

\begin{proof}
 This follows from Proposition \ref{prop:liftingmapsTpi0}.
\end{proof}




Let $(a_{1},\ldots,a_{n}):\mathbb{I}^{\oplus n}\rightarrow A$ be $\mathrm{T}$-extendable. Since $\mathrm{Sym}(R^{\oplus n})\rightarrow\mathrm{F}(\mathrm{Free}_{T}(\underline{\lambda}))$ is a homotopy epimorphism We have

\begin{align*}
A\big\slash\big\slash (a_{1},\ldots,a_{n})&\cong A\otimes_{\mathrm{Sym}(\mathbb{I}^{\oplus n})}^{\mathbb{L}}\mathbb{I}\\
&\cong A\otimes_{\mathrm{F}(\mathrm{Free}_{T}(\underline{\lambda}))}^{\mathbb{L}}\mathbb{I}\\
&\cong A\otimes_{\mathrm{F}(\mathrm{Free}_{T}(\underline{\lambda}))}B^{A}(n,\mathrm{Free}_{T}(f_{0}),\mathbb{I})
\end{align*}

\begin{rem}
In particular if $\mathrm{T}$ is concretely of $\mathbf{C}^{\heart}$-polynomial type, of $\mathbf{C}^{\heart}$ homotopy polynomial type, and $A$ is in the essential image of $\mathbf{F}$, then $A\big\slash\big\slash (a_{1},\ldots,a_{n})$ is also in the essential image of $\mathbf{F}$.
\end{rem}
%


\subsection{Derived Algebras}

Let $\underline{\mathbf{C}}$ be a Koszul derived algebraic context and let $\mathrm{T}$ be a $\Gamma$-sorted Lawvere theory of homotopy $\mathbf{C}^{\heart}$-polynomial type. 

\subsubsection{Derived $\mathrm{T}$-Coherence}

\begin{defn}\label{defm:strongTcoh}
\begin{enumerate}
    \item
    $A\in\mathbf{DAlg}^{\heart}(\mathbf{C})$ is said to be \textit{strongly }$\mathrm{T}$-\textit{coherent} if 
    \begin{enumerate}
        \item
         the class of $A$-algebras of the form $A\otimes\mathrm{T}(\underline{\lambda})\big\slash I$ with $I$ a finitely generated ideal is closed under finite colimits.
         \item 
         $\pi_{0}(A)\otimes\mathrm{T}(\underline{\lambda})\big\slash I$ is coherent for any $\underline{\lambda}$ and any finitely generated ideal $I$.
    \end{enumerate}
    \item 
    $\mathrm{T}$ is said to be  \textit{strongly }$\mathbf{C}^{\heart}$-\textit{coherent}  if $\mathbb{I}$ is strongly $\mathrm{T}$-coherent as an algebra.
\end{enumerate}

\end{defn}

\begin{defn}\label{defn:Tfinpres}
Let $A\in\mathbf{Alg_{D}}^{cn}$.
\begin{enumerate}
\item
An object $B\in{}_{A\big\backslash}\mathbf{DAlg}$ is said to be a $\mathrm{T}$-\textit{finitely presented }$A$-\textit{algebra} if it can be written as a finite colimit of algebras of the form $A\coprod^{\mathbb{L}}\mathrm{F}(\mathrm{Free_{T}}(\underline{\lambda}))$.
\item
An object $B\in{}_{A\big\backslash}\mathbf{DAlg}$ is said to be an \textit{almost} $\mathrm{T}$-\textit{finitely presented }$A$-\textit{algebra} if for each $n\ge0$, $\tau_{\le n}B$ is a $\mathrm{T}$-finitely presented $\tau_{\le n}$ $A$-algebra.
\end{enumerate}

The class of all $\mathrm{T}$-finitely presented algebras is denoted $\mathbf{DAlg}^{\mathrm{T}-fp}(\mathbf{C})$, and the class of all almost $\mathrm{T}$-finitely presented algebras is denoted $\mathbf{DAlg}^{\mathrm{T}-afp}(\mathbf{C})$.
\end{defn}

\begin{lem}\label{lem:pi0fp}
Let $A\in\mathbf{Alg_{D}}^{cn}(\mathbf{C})$ be such that $\pi_{0}(A)$ is strongly $\mathrm{T}$-coherent, and let $B$ be an almost $\mathrm{T}$-finitely presented $A$-algebra. Then $\pi_{0}(B)$ is $\mathrm{T}$-finitely presented as a $\pi_{0}(A)$-algebra and each $\pi_{n}(B)$ is finitely presented as a $\pi_{0}(B)$-module.
\end{lem}

\begin{proof}
Clearly we may assume that $B$ is in fact $\mathrm{T}$-finitely presented as an $\mathrm{A}$-algebra. Since $\pi_{0}$ commutes with colimits, $\pi_{0}(B)$ is a finite colimit of algebras of the form $\pi_{0}(A)\otimes^{\mathbb{L}}\mathrm{F}(\mathrm{Free}_{\mathrm{T}}(\underline{\lambda}))$ so it is by definition $\mathrm{T}$-finitely presented as a $\pi_{0}(A)$-algebra. Moreover since $\pi_{0}(B)$ is strongly $\mathrm{T}$-Coherent as an $A$-algebra, it is in fact a quotient of some $\pi(A)\hat{\otimes}\mathrm{F}(\mathrm{Free_{T}}(\underline{\gamma}))$ by a finitely generated ideal, and hence is a coherent algebra. Using the Tor spectral sequence, and the fact that finitely presented modules over coherent rings are closed under finite colimits, we find that each $\pi_{n}(B)$ is finitely presented as a $\pi_{0}(B)$-module.
\end{proof}

In fact, the Tor-spectral sequence implies the following, which applies when the the algebras in question are not necessarily strongly $\mathrm{T}$-coherent.

\begin{lem}\label{lem:pi0fp2}
Let $A,B,C\in\mathbf{DAlg}^{cn}(\mathbf{C})$ be such that for each $n$, $\pi_{n}(A)$, $\pi_{n}(B)$, and $\pi_{n}(C)$ are finitely generated modules over $\pi_{0}(A)$, $\pi_{0}(B)$, and $\pi_{0}(C)$ respectively. If $\pi_{0}(B\otimes^{\mathbb{L}}_{A}C)$ is coherent, then $\pi_{n}(B\otimes^{\mathbb{L}}_{A}C)$ is finitely generated over $\pi_{0}(B\otimes^{\mathbb{L}}_{A}C)$. 
\end{lem}

This motivates a more general class of algebras.

\begin{defn}\label{defn:T-cohalg}
Let $A\in\mathbf{DAlg}^{cn}(\mathbf{C})$. An $A$-algebra $B$ is said to be $\mathrm{T}$-\textit{coherent} as an $A$-algebra if
\begin{enumerate}
\item
each $\pi_{0}(B)$ is a finite colimit of algebras of the form $A\otimes\mathrm{F}(\mathrm{Free_{T}}(\underline{\lambda}))\big\slash I$ with $I$ a finitely generated ideal.
\item
each $\pi_{n}(B)$ is finitely presented as a $\pi_{0}(B)$-module.
\end{enumerate}
If $A=\mathbb{I}$ we just say that $B$ is a $\mathrm{T}$-coherent algebra. The class of all $\mathrm{T}$-coherent algebras is denoted $\mathbf{DAlg}^{\mathrm{T}-coh}(\mathbf{C})$, and of all connective $\mathrm{T}$-coherent algebras is denoted $\mathbf{DAlg}^{cn,\mathrm{T}-coh}(\mathbf{C})$.
\end{defn}

\begin{defn}\label{defn:discT-cohalg}
    Let $A\in\mathbf{DAlg}^{\heart}(\mathbf{C})$. $A$ is said to be a \textit{discrete }$\mathrm{T}$-\textit{finitely presented algebra} if it of the form $\pi_{0}(A)$ for $A$ a $\mathrm{T}$-finitely presented algebra in $\mathbf{DAlg}^{cn}(\mathbf{C})$. The class of all discrete $\mathrm{T}$-finitely presented algebras is denoted $\mathbf{DAlg}^{\heart,\mathrm{T}-coh}$.
\end{defn}

In particular for strongly $\mathrm{T}$-coherent algebras $A$, algebras which are almost $\mathrm{T}$-finitely presented as $A$-algebras are also $\mathrm{T}$-coherent as $A$-algebras.

\begin{example}
If $\mathrm{T}$ is strongly $\mathbf{C}^{\heart}$-coherent, then any $\mathrm{T}$-finitely presented algebra in $\mathbf{DAlg}^{cn}(\mathbf{C})$ is $\mathrm{T}$-coherent by Lemma \ref{lem:pi0fp}.
\end{example}

\begin{prop}\label{prop:Lfp}
Let $\mathrm{T}$ be strongly $\mathbf{C}^{\heart}$ coherent. Let $A\in\mathbf{DAlg}^{cn,\mathrm{T}-coh}$. Then each $\pi_{n}(\mathbb{L}_{A})$ is finitely presented as a $\pi_{0}(A)$-module.
\end{prop}

\begin{proof}
    This follows immediately from Proposition \ref{prop:altLcoh}.
\end{proof}

For some of the Lawvere theories we consider in analytic geometry, such as the affinoid and dagger affinoind ones, coherence seems to be the correct notion. However this class of algebras will be poorly behaved for the disc Lawvere theory controlling Stein geometry. Hence we introduce a more general concept.

\begin{defn}\label{defn:Tfinemb}
    Let $A\in\mathbf{DAlg}^{cn}(\mathbf{C})$. An $A$-algebra $B$ is said to be $\mathrm{T}$-\textit{finitely embeddable} as an $A$-algebra if there is a map $A\otimes^{\mathbb{L}}\mathrm{F}(\mathrm{Free_{T}}(\lambda_{1},\ldots,\lambda_{n})\rightarrow B$ which induces a surjection on $\pi_{0}$. If in addition $\mathrm{Ker}(\pi_{0}(A\otimes^{\mathbb{L}}\mathrm{F}(\mathrm{Free_{T}}(\lambda_{1},\ldots,\lambda_{n})))\rightarrow\pi_{0}(B))$ is a finitely generated $\pi_{0}(A\otimes^{\mathbb{L}}\mathrm{F}(\mathrm{Free_{T}}(\lambda_{1},\ldots,\lambda_{n})))$-module we say that $B$ is $\mathrm{T}$-\textit{globally finitely embeddable}
    
    If $A=\mathbb{I}$ we just say that $B$ is a $\mathrm{T}$-(globally) finitely embeddable algebra. The class of all $\mathrm{T}$-finitely embeddable algebras is denoted $\mathbf{DAlg}(\mathbf{C})^{cn,\mathrm{T}-f}$, and of all globally finitely embeddable algebras is denoted $\mathbf{DAlg}(\mathbf{C})^{cn,\mathrm{T}-gf}$.
\end{defn}

    We denote by $\mathbf{DAlg}^{\heart,\mathrm{T}-f}(\mathbf{C})$ (resp. $\mathbf{DAlg}^{\heart,\mathrm{T}-gf}(\mathbf{C})$) the class of algebras of the form $\pi_{0}(A)$ where $A\in\mathbf{DAlg}(\mathbf{C})^{cn,\mathrm{T}-f}$ (resp. $A\in\mathbf{DAlg}(\mathbf{C})^{cn,\mathrm{T}-gf}$). We also denote by

In particular any $\mathrm{T}$-coherent is finitely embeddable. In fact if $B$ is $\mathrm{T}$-coherent there is a map $\mathrm{F}(\mathrm{Free_{T}}(\lambda_{1},\ldots,\lambda_{n})) \rightarrow B$ which is a surjection on $\pi_{0}$, and such that the fibre is coherent as a $\mathrm{F}(\mathrm{Free_{T}}(\lambda_{1},\ldots,\lambda_{n}))$-module. In this sense $\mathrm{T}$-coherent algebras can be considered geometrically as the class of homotopically globally finitely embeddable objets.

It is worth keeping the following useful fact in mind.

\begin{prop}\label{prop:liftingmapsTpi0}
Let $A\in\mathbf{DAlg}^{cn}(\mathbf{C})$. Any map $\mathrm{F}(\mathrm{Free_{T}}(\lambda_{1},\ldots,\lambda_{n}))\rightarrow\pi_{0}(A)$ lifts to a map of $\mathrm{Sym}(\mathbb{I}^{\oplus n})$-algebras
$$\mathrm{F}(\mathrm{Free_{T}}(\lambda_{1},\ldots,\lambda_{n}))\rightarrow A$$
\end{prop}

\begin{proof}
Since $\mathrm{Sym}(\mathbb{I}^{\oplus n})$ is projective as an object of $\mathbf{DAlg}^{cn}(\mathbf{C})$, the composite 
$$\mathrm{Sym}(\mathbb{I}^{\oplus n})\rightarrow\mathrm{F}(\mathrm{Free_{T}}(\lambda_{1},\ldots,\lambda_{n}))\rightarrow\pi_{0}(A)$$
lifts to a map $\mathrm{Sym}(\mathbb{I}^{\oplus n})\rightarrow A$
The result now follows from Corollary \ref{cor:liftingetalemaps} and the fact that $\mathrm{Sym}(\mathbb{I}^{\oplus n})\rightarrow\mathrm{F}(\mathrm{Free_{T}}(\lambda_{1},\ldots,\lambda_{n}))$ is formally \'{e}tale.
\end{proof}

%

\begin{lem}[c.f. \cite{HA} Theorem 7.4.3.18]
Let $\pi_{0}(A)$ be strongly $\mathrm{T}$-coherent. If $B$ is almost $\mathrm{T}$-finitely presented as an $A$-algebra then each $\pi_{n}(\mathbb{L}_{B\big\slash A})$ is finitely presented as a $\pi_{0}(B)$-module. 
\end{lem}

\begin{proof}
Let $\mathbf{L}'$ be the full subcategory of ${}_{A\big\backslash}\mathbf{DAlg}^{cn}$ consisting of those $B$ such that $\pi_{n}(\mathbb{L}_{B\big\slash A})$ is finitely generated as a $B$-module and $\pi_{0}(B)$ is finitely $\mathrm{T}$-presented as a $\pi_{0}(A)$-algebra. We claim that $\mathbf{L}'$ is closed under finite colimits. It has an initial object, namely $A$. Since $\pi_{0}(A)$ is strongly $\mathrm{T}$-coherent, if $C$ is a such a finite colimit then $\pi_{0}(C)$ is of the form $\pi_{0}(C)\cong\pi_{0}(A)\otimes \mathrm{F}(\mathrm{Free_{T}}(\lambda_{1},\ldots,\lambda_{n}))\big\slash I$ with $I$ finitely generated. We have an exact sequence
$$I\big\slash I^{2}\rightarrow\pi_{0}(B)^{|n|}\rightarrow\Omega_{\pi_{0}(B)\big\slash\pi_{0}(A)}\rightarrow 0$$
Hence $\pi_{0}(\mathbb{L}_{B\big\slash A})\cong\Omega_{\pi_{0}(B)\big\slash\pi_{0}(A)}$ is finitely generated, and hence finitely presented, as a $\pi_{0}(B)$-module. The $\mathrm{Tor}$ spectral sequence implies that any pushout of objects in $\mathbf{L}'$ is still in $\mathbf{L}'$.

Now $\mathbb{L}_{A\coprod^{\mathbb{L}}\mathrm{F}(\mathrm{Free_{T}}(\underline{\lambda}))\big\slash A}\cong A\coprod^{\mathbb{L}}\mathrm{Free_{T}}(\underline{\lambda})^{\oplus |n|}$. Thus $\mathbf{L}'$ contains $\mathrm{T}$-finitely presented $A$-algebras. Since for any $A\rightarrow B$ and any $n<0$ we have that 
$$\tau_{\le n}\mathbb{L}_{B\big\slash A}\cong\tau_{\le n}\mathbb{L}_{A_{\le n}\big\slash B_{\le n}}$$
we also have that $\mathbf{L}'$ contains any almost $\mathrm{T}$-finitely presented $A$-algebra.
\end{proof}

\begin{lem}
Let $A\rightarrow B\rightarrow C$ be maps in $\mathbf{Alg_{D}}^{cn}(\mathbf{C})$ with $B$ being an almost $\mathrm{T}$-finitely presented $\mathrm{A}$-algebra. If $C$ is almost $\mathrm{T}$-finitely presented as a $B$-algebra then it is $\mathrm{T}$-almost finitely presented as an $A$-algebra.
\end{lem}

\begin{proof}
By truncating everything we may replace `almost $\mathrm{T}$-finitely presented' everywhere in the statement by `$\mathrm{T}$-finitely presented'. In this case we have $C$ is a finite colimit of algebras of the form $B\otimes^{\mathbb{L}}\mathrm{F}(\mathrm{Free_{\mathrm{T}}}(\underline{\lambda}))$, and $B$ is a finite colimit of algebras of the form $A\otimes^{\mathbb{L}}\mathrm{F}(\mathrm{Free_{\mathrm{T}}}(\underline{\gamma}))$. Since $\mathrm{F}(\mathrm{Free_{\mathrm{T}}}(\underline{\gamma}))\otimes^{\mathbb{L}}\mathrm{F}(\mathrm{Free_{\mathrm{T}}}(\underline{\lambda}))\cong\mathrm{F}(\mathrm{Free_{\mathrm{T}}}(\underline{\gamma},\underline{\lambda}))$ the claim is clear.
\end{proof}

\begin{lem}\label{lem:finpres}
If $B$ and $C$ are both (almost) $\mathrm{T}$-finitely presented $A$-algebras, Then $B\otimes^{\mathbb{L}}_{A}C$ is an (almost) $\mathrm{T}$-finitely presented $A$-algebra.
\end{lem}

\begin{proof}
The proof for $\mathrm{T}$-finitely presented $A$-algebras is clear.  Suppose $B$ and $C$ are almost $\mathrm{T}$-finitely presented. We need to show that $\tau_{\le n}(B\otimes_{A}^{\mathbb{L}}C)$ is  almost $\mathrm{T}$-finitely presented as a $\tau_{\le n}A$-algebra. Now $\tau_{\le n}(B)$ and $\tau_{\le n}(B)$ are almost $\mathrm{T}$-finitely presented. Thus $\tau_{\le n}(B)\otimes^{\mathbb{L}}_{\tau_{\le n}(A)}\tau_{\le n}(B)$ is $\mathrm{T}$-finitely presented as a $\tau_{\le n}(A)$-algebra, and $\tau_{\le n}(B\otimes_{A}^{\mathbb{L}}C)\cong\tau_{\le n}(\tau_{\le n}(B)\otimes^{\mathbb{L}}_{\tau_{\le n}(A)}\tau_{\le n}(B))$ is $\mathrm{T}$-finitely presented as a $\tau_{\le n}$-algebra, as required.
\end{proof}

\begin{cor}\label{cor:tensprodafp}
Let $A\in\mathbf{DAlg}^{cn}(\mathbf{C})$. Let $B$ be an $A$-algebra, and $C$ and $D$ be $B$-algebras. Suppose that $B$ is $\mathrm{T}$-finitely presented as an $A$-algebra, and that $C$ and $D$ are $\mathrm{T}$-finitely presented as $B$-algebras. Then $C\otimes^{\mathbb{L}}_{B}D$ is $\mathrm{T}$-finitely presented as a $\mathrm{A}$-algebra.
\end{cor}

\begin{prop}\label{prop:Tcohtens}
If $\mathrm{T}$ is strongly $\mathbf{C}^{\heart}$-coherent, and $A,B,C$ are all $\mathrm{T}$-coherent as $\mathbb{I}$-algebras then $B\otimes^{\mathbb{L}}_{A}C$ is $\mathrm{T}$-coherent as a $\mathbb{I}$-algebra.
\end{prop}

\begin{proof}
By the $\mathrm{T}$-coherence assumption on $A,B$, and $C$, $\pi_{0}(B)\otimes_{\pi_{0}(A)}\pi_{0}(C)$ is a finite colimit of algebras of the form $\mathrm{F}(\mathrm{Free_{T}}(\underline{\lambda}))\big\slash I$ with $I$ finitely generated. Moreover since $\mathrm{T}$ is strongly $\mathbf{C}^{\heart}$-coherent, $\pi_{0}(B)\otimes_{\pi_{0}(A)}\pi_{0}(C)$ is itself of the form $\mathrm{F}(\mathrm{Free_{T}}(\underline{\lambda}))\big\slash I$ with $I$ finitely generated. Thus $\pi_{0}(B)\otimes_{\pi_{0}(A)}\pi_{0}(C)$ is coherent. That $\pi_{n}(B\otimes^{\mathbb{L}}_{A}C)$ is finitely generated over $\pi_{0}(B)\otimes_{\pi_{0}(A)}\pi_{0}(C)$ now follows using an easy spectral sequence argument.
\end{proof}

\subsection{The Essential Image of $\mathbf{F}$}

Let $\mathrm{T}$ be of homotopy $\mathbf{C}^{\heart}$-polyomial type, and
suppose that the corresponding functor $\mathbf{F}:\mathbf{sAlg}_{\mathrm{T}}\rightarrow$ is fully faithful. Here we identify the essential image of $\mathbf{F}$

We first need a preparatory comment on Hochschild Homology.
$$\Sigma\mathrm{T}(\underline{\lambda})\cong\mathbb{I}\otimes^{\mathbb{L}}_{\mathrm{F}(\mathrm{Free}_{\mathrm{T}}(\underline{\lambda}))}\mathbb{I}\cong\mathbb{I}\otimes^{\mathbb{L}}_{\mathrm{Sym}(\mathbb{I}^{\oplus m})}\mathbb{I}$$
where we have used that $\mathrm{Sym}(\mathbb{I}^{\oplus n})\rightarrow\mathrm{F}(\mathrm{Free}_{\mathrm{T}}(\underline{\lambda}))$ is a homotopy epimorphism.
This in fact proves that the map $\Sigma^{q}\mathrm{Sym}(\mathbb{I}^{\oplus n})\rightarrow\Sigma^{q}\mathrm{T}(\underline{\lambda})$ is an equivalence for each $q>0$.

\begin{lem}\label{lem:Tconstructsequence}
Let $F:\mathrm{Alg_{T}^{ffg}}\rightarrow\mathbf{C}^{\heart}$ be a functor realising $\mathrm{T}$ as a $\Gamma$-sorted Fermat theory of homotopy $\mathbf{C}^{\heart}$-polynomial type.
\begin{enumerate}
    \item 
    Let $A\in\mathbf{DAlg}^{cn}(\mathbf{C})$ be in the essential image of $\mathbf{F}$. Then 
    \begin{enumerate}
        \item 
        $\pi_{0}(A)$ is in the essential image of $\mathrm{F}$
        \item 
        $\pi_{n}(A)$ is in the essential image of $j_{\pi_{0}(A)}{}_{\pi_{0}(A)}\mathrm{Mod}^{alg}(\mathbf{C}^{\heart})$.
    \end{enumerate}
    Morever, we have $\pi_{n}(\mathbf{F}(A))\cong j_{\pi_{0}(A)}(\pi_{n}(A))$. 
    \item 
Let $B$ be such that 
    \begin{enumerate}
        \item 
        $\pi_{0}(\mathbf{F}(B))\cong\mathrm{F}(\pi_{0}(B))$ is coherent in $\mathbf{C}^{\heart}$
        \item 
        each $\pi_{n}(B)$ is a finitely presented filtered $\pi_{0}(B)$ module.
    \end{enumerate}
Then $\mathbf{F}(B)$ is coherent.
\end{enumerate}

\end{lem}

\begin{proof}
\begin{enumerate}
    \item
     suppose that $A\cong\mathbf{F}(X)$ with $X\in\mathbf{sAlg}_{\mathrm{T}}$. We can write $X\cong\limind_{\alpha\in\mathrm{A}}X_{\alpha}$ where $X_{0}\cong\mathrm{T}(g)$ for some $g:Y\rightarrow\Lambda$, and for each successor ordinal $\alpha+1$ there is a pushout diagram
$$X_{\alpha+1}\cong X_{\alpha}\coprod^{\mathbb{L}}_{\Sigma^{n_{\alpha}}\mathrm{Free_{T}}(\gamma_{\alpha})}\mathrm{Free_{T}}(\gamma_{\alpha})(0)$$
or 
$$X_{\alpha+1}\cong X_{\alpha}\coprod^{\mathbb{L}}\Sigma^{n_{\alpha}}\mathrm{Free_{T}}(\gamma_{\alpha})$$
for some $n_{\alpha}>0$
Then
$$\mathbf{F}(X)\cong\limind_{\alpha\in\lambda}\mathbf{F}(X_{\alpha})$$
and 
$$\mathbf{F}(X_{\alpha+1})\cong\mathbf{F}(X_{\alpha})\otimes^{\mathbb{L}}_{\mathbf{F}(\Sigma^{n_{\alpha}}\mathrm{Free_{T}}(\gamma_{\alpha}))}\mathbf{F}(\mathrm{T}(0))\cong\mathbf{F}(X_{\alpha})\otimes^{\mathbb{L}}_{\Sigma^{n_{\alpha}}(\mathrm{Sym}(\mathbb{I}))}\mathbb{I}$$
or
$$\mathbf{F}(X_{\alpha+1})\cong\mathbf{F}(X_{\alpha})\otimes\Sigma^{n_{\alpha}}(\mathrm{Sym}(\mathbb{I}))$$
for some $n_{\alpha}>0$. In the former all that happens is that we take the quotient of $\pi_{n_{\alpha}}(X_{\alpha})$ by the image of a map $\pi_{0}(A)\rightarrow\pi_{n_{\alpha}}(X_{\alpha})$, and in the latter case we just add a copy of $\mathbb{I}$ to the homotopy group in degree $n_{\alpha}$. Since ${}_{\pi_{0}(A)}\mathrm{Mod}^{alg}(\mathbf{C}^{\heart})$ is closed under colimits we remain in this subcategory. 
\item
This follows from the first part, noting that the functor $j_{\pi_{0}(A)}$ is right exact,
\end{enumerate}
\end{proof}

We let $\mathbf{sAlg}_{\mathrm{T}}^{\pi_{0}-fp}$ denote the full subcagegory of $\mathbf{sAlg}_{\mathrm{T}}$ consisting of those algebras $A$ such that $\pi_{0}(A)$ is finitely presented as a discrete $\mathrm{T}$-algebra.

\begin{thm}\label{thm:essimageF}
Let $\mathrm{T}$ be a $1$-sorted Lawvere concretely of homotopy $\mathbf{C}^{\heart}$-polyomial type, and concretely of $\mathbf{C}^{\heart}$-polynomial type.
In particular the corresponding functors
$$\mathbf{F}:\mathbf{sAlg}_{\mathrm{T}}\rightarrow\mathbf{DAlg}^{cn}(\mathbf{C}^{\heart})$$
$$\mathrm{F}:\mathrm{Alg}_{\mathrm{T}}\rightarrow\mathrm{Comm}(\mathbf{C}^{\heart})$$ 
are fully faithful. Finally suppose that for any finitely presented $C\in\mathrm{Alg_{T}}$ the functor
$${}_{C}\mathrm{Mod}^{alg}\rightarrow{}_{\mathrm{F(C)}}\mathrm{Mod}(\mathbf{C}^{\heart})$$
is exact.
Then the essential image of $\mathbf{F}|_{\mathbf{sAlg}_{\mathrm{T}}^{\pi_{0}-fp}}$ consists of all objects $A$ of $\mathbf{DAlg}^{cn}(\mathbf{C})$ such that 
\begin{enumerate}
\item
$\pi_{0}(A)$ is in the essential image of $\mathrm{F}|_{\mathrm{Alg_{T}}^{fp}}$
\item
$\pi_{n}(A)\in{}_{\pi_{0}(A)}\mathrm{Mod}^{alg}(\mathbf{C}^{\heart})\subset {}_{\pi_{0}(A)}\mathrm{Mod}(\mathbf{C}^{\heart})$.
\end{enumerate}
If $\mathrm{T}$ is strongly $\mathbf{C}^{\heart}$-coherent then the essential image of $\mathbf{F}|_{\mathbf{sAlg}_{\mathrm{T}}^{coh}}$ is precisely $\mathbf{DAlg}^{cn,\mathrm{T}-coh}$. 
\end{thm}

\begin{proof}
These conditions are clearly necessary. Suppose the conditions hold. Write $\pi_{0}(A)$ as a sifted colimit $\colim_{\mathcal{I}}\mathrm{F}(\mathrm{Free_{T}}(g_{i}))$. We can write $\pi_{0}(A)\cong\mathrm{F}(\mathrm{Free_{T}}(S))\big\slash I$ for some set $S$, and some finitely generated ideal $I\in{}_{\mathrm{F}(\mathrm{Free_{T}}(S))}\mathrm{Mod^{alg}}(\mathbf{C}^{\heart})$. The map $\mathrm{F}(\mathrm{Free_{T}}(S))\rightarrow\pi_{0}(A)$ lifts to a map $\mathrm{Free_{T}}(S)\rightarrow A$.
Then we may write $A\cong\limind_{n}A_{\alpha}$ where
$$A_{0}\cong\mathrm{T}(g)$$
and $\mathrm{T}(g)\rightarrow\pi_{0}(A)$ is a surjection, and for each successor ordinal $\alpha+1$ there is a pushout diagram
$$A_{\alpha+1}\cong A_{\alpha}\otimes^{\mathbb{L}}_{\Sigma^{n_{\alpha}}\mathrm{Sym}(\mathbb{I})}\mathbb{I}\cong A_{\alpha}\hat{\otimes}_{\Sigma^{n_{\alpha}}\mathrm{T}(\underline{\lambda})}\mathrm{T}(0)$$
or 
$$A_{\alpha+1}\cong A_{\alpha}\hat{\otimes}\Sigma^{n_{\alpha}}\mathrm{F}(\mathrm{Free_{T}}(\underline{\lambda}))$$
Thus $A_{\alpha+1}$ is in the essential image of $\mathbf{F}$ whenever $A_{\alpha}$ is. Note that here we are using the fact that exactness of the functor ${}_{C}\mathrm{Mod}^{alg}\rightarrow{}_{\mathrm{F(C)}}\mathrm{Mod}(\mathbf{C}^{\heart})$ to ensure that taking finite limits and arbitrary colimits remains in the category of algebraic $\mathrm{F}(C)$-modules.. An easy transfinite induction then completes the proof. 
\end{proof}

\subsubsection{The Functor $\mathbf{F}$ and Coherent Algebras}

\begin{defn}
    Let $\mathrm{T}$ be a $1$-sorted Fermat theory. $A\in\mathbf{sAlg}_{\mathrm{T}}$ is said to be \textit{algebrically coherent} if 
    \begin{enumerate}
        \item 
        $\pi_{0}(A)$ is finitely presented in $\mathrm{Alg_{T}}$
        \item 
        $\pi_{n}(A)$ is finitely presented as a $\pi_{0}(A)$-module.
    \end{enumerate}
    Denote by $\mathbf{sAlg}_{\mathrm{T}}^{coh}$ the category of algebraically coherent $\mathrm{T}$-algebras.
\end{defn}

\begin{defn}
    A $1$-sorted Fermat theory $\mathrm{T}$ is said to be \textit{algebraically coherent} if any finitely presented object of $\mathrm{Alg_{T}}$ is coherent as a unital commutative algebra. 
\end{defn}

\begin{cor}\label{cor:cohessimageF}
Let $\mathrm{T}$ be a $1$-sorted Lawvere concretely of homotopy $\mathbf{C}^{\heart}$-polyomial type, and concretely of $\mathbf{C}^{\heart}$-polynomial type.
In particular the corresponding functors
$$\mathbf{F}:\mathbf{sAlg}_{\mathrm{T}}\rightarrow\mathbf{DAlg}^{cn}(\mathbf{C}^{\heart})$$
$$\mathrm{F}:\mathrm{Alg}_{\mathrm{T}}\rightarrow\mathrm{Comm}(\mathbf{C}^{\heart})$$ 
are fully faithful. Further suppose that the underlying algebra of each $\mathrm{T}(n)$ is coherent and that for any finitely presented $C\in\mathrm{Alg_{T}}$ the functor
$${}_{C}\mathrm{Mod}^{fp}\rightarrow{}_{\mathrm{F(C)}}\mathrm{Mod}(\mathbf{C}^{\heart})$$
is exact.
Then the essential image of $\mathbf{F}$ consists of all objects $A$ of $\mathbf{DAlg}^{cn}(\mathbf{C})$ such that 
\begin{enumerate}
\item
$\pi_{0}(A)$ is in the essential image of $\mathrm{F}$.
\item
$\pi_{n}(A)\in{}_{\pi_{0}(A)}\mathrm{Mod}^{alg}(\mathbf{C}^{\heart})\subset {}_{\pi_{0}(A)}\mathrm{Mod}(\mathbf{C}^{\heart})$ is finitely presented.
\end{enumerate}
\end{cor}

\subsection{$\mathrm{T}$-Smooth Maps}
In this subsection we define various classes of smooth maps `relative to $\mathrm{T}$'.

\begin{defn}\label{defn:Topensmooth}
Let $P$ be a projective object of $\mathbf{C}_{\ge0}$. A map $f:A\rightarrow B$ is said to be a $\mathrm{T}$-$P$-\textit{smooth} if 
\begin{enumerate}
\item
its is formally $P$-smooth.
\item
for any map $A\rightarrow C$ with $C$ $\mathrm{T}$-finitely embeddable, $\pi_{0}(C\otimes^{\mathbb{L}}_{A}B)$ is also $\mathrm{T}$-finitely embeddable
\end{enumerate}
$f$ is said to be $\mathrm{T}$-\textit{smooth} if it is $\mathrm{T}$-$P$-smooth for some $P$, $\mathrm{T}$-\textit{\'{e}tale morphism} if it is $\mathrm{T}$-$0$-smooth, and a $\mathrm{T}$-\textit{open immersion} if in addition it is a homotopy epimorphism.   Define the class of $\mathrm{T}$-smooth maps, $\mathrm{T}$-\'{e}tale maps, and $\mathrm{T}$-open immersions by $\mathbf{Sm}^{\mathrm{T}}$, $\textbf{\'{E}t}^{\mathrm{T}}$, and $\mathbf{open}^{\mathrm{T}}$-respectively
\end{defn}

\subsubsection{Regular Sequences}

We will define important classes of smooth maps using certain derived quotients. It will be useful to know when such quotients are discrete. Let $\mathrm{T}$ be a $\Gamma$-sorted Lawvere theory of homotopy $\mathbf{C}^{\heart}$-polynomial type, and let $A\in\mathbf{DAlg}^{cn}(\mathbf{C}^{\heart})$. 

\begin{defn}[\cite{bambozzi2020sheafyness} Definition 4.5]
A sequence $(a_{1},\ldots,a_{n}):\mathrm{Sym}(\mathbb{I}^{\oplus n})\rightarrow A$  is said to be $\mathrm{T}$-\textit{Koszul regular} if 
\begin{enumerate}
\item
it factors through a map $\mathrm{Sym}(\mathbb{I}^{\oplus n})\rightarrow\mathrm{F}(\mathrm{Free}_{\mathrm{T}}(\lambda_{1},\ldots,\lambda_{n}))\rightarrow A$ for some $(\lambda_{1},\ldots,\lambda_{n})\in\Lambda^{n}$
\item
the map $A\big\slash\big\slash(a_{1},\ldots,a_{n})\rightarrow A\big\slash(a_{1},\ldots,a_{n})$
\end{enumerate}
is an equivalence.
\end{defn}

\subsubsection{Localisations and $\mathrm{T}$-Open Maps}
Let $\mathrm{T}$ be a $\Gamma$-sorted Lawvere theory of homotopy $\mathbf{C}^{\heart}$-polynomial type.


A particularly important class of $\mathrm{T}$-open maps come from rational localisations. 
 For $\lambda\in\Gamma$ let $y_{\lambda}:\mathbb{I}\rightarrow\mathrm{F}(\mathrm{Free_{T}}(\lambda))$ denote the composition
$$\mathbb{I}\rightarrow\mathrm{Sym}(\mathbb{I})\rightarrow \mathrm{F}(\mathrm{Free_{T}}(\lambda))$$
where the first map is the inclusion of the $i$th coordinate.

\begin{defn}
Let $A\in\mathrm{Comm}(\mathbf{C}^{\heart})$. A collection of maps \[(a_{1},\ldots,a_{n}):\mathbb{I}^{\oplus n}\rightarrow A\] is said to \textit{generate the unit ideal in }$A$ if the induced map of $A$-modules
$$A^{\oplus n}\rightarrow A$$
is an epimorphism. 
\end{defn}

\begin{defn}[\cite{ben2021analytification}, Definition 6.9]\label{defn:Trat}
Let $A\in\mathbf{DAlg}^{cn}(\mathbf{C})$, let $f_{0},f_{1},\ldots,f_{n}:\mathbb{I}\rightarrow A$ be $\mathrm{T}$-extendable maps generating the unit ideal in $\pi_{0}(A)$, and let $(\lambda_{1},\ldots,\lambda_{n})$ be a sequence of elements of $\Lambda$. The $\mathrm{T}$-\textit{rational localisation of} $A$ \textit{at } $(f_{0},f_{1},\ldots,f_{n},(\lambda_{1},\ldots,\lambda_{n}))$ is the map
$$A\rightarrow A\otimes^{\mathbb{L}}\mathrm{F}(\mathrm{Free_{T}}(\lambda_{1},\ldots,\lambda_{n}))\big\slash\big\slash(f_{0}y_{\lambda_{1}}-f_{1},\ldots,f_{0}y_{\lambda_{n}}-f_{n})$$
The class of all $\mathrm{T}$-rational localisations is denoted $\mathbf{rat}^{\mathrm{T}}$. We also write
$$\overline{\mathbf{rat}}^{\mathrm{T}}\defeq\mathbf{rat}^{\mathrm{T}}\cap\mathbf{Sm}^{\mathrm{T}}$$
\end{defn}

\begin{defn}
    A map $A\rightarrow B$ in $\mathrm{Comm}(\mathbf{C}^{\heart})$ is said to be an \textit{discrete} $\mathrm{T}$-\textit{rational localisation} if it is of the form $A\rightarrow \pi_{0}(B)$ where $A\rightarrow B$ is a $\mathrm{T}$-rational localisation.
\end{defn}

\begin{example}
If $f_{0}=1$ then the localisation is called a $C$-\textit{Weierstrass localisation}.
\end{example}

\begin{defn}[\cite{ben2021analytification}, Definition 6.6]
\begin{enumerate}
\item
Let $A\in\mathbf{DAlg}^{cn}(\mathbf{C})$, let $g_{1},\ldots,g_{m}:\mathbb{I}\rightarrow A$ be $\mathrm{T}$-extendable maps
and let $(\lambda_{1},\ldots,\lambda_{n})$ be a sequence of elements of $\Lambda$. The \textit{pure} $\mathrm{T}$-\textit{Laurent localisation of }$A$ \textit {at }$g_{1},\ldots,g_{n},(\lambda_{1},\ldots,\lambda_{n})$ is the map
$$A\rightarrow A{\otimes}^{\mathbb{L}}\mathrm{F}(\mathrm{Free_{T}}(\lambda_{1},\ldots,\lambda_{n}))\big\slash\big\slash(g_{1}y_{\lambda_{1}}-1,\ldots,g_{1}y_{\lambda_{n}}-1)$$
\item
An $\mathrm{T}$-\textit{Laurent localisation} is a composition of a pure $\mathrm{T}$-Laurent localisation with a $\mathrm{T}$-Weierstrass localisation.
\end{enumerate}
\end{defn}

The following is proven in Lemma 6.10 \cite{ben2021analytification} for the case $\mathrm{C}=\mathrm{Ind(Ban_{R})}$ with $R$ a Banach ring. However the proof is essentially formal and works whenever $\mathrm{C}$ is a monoidal elementary quasi-abelian (or even exact) category.

\begin{lem}\label{lem:rathepi}
Any $\mathrm{T}$-derived rational localisation is a homotopy epimorphism. 
\end{lem}



\begin{defn}\label{defn:opencoherent}
\begin{enumerate}
    \item 
    An algebra $A\in\mathbf{DAlg}^{\heart}(\mathbf{C})$ is said to be $\mathrm{T}$-\textit{open coherent} if it is $\mathrm{T}$-localisation regular and finitely presented $A$-modules are transverse to $\mathrm{T}$-open maps $A\rightarrow B$.
    \item
    $\mathrm{T}$ is said to be $\mathbf{C}$-\textit{open coherent} if any discrete finitely $\mathrm{T}$-present algebra is $\mathrm{T}$-open coherent.
    
\end{enumerate}
\end{defn}

\begin{defn}\label{defn:loccoherent}
\begin{enumerate}
    \item 
    An algebra $A\in\mathbf{DAlg}^{\heart}$ is said to be $\mathrm{T}$-\textit{localisation regular} if whenever $f_{1},\ldots, f_{n}:\mathbb{I}\rightarrow A$ are $\mathrm{T}$-extendable maps generating the unit ideal of $\pi_{0}(A)$, then the sequence
    $$(f_{0}y_{\lambda_{1}}-f_{1},\ldots,f_{0}y_{\lambda_{n}}-f_{n})$$
    is a regular sequence in 
    $$A\otimes^{\mathbb{L}}\mathrm{F}(\mathrm{Free_{T}}(\lambda_{1},\ldots,\lambda_{n}))$$
    \item
    $\mathrm{T}$ is said to be $\mathbf{C}$-\textit{localisation regular} if any discrete $\mathrm{T}$-finitely embeddable algebra is localisation regular.
    \item 
    An algebra $A\in\mathbf{DAlg}^{\heart}(\mathbf{C})$ is said to be $\mathrm{T}$-\textit{localisation coherent} if finitely presented $A$-modules are transverse to  $\mathrm{T}$-rational localisations $A\rightarrow B$.
    \item
    $\mathrm{T}$ is said to be $\mathbf{C}$-\textit{localisation coherent} if any discrete $\mathrm{T}$-finitely embeddable algebra  is $\mathrm{T}$-localisation coherent.
\end{enumerate}
\end{defn}

\begin{cor}
    If $\mathrm{T}$ is $\mathbf{C}$-localisation regular and $\mathbf{C}$-localisation coherent then any $\mathrm{T}$-rational localisation of an algebra in $\mathbf{DAlg}(\mathbf{C})^{\mathrm{T}-coh}$ is a $\mathrm{T}$-open map.
\end{cor}

\subsubsection{Standard $\mathrm{T}$-\'{E}tale and $\mathrm{T}$-Smooth Maps}\label{subsubsection:standardT}

Let $A\in\mathbf{DAlg}^{cn}(\mathbf{C})$ and $f_{1},\ldots,f_{n}:\mathbb{I}\rightarrow\pi_{0}(A)\otimes^{\mathbb{L}}\mathrm{F}(\mathrm{Free_{T}}(\lambda_{1},\ldots,\lambda_{n}))$. Consider $A\hat{\otimes}^{\mathbb{L}}\mathrm{F}(\mathrm{Free_{T}}(\lambda_{1},\lambda_{n}))$. We get a map 
$$\mathrm{Sym}_{A}(A^{\otimes n})\cong A\otimes^{\mathbb{L}}\mathrm{Sym}(\mathbb{I}^{\oplus n})\rightarrow A\otimes^{\mathbb{L}}\mathrm{F}(\mathrm{Free_{T}}(\lambda_{1},\ldots,\lambda_{n}))$$
Passing to cotangent complexes we get a map

\begin{align*}
A\otimes\mathrm{\mathrm{F}\circ Free_{T}}(\lambda_{1},\ldots,\lambda_{n})^{\oplus n}& \cong A\otimes\mathrm{\mathrm{F}\circ Free_{T}}(\lambda_{1},\ldots,\lambda_{n})\otimes_{A}\mathbb{L}_{\mathrm{Sym}(A^{\oplus n})}\\
& \rightarrow \mathbb{L}_{A\otimes\mathrm{T}(\lambda_{1},\ldots,\lambda_{n})}\\
&\cong A\otimes\mathrm{\mathrm{F}\circ Free_{T}}(\lambda_{1},\ldots,\lambda_{n})^{\oplus n}
\end{align*}
This gives an element of 
$$\mathrm{Hom}(\mathbb{I},\pi_{0}(A\otimes\mathrm{\mathrm{F}\circ Free_{T}}(\lambda_{1},\ldots,\lambda_{n}))^{\oplus n^{2}}$$
which we call the \textit{Jacobian of }$(f_{1},\ldots,f_{n})$, and denote it $J(f_{1},\ldots,f_{n})$.

\begin{prop}
   Let  $f_{1},\ldots,f_{n}:\mathbb{I}\rightarrow\pi_{0}(A)\otimes^{\mathbb{L}}\mathrm{F}(\mathrm{Free_{T}}(\lambda_{1},\ldots,\lambda_{n}))$ be such that $J(f_{1},\ldots,f_{n})$ be such that $J(f_{1},\ldots,f_{n})$ is a unit. Then
   $$A\rightarrow A\otimes^{\mathbb{L}}\mathrm{F}(\mathrm{Free_{T}}(\lambda_{1},\ldots,\lambda_{n}))\big\slash\big\slash(f_{1},\ldots,f_{n})$$
   is formally \'{e}tale.
\end{prop}

\begin{proof}
     We have a homotopy pushout diagram
\begin{displaymath}
\xymatrix{
\mathrm{Sym}_{A}(A^{\oplus n})\ar[d]\ar[r] & A\otimes\mathrm{T}(\lambda_{1},\ldots,\lambda_{n})\ar[d]\\
A\ar[r] & B
}
\end{displaymath}
It therefore suffices to prove that the top horizontal map is formally \'{e}tale. This follows from the argument above.
\end{proof}

\begin{defn}\label{defn:Tstandardsmooth}
 A map $f:A\rightarrow B$ in $\mathbf{DAlg}^{cn}(\mathbf{C}^{\heart})$ is said to be 
 \begin{enumerate}
 \item 
 $\mathrm{T}$-\textit{standard \'{e}tale} if it is of the form 
   $$A\rightarrow A\otimes^{\mathbb{L}}\mathrm{F}(\mathrm{Free_{T}}(\lambda_{1},\ldots,\lambda_{n}))\big\slash\big\slash(f_{1},\ldots,f_{n})$$
   with $J(f_{1},\ldots,f_{n})$ invertible.
 \item
 $\mathrm{T}$-\textit{standard smooth of relative dimension }$k$ if it factorises as
 $$A\rightarrow A\otimes^{\mathbb{L}}\mathrm{F}(\mathrm{Free_{T}}(\lambda_{c+1},\ldots\lambda_{c+k}))\rightarrow A\otimes^{\mathbb{L}}\mathrm{F}(\mathrm{Free_{T}}(\lambda_{1},\ldots\lambda_{c+k}))\big\slash\big\slash (f_{1},\ldots,f_{c})$$
 with 
 $$A\otimes^{\mathbb{L}}\mathrm{F}(\mathrm{Free_{T}}(\lambda_{c+1},\ldots\lambda_{c+k}))\rightarrow A\otimes^{\mathbb{L}}\mathrm{F}(\mathrm{Free_{T}}(\lambda_{1},\ldots\lambda_{c+k}))\big\slash\big\slash (f_{1},\ldots,f_{c})$$
$\mathrm{T}$-standard \'{e}tale.
\item 
$\mathrm{T}$-\textit{open smooth of relative dimension }$k$ if it factorises as 
 $\mathrm{T}$-\textit{standard smooth of relative dimension }$k$ if it factorises as
 $$A\rightarrow A\otimes^{\mathbb{L}}\mathrm{F}(\mathrm{Free_{T}}(\lambda_{c+1},\ldots\lambda_{c+k}))\rightarrow A\otimes^{\mathbb{L}}\mathrm{F}(\mathrm{Free_{T}}(\lambda_{1},\ldots\lambda_{c+k}))\big\slash\big\slash (f_{1},\ldots,f_{c})$$
 with 
 $$A\otimes^{\mathbb{L}}\mathrm{F}(\mathrm{Free_{T}}(\lambda_{c+1},\ldots\lambda_{c+k}))\rightarrow A\otimes^{\mathbb{L}}\mathrm{F}(\mathrm{Free_{T}}(\lambda_{1},\ldots\lambda_{c+k}))\big\slash\big\slash (f_{1},\ldots,f_{c})$$
a homotopy monomorphism. 
\end{enumerate}
The class of all $\mathrm{T}$-standard smooth maps is denoted $\mathbf{sm}^{\mathrm{T}-std}$, of all open $\mathrm{T}$-standard smooth maps is denoted $\mathbf{sm}_{o}^{\mathrm{T}-std}$, and of all standard \'{e}tale maps is denoted $\textbf{\'{e}t}^{\mathrm{T}}$. We also write
$$\overline{\mathbf{sm}}^{\mathrm{T}-std}\defeq\mathbf{sm}^{\mathrm{T}-std}\cap\mathbf{Sm}^{\mathrm{T}}$$
$$\overline{\mathbf{sm}}_{o}^{\mathrm{T}-std}\defeq\mathbf{sm}_{o}^{\mathrm{T}-std}\cap\mathbf{Sm}^{\mathrm{T}}$$
$$\overline{\textbf{\'{e}t}}^{\mathrm{T}}\defeq\textbf{\'{e}t}^{\mathrm{T}}\cap\mathbf{Sm}^{\mathrm{T}}$$
for the class of \textit{strong }$\mathrm{T}$-\textit{standard smooth}, \textit{strong open}$\mathrm{T}$-\textit{standard smooth}, and \textit{strong }$\mathrm{T}$-\textit{\'{e}tale maps} respectively.
\end{defn}

Note that  a map is $\mathrm{T}$-standard \'{e}tale precisely if it is $\mathrm{T}$-standard smooth of relative dimension $0$.

\begin{rem}
Let $A\in\mathbf{DAlg}^{cn}(\mathbf{C})$. A $\mathrm{T}$-derived rational localisation of $A$ is $\mathrm{T}$-standard \'{e}tale.
\end{rem}

\begin{lem}
A $\mathrm{T}$-standard smooth map of relative dimension $k$ is formally smooth of relative dimension $k$.
\end{lem}

\begin{proof}
Let $A\rightarrow B$ be $\mathrm{T}$-standard smooth of relative dimension $k$. We can factor it as 
$$A\rightarrow A\otimes^{\mathbb{L}}\mathrm{F}(\mathrm{Free_{T}}(\lambda_{c+1},\ldots\lambda_{c+k}))\rightarrow A\otimes^{\mathbb{L}}\mathrm{F}(\mathrm{Free_{T}}(\lambda_{1},\ldots\lambda_{c+k}))\big\slash\big\slash (f_{1},\ldots,f_{c})$$
The first map is clearly formally smooth of relative dimension $k$. 
\end{proof}

\begin{lem}\label{lem:dersmoothrel}
Let $A\rightarrow B$ be a derived strong map in $\mathbf{Alg_{D}}^{cn}$ such that the map $\pi_{0}(A)\rightarrow\pi_{0}(B)\cong\pi(A)\otimes^{\mathbb{L}}\mathrm{\mathrm{F}\circ Free_{T}}(\lambda_{1},\ldots,\lambda_{c})\big\slash\big\slash(f_{1},\ldots,f_{c})$ is $\mathrm{T}$-standard smooth of relative dimension $k$. Then $A\rightarrow B$ is $\mathrm{T}$-standard smooth of relative dimension $k$. In fact we have $B\cong A\otimes^{\mathbb{L}}\mathrm{\mathrm{F}\circ Free_{T}}(\lambda_{1},\ldots,\lambda_{c})\big\slash\big\slash(f_{1},\ldots,f_{c})$
\end{lem}

\begin{proof}
The map $A\rightarrow B$ is formally smooth of relative dimension $k$ by Proposition \ref{prop:discretestrongsmooth} 
Define $\tilde{B}\defeq A\otimes^{\mathbb{L}}\mathrm{Free_{T}}(\lambda_{1},\ldots,\lambda_{c})\big\slash\big\slash(f_{1},\ldots,f_{c})$. Consider the factorisation

$$A\rightarrow A\otimes^{\mathbb{L}}\mathrm{Free_{T}}(\lambda_{c+1},\ldots,\lambda_{c+k})\rightarrow\tilde{B}$$
Now we the map $\mathrm{\mathrm{F}\circ Free_{T}(\lambda_{c+1},\ldots,\lambda_{c+k})}\rightarrow\pi_{0}(B)\cong\pi_{0}(A)\otimes\mathrm{Free_{T}}(\lambda_{1},\ldots,\lambda_{c})\big\slash(f_{1},\ldots,f_{c})$ lifts to a map $\mathrm{Free_{T}}(\lambda_{c+1},\rightarrow,\lambda_{c+k})\rightarrow A$ by Proposition \ref{prop:liftingmapsTpi0}. Thus we get a map
$$A\otimes^{\mathbb{L}}\mathrm{\mathrm{F}}(\mathrm{Free_{T}}(\lambda_{c+1},\ldots,\lambda_{c+k}))\rightarrow B$$
which is formally \'{e}tale. Now $A\otimes^{\mathbb{L}}\mathrm{\mathrm{F}\circ Free_{T}(\lambda_{c+1},\ldots,\lambda_{c+k})}\rightarrow\tilde{B}$ is also formally \'{e}tale.
By Corollary \ref{cor:liftingetalemaps} we get an isomorphism over $A\otimes^{\mathbb{L}}\mathrm{\mathrm{F}\circ Free_{T}(\lambda_{c+1},\ldots,\lambda_{c+k})}$, and hence over $A$, $\tilde{B}\cong B$.
\end{proof}

The same proof in fact shows the following.

\begin{lem}\label{lem:derratrel}
Let $A\rightarrow B$ be a derived strong map in $\mathbf{Alg_{D}}^{cn}$
\begin{enumerate}
    \item
    If  the map $\pi_{0}(A)\rightarrow\pi_{0}(B)$ is a $\mathrm{T}$-rational/ $\mathrm{T}$-Laurent/ $\mathrm{T}$-Weierstrass localisation, then $A\rightarrow B$ is a $\mathrm{T}$-rational/ $\mathrm{T}$-Laurent/ $\mathrm{T}$-Weierstrass localisation.
    \item 
      If  the map $\pi_{0}(A)\rightarrow\pi_{0}(B)$ is $\mathrm{T}$-open, then $A\rightarrow B$ is $\mathrm{T}$-open.
\end{enumerate}
\end{lem}

\begin{defn}
    A map $A\rightarrow B$ in $\mathrm{Comm}(\mathbf{C}^{\heart})$ is said to be an \textit{discrete} $\mathrm{T}$-\textit{standard \'{e}tale (resp. smooth map} if it is of the form $A\rightarrow \pi_{0}(B)$ where $A\rightarrow B$ is a $\mathrm{T}$-standard \'{e}tale (resp. smooth) map.
\end{defn}

\begin{defn}\label{defn:Tetaleregular}
\begin{enumerate}
    \item 
    An algebra $A\in\mathbf{DAlg}^{\heart}$ is said to be \textit{\'{e}tale regular} if whenever $f_{1},\ldots, f_{n}$ is a collection of elements of $A\otimes\mathrm{Free_{T}}(\underline{\lambda})$ whose Jacobian is a unit, then they define a regular sequence in  $A\otimes\mathrm{Free_{T}}(\underline{\lambda})$.
    \item
    $\mathrm{T}$ is said to be \textit{\'{e}tale regular} if any $\mathrm{T}$-finitely embeddable algebra is \'{e}tale regular.
    \item 
    An algebra $A\in\mathbf{DAlg}^{\heart}(\mathbf{C})$ is said to be \textit{\'{e}tale coherent} if finitely presented $A$-modules are transverse to standard $\mathrm{T}$-\'{e}tale maps $A\rightarrow B$.
    \item
    $\mathrm{T}$ is said to be \textit{\'{e}tale coherent} if any discrete finitely $\mathrm{T}$-embeddable algebra is \'{e}tale coherent.
\end{enumerate}
\end{defn}

\begin{rem}
    If $\mathrm{T}$ is $\mathrm{T}$-open regular/ $\mathrm{T}$-rational regular/ $\mathrm{T}$-\'{e}tale regular then $\overline{\mathbf{open}}^{\mathrm{T}}=\mathbf{open}^{\mathrm{T}}$/ $\overline{\mathbf{rat}}^{\mathrm{T}}=\mathbf{rat}^{\mathrm{T}}$/ $\overline{\textbf{\'{e}t}}^{\mathrm{T}}=\textbf{\'{e}t}^{\mathrm{T}}$.
\end{rem}

\begin{prop}\label{prop:etcohsmcoh}
 Let $A$ be a discretely finitely $\mathrm{T}$-presented algebra, let $\lambda_{1},\ldots,\lambda_{c+k}\in\Gamma$, and let $f_{1},\ldots,f_{c}$ be a sequence in $\mathrm{F}(\mathrm{Free_{T}}(\lambda_{1},\ldots,\lambda_{c+k}))$ such that the map
 $$A\rightarrow A\otimes^{\mathbb{L}}\mathrm{F}(\mathrm{Free_{T}}(\lambda_{1},\ldots,\lambda_{c+k}))\big\slash\big\slash (f_{1},\ldots,f_{c})$$
 is standard $\mathrm{T}$-smooth.
\begin{enumerate}
    \item 
    Suppose that $\mathrm{T}$ is $\mathbf{C}$-\'{e}tale regular. Then $\overline{\mathbf{sm}}^{\mathrm{T}-std}=\mathbf{sm}^{\mathrm{T}-std}$.
    \item 
    Suppose that $\mathrm{T}$ is $\mathbf{C}$-\'{e}tale regular and $\mathbf{C}$-\'{e}tale coherent. Let $M$ be a finitely presented $A$-module. Then the map
    $$M\otimes_{A}^{\mathbb{L}}A\otimes\mathrm{F}(\mathrm{Free_{T}}(\lambda_{1},\ldots,\lambda_{c+k}))\big\slash (f_{1},\ldots,f_{c})\rightarrow M\otimes_{A}\otimes\mathrm{F}(\mathrm{Free_{T}}(\lambda_{1},\ldots,\lambda_{c+k}))\big\slash (f_{1},\ldots,f_{c})$$
    is an equivalence.
\end{enumerate}
\end{prop}

\begin{proof}
    Consider the factorisation
     $$A\rightarrow A\otimes^{\mathbb{L}}\mathrm{F}(\mathrm{Free_{T}}(\lambda_{c+1},\ldots\lambda_{c+k}))\rightarrow A\otimes^{\mathbb{L}}\mathrm{F}(\mathrm{Free_{T}}(\lambda_{1},\ldots\lambda_{c+k}))\big\slash\big\slash (f_{1},\ldots,f_{c})$$
\begin{enumerate}
    \item 
         We have 
     $$A\otimes^{\mathbb{L}}\mathrm{F}(\mathrm{Free_{T}}(\lambda_{c+1},\ldots\lambda_{c+k}))\cong A\otimes\mathrm{F}(\mathrm{Free_{T}}(\lambda_{c+1},\ldots\lambda_{c+k}))$$
     since $\mathrm{F}(\mathrm{Free_{T}}(\lambda_{c+1},\ldots\lambda_{c+k}))$ is flat. By \'{e}tale regularity 
     $$A\otimes\mathrm{F}(\mathrm{Free_{T}}(\lambda_{1},\ldots\lambda_{c+k}))\big\slash\big\slash (f_{1},\ldots,f_{c})\cong A\otimes\mathrm{F}(\mathrm{Free_{T}}(\lambda_{1},\ldots\lambda_{c+k}))\big\slash (f_{1},\ldots,f_{c})$$
    \item
      $A\rightarrow\pi_{0}(A)\otimes\mathrm{F}(\mathrm{Free_{T}}(\lambda_{c+1},\ldots\lambda_{c+k}))$ is flat, and  by \'{e}tale coherence 
     $$A\otimes\mathrm{F}(\mathrm{Free_{T}}(\lambda_{c+1},\ldots\lambda_{c+k}))\rightarrow A\otimes\mathrm{F}(\mathrm{Free_{T}}(\lambda_{1},\ldots\lambda_{c+k}))\big\slash (f_{1},\ldots,f_{c})$$
     is transverse to finitely presented modules.
     Thus the composition 
     $$A\rightarrow A\otimes\mathrm{F}(\mathrm{Free_{T}}(\lambda_{c+1},\ldots\lambda_{c+k}))\rightarrow A\otimes\mathrm{F}(\mathrm{Free_{T}}(\lambda_{1},\ldots\lambda_{c+k}))\big\slash (f_{1},\ldots,f_{c})$$
     is transverse to finitely presented $A$-modules.
\end{enumerate}
     \end{proof}

\subsubsection{Coherence and $RR$-Quasicoherence for Modules}

The $\mathrm{T}$-rational localisations of interst to us in analytic geometry will in general not be flat. It is useful to consider those modules which are transverse to localisations.

\begin{defn}
 Fix a class of maps $\mathbf{P}$ in $\mathrm{Comm}(\mathbf{C}^{\heart})$. Let $A\in\mathrm{Comm}(\mathbf{C}^{\heart})$. An $A$-module $M$ is said to be $\mathbf{P}$-\textit{quasicoherent} if for any map $A\rightarrow B$ in $\mathbf{P}$, $B$ is transverse to $M$ as an $A$-module.
\end{defn}

The following is clear.

\begin{prop}
The class of $\mathbf{P}$-quasicoherent $A$-modules is closed under direct sums and filtered colimits.
\end{prop}

\begin{defn}\label{defn:Tqcoh}
An $A$-module $M$ is said to be 
\begin{enumerate}
    \item
    $(\mathrm{T},RR)$-\textit{quasicoherent} if it is $\mathbf{P}$-quasicoherent for $\mathbf{P}$ the class of discrete $\mathrm{T}$-rational localisations. 
    \item 
    $(\mathrm{T},\textrm{\'{e}t})$-\textit{quasicoherent} if it is $\mathbf{P}$-quasicoherent for $\mathbf{P}$ the class of discrete $\mathrm{T}$-standard \'{e}tale maps.. 
           \item
    $(\mathrm{T},\textrm{sm})$-\textit{quasicoherent} if it is $\mathbf{P}$-quasicoherent for $\mathbf{P}$ the class of discrete $\mathrm{T}$-standard smooth maps.
\end{enumerate}

\end{defn}

The terminology is inspired by Ramis-Ruget \cite{MR0352522}, and has also appeared in the context of analytic geometry in \cite{ben2020fr} Definition 7.4.

\begin{prop}\label{prop:colimloc}
Let $A=\limind_{\mathcal{I}} A_{i}$ with each $A_{i}\in\mathrm{Comm}(\mathbf{C}^{\heart})$ and let $M$ be an $A$-module which has a presentation $M\cong\colim M_{i}$ with each $M_{i}$ an $A_{i}$-module. If each $M_{i}$ is $(\mathrm{T},RR)$-/$(\mathrm{T},\textrm{\'{e}t})$-/$(\mathrm{T},\textrm{sm})$-quasicoherent as an $A_{i}$-module then $M$ is  $(\mathrm{T},RR)$-/$(\mathrm{T},\textrm{\'{e}t})$-/$(\mathrm{T},\textrm{sm})$-quasicoherent as an $A$-module.
\end{prop}

\begin{proof}
Let $A\rightarrow B$ be a $\mathrm{T}$-rational localisation. We may write $B\cong\colim B_{i}$, where $A\rightarrow B$ is the colimit of maps $A_{i}\rightarrow B_{i}$, each $A_{i}\rightarrow B_{i}$ being a $\mathrm{T}$-rational localisation. Then 
$$B\otimes_{A}^{\mathbb{L}}M\cong\limind_{\mathcal{I}}B_{i}\otimes_{A_{i}}^{\mathbb{L}}M_{i}\cong\limind_{\mathcal{I}}B_{i}\otimes_{A_{i}}M_{i}\cong B\otimes_{A}M$$
\end{proof}

By Proposition \ref{prop:colimloc} and Lemma \ref{lem:factscoherent} we have the following.

\begin{prop}
Let $A$ be coherent. Suppose that all finitely presented $A$-modules are $\mathbf{P}$-quasicoherent. Then any $M\in{}_{A}\mathrm{Mod}^{alg}$ is $\mathbf{P}$-quasicoherent.
\end{prop}

\subsubsection{Derived Coherent and Quasicoherent Algebras}

Let \[(\mathbf{C},\mathbf{C}_{\ge0},\mathbf{C}_{\le0},\mathbf{C}^{0})\] be a derived algebraic context, and $\mathrm{T}$ a $\Gamma$-sorted Lawvere of $\mathbf{C}^{\heart}$-polynomial type.

\begin{defn}
Let $A\in\mathbf{DAlg}(\mathbf{C})$. An $A$-module $M$ is said to be $\mathrm{T}$-\textit{quasi-coherent} if each $\pi_{n}(M)$ is a $(\mathrm{T},RR)$-quasicoherent $\pi_{0}(A)$-module.
\end{defn}

\begin{defn}
An algebra $A\in\mathbf{DAlg}(\mathbf{C})$ is said to be
 $\mathrm{T}$-\textit{quasi-coherent} if $\pi_{0}(A)$ is $\mathrm{T}$-finitely presented, and $A$ is $(\mathrm{T},RR)$-quasicoherent as an module over itself.
\end{defn}

\begin{example}
If $\mathrm{T}$ is localisation coherent then $\mathrm{T}$-coherent algebras are also $\mathrm{T}$-quasi-coherent.
\end{example}

\begin{lem}
    Let $\mathrm{T}$ be $\mathbf{C}$-open coherent/$\mathbf{C}$-rational coherent/ $\mathbf{C}$-\'{e}tale coherent, and let $A\in\mathbf{DAlg}^{cn,\mathrm{T}-coh}$. A map $A\rightarrow B$ is in $\overline{\mathbf{open}}_{\mathrm{T}}$/  $\overline{\mathbf{rat}}^{\mathrm{T}}$/ $\overline{\mathbf{sm}}_{\mathrm{T}}$ if and only if it is derived strong and $\pi_{0}(A)\rightarrow\pi_{0}(B)$ is in $\overline{\mathbf{open}}_{\mathrm{T}}$/ $\overline{\mathbf{rat}}^{\mathrm{T}}$/ $\overline{\mathbf{sm}}_{\mathrm{T}}$. In particular $B$ is also in $\mathbf{DAlg}^{cn,\mathrm{T}-coh}$.
\end{lem}

\begin{proof}
We prove the smooth case, the others being similar. 
    Suppose $f:A\rightarrow B$ is derived strong and $\pi_{0}(A)\rightarrow\pi_{0}(B)$ is in $\overline{\mathbf{sm}}_{\mathrm{T}}$, and is standard $\mathrm{T}$-smooth of relative dimension $k$. By Lemma \ref{lem:dersmoothrel} $f$ is standard $\mathrm{T}$-smooth of relative dimension $k$. Moreover if $A\rightarrow C$ is a map with $C$ discretely finitely $\mathrm{T}$-presented, then the map factors through $\pi_{0}(A)\rightarrow C$. We then have
    $$C\otimes_{A}^{\mathbb{L}}B\cong C\otimes^{\mathbb{L}}_{\pi_{0}(A)}\pi_{0}(A)\otimes_{A}^{\mathbb{L}}B\cong C\otimes_{\pi_{0}(A)}^{\mathbb{L}}\pi_{0}(B)$$
    which is discrete finitely $\mathrm{T}$-presented by assumption.

    Conversely suppose $f:A\rightarrow B$ is in $\overline{\mathbf{sm}}_{\mathrm{T}}$. $\pi_{0}(A)$ is discretely $\mathrm{T}$-finitely presented. Thus $\pi_{0}(A)\otimes_{A}^{\mathbb{L}}B$ is concentrated in degree $0$. Thus it must be equivalent to $\pi_{0}(B)$. Moreover $\pi_{0}(B)$ is itself discretely $\mathrm{T}$-finitely presented. Thus $\pi_{0}(A)\rightarrow\pi_{0}(B)$ is itself in $\overline{\mathbf{sm}}_{\mathrm{T}}$. Now by Proposition \ref{prop:etcohsmcoh}
     $\pi_{0}(A)\rightarrow\pi_{0}(B)$ is transverse to each $\pi_{n}(A)$. Thus $f$ is derived strong by Proposition \ref{prop:discretestrongsmooth}.

     Finally the above clearly shows that $B\in\mathbf{DAlg}^{cn,\mathrm{T}-coh}$.
     \end{proof}
 
In fact an idential proof gives the following.

\begin{lem}
    Let $f:A\rightarrow B$ be a map with $\pi_{n}(A)$ transverse to $\pi_{0}(B)$ over $\pi_{0}(A)$ for all $n$ and with $\pi_{0}(A)$ discretely finitely $\mathrm{T}$-presented. Then $f$ is in $\overline{\mathbf{open}}_{\mathrm{T}}$/ $\overline{\mathbf{rat}}^{\mathrm{T}}$/ $\overline{\textbf{\'{e}t}}_{\mathrm{T}}$/ $\overline{\mathbf{sm}}_{\mathrm{T}}$ if and only if $f$ is derived strong $\pi_{0}(f)$ is in $\overline{\mathbf{open}}_{\mathrm{T}}$/ $\overline{\mathbf{rat}}^{\mathrm{T}}$/ $\overline{\textbf{\'{e}t}}_{\mathrm{T}}$/ $\overline{\mathbf{sm}}_{\mathrm{T}}$.
\end{lem}

\chapter{Bornological Algebras}\label{BA}

\section{Bornological Modules}\label{sec:bornsp}

In this section we introduce our main derived algebraic context of interest, namely, the $(\infty,1)$-category of bornological spaces.

\subsection{Banach Rings and Modules}

We begin by defining semi-normed, normed, and Banach rings and their modules.

\begin{defn}
A \textit{semi-normed ring} is a unital commutative ring $R$, equipped with a function $|-|:R\rightarrow\mathbb{R}_{\ge0}$ such that there is $C_{R}\ge0$, and for all $s,r\in R$
\begin{enumerate}
\item
$|s+r|\le |s|+|r|$
\item
$|sr|\le |s||r|$
\item
$|0|=0$. 
\end{enumerate}
$R$ is said to be \textit{indiscrete} if $|1|=0$. $R$ is said to be a \textit{normed ring} if $|r|=0\Rightarrow r=0$. A normed ring is said to be a \textit{Banach ring} if it is complete as a metric space, with the metric induced by $|-|$. 
\end{defn}

\begin{defn}
A semi-normed ring $R$ is said to be \textit{non-Archimedean} if $|s+r|\le\mathrm{max}\{|r|,|s|\}$ fo r all $s,r\in R$.
\end{defn}

\begin{defn}
Let $R$ be a semi-normed ring. A \textit{semi-normed } $R$-\textit{module} is a $R$-module $M$, equipped with a function $||-||:M\rightarrow\mathbb{R}_{\ge0}$ such that there is $C_{M}\ge0$, such that for all $r \in R, m, n \in M$
\begin{enumerate}
\item
$||rm||\le C_{M}|r|||m||$
\item
$||m+n||\le  ||m||+||n||$
\item
$||0||=0$
\end{enumerate}
A map $f:M\rightarrow N$ of semi-normed $R$-modules is said to be \textit{bounded} if there is $C_{f}\ge 0$ such that for all $m\in M$,
$$||f(m)||\le C_{f}||m||$$
\end{defn}
Semi-normed $R$-modules arrange into a category $\mathrm{Norm}^{\frac{1}{2}}_{R}$.  

\begin{defn}
\begin{enumerate}
\item
A semi-normed $R$-module $M$ is said to be \textit{normed} if for $m\in M$, $||m||=0\Rightarrow m=0$. The full subcategory of $\mathrm{Norm}^{\frac{1}{2}}_{R}$ consisting of normed modules is denoted $\Norm_{R}$.
\item
A normed $R$-module $M$ is said to be \textit{Banach} if the metric on $M$, $d(m,n)\defeq||m-n||$ is complete. The full subcategory of $\Norm_{R}$ consisting of Banach modules is denoted $\Ban_{R}$. 
\item
If $R$ is non-Archimedean then a semi-normed/ normed/ Banach $R$-module $M$ is said to be \textit{non-Archimedean} if $||m+n||\le\mathrm{max}\{||m||,||n||\}$ for all $m,n\in M$. The full subcategories of $\mathrm{Norm}^{\frac{1}{2}}_{R},\mathrm{Norm}_{R},\mathrm{Ban}_{R}$ consisting of non-Archimedean modules will be denoted $\mathrm{Norm}^{\frac{1}{2},nA}_{R},\mathrm{Norm}^{nA}_{R},\mathrm{Ban}^{nA}_{R}$ respectively.
\end{enumerate}
\end{defn}
The inclusions
$$\mathrm{Norm}^{\frac{1}{2},nA}_{R}\rightarrow\mathrm{Norm}^{\frac{1}{2}}_{R}\;\;\; \Norm^{nA}_{R}\rightarrow\Norm_{R}\;\;\; \Ban^{nA}_{R}\rightarrow\Ban_{R}$$
commute with finite limits and colimits.

\begin{prop}
The categories $\mathrm{Norm}^{\frac{1}{2}}_{R},\Norm_{R},\Ban_{R}$ are all quasi-abelian. If $R$ is non-Archimedean then $\mathrm{Norm}^{\frac{1}{2},nA}_{R},\Norm^{nA}_{R},\Ban^{nA}_{R}$ are also quasi-abelian.
\end{prop}

\begin{proof}
For $\mathrm{Norm}^{\frac{1}{2}}_{R}$ and $\Norm_{R}$ this is shown for the case that $R$ is $\mathbb{C}$ in \cite{qacs} Section 3.2. For an arbitrary semi-normed ring the proof works identically. The proof for $\Ban_{R}$ is similar, and is mentioned for the case that $R=\mathbb{C}$ in e.g.  \cite{dcfapp} Section 3.3. The claim for the non-Archimedean versions follows immediately from the fact that the non-Archimedean subcategories are closed under finite limits and finite colimits.
\end{proof}

As for the case $R=\mathbb{C}$ considered in \cite{dcfapp}, for $R$ arbitrary the inclusion $\Norm_{R}\rightarrow\mathrm{Norm}^{\frac{1}{2}}_{R}$ has a left adjoint, 
$$\mathrm{Sep}:\mathrm{Norm}^{\frac{1}{2}}_{R}\rightarrow\Norm_{R}$$
which sends a semi-normed $R$-module $(M,\rho)$ to $M\big\slash\{m:\rho(m)=0\}$ equipped with the quotient norm. The functor $\mathrm{Sep}$ is exact.


The inclusion $\Ban_{R}\rightarrow\Norm_{R}$ also has a left adjoint
$$\mathrm{Cpl}:\Norm_{R}\rightarrow\Ban_{R}$$
It sends a normed $R$-module $M$ to equivalence classes of Cauchy sequences in $M$, with the limit norm.
The following is shown in Proposition 3.1.13 \cite{dcfapp}. We shall repeat the proof here for arbitrary semi-normed rings.

\begin{prop}\label{prop:cplexact}
The functor
$$\mathrm{Cpl}:\Norm_{R}\rightarrow\Ban_{R}$$
is exact.
\end{prop}

\begin{proof}
As a left adjoint it commutes with cokernels, i.e. it is strongly right exact in the sense of \cite{qacs} Definition 1.12, and therefore exact. It therefore suffices to show that it sends strict monomorphisms to strict monomorphisms. Let $f:E\rightarrow F$ be a strict monomorphism. Let $(x_{n})$ be a Cauchy sequence in $E$ such that $\rho_{F}(f(x_{n}))$ converges to $0$. Since $f$ is an admissible monomorphism $\rho_{E}(x_{n})$ converges to $0$. Thus the map $Cpl(E)\rightarrow Cpl(F)$ is injective. It remains to show that the map is admissible. Let $[(x_{n})]$ be an equivalence class in $Cpl(E)$, and let $C$ be such that for all $e\in E$, $\rho_{E}(e)\le C\rho_{F}(f(e))$. But $||[(x_{n})]||=\lim_{n\rightarrow\infty}||x_{n}||\le C \lim_{n\rightarrow\infty}||f(x_{n})||=||[f(x_{n})]||$. 
\end{proof}

\begin{rem}
For $R$ non-Archimedean the functors $\mathrm{Sep}$ and $\mathrm{Cpl}$ restrict to exact functors on the non-Archimedean categories, again giving left adjoints to the natural inclusions.
\end{rem}

Note: throughout we shall assume that our semi-normed/ normed/ Banach spaces live in some fixed Grothendieck universe - essentially their size will be bounded by some cardinal - which we shall suppress in notation. When we introduce ind-Banach spaces later, As in the condensed world, one can overcome this technicality by taking a colimit of categories over all cardinals.

\subsubsection{The Non-Expanding Categories and Projectives}

Let $\mathrm{Norm}^{\frac{1}{2},\le1}_{R}$ denote the wide subcategory of $\mathrm{Norm}^{\frac{1}{2}}_{R}$ consisting of maps $f:X\rightarrow Y$ such that $\rho_{Y}(f(x))\le\rho_{X}(x)$ for all $x\in X$. Similarly one defines the categories $\Norm^{\le1}_{R}$ and $\Ban^{\le1}_{R}$ These categories are studied in detail in Appendix A of \cite{koren}, and more details can be found there. In contrast to $\mathrm{Norm}^{\frac{1}{2}}_{R}$, $\Norm_{R}$, and $\Ban_{R}$, the non-expanding categories are complete and cocomplete. For example $\{X_{i}\}_{i\in\mathcal{I}}$ is a family of semi-normed $k$-modules, then $\coprod^{\le1}_{i\in\mathcal{I}}X_{i}$ is the semi-normed $k$-module whose underlying $k$-module is the direct sum $\bigoplus_{i\in\mathcal{I}}X_{i}$, with seminorm 
$$\rho_{\mathcal{I}}((x_{i})_{i\in\mathcal{I}})=\sum_{i\in\mathcal{I}}\rho_{X_{i}}(x_{i})$$
$\coprod^{\le1}_{i\in\mathcal{I}}X_{i}$  is called the \textit{non-expanding coproduct}. If each $X_{i}$ is normed, then the non-expanding coproduct is also normed, so this is also the coproduct in $\Norm^{\le1}_{R}$. In $\Ban^{\le1}_{R}$, the coproduct is given by
$$\{(x_{i})\in\prod_{i\in\mathcal{I}}X_{i}:\sum_{i\in\mathcal I}||x_{i}||<\infty\}$$
with norm $||(x_{i})||=\sum_{i\in\mathcal I}||x_{i}||<\infty$.

If $R$  is non-Archimedean then one can also define the subcategories  \[\mathrm{Norm}^{\frac{1}{2},nA,\le1}_{R},\Norm_R^{nA,\le1},\mathrm{Ban}_R^{nA,\le1}.\] The discussion above works mutatis mutandis for these categories with one caveat: the coproduct in this case is different. For $\{X_{i}\}_{i\in\mathcal{I}}$ is a family of non-Archimedean semi-normed $R$-modules, then $\coprod^{\le1}_{i\in\mathcal{I}}X_{i}$ is the semi-normed $R$-module whose underlying $R$-module is the direct sum $\bigoplus_{i\in\mathcal{I}}X_{i}$, with seminorm 
$$\rho_{\mathcal{I}}((x_{i})_{i\in\mathcal{I}})=\mathrm{sup}_{i\in\mathcal{I}}\rho_{X_{i}}(x_{i})$$
 In $\mathrm{Ban}_{R}^{nA,\le1}$ the coproduct is again given by completion.
 
 Coproducts are exact in the following precise sense. Let $\xNorm^{\le1}_{R}$ denote any of the categories $\mathrm{Norm}^{\frac{1}{2},\le1}_{R},\Norm_{R}^{\le1},\Ban_{R}^{\le1}$. We use similar notation $\xNorm^{nA,\le1}_{R}$ for the non-Archimedean versions. For $\mathcal{I}$ a set, let $\mathrm{Fun}^{ub}(\mathcal{I},\xNorm_{R}^{\le1})$ denote the wide subcategory of $\mathrm{Fun}^{ub}(\mathcal{I},\xNorm_{R}^{\le1})$ consisting of maps $(f_{i}):(X_{i})\rightarrow(Y_{i})$ such that there exists some $\epsilon>0$, such that for all $i\in\mathcal{I}$ and $x_{i}\in X_{i}$, $\rho_{Y_{i}}(f_{i}(x_{i}))\le \epsilon\rho_{X_{i}}(x_{i})$. The superscript $ub$ stands for \textit{uniformly bounded}. $\coprod_{i\in\mathcal{I}}^{\le1}:\mathrm{Fun}^{ub}(\mathcal{I},\xNorm_{k}^{\le1})\rightarrow\xNorm_{R}$ is a well-defined functor. Moreover the category $\mathrm{Fun}^{ub}(\mathcal{I},\xNorm_{R}^{\le1})$ is closed under computing kernels and cokernels in  $\mathrm{Fun}(\mathcal{I},\xNorm_{R}^{\le1})$. Finally, we say that a sequence
 $$0\rightarrow (X_{i})\rightarrow(Y_{i})\rightarrow (Z_{i})\rightarrow 0$$
 in $\mathrm{Fun}^{ub}(\mathcal{I},\xNorm_{R}^{\le 1})$ is exact if each sequence
 $$0\rightarrow X_{i}\rightarrow Y_{i}\rightarrow Z_{i}\rightarrow 0$$
 is exact in $\xNorm_{R}$.

\begin{prop}\label{prop:nonexpandingexact}
Here we do not distinguish between Archimedean and non-Archimedean categories.
\begin{enumerate}
\item
Let $\mathcal{I}$ be a set. The functor $\coprod_{i\in\mathcal{I}}^{\le1}:\mathrm{Fun}^{ub}(\mathcal{I},\mathsf{xNorm}_{R}^{\le1})\rightarrow\mathsf{xNorm}_{R}$ sends exact sequences to exact sequences, and commutes with kernels and cokernels. 
\item
Let $\mathcal{I}$ be a filtered category. The functor $\colim_{i\in\mathcal{I}}^{\le1}:\mathrm{Fun}^{ub}(\mathcal{I},\mathsf{xNorm}_{R}^{\le1})\rightarrow\mathsf{xNorm}_{R}$ sends exact sequences to exact sequences. It commutes with kernels and cokernels in the cases $\mathrm{Norm}^{\frac{1}{2},\le1}_{R},\Norm_{R}^{\le1}$. 
\end{enumerate}
\end{prop}

\begin{proof}
\item
In this case commuting with kernels and cokernels implies exactness. The fact that the non-expanding direct sum commutes with kernels and cokernels is clear.
\item
$\mathrm{Norm}^{\frac{1}{2},\le1}_{R},\Norm_{R}^{\le1}$ have a set of tiny generators, namely $\{R_{\delta}:\delta\in\mathbb{R}_{>0}\}$. Thus in this case filtered colimits commute with kernels and cokernels, and so exactness of filtered colimits is clear here. Exactness in $\Ban_{R}^{\le1}$ then follows from exactness in $\mathrm{Norm}^{\frac{1}{2},\le1}_{R},\Norm_{R}^{\le1}$ and Proposition \ref{prop:cplexact}.
\end{proof}

For $M\in\mathrm{Norm}^{\frac{1}{2}}_{R}$ and $\epsilon>0$, denote by $M_{\epsilon}$ the semi-normed $R$-module with the same underlying $R$-module as $M$, but with $\rho_{M_{\epsilon}}\defeq\epsilon\rho_{M}$. Note that $M\cong M_{\epsilon}$ in $\mathrm{Norm}^{\frac{1}{2}}_{R}$ but not in $\mathrm{Norm}^{\frac{1}{2},\le1}_{R}$. Finally, we denote by $M^{\rightarrow 0}$ the module $\coprod^{\le1}_{\epsilon\in\mathbb{R}_{>0}}M_{\epsilon}$.

\begin{prop}\label{prop:non-expandingproj}
For each $i\in\mathcal{I}$ let $\{\epsilon_{i}\}$ be a collection of non-negative real numbers. Then $\coprod^{\le1}_{i\in\mathcal{I}}R_{\epsilon_{i}}$ is projective.
\end{prop}

\begin{proof}
Let $f:X\rightarrow\coprod^{\le1}_{i\in\mathcal{I}}R_{\epsilon_{i}}$ be an admissible epimorphism. For each $i$ let $f_{i}:X\rightarrow R_{\epsilon_{i}}$ denote the composite of $f$ with the projection to $R_{\epsilon_{i}}$. This is clearly an admissible epimorphism. For any $\delta>0$ there is $x_{\delta}^{i}\in X$ such that $f_{i}(x^{i}_{\delta})=1$ and $\rho_{X}(x^{i}_{\delta})<(\delta+1)\epsilon_{i}$. Let $g_{i}:R_{\epsilon_{i}}\rightarrow X$ be the map sending $\lambda$ to $\lambda x^{i}_{\delta}$. Then $\rho_{X}(g_{i}(\lambda))\le (C+1)\epsilon_{i}|\lambda|$. Thus the collection of maps $g_{i}$ is uniformly bounded by $(C+1)$, so defines a bounded map $g:\coprod^{\le1}_{i\in\mathcal{I}}R_{\epsilon_{i}}\rightarrow X$ which is a right inverse for $f$.
\end{proof}

The proof of the next claim is similar to Proposition 3.2.11 of \cite{qacs} which deals with $k=\mathbb{C}$, but one has to be a bit more careful for an arbitrary semi-normed ring. 

\begin{prop}
The categories $\mathrm{Norm}^{\frac{1}{2}}_{R},\Norm_{R},$ and $\Ban_{R}$ have enough functorial projectives. If $R$ is non-Archimedean then $\mathrm{Norm}^{\frac{1}{2},nA}_{R},\Norm^{nA}_{R},$ and $\Ban^{nA}_{R}$ have enough functorial projectives.
\end{prop}

\begin{proof}
 First we consider the semi-normed case. Let $M$ be a semi-normed $k$-module. Consider the $R$-module $P(M)$ defined by
$$p:P(M)\defeq\Bigr(\coprod^{\le1}_{m:\rho_{M}(m)\neq0}R_{\rho_{M}(m)}\Bigr)\oplus\coprod^{\le1}_{m:\rho_{M}(m)=0}R^{\rightarrow 0}$$
By Proposition \ref{prop:non-expandingproj} this is a projective object. We claim that $P(M)\rightarrow M$ is an admissible epimorphism. Let $m\in M$ be such that $\rho_{M}(m)\neq0$. Consider the element $\delta_{m}$ of $P(M)$, which is the tuple with value $1$ corresponding to the copy of $R_{\rho_{M}(m)}$ indexed by $m$, and value $0$ in all other entries. Then $p(\delta_{m})=m$, and $\rho_{M}(m)=\rho_{P(M)}(\delta_{m})$. Now let $m\in M$ be such that $\rho_{M}(m)=0$. Consider the copy of $R^{\rightarrow0}$ indexed by $m$, and for each $\epsilon>0$ the summand $R_{\epsilon}$. Denote by $\delta_{m,\epsilon}$ the element of $P(M)$ whose value is $1$ on this copy of $R_{\epsilon}$ and zero elsewhere. Then $p(\delta_{m,\epsilon})=m$, and $\rho_{P(M)}(\delta_{m,\epsilon})=\epsilon$. This shows that the map $p$ is an admissible epimorphism. 

Now consider $\Norm_{R}$ and $\Ban_{R}$. In both cases the inclusion functors to $\mathrm{Norm}^{\frac{1}{2}}_{R}$ are exact. Thus the left adjoints $\mathrm{Sep}$ and $\mathrm{Cpl}=$ send projectives to projectives. Moreover if $M$ is a normed (resp. Banach) $R$-module, then $Sep(P(R_{N}^{N^{\frac{1}{2}}}(M)))\rightarrow M$ (resp. $Cpl(P(R_{N}^{N^{\frac{1}{2}}}(M)))\rightarrow M$) is an admissible epimorphism with $Sep(P(R_{N}^{N^{\frac{1}{2}}}(M)))$ (resp. $Cpl(P(R_{N}^{N^{\frac{1}{2}}}(M)))$) projective.
\end{proof}

\subsubsection{The Closed Monoidal Structure}

As explained in e.g. \cite{MR1888309}, the categories $\mathrm{Norm}^{\frac{1}{2}}_{R},\Norm_{R},\Ban_{R}$ have many monoidal structures. We are interested in the projective one. 

\begin{defn}
Let $X$ and $Y$ be semi-normed $R$-modules. Define the \textit{projective tensor semi-norm} $||-||_{\pi}$ on $X\otimes_{R}Y$ as follows.
$$||v||_{\pi}=\textrm{inf}\{\sum_{i=1}^{n}||x_{i}||_{X}||y_{i}||_{Y}:n\in\mathbb{N}, v=\sum_{i=1}^{n}x_{i}\otimes y_{i}\}$$
The \textit{projective tensor product} of $X$ and $Y$, denoted $X\otimes_{R,\pi}Y$ is the algebraic tensor product $X\otimes_{R}Y$ equipped with the projective tensor semi-norm.
\end{defn}

\begin{defn}
Let $R$ be non-Archimedean and let $X$ and $Y$ be non-Archimedean  semi-normed $R$-modules. Define the \textit{projective tensor semi-norm}, $||-||_{\pi}$ on $X\otimes_{R}Y$ as follows.
$$||v||_{\pi}=\textrm{inf}\{\mathrm{max}\{||x_{i}||_{X}||y_{i}||_{Y}\}_{i=1}^{n}, v=\sum_{i=1}^{n}x_{i}\otimes y_{i}\}$$
The \textit{projective tensor product} of $X$ and $Y$, denoted $X\otimes^{nA}_{R,\pi}Y$ is the algebraic tensor product $X\otimes_{R}Y$ equipped with the projective tensor semi-norm.
\end{defn}

The projective tensor product is closed.

\begin{defn}
Let $(X,\rho_{X})$ and $(Y,\rho_{Y})$ be semi-normed $R$-modules. Define $\underline{\Hom}(X,Y)$ to be the semi-normed $R$-module, whose underlying $R$-module is the space of bounded $R$-linear morphisms from $X$ to $Y$. The semi-norm is 
$$\rho_{\underline{\Hom}(X,Y)}(T)=\inf\{C:\rho_{Y}(T(x))\le C\rho_{X}(x)\;\forall x\in X\}$$
\end{defn}

\begin{prop}
Let $X$ be any semi-normed $R$-module.
\begin{enumerate}
\item
If $Y$ is a normed $R$-module then $\underline{\Hom}(X,Y)$ is normed.
\item
If $Y$ is a Banach $R$-module then $\underline{\Hom}(X,Y)$ is Banach.
\end{enumerate}
\end{prop}

\begin{proof}
\begin{enumerate}
\item
Let $T:X\rightarrow Y$ be a bounded map, and suppose that $\rho_{\underline{\Hom}(X,Y)}(T)=0$. Then for any $x\in X$, $\rho_{Y}(T(x))=0$. Since $Y$ is normed, $T(x)=0\;\;\forall x$. Thus $T=0$.
\item
Suppose that $Y$ is Banach. Let $(T_{n})$ be a Cauchy sequence in $\underline{\Hom}(X,Y)$. For any $x\in X$ $T_{n}(x)$ is a Cauchy sequence in $Y$. It therefore converges to some $T(x)$ in $Y$. It is easy to verify that $x\mapsto T(x)$ is linear and bounded, and that the sequence $(T_{n})$ converges to $\mathrm{T}$. 
\end{enumerate}
\end{proof}

\begin{cor}
\begin{enumerate}
\item
$\Norm_{R}$ inherits a closed monoidal structure from $\mathrm{Norm}^{\frac{1}{2}}_{R}$ with $X\overline{\otimes}_{R,\pi}T\defeq \mathrm{Sep}(X\otimes_{R,\pi}Y)$.
\item
$\Ban_{R}$ inherits a closed monoidal structure from $\mathrm{Norm}^{\frac{1}{2}}_{R}$ with $X\otimes_{R,\pi}Y\defeq \mathrm{Cpl}(X\otimes_{R,\pi}Y)$.
\end{enumerate}
\end{cor}

It is easily verified that the closed monoidal structure $\otimes,\underline{\Hom}$ on $\mathrm{Norm}^{\frac{1}{2}}_{k}$ restricts to one on $\xNorm_{R}^{\le1}$. This allows us to prove that projectives in $\xNorm_{R}^{\le1}$ are strongly flat, and that the tensor product of two projectives is projective.

\begin{prop}\label{prop:nonexpandingflat}
Let $\{F_{i}\}_{i\in\mathcal{I}}$ be a collection of (strongly) flat objects in $\mathsf{xNorm}_{R}$. Then $\coprod^{\le1}_{i\in\mathcal{I}}F_{i}$ is (strongly) flat. 
\end{prop}

\begin{proof}
There is a well-defined (strongly) exact functor $\xNorm_{R}\rightarrow\mathrm{Fun}^{ub}(\mathcal{I},\xNorm_{R})$ sending $X$ to $(F_{i}\otimes_{R,\pi}X)$. The functor $\xNorm_{R}\rightarrow\mathrm{Fun}^{ub}(\mathcal{I},\xNorm_{R})\rightarrow\xNorm_{R}$ is strongly exact, so the composite functor
$$X\mapsto\coprod^{\le1}_{i\in\mathcal{I}}(F_{i}\otimes_{R,\pi}X)$$
is (strongly) exact.  Moreover since $\otimes$ commutes with non-expanding colimits, this functor is isomorphic to $X\mapsto (\coprod^{\le1}_{i\in\mathcal{I}}F_{i})\otimes_{R,\pi}X$. 
\end{proof}

\begin{prop}
Let $P$ and $P'$ be projective objects in $\mathsf{xNorm}_{R}$. Then $P\otimes_{R,\pi}P'$ is projective. If $P$ is projective then $\mathrm{Sym}^{n}(P)$ is projective. Moreover projectives are strongly flat in $\mathsf{xNorm}_{R}$. 
\end{prop}

\begin{proof}
 Let $\tilde{R}$ denote $R$ in $\mathrm{Norm}^{\frac{1}{2}}_{R}$, $Sep(R)$ in $\Norm_{R}$, and $Cpl(R)$ in $\Ban_{R}$. $P$ is a retract of a projective of the form $\coprod^{\le1}_{i\in\mathcal{I}}\tilde{R}_{\epsilon_{i}}$ and $P'$ is a retract of $\coprod^{\le1}_{j\in\mathcal{J}}\tilde{R}_{\epsilon_{j}}$ with all $\epsilon_{i}>0,\epsilon_{j}>0$ . Thus $P\otimes P'$ is a retract of $\coprod^{\le1}_{i\in\mathcal{I}}\tilde{R}_{\epsilon_{i}}\otimes \coprod^{\le1}_{j\in\mathcal{J}}\tilde{R}_{\epsilon_{j}}$. But since $\otimes$ is a closed monoidal structure on $\xNorm_{R}^{\le1}$ it commutes with coproducts, so 
$$\coprod^{\le1}_{i\in\mathcal{I}}\tilde{R}_{\epsilon_{i}}\otimes \coprod^{\le1}_{j\in\mathcal{J}}\tilde{R}_{\epsilon_{j}}\cong\coprod^{\le1}_{(i,j)\in\mathcal{I}\times\mathcal{J}}\tilde{R}_{\epsilon_{i}}\otimes\tilde{R}_{\epsilon_{j}}\cong\coprod^{\le1}_{(i,j)\in\mathcal{I}\times\mathcal{J}}\tilde{R}_{\epsilon_{i}\epsilon_{j}}$$
and $\epsilon_{i}\epsilon_{j}>0$. 

If $P$ is a retract of $\coprod^{\le1}_{i\in\mathcal{I}}\tilde{R}_{\epsilon_{i}}$ then $\mathrm{Sym}^{n}(P)$ is a retract of $\mathrm{Sym}^{n}(\coprod^{\le1}_{i\in\mathcal{I}}\tilde{R}_{\epsilon_{i}})$. But this is just $\coprod^{\le1}_{[i_{1},\ldots,i_{n}]\in\mathcal{I}^{n}\big\slash\Sigma_{n}}R_{\epsilon_{i_{1}}\ldots\epsilon_{i_{n}}}$.

Each $R_{\epsilon_{i}}$ is strongly flat, so by Proposition \ref{prop:nonexpandingflat} $P=\coprod^{\le1}_{i\in\mathcal{I}}\tilde{R}_{\epsilon_{i}}$ is strongly flat. A retract of a strongly flat object is strongly flat, so any projective is strongly flat.
\end{proof}

\subsection{Ind-Banach and Bornological Modules}

 $\xNorm_{R}$ is not a particularly well-behaved category. In particular it only has finite limits and colimits in general. 
In this section we study bornological modules over an arbitrary semi-normed ring $k$, defined in \cite{bambozzi} Section 1.4. 
Fix a a semi-normed ring $R$, and consider the category
$$\mathrm{Ind}(\xNorm_{R})$$
By the results of the previous section we have the following

\begin{cor}\label{cor:IndBangood}
$\mathrm{Ind}(\xNorm_{R})$ is a monoidal elemenetary quasi-abelian category with symmetric projectives when equipped with the projective tensor product.
\end{cor}

In particular we get a derived algebraic context
$$(\mathbf{Ch}(\mathrm{Ind}(\xNorm_{R})),\mathbf{Ch}_{\ge0}(\mathrm{Ind}(\xNorm_{R})),\mathbf{Ch}_{\le0}(\mathrm{Ind}(\xNorm_{R})),\mathcal{P}^{0})$$
where $\mathcal{P}^{0}$ is the collection of projective objects of $\xNorm_{k}$. 

\begin{thm}\label{thm:bornillusie}
    The derived algebraic context 
    $$(\mathbf{Ch}(\mathrm{Ind}(\xNorm_{R})),\mathbf{Ch}_{\ge0}(\mathrm{Ind}(\xNorm_{R})),\mathbf{Ch}_{\le0}(\mathrm{Ind}(\xNorm_{R})),\mathcal{P}^{0})$$
    is Illusie.
\end{thm}

\begin{proof}
  First we prove that exterior powers and symmetric of projectives are projective. First assume $P=\coprod^{\le1}_{1\le j\le n}R_{\delta_{j}}$. Picking the standard basis of $P$ we find that $\Lambda^{n}(P)\cong\coprod^{\le1}_{1\le i_{1}<i_{2}<\ldots<i_{k}\le n}R_{\delta_{i_{1}}\ldots\delta_{i_{k}}}$ and is therefore projective. Now consider a projective of the form
 
 $$Q\cong\coprod_{i\in\mathcal{I}}^{\le 1}P_{i}$$
 with each $P_{i}$ of the form $R_{\delta_{i}}$. Then
 $$Q^{\otimes n}\cong\coprod^{\le1}_{[(i_{1},\ldots,i_{n})]\in(\mathcal{I}^{n}\big\slash\Sigma_{n})}\coprod^{\le1}_{\sigma\in\Sigma_{n}}P_{i_{\sigma_{1}}}\otimes\ldots\otimes P_{i_{\sigma(n)}}$$
 $\Lambda^{n}(Q)$ is the coequaliser of the maps 
 $$\coprod^{\le1}_{[i]=[(i_{1},\ldots,i_{n})]\in(\mathcal{I}^{n}\big\slash\Sigma_{n})}\;\coprod^{\le1}_{[\sigma]\in\Sigma_{n}\big\slash\mathrm{Stab}(i)}P_{i_{\sigma_{1}}}\otimes\ldots\otimes P_{i_{\sigma(n)}}\rightarrow \coprod^{\le1}_{[(i_{1},\ldots,i_{n})]\in(\mathcal{I}^{n}\big\slash\Sigma_{n})}\;\coprod^{\le1}_{[\sigma]\in\Sigma_{n}\big\slash\mathrm{Stab}(i)}P_{i_{\sigma_{1}}}\otimes\ldots\otimes P_{i_{\sigma(n)}}$$
 which is the same as the coproduct of the coequaliser of the maps
 $$\coprod^{\le1}_{[\sigma]\in\Sigma_{n}\big\slash\mathrm{Stab}(i)}P_{i_{\sigma_{1}}}\otimes\ldots\otimes P_{i_{\sigma(n)}}\rightarrow\coprod^{\le1}_{[\sigma]\in\Sigma_{n}\big\slash\mathrm{Stab}(i)}P_{i_{\sigma_{1}}}\otimes\ldots\otimes P_{i_{\sigma(n)}}$$
 This is precisely $\Lambda^{n}(\coprod^{\le1}_{\sigma\in\Sigma_{n}}P_{i_{\sigma_{1}}}\otimes\ldots\otimes P_{i_{\sigma(n)}})$ which as we have seen is projective.   Finally let $Q\cong\coprod_{j\in\mathcal{J}}P_{j}$ with each $P_{j}$ of the form $\coprod_{i_{j}\in\mathcal{I_{j}}}^{\le 1}R_{\delta_{i_{j}}}$. A similar argument to the above without the non-expanding coproduct shows again that $\Lambda^{n}(P)$ is projective for any projective $P$ in $\mathrm{Ind(xNorm_{R})}$.

For divided powers, again assume $P=\coprod^{\le1}_{1\le j\le n}R_{\delta_{j}}$. We have
  $$\Gamma^{n}(P^{\vee})\cong\mathrm{Sym}^{n}(P)^{\vee}$$
  is the dual in $\mathrm{Ban_{R}}$ of a free-module of finite rank and is therefore projective.  Now consider a projective of the form
 
 $$Q\cong\coprod_{i\in\mathcal{I}}^{\le 1}P_{i}$$

 with each $P_{i}$ of the form $R_{\delta_{i}}$. Then
 $$Q^{\otimes n}\cong\coprod^{\le1}_{[i]=[(i_{1},\ldots,i_{n})]\in(\mathcal{I}^{n}\big\slash\Sigma_{n})}\;\coprod^{\le1}_{[\sigma]\in\Sigma_{n}\big\slash\mathrm{Stab}(i)}P_{i_{\sigma_{1}}}\otimes\ldots\otimes P_{i_{\sigma(n)}}$$
 $\Gamma^{n}(Q)$ is the equaliser of the maps 
 $$\coprod^{\le1}_{[i]=[(i_{1},\ldots,i_{n})]\in(\mathcal{I}^{n}\big\slash\Sigma_{n}}\coprod^{\le1}_{[\sigma]\in\Sigma_{n}\big\slash\mathrm{Stab}(i)}P_{i_{\sigma_{1}}}\otimes\ldots\otimes P_{i_{\sigma(n)}}\rightarrow \coprod^{\le1}_{[i]=[(i_{1},\ldots,i_{n})]\in(\mathcal{I}^{n}\big\slash\Sigma_{n}}\coprod^{\le1}_{[\sigma]\in\Sigma_{n}\big\slash\mathrm{Stab}(i)}P_{i_{\sigma_{1}}}\otimes\ldots\otimes P_{i_{\sigma(n)}}$$
 Since contracting coproducts commute with kernels, this
  is the same as the contracting coproduct of the equaliser of the maps
 $$\coprod^{\le1}_{[\sigma]\in\Sigma_{n}\big\slash\mathrm{Stab}(i)}P_{i_{\sigma_{1}}}\otimes\ldots\otimes P_{i_{\sigma(n)}}\rightarrow\coprod^{\le1}_{[\sigma]\in\Sigma_{n}\big\slash\mathrm{Stab}(i)}P_{i_{\sigma_{1}}}\otimes\ldots\otimes P_{i_{\sigma(n)}}$$
  This is precisely $\Gamma^{n}(\coprod^{\le1}_{\sigma\in\Sigma_{n}}P_{i_{\sigma_{1}}}\otimes\ldots\otimes P_{i_{\sigma(n)}})$ which as we have seen is projective. Moreover this also shows that $\Gamma^{n}(Q)\cong\coprod_{[(i_{1},\ldots,i_{n})]\in\mathcal{I}^{n}\big\slash\Sigma^{n}}\Gamma^{n}(P_{i_{\sigma_{1}}}\otimes\ldots\otimes P_{i_{\sigma(n)}})$. A similar argument to the above without the non-expanding coproduct shows again that $\Gamma^{n}(P)$ is projective.

  Now we prove exactness of the Koszul complex. Consider first a projective of the form $R_{\delta}$. The Koszul comple in this case is just the map
  $$R_{\delta}\otimes\mathrm{Sym}(R_{\delta})\rightarrow\mathrm{Sym}(R_{\delta})$$
  which is clearly quasi-isomorphic to $\mathbb{I}$. Now consider a projecitve of the form $P=\coprod^{\le1}_{1\le i\le n}R_{\delta_{i}}$. The Koszul complex
  $$\Lambda(P)\otimes\mathrm{Sym}(P)$$
  is just given by tensoring up the Koszul complexes for each $R_{\delta_{i}}$. Since everything is flat, by K\"{u}nneth we again get a graded quasi-isomorphism
 $$\Lambda(P)\otimes\mathrm{Sym}(P)\rightarrow\mathbb{I}$$
Now for fixed $n$, in each homogeneous graded piece
  $$\Lambda^{(n-\bullet)}(P)\otimes\mathrm{Sym}^{\bullet}(P)$$
  is acyclic.
The formula for the differential implies that their norms are all bounded by $n$. For general $P\cong\coprod_{i\in\mathcal{I}}^{\le1}R_{\delta_{i}}$ we can write $P$ as a non-expanding filtered colimit of projectives of the form $\coprod^{\le1}_{1\le j\le n}R_{\epsilon_{j}}$. Then we use Proposition \ref{prop:nonexpandingexact}.

Next we prove exactness of the divided powers complex. 
  It is also useful to observe here that the natural isomorphism
  $$\Lambda^{n}(R^{\oplus m})^{\vee}\cong\Lambda^{n}((R^{\oplus m})^{\vee})$$
  is in fact an isomorphism of semi-normed/ normed/ Banach $R$-modules. 
For $P$ projective of finite rank the Koszul complex
 $$\Lambda^{(n-\bullet)}(P)\otimes\mathrm{Sym}^{\bullet}(P)$$
 is an acyclic complex of projectives, so taking its dual
 $$\Lambda^{(n-\bullet)}(P^{\vee})\otimes(\mathrm{Sym}^{\bullet}(P))^{\vee}$$
 is also an acyclic complex of projectives. But we have $(\mathrm{Sym}^{\bullet}(P))^{\vee}\cong\Gamma^{\bullet}(P^{\vee})$. Let $Q$ be a projective of finite rank whose predual is $P$, then 
 $$\Lambda^{(n-\bullet)}(Q)\otimes(\Gamma^{n}(Q))$$
 is also acyclic. Again the differentials have norm at most $n$. By the same argument using non-expanding colimits as for the Koszul comple we get the result for any projective of the form $P\cong\coprod^{\le1}_{i\in\mathcal{I}}R_{\delta_{i}}$. 
 \end{proof}

\begin{defn}
\begin{enumerate}
\item
The category of \textit{bornological }$R$-\textit{modules of convex type} is $\mathrm{Born}_{R}\defeq Ind^{m}(\mathrm{Norm}^{\frac{1}{2}}_{R})$.
\item
The category of \textit{separated bornological }$R$-\textit{modules} is $\mathrm{SBorn}_{R}\defeq Ind^{m}(\Norm_{R})$.
\item
The category of \textit{complete bornological }$R$-\textit{modules} is $\mathrm{CBorn}_{R}\defeq Ind^{m,}(\Ban_{R})$.
\end{enumerate}
\end{defn}

These categories are all concrete - in fact they are subcategories of the category of bornological modules. Let us recall how to describe the objects of these categories.

\begin{defn}\label{def:bornology}
A \textit{bornology} on a set \(X\) is a collection \(\mathfrak{B}\) of its subsets, satisfying the following:

\begin{itemize}
\item for every \(x \in X\), \(\{x\} \in \mathfrak{B}\);
\item if \(S\), \(T \in \mathfrak{B}\), then \(S \cup T \in \mathfrak{B}\);
\item if \(S \in \mathfrak{B}\), and \(T \subseteq S\), then \(T \in \mathfrak{B}\). 
\end{itemize}

\end{defn}

We call members of a bornology on a set its \textit{bounded subsets}. Let \((X,\mathfrak{B}_X)\) and \((Y,\mathfrak{B}_Y)\) be bornological sets. A function \(f \colon X \to Y\) is called \textit{bounded} if $f(U)\in\mathfrak{B}_{Y}$ for any $U\in\mathfrak{B}_{X}$. Bornological sets and bounded maps arrange into a category, $\mathsf{bSet}$.  

If \((X,\mathfrak{B}_X)\) and \((Y,\mathfrak{B}_Y)\) are bornological sets, then their product is the bornological set whose underlying set is $X\times Y$ such that a set $B\subset X\times Y$ is bounded if there are bounded subsets $B_{1}\in\mathfrak{B}_{X}, B_{2}\in\mathfrak{B}_{Y}$ such that $B\subset B_{1}\times B_{2}$.

If $R$ is a Banach ring and $M$ a semi-normed $R$-module, then $M$ can be equipped with a natural bornology called the \textit{von Neumann bornology}, whereby a subset $B\subset M$ is bounded if the restriction of the semi-norm on $M$ to $B$ is a bounded function. In this way, $R$ can itself be equipped with a bornology. 

\begin{defn}\label{def:bornological_module}
Let \(R\) be a Banach ring. A \textit{bornological \(R\)-module} \(M\) is an \(R\)-module with a bornology such that the scalar multiplication and addition maps
$$R\times M\rightarrow M$$
$$M\times M\rightarrow M$$
are bounded functions.
\end{defn}

We denote by $\mathsf{bMod}_{R}$ the subcategory of $\mathsf{bSet}$ consisting of bornological $R$-modules and bounded $R$-linear maps. To relate inductive systems with bornologies, first note that equipping a semi-normed $R$-module with the von Neumann bornology defines a canonical functor
$$b:\mathsf{sNMod}_R\rightarrow\mathsf{bMod}_{R}.$$

\begin{prop}
Let \(R\) be a Banach ring. There is a faithful functor $\mathsf{Ind}^m(\mathsf{sNMod}_R)\rightarrow\mathsf{bMod}_{R}$. If $R=k$ is a non-trivially valued Banach field, then this functor is fully faithful. The same claims hold for the categories \(\mathsf{Ind}^m(\mathsf{NMod}_k)\) and \(\mathsf{Ind}^m(\mathsf{Ban}_k)\), with the categories of separated and complete bornological \(k\)-vector spaces \(\mathsf{Born}_k\) and \(\mathsf{CBorn}_k\) on the right hand sides, respectively.
\end{prop}

For arbitrary $R$ it can be shown that the projective generators of $\mathrm{xNorm}_{R}$ are also projective generators for $\mathrm{Ind}^{m}(\mathrm{xNorm}_{R})$ rendering it a monoidal $\mathbf{AdMon}$-elementary quasi-abelian category. Thus there is induced projective model structure on $\mathrm{s}\mathrm{Ind}^{m}(\mathrm{xNorm}_{R})$ making it a combinatorial monoidal model category. Finally there is a  monoidal Quillen equivalence
$$\mathrm{s}\mathrm{Ind}^{m}(\mathrm{xNorm}_{R})\cong\mathrm{s}\mathrm{Ind}(\mathrm{xNorm}_{R})$$
Thus for the purposes of derived geometry, one may also work with the `concrete' bornological model.

\subsubsection{Locally Convex Spaces}

 For $k$ a Banach field the category $\mathrm{Pro}(\mathrm{Ban}_{k})$ contains the category $\hat{\mathcal{T}}_{c,k}$ of complete locally convex topological $k$-vector spaces as a full subcategory. Indeed if $E$ is an object of $\hat{\mathcal{T}}_{c,k}$ defined by a family of seminorms $\mathcal{P}$ then define $PB(E)=``\lim_{\leftarrow_{p\in\mathcal{P}}}"\hat{E}_{p}$ where $\hat{E}_{p}$ is the completion of $E$ with respect to the metric defined by the semi-norm $p$. This construction is functorial, lax monoidal, and $PB:\hat{\mathcal{T}}_{c,k}\rightarrow \mathrm{Pro}(\mathrm{Ban}_{k})$ is fully faithful. By an obvious Kan extension there is a canonical functor
$$PI:Pro(\mathpzc{E})\rightarrow Ind(\mathpzc{E})$$
Again this is lax monoidal. 

To a (complete) locally convex space $E$ one can functorially assign both the von Neumann bornology $vN(E)$ and the compact bornology $Cpt(E)$. The von Neumann bornology consists of the subsets of $E$ absorbed by all zero neighbourhoods. The compact bornology consists subsets with compact closure. These constructions define functors
$$\mathcal{T}_{c,k}\rightarrow\mathrm{sBorn}_{k}$$
where $\mathcal{T}_{c,k}$ is the category of complete locally convex topological $k$-vector spaces. 
 There is a natural transformation of functors $Cpt\rightarrow vN$. For details see \cite{MR2337277} Section 1.1.4. If $V$ is a nuclear locally convex space then the map $Cpt(V)\rightarrow vN(V)$ is an isomorphism by \cite{bambozzi2015stein} Lemma 3.67T

There is a natural isomorphism of functors $PI\circ PB\cong diss\circ vN$ (see \cite{bambozzi2015stein}), and therefore a natural transformation $diss\circ Cpt\rightarrow PI\circ PB$. Let $E$ be a Banach space and $F$ a complete locally convex space. Then
\begin{align*}
\Hom_{Ind(Ban)}(E,PI\circ PB(F))&\cong \Hom_{Pro(Ban)}(E,PB(F))\\
&\cong \Hom_{\hat{\mathcal{T}}_{c,k}}(E,F)\\
&\cong \Hom_{CBorn}(E,vN(F))\\
&\cong \Hom_{Ind(Ban)}(E,diss\circ vN(F))
\end{align*}
where the isomorphisms just arise from the fact that $Hom(E,-)$ commutes with projective limits in both $\hat{\mathcal{T}_{c}}$ and $CBorn$. When restricted to the category of nuclear Fr\'{e}chet spaces the functor $diss\circ vN$ is fully faithful by \cite{bambozzi2015stein} Example 3.22. Moreover since map $diss\circ Cpt\rightarrow PI\circ PB$ is a natural isomorphism this restriction of $PI\circ PB$ is also strong monoidal (\cite{MR2337277} Theorem 1.87). In particular the category of nuclear Fr\'{e}chet algebras over $\C$ embeds fully faithfully in the category of commutative complete bornological algebras. Since the category $CBorn_{\C}$ has good categorical properties, in particular it is closed monoidal \textbf{AdMon}-elementary, this is evidence that it provides a convenient setting in which to study analytic algebra.

\subsubsection{Non-Archimedean Modules and the Strong Context}

Although for our model of derived analytic geometry we will work with the context of ind-Banach/ Bornological modules, let us briefly mention another context which appears in the non-Archimedean setting. Let $k$ be a non-Archimedean semi-normed ring and consider the category $\xNorm_{k}^{nA,\le1}$. This is an additive, and in fact a quasi-abelian, category. However in \cite{kelly2016homotopy} there is another exact structure on this category called the \textit{strong exact structure} whereby a kernel-cokernel pair
\begin{displaymath}
\xymatrix{
0\ar[r] & K\ar[r]^{f}& L\ar[r]^{g} & M\ar[r] & 0
}
\end{displaymath}
is exact if for each $m\in M$ there is an $l\in L$ such that $g(l)=m$ and $||l||_{L}=||m||_{M}$. It is moreover shown in \cite{kelly2016homotopy} that, when equipped with the non-Archimedean projective monoidal structure, $\xNorm_{k}^{nA,\le1}$ is a projectively monoidal, weakly $\mathbf{AdMon}$-elementary, locally presentable exact category. Moreover in the cases $\Norm_{k}^{nA,\le1}$ and $\mathrm{Norm}^{\frac{1}{2},nA,\le1}_{k}$ it is in fact monoidal $\mathbf{AdMon}$-elementary. Moreover it has symmetric projectives. Indeed modules of the form $k_{r}$ for $r\in\mathbb{R}_{>0}$ form a projective generating set. Any symmetric power of such an object will just be a rescaling of $k$, and hence still projective. Tensoring with $k_{r}$ just amounts to rescaling, so is clearly exact and kernel preserving. Since tensor products commute with coproducts, are coproducts commute with kernels, we find that projectives are strongly flat. These context will be studied further in forthcoming work between the first and second authors and Jeroen Hekking.

\subsection{Fr\'{e}chet, Metrisable, and Nuclear Objects}

\subsubsection{Nuclear Objects}

Let $R$ be a Banach ring. By Corollary \ref{cor:strongnucnuc}, nuclear objects in $\mathbf{Ch}(\mathrm{CBorn}_{R})$ coincide with strongly nuclear objects. Moreover we have the following.

\begin{lem}[\cite{ben2020fr} Lemma 4.21]
   Discrete strongly nuclear objects in $\mathrm{Ind(Ban}_{R}\mathrm{)}$ are flat.
\end{lem}

Over a non-trivially valued Banach field, dual nuclear Fr\'{e}chet spaces are strongly nucler, and in fact are the basic nuclear objects. Nuclear Fr\'{e}chet spaces are also strongly nuclear objects.

\subsubsection{Fr\'{e}chet Spaces}

\begin{defn}
Let $R$ be a Banach ring. An object $F$ of $\mathrm{Ind(Ban}_{R}\mathrm{)}$ (resp. $\mathrm{Ind(Ban}^{nA}_{R}\mathrm{)}$ ) is said to be a \textit{bornological Fr\'{e}chet space} (resp. a \textit{non-Archimedean bornological Fr\'{e}chet space} if there exists a sequence
$$\cdots \rightarrow F_{n+1}\rightarrow F_{n}\rightarrow\cdots\rightarrow F_{1}\rightarrow F_{0}$$
with each $F_{i}\in \mathrm{Ban}_{R}$ (resp. each $F_{i}\in\mathrm{Ban}^{nA}_{R}$) such that $F\cong\lim_{\leftarrow}F_{i}$.
\end{defn}

\begin{lem}[\cite{ben2020fr}]
Bornological Fr\'{e}chet modules are $\aleph_{1}$-metrisable. 
\end{lem}


\subsubsection{The Open Mapping Property, Buchwalter's Theorem, and Derived Quotients}\label{subsubsec:derquotom}

\begin{defn}
A Banach ring $R$ is said to \textit{have the (non-Archimedean) open mapping property} if 
\begin{enumerate}
\item
a map of bornological (non-Archimedean) Fr\'{e}chet $R$-modules is an admissible epimorphism if and only if it is a surjection
\item
the kernel of a map of bornological (non-Archimedean) Fr\'{e}chet $R$-modules is a Fr\'{e}chet $R$-module.
\end{enumerate}
\end{defn}

\begin{example}
Any non-trivially valued (non-Archimedean) Banach field has the (non-Archimedean) open mapping property. Any algebra over such a field also has the open mapping property.
\end{example}

A related and useful result is Buchwalter's Theorem.

\begin{defn}
\begin{enumerate}
    \item
    An object $F$ of $\mathrm{CBorn}_{R}$ is said to have a \textit{bornology with a countable basis} if, when regarded as a bornological $R$-module, there is a countable collection $\{B_{n}\}_{n\in\mathbb{N}}$ of bounded subsets such that every bounded subset is contained in some $B_{n}$. 
    \item 
    $R$ is said to be \textit{have Buchwalter's property} if whenever $f:F\rightarrow G$ is a surjective map in $\mathrm{CBorn}_{R}$ and $F$ has a bornology with a countable basis, then $f$ is a strict epimorphism.
\end{enumerate}
\end{defn}

In particular if $k$ is a non-trivially valued Banach field then Buchwalter's Theorem (\cite{MR3577218} Theorem 4.9) implies that it has Buchwalter's property. The following is immediate.

\begin{lem}
Let $R$ have the (non-Archimedean) open mapping property (resp. Buchwalter's property). A complex $F_{\bullet}$ of Fr\'{e}chet $R$-modules (resp. of $R$-modules whose bornology has a countable basis) is acyclic if and only if it is algebraically exact.
\end{lem}

Let $R$ be a Banach ring satisfying the open mapping property (resp. Buchwalter's property). Let $A\in\mathrm{Comm}(\mathrm{Ind(Ban}_{R}\mathrm{)})$ be a bornological Fr\'{e}chet $R$-algebra (resp. an $R$-algebra whose bornology has a countable basis). Let $a_{1},\ldots,a_{n}\in A^{\oplus n}$ be an algebraically regular sequence  such that $(a_{0},\ldots,a_{n})\subset A$ is an adimssible ideal. Consider the Koszul complex
$$\wedge^{n}A^{n}\rightarrow\cdots\rightarrow\wedge^{2}A^{n}\rightarrow R\rightarrow A\big\slash(a_{0},\ldots,a_{n})$$
This is a sequence of Fr\'{e}chet modules (resp. a sequence of modules whose bornology have a countable basis). Since the sequence is regular it is algebraically exact, and by the open mapping property (resp. Buchwalter's property) it is admissibly exact. Now $\wedge^{n}A^{n}\rightarrow\cdots\rightarrow\wedge^{2}A^{n}\rightarrow R$ computes $A\big\slash\big\slash(a_{0},\ldots,a_{n})$. Thus in this case we have
$$A\big\slash\big\slash(a_{0},\ldots,a_{n})\rightarrow A\big\slash(a_{0},\ldots,a_{n})$$
is an equivalence.

\section{Some General Bornological Algebra}\label{sec:genbornalg}

Let us give some useful general results concerning bornological algebra. We start with the following somewhat trivial, but very useful, fact and its corollary.

\begin{lem}
Let $A$ be a bornological algebra, and let $M$, $M'$, and $N$ be bornological $A$-modules. Let $\pi:N\rightarrow M$ be an admissible epimorphism of $A$-modules, and $\phi:M\rightarrow M'$ be a map of $A$-modules such that $\phi\circ\pi$ is bounded. Then $\phi$ is bounded. 
\end{lem}

\begin{proof}
Let $B\subset M$ be bounded. There is a $\tilde{B}\in N$ such that $\pi(\tilde{B})=B$. Then $\phi(B)=\phi\circ\pi(\tilde{B})$ which is bounded.
\end{proof}

\begin{cor}
Let $A$ be a bornological algebra. Let $\phi:M\rightarrow M'$ be any map of finitely generated left $A$-modules. Then $\phi$ is bounded.
\end{cor}

\begin{proof}
Let $A^{m}\rightarrow M$ and $A^{m'}\rightarrow M'$ be admissible epimorphisms. By projectivity of $A^{m}$, the map $M\rightarrow M'$ lifts to a map $\overline{\phi}:A^{m}\rightarrow A^{m'}$. This will just correspond to a matrix of elements of $A$, multiplication by which is bounded since addition and multiplication in $A$ are bounded.
\end{proof}

\subsection{Semi-Completeness and Topological Nilpotence}

For this section fix a non-trivially valued Banach field $k$.

\begin{defn}[\cite{bambozzi}, Definition 2.1.40]
Let $E$ be a bornological $k$-vector space and $\Phi$ a filter of subsets of $E$. $\Phi$ is said to \textit{converge to }$0$\textit{ in the sense of Mackey} if there exists a bounded subset $B\subset E$ such that for every $\lambda\in k^{*}$ we have $\lambda B\in\Phi$. $\Phi$ is said to \textit{converge to }$a\in E$ if $\Phi-a$ converges to $0$. If $U\subset E$ then a sequence $\{x_{n}\}_{n\in\mathbb{N}}$ is said to be \textit{bornologically convergent in }$U$ if it converges in the sense of Mackey to an element of $U$.
\end{defn}

\begin{defn}[\cite{bambozzi}, Definition 2.1.41]
Let $E$ be a bornological $k$-vector space. A sequence $\{x_{n}\}_{n\in\mathbb{N}}\subset E$ is said to be \textit{Cauchy-Mackey} if the double sequence sequence $\{x_{n}-x_{m}\}_{n,m\in\mathbb{N}}$ converges to $0$ in the sense of Mackey. $E$ is said to be \textit{semi-complete} if all Cauchy-Mackey sequences in $E$ converge to some limit in the sense of Mackey.
\end{defn}

As mentioned in \cite{bambozzi} complete bornological spaces are semi-complete.

\begin{defn}[\cite{bambozzi2015stein} Definition 3.2]
Let $E$ be a bornological $k$-vector space. A subset $U\subset E$ is said to be \textit{bornologically closed} if whenever $\{x_{n}\}_{n\in\mathbb{N}}$ is a sequence of elements in $U$ which is bornologically convergent in $E$, is bornologically convergent in $U$.
\end{defn}

\begin{defn}[\cite{bambozzi2015stein} Definition 3.9]
A bornological vector space of convex type is called proper if its bornology has a basis of
bornologically closed subsets.
\end{defn}

The utility of proper spaces to us is the following result.

\begin{prop}[\cite{bambozzi2015stein} Proposition 3.15]
Let $E$ be a complete proper bornological space of convex type and $F\subset E$ a subspace. Then the bornological closure of $F$ coincides with the set of limit points of bornologically convergent sequences in $F$. 
\end{prop}

It is useful to make the following definition.

\begin{defn}
Let $E$ be a bornological $k$-vector space and $U\subset E$ a subset. 
\begin{enumerate}
\item $U$ is said to be \textit{bornologically dense in }$E$ if $E$ is the smallest bornologically closed subset of $E$ containing $U$.
    \item 
    $U$ is said to be \textit{strongly dense in }$E$ if for any element $e\in E$ there is a sequence $\{x_{n}\}_{n\in\mathbb{N}}$ which converges to $e$ in the sense of Mackey.
\end{enumerate}

\end{defn}

If $U\subset E$ is strongly dense then it is bornologically dense. The converse is not true in general, but it is for proper spaces.

\begin{defn}[\cite{bambozzi2015stein}, Definition 3.9]
    A bornological vector space $E$ is said to be \textit{proper} if its bornology a basis consisting of bornologically closed subsets.
\end{defn}

\begin{defn}[\cite{bambozzi}, Definition 2.1.42]
Let $A$ be a bornological algebra. $a\in A$ is said to be \textit{topologically nilpotent} if the sequence $\{a^{n}\}_{n\in\mathbb{N}}$ converges to $0$ in the sense of Mackey.
\end{defn}

The set of topologically nilpotent elements of $A$ is denote by $A^{\circ\circ}$. 

\begin{prop}[\cite{bambozzi}, Proposition 2.1.47]\label{prop:topnilpinvert}
Let $A$ be a semi-complete bornological algebra and let $x=1-y$ with $y\in A^{\circ \circ}$. Then $x$ is a unit in $A$ and 
$$x^{-1}=\sum_{n=0}^{\infty}y^{n}$$
\end{prop}

\begin{lem}[c.f. \cite{BGR} Lemma 6]
Let $A$ be a semi-complete bornological algebra, and $M$ an $A$-module. Let $N$ be a submodule of $M$ such that there are elements $x_{1},\ldots,x_{n}$ in $M$ with the property that $M\subset N+\Sigma_{i=1}^{n}A^{\circ\circ}x_{i}$. Then $N=M$. 
\end{lem}

\begin{proof}
The proof works exactly as in \cite{BGR} Lemma 6 using Proposition \ref{prop:topnilpinvert}.
\end{proof}

\begin{rem}
Let $E$ be a bornological $k$-vector space of convex type. Let $\{x_{n}\}_{n\in\mathbb{N}}$ be a sequence in $E$ converging to $a\in E$ in the sense of Mackey. Then $\{x_{n}\}_{n\in\mathbb{N}}$ converges to $a$ topologically in $E^{t}$. This follows immediately from \cite{bambozzi2015stein} Remark 3.4. In particular if $U\subset E$ is strongly dense, then $U\subset E^{t}$ is topologically dense.
\end{rem}

\begin{prop}
Let $E$ be a separated bornological $k$-vector space of convex type. Then the map $E\rightarrow\hat{E}$ where $\hat{E}$ is the completion of $E$, is bornologically dense. 
\end{prop}

\begin{proof}
Write $E\cong\limind_{\mathcal{I}} E_{i}$ with each $E_{i}$ being a normed space. Then $\hat{E}\cong\limind_{\mathcal{I}}\hat{E}_{i}$. Now $\limind_{\mathcal{I}}\hat{E}_{i}$ is a quotient of $\bigoplus_{i\in\mathcal{I}}\hat{E}_{i}$. The image of the map $\bigoplus_{i\in\mathcal{I}}E_{i}\rightarrow\bigoplus_{i\in\mathcal{I}}\hat{E}_{i}$ clearly has bornologically dense image, and so descends to a map on quotients
$$E\rightarrow\hat{E}$$
which has bornologically dense image. 
\end{proof}

\subsubsection{Sub-Multiplicative Algebras}


\begin{defn}
\begin{enumerate}
    \item 
    A bornological algebra of convex type $A$ is said to \textit{have a multiplicative bornology} if it can be written as a monomorphic colimit $A\cong\limind_{\mathcal{I}}A_{i}$ where each $A_{i}$ is a Banach algebra.
    \item 
    A Banach algebra $A$ is said to be \textit{sub-multiplicative} if $||ab||\le ||a|| ||b||$ for all $a,b\in A$.
    \item 
    A bornological algebra of convex type $A$ is said to be \textit{sub-multiplicative} if it has a presentation 
    $$A\cong\limind_{\mathcal{I}}A_{i}$$
    with each $A_{i}$ a sub-multiplicative Banach algebr and each map $A_{i}\rightarrow A_{i+1}$ is non-expanding.
\end{enumerate}

\end{defn}

\begin{example}
    Let $F$ be a bornological Fr\'{e}chet algebra of the form 
    $$F\cong\limpro_{\mathbb{N}}F_{n}$$
    with each $F_{n}$ a sub-multiplicative Banach algebra and each $F_{n+1}\rightarrow F_{n}$ being a non-expanding map of Banach algebras. Then $F$ is a sub-multiplicative bornological algebra by \cite{ben2020fr} Corollary 5.13.
\end{example}

\begin{prop}
Let $A$ be sub-mutliplicative. Then $A^{\circ\circ}$ is a neighbourhood of $0$ in $A^{t}$.
\end{prop}

\begin{proof}
Consider the ball $B_{A_{i}}(0,1)$ and $V=\bigcup_{i\in\mathcal{I}}B_{A_{i}}(0,1)\subset A$. This is open and clearly consists of topologically nilpotent elements.
\end{proof}

\begin{lem}\label{lem:completefg}
Let $A$  be submultiplicative. Let $M$ be a separated left bornological $A$-module whose completion $\hat{M}$ is a finitely generated $A$-module. Then $M$ is complete.
\end{lem}

\begin{proof}
Let $\pi:A^{\oplus n}\rightarrow\hat{M}$ be an admissible epimorphism. Then $\pi^{t}:(A^{t})^{\oplus n}\rightarrow\hat{M}^{t}$ is an admissible epimorphism. Denote the corresponding generators by $x_{1},\ldots,x_{m}$.  $(A^{\circ\circ})^{\oplus n}$ is a non-empty open neighbourhood of $0$ and $M^{t}\rightarrow\hat{M}^{t}$ has dense image. Therefore we must have $x_{j}\in M+\sum_{i=1}^{n}A^{\circ\circ}x_{i}$ for each $j$, so $M=\hat{M}$.
\end{proof}

\begin{cor}\label{cor:quotnoetherianproper}
Let $A$ be submultiplicative, and let $M$ be a complete bornological $A$-module. Suppose $M\cong E\big\slash F$ where $E$ is a complete, proper bornological $A$-module, and $F\subset E$ a bornologically closed sub-module. If $E$ is Noetherian then all submodules of $M$ are bornologically closed.
\end{cor}

\begin{proof}
Let $N\subset M$ be a submodule Denote by $\pi:E\rightarrow M$ the projection map, and consider $N'\defeq\pi^{-1}(N)$. Its completion $\hat{N'}$ is the bornological closure of $N'$ in $E$. In particular it is a submodule of $E$ and is therefore finitely generated. Thus we have $N'=\hat{N'}$. $N'$ is therefore bornologically closed in $E$. Now we have $N\cong N'\big\slash F\cap N'$ which is complete, and therefore closed in $M$ 
\end{proof}

\begin{defn}
A bornological left $A$-module $M$ is said to be \textit{strongly left Noetherian} if its underlying $A_{alg}$-module is Noetherian, and all submodules of $M$ are closed.
\end{defn}

\begin{lem}
Let $A$ be submultiplicative. Suppose further that $A$ is proper and Noetherian. Then any finitely generated $A$-module is strongly Noetherian. 
\end{lem}

\begin{proof}
Since $A$ is proper (resp. Noetherian), $A^{\oplus n}$ is proper (resp. Noetherian) for any $n$.  Now the claim follows immediately from Corollary \ref{cor:quotnoetherianproper}.
\end{proof}

\subsection{Fr\'{e}chet-Stein Algebras}

The formalism of Fr\'{e}chet-Stein algebras has proved very useful in non-Archimedean geometry \cite{MR1990669}. Here we give a definition of bornological Fr\'{e}chet-Stein algebras, in both the Archimedean and non-Archimedean settings. Fix a non-trivially valued Banach field $k$.

%
%

\begin{defn}
A complete bornological algebra $A$ is said to be a \textit{bornological Fr\'{e}chet-Stein} algebra if there is a sequence
$$A\rightarrow\cdots \rightarrow \overline{A}_{n+1}\rightarrow A_{n+1}\rightarrow\overline{A}_{n}\rightarrow A_{n}\rightarrow\dots\rightarrow\overline{A}_{0}\rightarrow A_{0}$$
and isomorphisms
$$A\cong\lim_{\leftarrow_{n\in\mathbb{N}}}A_{n}$$
$$A\cong\lim_{\leftarrow_{n\in\mathbb{N}}}\overline{A}_{n}$$
where
\begin{enumerate}
\item
$A$ is a bornological nuclear Fr\'{e}chet space. 
\item
each $\overline{A}_{n}$ is a strongly left Noetherian algebra, whose bornology has a countable basis.
\item 
finitely presented left-$\overline{A}_{n}$-modules $M$ are proper and nuclear.
\item
each $A_{n}$ is a Fr\'{e}chet algebra which is submultiplicative as a bornological algebra.
\item
for all $n$, $A\rightarrow A_{n}$ and $A\rightarrow\overline{A}_{n}$ are epimorphisms of complete bornological algebras. Moreover, the image of $A_{n+1}$ in $\overline{A}_{n}$ is strongly dense. 
\item
$\overline{A}_{n}$ is transverse over $\overline{A}_{n+1}$ to left finitely generated $\overline{A}_{n+1}$-modules.
\end{enumerate}
The data of a bornological Fr\'{e}chet-Stein algebra will include the projective systems. We denote this data by $(A,A_{n},\overline{A}_{n})$.
\end{defn}

This definition may seem overly complicated compared to the one in \cite{MR1990669}. The issue with the definition in loc. cit appears in the case of alegbras of analytic functions in the Archimedean setting, such as the algebra of entire functions on $\mathbb{C}$, $\mathcal{O}(\mathbb{C})$. As we will discuss later, this algebra can also be written as a limit of Banach algebras (`Archimedean Tate algebras'). Unfortunately unlike in the non-Archimedean case these are not Noetherian. However it can also be written as a projective limit of `overconvergent functions', which will be strongly Noetherian.

\begin{defn}
Let $(A,A_{n},\overline{A}_{n})$ be a bornological Fr\'{e}chet Stein algebra. A \textit{left pre-quasi-coherent sheaf for } $(A,A_{n},\overline{A}_{n})$ is a tuple $(M_{n},\overline{M}_{n},\alpha_{n},\beta_{n})$ where
\begin{enumerate}
\item
Each $M_{n}$ is a bornological Fr\'{e}chet space.
\item
$M_{n} \in {}_{A_{n}}\mathrm{Mod}(\mathrm{LH(Ind(Ban}_{k}\mathrm{)})$,
\item
$\overline{M}_{n}\in{}_{\overline{A}_{n}}\mathrm{Mod}(\mathrm{LH(Ind(Ban}_{k}\mathrm{)})$,
\item
$\alpha_{n}:A_{n}\hat{\otimes}_{\overline{A}_{n}}\overline{M}_{n}\rightarrow M_{n}$ is an isomorphism.
\item
$\beta_{n}:\overline{A}_{n}\hat{\otimes}_{A_{n+1}} M_{n+1}\rightarrow \overline{M}_{n}$ is an isomorphism.
\end{enumerate}
such that the compositions
\begin{displaymath}
\xymatrix{
\overline{A}_{n}\otimes_{\overline{A}_{n+1}}\overline{M}_{n+1}\ar[r] &\overline{A}_{n}\otimes_{A_{n+1}}A_{n+1}\otimes_{\overline{A}_{n+1}}\overline{M}_{n+1}\ar[rrr]^{\overline{A}_{n}\otimes_{A_{n+1}}\alpha_{n+1}} &&&\overline{A}_{n}\otimes_{A_{n+1}}M_{n+1}\rightarrow\overline{M}_{n}
}
\end{displaymath}
induce equivalences
$$\overline{A}_{n}\hat{\otimes}^{\mathbb{L}}_{\overline{A}_{n+1}}\overline{M}_{n+1}\rightarrow\overline{M}_{n+1}$$
\end{defn}

Pre-quasi-coherent sheaves arrange into an obvious category, which we denote by $\mathrm{PQCoh}_{(A,A_{n},\overline{A}_{n})}$. Note there are  natural maps 
$$M_{n+1}\rightarrow A_{n}\hat{\otimes}_{A_{n+1}}M_{n+1}\rightarrow M_{n}$$
$$\overline{M}_{n+1}\rightarrow A_{n}\hat{\otimes}_{A_{n+1}}\overline{M}_{n+1}\rightarrow \overline{M}_{n}$$
Denote by 
$$\Gamma:\mathrm{QCoh}_{(A,A_{n},\overline{A}_{n})}\rightarrow{}_{A}\mathrm{Mod}(\mathrm{LH(Ind(Ban}_{k}\mathrm{)})$$
the functor given by
$$(M_{n},\overline{M}_{n},\alpha_{n},\beta_{n})\mapsto\lim_{\leftarrow_{n}}M_{n}$$
Denote the essential image of this functor by $\mathcal{PQ}_{(A,A_{n},\overline{A}_{n})}$.

\begin{defn}
A pre-quasi-coherent sheaf $(M_{n},\overline{M}_{n},\alpha_{n},\beta_{n})$ is said to be \textit{quasicoherent} if 
\begin{enumerate}
\item
$\Gamma(M_{n},\overline{M}_{n},\alpha_{n},\beta_{n})$ is a bornological nuclear Fr\'{e}chet space. 
\item
The map $\overline{M}_{n+1}\rightarrow \overline{M}_{n}$ has strongly dense image.
\end{enumerate}
We denote the full subcategory of $\mathrm{PQCoh}_{(A,A_{n},\overline{A}_{n})}$ consisting of quasi-coherent sheaves by $\mathrm{QCoh}_{(A,A_{n},\overline{A}_{n})}$, and the essential image of the functor $\Gamma$ restricted to $\mathrm{QCoh}_{(A,A_{n})}$ by $\mathcal{Q}_{(A,A_{n},\overline{A}_{n})}$
\end{defn}

\begin{prop}
 Let $(M_{n},\overline{M}_{n},\alpha_{n},\beta_{n})$ be a pre-quasicoherent sheaf with each $\overline{M}_{n}$ finitely presented. Then $(M_{n},\overline{M}_{n},\alpha_{n},\beta_{n})$ is quasicoherent.
\end{prop}

\begin{proof}
First we show that $\overline{M}_{n+1}\rightarrow \overline{M}_{n}$ has bornologically dense image, and hence strongly dense image. Pick an admissible epimorphism $\overline{A}^{r}_{n+1}\rightarrow\overline{M}_{n+1}$. We have a commutative diagram
\begin{displaymath}
\xymatrix{
\overline{A}^{r}_{n+1}\ar[d]\ar[r] & \overline{M}_{n+1}\ar[d]\\
\overline{A}^{r}_{n}\ar[r] & \overline{M}_{n}
}
\end{displaymath}
The left vertical map is an epimorphism, and the bottom horizontal map is an admissible epimorphism. Thus the right vertical map is also an epimorphism.

Now each $\overline{M}_{n}$ is assumed nuclear. By \cite{bambozzi2015stein} Proposition 3.51, $\Gamma(M_{n},\overline{M}_{n},\alpha_{n},\beta_{n})$ is nuclear, as required.
\end{proof}

\begin{defn}
A pre-quasicoherent sheaf $(M_{n},\overline{M}_{n},\alpha_{n},\beta_{n})$ is said to be \textit{coadmissible} if each $M_{n}$ is finitely generated as a left $A_{n}$-module. The full subcategory of $\mathrm{QCoh}_{(A,A_{n},\overline{A}_{n})}$ is denoted $\mathrm{Coad}_{(A,A_{n},\overline{A}_{n})}$
\end{defn}

By the bornological Mittag-Leffler theorem, specifically \cite{bambozzi2015stein} Corolalry 3.80, we have the following.

\begin{cor}\label{cor:limitunder}
Let $(M_{n},\overline{M}_{n},\alpha_{n},\beta_{n})$ be a quasi-coherent sheaf for $(M_{n},\overline{M}_{n},\alpha_{n},\beta_{n})$ and write $M=\Gamma(M_{n},\overline{M}_{n},\alpha_{n},\beta_{n})$. The map $\mathbb{R}\lim_{\leftarrow_{n}}M_{n}\rightarrow\lim_{\leftarrow_{n}}M_{n}$ is an equivalence. In particular
$$\mathbb{R}\lim_{\leftarrow_{n}}\overline{M}_{n}\cong\mathbb{R}\lim_{\leftarrow_{n}}M_{n}\rightarrow\lim_{\leftarrow_{n}}M_{n}\cong\lim_{\leftarrow_{n}}\overline{M}_{n}$$
is also an equivalence.
\end{cor}

\begin{prop}\label{prop:coadexact}
    The category $\mathrm{Coad}_{(A,A_{n},\overline{A}_{n})}$ is abelian. Moreover the functor $\Gamma:\mathrm{Coad}_{(A,A_{n},\overline{A}_{n})}\rightarrow{}_{A}\mathrm{Mod}(\mathrm{CBorn}_{k})$ is exact. 
\end{prop}

\begin{proof}
    Let $f:(M_{n},\overline{M}_{n},\alpha_{n},\beta_{n})\rightarrow (M'_{n},\overline{M}'_{n},\alpha'_{n},\beta'_{n})$ be a map in $\mathrm{Coad}_{(A,A_{n},\overline{A}_{n})}$. Define
    $$\mathrm{coker}(f)\defeq(\mathrm{coker}(M_{n}\rightarrow M'_{n}),\mathrm{coker}(\overline{M}_{n}\rightarrow\overline{M}'_{n}),\tilde{\alpha}_{n},\tilde{\beta}_{n})$$
    where $\tilde{\alpha}_{n}$ and $\tilde{\beta}_{n}$ are the obvious induced maps on quotients. Each $\mathrm{coker}(\overline{M}_{n}\rightarrow\overline{M}'_{n})$ is finitely generated, so $\mathrm{coker}(f)$ is coherent.

    Define
     $$\mathrm{ker}(f)\defeq(A_{n}\hat{\otimes}_{\overline{A}_{n}}\mathrm{ker}(\overline{M}_{n}\rightarrow\overline{M}'_{n}),\mathrm{ker}(\overline{M}_{n}\rightarrow\overline{M}'_{n}),{\alpha}_{n},\tilde{\beta}_{n})$$
     By the transversality asusmption for quasi-coherent sheaves, we have
     $$\overline{A}_{n}\hat{\otimes}_{\overline{A}_{n+1}}\mathrm{ker}(\overline{M}_{n+1}\rightarrow\overline{M}'_{n+1})\cong \mathrm{ker}(\overline{M}_{n}\rightarrow\overline{M}'_{n})$$
     It is straightforward to check that this is a kernel. Moreover the strongly Noetherian assumption on $\overline{A}_{n}$ imply that it is an abelian category.

    Now since each $\overline{A}_{n}$ is strongly left Noetherian, the category ${}_{\overline{A}_{n}}\mathrm{Mod}^{fp}(\mathrm{CBorn})\cong{}_{\overline{A}_{n}}\mathrm{Mod}^{fp}(\mathrm{Ab})$ of finitely presented $\overline{A}_{n}$-modules is an exact abelian subcategory of ${}_{\overline{A}_{n}}\mathrm{Mod}(\mathrm{CBorn})$. The claim then follows immediately from Corollary \ref{cor:limitunder}. 
\end{proof}

\begin{lem}
Let $(M_{n},\alpha_{n})$ be a quasi-coherent sheaf for $(A,A_{n})$ and write $M=\Gamma(M_{n},\alpha_{n})$. The map $\Gamma(M)\hat{\otimes}_{A}\overline{A}_{i}\rightarrow M_{i}$ is an epimorphism. If $(M_{n},\alpha_{n})$ is coadmissible then the map is an isomorphism.
\end{lem}

\begin{proof}
Since $M$ is nuclear, by the Mittag-Leffler theorem in \cite{bambozzi2015stein}, we have an exat sequence
$$0\rightarrow M\rightarrow\prod_{n} M_{n}\rightarrow\prod_{n}M_{n}\rightarrow 0$$
Tensoring with $A_{i}$, and using the fact that both $A_{i}$ and $A$ are Fr\'{e}chet, we get a right exact sequence
$$A_{i}\hat{\otimes}_{A}M\rightarrow\prod_{n}A_{i}\hat{\otimes}_{A}M_{n}\rightarrow\prod_{n}A_{i}\hat{\otimes}_{A}M_{n}\rightarrow 0$$
Now the kernel of $\prod_{n}A_{i}\hat{\otimes}_{A}M_{n}\rightarrow\prod_{n}A_{i}\hat{\otimes}_{A}M_{n}$ computes the limit of the diagram which is eventually the constant diagram on $M_{i}$ since $A\rightarrow A_{i}$ is an epimorphism. Thus the map $A_{i}\hat{\otimes}_{A}M\rightarrow M_{i}$ is an epimorphism.

Now suppose that each $M_{i}$ is finitely generated. In particular it is Banach. Thus $M$ is Fr\'{e}chet and the bornology a countable basis. First we claim that $\overline{A}_{i}\hat{\otimes}_{A}M\rightarrow \overline{M}_{i}$ is an admissible epi. Now $\overline{A}_{i}\hat{\otimes}_{A}M\rightarrow \overline{M}_{i}$ is an epimorphism.  The set theoretic image of $\overline{A}_{i}\otimes_{A}M$ is a dense $\overline{A}_{i}$-submodule. Any $\overline{A}_{i}$-submodule is closed by assumption. Thus the set-theoretic image must be all of $\overline{M}_{i}$, so we get a surjection The bornology of $\overline{A}_{i}\hat{\otimes}_{A}M\rightarrow \overline{M}_{i}$. $\overline{A}_{i}\hat{\otimes}_{A}M$ has a countable basis. By Buchwalter's theorem (see e.g. \cite{MR3577218} Theorem 4.9), it is a strict epimorphism. Now an argument identical to \cite{MR1990669} Corollary 3.1 shows that it is in fact an isomorphism. 
\end{proof}

Immediately we get the following.

\begin{cor}\label{cor:ffcoh}
The functor
$$\Gamma:\mathrm{QCoh}_{(A,A_{n},\overline{A}_{n})}\rightarrow\mathcal{Q}_{(A,A_{n},\overline{A}_{n})}$$
is faithful. The restriction to $\mathrm{Coad}_{(A,A_{n},\overline{A}_{n})}$ is fully faithful. In particular $\mathcal{C}_{(A,A_{n},\overline{A}_{n})}$ is an abelian category closed under kernels and cokernels in ${}_{A}\mathrm{Mod}(\mathrm{CBorn}_{k})$.
\end{cor}

The inverse functor $\mathcal{C}_{(A,A_{n},\overline{A}_{n})}\rightarrow\mathrm{Coad}_{(A,A_{n},\overline{A}_{n})}$ sends $M$ to $(M,\overline{A}_{n}\otimes_{A}M,A_{n}\otimes_{A}M)$ with the obvious structure maps. In particular we get the following.

\begin{cor}
Let
$$0\rightarrow K\rightarrow M\rightarrow N\rightarrow 0$$
be an exact sequence in $\mathcal{C}_{(A,A_{n},\overline{A}_{n})}$. Then
$$0\rightarrow \overline{A}_{n}\hat{\otimes}_{A}K\rightarrow\overline{A}_{n}\hat{\otimes}_{A}M\rightarrow\overline{A}_{n}\hat{\otimes}_{A}N\rightarrow0$$
is an exact sequence in $\mathrm{Coad}_{(A,A_{n},\overline{A}_{n})}$. 
\end{cor}

\begin{lem}[c.f. \cite{MR1990669} Lemma 3.6]\label{lem:cohfg}
Let $M\in\mathcal{C}_{(A,A_{n},\overline{A}_{n})}$, and $N\subseteq M$. The following are equivalent.
\begin{enumerate}
\item
$N$ is in $\mathcal{C}_{(A,A_{n},\overline{A}_{n})}$.
\item
$M\big\slash N$ is in $\mathcal{C}_{(A,A_{n},\overline{A}_{n})}$.
\item
$N$ is bornologically closed in $M$.
\end{enumerate}

Moreover if each map $A\rightarrow A_{n}$ is transverse to $\mathcal{C}_{(A,A_{n},\overline{A}_{n})}$ then $\mathcal{C}_{(A,A_{n},\overline{A}_{n})}$ is thick in $\mathrm{CBorn}_{k}$
\end{lem}

\begin{proof}
$1)$ and $2)$ are equivalent by Corollary \ref{cor:ffcoh}. Suppose $N$ is bornologically closed in $M$. Write $M=\Gamma(M_{n},\overline{M}_{n},\alpha_{n},\beta_{n})$. For each $n$ let $\overline{N}_{n}\subseteq \overline{M}_{n}$ denote the submodule of $\overline{M}_{n}$ generate by the image of $N$ in $\overline{M}_{n}$. Define $N_{n}\defeq A_{n}\hat{\otimes}_{\overline{A}_{n}}N_{n}$. $\alpha$ and $\beta$ restrict to give well-defined maps such that $(N_{n},\overline{N}_{n},\alpha_{n},\beta_{n})$ is in $\mathrm{Coad}_{(A,A_{n},\overline{A}_{n})}$. Consider the exact sequence $0\rightarrow \overline{N}_{n}\rightarrow \overline{M}_{n}\rightarrow \overline{M}_{n}\big\slash \overline{N}_{n}\rightarrow 0$. Since everything in this sequence is coherent, we get an exact sequence
$$0\rightarrow \Gamma(N_{n},\overline{N}_{n},\alpha_{n},\beta_{n})\rightarrow M\rightarrow\Gamma(M_{n}\big\slash N_{n},\overline{M}_{n}\big\slash\overline{N}_{n},\overline{\alpha}_{n},\overline{\beta}_{n})\rightarrow 0$$
In particular $\Gamma(N_{n},\overline{N}_{n},\alpha_{n},\beta_{n})$ is a bornologically closed subset of $M$. There is a natural map $N\rightarrow\Gamma(N_{n},\overline{N}_{n},\alpha_{n},\beta_{n})$. Now $N$ is topologically dense in $\Gamma(N_{n},\overline{N}_{n},\alpha_{n},\beta_{n})^{t}$. Indeed let $(k_{n})_{n\in\mathbb{N}}\in\Gamma(N_{n},\overline{N}_{n},\alpha_{n},\beta_{n})$ with each $k_{n}$ in $\overline{N}_{n}$. Since every submodule of $\overline{M}_{n}$ is bornologically closed, there is $x_{n}\in N$ with $k_{n}=\overline{\rho}_{n}(x_{n})$, where $\overline{\rho}_{n}:M\rightarrow\overline{M}_{n}$ is the canonical map. The sequence $\{x_{n}\}$ converges to $(k_{n})_{n\in\mathbb{N}}$.

Since $\Gamma(N_{n},\overline{N}_{n},\alpha_{n},\beta_{n})$ is a bornological nuclear Fr\'{e}chet space, $N$ is also bornologically dense in $\Gamma(N_{n},\overline{N}_{n},\alpha_{n},\beta_{n})$ by \cite{bambozzi2015stein} Lemma 3.68. Thus $\Gamma(N_{n},\alpha_{n})$ is the bornological closure of $N$.

So if $N$ is closed it is in $\mathcal{C}_{(A,A_{n},\overline{A}_{n})}$. 

Conversely suppose that $N$ is in $\mathcal{C}_{(A,A_{n},\overline{A}_{n})}$. Then $N_{n}\cong A_{n}\hat{\otimes}_{A}N$ and $N\cong\Gamma(A_{n}\hat{\otimes}_{A}N)$ is closed. 

Now suppose each map $A\rightarrow A_{n}$ is transverse to $\mathcal{C}_{(A,A_{n},\overline{A}_{n})}$. Let
$$0\rightarrow N\rightarrow M\rightarrow M\big\slash N\rightarrow 0$$
be an exact sequence in $\mathrm{CBorn}_{k}$ with $N$ and $M\big\slash N$ in $\mathcal{C}_{(A,A_{n},\overline{A}_{n})}$. Define $\overline{M}_{n}\defeq\overline{A}_{n}\hat{\otimes}_{A}M$ and $M_{n}\defeq A_{n}\hat{\otimes}_{A}M$. By the transversality assumption we have exact sequences
$$overline{N}_{n}\rightarrow\overline{M}_{n}\rightarrow (M\big\slash N_{n})\rightarrow 0$$
Thus $\overline{M}_{n}$ is finitely generated, and hence finitely presented. In particular each $M_{n}$ is also finitely presented. Moreover taking the (derived) limit gives a commutative diagram with exact top and bottom rows
\begin{displaymath}
\xymatrix{
0\ar[r] &N\ar[d]\ar[r] & M\ar[d]\ar[r] & M\big\slash N\ar[d]\ar[r] & 0\\
0\ar[r] & \lim_{n} N_{n}\ar[r] & \lim_{n}M_{n}\ar[r] & \lim_{n}M_{n}\big\slash N_{n}\ar[r] & 0
}
\end{displaymath}
The first and second vertical maps are isomorphisms, so the middle one is as well, as required. 
\end{proof}

\begin{cor}
The forgetful functor $|-|:\mathcal{C}_{(A,A_{n},\overline{A}_{n})}\rightarrow{}_{k}\mathrm{Vect}$ commutes with kernels and cokernels. 
\end{cor}



\begin{cor}
Let $(A,A_{n},\overline{A}_{n})$ be a bornological Fr\'{e}chet-Stein algebra such that the underlying algebra $A$ is left coherent as a ring. Then $A$ is strongly left coherent.
\end{cor}

\begin{cor}[c.f. \cite{MR1990669} Proposition 3.7]
Let $(A,A_{n},\overline{A}_{n})$ be a bornological Fr\'{e}chet-Stein algebra.  If $I$ is a closed two-sided ideal of $A$ then $A\big\slash I$ is a bornological Fr\'{e}chet-Stein algebra.
\end{cor}

\begin{proof}
Let $I\subset A$ be a closed two-sided ideal. Both $I$ and $A\big\slash I$ are in $\mathcal{C}_{(A,A_{n},\overline{A}_{n})}$. We have an exact sequence
$$0\rightarrow I\rightarrow A\rightarrow A\big\slash I\rightarrow0$$
and therefore for each $n$ exact sequences
$$0\rightarrow \overline{A}_{n}\hat{\otimes}_{A}I\rightarrow \overline{A}_{n}\rightarrow\overline{A}_{n}\big\slash\overline{A}_{n}I\rightarrow 0$$
where $\overline{A}_{n}I$ is the left ideal of $\overline{A}_{n}$ generated by the image of $I$ under the map $A\rightarrow\overline{A}_{n}$. Now the image of $I$ in $\overline{A}_{n}$ is bornologically dense in $\overline{A}_{n}I$ Since $\overline{A}_{n}I$ is proper it is topologically dense, so its closure is all of $\overline{A}_{n}I$. Thus $\overline{A}_{n}I$ is in fact a two-sided ideal, and the quotient $\overline{A}_{n}\big\slash\overline{A}_{n}I$ is a complete bornological ring. $\overline{A}_{n}\big\slash\overline{A}_{n}I$ is also Noetherian. Similarly $A_{n}I$ is a two-sided ideal in $A_{n}$, and $A_{n}\big\slash A_{n}I$ is a Banach algebra. The quotient of a bornological nuclear Fr\'{e}chet space by a closed subspaces is still a bornological nuclear Fr\'{e}chet space.
\end{proof}

\begin{prop}
    Let $(A,A_{n},\overline{A}_{n})$ be a bornological Fr\'{e}chet Stein algebra and let $I\subseteq A$ be a closed two-sided ideal. Consider the bornological Fr\'{e}chet Stein $(A\big\slash I,A_{n}\big\slash I,\overline{A}_{n}\big\slash I)$. Let $M\in\mathcal{C}_{(A\big\slash I,A_{n}\big\slash I,\overline{A}_{n}\big\slash I)}$. Then regarding $M$ as an $A$-module, $M\in\mathcal{C}_{(A\big\slash I,A_{n},\overline{A}_{n})}$
\end{prop}

\begin{proof}
Write $M\cong\lim_{n}\overline{M}_{n}$ with each $\overline{M}_{n}$ an $\overline{A}_{n}\big\slash I$-module. Then this is also a limit of $\overline{A}_{n}$-modules.  It simply suffices to observe that $$\overline{A}_{n+1}\hat{\otimes}_{A_{n}}\overline{M}_{n+1}\cong(\overline{A}_{n+1}\big\slash I)\hat{\otimes}_{A_{n}}\overline{M}_{n+1}\cong(\overline{A}_{n+1}\big\slash I)\hat{\otimes}_{A_{n}\big\slash I}\overline{M}_{n+1}\cong \overline{M}_{n}$$
\end{proof}

\subsubsection{Derived Fr\'{e}chet-Stein Algebras}

\begin{defn}
An object $A\in\mathbf{DAlg}(\mathrm{Ind(Ban}_{k}\mathrm{)})$ is said to be a \textit{derived bornological Fr\'{e}chet-Stein} algebra if 
\begin{enumerate}
\item
$\pi_{0}(A)$ is a bornological Fr\'{e}chet-Stein algebra.
\item
each $\pi_{n}(A)$ is a coadmissible left $\pi_{0}(A)$-module.
\end{enumerate}
\end{defn}

\begin{lem}\label{lem:frstpres}
An object $A\in\mathbf{DAlg}^{cn}(\mathrm{Ind(Ban}_{k}\mathrm{)})$ is derived bornological Fr\'{e}chet-Stein algebra if and only if there is a sequence
$$A\rightarrow\dots \rightarrow \overline{A}_{n+1}\rightarrow A_{n+1}\rightarrow\overline{A}_{n}\rightarrow A_{n}\rightarrow\dots\rightarrow\overline{A}_{0}\rightarrow A_{0}$$
and equivalences
$$A\cong\mathbb{R}\lim_{\leftarrow_{n\in\mathbb{N}}}A_{n}$$
$$A\cong\mathbb{R}\lim_{\leftarrow_{n\in\mathbb{N}}}\overline{A}_{n}$$
where
\begin{enumerate}
\item
$(\pi_{0}(A),\pi_{0}(A_{n}),\pi_{0}(\overline{A}_{n}))$ is a presentation of $\pi_{0}(A)$ as a bornological Fr\'{e}chet-Stein algebra.
\item
each map $\overline{A}_{n+1}\rightarrow\overline{A}_{n}$ is derived strong.
\end{enumerate}
Moreover the map $A\rightarrow\overline{A}_{n}$ is strong.
\end{lem}

\begin{proof}
Suppose $A$ is a derived bornological Fr\'{e}chet-Stein algebra. 
We shall prove the following claim inductively: for each $n$ there is a sequence
$$A_{\le n}\rightarrow\mathbb{R}\lim_{\leftarrow}(\overline{A}_{k})_{\le n}\rightarrow\cdots\rightarrow (\overline{A}_{m+1})_{\le n}\rightarrow (\overline{A}_{m})_{\le n}\rightarrow\cdots\rightarrow (\overline{A}_{0})_{\le n}$$
For $\pi_{0}(A)=A_{\le 0}$ the result holds by definition. We write
$$\pi_{0}(A)=\mathbb{R}\lim_{\leftarrow}\ldots (\overline{A}_{i})_{\le0}$$
$$\pi_{0}(A)=\mathbb{R}\lim_{\leftarrow}\ldots (A_{i})_{\le0}$$
Now suppose we have proven the claim for some $n$. We have
$$A_{\le n+1}\cong A_{\le n}\oplus_{d_{n}}\Omega\pi_{n+1}(A)[n]$$
Write
$$A_{\le n}\rightarrow\lim_{\leftarrow}(\overline{A}_{k})_{\le n}\rightarrow\cdots\rightarrow (\overline{A}_{m+1})_{\le n}\rightarrow (\overline{A}_{m})_{\le n}\rightarrow\cdots\rightarrow (\overline{A}_{0})_{\le n}$$
$$A_{\le n}\rightarrow\lim_{\leftarrow}(A_{k})_{\le n}\rightarrow\cdots\rightarrow (A_{m+1})_{\le n}\rightarrow (A_{m})_{\le n}\rightarrow\cdots\rightarrow (A_{0})_{\le n}$$
Note that the module $(\overline{A}_{i})_{\le 0}\otimes_{\pi_{0}(A)}\pi_{n+1}(A)$ is a finitely generated $(A_{i})_{\le0}$ module. We get an induced derivation
$$d_{i,n}:(A_{i})_{\le n}\rightarrow  (A_{i})_{\le 0}\otimes_{\pi_{0}(A)}\pi_{n+1}(A)[n+2]$$
Define 
$$(\overline{A}_{i})_{\le n+1}\defeq(\overline{A}_{i})_{\le n}\oplus_{d_{i,n}}((\overline{A}_{i})_{\le 0}\otimes_{\pi_{0}(A)}\pi_{n+2}(A)[n+1])$$
$$(A_{i})_{\le n+1}\defeq (A_{i})_{\le n}\oplus_{d_{i,n}}((A_{i})_{\le 0}\otimes_{\pi_{0}(A)}\pi_{n+2}(A)[n+1])$$
Now define 
$$\overline{A}_{i}\defeq\mathbb{R}\lim_{\leftarrow_{n}}(\overline{A}_{i})_{\le n}$$
$$A_{i}\defeq\mathbb{R}\lim_{\leftarrow_{n}}(A_{i})_{\le n}$$
By construction the map $A\rightarrow\overline{A}_{i}$ is derived strong, and $(\pi_{0}(A),\pi_{0}(A_{n}),\pi_{0}(\overline{A}_{n}))$ is a presentation of $\pi_{0}(A)$ as a bornological Fr\'{e}chet-Stein algebra.

 Conversely suppose $A$ has a presentation 
 $$A\cong\mathbb{R}\lim_{\leftarrow}A_{i}$$
 $$A\cong\mathbb{R}\lim_{\leftarrow_{n\in\mathbb{N}}}\overline{A}_{n}$$
The $\pi_{m}(\overline{A}_{n})$ defines an object of $\mathrm{Coad}_{(\pi_{0}(A),\pi_{0}(A_{n}),\pi_{0}(\overline{A}_{n}))}$.
Thus the map 
 $$\lim_{\leftarrow}\pi_{n}(A_{i})\rightarrow\mathbb{R}\lim_{\leftarrow}\pi_{n}(A_{i})$$
 is an equivalence. A spectral sequence argument then shows that we have
 $$\pi_{n}(\mathbb{R}\lim_{\leftarrow}A_{i})\cong\lim_{\leftarrow}\pi_{n}(A_{i})$$
 In particular $\pi_{n}(\lim_{\leftarrow}A_{i})\in\mathcal{C}_{(\pi_{0}(A),\pi_{0}(A_{n}),\pi_{0}(\overline{A}_{n}))}$.

 Note the map $A\rightarrow\overline{A}_{n}$ is strong by construction.
\end{proof}

\section{Bornological Lawvere Theories}\label{sec:bornlawvere}

We know introduce the various Lawvere theories which will determine our geometries. They will all be generated by algebras in the sense of Subsubsection \ref{sec:Lawveregenalg}. First we give some general ways to check certain properties of these Lawvere theories. 

\subsection{Rigidification}

Here we do not distinguish between Archimedean and non-Archimedean, and suppress the `$nA$' decoration.
Let $\mathrm{R}$ be a Banach ring which has the open mapping property/ Buchwalter's prroperty. Let $\mathrm{T}$ be a $\underline{\Lambda}$-sorted Lawvere theory concretely of homotopy $\mathrm{Ind(Ban}_{R}\mathrm{)}$-polynomial type. Suppose that each $\mathrm{F}(\mathrm{Free_{T}}(\underline{\lambda}))$ is a Fr\'{e}chet $R$-module/ has a bornology with countable basis for each $\underline{\lambda}$. We claim that the functor
$$\mathrm{F}:\mathrm{sAlg_{T}}\rightarrow\mathrm{sComm}(\mathrm{Ind(Ban}_{R}\mathrm{)})$$
(resp.
$$\mathrm{F}:\mathrm{sAlg_{T}}\rightarrow\mathrm{sComm}(\mathrm{Ind(Ban}^{nA}_{R}\mathrm{)}))$$
preserves weak equivalences between fibrant-cofibrant objects.

The proof is essentially the same as \cite{borisov2017quasi} Lemma 10. First observe that filtered colimits in $\mathrm{sComm}(\mathrm{Ind(Ban}_{R}\mathrm{)})$ preserve weak equivalences.  Moreover $\mathrm{F}$ commutes with filtered colimits. It follows that we may assume that $A_{\bullet}$ is free and finitely generated in each degree. By Lemma \ref{lem:funcsimp} it suffices to prove that the map
$$\mathrm{F}(A_{\bullet}\otimes\Delta[1])\rightarrow\mathrm{F}(A_{\bullet})$$
is an equivalence. The kernel is itself a simplicial Fr\'{e}chet space. It is acyclic if and only if it is algebraically exact. Thus the map is a weak equivalence if and only if it is an algebraic equivalence. But each
$$\underline{\mathrm{Hom}}(\mathrm{Free_{T}}(\lambda),A_{\bullet}\otimes\Delta[1])\rightarrow\underline{\mathrm{Hom}}(\mathrm{Free_{T}}(\lambda),A_{\bullet})$$
is an equivalence sin $\mathrm{Free_{T}}(\lambda)$ is cofibrant in $\mathrm{sAlg_{T}}$. Taking the colimit over $\lambda$, and using the fact that $\mathrm{T}$ is concretely of homotopy $\mathrm{Ind(Ban}_{R}\mathrm{)}$-polynomial type, gives that 
$$\mathrm{F}(A_{\bullet}\otimes\Delta[1])\rightarrow\mathrm{F}(A_{\bullet})$$
is algebraically exact, as required.

\subsection{Homotopy Epimorphisms in Bornological Fermat Theories}

Let $\mathrm{T}$ be a $\Gamma$-filtered Fermat theory, and let \[F:\mathrm{T}^{op}\rightarrow\mathrm{Comm}(\mathrm{CBorn}_{R})\] be a fully faithful finite coproduct-preserving functor, where $R$ is a Banach ring. Let $\gamma\in\Gamma$, and consider the map $\mathrm{Sym}(R)\rightarrow \mathrm{F}(\mathrm{Free_{T}}(\gamma))$. Suppose that $\mathrm{F}(\mathrm{Free_{T}}(\gamma))$ is flat as an object of $\mathrm{Comm}(\mathrm{CBorn}_{R})$.
Consider the complex
$$(x-y)\times:\mathrm{F}(\mathrm{Free_{T}}(\gamma,\gamma))\rightarrow \mathrm{F}(\mathrm{Free_{T}}(\gamma,\gamma))$$
Note that since $\mathrm{T}$ is Fermat, the multiplication by $(x-y)$ map is automatically injective.
Suppose that the underlying algebra $\mathrm{Hom}(\mathrm{Sym}(R),\mathrm{F}(\mathrm{Free_{T}}(\gamma,\gamma)))\cong \mathrm{Free_{T}}(\gamma)^{ind}$. The map
$\pi:\mathrm{Free_{T}}(\gamma,\gamma)\rightarrow\mathrm{Free_{T}}(\gamma)$ is a split epimorphism, with splitting $s_{y}$ corresponding to the projection in the second factor. Thus $\mathrm{F}(\pi):\mathrm{F}(\mathrm{Free_{T}}(\gamma,\gamma))\rightarrow\mathrm{F}(\mathrm{Free_{T}}(\gamma))$ is also a split epimorphism. Let $f(x,y)\in\mathrm{Free_{T}}(\gamma,\gamma)(\delta)$ and consider $f(x,y)-f(y,y)\in\mathrm{Free_{T}}(\gamma,\gamma)(\delta)$. Note that $f(y,y)=s_{y}\circ\pi(f(x,y))$. Since $\mathrm{T}$ is a Fermat theory, there is a unique $g(x,y,z)\in\mathrm{Free_{T}}(\lambda,\lambda,\lambda)(\frac{\delta}{\lambda})$ such that 
$$f(x,y)-f(y,y)=(x-y)g(x,y,y)$$

Finally, suppose that 
\begin{enumerate}
    \item 
    $R$ satisfies the open mapping property (resp. Buchwalter's property) and each $\mathrm{F}(\mathrm{Free_{T}}(\underline{\lambda}))$ is a Fr\'{e}chet space (resp. has a bornology with a countable basis).
    \item 
    The induced assignment $\mathrm{Free}_{T}(\gamma)^{ind}\rightarrow\mathrm{Free}_{T}(\gamma,\gamma)^{ind}, f(x,y)\mapsto g(x,y,y)$ is bounded for the bornology determined by $\mathrm{F}$.
\end{enumerate}
Then the sequence 
$$\mathrm{F}(\mathrm{T}(\gamma,\gamma))\rightarrow \mathrm{F}(\mathrm{T}(\gamma,\gamma))\rightarrow\mathrm{F}(\mathrm{T}(\gamma))$$
is short exact. This implies then that $\mathrm{T}$ is of homotopy ind-$\mathrm{CBorn}_{R}$-polynomial type. Indeed we have a short exact sequence
$$\mathrm{Sym}(R\oplus R)\rightarrow\mathrm{Sym}(R\oplus R)\rightarrow\mathrm{R}$$
where the first map is again multiplication by $(x-y)$. Thus $\mathrm{Sym}(R\oplus R)\rightarrow\mathrm{Sym}(R\oplus R)$ is a flat resolution of $\mathrm{Sym}(R)$. Using that $\mathrm{Free_{T}}(\gamma)$ is flat, tensoring both sides by $\mathrm{F}\circ\mathrm{Free_{T}}(\gamma)$, we find that the complex

$$(x-y)\times:\mathrm{F}(\mathrm{Free_{T}})(\gamma,\gamma)\rightarrow \mathrm{F}(\mathrm{Free_{T}}(\gamma,\gamma))$$
computes $\mathrm{F}(\mathrm{Free_{T}}(\gamma))\otimes^{\mathbb{L}}_{\mathrm{Sym}(R)}\mathrm{F}(\mathrm{Free_{T}}(\gamma))$. But this is just $\mathrm{F}(\mathrm{Free_{T}}(\gamma))$ as required. 

\subsection{Bornological Algebras}

The Lawvere theories in which we are interested are generated by algebras (see Section \ref{sec:Lawveregenalg}). Fix a Banach ring $R$.

\subsubsection{Polynomial and Power Series Algebras}

\begin{defn}
Let $E$ be 
\begin{enumerate}
\item
The \textit{bornological polynomial algebra in }$n$-\textit{variables} is the symmetric algebra $\mathrm{Sym}(R^{\oplus n})=\coprod_{m=0}^{\infty}\mathrm{Sym}^{m}(R^{\oplus n})$ computed internal to $\mathrm{Ind(Ban}_{R}\mathrm{)}$.
\item
The \textit{bornological power series algebra in }$n$-\textit{variables} is $\widehat{\mathrm{Sym}(R^{\oplus n})}=\prod_{m=0}^{\infty}\mathrm{Sym}^{m}(R^{\oplus n})$ computed internal to $\mathrm{Ind(Ban}_{R}\mathrm{)}$.
\end{enumerate}
\end{defn}

Consider the Lawvere theory $\mathrm{FPS}_{R}\defeq\mathrm{FPS}_{\mathrm{CBorn}_{R}}$. It turns out that in this case $\mathrm{FPS}_{R}$ is generated by the algebra 
$$\mathrm{FPS}_{R}(1)=\widehat{\mathrm{Sym}}(R)$$
Indeed $\widehat{\mathrm{Sym}}(R^{\oplus m})$ is a Fr\'{e}chet $R$-module, so by \cite{ben2020fr} $\widehat{\mathrm{Sym}}(R^{\oplus m)})\hat{\otimes}_{\pi,R}(-)$ commutes with countable products. Thus
\begin{align*}
\widehat{\mathrm{Sym}}(R^{\oplus m})\hat{\otimes}_{\pi,R}\widehat{\mathrm{Sym}}(R^{\oplus n})&\cong\prod_{p,q}\mathrm{Sym}^{p}(R^{\oplus m})\hat{\otimes}_{\pi,R}\mathrm{Sym}^{q}(R^{\oplus n})\\
&\cong\prod_{k}\mathrm{Sym}^{k}(R^{\oplus m}\oplus R^{\oplus n})\\
&\cong\widehat{\mathrm{Sym}}(R^{\oplus (m+n)})
\end{align*}
as required.

 $\mathrm{FPS}_{R}$ is a Fermat theory. Let 
$$f(x,z_{1},\ldots,z_{n})=\sum_{i,j_{1},\ldots,j_{n}}a_{i,j_{1},\ldots,j_{n}}x^{i}z_{1}^{j_{1}}\ldots z_{n}^{j_{n}}\in\mathrm{FPS}_{R}(n+1)$$
 Write 
 $$g(x,y,z_{1},\ldots,z_{n})=\sum_{i,j_{1},\ldots,j_{n}}a_{i,j_{1},\ldots,j_{n}}\sum_{m=0}^{i-1}x^{i}y^{i-1-m}z_{1}^{j_{1}}\ldots z_{n}^{j_{n}}$$
 A straightforward computation shows that $g(x,y,z_{1},\ldots,z_{n})$ gives the required function, and that the assigment $f\mapsto g$ is bounded.

\begin{lem}
The map
$$\mathrm{Sym}(R^{\oplus n})\rightarrow\widehat{\mathrm{Sym}(R^{\oplus n})}$$
is a homotopy epimorphism. In particular $\mathrm{FPS}_{R}$ is a Lawvere theory of homotopy $\mathrm{LH(}\mathrm{Ind(Ban}_{R}\mathrm{))}$-polynomial type. In fact it is concretely of $\mathrm{LH(}\mathrm{Ind(Ban}_{R}\mathrm{))}$-polynomial type. 
\end{lem}

\begin{proof}
By \cite{ben2021analytification} Lemma 4.8 it is of homotopy polynomial type. The fact that it is concretely of homotopy polynomial type simply follows from the fact that $\mathrm{Hom}(R,-)$ commutes with products.
\end{proof}

\begin{cor}
The functor
$$\mathbf{sAlg}_{\mathrm{FPS}_{R}}\rightarrow\mathbf{DAlg}^{cn}(\mathrm{CBorn}_{R})$$
is fully faithful.
\end{cor}

Moreover we have the following.

\begin{prop}
    If $R=k$ is a non-trivially valued Banach field of characteristic $0$ then $\mathbf{FPS}_{k}$ is a Fermat theory admitting integration.
\end{prop}

\begin{proof}
    Most of this is clear. Integration is defined formally in the obvious way.
\end{proof}

\subsection{`Analytic' Algebras}

We now define some classes of algebras which sit strictly between polynomials and power series. roughly speaking, these should be functions on analytic neighbourhoods, possibly with boundary. For a Banach ring $R$ whose multiplication map is of semi-norm $C$ we denote by $\underline{R}$ the $\mathbb{R}_{>0}$-filtered ring with
$$F_{\lambda}R=\{r\in R:|r|\le\frac{\lambda}{C}\}$$

\subsubsection{The Contracting Algebras and Affinoid Algebras}

\begin{defn}
Let $E$ be a Banach $R$-module
\begin{enumerate}
\item
Define
$$R\{E\}\defeq\mathrm{Sym}^{\le1}(E)$$
to be the free symmetric algebra on $E$ commuted in the symmetric monoidal category $\mathrm{Ban}_{R}^{\le1}$.
\item
If $R$ and $E$ are both non-Archimedean define
$$R\Bigr<E\Bigr>\defeq\mathrm{Sym}^{nA,\le1}(E)$$
to be the free symmetric algebra on $E$ commuted in the symmetric monoidal category $\mathrm{Ban}_{R}^{nA,\le1}$.
\end{enumerate}
\end{defn}

\begin{rem}\label{rem:afndtensorsum}
If $E$ and $F$ are Banach $R$-modules then we have
$$R\{E\}\hat{\otimes}_{\pi,R}R\{F\}\cong R\{E\oplus F\}$$
and in the non-Archimedean case,
$$R\Bigr<E\Bigr>\hat{\otimes}_{\pi,R}R\Bigr<F\Bigr>\cong R\Bigr<E\oplus F\Bigr>$$
\end{rem}

\begin{defn}
\begin{enumerate}
Let $R$ be a Banach ring and $r=(r_{1},\ldots,r_{n})\in\mathbb{R}_{>0}^{n}$ a polyradius. Define
\item
$$R\{\frac{x_{1}}{r_{1}},\cdots,\frac{x_{n}}{r_{n}}\}\defeq R\{R_{r_{1}}\oplus\ldots \oplus R_{r_{n}}\}$$
This can be equivalently described as the algebra
 $$\{\sum_{I\in\mathbb{N}^{n}}a_{I}X^{I}:\sum_{i\in\mathbb{N}^{I}}|a_{I}|r^{I}<\infty\}$$
with norm  $||\sum_{I\in\mathbb{N}^{n}}a_{I}X^{I}||\defeq \sum_{I\in\mathbb{N}^{n}}|a_{I}|r^{I}$.
\item
If $R$ is non-Archimedean define
$$R\Big<\frac{x_{1}}{r_{1}},\ldots,\frac{x_{n}}{r_{n}}\Bigr>\defeq R\Big<R_{r_{1}}\oplus\cdots \oplus R_{r_{n}}\Bigr>$$
$$\{\sum_{I\in\mathbb{N}^{n}}a_{I}X^{I}:\lim_{I\rightarrow\infty}|a_{I}|r^{I}=0\}$$
with norm  $||\sum_{I\in\mathbb{N}^{n}}a_{I}X^{I}||\defeq \underset{I\in\mathbb{N}^{n}}\sup |a_{I}|r^{I}$
\end{enumerate}
\end{defn}

In the non-Archimedean case we will also use the notation $T_{R}^{n}(r)=R\Big<\frac{x_{1}}{r_{1}},\ldots,\frac{x_{n}}{r_{n}}\Bigr>$. We call this the \textit{Tate algebra of polyradius }$r$. 
More generally, an \textit{affinoid algebra} is an algebra of the form $T^{n}_{R}(r)\big\slash I$, for some $n$ and $r$, where $I$ is a finitely generated closed ideal.

\begin{defn}
For $R$ non-Archimedean the $R$-\textit{Tate Lawvere theory} is the $\mathbb{R}_{>0}$-sorted $\underline{R}$-Lawvere theory $\mathrm{Tate}_{R}$, defined by
$$\mathrm{Hom}_{\mathrm{Tate}_{R}}(\underline{\lambda},\underline{\gamma})\defeq\mathrm{Hom}_{\mathrm{Comm(CBorn^{nA}_{R})}}(T^{n}_{R}(\underline{\gamma}),T^{m}_{R}(\underline{\lambda}))$$
\end{defn}

Thanks to Remark \ref{rem:afndtensorsum} this Lawvere theory is generated by algebras.

\begin{rem}
One can also define an Archimedean `Tate' Lawvere theory using the algebras $R\{\frac{x_{1}}{r_{1}},\ldots,\frac{x_{n}}{r_{n}}\}$. However apart from briefly as an ancillary tool in this work we will not have any use for this. In fact homotopy-theoretically it is not very well-behaved. We will denote this by $\mathrm{Tate}^{arch}_{R}$. 
\end{rem}

Let us now show that $\mathrm{Tate}_{R}$ is an $\underline{\mathbb{R}_{>0}}$-filtered Fermat theory in the non-Archimedean case. Let 
$$f(x,z_{1},\ldots,z_{n})=\sum_{i,j_{1},\ldots,j_{n}}a_{i,j_{1},\ldots,j_{n}}x^{i}z_{1}^{j_{1}}\ldots z_{n}^{j_{n}}\in\mathrm{Tate}_{R}(\lambda,\gamma_{1},\ldots,\gamma_{n})(\delta)$$
 Write 
 $$g(x,y,z_{1},\ldots,z_{n})=\sum_{i,j_{1},\ldots,j_{n}}a_{i,j_{1},\ldots,j_{n}}\sum_{m=0}^{i-1}x^{i}y^{i-1-m}z_{1}^{j_{1}}\ldots z_{n}^{j_{n}}$$
 A straightforward computation shows that $g(x,y,z_{1},\ldots,z_{n})$ gives the required function, and that the assigment $f\mapsto g$ is bounded.

By \cite{ben2021analytification} Corollary 5.6 we get the following.

\begin{lem}
$\mathrm{Tate}_{R}$ is of homotopy $\mathrm{Ind(Ban^{nA}_{R})}$-polynomial type, and concretely of homotopy $\mathrm{Ind(Ban^{nA}_{R})}$-polynomial type.
\end{lem}

\begin{proof}
By \cite{ben2021analytification} Corollary 5.6 it is of homotopy polynomial type. Now $\mathrm{Hom}(\mathrm{Sym}(R),R<\lambda_{1}^{-1}x_{1},\ldots,\lambda_{n}^{-1}x_{n}>)$ is the underlying algebra of the Banach algebra $R<\lambda_{1}^{-1}x_{1},\ldots,\lambda_{n}^{-1}x_{n}>$. On the other hand, any element $f\in R<\lambda_{1}^{-1}x_{1},\ldots,\lambda_{n}^{-1}x_{n}>$ induces a map $R<\gamma>\rightarrow R<\lambda_{1}^{-1}x_{1},\ldots,\lambda_{n}^{-1}x_{n}>$ sending the coordinate function $y$ to $f$, as long as $\gamma>||f||$. Thus $\limind_{\gamma}\mathrm{Hom}(R<\gamma>,R<\lambda_{1}^{-1}x_{1},\ldots,\lambda_{n}^{-1}x_{n}>)$ is also the underlying algebra of $R<\lambda_{1}^{-1}x_{1},\ldots,\lambda_{n}^{-1}x_{n}>$.
\end{proof}

Moreover we have the following.

\begin{prop}
    If $R=k$ is a non-trivially valued non-Archimedean Banach field of characteristic $0$ then $\mathrm{Tate}_{k}$ is a strict $\mathbb{R}_{\ge0}$-sorted $\underline{k}$-Fermat theory admitting integration. Moreover it is strongly $\mathrm{Ind(Ban}_{k}\mathrm{)}$-coherent.
\end{prop}

The reason this theory admits integration while the Tate one doesn't is as follows. One can try and define integration formally on Tate algebras. However the resulting power series may not converge at the boundary. Using the standard formula for the radius of convergence of a non-Archimedean power series however does show that formally integrating does not change the radius of convergence.

\subsubsection{Dagger Affinoid Algebras}

Let $E$ be a Banach $R$-module. 

\begin{defn}
\begin{enumerate}
\item
Let $R$ be any Banach ring. If $\rho=(\rho_{1},\ldots,\rho_{n})\in\mathbb{R}_{>0}^{n}$ is a polyradius, we define 
$$W^{n}_{R}(\rho)=R\{\rho_{1}^{-1}X_{1},\ldots,\rho_{n}^{-1}X_{n}\}^{\dagger}=``\limind_{r>\rho}"R\{r_{1}^{-1}X_{1},\ldots,r^{-1}X_{n}\}$$
$W^{n}_{R}(\rho)$ is called a \textit{strict dagger affinoid}. A \textit{dagger affinoid algebra} is an algebra of the form $W^{n}_{R}(r)\big\slash I$, for some $n$ and $r$, where $I$ is a finitely generated closed ideal. 
\item
Let $R$ be non-Archimedean. If $\rho=(\rho_{1},\ldots,\rho_{n})\in\mathbb{R}_{>0}^{n}$ is a polyradius, we define 
$$W^{nA,n}_{R}(\rho)=R\Bigr<\rho_{1}^{-1}X_{1},\ldots,\rho_{n}^{-1}X_{n}\Bigr>^{\dagger}=``\limind_{r>\rho}"T^{n}_{R}(r)$$
$W^{n}_{R}(\rho)$ is called a \textit{strict non-Archimedean dagger affinoid}. A \textit{non-Archimedean dagger affinoid algebra} is an algebra of the form $W^{nA,n}_{R}(r)\big\slash I$, for some $n$ and $r$, where $I$ is a finitely generated closed ideal. 
\end{enumerate}
\end{defn}

By the results for Tate algebras, and the fact that tensor products commute with colimits we have 
$$W^{m}_{R}(\delta)\hat{\otimes}_{\pi,R}W^{n}_{R}(\rho)\cong W^{m+n}_{R}(\delta,\rho)$$
and if $R$ is non-Archimedean
$$W^{nA,m}_{R}(\delta)\hat{\otimes}^{nA}_{\pi,R}W^{nA,n}_{R}(\rho)\cong W^{nA,m+n}_{R}(\delta,\rho)$$

Thus they give rise to Lawvere theories generated by algebras in $\mathrm{Ind(Ban}_{R}\mathrm{)}$ and in $\mathrm{Ind(Ban}^{nA}_{R}\mathrm{)}$ respectively. We denote these Lawvere theories by $\mathbf{Afnd}^{\dagger}_{R}$ and $\mathbf{Afnd}^{nA,\dagger}$ respectively.  Some of the geometry in this setting was done in \cite{MR3448274}. 
 
 \begin{lem}
$\mathrm{Afnd}^{\dagger}_{R}$ is of homotopy $\mathrm{Ind(Ban_{R})}$-polynomial type, and concretely of homotopy $\mathrm{Ind(Ban_{R})}$-polynomial type. For $R$ non-Archimedean $\mathrm{Afnd}^{nA,\dagger}_{R}$ is of homotopy $\mathrm{Ind(Ban^{nA}_{R})}$-polynomial type, and concretely of homotopy $\mathrm{Ind(Ban^{nA}_{R})}$-polynomial type. 
\end{lem}

\begin{proof}
    The fact that it is of homotopy $\mathrm{Ind(Ban_{R})}$-polynomial type is proved in \cite{ben2021analytification} Lemma 4.14. The fact that it is concretely of homotopy polynomial type is similar to the affinoid case.
\end{proof}

Moreover we have the following.

\begin{prop}
    If $R=k$ is a non-trivially valued Banach field of characteristic $0$ then $\mathrm{Afnd}^{\dagger}_{R}$ is a Fermat theory. Moreover it is strongly $\mathrm{Ind(Ban}_{k}\mathrm{)}$-coherent.
\end{prop}

%

\subsubsection{Disc or Dagger Stein Algebras}

\begin{defn}
\begin{enumerate}
\item
Let $R$ be any Banach ring. If $\rho=(\rho_{1},\ldots,\rho_{n})\in(\mathbb{R}_{>0},\cup\{\infty\})^{n}$ is a polyradius, we define 
$$\mathcal{O}(D^{n}_{<\rho,R})=\lim_{r<\rho}R\Bigr\{r_{1}^{-1}X_{1},\ldots,r_{n}^{-1}X_{n}\Bigr\}$$
$\mathcal{O}(D^{n}_{<\rho,R})$ is called a \textit{disc algebra}. 
\item
Let $R$ be a non-Archimedean Banach ring. If $\rho=(\rho_{1},\ldots,\rho_{n})\in(\mathbb{R}_{>0},\cup\{\infty\})^{n}$ is a polyradius, we define 
$$\mathcal{O}^{nA}(D^{n}_{<\rho,R})=\lim_{r<\rho}R\Bigr<r_{1}^{-1}X_{1},\ldots,r_{n}^{-1}X_{n}\Bigr>$$
$\mathcal{O}^{nA}(D^{n}_{<\rho,R})$ is called a \textit{non-Archimedean disc algebra}. 
\end{enumerate}
\end{defn}

By \cite{ben2020fr} Lemma 6.10 we have 
$$\mathcal{O}(D^{m}_{<\delta,R})\hat{\otimes}_{\pi,R}\mathcal{O}(D^{n}_{<\rho,R})\cong \mathcal{O}(D^{m+n}_{<(\delta,\rho),R})$$
and if $R$ is non-Archimedean
$$\mathcal{O}^{nA}(D^{m}_{<\delta,R})\hat{\otimes}_{\pi,R}\mathcal{O}^{nA}(D^{n}_{<\rho,R})\cong \mathcal{O}^{nA}(D^{m+n}_{<(\delta,\rho),R})$$

Thus they give rise to Lawvere theories generated by algebras in $\mathrm{Ind(Ban}_{R}\mathrm{)}$ and in $\mathrm{Ind(Ban}^{nA}_{R}\mathrm{)}$ respectively .We denote these Lawvere theories by $\mathrm{Disc}_{R}$ and $\mathrm{Disc}_{R}^{nA}$ respectively. As before we get the following. 

 \begin{lem}
$\mathrm{Disc}_{R}$ is of homotopy $\mathrm{Ind(Ban_{R})}$-polynomial type, and concretely of homotopy $\mathrm{Ind(Ban_{R})}$-polynomial type. For $R$ non-Archimedean $\mathrm{Disc}^{nA}_{R}$ is of homotopy $\mathrm{Ind(Ban_{R})}$-polynomial type, and concretely of homotopy $\mathrm{Ind(Ban_{R})}$-polynomial type. 
\end{lem}

\begin{prop}
    If $R=k$ is a non-trivially valued Banach field of characteristic $0$, either Archimedean or non-Archimedean then $\mathrm{Disc}_{k}$ is a Fermat theory admitting integration. 
\end{prop}

\begin{lem}\label{lem:DiscFrSTn}
Let $k$ be a non-trivially valued Banach field, and let $A$ be a discrete finitely $\mathrm{Disc}_{k}$-presented algebra (resp. a finitely $\mathrm{Disc}^{nA}_{k}$-presented algebra). Then $A$ is a bornological Fr\'{e}chet Stein algebra. 
\end{lem}

\begin{proof}
We prove the Archimedean case, the non-Archimedean case being identical. It suffices to prove the claim for $A=\mathcal{O}(D^{n}_{k,<\lambda})$, since quotients of bornological Fr\'{e}chet-Stein algebras are bornological Fr\'{e}chet-Stein algebras. Let $(\rho^{m})$ be an increasing sequence converging to $\lambda$. We may write
\begin{equation}
    \begin{split}
\mathcal{O}(D^{n}_{k,<\lambda}) & \cong \lim_{m}k<(\rho^{m}_{1})^{-1}x_{1},\ldots,(\rho^{m}_{n})^{-1}x_{n}>^{\dagger} \\ & \cong\lim_{m}k<(\rho^{m}_{1})^{-1}x_{1},\ldots,(\rho^{m}_{n})^{-1}x_{n}>
\end{split}
\end{equation}
Moreover there are natural dense maps

\begin{equation}
    \begin{split}
k<(\rho^{m}_{1})^{-1}x_{1},\ldots,(\rho^{m}_{n})^{-1}x_{n}>^{\dagger} & \rightarrow k<(\rho^{m}_{1})^{-1}x_{1},\ldots,(\rho^{m}_{n})^{-1}x_{n}> \\ & \rightarrow k<(\rho^{m-1}_{1})^{-1}x_{1},\ldots,(\rho^{m-1}_{n})^{-1}x_{n}>^{\dagger} \\ & \rightarrow k<(\rho^{m-1}_{1})^{-1}x_{1},\ldots,(\rho^{m-1}_{n})^{-1}x_{n}>
\end{split}
\end{equation}
The algebras $k<(\rho^{m}_{1})^{-1}x_{1},\ldots,(\rho^{m}_{n})^{-1}x_{n}>^{\dagger}$ are strongly Noetherian. This gives the presentation as a bornological Fr\'{e}chet-Stein algebra.
\end{proof}

This brings us to the dagger Stein algebras of \cite{bambozzi2015stein}

\subsubsection{Algebras of Entire Functions}

Let $\mathbf{EFC}_{R}$ denote the Lawvere theory generated by $\mathcal{O}(R)\defeq\mathcal{O}(D^{1}_{<\infty,R})$. If $R$ is non-Archimedean let $\mathbf{EFC}^{nA}_{R}$ denote the Lawvere theory generated by $\mathcal{O}^{nA}(R)\defeq\mathcal{O}^{nA}(D^{1}_{<\infty,R})$. These are $1$-sorted Lawvere theories of homotopy $\mathrm{Ind(Ban}_{R}\mathrm{)}$-polynomial type and homotopy $\mathrm{Ind(Ban}^{nA}_{R}\mathrm{)}$-polynomial type respectively (again by \cite{ben2021analytification} Section 5). Moreover they are in fact of concrete polynomial type. In particular, we have the following.

\begin{lem}\label{lem:embedEFC}
The functor
$$\mathbf{sAlg}_{\mathbf{EFC}_{R}}\rightarrow\mathbf{DAlg}^{cn}(\mathrm{Ind(Ban}_{R}\mathrm{)})$$
is fully faithful. If $R$ is non-Archimedean then the functor
$$\mathbf{sAlg}_{\mathbf{EFC}^{nA}_{R}}\rightarrow\mathbf{DAlg}^{cn}(\mathrm{Ind(Ban}^{nA}_{R}\mathrm{)})$$
is fully faithful.
\end{lem}

$\mathbf{EFC}$ here stands for `Entire Functional Calculus' as in \cite{MR4036665}.

\begin{prop}
    If $R=k$ is a non-trivially valued Banach field of characteristic $0$, either Archimedean or non-Archimedean then $\mathrm{EFC}_{k}$ is a Fermat theory admitting integration. 
\end{prop}


\subsubsection{Germs}

For $n\in\mathbb{N}_{0}$ let $\mathrm{Germ}_{R}$ be the $1$-sorted Lawvere theory generated by $\mathrm{Germ}_{R}(n)\defeq\limind_{r\in\mathbb{R}^{n}_{>0}}R\{\frac{x_{1}}{r_{1}},\ldots,\frac{x_{n}}{r_{n}}\}$ and, in the non-Archimedean case, let $\mathrm{Germ}_{R}^{nA}$ be the Lawvere theory generated by $\limind_{r\in\mathbb{R}^{n}_{>0}}R\Bigr<\frac{x_{1}}{r_{1}},\ldots,\frac{x_{n}}{r_{n}}\Bigr>$. Since colimits of homotopy epimorphisms are homotopy epimorphisms we get the following.

\begin{lem}
$\mathrm{Germ}_{R}$ is of homotopy $\mathrm{Ind(Ban}_{R}\mathrm{)}$-polynomial type, and in the non-Archimedean setting,  $\mathrm{Germ}^{nA}_{R}$ is of homotopy $\mathrm{Ind(Ban^{nA}}_{R}\mathrm{)}$-polynomial type. Moreover it is strongly $\mathrm{Ind(Ban}_{k}\mathrm{)}$-coherent, is Fermat, and admits integration. 
\end{lem}

\subsubsection{$C^{\infty}$-Algebras}

As our final example we consider algebras of smooth functions following work from \cite{borisov2017quasi}. Here we fix $R$ to be the Banach field $\mathbb{R}$ of real numbers with the usual Euclidean norm. For each $m,n\in\mathbb{N}$ let $C^{\infty}(\mathbb{R}^{m},\mathbb{R}^{n})$ denote the space of smooth function from $\mathbb{R}^{m}$ to $\mathbb{R}^{n}$. $C^{\infty}(\mathbb{R}^{m},\mathbb{R}^n)$ can be endowed with a locally convex topology with basis of open sets parametrised by $r,\epsilon\in\mathbb{R}_{>0},q\in\mathbb{R}^{m},n\in\mathbb{Z}_{\ge0}$:
$$B_{q,r}^{\epsilon,m}\defeq\{f\in\mathcal{C}^{\infty}(\mathbb{R}^{m},\mathbb{R}):\mathrm{sup}_{||p-q||\le r}||J^{n}(f)(p)||<\epsilon\}$$
where $||J^{m}(f)(p)||\defeq\sum_{||\underline{k}||\le m}||D^{\underline{k}}(f)(p)||$. 
$$C^{\infty}(\mathbb{R}^{m}\mathbb{R}^{n})\cong C^{\infty}(\mathbb{R}^{m}\mathbb{R})^{n}$$
can then be endowed with the product locally convex topology. By equipping this with the associated pre-compact bornology, this becomes a bornological Fr\'{e}chet space over $\mathbb{R}$.

\begin{defn}\label{defn:lawvcart}
The $C^{\infty}$ \textit{Lawvere theory}, denoted $\mathrm{CartSp}_{\mathrm{smooth}}$ is the $1$-sorted Lawvere theory with 
$$\mathrm{Hom}_{\mathrm{CartSp}_{\mathrm{smooth}}}(m,n)\defeq\mathrm{Hom}_{C^{\infty}\mathrm{Man}}(\mathbb{R}^{m},\mathbb{R}^{n})$$
where $C^{\infty}\mathrm{Man}$ is the category of smooth manifolds and smooth functions, and $\mathbb{R}^{k}$ is equipped with its usual smooth structure. 
\end{defn}

The notation $\mathrm{CartSp}_{\mathrm{smooth}}$ is from \cite{spivak2010derived}.

\begin{lem}
The map
$$\mathrm{Sym}(\mathbb{R}^{n})\rightarrow\mathrm{C}^{\infty}(\mathbb{R}^{ n},\mathbb{R})$$
in $\mathrm{Ind(Ban_{\mathbb{R}})}$ is a homotopy epimorphism. In particular the functor
$$\mathbf{sAlg}_{\mathrm{CartSp}_{\mathrm{smooth}}}\rightarrow\mathbf{DAlg}^{cn}(\mathrm{Ind(Ban}^{nA}_{\mathbb{R}}\mathrm{)}$$
is fully faithful.
\end{lem}

\begin{proof}
First observe that polynomials are bornologically dense in $\mathrm{C}^{\infty}(\mathbb{R}^{ n},\mathbb{R})$. Thus the map is an epimorphism. Consider the map $\mathrm{C}^{\infty}(\mathbb{R}^{ 2},\mathbb{R})\rightarrow \mathrm{C}^{\infty}(\mathbb{R}^{ 2},\mathbb{R})$ which sends $f$ to
$$\frac{f(x,y)-f(y,y)}{x-y}$$
This is clearly well-defined as a map of algebras, and is easily seen to be continuous. This gives a splitting of the Koszul complex.
\end{proof}

Although formulated and proved in a slightly different fashion full faithfulness of $\mathbf{sAlg}_{\mathrm{CartSp}_{\mathrm{smooth}}}\rightarrow\mathbf{DAlg}^{cn}(\mathrm{Ind(Ban}^{nA}_{\mathbb{R}}\mathrm{)}$ was previously established in \cite{borisov2017quasi}.

The following is proven in \cite{carchedi2012homological}.

\begin{prop}
$\mathrm{CartSp}$ is a Fermat theory admitting integration. 
\end{prop}

\subsubsection{Exactness of the Functors $\mathbf{F}_{B}$}

\begin{prop}
    Let $k$ be a non-trivially valued Banach field, and let $\mathrm{T}\in\{\mathrm{FPS}_{k},\mathrm{Disc}_{k},\mathrm{Afnd}^{\dagger}_{k}\}$ for $k$ Archimedean, and $\mathrm{T}\in\{\mathrm{FPS}_{k},\mathrm{Disc}^{nA}_{k},\mathrm{Afnd}^{nA,\dagger}_{k},\mathrm{Tate}_{k}\}$ in the case that $k$ is non-Archimedean. Let $B\in\mathbf{sAlg}_{\mathrm{T}}$ be such that $\pi_{0}(B)$ is of the form $\mathrm{Free_{T}}(\underline{\lambda})\big\slash I$ with $I$ finitely generated, and such that $I$ is closed in $\mathrm{F}(\mathrm{Free_{T}}(\underline{\lambda}))$. Then the functor
    $$\mathbf{F}_{B}:{}_{B}\mathbf{Mod}^{alg}\rightarrow{}_{\mathbf{F}(B)}\mathbf{Mod}(\mathrm{Ind(Ban}_{k}\mathrm{)}$$
    is $t$-exact.
\end{prop}

\begin{proof}
    This is a consequence of the open mapping theorem and Buchwalter's Theorem.
\end{proof}

By Theorem \ref{thm:essimageF} we then get the following.

\begin{cor}
    For $\mathrm{T}\in\{\mathrm{FPS}_{k},\mathrm{EFC}_{k}\}$ in both the Archimedean and non-Archimedean cases, the essential image of the functor
    $$\mathbf{F}:\mathbf{sAlg}^{\pi_{0}-gf}_{\mathrm{T}}\rightarrow\mathbf{DAlg}^{cn}(\mathrm{Ind(Ban}_{k}\mathrm{)})$$
    consists of those algebras $A$ such that $\pi_{0}(A)$ is of the form $\mathrm{F}(\mathrm{Free_{T}}(\underline{\lambda})\big\slash I)$ for some closed finitely generated ideal $I$, and $\pi_{n}(A)\in{}_{\pi_{0}(A)}\mathrm{Mod}(\mathrm{Ind(Ban}_{k}\mathrm{)}$. 
\end{cor}

\subsection{The Cotangent Complex}\label{sec:cotangent}

In this section we verify that the cotangent complex, and the canonical differential, of the algebras we have defined, give us what we would hope. Let
$$\mathrm{T}\in\{\mathbb{A}_{\mathrm{Ind(Ban}_{R}\mathrm{)}},\mathrm{FPS}_{R},\mathrm{Afnd}_{R}^{\dagger},\mathrm{EFC}_{R},\mathrm{Disc},\mathrm{Germ}_{R}\}$$
and if $R$ is non-Archimedean, let
$$\mathrm{T}\in\{\mathbb{A}_{\mathrm{Ind(Ban}^{nA}_{R}\mathrm{)}},\mathrm{Tate}_{R},\mathrm{FPS}^{nA}_{R},\mathrm{Afnd}_{R}^{nA,\dagger},\mathrm{EFC}^{nA}_{R},\mathrm{Disc}^{nA},\mathrm{Germ}^{nA}_{R}\}$$

Write $A_{(\lambda_{1},\ldots,\lambda_{n})}$ for the image of $\mathrm{T}(\lambda_{1},\ldots,\lambda_{n})$ in $\mathrm{Comm}(\mathrm{Ind(Ban}_{R}\mathrm{)})$ (or $\mathrm{Comm}(\mathrm{Ind(Ban}^{nA}_{R}\mathrm{)}$). Since $\mathrm{T}$ is of homotopy polynomial type, we have
$$\mathbb{L}_{A_{(\lambda_{1},\ldots,\lambda_{n})}}\cong A_{(\lambda_{1},\ldots,\lambda_{n})}^{\oplus n}$$
regarded as a module over itself. Consider the canonical derivation
$$d:A_{(\lambda)}\rightarrow A_{(\lambda_{1},\ldots,\lambda_{n})}^{\oplus n}$$
$A_{(\lambda_{1},\ldots,\lambda_{n})}^{\oplus n}$ is a subalgebra of the ring of formal power series $R[[x_{1},\ldots,x_{n}]]$. This algebra is equipped with a bounded derivation $\tilde{d}$ sending
$$\sum_{(i_{1},\ldots,i_{n})\in\mathbb{N}_{0}}a_{(i_{1},\ldots,i_{n})}x_{1}^{i_{1}}\ldots x_{n}^{i_{n}}$$
to 
$$(\sum_{i\in\mathbb{N}_{0}}i_{1}a_{(i_{1},\ldots,i_{n})}x^{i_{1}-1}x_{2}^{i_{2}}\ldots x_{n}^{i_{n}},\ldots,\sum_{i\in\mathbb{N}_{0}}i_{n}a_{(i_{1},\ldots,i_{n})}x^{i_{1}}x_{2}^{i_{2}}\ldots x_{n}^{i_{n}-1})$$
Moreover it restricts to a bounded derivation $\tilde{d}$ on $A_{(\lambda_{1},\ldots,\lambda_{n})}$. $\tilde{d}$ coincides with $d$. Indeed this is certainly true when we restrict to the subalgebra
$$R[x_{1},\ldots,x_{n}]\subset A_{(\lambda)}$$
Since $R[x_{1},\ldots,x_{n}]$ is bornologically dense in $A_{(\lambda)}$, it follows that $d=\tilde{d}$. 

Consider also the case $R=\mathbb{R}$ and the Lawvere theory $\mathrm{CartSp_{smooth}}$. Again the algebra $C^{\infty}(\mathbb{R}^{n})$ is equipped with a derivation $\tilde{d}:C^{\infty}(\mathbb{R})\rightarrow C^{\infty}(\mathbb{R})^{\oplus n}$,
$$f(x)\mapsto\Bigr(\frac{\partial f(x_{1},\ldots,x_{n}))}{\partial x_{1}},\ldots,\frac{\partial f(x_{1},\ldots,x_{n}))}{\partial x_{n}}\Bigr)$$
$$R[x_{1},\ldots,x_{n}]\subset C^{\infty}(\mathbb{R}^{n})$$
is again bornologically dense, so once more we must have $\tilde{d}=d$. 

Now let 
$$f:A_{(\lambda_{1},\ldots,\lambda_{m})}\rightarrow A_{(\gamma_{1},\ldots,\gamma_{n})}$$
 be a map of algebras. This corresponds to an $m$-tuple $(f_{1},\ldots,f_{m})$ of elements of $A_{(\gamma_{1},\ldots,\gamma_{n})}$. The map
 $$ A_{(\gamma_{1},\ldots,\gamma_{n})}^{\oplus m}\cong A_{(\gamma_{1},\ldots,\gamma_{n})}\hat{\otimes}_{A_{(\lambda_{1},\ldots,\lambda_{m})}}\mathbb{L}_{A_{(\lambda_{1},\ldots,\lambda_{m})}}\rightarrow\mathbb{L}_{A_{(\gamma_{1},\ldots,\gamma_{n})}}\cong A_{(\gamma_{1},\ldots,\gamma_{n})}^{\oplus n}$$
 is given by  the $n\times m$ matrix

\[ \left( \begin{array}{ccc}
\frac{\partial f_{1}}{\partial x_{1}} & \ldots & \frac{\partial f_{m}}{\partial x_{1}}  \\
\vdots & \ddots & \vdots \\
\frac{\partial f_{1}}{\partial x_{n}}  & \ldots &  \frac{\partial f_{m}}{\partial x_{n}}  \end{array} \right)\]

\section{Morphisms of Bornological Algebras}

Let us now discuss various types of maps introduced in Chapter 2 in the context of the various bornological algebras introduced above.

\subsection{Formally \'{E}tale and Smooth Morphisms}

Let us begin with formally \'{e}tale morphisms. Consider first the case $R=\mathbb{C}$ and $\mathrm{T}=\mathrm{Disc}_{R}$. In this case a map $f:A_{(\lambda_{1},\ldots,\lambda_{m})}\rightarrow A_{(\gamma_{1},\ldots,\gamma_{n})}$ between reduced globally finitely embeddable $\mathrm{Disc}_{\mathbb{C}}$ algebras is formally \'{e}tale precisely if it corresponds to a local biholomorphism. Any $\mathrm{Disc}_{\mathbb{C}}$-\'{e}tale map between any finitely $\mathrm{Disc}_{\mathbb{C}}$-presented algebras corresponds to a local biholomorphism If $R=\mathbb{R}$ and $\mathrm{T}=\mathrm{Cart_{smooth}}$ then any formally \'{e}tale map between algebras of functions on smooth manifolds corresponds to a local diffeomorphism. 

Let us establish some more general results.

\begin{lem}\label{lem:discreteTsmooth}
Let $k$ be a non-trivially valued Banach field. Let $\mathrm{T}$ be a $\Gamma$-sorted Lawvere theory of homotopy $\mathrm{Ind(Ban}_{k}\mathrm{)}$-polynomial type such
    that the bornology of each $\mathrm{F}(\mathrm{T}(\lambda_{1},\ldots,\lambda_{n}))$ has a countable basis
    for any $\underline{\lambda}$, $\mathrm{F}(\mathrm{Free_{T}}(\underline{\lambda}))$ is strongly (left) coherent.
    Let $A$ be a discrete finitely $\mathrm{T}$-presented algebra and let $A\rightarrow B\cong A\hat{\otimes}^{\mathbb{L}}\mathrm{T}(\lambda_{1},\ldots,\lambda_{c+k})\big\slash\big\slash (f_{1},\ldots,f_{c})$ with the sequence $f_{1},\ldots,f_{c}$ Koszul regular in the algebraic sense. Then $B$ is a discrete finitely $\mathrm{T}$-presented algebra.
\end{lem}

\begin{proof}
The Koszul complex is exact by the open mapping theorem/ Buchwalter's Theorem
\end{proof}

\begin{cor}
Let $k$ be a non-trivially valued Banach field. Let $\mathrm{T}$ be a $\Gamma$-sorted Lawvere theory of homotopy $\mathrm{Ind(Ban}_{k}\mathrm{)}$-polynomial type such
\begin{enumerate}
    \item
    that the bornology of each $\mathrm{F}(\mathrm{Free_{\mathrm{T}}}(\lambda_{1},\ldots,\lambda_{n}))$ has a countable basis
    \item 
    for any $\underline{\lambda}$, $\mathrm{F}(\mathrm{Free_{T}}(\underline{\lambda}))$ is strongly (left) coherent.
    \item 
    finitely presented modules over discrete finitely $\mathrm{T}$-presented algebras are transverse to standard $\mathrm{T}$-\'{e}tale morphisms.
\end{enumerate}
Then any standard $\mathrm{T}$-smooth morphism $f:A\rightarrow B$, where $A$ is a coherent $\mathrm{T}$-presented $k$-algebra is derived strong whenever the corresponding sequence is Koszul regular.
\end{cor}

\begin{proof}
By Lemma \ref{lem:discreteTsmooth} and Proposition \ref{prop:discretestrongsmooth} it suffices to prove that finitely presented $\pi_{0}(A)$-modules are transverse to standard $\mathrm{T}$-smooth maps. But using the factorisation
$$\pi_{0}(A)\rightarrow\pi_{0}(A)\hat{\otimes}\mathrm{T}(\lambda_{c+1},\ldots,\lambda_{c+k})\rightarrow\pi_{0}(A)\hat{\otimes}\mathrm{T}(\lambda_{1},\ldots,\lambda_{c+k})\big\slash (f_{1},\ldots,f_{c})$$
it clearly suffices to prove it for standard $\mathrm{T}$-\'{e}tale maps, and we have assumed this.
\end{proof}

In analytic geometry we can get infinite products of fields as \'{e}tale extensions.

\begin{rem}
Let $A$ be a finitely $\mathrm{Disc}_{\mathbb{C}}$-presented algebra which is reduced. If $A\rightarrow B$ is an $\mathrm{Disc}_{\mathbb{C}}$-\'{e}tale extension then $B$ is reduced. This is essentially a consequence of the open mapping theorem.
\end{rem}

\begin{lem}
A $\mathrm{Disc}_{\mathbb{C}}$-\'{e}tale extension of $\mathbb{C}$ is a countably infinite products of fields.
\end{lem}

\begin{proof}
Consider the space $A=\mathcal{O}(D^{n}_{k<\lambda})\big\slash(f_{1},\ldots,f_{n})$ where the determinant of the Jacobian of $(f_{1},\ldots,f_{n})$ is invertible in $A$. $A$ is the algebra of functions on the common vanishing locus of the $f_{i}$. Now by the condition on the Jacobian this must be a discrete, and hence either finite if $\lambda<\infty$, or countable otherwise. Thus $A$ is a product of the algebra of functions on each point, i.e. a product of copies of $\mathbb{C}$.

Let us show that we can get any countable product. Let $\{z_{k}\}_{k\in N}$ in $\mathbb{C}$. As a consequence of the Mittag-Leffler Theorem, there is an analytic function $f(x)\in\mathcal{O}(\mathbb{C})$ such that 
\begin{enumerate}
\item
$\{z_{k}\}_{k\in N}$ are precisely the zeroes of $f(x)$.
\item
each $z_{k}$ is a zero of $f$ of order $1$.
\end{enumerate}
Then $\mathcal{O}(\mathbb{C})\big\slash (f)$ is an \'{e}tale extension of $\mathbb{C}$. The map $\mathcal{O}(\mathbb{C})\big\slash (f)\rightarrow\prod_{k\in N}\mathbb{C}$ which sends $[g]$ to $g(z_{k})$ is an isomorphism.
\end{proof}

Often the diagonal of an \'{e}tale map will be a Zariski open immersion, as in the algebraic case.

\begin{lem}
Let $A$ be a complete bornological $R$-algebra and $I=(f_{1},\ldots,f_{n})$ a finitely generated admissible ideal such that
\begin{enumerate}
\item
$(f_{i}f_{j})_{1\le i,j\le n}$ is an admissible ideal.
\item
$I^{2}=I$.
\end{enumerate}
Then $A\rightarrow A\big\slash I$ is a Zariski open immersion.
\end{lem}

\begin{proof}
By the closedness assumption and Nakayama's lemma we have $I=(e)$ where $e\in A$. Then $A\big\slash I\cong A_{(1-e)}$ by Corollary \ref{cor:idempotentloc}.
\end{proof}

\subsection{Finite and Flat Morphisms}
Let
$$\mathrm{T}\in\{\mathbb{A}_{\mathrm{Ind(Ban}_{R}\mathrm{)}},\mathrm{FPS}_{R},\mathrm{Afnd}_{R},\mathrm{Afnd}_{R}^{\dagger},\mathrm{EFC}_{R},\mathrm{Disc}_{R},\mathrm{Germ}_{R}\}$$
and if $R$ is non-Archimedean, let
$$\mathrm{T}\in\{\mathbb{A}_{\mathrm{Ind(Ban}^{nA}_{R}\mathrm{)}},\mathrm{Tate}_{R},\mathrm{FPS}^{nA}_{R},\mathrm{Afnd}_{R}^{nA,\dagger},\mathrm{EFC}^{nA}_{R},\mathrm{Disc}^{nA}_{R},\mathrm{Germ}^{nA}_{R}\}$$
We also allow $\mathrm{T}=\mathrm{CartSp_{smooth}}$ when $R=\mathbb{R}$.

Let $S\in\mathrm{Comm(Ind(Ban_{R}))}$, (or $\mathrm{Comm(Ind(Ban^{nA}_{R}))}$ if $R$ is non-Archimedean, and let $A_{(\lambda)}$ be a generating algebra of one of these Lawvere theories. Consider a map of the form
$$S\rightarrow S\hat{\otimes}_{R}A_{(\lambda)}\big\slash (f)$$
where $f\in A_{(\lambda)}$. We want to determine when this map is finite, and in particular faithfully flat.

\begin{lem}
Let $R$ be a non-Archimedean Banach ring, and let $f=\sum_{i=0}^{n}a_{i}x^{i}\in R[x]$ be a monic polynomial with $|a_{i}|\le 1$ for all $i$. Then there is an isomorphism of Banach $R$-modules
$$R\Bigr<x\Bigr>\big\slash\overline{(f})\cong R^{\oplus n}$$
\end{lem}

\begin{proof}.
Let $\{e_{1},\ldots,e_{n}\}$ be the standard basis of $R^{\oplus n}$. The map $R^{\oplus n}\rightarrow R\Bigr<x\Bigr>$, $e_{i}\mapsto x^{i}$ is clearly bounded. Therefore the composite map 
$$R^{\oplus n}\rightarrow R\Bigr<x\Bigr>\big\slash(f)$$
is bounded. Consider the normed polyonomial subalgebra $R[x]\subset R\Bigr<x\Bigr>$ with norm being the restriction of the norm on $R\Bigr<x\Bigr>$. We therefore get an exact sequence. There is a map $R[x]\rightarrow R^{\oplus n}$ induced by sending $x^{i}$ to $e_{i}$ for $0\le i\le n-1$. The map on $x^{k}$ for $k\ge n$ is defined by writing it in the basis $x^{i}+(f(x))$. We get an exact sequence.
$$0\rightarrow (f)\rightarrow R[x]\rightarrow R^{\oplus n}\rightarrow0$$
and by completion an exact sequence
$$0\rightarrow \overline{(f)}\rightarrow R\Bigr<x\Bigr>\rightarrow R^{\oplus n}\rightarrow0$$

\end{proof}

Similarly one can prove the following. 

\begin{lem}
Let $R$ be an Archimedean Banach ring, and let $f=\sum_{i=0}^{n}a_{i}x^{i}\in R[x]$ be a monic polynomial with $\sum_{i=0}^{n}|a_{i}|\le 1$. Then there is an isomorphism of Banach $R$-modules
$$R\Bigr\{x\Bigr\}\big\slash\overline{(f)}\cong R^{\oplus n}$$
\end{lem}


\begin{cor}\label{cor:quotient free}
Let $R$ be a Banach ring and
let $f\in R[x]$. Then there is an isomorphism
$$\mathcal{O}(R)\big\slash\overline{(f)}\cong R^{\oplus n}$$
\end{cor}

\begin{proof}
We prove the Archimedean claim, the non-Archimedean claim being identical.
Write
$$\mathcal{O}(R)\cong\lim_{\lambda>0}R_{i}\{\frac{1}{\lambda}x\}$$
Once again we need to show that the induced map
$$R[x]\subset\mathcal{O}(R)\rightarrow R^{\oplus n}$$
is bounded, where the bornology on $R[x]$ is induced from $\mathcal{O}(R)$.  However for $\lambda$ sufficiently large we have a bounded map $R\{\frac{1}{\lambda}x\}\rightarrow R^{\oplus n}$. This suffices to prove the claim.
\end{proof}

\begin{cor}
Let $R$ be an arbitrary Banach ring and let $f\in R[x]$. Let $\mathrm{T}\in\{\mathrm{Afnd}_{R},\mathrm{EFC}_{R}\}$ (or $\mathrm{T}\in\{\mathrm{Tate}_{R},\mathrm{EFC}^{nA}_{R}\}$ if $R$ is non-Archimedean). Then there is a finite map of multiplicatively convex complete bornological rings
$$R\rightarrow\overline{R}$$
with $\overline{R}\in{}_{R}\mathrm{Alg_{T}}$
such that in $\overline{R}$, $f\cong\prod_{i=0}^{n}(x-\alpha_{i})$ for some $\alpha_{i}\in\overline{R}$. 
\end{cor}

\subsection{Localisations}

Let us now discuss $\mathrm{T}$-localisations for the specified Lawvere theories. As we will see pure Laurent localisations are sometimes flat (at least in the non-Archimedean case), but even in simple cases, general localisations are not flat.

\begin{lem}
Let $R$ be a non-Archimedean Banach ring, and let $g\in R$, and $\lambda\in\mathbb{R}_{>0}$. The map
$$R\rightarrow R\Bigr<\lambda x\Bigr>\big\slash(1-gx)$$
is transverse to filtered colimits of Banach $R$-modules.
\end{lem}

\begin{proof}
(c.f. \cite{bambozzi2020sheafyness} Proposition 4.23). Let $M$ be a Banach $R$-module. We claim that the multiplication by $(1-xg)$ map
$$M\Bigr<x\Bigr>\rightarrow M\Bigr<x\Bigr>$$
is an admissible monomorphism. First we prove that it is a monomorphism. Suppose $(1-xg)m(x)=0$. Solving the equations gives 
$$m_{0}=0,\;\; m_{i}=fm_{i-1},i\ge 1$$
which clearly implies $m=0$. Now write $m(x)=\sum_{i=0}^{\infty}m_{i}x^{i}$. Defining $m_{-1}=0$ for simplicity of notation, we have
\begin{align*}
\mathrm{sup}{||m_{i}||\lambda^{i}}&=\mathrm{sup}||m_{i}-gm_{i-1}+gm_{i-1}||\lambda^{i}\\
&\le\mathrm{sup}\{\mathrm{max}\{||m_{i}||\lambda^{i},||gm_{i-1}||\lambda^{i}\}\}\\
&\le\mathrm{max}\{\mathrm{sup}||m_{i}||\lambda^{i},\mathrm{sup}\lambda||gm_{i}||\lambda^{i}\}\\
&\le\mathrm{max}\{1,|g|\lambda\}\mathrm{sup}||m_{i}||\lambda^{i}
\end{align*}
Hence the map is open. It remains to show that the sub-module $(1-gx)M\Bigr<\lambda x\Bigr>$ is closed. Let $(1-xg)m_{i}(x)$ be convergent. By looking at coefficients of the $x_{i}$, it is clear that $m_{i}(x)$ converges to some $m(x)$, and therefore $(1-xg)m_{i}(x)$ converges to $(1-xg)m(x)$. 

In particular taking $M=R$ gives that  
$$R\Bigr<\lambda x\Bigr>\rightarrow R\Bigr<\lambda x\Bigr>$$
is a flat resolution of $R\Bigr<\lambda x\Bigr>\big\slash(1-gx)$ as an $R$-module. We have also shown that for $M$ a Banach $R$-module we have 
$$M\otimes_{R}^\mathbb{L}R\Bigr<\lambda x\Bigr>\big\slash(1-gx)\cong M\otimes_{R}R\Bigr<\lambda x\Bigr>\big\slash(1-gx)$$
For a general Ind-Banach $R$-module $M$, we may write $M=\mathrm{colim}_{\mathcal{I}}M_{i}$ where $M_{i}$ is a Banach $R$-module and the category $\mathcal{I}$ is filtered. Since filtered colimits are exact, for such $M$ we also have 
$$M\otimes_{R}^\mathbb{L}R\Bigr<\lambda x\Bigr>\big\slash(1-gx)\cong M\otimes_{R}R\Bigr<\lambda x\Bigr>\big\slash(1-gx)$$
Therefore $R\Bigr<\lambda x\Bigr>\big\slash(1-gx)$ is flat as an $R$-module. 
\end{proof}

\begin{lem}
Let $R$ be a Banach ring satisfying the open mapping property (resp. Buchwalter's property), and let $S$ be a complete bornological $R$-algebra. Let $g\in S$. The map
$$M\hat{\otimes}^{\mathbb{L}}_{S}S\hat{\otimes}\mathcal{O}(R)\big\slash(1-gx)\rightarrow M\hat{\otimes}_{S}S\hat{\otimes}\mathcal{O}(R)\big\slash(1-gx)$$
is a equivalence whenever $M$ is an $S$-module which is a filtered colimit of $S$-modules which are flat and Fr\'{e}chet as $R$-modules (or flat and with a bornology which has a countable basis as $R$-modules). The same result is true in the non-Archimedean setting.
\end{lem}

\begin{proof}
Since filtered colimits are exact, we may assume that $M$ is Fr\'{e}chet (or has a bornology with a countable basis). As in the previous result, the derived tensor product is computed by the two-term sequence
$$M\hat{\otimes}\mathcal{O}(R)\rightarrow M\hat{\otimes}\mathcal{O}(R)$$
Also as before (by comparing coefficients), algebraically there is an exact sequence
$$0\rightarrow M\hat{\otimes}\mathcal{O}(R)\rightarrow M\hat{\otimes}\mathcal{O}(R)\rightarrow  M\hat{\otimes}_{S}S\hat{\otimes}\mathcal{O}(R)\big\slash(1-gx)\rightarrow 0$$
Using the open mapping property (resp. Buchwalter's property), this is also admissibly exact.
\end{proof}



\begin{lem}\label{lem:transverse}
    Let $k$ be a non-trivially valued Banach field. Let $\mathrm{T}$ be one of $\{\mathrm{Disc}_{k},\mathrm{Afnd}^{\dagger}_{k}\}$ in the Archimedean case, or one of $\{\mathrm{Disc}^{nA}_{k},\mathrm{Afnd}^{\dagger,nA}_{k},\mathrm{Tate}_{k}^{nA}\}$ in the non-Archimedean case. Let $A$ be a discrete finitely $\mathrm{T}$-presented algebra and $M$ a finitely presented $A$-module. Let $A\rightarrow B$ be a $\mathrm{T}$-rational localisation, or a $\mathrm{T}$-open immersion in the disc cases. Then the map
    $$B\hat{\otimes}_{A}^{\mathbb{L}}M\rightarrow B\hat{\otimes}_{A}M$$
    is an equivalence.
\end{lem}

\begin{proof}
Since the algebras involved are coherent, there is a resolution of $M$ by degree-wise free $A$-algebras of finite rank. By Buchwalter's Theorem, it suffices to prove that 
$$B\otimes_{A}^{\mathbb{L}}M\rightarrow B\otimes_{A}M$$
is an algebraic equivalence. For Archimedean Steins this is \cite{MR1420618} Proposition 4.3.3 and Theorem 4.3.6. The non-Archimedean Stein case can be proven as in loc. cit. For affinoids this is \cite{BGR} Corollary 6. The dagger affinoid case is \cite{bambozzi2014generalization} Proposition 6.1.13.
\end{proof}

\subsubsection{Localisations are Not All Flat}

In \cite{ben2020fr} the first and third authors gave a simple, geometrically intutive, example of a localisation of Tate $\mathbb{Q}_{p}$-algebras which is not flat. Let us briefly recall the construction. Let $A=\mathbb{Q}_{p}<x>$. Consider the algebras 
$$A_{V}\defeq A<3y>\big\slash(y-x)$$
$$A_{W}\defeq A\bigr<\frac{z}{2}\Bigr>(xz-1)$$
Geometrically, these are the algebras of functions on $V=\{t\in\mathbb{Q}_{p}:|t|\le\frac{1}{3}\}$ and $W=\{t\in\mathbb{Q}_{p}:\frac{1}{2}\le |t|\le 1\}$ respectively. The map $A\rightarrow A_{W}$ is a closed embedding so we get an admissible exact sequence
$$0\rightarrow A\rightarrow A_{W}\rightarrow A_{W}\big\slash A\rightarrow 0$$
Now we have $A_{V}\hat{\otimes}_{A}A_{W}\cong0$ (geometrically this is because $V\cap W=\emptyset$). Thus tensoring the above sequence with $A_{V}$ gives the sequence
$$A_{V}\rightarrow 0\rightarrow 0$$
which is clearly not short exact.

Due to work of Frunz\v{a} \cite{frunzua1975taylor}, Putinar \cite{MR1420618}, and Finch \cite{finch1975single} we can also give an example of a natural localisation of Stein algebras over $\mathbb{C}$ which is not flat. Consider the incluson of open ball $B(0,k)\subset B(0,r)$ where $0<k<r\le\infty$. The inclusion map corresponds to the restriction map map of Stein algebras
$$\mathcal{O}(B(0,r))\rightarrow \mathcal{O}(B(0,k))$$
Recall the Koszul resolution
$$\mathcal{O}(B(0,r)\times B(0,r))\rightarrow\mathcal{O}(B(0,r)\times B(0,r))\rightarrow\mathcal{O}(B(0,r))$$
where the first map is multiplication by $(z-w)$, $z$ being the coordinate in the first variable, and $w$ the coordiante in the second. If $M$ is a $\mathcal{O}(B(0,r))$-module which is flat as a $\mathbb{C}$-module, then 
$$(z-T):\mathcal{O}(B(0,r))\hat{\otimes}M\rightarrow\mathcal{O}(B(0,r))\hat{\otimes}M$$
gives a flat resolution of $M$ as a $\mathcal{O}(B(0,r))$-module. Here $T$ is the operator $T:M\rightarrow M$ given by the action of $w$ on $M$. $\mathcal{O}(B(0,k))\hat{\otimes}^{\mathbb{L}}_{\mathcal{O}(B(0,r))}M$ is then computed by the complex
$$(z-T):\mathcal{O}(B(0,k))\hat{\otimes}M\rightarrow\mathcal{O}(B(0,k))\hat{\otimes}M$$
Just as a chain complex, we have $H^{1}(\mathcal{O}(B(0,k))\hat{\otimes}^{\mathbb{L}}_{\mathcal{O}(B(0,r))}M)\cong 0$ if and only if the map
$$\mathcal{O}(B(0,k))\hat{\otimes}M\rightarrow\mathcal{O}(B(0,k))\hat{\otimes}M$$
is injective with closed image. This is the same as saying that $T$ has \textit{Bishop's property } $(\beta)$ or \textit{the single valued extension property}. The translation of this property into the homological language of the Koszul complex is due to work of \cite{frunzua1975taylor} and \cite{MR1420618}. Finch \cite{finch1975single} showed that if an operator $T$ is surjective but not injective then it does not satisfy the single-valued extension property.

Let $E$ be a flat Banach space and $T:E\rightarrow E$ a linear operator such that the spectrum $\sigma(T)$ is contained in $B(0,r)$. Let $\Gamma$ be the boundary of an open disc in $B(0,r)$ containing $\sigma(T)$. $E$ becomes a $\mathcal{O}(B(0,r))$-module by defining
$$f\mapsto f(T)\defeq\int_{\Gamma}f(z)(z-T)^{-1}dz$$
Suppose $T:E\rightarrow E$ is such an operator, and moreover that it is onto but not injective. For example, one can take $E=\bigoplus^{\le1}_{n\in\mathbb{N}_{0}}\mathbb{C}_{\delta}$ with $\delta>\frac{1}{r}$, and $T:E\rightarrow E$ to be the map $(r_{1},r_{2},r_{3},r_{4},\ldots)\mapsto(r_{1}+r_{2},r_{3},r_{4},r_{5},\ldots)$. 

Let $x_{0}\in\mathrm{Ker}(T)$ be such that $||x_{0}||=1$. By the open mapping theorem we may inductively define $x_{n}$ such that there exists a $\epsilon>0$ with $Tx_{n}=x_{n-1}$ and $||x_{n}||\le \epsilon^{n}$. Then 
$$f(z)\defeq\sum x_{n}z^{n}$$
is analytic on the ball $B(0,\frac{1}{\epsilon})$.  Then we have
$$(z-T)f=\lim x_{n}z^{n+1}=0$$
In particular $(z-T)$ is not injective, and $\mathcal{O}(B(0,r))\rightarrow\mathcal{O}(B(0,\frac{1}{\epsilon}))$ is not flat.

In our example one can take $\epsilon=\delta$ for any $\delta>\frac{1}{r}$, i.e. for any $\mathcal{O}(B(0,r))\rightarrow \mathcal{O}(B(0,k))$ is not flat for any $k<r$. 

\subsubsection{$C^{\infty}$-Localisations Are Flat}

In the smooth case we find that in fact many maps are flat, essentially because of the existence of partitions of unity.

\begin{prop}[\cite{ogneva} Theorem 2]\label{thm:localisationflat}
Let $M$ be a finite-dimensional smooth manifold, and $U\subset M$ an open subset which is contained in a coordinate chart. Then $C^{\infty}(U)$ is a flat bornological $C^{\infty}(M)$-module.
\end{prop}

\begin{proof}
We prove in exactly the same way as \cite{ogneva} Theorem 2, that $C^{\infty}(U)$ is a retract of the complete boronological $C^{\infty}(M)$-module, $C^{\infty}(M\times U)\cong C^{\infty}(M)\wotimes C^{\infty}(U)$. As a nuclear (Fr\'{e}chet) space, $C^{\infty}(U)$ is flat as a complete bornological $\mathbb{R}$-vector space. Thus being a free $C^{\infty}(M)$-module on a flat space, $C^{\infty}(M\times U)$ is a flat module. 

Consider the canonical projection map
$$\pi_{U}:C^{\infty}(M\times U)\rightarrow C^{\infty}(U)$$
Ogneva shows that this map has a left inverse, and we shall explain why it is bounded. Let $U\subset W$ with $(W,\omega)$ being a chart. Consider 
$$\psi(x)\defeq\mathrm{min}(1,\mathrm{dist}(x,\partial\omega(U)))$$
defined on $\omega(U)$. Let $\phi(x)$ be any smooth function on $\omega(U)$ with $0<\phi(X)<\psi(x)$, and define the function $\theta$ on $\omega(U)\times\omega(U)$ by
\[
\theta(x,y)\defeq
\begin{cases}
\mathrm{exp}\frac{|x-y|^{2}}{|x-y|^{2}-\phi(y)^{2}},&\;\;\mathrm{if }|x-y|\le\phi(y)\\
0&\;\;\mathrm{if }|x-y|>\phi(y)
\end{cases}
\]
For $(s,t)\in M\times U$ define
\[
F(s,t)\defeq
\begin{cases}
\theta(\omega(s),\omega(t))&\;\;s\in U\\
0 & s\notin U
\end{cases}
\]
$F(s,t)$ is a well-defined element of $C^{\infty}(M\times U)$, and $F(s,s)=1$ for $s\in U$. Define $\rho:C^{\infty}(U)\rightarrow C^{\infty}(M\times U)$ by
\[
(\rho f)(s,t)\defeq
\begin{cases}
f(s)F(s,t)&\;\;s\in U\\
0 & s\notin U
\end{cases}
\]
It is easy to see that this map is continuous, and therefore bornologically bounded. 
\end{proof}

\subsection{Derived Algebras}

Let us now consider derived algebras over the various analytic Lawvere theories.

\begin{prop}\label{prop:NoetherianStrong}
    Let $k$ be a non-trivially valued Banach field.
\begin{enumerate}
    \item 
    For $k$ Archimedean $\mathrm{Afnd}_{k}^{\dagger}$ is strongly  $\mathrm{CBorn}_{k}$-Noetherian.
    \item 
    For $k$ non-Archimedean both $\mathrm{Afnd}_{k}^{nA,\dagger}$ and $\mathrm{Tate}_{k}$ are strongly  $\mathrm{CBorn}^{nA}_{k}$-Noetherian
\end{enumerate}
\end{prop}

\begin{proof}
    \begin{enumerate}
        \item 
      The dagger affinoid case is \cite{bambozzi2014generalization} Theorem 3.1.11.
        \item
         The dagger affinoid case is again \cite{bambozzi2014generalization} Theorem 3.1.11. The affinoid case is \cite{BGR} Chapter 6 Proposition 3.
    \end{enumerate}
\end{proof}

  

\begin{lem}\label{loc:cohTan}
 Let $k$ be a non-trivially valued non-Archimedean Banach field, and let $A\rightarrow B$ be a $\mathrm{T}$-rational localisation of $\mathrm{T}$-coherent algebras where $\mathrm{T}\in\,\mathrm{Afnd}^{\dagger}_{k}\}$ in the Archimedean case, or $\mathrm{T}\in\{\mathrm{Disc}^{nA}_{k},\mathrm{Afnd}^{nA,\dagger}_{k},\mathrm{Tate}_{k}\}$ in the non-Archimedean case. Then $\pi_{0}(A)\rightarrow\pi_{0}(B)$ is a (derived) $\mathrm{T}$-rational localisation and $A\rightarrow B$ is derived strong. Conversely if $A\rightarrow B$ is derived strong and $\pi_{0}(A)\rightarrow\pi_{0}(B)$ is a $\mathrm{T}$-rational localisation, then $A\rightarrow B$ is a $\mathrm{T}$-rational localisation.
\end{lem}

\begin{proof}
  Suppose that $A\rightarrow B$ is a $\mathrm{T}$-rational localisation of $\mathrm{T}$-coherent algebras. By Lemma \ref{lem:discreteTsmooth} and Lemma \ref{lem:transverse} it is derived strong and $\pi_{0}(A)\rightarrow\pi_{0}(B)$ is a $\mathrm{T}$-rational localisation. 

The converse is immediate from Lemma \ref{lem:derratrel}
  \end{proof}

\subsubsection{Derived Dagger Stein Algebras}

Let $k$ be a non-trivially valued Banach field.

\begin{defn}[\cite{bambozzi2015stein}]\label{defn:daggerstein}
    A \textit{dagger Stein algebra over }$k$ is a complete bornological algebra $A$ which is isomorphic to an inverse limit of $k$-dagger affinoids
    $$A\cong\lim_{\mathbb{N}^{op}} A_{i}$$
    where each $A\rightarrow A_{i+1}$ is a Weierstrass localisation and $\mathcal{M}(A_{i})$ is contained in the interior of $\mathcal{M}(A_{i+1})$. The full subcategory of $\mathrm{Comm}(\mathrm{Ind(Ban}_{k}\mathrm{)})$ consisting of dagger Stein algebras is denoted $\mathrm{StnAlg}^{\dagger}$. 
\end{defn}

Here $\mathcal{M}(A)$ is the Berkovich spectrum of $A$, defined in \cite{bambozzi2015stein} Definition 4.1, and recalled in the present work in Chapter 9. By Lemma 4.15 in \cite{bambozzi2015stein}, and the fact that dagger affinoids are strongly Noetherian, dagger Stein algebras are bornological Fr\'{e}chet-Stein algebras. In fact by \cite{bambozzi2015stein} Theorem 4.25, the assignment $A\mapsto \mathcal{M}(A)$ gives an antiequivalence between $\mathrm{StnAlg}^{\dagger}$ and the category $\mathrm{Stn}^{\dagger}$ of dagger Stein spaces defined in \cite{bambozzi2015stein} Definition 4.13.

\begin{defn}
    An object $A$ of $\mathbf{DAlg}^{cn}(\mathrm{Ind(Ban}_{k}\mathrm{)})$ is said to be a \textit{derived dagger Stein algebra} if 
    \begin{enumerate}
        \item 
        $\pi_{0}(A)$ is a dagger Stein algebra.
        \item 
        each $\pi_{n}(A)$ is a coadmissible $\pi_{0}(A)$-module.
    \end{enumerate}
    The full subcategory of $\mathbf{DAlg}^{cn}(\mathrm{Ind(Ban}_{k}\mathrm{)})$ consisting of derived dagger Stein algebras is denoted $\mathbf{dStnAlg}^{\dagger}$, 
\end{defn}

Let $A\in \mathbf{dStnAlg}^{\dagger}$. We can construct a presentation $(A,A_{n},\overline{A}_{n})$ such that 
$$\pi_{0}(A)\cong\lim_{\mathbb{N}^{op}}\pi_{0}(\overline{A}_{n})$$
is a presentation as in Definition \ref{defn:daggerstein}. We will call such a presentation a \textit{geometric presentation}. In fact the proof of \cite{bambozzi2015stein} Lemma 4.19 implies that we can construct a diagram
$$A\rightarrow\ldots\rightarrow A_{n+1}\rightarrow\overline{A}_{n+1}\rightarrow A_{n}\rightarrow\ldots$$
where $(A,A_{n}\overline{A}_{n})$ is a geometric presentation, each $A_{n}$ is itself a dagger Stein algebra, such that the maps $A\rightarrow A_{n}$ are homotopy monomorphisms. Equivalently by \cite{bambozzi2015stein} Section 5 they correspond to open immersions $\mathcal{M}(A_{n})\rightarrow\mathcal{M}(A)$.

\begin{prop}
    Let $(A,A_{n},\overline{A}_{n})$ and $(B,B_{n},\overline{B}_{n})$ be geometric presentations of dagger Stein algebras. Let $F:A\rightarrow B$ be a map of algebras. Possibly after re-indexing there are inverse systems of maps 
    $$f_{n}:\overline{A}_{n}\rightarrow\overline{B}_{n}$$
    such that 
    $$f=\mathbf{lim}_{\mathbb{N}^{op}}f_{n}$$
\end{prop}

\begin{proof}
We can write
$$\pi_{0}(f)\cong\lim_{\mathbb{N}^{op}}\tilde{f}_{i}$$
using \cite{bambozzi2015stein} Lemma 4.24. Since each $A\rightarrow\overline{A}_{n+1}$ is a homotopy epimorphism and hence formally \'{e}tale, we can find an essentially unique lift of the map $\pi_{0}(A_{0})\rightarrow \pi_{0}(B_{0})$ to $A_{0}\rightarrow B_{0}$ as a map of $A$-algebras. We can then sequentially find lifts in the diagrams
\begin{displaymath}
    \xymatrix{
A\ar[d]\ar[rr] & & B_{n+1}\ar[d]\\
A_{n+1}\ar[r] & A_{n}\ar[r] & B_{n}
    }
\end{displaymath}
\end{proof}

\begin{defn}
\begin{enumerate}
    \item 
     Let $A$ be a dagger Stein algebra. We define $\mathrm{Coad}(A)$ to be the full subcategory of ${}_{A}\mathrm{Mod}(\mathrm{Ind(Ban}_{k}\mathrm{)})^{\heart}$ consisting of those modules $M$ such that there exists a geometric presentation $(A,A_{n},\overline{A}_{n})$ with $M\in\mathrm{Coad}(A,A_{n},\overline{A}_{n})$ for each $m$.
     \item 
          Let $A$ be a derived dagger Stein algebra. We define $\mathbf{Coad}(A)$ to be the full subcategory of ${}_{A}\mathbf{Mod}$ consisting of those modules $M$ such that each $\pi_{m}(M)$ is in $\mathrm{Coad}(\pi_{0}(A))$.
\end{enumerate}
\end{defn}

It follows from descent for coherent sheaves, and the bornological versions of Cartan's/ Kiehl's theorems \cite{MR3975530}, that the category $\mathbf{Coad}(A)$ is independent of the choice of geometric presentation. Moreover any two geometric presentations have a common refinement. In fact Theoerems A and B imply that $\mathrm{Coad}(A)$ is equivalent to the category of ($\mathrm{CBorn}_{\mathbb{C}}$-valued) sheaves on $\mathcal{M}(A)$ which are coherent. 


\begin{rem}
We may always pick a geometric presentation $(A,A_{n},\overline{A}_{n})$ of a dagger Stein algebra such that each $A\rightarrow A_{n}$ is transverse to coadmissible $A$-modules. Indeed topologically we may write $\mathcal{M}(\overline{A}_{n})$ as an intersection of dagger Stein neighbourhoods $\mathcal{M}(A_{r})$ of decreasing radius in $\mathcal{M}(A)$, and then we have that $\overline{A}_{n}\cong\colim_{r} A_{r}$. This can easily be made formal using the proof of Lemma 4.19 in \cite{bambozzi2015stein}, for example. By \cite{MR1420618} Corollary 4.2.5 -proved in the ARchimedean case there, but the proof also works in the non-Archimedean case, the map $A\rightarrow A_{r}$ is transverse to objects in $\mathrm{Coad}(A)$, so is the map $A\rightarrow\colim_{r}A_{r}$.

 In particular if $A$ is ae derived dagger Stein algebra, we can and will assume our presentation $(A,A_{n},\overline{A}_{n})$ is such that the map $A\rightarrow\overline{A}_{n}$ is derived strong.
\end{rem}

\begin{prop}\label{prop:coad}
Let $A$ be a dagger Stein algebra. $\mathrm{Coad}(A)$  is thik, closed under finite limits and finite colimits.
\end{prop}

\begin{proof}
    Let 
    $$M\rightarrow N$$
    be a map in $\mathrm{Coad}(A)$. By picking common refinements, we may assume there is a single geometric presentation $(A_{n},\overline{A}_{n})$ such that $M$ and $N$ are both in $\mathrm{Coad}(A,A_{n},\overline{A}_{n})$. Then closure under extensions, finite limits, and finite colimits follows from Proposition \ref{prop:coadexact}.

\end{proof}

\begin{cor}
    Let $A$ be a derived dagger Stein algebra. Then $\mathbf{Coad}_{+}(A)$ is thick and closed under finite limits and finite colimits.
\end{cor}



\begin{thm}\label{thm:coadpresheaf}
    Let $f:A\rightarrow B$ be a map of derived dagger Stein algebras, and let $M\in\mathbf{Coad}_{+}(A)$. Then $B\hat{\otimes}_{A}^{\mathbb{L}}M$ is in $\mathbf{Coad}_{+}(B)$.
\end{thm}

\begin{proof}
The proof can be seen as a derived bornological version of \cite{fan2017d} Lemma 2.2.10. Suppose we write $f$ as $\mathbf{lim}_{\mathbb{N}^{op}}f_{i}$ where $f_{i}:\overline{A}_{i}\rightarrow \overline{B}_{i}$ is a map and $(A,A_{n},\overline{A}_{n})$, $(B,B_{n},\overline{B}_{n})$ are presentations. Now 
$$\overline{B}_{i}\hat{\otimes}_{B}^{\mathbb{L}}(B\hat{\otimes}_{A}^{\mathbb{L}}M)\cong \overline{B}_{i}\hat{\otimes}_{A}^{\mathbb{L}}M\cong \overline{B}_{i}\hat{\otimes}_{\overline{A}_{i}}^{\mathbb{L}}\overline{A}_{i}\hat{\otimes}_{A}^{\mathbb{L}}M$$
Since $A\rightarrow \overline{A}_{i}$ is derived strong we have $\pi_{n}(\overline{A}_{i}\otimes_{A}^{\mathbb{L}}M)\cong\pi_{0}(\overline{A}_{i})\otimes_{\pi_{0}(A)}^{\mathbb{L}}\pi_{n}(M)$. In particular each $\pi_{n}(\overline{A}_{i}\otimes_{A}^{\mathbb{L}}M)$ is finitely presented as a $\pi_{0}(A_{i})$-module. By Lemma \ref{lem:cohfilt} each $\pi_{n}(\overline{B}_{i}\hat{\otimes}^{\mathbb{L}}_{\overline{A}_{i}}\overline{A}_{i}\hat{\otimes}^{\mathbb{L}}_{A}M)$ is finitely presented as a $\pi_{0}(\overline{B_{i}})$-module. Moreover each map $\overline{B}_{i}\rightarrow B_{i}$ is also derived strong, so we have that $\pi_{n}(B_{i}\hat{\otimes}_{\overline{B}_{i}}^{\mathbb{L}}\overline{B}_{i}\hat{\otimes}^{\mathbb{L}}_{\overline{A}_{i}}\overline{A}_{i}\hat{\otimes}^{\mathbb{L}}_{A}M)\cong\pi_{n}(B_{i}\hat{\otimes}_{A}^{\mathbb{L}}M)$ is also finitely presented as a $\pi_{0}(B_{i})$-module. Thus it is a nuclear Fr\'{e}chet space. Next we claim that the maps

$$\pi_{n}(B_{i+1}\hat{\otimes}_{A}^{\mathbb{L}}M)\rightarrow \pi_{n}(B_{i}\hat{\otimes}_{A}^{\mathbb{L}}M)$$
have dense image. Now we have that $B_{i+1}\rightarrow B_{i}$ is derived strong. So by Proposition \ref{prop:derstronghtpy} we have that 
$$\pi_{*}(B_{i}\otimes_{A}^{\mathbb{L}}M)\cong\pi_{0}(B_{i})\hat{\otimes}^{\mathbb{L}}_{\pi_{0}(B_{i+1})}\pi_{*}(B_{i+1}\otimes_{A}^{\mathbb{L}}M)$$
Since $\pi_{0}(B_{i+1})\rightarrow\pi_{0}(B_{i})$ has dense image, the tensor product also has dense image. By the bornological Mittag-Leffler theorem, \cite{bambozzi2015stein} Corollary 3.80, we have that the map
$$\mathbf{lim}\;\pi_{n}(B_{i}\hat{\otimes}_{A}^{\mathbb{L}}M)\cong\lim\pi_{n}(B_{i}\hat{\otimes}_{A}^{\mathbb{L}}M)$$
is an equivalence. Write $M_{i}\defeq A_{i}\otimes_{A}^{\mathbb{L}}B$, which is a nuclear Fr\'{e}chet space, so that the above computation in fact gives 
 $$\mathbf{lim}\;\pi_{n}(B_{i}\hat{\otimes}_{A_{i}}^{\mathbb{L}}M_{i})\cong\lim\pi_{n}(B_{i}\hat{\otimes}_{A_{i}}^{\mathbb{L}}M_{i})$$
 Now consider the bar complex
$$|B_{i}\hat{\otimes}^{\mathbb{L}}A_{i}^{\hat{\otimes}^{\mathbb{L}}n}\hat{\otimes}^{\mathbb{L}}M_{i}|$$
which computes $B_{i}\otimes_{A_{i}}^{\mathbb{L}}M_{i}$. The homology of each term of this complex is lim-acyclic. Indeed since everything is homotopy flat we have 
$$\pi_{*}(B_{i}\hat{\otimes}^{\mathbb{L}}A_{i}^{\hat{\otimes}^{\mathbb{L}}n}\hat{\otimes}^{\mathbb{L}}M_{i})\cong\pi_{*}(B_{i})\hat{\otimes}\pi_{*}(A_{i})^{\hat{\otimes}n}\hat{\otimes}\pi_{*}(M_{i})$$
This is a sequence of nuclear Fr\'{e}chet spaces whose transition maps are dense. Thus by the bornological Mittag-Leffler theorem we again that it is lim-acyclic. We have shown that for fixed $m$ the sequence $\pi_{m}(|B_{i}\hat{\otimes}^{\mathbb{L}}A_{i}^{\hat{\otimes}^{\mathbb{L}}n}\hat{\otimes}^{\mathbb{L}}M_{i}|)$ is lim-acyclic. It follows that we have 
$$\mathbf{lim}_{\mathbb{N}^{op}} B_{i}\hat{\otimes}_{A_{i}}^{\mathbb{L}}M_{i}\cong\mathbf{lim}_{\mathbb{N}^{op}}|B_{i}\hat{\otimes}^{\mathbb{L}}A_{i}^{\hat{\otimes}^{\mathbb{L}}n}\hat{\otimes}^{\mathbb{L}}M_{i}|\cong |B\hat{\otimes}^{\mathbb{L}}A^{\hat{\otimes}^{\mathbb{L}}n}\hat{\otimes}^{\mathbb{L}}M|\cong B\otimes_{A}^{\mathbb{L}}M$$
as rquired.



\end{proof}

%
%

A similar proof can be used to show the following.

\begin{prop}
    Let $A$ be a derived dagger Stein and $M,N\in\mathbf{Coad}_{+}(A)$. Then $M\otimes_{A}^{\mathbb{L}}N\in\mathbf{Coad}_{+}(A)$.
\end{prop}

\begin{defn}
    A derived dagger Stein algebra $A$ is said to be \textit{finitely embeddable} if $\pi_{0}(A)\cong\mathcal{O}(D^{n}_{<\rho,k})\big\slash I$ for some $0<\rho\le\infty$ and some $n$
\end{defn}

\begin{lem}
    Let $k$ be a non-trivially valued Banach ring and let $A$ be a finitely embeddable derived dagger Stein algebra $A$ is such that $\pi_{0}(A)\cong\mathcal{O}(D^{m}_{k,<\infty})\big\slash J$ for some closed ideal $J$.
\end{lem}

\begin{proof}
For both claims it suffices to observe that by \cite{MR4036665} Proposition 1.13 and Proposition 1.24 any quotient of $\mathcal{O}(D^{n}_{k,<r})$ by a finitely generated ideal is equivalent to a quotient of $\mathcal{O}(D^{m}_{k,<\infty})$ by a closed ideal, for some $n$.
\end{proof}

\begin{lem}\label{lem:cotangentfembed}
    Let $k$ be a non-trivially valued Banach field of characteristic $0$, and $A\rightarrow B$ be a map of finitely embeddable derived dagger Stein algebras. Then $\mathbb{L}_{B\big\slash A}\in\mathbf{Coad}_{\ge0}(B)$.
\end{lem}

\begin{proof}
    It suffices to prove this in the absolute case. Thus we consider $\mathbb{I}\rightarrow A$ with $A$ derived dagger Stein. We have a surjection $\mathcal{O}(D^{n}_{<\rho,k})\rightarrow\pi_{0}(A)$ which lifts to a map $\mathcal{O}(D^{n}_{<\rho,k})\rightarrow A$. Since $\mathbb{L}_{\mathcal{O}(D^{n}_{<\rho,k})}\cong \mathcal{O}(D^{n}_{<\rho,k})^{\oplus n}$ is in $\mathrm{Coad}(\mathcal{O}(D^{n}_{<\rho,k}))$, it suffices to prove that $\mathbb{L}_{A\big\slash \mathcal{O}(D^{n}_{<\rho,k})}$ is coadmissible. This follows immediately from Lemma \ref{lem:cotangentsubcat} 

    \end{proof}

\subsection{Change of Base}\label{subsec:changeofbase}

In this subsection we discuss change of base of Banach rings. Let $f:R\rightarrow S$ be a map of Banach rings. The tensor-hom adjunction gives an adjunction
$$\adj{S\hat{\otimes}^{\mathbb{L}}_{R}(-)}{\mathbf{Ch}(\mathrm{Ind(Ban}_{R}\mathrm{)})}{\mathbf{Ch}(\mathrm{Ind(Ban}_{S}\mathrm{)})}{|-|}$$
where $S$ $|-|$ denotes the forgetful functor. If $R$ and $S$ are both non-Archimedean, we similarly get an adjunction
$$\adj{S\hat{\otimes}^{\mathbb{L}}_{R}(-)}{\mathbf{Ch}(\mathrm{Ind(Ban^{nA}}_{R}\mathrm{)})}{\mathbf{Ch}(\mathrm{Ind(Ban^{nA}}_{S}\mathrm{)})}{|-|}$$
In both cases left adjoint is strong monoidal and left $t$-exact, so they determines a transformation of derived algebraic contexts. However if perhaps only $S$ is non-Archimedean, there is still a transformation of derived algebraic contexts from
$\underline{\mathbf{Ch}}(\mathrm{Ind(Ban}_{R}\mathrm{)})$ to $\underline{\mathbf{Ch}}(\mathrm{Ind(Ban^{nA}}_{S}\mathrm{)})$. First consider the underived categories $\mathrm{Ch}(\mathrm{Ind(Ban}_{S}\mathrm{)})$ and $\mathrm{Ch}(\mathrm{Ind(Ban^{nA}}_{S}\mathrm{)})$. Now the inclusion
$$i:\mathrm{Ch}(\mathrm{Ind(Ban^{nA}}_{S}\mathrm{)})\subseteq \mathrm{Ch}(\mathrm{Ind(Ban}_{S}\mathrm{)})$$
is fully faithful, and by \cite{MR3448274} Theorem 7.4, this inclusion has a strongly monoidal left adjoint $\pi$. 

\begin{lem}
The adjunction
$$\adj{\pi_{S}}{\mathrm{Ch}(\mathrm{Ind(Ban}_{S}\mathrm{)})}{\mathrm{Ch}(\mathrm{Ind(Ban^{nA}}_{S}\mathrm{)})}{i}$$
is a strongly monoidal Quillen adjunction which is left $t$-exact.
\end{lem}

\begin{proof}
Left $t$-exactness is clear by construction. The rest follows immediately from the fact that by the proof of \cite{MR3448274} Theorem 7.4, for $V$ a normed set $\pi(\coprod^{\le1}_{v\in V}S_{\rho(v)})\cong \coprod^{\le1,nA}_{v\in V}S_{\rho(v)}$.
\end{proof}

In particular if $R\rightarrow S$ is a map of rings with $S$ non-Archimedean, we get an adjunction

$$\adj{\pi_{S}\circ S\hat{\otimes}^{\mathbb{L}}_{R}(-)}{\mathbf{Ch}(\mathrm{Ind(Ban_{R}}\mathrm{)})}{\mathbf{Ch}(\mathrm{Ind(Ban^{nA}}_{S}\mathrm{)})}{|-|\circ i}$$
which is a transformation of derived algebraic contexts.

Let $\mathrm{T}_{R}$ be one of the Lawvere theories \[\{\mathrm{Disc}_{R},\mathrm{EFC}_{R},\mathrm{Afnd}_{R}^{\dagger},\mathrm{Tate}_{R},\mathrm{Germ}_{R},\mathbb{A}_{\mathrm{CBorn}_{R}},\mathrm{FPS}_{R}\}\] with either $R$ and $S$ both Archimedean, or $R$ and $S$ both Archimedean and let $R\rightarrow S$ be a map of rings. Then we have
$$S\hat{\otimes}^{\mathbb{L}}_{R}\mathrm{Free_{\mathrm{T}_{R}}}\cong\mathrm{Free_{\mathrm{T}_{S}}}$$
If $R$ is Archimedean and $S$ is non-Archimedean, then for $T$ one of $\{\mathrm{Afnd}_{R}^{\dagger},\mathrm{Tate}_{R},\mathrm{Germ}_{R},\mathbb{A}_{\mathrm{CBorn}_{R}}\}$ we have

$$\pi_{S}(S\hat{\otimes}^{\mathbb{L}}_{R}\mathrm{Free_{\mathrm{T}_{R}}})\cong\mathrm{Free_{\mathrm{T}_{S}}}$$

\chapter{Stacks and Sheaves}\label{StSh}

Now that we've developed the necessary commutative algebra (i.e. affine geometry), in this section we will introduce our conventions for more general geometric stacks. We will also establish some very general descent results for sheaves on stacks.

\section{Topologies}

We begin by discussing topologies. Recall that a Grothendieck pre-topology on an $(\infty,1)$-category is essentially equivalent to the data of a pre-topology on its homotopy category.

\begin{defn}[\cite{toen2008homotopical}*{Proposition 4.3.5}]
Let $\mathbf{M}$ be a finitely complete $(\infty,1)$-category. A \textit{Grothendieck pre-topology} on $\mathbf{M}$, denoted $\tau$, is a Grothendieck pre-topology on the homotopy category $\Ho(\mathbf{M})$.
\end{defn}

\subsection{Local Maps}\label{subsubsec:local}

Particularly when we discuss \'{e}tale and smooth maps of later, we will need techniques to make classes of maps `local'.

\begin{defn}
Let $\tau$ be a pre-topology on and $(\infty,1)$-category $\Ho(\mathbf{M})$. A class of maps $\mathbf{P}$ in $\Ho(\mathbf{M})$ is said to be $\tau$-\textit{local} if whenever $f:X\rightarrow Y$ is a map such that there is a cover $g:C\rightarrow X$ in $\tau$ with $f\circ g\in\mathbf{P}$ then $f\in\mathbf{P}$.
\end{defn}

For example, in classical (derived) algebraic geometry over a unital commutative ring $R$, smooth maps are local for the \'{e}tale pre-topology. In our contexts it is not clear that this will always be the case. However there is a way to turn a class of maps into a $\tau$-local class.

\begin{defn}
Let $\mathbf{P}$ be a class of maps in $\Ho(\mathbf{M})$ and $\tau$ a pre-topology on $\Ho(\mathbf{M})$. Denote by $\mathbf{P}^{\tau}$ the class of maps $f:X\rightarrow Y$, such that there is a cover $g:C\rightarrow X$ in $\tau$ with $f\circ g\in\mathbf{P}$.
\end{defn}

Clearly $\mathbf{P}\subset\mathbf{P}^{\tau}$ with equality if and only if $\mathbf{P}$ is $\tau$-local.

\begin{prop}
Let $\mathbf{P}$ be a class of maps in $\Ho(\mathbf{M})$ and $\tau$ a  pre-topology on $\mathbf{M}$. Then
\begin{enumerate}
\item
$\mathbf{P}^{\tau}$ is $\tau$-local.
\item
If $\mathbf{P}$ is stable under composition, pullback and equivalences then so is $\mathbf{P}^{\tau}$
\end{enumerate}
\end{prop}

\begin{proof}
\begin{enumerate}
\item
Let $f:X\rightarrow Y$ be a map, and let $g:C\rightarrow X$ be a cover in $\tau$ such that $f\circ g$ is in $\mathbf{P}^{\tau}$. Then there is a cover $q:D\rightarrow C$ in $\tau$ such that $f\circ g\circ q$ is in $\mathbf{P}$. Then $g\circ q:D\rightarrow X$ is a cover in $\tau$ such that $f\circ(g\circ q)\in\mathbf{P}$.
\item
First suppose that $\mathbf{P}$ is stable under homotopy pullbacks. Let $f:X\rightarrow Y$ be a map in $\mathbf{P}^{\tau}$, and let $g:C\rightarrow X$ be a cover in $\tau$ such that $f\circ g$ is in $\mathbf{P}$. Let $h:Z\rightarrow Y$ be a map. Consider the pullback $f\times_{Y}Z:X\times_{Y}Z\rightarrow Z$, and the cover $g\times_{Y}Z:C\times_{Y}Z\rightarrow X\times_{Y}Z$ in $\tau$. Now $f\times_{Y}Z\circ g_{Y}Z=(f\circ g)\times_{Y}Z$. Each of these maps is in $\mathbf{P}$ since $\mathbf{P}$ is stable under homotopy pullback. 

Now suppose that $\mathbf{P}$ is also stable by equivalences. Let
\begin{displaymath}
\xymatrix{
X\ar[r]^{f}\ar[d]^{u} & Y\ar[d]^{w}\\
\tilde{X}\ar[r]^{\tilde{f}} & \tilde{Y}
}
\end{displaymath}
be a commutative diagram where $u$ and $w$ are equivalences, and that $f$ is in $\mathbf{P}^{\tau}$. Let $g:C\rightarrow X$ be a cover such that $f\circ g$ is in $\mathbf{P}$. Since $u$ is an equivalence $u\circ g:C\rightarrow\tilde{X}$ is a cover in $\tau$. Moreover $\tilde{f}\circ u\circ g=w\circ f\circ g$. Now $w$ is an equivalence and $f\circ g\in\mathbf{P}$, so $w\circ f\circ g\in\mathbf{P}$. Therefore $w\circ f\in\mathbf{P}^{\tau}$. 

Suppose now that $\tilde{f}$ is in $\mathbf{P}^{\tau}$. Let $\tilde{g}:\tilde{C}\rightarrow\tilde{X}$ be a cover such that $\tilde{f}\circ\tilde{g}\in\mathbf{P}$. Consider the cover $\tilde{g}\times_{\tilde{X}}X:\tilde{C}\times_{\tilde{X}}X\rightarrow X$. Now we have a commutative diagram
\begin{displaymath}
\xymatrix{
\tilde{C}\times_{\tilde{X}}X\ar[d]^{\tilde{C}\times^{h}_{\tilde{X}}u}\ar[rr]^{f\circ (\tilde{g}\times_{\tilde{X}}X)} & & Y\ar[d]^{w}\\
\tilde{C}\ar[rr]^{\tilde{f}\circ\tilde{g}}& & \tilde{Y}
}
\end{displaymath}
Since $\tilde{C}\times_{\tilde{X}}u$ and $w$ are equivalences, $\tilde{f}\circ\tilde{g}$ is in $\mathbf{P}$, and $\mathbf{P}$ is stable under equivalences, $f\circ (\tilde{g}\times_{\tilde{X}}X)$ is in $\mathbf{P}$. Hence $f$ is in $\mathbf{P}^{\tau}$. 
%

Finally suppose that $\mathbf{P}$ is stable under composition. Let $f:X\rightarrow Y$ and $g:Y\rightarrow Z$ be in $\mathbf{P}^{\tau}$. Let $C\rightarrow Y$ be a cover in $\tau$ such that $C\rightarrow Y\rightarrow Z$ is in $\mathbf{P}$. Then $C\times_{Y}X\rightarrow X$ is a cover in $\tau$. Let $D\rightarrow X$ be a cover in $\tau$ such that $D\rightarrow X\rightarrow Y$ is in $\mathbf{P}$. Then $C\times_{Y}D\rightarrow X$ is in $\tau$, and $C\times_{Y}D\rightarrow Z$ is in $\mathbf{P}$. This argument is summed up in the diagram below.
\begin{displaymath}
\xymatrix{
C\times_{Y}D\ar[d]\ar[r] & C\times_{Y}X\ar[d]\ar[r] & C\ar[d]\ar[dr]\\
D\ar[r] & X\ar[r] & Y\ar[r] & Z
}
\end{displaymath}
where both squares, and hence the rectangle are pullbacks, the map $D\rightarrow Y$ is in $\mathbf{P}$, all maps in the left-hand square are in $\tau$, the map $C\rightarrow Y$ is in $\tau$, and the map $C\rightarrow Z$ is in $\mathbf{P}$. 

\end{enumerate}
\end{proof}

\section{Stacks}\label{sec:stacks}

Let us now introduce our formalism for stacks.

\subsection{A Relative Context For Geometry}

\subsubsection{General Sites}
Let us give a very general setup for geometry in the style of \cite{toen2008homotopical}. The material in this subsection also appeared in \cite{kelly2021analytic}.


Let $\mathbf{M}$ be an $(\infty,1)$-category and $\tau$ a pre-topology on $\mathbf{M}$. Recall that the category of \textit{pre-stacks} (in $\infty$-groupoids) on $\mathbf{M}$ is the category of functors $\mathbf{PrStk}(\mathbf{M})\defeq\mathbf{Fun}(\mathbf{M}^{op},\mathbf{Grpd})$ where $\mathbf{Grpd}$ is the category of $\infty$-groupoids, which is presented by the model category of simplicial sets (with its Quillen model structure).
We define
$$\mathbf{Stk}(\mathbf{M},\tau)$$
to be the full subcategory of $\mathbf{PrStk}(\mathbf{M})$ consisting of stacks which satisfy descent for hypercovers in $\tau$.
$\mathbf{Stk}(\mathbf{M},\tau)$ is a reflective subcategory of $\mathbf{PreStk}(\mathbf{M})$. Thus it is an $(\infty,1)$-topos (see \cite{HA}*{Chapter 6}). It is also closed monoidal for the cartesian monoidal structure, in particular one can form the mapping stack $\underline{\mathbf{Map}}(\mathcal{X},\mathcal{Y})$ between any two-stacks. 

\begin{defn}
Let $(\mathbf{M},\tau)$ and $(\mathbf{M}',\tau')$ be sites, and let $F:\mathbf{M}\rightarrow\mathbf{M}'$ be a functor. $F$ is said to be a \textit{continuous functor of sites} if whenever $\{U_{i}\rightarrow U\}$ is a cover in $\tau$, then $\{F(U_{i})\rightarrow F(U)\}$ is a cover in $\tau'$. 
\end{defn}
If $F:\mathbf{M}\rightarrow\mathbf{M}'$ is any functor then the pre-composition functor
$$F^{*}:\mathbf{PreStk}(\mathbf{M}')\rightarrow\mathbf{PreStk}(\mathbf{M})$$
has a left adjoint
$$F_{!}:\mathbf{PreStk}(\mathbf{M})\rightarrow\mathbf{PreStk}(\mathbf{M}')$$
given by left Kan extension.

If $F:(\mathbf{M},\tau)\rightarrow(\mathbf{M}',\tau')$ is a continuous functor of sites then pre-composition restricts to a functor
$$F^{*}:\mathbf{Stk}(\mathbf{M}',\tau')\rightarrow\mathbf{Stk}(\mathbf{M},\tau)$$
This functor has a left adoint
$$F^{\#}:\mathbf{Stk}(\mathbf{M},\tau)\rightarrow\mathbf{Stk}(\mathbf{M}',\tau')$$
given by composing of the stackification functor with $F^{*}$. Exactly as in \cite{kelly2021analytic} Proposition 6.7 (1), we have the following.

\begin{prop}
If $\mathcal{X}\rightarrow\mathcal{Y}$ is an epimorphism in $\mathbf{Stk}(\mathbf{M},\tau)$ then $F^{\#}(\mathcal{X})\rightarrow F^{\#}(\mathcal{Y})$ is an epimorphism in $\mathbf{Stk}(\mathbf{M}',\tau')$. 
\end{prop}

\subsubsection{\v{C}ech Stacks}

Let $(\mathbf{M},\tau)$ be a site. A \v{C}ech \textit{cover} $K_{\bullet}\rightarrow U$ in $\mathbf{M}$ is a hypercover which is the nerve of some collection $\{U_{i}\rightarrow U\}$ in $\tau$. 

\begin{defn}
A \textit{\v{C}ech stack} is a presack $\mathcal{X}$ which satisfies descent for \v{C}ech covers. 
\end{defn}

The full subcategory of $\mathbf{PreStk}(\mathbf{M})$ consisting of \v{C}ech stacks is denoted $\check{\mathbf{Stk}}(\mathbf{M},\tau)$. Again the inclusion 
$$\check{\mathbf{Stk}}(\mathbf{M},\tau)\rightarrow\mathbf{PreStk}(\mathbf{M})$$
has a left adjoint $\mathcal{X}\mapsto\check{\mathcal{X}}$ which we call \textit{\v{C}ech stackification}. It is not clear that some of the topologies of interest will be hyper-subcanonical, only subcanonical for \v{C}ech covers. Thus we have need to work with \v{C}ech stacks. In what follows in this section and beyond, unless specified, `stack' will refer to both (hyper)stacks and \v{C}ech stack, as will the category $\mathbf{Stk}(\mathbf{M},\tau)$.

\subsubsection{(Relative) Geometry Tuples}

We now recall the definition of $(\infty,1)$-geometry tuples from \cite{kelly2021analytic}.

\begin{defn}
Let $(\mathbf{M},\tau)$ be a site. An object $X$ of $\mathbf{M}$ is said to be \textit{admissible} if $\Map(-,X)$ is a stack.
\end{defn}

\begin{defn}
An $(\infty,1)$-\textit{pre-geometry triple} is a triple $(\mathbf{M},\tau,\mathbf{P})$ where
\begin{enumerate}
\item
$\mathbf{M}$ is an $(\infty,1)$-category.
\item
$\tau$ is a Grothendieck pre-topology on $\mathbf{M}$.
\item
$\mathbf{P}$ is a class of maps in $\mathbf{M}$
\end{enumerate}
such that
\begin{enumerate}
\item
If $\{U_{i}\rightarrow U\}$ is a cover in $\tau$, then each $U_{i}\rightarrow U$ is in $\mathbf{P}$.
\item
$\mathbf{P}$ is local for the pre-topology $\tau$.  
\item
The class $\mathbf{P}$ contains isomorphisms, and is closed under pullbacks and compositions.
\end{enumerate}
An $(\infty,1)$-pre-geometry triple is said to be an $(\infty,1)$-\textit{geometry triple} if every object of $\mathbf{M}$ is admissible.
\end{defn}

We will need a relative version of this.

\begin{defn}\label{defn:strongtup}
\begin{enumerate}
\item
A \textit{relative} $(\infty,1)$-\textit{pre-geometry tuple} is a tuple $(\mathbf{M},\tau,\mathbf{P},\mathbf{A})$ where $(\mathbf{M},\tau,\mathbf{P})$ is an $(\infty,1)$-pre-geometry triple, and $\mathbf{A}\subset\mathbf{M}$ is a full subcategory such that if $f:X\rightarrow Y$ is a map in $\mathbf{P}\cap\mathbf{A}$, and $Z\rightarrow Y$ is any map with $Z$ in $\mathbf{A}$, then $X\times_{Y}Z$ is in $\mathbf{A}$. 
\item
A relative $(\infty,1)$-pre-geometry tuple $(\mathbf{M},\tau,\mathbf{P},\mathbf{A})$ is said to be \textit{strong} if whenever $\{U_{i}\rightarrow U\}$ is a cover in $\mathbf{M}$ (not necessarily in $\mathbf{A}$), and $Z\rightarrow U$ is any map with $Z$ in $\mathbf{A}$, then $\{U_{i}\times_{U}Z\rightarrow Z\}$ is a cover in $\tau|_{\mathbf{A}}$.
%

\end{enumerate}
\end{defn}

The idea behind this definition is that $\mathbf{A}$ is the collection of (affine) geometric objects in which we are primarily interested, and $\mathbf{M}$ is a larger class of spaces which permits more categorical constructions than are available in $\mathbf{A}$. The last condition essentially means that for a map in $\mathbf{P}$, the fibres are in $\mathbf{A}$.

Note that in particular if $(\mathbf{M},\tau,\mathbf{P},\mathbf{A})$ is a relative $(\infty,1)$-pre-geometry tuple then $(\mathbf{A},\tau|_{\mathbf{A}},\mathbf{P}|_{\mathbf{A}})$ is an $(\infty,1)$-pre-geometry triple, where $\mathbf{P}|_{\mathbf{A}}$ denotes the class of maps $f:X\rightarrow Y$ in $\mathbf{A}$ which are in $\mathbf{P}$ when regarded as maps in $\mathbf{M}$, and $\tau|_{\mathbf{A}}$ denotes the restriction of $\tau$ to $\mathbf{A}$. 

\begin{rem}
If $(\mathbf{M},\tau,\mathbf{P},\mathbf{A})$ is a relative $(\infty,1)$-geometry tuple, then there is an associated strong relative $(\infty,1)$-geometry tuple defined as follows. Let $\mathbf{P}_{\mathbf{A}}\subset\mathbf{P}$ denote the class of maps $f:X\rightarrow Y$ such that whenever $Z\rightarrow Y$ is a map with $Z\in\mathbf{A}$ then $X\times_{Y}Z$ is in $\mathbf{A}$. Let $\tau_{\mathbf{A}}$ denote the class of covers $\{U_{i}\rightarrow V\}$ in $\tau$ such that each $U_{i}\rightarrow V$ is in $\mathbf{P}_{\mathbf{A}}$. Then $(\mathbf{M},\tau_{\mathbf{A}},\mathbf{P})$ is a strong relative $(\infty,1)$-geometry tuple. 
\end{rem}

\begin{defn}\label{defn:strongtupdesc}
Let $(\mathbf{M},\tau,\mathbf{P},\mathbf{A})$ be a (strong) relative $(\infty,1)$-pre-geometry tuple .
\begin{enumerate}
\item
An object $X$ of $\mathbf{M}$ is said to be (hyper)-$\mathbf{A}$-admissible if the restriction of $\Map(-,X)$ to $\mathbf{A}$ is a (hyper)stack for $\tau|_{\mathbf{A}}$.
\item
$(\mathbf{M},\tau,\mathbf{P},\mathbf{A})$  is said to be a \textit{(strong) relative} $(\infty,1)$-\textit{geometry tuple} if each $X\in\mathbf{A}$ is $\mathbf{A}$-admissible. 
\item 
$(\mathbf{M},\tau,\mathbf{P},\mathbf{A})$  is said to be a \textit{hyper (strong) relative} $(\infty,1)$-\textit{geometry tuple} if each $X\in\mathbf{A}$ is hyper-$\mathbf{A}$-admissible. 
\end{enumerate}
\end{defn}

%
%
For a  strong relative $(\infty,1)$-geometry tuple $(\mathbf{M},\tau,\mathbf{P},\mathbf{A})$ $\mathbf{A}$-stacks embed fully faithfully in $\mathbf{M}$-stacks.

\begin{prop}
Let $(\mathbf{M},\tau,\mathbf{P},\mathbf{A})$ be a relative $(\infty,1)$-geometry tuple.
\begin{enumerate}
\item
If $\mathcal{X}\rightarrow\mathcal{Y}$ is an epimorphism in $\mathbf{Stk}(\mathbf{A},\tau|_{\mathbf{A}})$ then $i^{\#}(\mathcal{X})\rightarrow i^{\#}(\mathcal{Y})$ is an epimorphism in $\mathbf{Stk}(\mathbf{M},\tau)$.
\item
If the relative $(\infty,1)$-geometry tuple is strong then $i^{\#}$ is fully faithful. 
\end{enumerate}
\end{prop}

\subsubsection{Geometric Stacks}

An important subclass of stacks is furnished by geometric stacks. They are the stacks which admit an $n$-coskeletal hypercover for some $n$. Most examples of interest fall into this subclass. We fix  a relative $(\infty,1)$-pre-geometry tuple $(\mathbf{M},\tau,\mathbf{P},\mathbf{A})$.

\begin{defn}[\cite{toen2008homotopical} Definition 1.3.3.1]
\begin{enumerate}
\item
A prestack $\mathcal{X}$ in $\mathbf{Stk}(\mathbf{A},\tau|_{\mathbf{A}})$ is $(-1)$-\textit{geometric} if it is of the form $\mathcal{X}\cong\Map(-,M)$ for some $\mathbf{A}$-admissible $M\in\mathbf{M}$.
\item
A morphism $f:\mathcal{X}\rightarrow\mathcal{Y}$ in $\mathbf{Stk}(\mathbf{A},\tau|_{\mathbf{A}})$ is $(-1)$-\textit{representable} if for any map $U\rightarrow\mathcal{Y}$ with $U\in\mathbf{A}$ the pullback $\mathcal{X}\times_{\mathcal{Y}}U$ is $(-1)$-geometric.
\item
A morphism $f:\mathcal{X}\rightarrow\mathcal{Y}$ in  $\mathbf{Stk}(\mathbf{A},\tau|_{\mathbf{A}})$ is in $(-1)-\mathbf{P}$ if it is $(-1)$-representable and for any map $U\rightarrow\mathcal{Y}$ with $U\in\mathbf{A}$, the map $\mathcal{X}\times_{\mathcal{Y}}U\rightarrow U$ is in $\mathbf{P}$.
\end{enumerate}
\end{defn}

\begin{defn}[\cite{toen2008homotopical} Definition 1.3.3.1]\label{defn:ngeom}
Let $n\ge0$.
\begin{enumerate}
\item
Let $\mathcal{X}$ be in  $\mathbf{Stk}(\mathbf{A},\tau|_{\mathbf{A}})$. An $n$\textit{-atlas} for $\mathcal{X}$ is a set of morphisms $\{U_{i}\rightarrow \mathcal{X}\}$ such that 
\begin{enumerate}
\item
Each $U_{i}$ is $(-1)$-geometric.
\item
Each map $U_{i}\rightarrow\mathcal{X}$ is in $(n-1)$-$\mathbf{P}$ and is an epimorphism of stacks.
\end{enumerate}
\item
$\mathcal{X}$ is $n$-\textit{geometric}  if
\begin{enumerate}
\item
The map $\mathcal{X}\rightarrow\mathcal{X}\times\mathcal{X}$ is $(n-1)$-representable.
\item
$\mathcal{X}$ admits an $n$-atlas.
\end{enumerate}
\item
A morphism of prestacks $f:\mathcal{X}\rightarrow\mathcal{Y}$ is $n$-\textit{representable} if for any map  $U\rightarrow \mathcal{Y}$ with $U\in\mathbf{A}$ the pullback $\mathcal{X}\times_{\mathcal{Y}}U$ is $n$-geometric.
\item
A morphism of prestacks $\mathcal{X}\rightarrow\mathcal{Y}$ is in $n$-$\mathbf{P}$ if it is $n$-representable, and for any map $U\rightarrow\mathcal{Y}$ with $U\in\mathbf{A}$, there is an $n$-atlas $\{U_{i}\rightarrow\mathcal{X}\times_{\mathcal{Y}}U\}_{i\in\mathcal{I}}$ such that each map $U_{i}\rightarrow U$ is in $\mathbf{P}$.
\end{enumerate}
\end{defn}

Consider now the absolute case, i.e. an $(\infty,1)$-pre-geometry triple $(\mathbf{M},\tau,\mathbf{P})$. Let $\mathbf{Stk}_{n}$ denote the class of $n$-geometric stacks. It is closed under finite limits.


For a class of maps $\mathbf{P}$ in $\mathbf{A}$ let $\mathbf{P}^{\coprod}$ denote the class of maps between disjoint unions
$$f:\coprod_{i\in\mathcal{I}}X_{i}\rightarrow\coprod_{j\in\mathcal{J}}Y_{j}$$
such that there is a map $\alpha:\mathcal{I}\rightarrow\mathcal{J}$, and for each $i\in\mathcal{I}$ a map $p_{i}: X_{i}\rightarrow Y_{\alpha(i)}$ in $\mathbf{P}$ inducing the map $f$.

\begin{thm}
A stack $\mathcal{X}$ is $n$-geometric for some $n$ if and only if for some $m$ there is an $m$-skeletal hypercover 
$$K_{\bullet}\rightarrow\mathcal{X}$$
where each $K_{n}$ is a disjoint union of representables, and all face and degeneracy maps are in $\mathbf{P}^{\coprod}$. Moreover in this case, with $K_{0}=\coprod U_{i_{0}}$ we have that each $U_{i_{0}}\rightarrow\mathcal{X}$ is in $\mathbf{P}^{rep}$. We call such a cover \textit{pseudo-representable}.
\end{thm}

\begin{proof}
This is the content of Proposition 4.1 and Theorem 4.7 in \cite{MR3033634}. There it is proven in the language of model categories, but the proof works for $(\infty,1)$-categories.
\end{proof}

\subsubsection{Functors Between Categories of Geometric Stacks}

Let $(\mathbf{M},\tau,\mathbf{P},\mathbf{A})$ and $(\mathbf{M}',\tau',\mathbf{P}',\mathbf{A}')$ be relative $(\infty,1)$-geometry tuples. 

\begin{defn}
A morphism $(\mathbf{M},\tau,\mathbf{P},\mathbf{A})\rightarrow(\mathbf{M}',\tau',\mathbf{P}',\mathbf{A}')$ of relative $(\infty,1)$-geometry tuples is a functor $F:\mathbf{M}\rightarrow\mathbf{M}'$ such that
\begin{enumerate}
\item
$F:(\mathbf{M},\tau)\rightarrow(\mathbf{M}',\tau')$ is a continuous morphism of sites.
\item
If $f\in\mathbf{P}|_{\mathbf{A}}$ then $F(f)\in\mathbf{P}|_{\mathbf{A}}$. 
\item
If
\begin{displaymath}
\xymatrix{
U\times_{W}V\ar[r]\ar[d] & U\ar[d]^{f}\\
V\ar[r] & W
}
\end{displaymath}
is a pullback diagram in $\mathbf{M}$ with $f$ in $\mathbf{P}|_{\mathbf{A}}$, then 
\begin{displaymath}
\xymatrix{
F(U\times_{W}V)\ar[r]\ar[d] & F(U)\ar[d]^{F(f)}\\
F(V)\ar[r] & F(W)
}
\end{displaymath}
is a pullback diagram in $\mathbf{M}'$.
\end{enumerate}
\end{defn}

Let  $F:(\mathbf{M},\tau,\mathbf{P},\mathbf{A})\rightarrow(\mathbf{M}',\tau',\mathbf{P}',\mathbf{A}')$ be a morphism of relative $(\infty,1)$-geometry tuples. If $\mathcal{Y}\in\mathbf{Stk}(\mathbf{A}',\tau'|_{\mathbf{A}'})$, then $\mathcal{Y}\circ F\in\mathbf{Stk}(\mathbf{A},\tau|_{\mathbf{A}})$. Recall the restriction functor
$$F^{*}:\mathbf{Stk}(\mathbf{A}',\tau'|_{\mathbf{A}'})\rightarrow\mathbf{Stk}(\mathbf{A},\tau|_{\mathbf{A}})$$
and its left adjoint
$$F^{\#}:\mathbf{Stk}(\mathbf{A},\tau|_{\mathbf{A}})\rightarrow\mathbf{Stk}(\mathbf{A}',\tau'|_{\mathbf{A}'})$$

Let $(\mathbf{M},\tau,\mathbf{P},\mathbf{A})$ be a relative $(\infty,1)$-geometry tuple such that $(\mathbf{M},\tau,\mathbf{P})$ is an $(\infty,1)$-geometry triple. Suppose further that $\mathbf{A}$ has is closed under finite limits in $\mathbf{M}$. In particular $(\mathbf{A},\tau|_{\mathbf{A}},\mathbf{P}|_{\mathbf{A}})$ is an $(\infty,1)$-geometry triple. Moreover the inclusion functor $\mathbf{A}\rightarrow\mathbf{M}$ determines

\begin{lem}[\cite{kelly2021analytic} Lemma 6.12]
Let $(\mathbf{M},\tau,\mathbf{P},\mathbf{A})$ be a relative $(\infty,1)$-geometry tuple such that $(\mathbf{M},\tau,\mathbf{P})$ is an $(\infty,1)$-geometry triple. Let $f:\mathcal{X}\rightarrow\mathcal{Y}$ be a map in $\mathbf{Stk}_{n}(\mathbf{A},\tau|_{\mathbf{A}},\mathbf{P}|_{\mathbf{A}})$ which is in $n-\mathbf{P}|_{\mathbf{A}}$. If $\mathbf{A}$ has all finite limits which are preserved by the inclusion $\mathbf{A}\rightarrow\mathbf{M}$ then $i^{\#}(f)$ is in $n-\mathbf{P}$.
\end{lem}

In forthcoming work \cite{RhiannonRepresentability}, Savage extends this result to prove the following.

\begin{lem}[\cite{RhiannonRepresentability} 2.4.4]
Let $F:(\mathbf{M},\tau,\mathbf{P},\mathbf{A})\rightarrow(\mathbf{M}',\tau',\mathbf{P}',\mathbf{A}')$ be a morphism of relative $(\infty,1)$-geometry tuples. Then $F^{\#}$ sends $n$-geometric stacks to $n$-geometric stacks,
and $n-\mathbf{P}$-representable morphisms to $n-\mathbf{P}'$-representable morphisms. 
\end{lem}


\subsubsection{Schemes}

Let $(\mathbf{M},\tau,\mathbf{P},\mathbf{A})$ be a strong $(\infty,1)$-pre-geometry tuple.

\begin{defn}
A map $f:\mathcal{X}\rightarrow\mathcal{Y}$ in an $(\infty,1)$-category $\mathbf{M}$ is said to be a \textit{monomorphism} if the map $\mathcal{X}\rightarrow\mathcal{X}\times_{\mathcal{Y}}\mathcal{X}$ is an isomorphism.
\end{defn}

It is easy to see that the class $\mathbf{Mon}$ of all monomorphisms is in fact itself closed under composition and pullbacks, and contains all isomorphisms. 

\begin{defn}
Let $(\mathbf{M},\tau,\mathbf{P},\mathbf{A})$ be a relative $(\infty,1)$-geometry tuple. An $n$-geometric stack $\mathcal{X}$ is said to be an $n$-\textit{scheme} if 
\begin{enumerate}
\item
it is a $(-1)$-geometric stack for $n=-1$.
\item
for $n\ge 0$ it has an  $n$-atlas $\{U_{i}\rightarrow\mathcal{X}\}_{i\in\mathcal{I}}$ such that 
\begin{enumerate}
\item
each $U_{i}\rightarrow\mathcal{X}$ is a monomorphism in $\mathbf{PreStk}(\mathbf{A},\tau|_{\mathbf{A}})$,
\item
for each $n\in\mathbb{N}$ and each $(i_{1},\ldots,i_{n})\in\mathcal{I}^{n}$, each $U_{i_{1}}\times_{\mathcal{X}}U_{i_{2}}\times_{\mathcal{X}}\ldots\times_{\mathcal{X}}U_{i_{n}}$ is an $(n-1)$-scheme.
\end{enumerate}
\end{enumerate}
\end{defn}

In \cite{kelly2021analytic} we had the additional condition that diagonal map $\mathcal{X}\rightarrow\mathcal{X}\times \mathcal{X}$ is a monomorphism in $\mathbf{PreStk}(\mathbf{A},\tau|_{\mathbf{A}})$, which is also present in the definition of scheme in \cite{toen2008homotopical} - however we believe this condition is unnecessary, and in fact virtually never satisfied.

The full subcategory of $\mathbf{Stk}(\mathbf{A},\tau|_{\mathbf{A}})$ consisting of $n$-schemes is denoted $\mathbf{Sch}_{n}(\mathbf{M},\tau,\mathbf{P},\mathbf{A})$. We also write \[\mathbf{Sch}(\mathbf{M},\tau,\mathbf{P},\mathbf{A})\defeq\bigcup_{n=-1}^{\infty}\mathbf{Sch}_{n}(\mathbf{M},\tau,\mathbf{P},\mathbf{A}).\]

\section{Sheaves}\label{sec:sheaves} 

In this section we mostly introduce conventions for discussing sheaves (of categories) on sites and stacks.

\subsection{Sheaves of categories}

We begin with some notation.

\begin{defn}
Let $\kappa$ be an ordinal. Denote by 
\begin{enumerate}
\item
$\mathbf{Cat}$ the $(\infty,1)$-category of $(\infty,1)$-categories with morphisms being functors.
\item
$\mathbf{Cat^{L}}$ the $(\infty,1)$-category of $(\infty,1)$-categories with morphisms being left-adjoint functors.
\item
$\mathbf{Cat^{L,\kappa}}$ the $(\infty,1)$-category of $(\infty,1)$-categories with morphisms being left-adjoint functors whose right adjoints commute with $\kappa$-filtered colimits.
\item
$\mathbf{Pr^{L}}\subset\mathbf{Cat^{L}}$/ $\mathbf{Pr^{L,\kappa}}\subset\mathbf{Cat^{L,\kappa}}$ the full subcategories consisting of locally presentable categories.
\end{enumerate}
We will write $\mathbf{Cat}^{*}$ when we are referring to one of the above categories.
\end{defn}

\begin{defn}
Let $\mathbf{A}$ be an $(\infty,1)$-category. A $\mathbf{Cat}^{*}$\textit{-presheaf on }$\mathbf{A}$ is an $(\infty,1)$-functor
$$\mathbf{Q}:\mathbf{A}^{op}\rightarrow\mathbf{Cat}^{*}$$
\end{defn}

Let $\mathbf{Cat}$ denote the $(\infty,1)$-category of $(\infty,1)$-categories. For each $\mathbf{Cat}^{*}$ there is an inclusion functor $\mathbf{Cat}^{*}\rightarrow\mathbf{Cat}$ which by \cite{HTT} Proposition 5.5.3.13, in the case of $\mathbf{Pr^{L}}$, commutes with limits. For $\mathbf{Q}$ a $\mathbf{Cat}^{*}$-presheaf on $\mathbf{A}$, we denote by $\mathbf{Q}^{o}$ the functor $\mathbf{A}^{op}\rightarrow\mathbf{Cat}$ given by post-composition with this inclusion. 

\begin{defn}
Let $\mathbf{A}$ be an $(\infty,1)$-category and $\mathbf{Q}$ a $\mathbf{Cat}^{*}$-presheaf on $\mathbf{A}$. A \textit{sub-presheaf of }$\mathbf{Q}$, denoted $\mathbf{N}\subset\mathbf{Q}$, is a functor
$$\mathbf{N}:\mathbf{A}^{op}\rightarrow\mathbf{Cat}$$
together with an object-wise fully faithful natural transformation $\mathbf{N}\rightarrow\mathbf{Q}^{o}$.
\end{defn}

We will often restrict to subcategories which are not preserved by the right adjoints. Let us give our first important examples of sub-presheaves.

\begin{example}\label{ex:compactsubpresheaf}
Let $\mathbf{Q}$ be a $\mathbf{Cat^{L,\kappa}}$ presheaf on $\mathbf{Q}$. For $\kappa$ a regular cardinal, let $\mathbf{Q}^{\kappa-cpct}$ denote the full subcategory of $\mathbf{Q}$ consisting of $\kappa$-compact objects. Then this defines a sub-presheaf of $\mathbf{Q}$.
\end{example}

\subsubsection{Presheaves on Stacks}

Consider the $(\infty,1)$-category \[\mathbf{PreStk}(\mathbf{A})=\mathbf{Fun}(\mathbf{A}^{op},\mathbf{sSet})\] of pre-stacks on $\mathbf{A}$. Via the Yoneda embedding, $\mathbf{PreStk}(\mathbf{A})$ is the free cocompletion of $\mathbf{A}$. Since $\mathbf{Pr^{L}}$ and $\mathbf{Cat}$ are complete, for $\mathbf{Cat}^{*}\in\{\mathbf{Cat},\mathbf{Pr^{L}}\}$, any $\mathbf{Cat}^{*}$-presheaf of $(\infty,1)$-categories on $\mathbf{A}$ uniquely extends to a limit preserving functor, which we also denote by $\mathbf{Q}$,
$$\mathbf{Q}:\mathbf{PreStk}(\mathbf{A})^{op}\rightarrow\mathbf{Cat}^{*}$$
In particular $\mathbf{Q}$ extends to a $\mathbf{Cat}^{*}$-presheaf on $\mathbf{PreStk}(\mathbf{A})$, and in the $\mathbf{Pr^{L}}$ case if $f:\mathcal{X}\rightarrow\mathcal{Y}$ is a map of pre-stacks, then we get an adjunction
$$\adj{f_{\mathbf{Q}}^{*}}{\mathbf{Q}(\mathcal{Y})}{\mathbf{Q}(\mathcal{X})}{f_{\mathbf{Q},*}}$$
The category $\mathbf{Q}(\mathcal{X})$ is computed as the limit
$$\mathbf{Q}(\mathcal{X})=\mathbf{lim}_{\mathbf{Map}(-,X)\rightarrow\mathcal{X}}\mathbf{Q}(X)$$
This limit can be computed in either $\mathbf{Pr^{L}}$ or $\mathbf{Cat}$. In particular if $\mathbf{N}\subset\mathbf{Q}$ is a subpresheaf, we can define $\mathbf{N}:\mathbf{PreStk}(\mathbf{A})\rightarrow\mathbf{Cat}$ by the same formula, and we still have a natural transformation $\mathbf{N}(\mathcal{X})\rightarrow\mathbf{Q}(\mathcal{X})$. Using the construction of the mapping space in the limit category in \cite{dhillonzsamboki} Proposition 4.12, for each $\mathcal{X}$, $\mathbf{N}(\mathcal{X})\rightarrow\mathbf{Q}(\mathcal{X})$ is fully faithful. Thus we get a sub-presheaf $\mathbf{N}\subset\mathbf{Q}$ on $\mathbf{PreStk}(\mathbf{A})$. In what follows we will usually suppress the notation $f_{\mathbf{Q}}^{*}$ (resp. $f_{\mathbf{Q},*}$) to $f^{*}$ (resp. $f_{*}$).

\subsection{Conservative, $\kappa$-Filtered, and Base-Change Maps}

In this section we discuss some important technical local properties of presheaves of categories. 

\subsubsection{Conservative Maps}
Let $\mathbf{Q}:\mathbf{A}^{op}\rightarrow\mathbf{Cat^{L}}$ be a presheaf of $(\infty,1)$-categories on an $(\infty,1)$-category $\mathbf{A}$.

\begin{defn}
A map $f:X\rightarrow Y$ is said to be  
\begin{enumerate}
\item
$\mathbf{Q}$-\textit{conservative} if $f^{*}$ is a conservative functor.
\item
$\mathbf{Q}$-\textit{pushforward conservative} if $f_{*}$ is a conservative functor.
\end{enumerate}
\end{defn}

\subsubsection{Formally $\kappa$-Filtered Maps}

\begin{defn}
Let $L:\mathbf{C}\rightarrow\mathbf{D}$ be a left adjoint functor of $(\infty,1)$-categories where both $\mathbf{C}$ and $\mathbf{D}$ have $\kappa$-small sifted limits. $L$ is said to be  
\begin{enumerate}
\item
\textit{weakly} $\kappa$-\textit{filtered} if $L$ commutes with $\kappa$-small products.
\item
\textit{formally} $\kappa$-\textit{filtered} if for any $\kappa$-small  diagram $\kappa:\mathcal{I}\rightarrow\mathbf{C}$ the natural transformation
$$L(\lim_{i}C_{i})\rightarrow \lim_{i}L(C_{i})$$
is an equivalence, where the left-hand limit is computed in $\mathbf{C}$, and the right-hand limit in $\mathbf{D}$.
\end{enumerate}


\end{defn}

\begin{rem}
If $\mathbf{C}$ and $\mathbf{D}$ are stable then $F$ is formally $\kappa$-filtered if and only if it is weakly $\kappa$-filtered. Moreover if $\kappa=\aleph_{0}$ then all left exact functors between stable categories are $\kappa$-filtered, since in this case the limits are finite and left exact functor between stable categories commute with all limits.
\end{rem}

\begin{defn}
Let $\mathbf{Q}:\mathbf{A}^{op}\rightarrow\mathbf{Cat}^{\mathbf{L}}$ be a presheaf of categories. A map $f:\mathcal{X}\rightarrow\mathcal{Y}$ in $\mathbf{PreStk}(\mathbf{A})$ is said to be $\mathbf{Q}$-\textit{formally }$\kappa$-filtered if the functor $f^{*}:\mathbf{Q}(\mathcal{Y})\rightarrow\mathbf{Q}(\mathcal{X})$ is formally $\kappa$-filtered.
\end{defn}

\subsubsection{Base-Change Maps}
Most of this subsubsection is a straightforward adaptation of \cite{MR3037900} Section 1.3. 
Let 
\begin{displaymath}
\xymatrix{
X\times_{Z}Y\ar[d]^{f'}\ar[r]^{g'} & X\ar[d]^{f}\\
Y\ar[r]^{g} & Z
}
\end{displaymath}
be a pullback diagram in $\mathbf{A}$. Let $\mathcal{F}\in \mathbf{Q}(X)$ By adjunction we have the following equivalences
\begin{align*}
\mathbf{Map}(g^{*}f_{*}\mathcal{F}_,f'_{*}g'^{*}\mathcal{F})&\cong\mathbf{Map}(f_{*}\mathcal{F},g_{*}f'_{*}g'^{*}\mathcal{F})\\
&\cong\mathbf{Map} (f_{*}\mathcal{F},f_{*}g'_{*}g'^{*}\mathcal{F})
\end{align*}
Thus the unit $\mathcal{F}\rightarrow g'_{*}g'^{*}\mathcal{F}$ induces a map $g^{*}f_{*}\mathcal{F}\rightarrow f'_{*}g'^{*}\mathcal{F}$

\begin{defn}
\begin{enumerate}
\item
$(g,\mathcal{F},f)$ is said to \textit{have the} $\mathbf{Q}$-\textit{base change property} if the map $g^{*}f_{*}\mathcal{F}\rightarrow f'_{*}g'^{*}\mathcal{F}$ is an equivalence. 
\item
$(\mathcal{F},f)$ is said  to \textit{have the} $\mathbf{Q}$-\textit{base change property} if the map $g^{*}f_{*}\mathcal{F}\rightarrow f'_{*}g'^{*}\mathcal{F}$ is an equivalence for any $g:Y\rightarrow Z$.
\item
$(g,f)$ is said to \textit{have the} $\mathbf{Q}$-\textit{base change property} if the map $g^{*}f_{*}\mathcal{F}\rightarrow f'_{*}g'^{*}\mathcal{F}$ is an equivalence for any $\mathcal{F}\in\mathbf{Q}(X)$
\item
$f$ is said  to \textit{have the} $\mathbf{Q}$-\textit{base change property} $(g,f)$ has the $\mathbf{Q}$-base change property for any $g$.
\end{enumerate}
\end{defn}

\begin{defn}
Let 
\begin{displaymath}
\xymatrix{
X\times_{Z}Y\ar[d]^{f'}\ar[r]^{g'} & X\ar[d]^{f}\\
Y\ar[r]^{g} & Z
}
\end{displaymath}
be a pullback square. It is said to have the \textit{vertical base-change property} if $(g,f)$ has the base-change property, the  \textit{horizontal base-change property} if $(f,g)$ has the base-change property, and the \textit{symmetric base-change property} if it has both the vertical base-change property and the horizontal base-change property.
\end{defn}

The base-change property can be tested on affines. Indeed the following proof works identically to \cite{MR3037900} Proposition 1.3.6.

\begin{lem}
Let $f:X\rightarrow Z$ be a map and let $\mathcal{F}\in\mathbf{Q}(X)$. The following are equivalent.
\begin{enumerate}
\item
 $(\mathcal{F},f)$ has the $\mathbf{Q}$-base-change property
 \item
 for any $g:U\rightarrow Z$ with $U$ affine $(g,\mathcal{F},f)$ has the $\mathbf{Q}$-base-change property.
 \item
 For any $h:V\rightarrow U$ and $g:U\rightarrow Z$ with $U,V$ affine, $(h,\mathcal{F}_{V}',\pi_{V}')$ satisfies the $\mathbf{Q}$-base change property, where $\mathcal{F}_{V}'$ is the pullback of $\mathcal{F}$ to $U\times_{Y}X$, and $\pi_{V}':U\times_{Y}X\rightarrow U$ is the projection.
 \end{enumerate}
\end{lem}



\subsubsection{$\mathbf{Q}$-Embeddings}

\begin{defn}
A map $f:X\rightarrow Y$ is said to be a
\begin{enumerate}
\item
 $\mathbf{Q}$-\textit{embedding} if the adjunction
$$\adj{f^{*}}{\mathbf{Q}(Y)}{\mathbf{Q}(X)}{f_{*}}$$
realises $\mathbf{Q}(X)$ as a reflective subcategory of $\mathbf{Q}(Y)$. 
\item
 $\mathbf{Q}$-\textit{coembedding} if the adjunction
$$\adj{f^{*}}{\mathbf{Q}(Y)}{\mathbf{Q}(X)}{f_{*}}$$
realises $\mathbf{Q}(Y)$ as a coreflective subcategory of $\mathbf{Q}(X)$. 
\item
$\mathbf{Q}$-equivalence if the adjunction
$$\adj{f^{*}}{\mathbf{Q}(Y)}{\mathbf{M}(X)}{f_{*}}$$
is an equivalence.
\end{enumerate}
\end{defn}

Let $f:\mathcal{X}\rightarrow\mathcal{Y}$ be a map of pre-stacks. The \textit{image of }$f$, denoted $\mathrm{im}(f)$, is the geometric realisation of the nerve of $f$
$$\mathrm{im}(f)\defeq|N(f)|$$
where $N(f)$ is the simplicial object given in degree $n$ by $\mathcal{X}^{\times_{\mathcal{Y}}n}$. Note that $f$ factorises through a natural map $d_{f}:\mathrm{im}(f)\rightarrow\mathcal{Y}$.

\begin{defn}
A map of pre-stacks $f:\mathcal{X}\rightarrow\mathcal{Y}$ is said to be 
\begin{enumerate}
\item
a $\mathbf{Q}$-\textit{descent morphism} if $d_{f}:\mathrm{im}(f)\rightarrow\mathcal{Y}$ is a $\mathbf{Q}$-embedding
\item
a $\mathbf{Q}$-\textit{codescent morphism} if $d_{f}:\mathrm{im}(f)\rightarrow\mathcal{Y}$ is a $\mathbf{Q}$-embedding
\item
an \textit{effective }$\mathbf{Q}$-\textit{descent morphism} $d_{f}:\mathrm{im}(f)\rightarrow\mathcal{Y}$ is a $\mathbf{Q}$-equivalence.
\end{enumerate}
\end{defn}

\begin{defn}
A map $f:\mathcal{X}\rightarrow \mathcal{Z}$ be a map of prestacks is said to be a 
\begin{enumerate}
\item
\textit{universal }$\mathbf{Q}$-\textit{embedding} if for any map $Y\rightarrow\mathcal{Z}$ with $Y$ affine, $f':Y\times_{\mathcal{Z}}\mathcal{X}\rightarrow Y$ is a $\mathbf{Q}$-embedding.
\item
\textit{universal }$\mathbf{Q}$-\textit{coembedding} if for any map $Y\rightarrow\mathcal{Z}$ with $Y$ affine, $f':Y\times_{\mathcal{Z}}\mathcal{X}\rightarrow Y$ is a $\mathbf{Q}$-coembedding.
\item
\textit{universal }$\mathbf{Q}$-\textit{descent map} if for any map $Y\rightarrow\mathcal{Z}$ with $Y$ affine, $f':Y\times_{\mathcal{Z}}\mathcal{X}\rightarrow Y$ is a $\mathbf{Q}$-descent map.
\item 
\textit{universal }$\mathbf{Q}$-\textit{codescent map} if for any map $Y\rightarrow\mathcal{Z}$ with $Y$ affine, $f':Y\times_{\mathcal{Z}}\mathcal{X}\rightarrow Y$ is a $\mathbf{Q}$-effective codescent map.
\item
\textit{universal }$\mathbf{Q}$-\textit{effective descent map} if for any map $Y\rightarrow\mathcal{Z}$ with $Y$ affine, $f':Y\times_{\mathcal{Z}}\mathcal{X}\rightarrow Y$ is a $\mathbf{Q}$-effective descent map.
\end{enumerate}
\end{defn}

$\mathbf{Q}$-embeddings are clearly push-forward conservative.

\begin{lem}\label{lem:pullbackembedding}
Let
\begin{displaymath}
\xymatrix{
X\times_{Z}Y\ar[d]^{f'}\ar[r]^{g'} & X\ar[d]^{f}\\
Y\ar[r]^{g} & Z
}
\end{displaymath}
be a fibre-product diagram which is horizontally base change.
\begin{enumerate}
\item
Suppose that $f$ is a $\mathbf{Q}$-coembedding and that $g$ is push-forward conservative. Then $f'$ is a $\mathbf{Q}$-coembedding.
\item
Suppose that $f$ is a $\mathbf{Q}$-embedding and that $g'$ is push-forward conservative. Then $f'$ is a $\mathbf{Q}$-embedding.
\item
Suppose that $f$ is a $\mathbf{Q}$-equivalence and that $g$ is push-forward conserivative. Then the square is symmetrically base-change..
\end{enumerate}
\end{lem}

\begin{proof}
\begin{enumerate}
\item
$f$ being an $\mathbf{Q}$-coembedding amounts to the unit $\mathrm{Id}\rightarrow f_{*}f^{*}$ being an equivalence. We need to show that $\mathrm{Id}\rightarrow (f')_{*}(f')^{*}$ is an equivalence. Since $g$ is push-forward conservative it suffices to show that $g_{*}\rightarrow g_{*}(f')_{*}(f')^{*}$ is an equivalence. Now we have
$$g_{*}(f')_{*}(f')^{*}\cong f_{*}(g')_{*}(f')^{*}\cong f_{*}f^{*}g_{*}\cong g_{*}$$
using base-change and the fact that $\mathrm{Id}\rightarrow f_{*}f^{*}$ is an equivalence.
\item
We need to show that $(f')^{*}(f')_{*}\rightarrow Id$ is an equivalence. The proof works by applying $(g')_{*}$ to both sides, and unravelling the base-change formula.
\item
We need to show that
$$g^{*}f_{*}\rightarrow(f')_{*}(g')^{*}$$
is an equivalence. It suffices to show that
$$g_{*}g^{*}f_{*}\rightarrow g_{*}f'_{*}g'_{*}\cong f_{*}(g')_{*}(g')^{*}$$
is an equivalence. Now $f^{*}$ is conservative, so it suffices to prove that 
$$f^{*}g_{*}g^{*}f_{*}\rightarrow g_{*}f'_{*}g'_{*}\cong(g')_{*}(g')^{*}$$
is an equivalence. But using base change for $g$ we have
$$f^{*}g_{*}g^{*}f_{*}\cong (g')_{*}(f')^{*}g^{*}f_{*}\cong(g')_{*}(g')^{*}f^{*}f_{*}\cong(g')_{*}(g')^{*}$$
as required.
\end{enumerate}
\end{proof}



\begin{lem}\label{lem:reversepullback}
Let
\begin{displaymath}
\xymatrix{
X\times_{Z}Y\ar[d]^{f'}\ar[r]^{g'} & X\ar[d]^{f}\\
Y\ar[r]^{g} & Z
}
\end{displaymath}
be a fibre-product diagram in which $g$ is conservative.
\begin{enumerate}
\item
If the square is vertically base change and $f'$ is a $\mathbf{Q}$-coembedding then $f$ is a $\mathbf{Q}$-coembedding.
 \item
If the square is vertically base change, $g'$ is conservative, and $f'$ is push-forward conservative, then $f$ is push-forward conservative.
 \end{enumerate}
 In particular if the square is vertically base change, $f'$ is a $\mathbf{Q}$-equivalence, and $g'$ is conservative, then $f$ is an $\mathbf{Q}$-equivalence.
\end{lem}

\begin{proof}
\begin{enumerate}
\item
We need to show that $\mathrm{Id}\rightarrow f_{*}f^{*}$ is an equivalence. This is equivalent to showing $g^{*}\rightarrow g^{*}f_{*}f^{*}$ is an equivalence. By vertical base change we have
$$g^{*}f_{*}f^{*}\cong(f')_{*}(g')^{*}f^{*}\cong (f')_{*}(f')^{*}g^{*}\cong g^{*}$$
as required.
\item
Let $h$ be a map in $\mathbf{Q}(X)$ such that $f_{*}(h)$ is an equivalence. Then $g^{*}f_{*}(h)\cong (f')_{*}(g')^{*}(h)$ is an equivalence. Since $f'$ is push-forward conservative and $g'$ is conservative by assumption, $h$ is an equivalence.
\end{enumerate}
\end{proof}



\subsubsection{Presheaves of Affine Type}

\begin{defn}
$\mathbf{Q}$ is said to be \textit{of affine type} if all maps in $\mathbf{A}$ have the $\mathbf{Q}$-base change property, and are $\mathbf{Q}$-pushforward conservative. We denote by
$$\mathbf{PreShv}^{\mathrm{aff}}(\mathbf{A},\mathbf{Cat^{L}})\subset\mathbf{PreShv}(\mathbf{A},\mathbf{Cat^{L}})$$
the full subcategory consisting of presheaves of affine type.
\end{defn}

\begin{example}
Let $(\mathbf{C},\mathbf{D},\theta)$ be an $(\infty,1)$-algebra context. Consider the $\mathbf{Pr^{L}}$-presheaf $\mathbf{QCoh}$ on $\mathbf{Alg_{D}}(\mathbf{C})$ which, roughly speaking, assigns to $A$ the $(\infty,1)$-category ${}_{A}\mathbf{Mod}$, and to $f:A\rightarrow B$ the adjunction
$$\adj{B\otimes^{\mathbb{L}}_{A}(-)}{{}_{A}\mathbf{Mod}}{{}_{B}\mathbf{Mod}}{|-|_{A}}$$
This is of affine type. 
\end{example}

\begin{lem}
If $\mathbf{Q}$ is of affine type then universal $\mathbf{Q}$-embeddings are $\mathbf{Q}$-base change.
\end{lem}

\begin{proof}
This follows immediately from Lemma \ref{lem:pullbackembedding}.
\end{proof}

\section{Presheaves with $t$-structures}\label{sec:presheavet}

Let $\mathbf{Cat}^{\mathbf{L},t}$ the category whose objects are stable $(\infty,1)$-categories equipped with a $t$-structure, and whose morphisms are left adjoint functors which are right exact for the $t$-structure with left exact right adjoints. We also denote by $\mathbf{Cat}^{\mathbf{L},t,re}$ the wide subcategory of $\mathbf{Cat}^{\mathbf{L},t}$ consisting of functors whose right adjoints are exact for the $t$-structures. We will typically generally work with presheaves of categories with $t$-structures with exact right adjoints. 

\begin{defn}
A \textit{presheaf with }$t$-\textit{structure} on an $(\infty,1)$-category $\mathbf{A}$ is a functor
$$\mathbf{Q}:\mathbf{A}^{op}\rightarrow\mathbf{Cat}^{\mathbf{L},t}$$
\end{defn}

Note that if
$$\mathbf{Q}:\mathbf{A}^{op}\rightarrow\mathbf{Cat}^{\mathbf{L},t}$$
is a presheaf with $t$-structure then we may also regard it as a $\mathbf{Cat}^{\mathbf{L}}$-presheaf. Consider the extended presheaf
$$\mathbf{Q}:\mathbf{PreStk}(\mathbf{A})^{op}\rightarrow\mathbf{Cat^{L}}$$
For $\mathcal{X}\in\mathbf{PreStk}(\mathbf{A})^{op}$ denote by $\mathbf{Q}_{\ge0}(\mathcal{X})$ the full subcategory of $\mathbf{Q}(\mathcal{X})$ consisting of objects $F$ such that for any map
$$f:\mathrm{Spec}(A)\rightarrow\mathcal{X}$$
with $A\in\mathbf{A}$, $f^{*}_{\mathbf{Q}}(F)\in\mathbf{Q}_{\ge0}(A)$. By right exactness of $f^{*}_{\mathbf{Q}}$ and \cite{HA} Proposition 1.4.4.11 this defines a $t$-structure on $\mathbf{Q}(\mathcal{X})$, such that we get a presheaf with $t$-structure
$$\mathbf{Q}:\mathbf{PreStk}(\mathbf{A})^{op}\rightarrow\mathbf{Cat}^{\mathbf{L},t}$$

\subsection{Basics of Presheaves With $t$-Structure}

Let us discuss some basic constructions with presheaves with $t$-structures and their properties.

\subsubsection{Truncated Objects}
 For each $n\in\mathbb{Z}$ there is a subpresheaf $\mathbf{Q}_{\ge n}$ of $\mathbf{Q}$. This in fact extends to a presheaf
 $$\mathbf{Q}_{\ge n}:\mathbf{A}^{op}\rightarrow\mathbf{Cat}^{\mathbf{L}}$$
  We also have a sub-presheaf $\mathbf{Q}_{+}$. If $\mathbf{Q}$ is a $\mathbf{Cat}^{\mathbf{L},t,re}$-presheaf this in fact defines a presheaf
$$\mathbf{Q}_{+}:\mathbf{A}^{op}\rightarrow\mathbf{Cat}^{\mathbf{L}}$$

There is also a presheaf  
$$\mathbf{Q}_{\le n}:\mathbf{A}^{op}\rightarrow\mathbf{Cat}^{\mathbf{L}}$$
but it is not a sub-presheaf of $\mathbf{Q}$, as the left adjoint is not the restriction of $f^{*}$, but rather the restriction of $f^{*}$ followed by truncation.

Finally there is a presheaf
$$\mathbf{Q}^{\heart}:\mathbf{A}^{op}\rightarrow\mathbf{Cat}^{\mathbf{L}}$$

\subsubsection{Completeness}

Denote by $\mathbf{Cat}^{\mathbf{L},t,lc}$ (resp. $\mathbf{Cat}^{\mathbf{L},t,rc}$/ $\mathbf{Cat}^{\mathbf{L},t,c}$ ) the full subcategory of $\mathbf{Cat}^{\mathbf{L},t,re}$ consisting of stable categories with $t$-structures whose $t$-structure is left complete (resp. right complete/ complete). Note that here we are implicitly assuming that the right adjoint functors are $t$-exact.

\subsubsection{$\mathbf{P}$-Quasicoherence}

Here we consider a useful subcategory of the heart, which under some restrictions defines a sub-presheaf. 

\begin{defn}\label{defn:Qflatsub}
Let $\mathbf{Q}:\mathbf{A}^{op}\rightarrow\mathbf{Pr}^{\mathbf{L},t}$ be a presheaf with $t$-structure. A map $f:X\rightarrow Y$ in $\mathbf{A}$ is said to be $\mathbf{Q}$-\textit{flat} if the map
$$f^{*}:\mathbf{Q}(Y)\rightarrow\mathbf{Q}(X)$$
is exact for the $t$-structures.
\end{defn}

\begin{defn}
Let $\mathbf{P}$ be a class of maps in $\mathbf{A}$ and $Y$ an object of $\mathbf{A}$. An object $F$ of $\mathbf{Q}_{\le n}(Y)$ is said to be $(n,\mathbf{P})$-\textit{quasicoherent} if for any map $f:X\rightarrow Y$ in $\mathbf{P}$, $f^{*}_{\mathbf{Q}}(F)\in\mathbf{Q}_{\le n}(X)$.
\end{defn}

We let $\mathbf{Q}^{\heart}_{\mathbf{P}}(X)\subset\mathbf{Q}(X)^{\heart}$ denote the full subcategory consisting of those objects in the heart which are $(0,\mathbf{P})$-quasicoherent. 

\begin{rem}
Note that if $\mathbf{P}$ consists of $\mathbf{Q}$-flat maps then any object of $\mathbf{Q}_{\le n}(A)$ is $(n,\mathbf{P})$-quasicoherent. In particular $\mathbf{Q}^{\heart}_{\mathbf{P}}(A)=\mathbf{Q}(A)^{\heart}$. However for many of the topologies we are considering, maps in $\mathbf{P}$ will not in general be flat. 
\end{rem}

\begin{lem}
Let $F\in\mathbf{Q}_{b}(X)$ be such that for each $m$, $H_{m}(F)$ is $(0,\mathbf{P})$-quasicoherent. Then for each $m\le0$, and any map $f:Y\rightarrow X$ in $\mathbf{P}$, the map
$$f^{*}H_{m}(F)\rightarrow H_{m}(f^{*}F)$$
is an isomorphism. If $\mathbf{Q}$ is a
\begin{enumerate}
\item
$\mathbf{Cat}^{\mathbf{L},t,lc}$-presheaf then this is also true for complexes in $\mathbf{Q}_{+}(X)$
\item
$\mathbf{Cat}^{\mathbf{L},t,rc}$-presheaf and filtered colimits are $t$-exact then this is also true for complexes in $\mathbf{Q}_{-}(X)$.
\item
$\mathbf{Cat}^{\mathbf{L},t,c}$-presheaf and filtered colimits are $t$-exact then this is also true for complexes in $\mathbf{Q}(\mathcal{X})$.
\end{enumerate}
\end{lem}

\begin{proof}
 We may assume that $F$ is concentrated in non-negative degrees. Let $l$ denote the length of $F$, i.e. $l$ is the maximum $m$ such that $H_{m}(F)\neq 0$. The proof is by induction on $l$. If $l=0$ then $F$ is equivalent to a $(0,\mathbf{P})$-quasicoherent sheaf concentrated in degree $0$, for which the result is clear. Suppose we have proven the claim for some $l$, and let $F$ have length $l+1$ There is a cofibre sequence
$$F_{\ge l+1}\rightarrow F\rightarrow F_\le{l}$$
where $F_{\le l}$ has length $l$ and $F_{\ge l+1}$ has homology only concentrated in degree $l+1$, with $H_{l+1}(F_{\ge l+1})\cong H_{l+1}(F)$. By shifting, it follows from the $l=0$ case that the claim is true for $F_{\ge l+1}$. By the inductive hypothesis it is true for $F_{\le l}$. Now an easy long exact sequence argument gives the claim for $F$. 

Suppose now that $\mathbf{Q}$ is a $\mathbf{Cat}^{\mathbf{L},t,lc}$-presheaf. Let $f:\mathcal{X}\rightarrow\mathcal{Y}$ be a map. We have
$$\tau_{\le n}f^{*} (M)\cong \tau_{\le n}f^{*} (\tau_{\le n}M)$$
Thus for $l\le m$
$$H_{l}(f^{*}(M))\cong H_{l}(f^{*} (\tau_{\le n}M))\cong f^{*}H_{l}(\tau_{\le n}M)\cong f^{*}H_{l}(M)$$
where we have used the result for bounded complexes.

The right completeness claim is similar. We write $X\cong\limind\tau_{\ge n}X$ with each $\tau_{\ge n}$ bounded below. Then $f^{*}X\cong\limind f^{*}(\tau_{\ge n}X)$. Now by assumption $f^{*}H_{m}(f)\cong H_{m}(f^{*}(\tau_{\ge n}X))$ for $n$ sufficiently large. By exactness of filtered colimits $H_{m}(f^{*}X)\cong f^{*}H_{m}(X)$. 

The final claim follows from the previous ones.


\end{proof}

\begin{notation}
We denote by $\mathbf{Q}_{\mathbf{P}-QC}(A)$ the full subcategory of $\mathbf{Q}(A)$  consisting of objects $F$ such that $H_{m}(F)$ is $(0,\mathbf{P})$-quasicoherent for each $m\in\mathbb{Z}$. For $*\in\{\ge0,\le0,\emptyset,+,-\}$ we write $\mathbf{Q}_{\mathbf{P}-QC}^{*}(A)\defeq\mathbf{Q}^{*}(A)\cap\mathbf{Q}_{\mathbf{P}-QC}$.
\end{notation}

Note that since not all maps are in $\mathbf{P}$, $\mathbf{Q}_{\mathbf{P}-QC}^{*}$ does not define a sub-presheaf of $\mathbf{Q}^{*}$. Instead we denote by $\mathbf{A}_{\mathbf{P}}$ the wide subcategory of $\mathbf{A}$ consisting of maps in $\mathbf{P}$. $\mathbf{Q}$ restricts to a $\mathbf{Cat}^{*}$-presheaf on $\mathbf{A}_{\mathbf{P}}$, and $\mathbf{Q}_{\mathbf{P}-QC}^{b}$ \textit{does} define a sub-presheaf on $\mathbf{A}_{\mathbf{P}}$. In the presence of various completeness assumptions on the $t$-structures, this can also be extended to various categories of unbounded (or bounded in one direction) complexes.

Now let $\mathcal{X}$ be a prestack on $\mathbf{A}$ and denote its restriction to $\mathbf{A}_{\mathbf{P}}$ by $\mathcal{X}_{\mathbf{P}}$. Then 
$$\mathbf{Q}(\mathcal{X}_{\mathbf{P}})=\lim_{\mathbf{A}_{\mathbf{P}}\big\slash\mathcal{X}}\mathbf{Q}(A)$$
There is an obvious map
$$\mathbf{Q}(\mathcal{X})\rightarrow\mathbf{Q}(\mathcal{X}_{\mathbf{P}})$$
We will compare these categories later.

\section{Monoidal Presheaves}\label{sec:monoidalpresheaves}

Let us now discuss presheaves valued in monoidal categories.

\begin{defn}
Denote by 
\begin{enumerate}
\item
$\mathbf{Cat}^{\otimes}$ the $(\infty,1)$-category of symmetric monoidal $(\infty,1)$-categories with morphisms being strong monoidal functors.
\item
$\mathbf{Cat^{L,\otimes}}$ the $(\infty,1)$-category of symmetric monoidal $(\infty,1)$-categories with morphisms being strong monoidal left-adjoint functors.
\item
$\mathbf{Pr^{L,\otimes}}\subset\mathbf{Cat^{L,\otimes}}$ the full subcategory consisting of locally presentably symmetric monoidal categories.
\item
$\mathbf{Pr^{L,\otimes,t,re}}\subset\mathbf{Pr^{L,\otimes}}\times_{\mathbf{Cat}}\mathbf{Cat}^{t,re}$ the full subcategory consisting of those monoidal, stable $(\infty,1)$-categories with $t$-structure $(\mathbf{C},\mathbf{C}_{\ge0},\mathbf{C}_{\le0})$ for which $\mathbf{C}_{\ge0}$ is closed under the monoidal product.
\end{enumerate}
We will write $\mathbf{Cat}^{*,\otimes}$ when we are referring to one of the above categories.
\end{defn}

\begin{defn}
A \textit{monoidal $\mathbf{Cat}^{*}$-presheaf} on an $(\infty,1)$-category $\mathbf{A}$ is a functor
$$\mathbf{Q}:\mathbf{A}^{op}\rightarrow\mathbf{Cat}^{*,\otimes}$$
\end{defn}

Let $\mathbf{Q}$ be a $\mathbf{Pr^{L,\otimes,t,re}}$-presheaf. $\mathbf{Q}$ extends to a functor

$$\mathbf{Q}:\mathbf{PreStk}(\mathbf{A})^{op}\rightarrow\mathbf{Pr^{L,\otimes}}$$
again by Kan extension.

Given a monoidal presheaf $\mathbf{Q}$, we denote by $\mathbf{Perf}^{\mathbf{Q}}$ the subpresheaf such that $\mathbf{Perf}^{\mathbf{Q}}(U)$ consists of finite colimits of the monoidal unit $\mathbb{I}_{\mathbf{Q}(U)}$. 

\subsection{The Projection Formula}[c.f. \cite{stacks-project} 20.54]
For $\mathcal{X}\in\mathbf{PreStk}(\mathbf{A})$ denote by $\otimes_{\mathcal{X}}$ the tensor product functor, $\mathcal{O}_{\mathcal{X}}$ the monoidal unit, and $\underline{\mathbf{Map}}_{\mathcal{O}_{\mathcal{X}}}$ the internal hom. We also write $(-)^{\vee}_{\mathcal{X}}$ for $\underline{\mathbf{Map}}_{\mathcal{O}_{\mathcal{X}}}(-\mathcal{O}_{\mathcal{X}})$ By an easy adjunction argument for any $p:\mathcal{X}\rightarrow\mathcal{Y}$ we have an isomorphism
$$p_{*}\underline{\mathbf{Map}}_{\mathcal{O}_{\mathcal{X}}}(p^{*}(-),-)\cong\underline{\mathbf{Map}}_{\mathcal{O}_{\mathcal{Y}}}((-),p_{*}(-))$$
By adjunction this gives a natural map
$$p^{*}\underline{\mathbf{Map}}_{\mathcal{O}_{\mathcal{Y}}}((-),p_{*}(-))\rightarrow\underline{\mathbf{Map}}_{\mathcal{O}_{\mathcal{X}}}(p^{*}(-),-)$$
and thus a natural map
$$\delta:p^{*}\underline{\mathbf{Map}}_{\mathcal{O}_{\mathcal{Y}}}(p^{*}(-),p^{*}(-))\rightarrow\underline{\mathbf{Map}}_{\mathcal{O}_{\mathcal{X}}}((-),(-))$$
We also have a natural map
$$\eta:p_{*}(-)\otimes_{\mathcal{O}_{\mathcal{Y}}}(-)\rightarrow (p_{*}(-)\otimes_{\mathcal{O}_{\mathcal{X}}}p^{*}(-))$$

It is important to understand for which $\mathcal{F}$ and $\mathcal{G}$ the components of these natural transformations, $\delta_{\mathcal{F},\mathcal{G}}$ and $\eta_{\mathcal{F},\mathcal{G}}$ are equivalences.

\begin{defn}
A triple $(\mathcal{F},\mathcal{G},p)$ with $p:\mathcal{X}\rightarrow\mathcal{Y}$, $\mathcal{F}\in\mathbf{Q}(\mathcal{X})$, and $\mathcal{G}\in\mathbf{Q}(\mathcal{Y})$ is said to \textit{satisfy the projection formula} if 
$$\eta_{\mathcal{F},\mathcal{G}}:f_{*}(\mathcal{F})\otimes_{\mathcal{O}_{\mathcal{Y}}}(\mathcal{G})\rightarrow p_{*}(\mathcal{F})\otimes_{\mathcal{O}_{\mathcal{X}}}f^{*}(\mathcal{G})$$
is an equivalence.
\end{defn}

\begin{lem}
Let $\mathcal{G}\in\mathbf{Q}(\mathcal{X})$ be perfect. Then for any map $p:\mathcal{X}\rightarrow\mathcal{Y}$ and any $\mathcal{F}\in\mathbf{Q}(\mathcal{X})$ the maps
$$\delta_{\mathcal{F},\mathcal{G}}:p^{*}\underline{\mathbf{Map}}_{\mathcal{O}_{\mathcal{Y}}}(p^{*}(\mathcal{F}),p^{*}(\mathcal{G}))\rightarrow\underline{\mathbf{Map}}_{\mathcal{O}_{\mathcal{X}}}(\mathcal{F},\mathcal{G})$$
$$\eta_{\mathcal{F},\mathcal{G}}:p_{*}(\mathcal{F})\otimes_{\mathcal{O}_{\mathcal{Y}}}(\mathcal{G})\rightarrow p_{*}(\mathcal{F})\otimes_{\mathcal{O}_{\mathcal{X}}}p^{*}(\mathcal{G})$$
are equivalences.
\end{lem}

\begin{proof}
Fix $\mathcal{G}$. Since we are in the stable setting the class of objects $\mathcal{F}$ for which $\delta_{\mathcal{F},\mathcal{G}}$ is an equivalence is closed under finite colimits. Thus it suffices to observe that the formulas hold for the strucure sheaf $\mathcal{O}_{\mathcal{X}}$.
\end{proof}

A straightforward computation using the base change formula gives the following. 

\begin{lem}
Let $p:\mathcal{X}\rightarrow\mathcal{Y}$ be a universal $\mathbf{Q}$-base change morphism. Then $(\mathcal{F},\mathcal{G},p)$ satisfies the projection formula if and only if for any $f:U\rightarrow\mathcal{Y}$ with $U\in\mathbf{A}^{op}$ the triple $((f')^{*}\mathcal{F},f^{*}\mathcal{G},p')$ satisfies the projection formula, where $f'$ and $p'$ are defined by the following pullback diagram
\begin{displaymath}
\xymatrix{
U\times_{\mathcal{Y}}\mathcal{X}\ar[r]^{f'}\ar[d]^{p'} & \mathcal{X}\ar[d]^{p}\\
U\ar[r]^{f} & \mathcal{Y}
}
\end{displaymath}
\end{lem}

In particular if $p$ is a base-change morphism then the projection formula can be verified locally.


\subsubsection{The Plus Construction}\label{subsubsec:plus}
(following 
\cite{porta2017representability})
Let $p:\mathcal{X}\rightarrow\mathcal{Y}$ be a map of prestacks. Define 
$$p_{+}\defeq(-)_{\mathcal{Y}}^{\vee}\circ p_{*}\circ(-)_{\mathcal{X}}^{\vee}$$

Consider the counit map
$$p^{*}p_{*}\mathcal{F}^{\vee}\rightarrow\mathcal{F}^{\vee}$$
Taking duals gives a map
$$\mathcal{F}^{\vee\vee}\rightarrow (p^{*}p_{*}\mathcal{F}^{\vee})^{\vee}$$
If $\mathcal{F}$ is perfect then $\mathcal{F}^{\vee\vee}\cong\mathcal{F}$, and $(p^{*}p_{*}\mathcal{F}^{\vee})^{\vee}\cong p^{*}p_{+}\mathcal{F}$. Thus we have a natural transformation
$$\eta_{\mathcal{F}}:\mathcal{F}\rightarrow p^{*}p_{+}\mathcal{F}$$

\begin{lem}
Let $p:\mathcal{X}\rightarrow\mathcal{Y}$ be a map such that $p_{*}\mathbf{Perf}^{\mathbf{Q}}(\mathcal{X})\subseteq\mathbf{Perf}^{\mathbf{Q}}(\mathcal{Y})$. Then for $\mathcal{F}\in\mathbf{Perf}^{\mathbf{Q}}(\mathcal{X})$, $\mathcal{G}\in\mathbf{Perf}^{\mathbf{Q}}(\mathcal{Y})$ there is an equivalence, natural in $\mathcal{F}$ and $\mathcal{G}$,
$$\underline{\mathbf{Map}}_{\mathcal{O}_{\mathcal{Y}}}(p_{+}\mathcal{F},\mathcal{G})\cong p_{*}\underline{\mathbf{Map}}_{\mathcal{O}_{\mathcal{X}}}(\mathcal{F},f^{*}\mathcal{G})$$
\end{lem}

\begin{proof}
We have 
\begin{align*}
\underline{\mathbf{Map}}_{\mathcal{O}_{\mathcal{Y}}}(p_{+}\mathcal{F},\mathcal{G})&\cong p_{*}(\mathcal{F}^{\vee})\otimes\mathcal{G}\\
&\cong p_{*}(\mathcal{F}^{\vee}\otimes p^{*}\mathcal{G})\\
&\cong p_{*}\underline{\mathbf{Map}}_{\mathcal{O}_{\mathcal{X}}}(\mathcal{F},f^{*}\mathcal{G})
\end{align*}
\end{proof}

\begin{defn}\label{defn:cohprop}
Let $f:\mathcal{X}\rightarrow \mathcal{Y}$ be a map of presstacks. $f$ is said to be \textit{strongly} $\mathbf{Q}$-\textit{cohomologically proper} if for any $g:U\rightarrow\mathcal{Y}$ with $U$ affine, the morphism
$$(f')_{*}:\mathbf{Q}(U\times_{\mathcal{Y}}\mathcal{X})\rightarrow\mathbf{Q}(U)$$
sends objects in $\mathbf{Perf}^{\mathbf{Q}}(U\times_{\mathcal{Y}}\mathcal{X})$ to objects in $\mathbf{Perf}^{\mathbf{Q}}(U)$. 
\end{defn}

\begin{defn}\label{defn:perfbase}
  Let  $f:\mathcal{X}\rightarrow \mathcal{Y}$ be strongly $\mathbf{Q}$-cohomologically proper. $f$ is said to be \textit{perfect base-change} if for any map $U\rightarrow\mathcal{Y}$ with $U$ affine and any fibre-product diagram
  \begin{displaymath}
\xymatrix{
\mathcal{X}\times_{\mathcal{Y}}U\ar[d]^{f'}\ar[r]^{g'} & \mathcal{X}\ar[d]^{f}\\
U\ar[r]^{g} & \mathcal{Y}
}
\end{displaymath}
the natural map
$$g^{*}f_{*}\mathcal{F}\rightarrow f'_{*}(g')^{*}\mathcal{F}$$
is an equivalence whenever $\mathcal{F}\in\mathbf{Perf}^{\mathbf{Q}}(\mathcal{X})$
\end{defn}

\begin{lem}\label{lem:Perfbasechange}
Let
\begin{displaymath}
\xymatrix{
\mathcal{X}\ar[r]^{f'}\ar[d]^{p'} & \mathcal{X}\ar[d]^{p}\\
\mathcal{Y}'\ar[r]^{f} & \mathcal{Y}
}
\end{displaymath}
be a pullback diagram in which the base-change formula is satisfied for all objects in $\mathbf{Perf}^{\mathbf{Q}}(\mathcal{X})$, and such that both $p'_{*}$ and $p_{*}$ preserve perfect objects. Let $\mathcal{F}\in\mathbf{Perf}^{\mathbf{Q}}(\mathcal{X})$. Then there is a natural equivalence
$$(p')_{+}(f')^{*}\mathcal{F}\rightarrow f^{*}p_{+}\mathcal{F}$$
\end{lem}

\begin{proof}
This follows simply from dualising the usual base change formula.
\end{proof}

\subsection{Pre-Stacks of Modules}

Let us restrict to the following setup, which is the crucial one for the purposes of derived geometry. Fix a presentably symmetric monoidal $(\infty,1)$-category $\mathbf{C}$, a monad $\mathbf{D}:\mathbf{C}\rightarrow\mathbf{C}$, and a map of monads $\theta:\mathbf{D}\rightarrow\mathbf{Comm}(-)$ such that the induced functor
$$\theta:\mathbf{Alg_{D}}(\mathbf{C})\rightarrow\mathbf{Comm}(\mathbf{C})$$
commutes with limits and colimits. We consider the presheaf $\mathbf{QCoh}$ on $\mathbf{Alg_{D}}^{op}$.
This determines a presheaf
$$\mathbf{QCoh}:\mathbf{Alg_{D}(C)}^{op}\rightarrow\mathbf{Pr}^{\mathbf{L},\aleph_{0},\otimes}$$

\subsubsection{Coherent Sheaves}\label{subsbsec:coherent}

Let $\mathbf{Coh}_{+}(\mathrm{Spec}(A))\subset\mathbf{QCoh}(\mathrm{Spec}(A))$ denote the subpresheaf consisting of objects $F$ such that $\pi_{m}(F)$ is finitely presented over $\pi_{0}(A)$ for all $m$. Let $A\rightarrow B$ be a map of algebras such that $\pi_{0}(A)$ and $\pi_{0}(B)$ are coherent. By Corollary \ref{cor:sthngleftcoh}, $B\otimes^{\mathbb{L}}_{A}F$ is in $\mathbf{Coh}_{+}(\mathrm{Spec}(A))$. Suppose for all $\mathrm{Spec}(A)\in\mathbf{A}$, $\pi_{0}(A)$ is coherent. Then we get subpresheaves
$$\mathbf{Coh}_{\ge n}|_{\mathbf{A}}\subset\mathbf{QCoh}|_{\mathbf{A}}$$
$$\mathbf{Coh}_{+}|_{\mathbf{A}}\subset\mathbf{QCoh}|_{\mathbf{A}}$$
\subsubsection{Compact, Perfect, Nuclear, and $\kappa$-Filtered Objects}\label{subsbsec:kappfilt}
Since $\mathbf{QCoh}$ is a $\mathbf{Pr}^{\mathbf{L},\aleph_{0},\otimes}$-presheaf we have a well-defined sub-presheaf
$$\mathbf{QCoh}^{\aleph_{0}-cpct}\subset\mathbf{QCoh}$$

In presheaves of module categories we may also define the sub-presheaves $\mathbf{Perf}$ and $\mathbf{Nuc}$ of perfect and nuclear objects respectively. In general we will have an inclusion
$$\mathbf{Perf}\subseteq\mathbf{QCoh}^{\aleph_{0}-cpct}$$
which may be strict, and in fact will be so in our main case of interest.

We also write $\mathbf{Perf}_{\ge n}$ for the sub-presheaf of perfect objects which are homologically concentrated in degrees $\ge n$.


\subsubsection{Modules in Spectral Algebraic Contexts}
Fix a spectral algebraic context
$$\underline{\mathbf{C}}\defeq(\mathbf{C},\mathbf{C}_{\ge0},\mathbf{C}_{\le0},\mathbf{C}^{0},\mathbf{D},\theta)$$

We write
$$\mathbf{Aff}^{cn}_{\underline{\mathbf{C}}}\defeq(\mathbf{Alg}_{\mathbf{D}}^{cn})^{op}$$

When $\underline{\mathbf{C}}$ is a derived algebraic context of the form $\underline{\mathbf{Ch}}(\mathpzc{E})$ for $\mathpzc{E}$ a monoidal elementary exact category we write 
$$\mathbf{Aff}^{cn}_{\mathpzc{E}}\defeq\mathbf{Aff}_{\underline{\mathbf{Ch}}(\mathpzc{E})}$$

 
 If $f:A\rightarrow B$ is a map in $\mathbf{Alg_{D}}^{cn}(\mathbf{C})$ then the left adjoint functor 
 $$\mathbf{QCoh}(\mathrm{Spec}(A))\cong{}_{A}\mathbf{Mod}\rightarrow\mathbf{QCoh}(\mathrm{Spec}(B))\cong\mathbf{Mod}(B)$$
 is given by $B\otimes^{\mathbb{L}}_{A}(-)$. The $t$-structure on ${}_{A}\mathbf{Mod}$ means that this in fact defines a functor
 $$\mathbf{QCoh}:(\mathbf{Aff}^{cn}(\mathbf{C}))^{op}\rightarrow\mathbf{Pr}^{\mathbf{L},t,re,c}$$
 Note that 
$$\mathbf{QCoh}:(\mathbf{Aff}^{cn}(\mathbf{C}))^{op}\rightarrow\mathbf{Pr}^{\mathbf{L},t}$$
is of affine type, but 
$$\mathbf{QCoh}:\mathbf{PreStk}(\mathbf{C})^{op}\rightarrow\mathbf{Pr}^{\mathbf{L},t}$$
will not be in general.

\subsubsection{Vector Bundles and Locally Free Sheaves}
We conclude this chapter with a discussion of vector bundles. Fix a derived algebraic context
$$(\mathbf{C},\mathbf{C}_{\ge0},\mathbf{C}_{\le0},\mathbf{C}^{0})$$
and full subcategories 
$\mathbf{A}\subset\mathbf{DAlg}^{cn}(\mathbf{C})^{op}$, $\mathbf{V}\subset\mathbf{C_{\ge0}}$.

\begin{defn}
Let $\mathrm{Spec}(A)\in\mathbf{A}$, and $M\in{}_{A}\mathbf{Mod}^{cn}$. $M$ is said to be \textit{locally free of type} $\mathbf{V}$ if there exists a cover $\{\mathrm{Spec}(A_{i})\rightarrow\mathrm{Spec}(A)\}$ and $M_{i}\in\mathbf{V}$ for each $i\in\mathcal{I}$ such that $A_{i}\otimes^{\mathbb{L}}_{A}F\cong A_{i}\otimes^{\mathbb{L}}M_{i}$. 
\end{defn}

Let $\mathbf{QCoh}^{lf,\mathbf{V}}(\mathrm{Spec}(A))\subset\mathbf{QCoh}(\mathrm{Spec}(A))$ denote the full subcategory consisting of locally free sheaves of type $\mathbf{V}$. Note that by construction
$$\mathbf{QCoh}^{lf,\mathbf{V}}\subset\mathbf{QCoh}$$
is a local sub-presheaf.

\begin{defn}
Fix a stack $\mathcal{Y}$. A stack $\mathcal{X}\in\mathbf{Stk}(\mathbf{A}\big\slash_{\mathcal{Y}})$ is said to be a \textit{vector bundle of type} $\mathbf{V}$ \textit{over} $\mathcal{Y}$ if for any map $U\rightarrow\mathcal{Y}$ there exists a cover $\coprod U_{i}\rightarrow U$ such that $U_{i}\times_{\mathcal{Y}}\mathcal{X}\cong U_{i}\times\mathrm{Spec}(\mathbf{LSym}(M))$ for some $M\in\mathbf{V}$. 
\end{defn}

We denote the category of vector bundles of type $\mathbf{V}$ over $\mathcal{Y}$ by $\mathbf{Vect}_{\mathbf{V}}(\mathcal{Y})$. When $\mathbf{V}=\mathbf{C}_{\ge0}$ we just write this as $\mathbf{Vect}(\mathcal{Y})$.

Let $\mathbf{DAlg}^{cn}(\mathbf{QCoh}(-))$ denote the presheaf of categories which assigns to $\mathrm{Spec}(A)$ the $(\infty,1)$-category ${}_{{}_{A}\big\backslash}\mathbf{DAlg}^{cn}(\mathbf{C})$. There is a forgetful natural transformation of presheaves
$$\mathbf{DAlg}^{cn}(\mathbf{QCoh}(-))\rightarrow\mathbf{QCoh}(-)$$
Its left adjoint sends $F\in\mathbf{QCoh}(\mathrm{Spec}(A))$ to $\mathbf{LSym}_{A}(F)$.

Let $\mathcal{Y}$ be a prestack and let $\mathcal{B}\in\mathbf{DAlg}^{cn}(\mathbf{QCoh}(\mathcal{Y}))$. Define a pre-stack on $\mathbf{A}_{\big\slash\mathcal{Y}}$, $\mathrm{Spec}_{\mathcal{Y}}(\mathcal{B})\rightarrow\mathcal{Y}$, by defining 

$$\mathrm{Spec}_{\mathcal{Y}}(\mathcal{B})(\mathrm{Spec}(A)\rightarrow\mathbf{Y})\defeq\mathbf{Map}_{{}_{A\big\backslash}\mathbf{DAlg}^{cn}(\mathbf{C})}(\mathcal{B}(\mathrm{Spec}(A)),A)$$

Consider the functor 
$$\mathrm{Spec}\circ\mathbf{LSym}_{\mathcal{Y}}\circ(-)^{\vee}:\mathbf{QCoh}(\mathcal{Y})^{op}\rightarrow\mathbf{Stk}(\mathbf{A}_{\big\slash \mathcal{Y}})$$

If $\mathbf{Q}$ consists of free modules of finite rank, then this functor sends a locally free sheaf to a vector bundle. Moreover it is fully faithful.

\chapter{Descent}\label{D}

In this chapter we discuss various topologies which satisfy descent. The proofs for different `types' topologies will be fundamentally different.

\section{Some Generalities}

 Let $\mathbf{A}$ be equipped with a Grothendieck pre-topology $\tau$, and fix a presheaf of categories
$$\mathbf{Q}:\mathbf{A}^{op}\rightarrow\mathbf{Pr^{L}}$$
and a subpresheaf $\mathbf{N}\subset\mathbf{Q}^{o}$.

\subsection{Descent and Sub-presheaves}
Let $K_{\bullet}$ be a simplicial pre-stack and $f:K_{\bullet}\rightarrow\mathcal{X}$ a morphism, where $\mathcal{X}$ is a pre-stack regarded as a constant simplicial pre-stack. 

\begin{defn}
\begin{enumerate}
\item
$K_{\bullet}\rightarrow\mathcal{X}$ is said to be a weak $\mathbf{N}$-\textit{cover} if the map
$$\mathbf{N}(|K_{\bullet}|)\cong\lim_{n}\mathbf{N}(K_{n})\rightarrow\mathbf{Q}(\mathcal{X})$$ 
induces an equivalence with the full subcategory $\mathbf{Q}(\mathcal{X})^{\mathbf{N}_{K}}$ of $\mathbf{Q}(\mathcal{X})$ consisting of those objects $F$ such that for each $f_{n}:K_{n}\rightarrow\mathcal{X}$, $f_{n}^{*}F\in\mathbf{N}(K_{n})$, with inverse given by the restriction of the map 
$$\mathbf{Q}(\mathcal{X})\rightarrow\lim_{n}\mathbf{Q}(K_{n})$$
to $\mathbf{Q}(\mathcal{X})^{\mathbf{N}_{K}}$
\item
$K_{\bullet}\rightarrow\mathcal{X}$ is said to be a $\mathbf{N}$-\textit{cover} if it is a $\mathbf{N}$-weak cover and the map $\mathbf{N}(\mathcal{X})\rightarrow \mathbf{Q}(\mathcal{X})^{\mathbf{N}_{K}}$ is an equivalence.
\end{enumerate}
\end{defn}

\begin{defn}
\begin{enumerate}
\item
$\mathbf{N}$ is said to \textit{satisfy (weak) hyperdescent} if for any pseudo-representable hypercover $K_{\bullet}\rightarrow U$, with $U$ representable is a (weak) $\mathbf{N}$-cover. 
\item
$\mathbf{N}$ is said to \textit{satisfiy (weak) \v{C}ech descent} if any \v{C}ech hypercover $\check{C}(\{U_{i}\rightarrow U\})\rightarrow U$ is a (weak) $\mathbf{N}$-cover.
\end{enumerate}
\end{defn}

In this chapter we will deduce several descent results for various kinds of pre-topologies. Before this though we will understand when descent for a pre-sheaf implies decsent for certain sub-pre-sheaves.


\subsection{Descent for Subcategories}

\begin{defn}
Let $\mathbf{Q}:\mathbf{A}^{op}\rightarrow\mathbf{Pr^{L}}$ be a $\mathbf{Cat}^{*}$-presheaf, and let $\mathbf{N}\subset\mathbf{Q}^{o}$ be a sub-presheaf. $\mathbf{N}$ is said to be \textit{local} if whenever $\{f_{i}:U_{i}\rightarrow U\}_{i\in\mathcal{I}}$ is a cover in $\tau$, and $m\in\mathbf{Q}(U)$ is such that $(f_{i})^{*}_{\mathbf{Q}}(m)$ is in $\mathbf{N}(U_{i})$ for all $i$, then $m\in\mathbf{Q}(U)$. 
\end{defn}

\begin{lem}
Suppose that $\mathbf{Q}$ satisfies hyperdescent (resp. \v{C}ech descent) and that $\mathbf{N}$ is local. Then $\mathbf{N}$ satisfies hyperdescent (resp. \v{C}ech descent).
\end{lem}

\begin{proof}
It suffices to prove that if $K_{\bullet}\rightarrow V$ is a hypercover, and $m_{n}\in\mathbf{N}(K_{n})$ for each $n$, then $R_{K}(m_{\bullet})\in\mathbf{N}(V)$, where $R_{K}:\lim_{n}\mathbf{N}(K_{n})$ is the limit functor. Write $K_{0}=\coprod U_{i}$. Then $\{f_{i}:U_{i}\rightarrow V\}$ is a cover. By assumption we have that $(f^{*}_{i})_{\mathbf{Q}}(R_{K}((m_{\bullet})))\in\mathbf{N}(U_{i})$ for each $i$. Since $\mathbf{N}$ is local this means $R_{K}((m_{\bullet}))\in\mathbf{N}(V)$. This completes the proof. 
\end{proof}

\begin{lem}
Let $\tau$ be a pre-topology on $\mathbf{A}^{op}$ such that 
$$\mathbf{Q}:\mathbf{A}^{op}\rightarrow\mathbf{Pr^{L}}$$
satisfies descent for $\tau$. 
Suppose that for each $f:U\rightarrow V$ the functor $f_{*,\mathbf{Q}}$ commutes with $\kappa$-filtered colimits.
Let $\tau$ be a topology such that whenever $\{f_{i}:U_{i}\rightarrow U\}_{i\in\mathcal{I}}$ is a cover then either
\begin{enumerate}
    \item 
    $\kappa\ge\aleph_{0}$ and $\mathcal{I}$ $\kappa$-small. 
    \item 
    $\kappa<\aleph_{0}$ and the $U_{i}\rightarrow U$ are monomorphisms.
\end{enumerate}

Let $m\in\mathbf{Q}(U)$ be such that $(f_{i})^{*}_{\mathbf{Q}}(m)$ is in $\mathbf{Q}(U_{i})^{\kappa-cpct}$ for all $i\in\mathcal{I}$. Then $m\in\mathbf{Q}(U)^{\kappa-cpct}$. In particular the functor subpresheaf
$\mathbf{Q}^{\kappa-cpct}$
satisfies descent if all covers have a $\kappa$-small refinement.
\end{lem}

\begin{proof}
Let $m$ be such that $f_{i}^{*}m$ is $\kappa$-compact for each $i$. Then $f^{*}(i_{1},\ldots,i_{n})(m)$ is $\kappa$-compact for each $(i_{1},\ldots,i_{n})\in\mathcal{I}^{n}$ since the push-forward functors commute with $\kappa$-filtered colimits. Let $F:\mathcal{J}\rightarrow\mathbf{Q}(U)$ be a $\kappa$-filtered diagram. We then have
\begin{align*}
\mathbf{Map}(m,\colim_{j\in\mathcal{J}})&\cong\mathbf{Map}(m,(f_{\underline{i}})_{*}f_{\underline{i}}^{*}\colim_{j\in\mathcal{J}}F(j))\\
&\cong\lim_{\underline{i}\in\mathcal{I}^{n}}\mathbf{Map}(m,f_{\underline{i}}^{*}\colim_{j\in\mathcal{J}}F(j))\\
&\cong\lim_{\underline{i}\in\mathcal{I}^{n}}\colim_{j\in\mathcal{J}}\mathbf{Map}(f_{\underline{i}}^{*}m,f_{\underline{i}}^{*}F(j))\\
&\cong\colim_{j\in\mathcal{J}}\lim_{\underline{i}\in\mathcal{I}^{n}}\mathbf{Map}(f_{\underline{i}}^{*}m,f_{\underline{i}}^{*}F(j))\\
&\cong\colim_{j\in\mathcal{J}}\mathbf{Map}(f_{\underline{i}}^{*}m,\lim_{\underline{i}\in\mathcal{I}^{n}}f_{\underline{i}}^{*}F(j))\\
&\cong\colim_{j\in\mathcal{J}}\mathbf{Map}(m,(\lim_{\underline{i}\in\mathcal{I}^{n}}(f_{\underline{i}})_{*}f_{\underline{i}}^{*}F(j))\\
&\cong\colim_{j\in\mathcal{J}}\mathbf{Map}(m,F(j))
\end{align*}
where we have used that $\kappa$-filtered colimits commute with $\kappa$-small limits
\end{proof}

\subsection{Descent for Sheaves on Stacks}
Let us now compare sheaves on prestacks with sheaves on their stackifications. Let $i:\mathcal{X}\rightarrow\tilde{X}$ be a stackification. Consider the induced adjunction
$$\adj{i^{*}_{\mathbf{Q}}}{\mathbf{Q}(\tilde{\mathcal{X}})}{\mathbf{Q}(\mathcal{X})}{i_{*,\mathbf{Q}}}$$
The following is tautological.

\begin{prop}
\begin{enumerate}
\item
If $\mathbf{Q}$ satisfies hyperdescent then the adjunction
$$\adj{i^{*}_{\mathbf{Q}}}{\mathbf{Q}(\tilde{\mathcal{X}})}{\mathbf{Q}(\mathcal{X})}{i_{*,\mathbf{Q}}}$$
is an equivalence.
\item
If $\mathbf{Q}$ satisfies \v{C}ech descent then the adjunction
$$\adj{i^{*}_{\mathbf{Q}}}{\mathbf{Q}(\check{\mathcal{X}})}{\mathbf{Q}(\mathcal{X})}{i_{*,\mathbf{Q}}}$$
is an equivalence.
\end{enumerate}
\end{prop}

In particular we get the following.

\begin{cor}
\begin{enumerate}
\item
Let $\mathcal{X}$ be a stack and $K_{\bullet}\rightarrow\mathcal{X}$ a hypercover. Suppose that $\mathbf{Q}$ satisfies hyperdescent. Then the map
$$\underset{n}\lim\mathbf{Q}(K_{n})\rightarrow\mathbf{Q}(\mathcal{X})$$
is an equivalence. 
\item
Let $\mathcal{X}$ be a \v{C}ech stack and $K_{\bullet}\rightarrow\mathcal{X}$ a \v{C}ech cover. Suppose that $\mathbf{Q}$ satisfies \v{C}ech descent. Then the map
$$\underset{n}\lim\mathbf{Q}(K_{n})\rightarrow\mathbf{Q}(\mathcal{X})$$
is an equivalence. 
\end{enumerate}
\end{cor}


Several of the topologies in which we are interested, often those  allowing infinite covers, do not appear to have descent unless we restrict to sub-presheaves. 
However, we can still relate sheaves on $\mathcal{X}$ to sheaves on its stackification.
Let $\mathbf{N}\subset\mathbf{Q}$ be a subpresheaf. 
\begin{lem}
If $\mathbf{N}$ satisfies hyper-descent (resp. \v{C}ech descent) then the functor
$$(i_{*})_{\mathbf{Q}}:\mathbf{N}(\mathcal{X})\rightarrow\mathbf{Q}(\tilde{\mathcal{X}})$$
(resp. 
$$(i_{*})_{\mathbf{Q}}:\mathbf{N}(\mathcal{X})\rightarrow\mathbf{Q}(\check{\mathcal{X}})$$)
is fully faithful.
\end{lem}
\begin{proof}
Write $\mathcal{X}^{(1)}\defeq \underset{(U\leftarrow K_{\bullet}\rightarrow\mathcal{X})} \colim U$ where the colimit is taken over the category $\mathrm{Span}^{aff,hyp,\mathcal{X}}$ of all spans $U\leftarrow K_{\bullet}\rightarrow\mathcal{X}$ where $K_{\bullet}\rightarrow U$ is a hypercover. Inductively define $\mathcal{X}^{(n+1)}\defeq \left(\mathcal{X}^{(n)}\right)^{(1)}$. Then $\tilde{X}\cong\limind_{n}\tilde{\mathcal{X}}^{(n)}$.
It suffices to show that the map
$$\mathbf{N}(\mathcal{X})\rightarrow\mathbf{N}(\tilde{\mathcal{X}}^{(1)})$$
is an equivalence.
Consider the functor
$$\mathbf{N}:\mathrm{Span}^{aff,hyp,\mathcal{X}}\rightarrow\mathbf{Cat}$$
sending $U\leftarrow K_{\bullet}\rightarrow\mathcal{X}$ to $\mathbf{N}(U)$. Now we have
$$\underset{(\mathrm{Span}^{aff,hyp,\mathcal{X}})^{op}}\lim\mathbf{N}(U\leftarrow K_{\bullet}\rightarrow\mathcal{X})\cong\underset{(\mathrm{Span}^{aff,hyp,\mathcal{X}})^{op}}\lim \lim_{n}\mathbf{N}(K_{n})\cong\mathbf{N}(\mathcal{X})$$
\end{proof}

\subsubsection{$\mathbf{P}$-Quasicoherent Descent}

Let $(\mathbf{A},\tau,\mathbf{P})$ be an $(\infty,1)$-geometry triple, and $\mathbf{Q}$ a $\mathbf{Cat}^{*}$-presheaf on $\mathbf{A}$. Essentially by definition we have the following.

\begin{prop}
$\mathbf{Q}_{\mathbf{P}-QC}\subset\mathbf{Q}|_{\mathbf{A}_{\mathbf{P}}}$ is a flat sub-presheaf (Definition \ref{defn:Qflatsub}).
\end{prop}

\begin{lem}
Let $\mathcal{X}$ be a pre-stack. Suppose there exists a pseudo-representable hypercover $K\rightarrow\mathcal{X}$ with each $K$ a disjoint union of objects of $\mathbf{A}$, and each face map $K_{n+1}\rightarrow K_{n}$ being in $\mathbf{P}^{\coprod}$. If $\mathbf{Q}$ has hyperdescent for $\tau$ then the functor
$$\mathbf{Q}_{\mathbf{P}-QC}(\mathcal{X})\rightarrow\mathbf{Q}(\mathcal{X})$$
is fully faithful. 
\end{lem}

\begin{proof}
We have $|\mathbf{Q}(K)|\cong \mathbf{Q}(\mathcal{X})$. Clearly $\mathbf{Q}|_{\mathbf{P}}$ also satisfies hyperdescent for $(\mathbf{A}_{\mathbf{P}},\tau)$ and $K\rightarrow\mathcal{X}|_{\mathbf{P}}$ is also a hypercover. Thus
$$\mathbf{Q}|_{\mathbf{P}}(\mathcal{X}|_{\mathbf{P}})\cong|\mathbf{Q}|_{\mathbf{P}}(K)|\cong|\mathbf{Q}(K)|\cong\mathbf{Q}(\mathcal{X})$$ 
Since $\mathbf{Q}_{\mathbf{P}-QC}\subset\mathbf{Q}|_{\mathbf{P}}$, we get the result.
\end{proof}

\subsection{The $\kappa$-Embedding pre-topology}

One very general pre-topology which satisfies descent, and which will be of use to us, is the $\kappa$-embedding pre-topology. Let $(\mathbf{A},\tau)$ be a site, $\mathbf{Q}$ a $\mathbf{Pr}^{\mathbf{L}}$-valued presheaf, and $\mathbf{N}\subseteq\mathbf{Q}|_{\mathbf{P}}$ a subpresheaf. 

\begin{defn}
    A map $\mathcal{X}\rightarrow\mathcal{Y}$ in an $(\infty,1)$-category $\mathbf{C}$ is said to be a \textit{homotopy monomorphism} if the natural map
    $$\mathcal{X}\rightarrow\mathcal{X}\times_{\mathcal{Y}}\mathcal{X}$$
    is an equivalence.
\end{defn}

The class of all homotopy monomorphisms contains all isomorphisms, and is stable by base change and composition. Fix a class of maps $\mathbf{P}$ in $\mathbf{A}$ which is closed under composition and pullback, and let $\mathbf{P}^{\kappa}_{epi}$ denote the class of maps $f:\mathrm{Spec}(B)\rightarrow \mathrm{Spec}(A)$ such that

\begin{enumerate}
    \item 
    $f$ is in $\mathbf{P}$
    \item 
    $f$ is a homotopy monomorphism
    \item 
    $f$ is a universal base-change morphism
    \item 
    $f$ is a universal $\mathbf{Q}$-embedding.
    \item 
    If $\mathrm{Spec}(C)\rightarrow\mathrm{Spec}(A)$ is a map, and $f'$ is the base-change of $f$ along this map, then
    $(f')^{*}$ is formally $\kappa$ filtered - i.e. it commutes with $\kappa$-small limits.
\end{enumerate}

Clearly the class $\mathbf{P}^{\kappa}_{epi}$ contains all equivalences, and is stable by pullback and composition.

\begin{defn}
Suppose $\mathbf{Q}$ is of affine type. The $(\kappa,\mathbf{P},\mathbf{N})$-embedding pre-topology on $\mathbf{A}$ is the pre-topology whose covers are of the form $\{f_{i}:\mathrm{Spec}(B_{i})\rightarrow\mathrm{Spec}(A)\}_{i\in\mathcal{I}}$ where
\begin{enumerate}
    \item 
    Each $\mathrm{Spec}(B_{i})\rightarrow\mathrm{Spec}(A)$ is in $\mathbf{P}^{\kappa}_{epi}$.
    \item 
    There exists a $\kappa$-small subset $\mathcal{J}$ of $\mathcal{I}$ such that whenever $\alpha:M\rightarrow N$ is a map in $\mathbf{N}(\mathrm{Spec}(A))$ is such that $f_{i}^{*}(\alpha)$ is an equivalence for all $j\in\mathcal{J}$, then $\alpha$ is an equivalence. 
\end{enumerate}
This pre-topology is denoted $\tau^{\mathbf{P},\mathbf{N},\kappa,epi}$
\end{defn}

\begin{prop}
  If $\mathbf{Q}$ is of affine type then $\tau^{\mathbf{P},\mathbf{N},\kappa,epi}$ does define a pre-topology on $\mathbf{A}$.
\end{prop}

\begin{proof}
The only point that is not clear is that if $\{f_{i}:\mathrm{Spec}(B_{i})\rightarrow\mathrm{Spec}(A)\}_{i\in\mathcal{I}}$ is a cover, and $g:\mathrm{Spec}(C)\rightarrow\mathrm{Spec}(A)$ is any map, then $\{f'_{i}:\mathrm{Spec}(C)\times_{\mathrm{Spec}(A)}\mathrm{Spec}(U_{i})\rightarrow\mathrm{Spec}(A)\}_{i\in\mathcal{I}}$ satisfies the second condition. So let $\alpha:M\rightarrow N$ be a map in $\mathbf{N}(\mathrm{Spec}(C))$ such that $(f'_{i})^{*}(\alpha)$ is an equivalence for all $i$. Since all maps are push-forward conservative it suffices to prove that $g_{*}(\alpha)$ is an equivalence. But we have 
$$f_{i}^{*}g_{*}(\alpha)\cong (g')_{*}(f'_{i})^{*}(\alpha)$$
Hence $f_{i}^{*}g_{*}(\alpha)$ is an equivalence for all $i$, so $g_{*}(\alpha)$ is an equivalence, as required.
\end{proof}

\begin{prop}
    Let $\mathbf{Q}$ be of affine type, and let $\kappa$ be a cardinal. Let $K_{\bullet}\rightarrow\mathrm{Spec}(A)$ be the nerve of a \v{C}ech cover in $\tau^{\mathbf{P},\mathbf{N},\kappa,epi}$. Then the map
    $$\mathbf{Q}(|K_{\bullet}|)\rightarrow\mathbf{Q}(\mathrm{Spec}(A))$$
    is fully faithful. Moreover $\mathbf{N}$ satisfies weak \v{C}ech descent for $\tau^{\mathbf{P},\mathbf{N},\kappa,epi}$.
 In particular, if $\mathbf{N}=\mathbf{Q}$, then $\mathbf{Q}$ satisfies \v{C}ech descent.
\end{prop}

\begin{proof}
The claim about conservativity (i.e. weak \v{C}ech descent) for
   $$\mathbf{N}(\mathrm{Spec}(A)))\rightarrow\mathbf{Q}(K_{\bullet})$$
   is clear, so we just need to prove that $\mathbf{Q}(|K_{\bullet}|)\rightarrow\mathbf{Q}(\mathrm{Spec}(A))$ is fully faithful. Let $F_{\bullet}\in\mathbf{Q}(|K_{\bullet}|)$ The functor $\mathbf{Q}|K_{\bullet}|\rightarrow\mathbf{Q}(\mathrm{Spec}(A))$ sends $E_{\bullet}\in\mathbf{Q}(|K_{\bullet}|)$ to $\mathbf{lim}_{n}\prod_{i_{n}\in\mathcal{I}_{n}}(f_{i_{n}})_{*}E_{i_{n}}$ where $K_{n}=\coprod_{i_{n}\in\mathcal{I}_{n}}\mathrm{Spec}(B_{i_{n}})$ and $E_{i_{n}}\in\mathrm{Spec}(B_{i_{n}})$. Let $f_{i_{n}}:\mathrm{Spec}(B_{i_{n}})\rightarrow\mathrm{Spec}(A)$ be the natural map. We need to show that for each $i_{0}\in\mathcal{I}_{0}$, 
   $$f_{i_{0}}^{*}(\mathbf{lim}_{n}\prod_{i_{n}\in\mathcal{I}_{n}}(f_{i_{n}})_{*}E_{i_{n}})\cong E_{i_{0}}$$
   Note that regardless of what $\kappa$ is $\mathbf{lim}_{n}\prod_{i_{n}\in\mathcal{I}_{n}}(f_{i_{n}})_{*}E_{i_{n}}$ is $\kappa$-small. Indeed if $\kappa\ge\aleph_{0}$ then this is clear. But even if $\kappa$ is finite then since maps in $\mathbf{P}$ are monomorphisms we may assume our \v{C}ech cover is finite, and then since all maps are monomorphisms, $\mathbf{lim}_{n}\prod_{i_{n}\in\mathcal{I}_{n}}(f_{i_{n}})_{*}E_{i_{n}}$ is a finite limit. Thus we have 
   
   $$f_{i_{0}}^{*}(\mathbf{lim}_{n}\prod_{i_{n}\in\mathcal{I}_{n}}(f_{i_{n}})_{*}E_{i_{n}})z\cong\mathbf{lim}_{n}\prod_{i_{n}\in\mathcal{I}_{n}}(f_{i_{0}})^{*}(f_{i_{n}})_{*}E_{i_{n}}$$
   This brings us to the situation of a cover $\{\mathrm{Spec}(B'_{i})\rightarrow\mathrm{Spec}(A')\}$ where some $\mathrm{Spec}(B'_{i})=\mathrm{Spec}(A)$, and this clearly satisfies descent.
\end{proof}

 For finite $\kappa$ and stable $(\infty,1)$-categories the statement is just Proposition 10.5 in \cite{scholze2019lectures}, and our proof is a minor modification. 

 \section{Descent for Modules}\label{sec:descentmod}
 
In this section we specialise to the presheaf $\mathbf{QCoh}$.
 
 \subsection{Descendable Maps}

 In \cite{MR3459022} Mathew introduced the notion of a \textit{descendable} map $f:A\rightarrow B$ in $\mathbf{Comm}(\mathbf{C})$ for $\mathbf{C}$ a closed symmetric monoidal presentable $(\infty,1)$-category. Consider the cosimplicial object in $\mathbf{Comm}(\mathbf{C})$ given by the bar complex 
 $$\mathrm{CB}_{n}(f)\defeq B^{\otimes_{A}n+1}$$
 and its augmented version $\mathrm{CB}^{aug}_{\bullet}(f)$ with $\mathrm{CB}^{aug}_{-1}(f)=A$. Consider the pro-cosimplicial object given the tower of partial totalisations
 $$\{\mathrm{Tot}_{n}(\mathrm{CB}_{\bullet}(f))\}_{n\ge0}$$
 
  Recall that an object of $\mathbf{Pro(C)}$ is constant (\cite{MR3459022} Definition 3.9) if it is equivalent to an object in the image of $\mathbf{C}\rightarrow\mathbf{Pro(C)}$. 
 
 \begin{defn}[\cite{MR3459022} Definition 3.18/ Proposition 3.20]
 $A\rightarrow B$ is said to be \textit{descendable} if it is a constant pro-object which converges to $A$, i.e. 
 $$A\rightarrow\mathbf{lim}_{n}\mathrm{Tot}_{n}(\mathrm{CB}_{\bullet}(f))$$
 is an equivalence.
 \end{defn}
 
 Descendability can also be phrased in terms of the augmented complex. For each $n$ write
 $$T_{aug}^{n}(f)\defeq\mathrm{cofib}(A\rightarrow\mathrm{Tot}_{n}(\mathrm{CB}_{\bullet}(f))$$
 
 Then $A\rightarrow B$ is descendable precisely if $\{T_{aug}^{n}(f)\}_{n\ge0}$ is pro-constant and converges to the zero object. 
 
 \begin{defn}
 A tower $\{X_{n}\}_{n\ge0}$ is said to be \textit{nilpotent} if there exists an $N$ such that for all $n\ge0$ the map $X_{n+N}\rightarrow X_{n}$ is null-homotopic. 
 \end{defn}
 
 As explained in \cite{MR3459022}, nilpotent objects are pro-constant and converge to $0$, but in general nilpotence is in fact much stronger than this. However it turns out that a map is descendable precisely if $\{T_{aug}^{n}(f)\}_{n\ge0}$ is nilpotent. It is nilpotence which we will primarily use to establish descendability. 
 
 \begin{prop}[\cite{MR3459022} Proposition 3.22]
If $f:A\rightarrow B$ is descendable then the adjunction
$$\adj{B\otimes^{\mathbb{L}}_{A}(-)}{{}_{A}\mathbf{Mod(C)}}{{}_{B}\mathbf{Mod(C)}}{|-|}$$
is comonadic. Thus the functor
$${}_{A}\mathbf{Mod(\mathbf{C})}\rightarrow\mathbf{lim}_{n}{}_{B^{\otimes_{A}^{\mathbb{L}}n+1}}\mathbf{Mod(C)}$$
is an equivalence. 
 \end{prop}
 
 \begin{rem}
 There are potentially maps of algebras which satisfiy descent, in that the functor
$${}_{A}\mathbf{Mod(\mathbf{C})}\rightarrow\mathbf{lim}_{n}{}_{B^{\otimes_{A}^{\mathbb{L}}n+1}}\mathbf{Mod(C)}$$
is an equivalence, but which are not descendable in the sense of Mathew. Indeed as explained in \cite{MR3459022} it is not clear that all faithfully flat maps of $E_{\infty}$ algebras are descendable. However as we will explain in the next section, \cite{lurie2018spectral} Section D.4 says that any such map does indeed satisfy descent
 \end{rem}
 
 The class of descendable morphisms of commutative monoid objects forms a pre-topology. This follows immediately from the following two results.
 
 \begin{prop}[\cite{MR3459022} Proposition 3.23]
 Let $A\rightarrow B$ and $B\rightarrow C$ be maps in $\mathbf{Comm(C)}$. Then
 \begin{enumerate}
 \item
 If $A\rightarrow B$ and $B\rightarrow C$ are descendable, then so is $A\rightarrow C$.
 \item
 If $A\rightarrow C$ is descendable then so is $A\rightarrow B$.
 \end{enumerate}
 \end{prop}
 
 \begin{prop}[\cite{MR3459022} Corollary 3.21]\label{prop:basechangedescendable}
 Let $F:\mathbf{C}\rightarrow\mathbf{C'}$ be a symmetric monoidal functor between symmetric monoidal, stable $(\infty,1)$-categories. Let $f:A\rightarrow B$ be a descendable map in $\mathbf{Comm(C)}$. Then $F(f):F(A)\rightarrow F(B)$ is a descendable map in $\mathbf{Comm(C')}$.
 \end{prop}

\begin{cor}
Let $(\mathbf{C},\mathbf{D},\theta)$ be a stable $(\infty,1)$-algebra context. The collection of maps $\mathrm{Spec}(B)\rightarrow\mathrm{Spec}(A)$ such that $\Theta(A)\rightarrow\Theta(B)$ is descendable defines a pre-topology on $\mathbf{Alg_{D}}(\mathbf{C})^{op}$. 
\end{cor}

\begin{example}\label{Example:htpymonodescendable}
Let $(\mathbf{C},\mathbf{D},\theta)$ be a stable $(\infty,1)$-algebra context, and let $\{\mathrm{Spec}(A_{i})\rightarrow\mathrm{Spec}(A)\}_{i\in\mathcal{I}}$ be a collection of homotopy monomorphisms with $\mathcal{I}$ finite. Then the map $A\rightarrow B\defeq \underset{i\in\mathcal{I}}\prod A_{i}$ is descendable precisely if 
$$A\cong\mathbf{lim}B^{\otimes_{n}^{\mathbb{L}}A}$$
is an equivalence.
\end{example}

\begin{prop}
    Let $f:A\rightarrow B$ be a map in $\mathbf{DAlg}^{cn}(\mathbf{C})$ such that $B$ is formally $\aleph_{1}$-filtered as an $A$-module. If
    $$A\rightarrow\mathbf{lim}_{n}B^{\otimes_{A}^{\mathbb{L}}n}$$
    is an equivalence then $f$ is descendable.
\end{prop}

\begin{proof}
We have
  $$A\rightarrow\mathbf{lim}_{n}B^{\otimes_{A}^{\mathbb{L}}n}$$
  Since $B$ is is formally $\aleph_{1}$-filtered as an $A$-module so is each $T^{aug}_{m}(f)$. Thus we have
    $$T^{aug}_{m}(f)\cong \underset{n}{\mathbf{lim}} B^{\otimes^{\mathbb{L}}_{A}n}\otimes T^{aug}_{m}(f)$$
    So it suffices to prove that $\{B\otimes_{A}T^{m}_{aug}(f)\}$ is nilpotent. But $\{B\otimes^{\mathbb{L}}_{A}T^{m}_{aug}(f)\}\cong\{T^{m}_{aug}(B\times_{A}^(\mathbb{L}f))$ and the result follows since $B\rightarrow B\otimes_{A}^{\mathbb{L}}B$ is split, using \cite{MR3459022} Example 3.11.
\end{proof}

 \begin{lem}\label{lem:strongdesc}
     Let $A\rightarrow B$ be a derived strong map. Then $\pi_{0}(A)\rightarrow\pi_{0}(B)$ is descendable. 
     
     Conversely if $\pi_{0}(A)\rightarrow\pi_{0}(B)$ is descendable then for any $A$-module $M$ the map
     $$M\rightarrow\underset{n}{\mathbf{lim}}M\otimes_{A}^{\mathbb{L}}B^{\otimes_{A}^{\mathbb{L}}n}$$
     is an equivalence. In particular if $B$ is formally $\aleph_{1}$-filtered as an $A$-module then $A\rightarrow B$ is descendable.
 \end{lem}

 \begin{proof}
 Suppose that $A\rightarrow B$ is a derived strong map. Then $\pi_{0}(A)\rightarrow\pi_{0}(A)\otimes_{A}^{\mathbb{L}}B$ is descendable. Since $A\rightarrow B$ is strong we have $\pi_{0}(A)\otimes_{A}^{\mathbb{L}}B\cong\pi_{0}(B)$.

 Conversely suppose $\pi_{0}(A)\rightarrow\pi_{0}(B)$ is descendable. Let $M$ be an $A$-module. By descendability for $\pi_{0}(A)\rightarrow\pi_{0}(B)$, we have for each $m$
 $$\pi_{m}(M)\cong\underset{n}{\mathbf{lim}}\;\pi_{m}(M)\otimes_{\pi_{0}(A)}^{\mathbb{L}}\pi_{0}(B)^{\otimes_{\pi_{0}(A)}^{\mathbb{L}}n}\cong\underset{n}{\mathbf{lim}}\;\pi_{m}(M\otimes_{A}^{\mathbb{L}}B^{\otimes_{A}^{\mathbb{L}}n})$$
 where we have used derived strongness. A spectral sequence argument now implies the claim.

 \end{proof}

\subsection{Faithfully Flat Descent}

We finish with a generalisation of faithfully flat descent, which is significantly more straightforward. The proof is an easy adaptation of \cite{lurie2018spectral} Section D.4. Let $(\mathbf{C},\mathbf{C}_{\ge0},\mathbf{C}_{\le0},\mathbf{D},\theta)$ be a flat spectral algebraic context.

 Note that since projectives are strong and homotopy flat, strongly flat modules are in particular flat.
 
  \begin{defn}
 A map $f:A\rightarrow B$ in $\mathbf{Alg_{D}}(\mathbf{C})$ is said to be \textit{faithfully flat} if
 \begin{enumerate}
 \item
$B$ is flat as an $A$-module.
 \item
If $B\otimes^{\mathbb{L}}_{A}M$ is a trivial $B$-module then $M$ is a trivial $A$-module.
 \end{enumerate}
 \end{defn}
 
  \begin{lem}\label{Lem:faithfullyflatcofib}
 Let $f:A\rightarrow B$ be a map in $\mathbf{Alg}_{\mathbf{D}}^{cn}(\mathbf{C})$. If $f$ is faithfully flat  then the cofibre $\mathrm{cofib}(f)$ is flat when regarded as an $A$-module.
 \end{lem}
 
 \begin{proof}
 The proof is a minor modification of the second half of Lemma D.4.4.3 in \cite{lurie2018spectral}. Suppose $f$ is faithfully flat. Note that $\mathrm{cofib}(f)$ is connective since $A$ and $B$ are. Let $M$ be discrete. We have the fibre sequence
 $$M\rightarrow M\otimes^{\mathbb{L}}_{A}B\rightarrow M\otimes^{\mathbb{L}}_{A}\mathrm{cofib}(f)$$
 By assumption $M\otimes^{\mathbb{L}}_{A}B$ is discrete. To prove that $M\otimes^{\mathbb{L}}_{A}\mathrm{cofib}(f)$ is discrete it remains to show that $\theta:M\rightarrow M\otimes^{\mathbb{L}}_{A}B$ is a monomorphism in $\mathbf{C}^{\heart}$. Write $K=\mathrm{Ker}(\theta)$, which itself is discrete. We get a monomorphism $K\otimes_{\pi_{0}(A)}\pi_{0}(B)\rightarrow M\otimes_{\pi_{0}(A)}\pi_{0}(B)$. But the map $K\rightarrow M\otimes_{\pi_{0}(A)}\pi_{0}(B)$ is zero, so $K\otimes_{\pi_{0}(A)}\pi_{0}(B)\rightarrow M\otimes_{\pi_{0}(A)}\pi_{0}(B)$ is also zero. In particular $K\otimes_{\pi_{0}(A)}\pi_{0}(B)\cong 0$. Since $K$ is discrete and $A\rightarrow B$ flat, $K\otimes^{\mathbb{L}}_{A}B\cong K\otimes_{\pi_{0}(A)}\pi_{0}(B)\cong0$.
%
  \end{proof}

In the algebraic case, the converse to Lemma \ref{Lem:faithfullyflatcofib} is true using a derived version of Lazard's Theorem. In our case we have a partial converse.

\begin{prop}
Let $f:A\rightarrow B$ be a map in $\mathbf{Alg_{D}}^{cn}(\mathbf{C})$. If $\mathrm{cofib}(f)$ is strongly flat as an $A$-module then $f$ is faithfully flat.
\end{prop}

\begin{proof}
 The proof is as in the second half of Lemma D.4.4.3 in \cite{lurie2018spectral}. Suppose that $\mathrm{cofib}(f)$ is a strongly flat $A$-module. Write $\mathrm{cofib}(f)\cong\lim_{\rightarrow_{\alpha}}N_{\alpha}$. Pulling back the fibre sequence 
 $$A\rightarrow B\rightarrow \mathrm{cofib}(f)$$ 
 along the maps $N_{\alpha}\rightarrow \mathrm{cofib}(f)$ we get a commutative diagrams of $A$-modules
 \begin{displaymath}
 \xymatrix{
 A\ar[r]^{f} & B\ar[r] & \mathrm{cofib}(f)\\
 A\ar[r]\ar[u] & B_{\alpha}\ar[r]\ar[u] & N_{\alpha}\ar[u]
 }
 \end{displaymath}
 bottom sequence is a fibre sequence, and hence a cofibre sequence. Since $ A$ is connective and $N_{\alpha}$ is projective, the fibre sequence splits.  
  For any  right $A$-module $M$,
$$M\rightarrow B_{\alpha}\otimes_A^{\mathbb{L}} M\rightarrow N_{\alpha}\otimes_A^{\mathbb{L}} M$$
is a (split) fibre sequence. Thus $B_{\alpha}\otimes_A M$ is discrete, since $M$ and $N_{\alpha}\otimes_A M$ are. Since $B\otimes_A^{\mathbb{L}} M\cong \underset{\alpha}\colim B_{\alpha}\otimes_A^{\mathbb{L}} M$, $B\otimes^{\mathbb{L}}_{A}M$ is discrete. So $B$ is flat over $A$. Moreover since the sequence is split, for any projective generator $P$,
$$\Hom(P,M)\rightarrow \Hom(P,B_{\alpha}\otimes_A^{\mathbb{L}} M)\rightarrow \Hom(N_{\alpha}\otimes_A M)$$
is split. Thus $\pi_{*}\Hom(P,M)\rightarrow\pi_{*}\Hom(P,M\otimes^{\mathbb{L}}_{A}B_{\alpha})$ is a monomorphism. Taking colimits, $\pi_{*}\Hom(P,M)\rightarrow\pi_{*}\Hom(P,M\otimes^{\mathbb{L}}_{A}B)$ is a monomorphism. Thus if $M\otimes^{\mathbb{L}}_{A}B$ is trivial, so is $M$. 
 \end{proof}
 
 \begin{defn} 
Let $\mathbf{C}$ be a stable $\infty$-category with a t-structure. A monad $\mathrm{T}$ on $\mathbf{C}$ is called faithfully flat if the underlying functor $T:\mathbf{C}\to \mathbf{C}$ is right t-exact and preserves small colimits, and for every $X\in \mathbf{C}_{\leq0}$, the cofibre of the unit map $X\to TX$ also belongs to $\mathbf{C}_{\leq0}$.
\end{defn} 

By Lemma \ref{Lem:faithfullyflatcofib} we have the following.

\begin{cor}
Let $f:A\rightarrow B$ be faithfully flat. The monad \[B\otimes^{\mathbb{L}}_{A}(-):{}_{A}\mathbf{Mod}\rightarrow{}_{A}\mathbf{Mod}\] is faithfully flat.
\end{cor}

Since ${}_{A}\mathbf{Mod}(\mathbf{C})$ is a complete stable $(\infty,1)$-category, Proposition D.6.4.6 in \cite{lurie2018spectral} gives the following.
 
 \begin{cor}\label{cor:ffdesc}
Let $f:A\rightarrow B$ be faithfully flat. The functor $B\otimes^{\mathbb{L}}_{A}(-):{}_{A}\mathbf{Mod}(\mathbf{C})\rightarrow{}_{B}\mathbf{Mod}(\mathbf{C})$ is comonadic. In particular it satisfies descent.
 \end{cor}

 \subsection{The $\kappa$-Embedding pre-topology for Modules}
  Let $(\mathbf{C},\mathbf{C}_{\ge0},\mathbf{C}_{\le0},\mathbf{D},\theta)$ be a stable $(\infty,1)$-algebra context. Fix a class of maps $\mathbf{P}$ in $\mathbf{Aff}^{cn}_{\mathbf{C}}$ and a subpresheaf $\mathbf{N}\subseteq\mathbf{QCoh}|_{\mathbf{P}}$. Note that a homotopy monomorphism
  $$\mathrm{Spec}(B)\rightarrow\mathrm{Spec}(A)$$
  in $\mathbf{Alg_{D}}^{cn}(\mathbf{C})$, corresponds precisely to a \textit{homotopy epimorphism} $A\rightarrow B$, i.e. a map such that
  $$B\otimes_{A}^{\mathbb{L}}B\rightarrow B$$
  is an equivalence.

\begin{defn}
The $(\kappa,\mathbf{P},\mathbf{N})$-\textit{embedding pre-topology on }$\mathbf{Alg_{D}}^{op}$ is $\tau^{\mathbf{P},\mathbf{N},\kappa,mono}$.
\end{defn}

When $\kappa=\aleph_{0}$, $\mathbf{P}=\mathbf{hm}$, and $\mathbf{N}=\mathbf{QCoh}$ is the class of all homotopy monomorphisms, we denote by $\tau^{hm}$ the $(\aleph_{0},\mathbf{hm})$-embedding topology. 

By Propositions \ref{prop:htypepistrong}, Example \ref{Example:htpymonodescendable}, and Lemma \ref{lem:strongdesc} we have the following.

\begin{lem}
    Let $\{A\rightarrow B_{i}\}$ be a finite collection of maps between connective algebras. Suppose that each $\pi_{0}(A)\rightarrow\pi_{0}(B_{i})$ is a homotopy epimorphism and that these maps are transverse to each $\pi_{m}(A)$. Then the following are equivalent.
    \begin{enumerate}
        \item 
          $\{\mathrm{Spec}(B_{i})\rightarrow\mathrm{Spec}(A)\}_{\in\mathcal{I}}$ is a cover in $\tau^{hm}$.
          \item 
          The map $A\rightarrow\prod B_{i}$ is descendable.
          \item 
          The map $\pi_{0}(A)\rightarrow\prod\pi_{0}(B_{i})$ is descendable and each map $A\rightarrow B_{i}$ is derived strong.
          \item 
                     $\{\mathrm{Spec}(\pi_{0}(B_{i}))\rightarrow\mathrm{Spec}(\pi_{0}(A))\}_{\in\mathcal{I}}$ is a cover in $\tau^{hm}$  and each map $A\rightarrow B_{i}$ is derived strong.
    \end{enumerate}
\end{lem}

When $\kappa\ge\aleph_{1}$ the situation becomes a little more subtle. We will generally restrict to a category $\mathbf{A}\subseteq\mathbf{Aff}^{cn}_{\mathbf{C}}$ such that maps in $\mathbf{P}|_{\mathbf{A}}$ are derived strong. We will consider a sub-presheaf $\mathbf{N}\subseteq\mathbf{QCoh}|_{\mathbf{A}}$ such that for any $M\in\mathbf{N}(\mathrm{Spec}(A))$, any integer $m$, and any cover $\{\mathrm{Spec}(\pi_{0}(B_{i}))\rightarrow\mathrm{Spec}(\pi_{0}(A))\}_{i\in\mathcal{I}}$, the map 

$$\pi_{m}(M)\rightarrow\mathbf{lim}_{n}(\prod_{(i_{1},\ldots,i_{n})\in\mathcal{I}^{n})}\pi_{m}(M)\otimes_{\pi_{0}(A)}^{\mathbb{L}}\pi_{0}(B_{i_{1}})\otimes_{\pi_{0}(A)}^{\mathbb{L}}\ldots\otimes_{\pi_{0}(A)}^{\mathbb{L}}\pi_{0}(B_{i_{n}}))$$
is an equivalence. We then extend $\mathbf{N}$ to a sub-presheaf on $\mathbf{Aff}^{cn}_{\mathbf{C}}$ by declaring $M\in\mathbf{N}(\mathrm{Spec}(C))$ if for every map $C\rightarrow A$ with $A\in\mathbf{A}$, $A\otimes_{C}^{\mathbb{L}}M\in\mathbf{N}(\mathrm{Spec}(A))$. 

\begin{lem}
    Let $\kappa$ be a cardinal and let $\mathbf{N}\subseteq\mathbf{QCoh}_{\mathbf{P}}$ be a sub-presheaf. Let $\{\mathrm{Spec}(B_{i})\rightarrow\mathrm{Spec}(A)\}$ be a $\kappa$-small collection of maps. 
    
        Suppose that each $\pi_{0}(A)\rightarrow\pi_{0}(B_{i})$ is a homotopy epimorphism, and that these maps are transverse to each $\pi_{m}(A)$. Suppose that for each $M\in\mathbf{N}(\mathrm{Spec}(A))$ we have $\pi_{m}(M)\in\mathbf{N}(\mathrm{Spec}(\pi_{0}(A)))$. Then the following are equivalent.
    \begin{enumerate}
        \item 
         $\{\mathrm{Spec}(B_{i})\rightarrow\mathrm{Spec}(A)\}$ is a cover in $\tau^{\mathbf{hm},\mathbf{N},\kappa,mono}$.
         \item 
            $\{\mathrm{Spec}(\pi_{0}(B_{i}))\rightarrow\mathrm{Spec}(\pi_{0}(A))\}$ is a cover in $\tau^{\mathbf{hm},\mathbf{N},\kappa,mono}$.
    \end{enumerate}
\end{lem}

\begin{proof}
    Let  $\{\mathrm{Spec}(B_{i})\rightarrow\mathrm{Spec}(A)\}$ is a cover in $\tau^{\mathbf{hm},\mathbf{N},\kappa,mono}$. Since $\pi_{0}(A)\rightarrow\pi_{0}(B_{i})$ is a homotopy monomorphism and each $\pi_{m}(A)$ is transverese to these maps, each map $A\rightarrow B_{i}$ is derived strong and the claim is immediate. 

    Conversely suppose   $\{\mathrm{Spec}(\pi_{0}(B_{i}))\rightarrow\mathrm{Spec}(\pi_{0}(A))\}$ is a cover in $\tau^{\mathbf{hm},\mathbf{N},\kappa,mono}$. The maps $A\rightarrow B_{i}$ are still derived strong. Let $M\in\mathbf{N}(\mathrm{Spec}(A))$. Then $\pi_{m}(B_{i}\otimes_{A}^{\mathbb{L}}M)\cong\pi_{0}(B_{i})\otimes_{A}^{\mathbb{L}}\pi_{n}(M)$. If $B_{i}\otimes_{A}^{\mathbb{L}}M$ is acyclic. then each $\pi_{0}(B_{i})\otimes_{A}^{\mathbb{L}}\pi_{n}(M)$ is zero. Thus $\pi_{n}(M)\cong0$ for all $m$ as required.
\end{proof}




\subsection{Transverse \v{C}ech Descent and Hyperdescent}

In general it is most definitely not the case that \v{C}ech descent implies hyperdescent. However it is true in certain situations.

\begin{lem}[c.f. \cite{toen2008homotopical} Lemma 2.2.2.13.]\label{lem:hyperdescentheart}
Let $(\mathbf{Aff}^{cn}_{\mathbf{C}},\tau,\mathbf{P},\mathbf{A})$ be a relative $(\infty,1)$-pre-geometry tuple. Let $\mathbf{Q}\subset\mathbf{QCoh}|_{\mathbf{A}}$ be a sub-presheaf. Write $\mathbf{A}^{\heart} \subset\mathbf{A}$ for the full subcategory consisting of objects $\mathrm{Spec}(A)$ with $A\in\mathbf{DAlg}^{\heart}(\mathbf{C})$. Finally denote by $\mathbf{Q}^{\heart}:(\mathbf{A}^{\heart})^{op}\rightarrow\mathrm{Cat}$ the assignment which sends $\mathrm{Spec}(A)$ to $\mathbf{Q}(A) \cap {}_{A}\mathbf{Mod}(\mathbf{C})^{\heart}$ Suppose that 
\begin{enumerate}
\item
$\tau|_{\mathbf{A}^{\heart}}$ defines a pre-topology on $\mathbf{A}^{\heart}$
\item 
objects of $\mathbf{Q}^{\heart}$ are transverse to maps in $\mathbf{P}|_{\mathbf{A}^{\heart}}$
\item
$\mathbf{Q}^{\heart}$ is a subpresheaf of $\mathbf{QCoh}^{\heart}$ on $\mathbf{A}^{\heart}$ which satisfies \v{C}ech descent for $\tau|_{\mathbf{A}^{\heart}}$
\item 
If $X\in\mathbf{Q}(\mathrm{Spec}(A))$ then $\pi_{n}(X)\in\mathbf{Q}^{\heart}(\mathrm{Spec}(A))$ for all $n$.
\item
if $\mathrm{Spec}(A)\in\mathbf{A}$ then $\mathrm{Spec}(\pi_{0}(A))\in\mathbf{A}$.
\item
whenever $\{\mathrm{Spec}(A_{i})\rightarrow\mathrm{Spec}(A)\}$ is a cover in $\tau|_{\mathbf{A}}$,  \[\{\mathrm{Spec}(\pi_{0}(A_{i}))\rightarrow\mathrm{Spec}(\pi_{0}(A))\}\] is a cover in $\tau|_{\mathbf{A}^{\heart}}$
\item
if $f:\mathrm{Spec}(B)\rightarrow\mathrm{Spec}(A)$ is a map in $\mathbf{P}|_{\mathbf{A}}$ then the corresponding map $A\rightarrow B$ is derived strong.
\end{enumerate}
Then $\mathbf{Q}|_{\mathbf{A}}$ satisfies hyperdescent.
%
\end{lem}

\begin{proof}
 Let $K_{\bullet}\rightarrow\mathrm{Spec}(A)$ be a pseudo-representable hypercover with $K_{n}=\coprod_{i_{n}\in\mathcal{I}_{n}} \mathrm{Spec}(B_{i_{n}})$. Let $E_{\bullet}\in\mathbf{QCoh}(|K|)$, whose value on $B_{i_{n}}$ is $E_{i_{n}}$. Since $A\rightarrow B_{i_{0}}$ is derived strong we have
$$\mathrm{Tor}^{p}_{\pi_{0}(A)}(\pi_{*}(\mathbf{lim}_{n}\prod E_{i_{n}}),\pi_{0}(B_{i_{0}}))_{q}\Rightarrow\pi_{*}((\mathbf{lim}_{n}\prod E_{i_{n}})\otimes_{A}^{\mathbb{L}}B_{i_{0}})$$
Consider the spectral sequence
$$\pi_{p}(\mathrm{Tot}(\pi_{q}(\prod E_{i_{n}})))\Rightarrow\pi_{q-p}(\mathbf{lim}_{n}\prod E_{i_{n}})$$
 Since each $\pi_{q}(E_{i_{n}})$ is transverse to $\pi_{0}(A)$ over $\pi_{0}(B_{i_{n}})$, we have that $\pi_{q}(E_{i_{n}})$  defines a well-defined object of $\mathbf{QCoh}(\pi_{0}(K_{\bullet}))$. Now \v{C}ech  for groupoids implies hyperdescent, so we find that $\pi^{p}(\mathrm{Tot}(\pi_{q}(\prod E_{i_{\bullet}})))\cong 0$ for $p\neq 0$. Thus the spectral sequence degenerates to give an isomorphism
$$\pi_{p}(\mathbf{lim}_{n}E_{i_{n}})\cong\mathrm{ker}(\pi_{p}(\prod E_{i_{0}})\rightarrow\pi_{p}(\prod E_{i_{1}}))$$
Therefore $\pi_{p}(\mathbf{lim}_{n}\prod E_{i_{n}})$ is the $\pi_{0}(A)$-module obtained by descent from $\pi_{p}(E_{i_{n}})$ on $t_{\le0}(K_{\bullet})$. Thus we get an isomorphism in ${}_{\pi_{0}(B_{i_{0}})}\mathbf{Mod}(\mathbf{C})^{\heart}$
$$\pi_{p}(\mathbf{lim}_{n}\prod E_{i_{n}})\otimes_{\pi_{0}(A)}^{\mathbb{L}}\pi_{0}(B_{i_{0}})\rightarrow\pi_{p}(E_{i_{0}})$$
Since $A\rightarrow B_{i_{0}}$ is strong this gives that the map
$$E_{i_{0}}\rightarrow (\mathbf{lim}_{n}\prod E_{i_{n}})\otimes_{A}^{\mathbb{L}}B_{i_{0}}$$
is an equivalence, as required.
\end{proof}

\subsubsection{Perfect Descent}

In quite general circumstances, $\mathbf{Perf}_{\ge n}$ and $\mathbf{Perf}$ are stacks.

\begin{lem}
Suppose that $\mathbf{N}\subset\mathbf{QCoh}(-)$ is such that
\begin{enumerate}
\item
$\mathbf{N}$ satisfies \v{C}ech descent.
\item
$\mathbf{N}(\mathrm{Spec}(A))$ contains $A$.
\item
$\mathbf{N}(\mathrm{Spec}(A))\subset{}_{A}\mathbf{Mod}$ is closed under the tensor product and internal hom.
\end{enumerate}
 Then $\mathbf{Perf}(-)$ is a local sub-presheaf of $\mathbf{N}$. In particular it satisfies \v{C}ech descent, and also satisfies hyperdescent whenever $\mathbf{N}$ does.
\end{lem}

\begin{proof}
The proof is as in in \cite{toen2008homotopical} Corollary 1.3.7.4. 
Let $\{f_{i}:U_{i}=\mathbf{Spec}(B_{i})\rightarrow \mathbf{Spec}(A)=U\}$ be a cover, and $P,Q\in{}_{A}\mathbf{Mod}$ where $Q\in\mathbf{N}(\mathrm{Spec}(A))$ and $P$ is such that $B_{i}\otimes^{\mathbb{L}}_{A}P\in\mathbf{Perf}(B_{i})$ for each $i$. Define
$$\mathbf{Map}_{A}(P,Q)_{(i_{1},\ldots,i_{n})}\defeq\mathbf{Map}_{A}(P,Q)\otimes^{\mathbb{L}}_{A}B_{(i_{1},\ldots,i_{n})}$$
$$\mathbf{Map}_{A}(P,Q_{(i_{1},\ldots,i_{n})})\defeq\mathbf{Map}_{A}(P,Q\otimes^{\mathbb{L}}_{A}B_{(i_{1},\ldots,i_{n})})$$
where \[B_{(i_{1},\ldots,i_{n})}\defeq B_{i_{1}}\otimes_{A}^{\mathbb{L}}B_{i_{2}}\otimes^{\mathbb{L}}_{A}\cdots\otimes^{\mathbb{L}}_{A}B_{i_{n}}.\] For $*=(i_{1},\ldots,i_{n})$ there is a natural map
$$\mathbf{Map}_{A}(P,Q)_{*}\rightarrow\mathbf{Map}_{A}(P,Q_{*})$$
$\mathbf{Map}_{A}(P,Q)_{*}$ is clearly a well-defined object of $\mathbf{QCoh}(\check{C}(\{U_{i}\rightarrow U\}))$.

 $\mathbf{Map}_{A}(P,Q_{*})$ is also a a well-defined object of \[\mathbf{QCoh}(\check{C}(\{U_{i}\rightarrow U\})).\] This follows from the identifications
\begin{equation}
    \begin{split}
    \mathbf{Map}_{A}(P,Q_{*})\otimes^{\mathbb{L}}_{B_{*}}B_{i}& \cong \mathbf{Map}_{B}(B_{*}\otimes^{\mathbb{L}}_{A}P,Q_{*})\otimes_{B_{*}}B_{i} \\ &\cong\mathbf{Map}_{B}(B_{*}\otimes^{\mathbb{L}}_{A}P,Q_{*}\otimes^{\mathbb{L}}_{B_{*}}B_{i})
    \end{split}
\end{equation}
where we have used that $B_{*}\otimes_{A}P$ is perfect. We then have
\begin{align*}
\lim_{n,(i_{1},\ldots,i_{n})\in\mathcal{I}^{n}}\mathbf{Map}_{A}(P,Q)_{(i_{1},\ldots,i_{n})}&\cong\mathbf{Map}_{A}(P,Q)\otimes^{\mathbb{L}}_{B_{*}}B_{(i_{1},\ldots,i_{n})}\\
&\cong\lim_{n,(i_{1},\ldots,i_{n})\in\mathcal{I}^{n}}\mathbf{Map}_{A}(P,Q\otimes^{\mathbb{L}}_{A}B_{(i_{1},\ldots,i_{n})})\\
&\cong\lim_{n,(i_{1},\ldots,i_{n})\in\mathcal{I}^{n}}\mathbf{Map}_{A}(P,Q_{(i_{1},\ldots,i_{n})})
\end{align*}
Thus by descent, the map
$$\mathbf{Map}_{A}(P,Q)\otimes^{\mathbb{L}}_{A}B_{i}\rightarrow\mathbf{Map}_{A}(P,Q\otimes^{\mathbb{L}}_{A}B_{i})$$
is an isomorphism. Setting $Q=A$ gives that
$$\mathbf{P}^{\vee}\otimes_{A}^{\mathbb{L}}B_{i}\rightarrow\mathbf{Map}_{B_{i}}(P\otimes_{A}^{\mathbb{L}}B_{i},B_{i})$$
is an isomorphism. Since $P\otimes_{A}^{\mathbb{L}}B$ is perfect, this implies that for any $A$-module $Q$ we have
$$P^{\vee}\otimes_{A}^{\mathbb{L}}Q\otimes_{A}^{\mathbb{L}}B_{i}\cong\mathbf{Map}_{A}(P,Q)\otimes^{\mathbb{L}}_{A}B_{i}$$
is an equivalence. Thus again by descent we have that
$$P^{\vee}\otimes_{A}^{\mathbb{L}}Q\rightarrow\mathbf{Map}_{A}(P,Q)$$
is an equivalence for any $A$-module $Q$, i.e. $P$ is perfect.
\end{proof}

In applications later $\mathbf{N}$ will be $\mathbf{Coh}_{+}$. We also then have the following.

\begin{cor}\label{cor:perfdescent}
    Suppose that $(\mathbf{A},\mathbf{P})$ satisfies Theorems A and B. Then $\mathbf{Perf}_{\ge n}$ satisfies hyperdescent for each $n$. If the pre-topology consists of finite covers then $\mathbf{Perf}$ satisfies hyperdescent.
\end{cor}

\subsubsection{Nuclear Descent}

Suppose that all nuclear objects in $\mathbf{C}$ are strongly nuclear. Let $\tau$ be a pre-topology on $\mathbf{A}\subseteq\mathbf{Alg_{D}}(\mathbf{C})$ such that if $\mathrm{Spec}(A)\in\mathbf{A}$, then $A$ is strongly nuclear. We then get a sub-sheaf $\mathbf{Nuc}\subseteq\mathbf{QCoh}|_{\mathbf{A}}$.

\begin{lem}
$\mathbf{Nuc}$ is a local subsheaf of $\mathbf{QCoh}$ for the finite homotopy monomorphism topology.
\end{lem}

\begin{proof}
    Let $M$ be an $A$-module such that each $B_{i}\otimes^{\mathbb{L}}_{A}M$ is a nuclear $B_{i}$-module. Let $c$ be a compact $A$-module. Consider the map
    $$\underline{\mathbf{Map}}_{A}(A\otimes^{\mathbb{L}}c,A)\otimes^{\mathbb{L}}_{A}M\rightarrow\underline{\mathbf{Map}}_{A}(A\otimes^\mathbb{L}c,M)$$
    which by the natural isomorphisms is equivalent to the map
    $$\underline{\mathbf{Map}}(c,|A|)\otimes_{A}^{\mathbb{L}}M\rightarrow\underline{\mathbf{Map}}(c,|M|)$$

 Now since $B_{i}\otimes_{A}^{\mathbb{L}}M$ is nuclear as a $B_{i}$-module and $B_{i}$ are nuclear as $A$-modules, we have 
\begin{align*}
\underline{\mathbf{Map}}(c,|A|)\otimes_{A}^{\mathbb{L}}B_{i}\otimes_{A}^{\mathbb{L}}M&\cong\underline{\mathbf{Map}}(c,B_{i})\otimes_{A}^{\mathbb{L}}M\\
&\cong\underline{\mathbf{Map}}(c,B_{i})\otimes_{B_{i}}^{\mathbb{L}}B_{i}\otimes_{A}^{\mathbb{L}}M\\
&\cong\underline{\mathbf{Map}}(c,B_{i}\otimes^{\mathbb{L}}_{A}M)\\
&\cong\underline{\mathbf{Map}}_{B_{i}}(B_{i}\otimes^{\mathbb{L}}c,B_{i}\otimes^{\mathbb{L}}_{A}M)\\
\end{align*}
    Similarly for any $B_{(i_{1},\ldots,i_{n})}$ with $(i_{1},\ldots i_{n})\in\mathcal{I}^{n}$ we also have
        $$|\underline{\mathbf{Map}}(c,A)\hat{\otimes}_{A}B_{(i_{1},\ldots,i_{n})}\otimes_{A}^{\mathbb{L}}M|\cong\underline{\mathbf{Map}}(B_{i}\otimes^{\mathbb{L}}c,B_{(i_{1},\ldots,i_{n})}\otimes_{A}^{\mathbb{L}}M)$$
        Taking limits over $\mathcal{I}^{n}$, and using that this is a finite limit, gives that 
         $$\underline{\mathbf{Map}}(c,|A|)\otimes_{A}^{\mathbb{L}}M\rightarrow\underline{\mathbf{Map}}(c,|M|)$$
         is an equivalence, which proves that $M$ is a nuclear $A$-module.
\end{proof}

\subsection{Descent for Properties of Maps}

We conclude with some results regarding descent for properties of maps.

\begin{lem}
Let $\mathbf{Q}:\mathbf{A}^{op}\rightarrow\mathbf{Pr^{L}}$ be a presheaf, $\tau$ a pre-topology conisting of flat maps which satisfies \v{C}ech descent, and $\{f_{i}:U_{i}\rightarrow V\}_{i\in\mathcal{I}}$ be a cover which satisfies descent for $\mathbf{Q}$, and $g:V\rightarrow W$ a map. If 
$$g_{I}:U_{i_{1}}\times_{V} U_{i_{2}}\times_{V}\cdots \times_{V} U_{i_{n}}\rightarrow W$$ 
is flat (in the sense that $g^{*}$ commutes with finite limits) for each $I=\{i_{1},\ldots,i_{n}\}\subset\mathcal{I}$, then $g$ is flat.
\end{lem}

\begin{proof}
We have $g^{*}=\lim_{I}g^{*}_{I}$. Each $g^{*}_{I}$ commutes with finite limits, so $g^{*}$ clearly does as well.
\end{proof}

\begin{lem}\label{lem:descentforsmoothmaps}
Let $\{f_{i}:U_{i}\rightarrow V\}_{i\in\mathcal{I}}$ be a cover in $(\mathbf{Alg_{D}}^{cn})^{op}$ consisting of formally \'{e}tale maps, and which satisfies \v{C}ech descent for $\mathbf{Perf}_{\ge0}$. Let $g:V\rightarrow W$ be a map. If for all $I\subset\mathcal{I}$
$$g_{I}:U_{i}\times_{V} U_{i_{2}}\times_{V}\cdots \times_{V} U_{i_{n}}\rightarrow W$$ 
is formally smooth of relative dimension $k$, then $g$ is formally smooth of relative dimension $k$
\end{lem}

\begin{proof}
Write $W=\mathrm{Spec}(B)$, $V=\mathrm{Spec}(A)$, $U_{i}=\mathrm{Spec}(A_{i})$
Observe that $\mathbb{L}_{V\big\slash W}$ is the limit of the effective descent data $(\mathbb{L}_{B\big\slash  (A_{i_{1}}\otimes^{\mathbb{L}}_{A}\ldots\otimes_{A}^{\mathbb{L}}A_{i_{n}})})\cong (A_{i_{1}}\otimes^{\mathbb{L}}_{A}\ldots\otimes_{A}^{\mathbb{L}}A_{i_{n}})^{\oplus k}$. Hence   $\mathbb{L}_{B\big\slash A}\cong A^{\oplus n}$. Moreover each $\mathbb{L}_{A_{i}}\otimes_{A_{i}}^{\mathbb{L}}B\rightarrow\mathbb{L}_{B}$ has a section. But $\mathbb{L}_{A_{i}}\cong\mathbb{L}_{A}\otimes_{A}^{\mathbb{L}}A_{i}$, so that $\mathbb{L}_{A_{i}}\otimes_{A_{i}}^{\mathbb{L}}B\cong\mathbb{L}_{A}\otimes_{A}^{\mathbb{L}}B$ has a section.
\end{proof}

\subsubsection{Descent of Transversality}

\begin{lem}\label{lem:decsenttransv}
Let $f:A\rightarrow B$ be a map in $\mathbf{Alg_{D}}(\mathbf{C}^{\heart})$. Let \[\{\mathrm{Spec}(B_{i})\rightarrow\mathrm{Spec}(B)\}_{i\in\mathcal{I}}\] be a cover with each $B_{i}\in\mathbf{Alg_{D}}(\mathbf{C}^{\heart})$. Let $M\in{}_{A}\mathrm{Mod}(\mathbf{C}^{\heart})$ be such that
\begin{enumerate}
\item
for each $B_{i}$ the map
$$B_{i}\otimes^{\mathbb{L}}_{A}M\rightarrow B_{i}\otimes_{A}M$$
is an equivalence.
\item
for each $i,j\in\mathcal{I}$ the map
$$B_{i}\otimes^{\mathbb{L}}_{A}B_{j}\rightarrow B_{i}\otimes_{A}B_{j}$$
is an equivalence.
\item
the map
$$B\otimes_{A}^{\mathbb{L}}M\rightarrow \underset{n\in\mathbb{N},(i_{1},\ldots,i_{n})\in\mathcal{I}^{n}}{\mathbf{lim}} B_{i_{1}}\otimes_{A}^{\mathbb{L}}\cdots \otimes_{A}^{\mathbb{L}}B_{i_{n}}\otimes_{A}^{\mathbb{L}}M$$
is an equivalence.
\end{enumerate}
Then the map
$$B\otimes^{\mathbb{L}}_{A}M\rightarrow B\otimes_{A}M$$
is an equivalence.
\end{lem}

\begin{proof}
$B\otimes_{A}^{\mathbb{L}}M$ can only have non-zero homology in non-negative degrees, and \[\underset{n\in\mathbb{N},(i_{1},\ldots,i_{n})\in\mathcal{I}^{n}}{\mathbf{lim}} B_{i_{1}}\otimes_{A}^{\mathbb{L}}\cdots \otimes_{A}^{\mathbb{L}}B_{i_{n}}\otimes_{A}^{\mathbb{L}}M\cong \underset{n\in\mathbb{N},(i_{1},\ldots,i_{n})\in\mathcal{I}^{n}}{\mathbf{lim}} B_{i_{1}}\otimes_{A}\cdots \otimes_{A} B_{i_{n}}\otimes_{A}M\]
can only have non-zero homology in non-positive degree. Thus  $B\otimes_{A}^{\mathbb{L}}M$ can only have non-zero homology in degree $0$, as required.
\end{proof}

\subsubsection{A Note on Descent and Subcanonicity of Topologies}

Descent for certain subcategories of $\mathbf{Qcoh}$ is often related to subcaonicity of topologies. This fact is also used extensively in \cite{RhiannonRepresentability}.

\begin{prop}\label{prop:descentsubcan}
    Let $(\mathbf{Aff}^{cn}_{\mathbf{C}},\tau,\mathbf{P},\mathbf{A})$ be a relative $(\infty,1)$-pre-geometry tuple. If whenever $K_{\bullet}\rightarrow\mathrm{Spec}(A)$ is a pseudo-representable hypercover (resp. a \v{C}ech cover) with $K_{n}\cong\coprod_{i_{n}\in\mathcal{I}_{n}}\mathrm{Spec}(B_{i_{n}})$ with $\mathrm{Spec}(A)\in\mathbf{A}$ then the map
    $$A\rightarrow\underset{n\in\mathbb{N}}{\mathbf{lim}}\prod_{i_{n}\in\mathcal{I}_{n}} B_{i_{n}}$$
    is an equivalence, then the pre-topology is hyper-subcanonical (resp. \v{C}ech subcanonical). In particular if $\mathbf{N}\subseteq\mathbf{QCoh}|_{\mathbf{A}}$ is a subpresheaf satisfying hyperdescent (resp. \v{C}ech descent) with $A\in\mathbf{N}(\mathrm{Spec}(A))$ for each $A\in\mathbf{A}$, then the pre-topology is hyper-subcanonical (resp. \v{C}ech subcanonical).
\end{prop}

\begin{proof}
    The proof is identical to \cite{toen2008homotopical} Corollary 1.3.2.5. 
\end{proof}

\begin{cor}
   Let $$(\mathbf{C},\mathbf{D},\theta)$$
   be an $(\infty,1)$-algebraic context, $\mathbf{P}$ a class of maps, and $\tau$ a pre-topology on $(\mathbf{Alg}_{D})^{op}$ which is contained in the descendable pre-topology and whose covers consist of maps in $\mathbf{P}$. Then
    $$(\mathbf{C},\tau,\mathbf{P})$$
    is an $(\infty,1)$-geometry tuple.
\end{cor}

\begin{cor}\label{cor:Postnikovdecsendablecontext}
   Let $$(\mathbf{C},\mathbf{C}_{\ge0},\mathbf{C}_{\le0},\mathbf{C}^{0},\mathbf{D},\theta)$$
   be a Postnikov algebraic context, $\mathbf{P}$ a class of maps, and $\tau$ a pre-topology on $(\mathbf{Alg}_{D}^{cn})^{op}$ which is contained in the descendable pre-topology and whose covers consist of maps in $\mathbf{P}$. Then
    $$(\mathbf{C}_{\ge0},\tau,\mathbf{P})$$
    is an $(\infty,1)$-geometry tuple.
\end{cor}

\chapter{Relative Geometry}\label{RA}

In this somewhat technical chapter we will apply work of the previous chapters to develop some aspects of geometry relative to contexts, including some basic obstruction theory, the rudiments of geometry relative to Lawvere theories, and analyitification.

 \section{Geometry Relative to an $(\infty,1)$-Algebra Context}\label{sec:geometryrelac}

Here we introduce a formalism for geometry relative to spectral algebraic contexts. Most of the material in this section is a straightforward generalisation of \cite{kelly2021analytic}. We fix an $(\infty,1)$-algebra context (Definition \ref{defn:infA}).
$$(\mathbf{C},\mathbf{D},\theta,\mathbf{G}_{0},\mathbf{S})$$

\begin{defn}
Let
\begin{enumerate}
\item
$$\underline{\mathbf{C}}=(\mathbf{C},\mathbf{D},\theta,\mathbf{G}_{0},\mathbf{S})$$

be an $(\infty,1)$-algebra context.
\item
$\mathbf{A}\subset\mathbf{Aff}_{\underline{\mathbf{C}}}$ be a full subcategory.
\item
$\tau$ be a pre-topology on $\mathbf{Aff}_{\underline{\mathbf{C}}}$
\item
$\mathbf{P}$ be a class of maps in $\mathbf{Aff}_{\underline{\mathbf{C}}}$
\end{enumerate}
such that
$$(\mathbf{Aff}_{\underline{\mathbf{C}}},\tau,\mathbf{P},\mathbf{A})$$
is a strong relative $(\infty,1)$-geometry tuple (Definitions \ref{defn:strongtupdesc} and \ref{defn:strongtup}). A \textit{good system of module categories relative to} $(\mathbf{A},\tau,\mathbf{P})$, denoted, $\mathbf{F}$ is a sub-pre-sheaf $\mathbf{F}\subset\mathbf{QCoh}|_{\mathbf{A}}$ such that
\begin{enumerate}
\item
$A\in\mathbf{F}_{\mathrm{Spec}(A)}$ for any $\mathrm{Spec}(A)\in\mathbf{A}$;
\item
$\mathbb{L}_{\mathrm{Spec}(A)}\in\mathbf{F}_{A}$ for any $\mathrm{Spec}(A)\in\mathbf{A}$. 
\item
$\mathbf{F}(X)$ is closed under equivalences, finite colimits, and retracts;
\item
$f^{*}M\in\mathbf{F}_{Y}$ for any map $f:Y\rightarrow X$ and any $M\in\mathbf{F}_{Y}$. 
\item
 $\mathbf{F}$ satisfies weak \v{C}ech descent;
\end{enumerate}
\end{defn}

Moreover for each $n\in\mathbb{Z}$ we also define
$$\mathbf{F}(X)_{\ge n}\defeq\mathbf{F}(X)\times_{\mathbf{QCoh}(X)}\mathbf{QCoh}(X)_{\ge n}$$
$$\mathbf{F}(X)^{cn}\defeq\mathbf{F}(X)\times_{\mathbf{QCoh}(X)}\mathbf{QCoh}(X)^{cn}$$

Let $X=\mathrm{Spec}(A)\in\mathbf{Aff}$ and let $M\in{}_{A}\mathbf{Mod}_{\ge 1}$. Consider the trivial square-zero extension $X[M]$, and let $d:A\rightarrow A\oplus M$ be a derivation. We define $X_{d}[M]$ to be $\mathrm{Spec}(A\oplus_{d}M)$, where $A\oplus_{d}\Omega M$ is defined by the pullback diagram below.
\begin{displaymath}
\xymatrix{
A\oplus_{d}\Omega M\ar[d]\ar[r] & A\ar[d]^{i}\\
A\ar[r]^{d} & A\oplus M
}
\end{displaymath}
Here $i:A\rightarrow A\oplus M$ is the trivial derivation. 

The following is motivated, as explained in \cite{kelly2021analytic}, by obstruction theory in the context of HKR theorems.

\begin{defn}
A \textit{weak relative \'{e}tale} $(\infty,1)$-\textit{AG context} is a tuple $(\mathbf{Aff}_{\underline{\mathbf{C}}},\tau,\mathbf{P},\mathbf{A},\mathbf{F})$ where
\begin{enumerate}
\item
$\underline{\mathbf{C}}=(\mathbf{C},\mathbf{D},\theta,\mathbf{G}_{0},\mathbf{S})$ is an $(\infty,1)$-algebra context;
\item
$\tau$ is a Grothendieck pre-topology on $\mathbf{Aff}_{\underline{\mathbf{C}}}$;
\item
$\mathbf{P}$ is a class of maps in $\mathsf{Ho}(\mathbf{Aff}_{\underline{\mathbf{C}}})$ which is stable by composition, pullback, and contains isomorphisms;
\item
$\mathbf{P}$ consists of formally \'{e}tale morphisms;
\item
if $\{U_{i}\rightarrow U\}$ is a cover in $\tau$ then each $U_{i}\rightarrow U$ is in $\mathbf{P}$.
\item
$\mathbf{F}$ is a good system of modules on $\mathbf{A}$,
\item
For $M\in{}_{A}\mathbf{Mod}^{\mathbf{G}}_{1}\cap\mathbf{F}$  $X[M]\defeq\mathrm{Spec}(A\oplus M)\in\mathbf{A}$. 
\end{enumerate}
If $\underline{\mathbf{C}}$ is either a Postnikov algebraic context or a connective Postnikov algebraic context then we call $(\mathbf{Aff}_{\underline{\mathbf{C}}},\tau,\mathbf{P},\mathbf{A},\mathbf{F})$ a weak \textit{weak relative \'{e}tale} $(\infty,1)$-\textit{SAG context}, where SAG stands for \textit{spectral algebraic geometry}.
\end{defn}


This somewhat cumbersome amount of data can be greatly simplified if we are working in a connective Postnikov algebraic context setting and/ or the pre-topology satisfies descent for all modules. In the connective setting we can take $\mathbf{S}=\mathbf{Alg_{D}}$ and $\mathbf{G}=\mathbf{C}_{\ge0}$, and in the case that we have descent for all modules we can remove take $\mathbf{F}=\mathbf{Mod}$ and $\mathbf{A}=\mathbf{Alg_{D}}^{op}$. However there are cases of interests where at least one, and possibly both, of these simplifications cannot be made, namely in the setting of complex analytic geometry. Here we can only hope to have descent for metrisable, and in particular coherent, sheaves on Stein spaces. Thus we must restrict $\mathbf{F}$ to be so-called $RR$-quasicoherent modules, and $\mathbf{S}$ to be derived/ coherent Stein spaces. On the other hand the restrictions on $\mathbf{S}$ and $\mathbf{G}$ will be important for obstruction theory, particularly if we want to work with non-connective objects and obtain some form of spectral analytic geometry, or `complicial' analytic geometry in the sense of \cite{toen2008homotopical} Chapter 2.3. 

When $\underline{\mathbf{C}}$ is a connective Postnikov alebraic context whose underlying spectral algebraic context is a derived algebraic context, then in \cite{kelly2021analytic} Definition 7.5, the second and third author, together with Mukherjee, called $(\mathbf{Aff}_{\underline{\mathbf{C}}},\tau,\mathbf{P},\mathbf{S},\mathbf{F})$ a \textit{weak relative DAG context}. In this setting the following HKR theorem holds.

\begin{thm}[\cite{kelly2021analytic} Theorem 9.1]
    Let $(\mathbf{Aff}_{\underline{\mathbf{C}}},\tau,\mathbf{P},\mathbf{S},\mathbf{F})$ be w aweak relative DAG context enriched over $\mathbb{Q}$. Let $\mathcal{X}\in\mathbf{Sch}(\mathbf{A},\tau|_{\mathbf{A}},\mathbf{P}|_{\mathbf{A}})$. Then there is an equivalence, natural in $\mathcal{X}$
    $$\mathcal{L}(\mathcal{X})\cong T\mathcal{X}[-1]$$
    between the loop space and the shifted tangent space.
\end{thm}

\subsection{The Cotangent Complex and Obstruction Theories}\label{subsec:cotangent}

In this section we show that the abstract setup of AG contexts lends itself to a general theory of cotangent complexes and obstruction theories. Fix a weak relative \'{e}tale $(\infty,1)$-AG context $(\mathbf{Aff}_{\underline{\mathbf{C}}},\tau,\mathbf{P},\mathbf{S},\mathbf{F})$. In the context of algebraic contexts arising from exact categories most of the material in this section appeared in  \cite{kelly2021analytic} Section 8, and much more along these lines is accomplished in greater detail in \cite{RhiannonRepresentability}. In the setting of derived analytic geometry similar results appear in \cite{porta2017representability}. In all of these sources, the results are typically reasonably straightforward modifications of analogous results in \cite{toen2008homotopical}.

\subsubsection{The Cotangent Complex}

\begin{defn}[\cite{toen2008homotopical}, Definition 1.4.1.4]
Let $\mathcal{X}\in\mathbf{PreStk}(\mathbf{Aff}_{\underline{\mathbf{C}}})$, $Y\in\mathbf{Aff}^{cn}$, and $x:Y\rightarrow\mathcal{X}$ a map. Define 
$$\mathbf{Der}_{\mathcal{X}}(\mathrm{Spec}(A),-):\mathbf{QCoh}(\mathrm{Spec}(A))\rightarrow\textbf{sSet}$$
to be the functor
$$\mathrm{Map}_{{}_{Y\big\backslash}\mathbf{Aff}}(Y[M],\mathcal{X})$$
\end{defn}

We can also consider the restriction of $\mathbf{Der}_{\mathcal{X}}(\mathrm{Spec}(A),-)$ to ${}_{A}\mathbf{Mod}^{\mathbf{G}}_{0}$, and we denote this functor by $\mathbf{Der}_{\mathcal{X},0}(\mathrm{Spec}(A),-)$

\begin{defn}[\cite{toen2008homotopical} Section 1.4.1]
Let $f:\mathcal{X}\rightarrow\mathcal{Y}$ be a map in $\mathbf{PreStk}(\mathbf{Aff})$ and $x:U\rightarrow\mathcal{X}$ a map. By composition we get a map
$$\mathbf{Der}_{0,\mathcal{X}}(U,-)\rightarrow \mathbf{Der}_{0,\mathcal{Y}}(U,-).$$
We define 
$$\mathbf{Der}_{0,\mathcal{X}\big\slash\mathcal{Y}}(U,-):{}_{A}\mathbf{Mod}^{\mathrm{cn}}\rightarrow\textbf{sSet}$$
to be the fibre of the map $\mathbf{Der}_{0,\mathcal{X}}(U,-)\rightarrow \mathbf{Der}_{0,\mathcal{Y}}(U,-)$.
\end{defn}

A cotangent complex at $x:Y\rightarrow\mathcal{X}$ is, in principle, an object of ${}_{A}\mathbf{Mod}$ copresenting $\mathbf{Der}_{\mathcal{X}}(\mathrm{Spec}(A),-)$. However we only expect such an object to exist with restrictions on $M$ and $A$.

\begin{defn}[\cite{toen2008homotopical}, Definition 1.4.1.14]
Let $\mathrm{Spec}(A)=U\in\mathbf{S}^{op}\cap\mathbf{A}$. A map $f:\mathcal{X}\rightarrow\mathcal{Y}$ in $\mathbf{PreStk}(\mathbf{Aff})$ is said to \textit{have a cotangent complex at }$x:U\rightarrow\mathcal{X}$ \textit{relative to} $\mathbf{F}$ if there is a $\mathbf{G}$-connective $\mathbb{L}_{\mathcal{X}\big\slash\mathcal{Y},x}\in\mathbf{F}^{cn}_{U}$ and an equivalence 
$$\mathbf{Der}_{0,\mathcal{X}\big\slash\mathcal{Y}}(U,-)\cong\mathrm{Map}_{\mathbf{QCoh}(U)}(\mathbb{L}_{\mathcal{X}\big\slash\mathcal{Y},x},-)$$
in $\mathbf{Fun}({}_{A}\mathbf{Mod}^{\mathbf{G}}_{0},\mathbf{sSet})$.
\end{defn}

In the setting of geometry relative to a weak relative $(\infty,1)$-AG tuples it is possible to meaningfully define cotangent complexes of stacks, and tangent stacks. Here we follow \cite{toen2008homotopical} Section 1.4.1.


Let

\begin{displaymath}
\xymatrix{
\mathrm{Spec}(B)\ar[dr]^{y}\ar[rr]^{u} & & \mathrm{Spec}(A)\ar[dl]^{x}\\
& \mathcal{X} & 
}
\end{displaymath}
be a commutative diagram with $u$ a morphism in $\mathbf{S}^{op}\cap\mathbf{A}$. Suppose that $f:\mathcal{X}\rightarrow\mathcal{Y}$ has a relative cotangent complex at both $x$ and $y$. Using the universal properties of $\mathbb{L}_{\mathcal{X}\big\slash\mathcal{Y},x}$ and $\mathbb{L}_{\mathcal{X}\big\slash\mathcal{Y},y}$, there is a natural map.
$$\tilde{u}:\mathbb{L}_{\mathcal{X}\big\slash\mathcal{Y},y}\rightarrow u^{*}(\mathbb{L}_{\mathcal{X}\big\slash\mathcal{Y},x})$$

\begin{defn}[\cite{toen2008homotopical} Definition 1.4.1.15]
A morphism $f:\mathcal{X}\rightarrow\mathcal{Y}$ of prestacks is said to have a \textit{global cotangent complex relative to }$(\mathbf{A},\mathbf{F})$ if
\begin{enumerate}
\item
for any $U\in\mathbf{S}^{op}\cap\mathbf{A}$ and any $x:U\rightarrow\mathcal{X}$, $f$ has a relative cotangent complex at $x$ relative to $\mathbf{F}$
\item
For any morphism $u:U\rightarrow V$ in $\mathbf{S}^{op}\cap\mathbf{A}$, and any commutative diagram
\begin{displaymath}
\xymatrix{
U\ar[dr]^{y}\ar[rr]^{u} & &V\ar[dl]^{x}\\
& \mathcal{X} & 
}
\end{displaymath}
the map 
$$\tilde{u}:\mathbb{L}_{\mathcal{X}\big\slash\mathcal{Y},y}\rightarrow u^{*}(\mathbb{L}_{\mathcal{X}\big\slash\mathcal{Y},x})$$
is an equivalence. 
\end{enumerate}
$f$ is said to \textit{have a global contangent complex} if it has a global cotangent complex relative to $(\mathbf{Aff},\mathbf{Mod})$.
\end{defn}

\begin{rem}
If $f:\mathcal{X}\rightarrow\mathcal{Y}$ has a global contangent complex, then this data determines an object $\mathbb{L}_{\mathcal{X}\big\slash\mathcal{Y}}\in\mathbf{QCoh}(\mathcal{X})$.
\end{rem}

\begin{defn}[ \cite{toen2008homotopical} Definition 1.4.2.1]
\begin{enumerate}
\item
A pre-stack $\mathcal{X}\in\mathbf{PreStk}(\mathbf{Aff}_{\mathbf{C}})$ is said to be \textit{inf-cartesian relative to }$(\mathbf{A},\mathbf{F})$ if for any $A\in\mathbf{S}$ with $\mathrm{Spec}(A)\in\mathbf{A}$, any $1$-connective $M\in\mathbf{F}_{A}\cap{}_{A}\mathbf{Mod}^{\mathbf{G}}_{\ge1}$, and any derivation $d\in\pi_{0}\mathbf{Der}(A,M)$, the square 
\begin{displaymath}
\xymatrix{
\mathcal{X}(A\oplus_{d}\Omega M)\ar[d]\ar[r] & \mathcal{X}(A)\ar[d]^{d}\\
\mathcal{X}(A)\ar[r]^{s} & \mathcal{X}(A\oplus M)
}
\end{displaymath}
is cartesian in $\textbf{sSet}$.
\item
A pre-stack $\mathcal{X}$ \textit{has an obstruction theory relative to }$(\mathbf{A},\mathbf{F})$ if it has a global cotangent complex relative to $(\mathbf{A},\mathbf{F})$, and it is inf-cartesian relative to $(\mathbf{A},\mathbf{F})$.
\end{enumerate}
\end{defn}

\begin{defn}[ \cite{toen2008homotopical} Definition 1.4.2.2]
\begin{enumerate}
\item
A map of pre-stacks $f:\mathcal{X}\rightarrow\mathcal{Y}$ is said to be \textit{inf-cartesian relative to }$(\mathbf{A},\mathbf{F})$ if for any $A\in\mathbf{A}$ with $\mathrm{Spec}(A)\in\mathbf{S}$, any $1$-connective $M\in\mathbf{F}_{A}\cap{}_{A}\mathbf{Mod}^{\mathbf{G}}_{1}$, and any derivation $d\in\pi_{0}\mathbf{Der}(A,M)$, the square 
\begin{displaymath}
\xymatrix{
\mathcal{X}(A\oplus_{d}\Omega M)\ar[d]\ar[r] & \mathcal{Y}(A\oplus_{d}\Omega M)\ar[d]^{d}\\
\mathcal{X}(A)\times_{\mathcal{X}(A\oplus M)}\mathcal{X}(A)\ar[r]^{s} & \mathcal{Y}(A)\times_{\mathcal{Y}(A\oplus M)}\mathcal{Y}(A)
}
\end{displaymath}
is cartesian in $\textbf{sSet}$.
\item
A map of pre-stacks $f:\mathcal{X}\rightarrow\mathcal{Y}$ \textit{has an obstruction theory relative to }$(\mathbf{A},\mathbf{F})$ if it has a global cotangent complex relative to $(\mathbf{A},\mathbf{F})$, and it is inf-cartesian relative to $(\mathbf{A},\mathbf{F})$.
\end{enumerate}
\end{defn}

\begin{prop}[\cite{kelly2021analytic}, Proposition 8.12]
Let $f:\mathcal{X}\rightarrow\mathcal{Y}$ be a representable map in $\mathbf{Stk}(\mathbf{S},\tau|_{\mathbf{S}},\mathbf{P}|_{\mathbf{S}})$. Then $f$ has an obstruction theory relative to $(\mathbf{A},\mathbf{F})$.
\end{prop}

Again, exactly as in \cite{toen2008homotopical} Proposition 1.4.2.5 one can prove the following.

\begin{prop}
Let $f:\mathcal{F}\rightarrow \mathcal{G}$ be a map of stacks which has an obstruction theory relative to $(\mathbf{A},\mathbf{F})$. Let $X\in\mathbf{S}^{op}\cap\mathbf{A}$, $M\in\mathbf{F}_{X}\cap{}_{A}\mathbf{Mod}^{\mathbf{G}}_{1}$, $d:A\rightarrow M$ a derivation, and $h:X\rightarrow \mathcal{F}$ a map. Write $g=h\circ f$ Suppose we are given a commutative diagram
\begin{displaymath}
    \xymatrix{
X\ar[d]^{q}\ar[r]^{g} & \mathcal{F}\ar[d]^{f}\\
X\oplus_{d}\Omega M\ar[r]^{p} & \mathcal{G}
    }
\end{displaymath}
\begin{enumerate}
\item
There exists a natural obstruction
$$\alpha(x)\in\pi_{0}(\mathbf{Map}_{{}_{\mathbf{QCoh}(X)}}(f^{*}\mathbb{L}_{\mathcal{F}\big\slash\mathcal{G}},M))$$
which vanishes if and only if there exists a lift in the diagram $x$, i.e. a map $g_{d}:X\oplus_{d}\Omega M\rightarrow\mathcal{F}$ such that $q_{d}\circ q=g$, $f\circ g_{d}=p$.

\item
If $\alpha(x)=0$ then the space of lifts is non-canonically eqivalent to 
$$\Omega_{\alpha(x),0}\mathbf{Map}_{\mathbf{QCoh}(X)}(f^{*}\mathbb{L}_{\mathcal{F}\big\slash\mathcal{G}},M)$$
\end{enumerate}
\end{prop}

\subsubsection{Geometry Relative to Postnikov Contexts}

Let us now specialise to the case of either a Koszul algebraic context, or a connective Koszul algebraic context. We fix such a context

$$(\mathbf{C},\mathbf{C}_{\ge0},\mathbf{C}_{\le0},\mathbf{C}^{0},\mathbf{D},\theta,\mathbf{G}_{0},\mathbf{S})$$

\subsubsection{Lifting Maps}

The fact that in Koszul Postnikov algebraic contexts, the Postnikov towers of a $\mathbf{D}$-algebra is a sequence of square-zero extensions allows us to lift maps from the zero-truncation. For $A\in\mathbf{Alg_{D}}$, $n\ge0$, and $X=\mathrm{Spec}(A)$, we denote
$$t_{\le n}X\defeq\mathrm{Spec}(\tau_{\ge n}A)$$
The material in this such is a straightforward modification of the analogous results in \cite{porta2017representability} Section 5.

\begin{defn}
Let $\mathcal{F}$ be a stack and $M\in\mathbf{QCoh}(\mathcal{F})$. $M$ is said to be \textit{projective of tor-dimension} $n$ \textit{ relative to } $\mathbf{A}$ if for any $n\in\mathbb{N}_{0}$, any $X=\mathrm{Spec}(A)\in\mathbf{S}^{op}\cap\mathbf{A}$ such that $t_{\le n}X\rightarrow X$ is an equivalence, and any map $f:X\rightarrow F$, one has that 
$$\mathbf{Map}_{\mathbf{QCoh}(\mathcal{X})}(f^{*}M,N)\cong 0$$
whenever $N\in{}_{A}\mathbf{Mod}_{\ge n+1}$. 
\end{defn}

\begin{cor}
Let $\mathcal{F}\rightarrow\mathcal{G}$ be a map of stacks which has an obstruction theory relative to $(\mathbf{A},\mathbf{F})$ and $f:X\rightarrow \mathcal{F}$ a map with $X=\mathrm{Spec}(A)\in\mathbf{S}^{op}\cap\mathbf{A}$. Suppose that  $\mathbb{L}_{\mathcal{F}\big\slash\mathcal{G}}$ is projective of tor-dimension $0$ relative to $\mathbf{A}$.
Then for any $M\in\mathbf{F}_{X}\cap{}_{A}\mathbf{Mod}^{\mathbf{G}}_{1}$ and any derivation $d:A\rightarrow M$, any commutative diagram
\begin{displaymath}
    \xymatrix{
X\ar[d]^{q}\ar[r]^{g} & \mathcal{F}\ar[d]^{f}\\
X\oplus_{d}\Omega M\ar[r]^{p} & \mathcal{G}
    }
\end{displaymath}
admits a lift, which is unique up to a contractible choice.
\end{cor}


\begin{cor}
Let $f:\mathcal{F}\rightarrow\mathcal{G}$ be a map of stacks which has an obstruction theory relative to $(\mathbf{A},\mathbf{F})$ and $f:X=\mathrm{Spec}(A)\rightarrow \mathcal{G}$ a map with $X\in\mathbf{S}^{op}\cap\mathbf{A}$.  Suppose that 
\begin{enumerate}
\item
for $n\ge0 $ we have $\pi_{n}(A)\in\mathbf{F}_{A}\cap\mathbf{G}_{0}$. 
\item
$\mathbb{L}_{\mathcal{F}\big\slash\mathcal{G}}$ is projective and in tor-amplitude $0$ relative to $(\mathbf{A},\mathbf{F})$. 
\end{enumerate}
Then any lift of $t_{\le 0}(X)\rightarrow\mathcal{G}$ to a map $t_{\le0}(X)\rightarrow\mathcal{F}$ extends uniquely, up to a contractible choice, to a lift $X\rightarrow\mathcal{F}$ of $X\rightarrow\mathcal{G}$.
\end{cor}

%

\subsubsection{Higher (Underived) Geometry}

Let  \[(\mathbf{C},\mathbf{C}_{\ge0},\mathbf{C}_{\le0},\mathbf{C}^{0},\mathbf{D},\mathbf{\theta},\tau,\mathbf{P},\mathbf{A},\mathbf{F})\] be a (connective) weak relative smooth SAG context. Let $\mathbf{A}^{\heart}$ denote the full subcategory of $\mathbf{Aff_{C}}^{cn}$ consisting of objects of the form $\mathrm{Spec}(A)$ where $\mathrm{Spec}(A)\in\mathbf{A}$, and $A\in\mathbf{C}^{\heart}$.

\begin{defn}
$(\mathbf{C},\mathbf{C}_{\ge0},\mathbf{C}_{\le0},\mathbf{C}^{0},\mathbf{D},\mathbf{\theta},\tau,\mathbf{P},\mathbf{A},\mathbf{F})$ is said to be a \textit{derived enhacement of the classical theory }$\mathbf{A}^{\heart}$ if
\begin{enumerate}
\item
$(\mathbf{A},\tau|_{\mathbf{A}},\mathbf{P}|_{\mathbf{A},}\mathbf{A}^{\heart})$ is a strong relative $(\infty,1)$-geometry tuple. 
\item
for any $\mathrm{Spec}(A)\in\mathbf{A}$, $\pi_{0}(A)\in\mathbf{A}^{\heart}$
\end{enumerate}
\end{defn}
Then we get a fully faithful functor
$$\mathbf{Stk}(\mathbf{A}^{\heart},\tau|_{\mathbf{A}^{\heart}})\rightarrow\mathbf{Stk}(\mathbf{S},\tau|_{\mathbf{A}})$$
which induces a fully faithful functor between $n$-geometric stacks. $\mathbf{Stk}(\mathbf{S}^{\heart},\tau|_{\mathbf{S}^{\heart}})$ is the category of \textit{higher stacks} on the $1$-categorical Grothendieck site $(\mathbf{S}^{\heart},\tau|_{\mathbf{S}^{\heart}})$. 

Producing such enhancements usually comes in the following flavour. We start with a subcategory $\mathbf{A}^{\heart}\subseteq\mathrm{Comm}(\mathbf{C}^{\heart})$, a sub-presheaf $\mathbf{Q}^{\heart}\subseteq\mathbf{QCoh}^{\heart}|_{\mathbf{A}^{\heart}}$, and an $(\infty,1)$-pre-geometry tuple $(\mathbf{Aff}^{cn}_{\underline{\mathbf{C}}},\tau,\mathbf{P})$. We suppose that 
\begin{enumerate}
    \item
    $$(\mathbf{Aff}^{cn}_{\underline{\mathbf{C}}},\tau,\mathbf{P},\mathbf{A}^{\heart})$$
is a strong relative $(\infty,1)$-pre-geometry tuple
\item 
If $\mathrm{Spec}(B)\rightarrow\mathrm{Spec}(A)$ is a map in $\mathbf{P}$, and $\mathrm{Spec}(C)\rightarrow\mathrm{Spec}(A)$ is any map with $\mathrm{Spec}(C)\in\mathbf{A}^{\heart}$, then $\mathrm{Spec}(C\otimes_{A}^{\mathbb{L}}B)$ is in $\mathbf{A}^{\heart}$.
\item 
 $A\in\mathbf{Q}^{\heart}(\mathrm{Spec}(A))$ for any $\mathrm{Spec}(A)\in\mathbf{A}^{\heart}$. 
\item 
$\mathbf{P}$ consists of formally smooth maps,
\item
if $\mathrm{Spec}(B)\rightarrow\mathrm{Spec}(A)$ is a map in $\mathbf{P}|_{\mathbf{A}^{\heart}}$ and $M\in\mathbf{Q}^{\heart}(\mathrm{Spec}(A))$ then $B\otimes_{A}^{\mathbb{L}}M\cong B\otimes_{A} M$. 
\end{enumerate}

We then let $\mathbf{A}^{\mathbf{Q}}$ denote the class of objects $\mathrm{Spec}(A)$ in $\mathbf{Aff}^{cn}_{\underline{\mathbf{C}}}$ such that
\begin{enumerate}
    \item 
    $\mathrm{Spec}(\pi_{0}(A))\in\mathbf{A}^{\heart}$
    \item 
    $\pi_{n}(A)\in\mathbf{Q}^{\heart}(\mathrm{Spec}(\pi_{0}(A))$ for each $n$,
\end{enumerate}

Now let $\mathrm{Spec}(B)\rightarrow\mathrm{Spec}(A)$ be a map in $\mathbf{P}$ and let $\mathrm{Spec}(C)\rightarrow\mathrm{Spec}(A)$ be any map with $\mathrm{Spec}(A)\in\mathbf{A}^{\mathbf{Q}}$. We claim that $\mathrm{Spec}(B\otimes_{A}^{\mathbb{L}}C)$ is in $\mathbf{A}^{\mathbf{Q}}$. Indeed $B\otimes^{\mathbb{L}}_{A}\pi_{0}(C)$ is by assumption in $\mathbf{A}^{\heart}$, and the map $\pi_{0}(C)\rightarrow B\otimes^{\mathbb{L}}_{A}\pi_{0}(C)\cong\pi_{0}(B\otimes_{A}^{\mathbb{L}}C)$ is in $\mathbf{P}$. In particular both $C\rightarrow B\otimes_{A}^{\mathbb{L}}C$ and  $\pi_{0}(C)\rightarrow B\otimes^{\mathbb{L}}_{A}\pi_{0}(C)\cong\pi_{0}(B\otimes_{A}^{\mathbb{L}}C)$ are formally smooth. each $\pi_{n}(C)$ is transverse to the latter map. Thus 
$$C\rightarrow B\otimes_{A}^{\mathbb{L}}C$$
is derived strong. Hence 
$$\pi_{n}(B\otimes_{A}^{\mathbb{L}}C)\cong\pi_{0}(B\otimes_{A}^{\mathbb{L}}C)\otimes^{\mathbb{L}}_{\pi_{0}(C)}\pi_{n}(C)$$
Now $\pi_{0}(B\otimes_{A}^{\mathbb{L}}C),\pi_{0}(C)\in\mathbf{A}^{\heart}$ and $\pi_{n}(C)\in\mathbf{Q}^{\heart}(\mathrm{Spec}(C))$, so by assumption $\pi_{0}(B\otimes_{A}^{\mathbb{L}}C)\otimes^{\mathbb{L}}_{\pi_{0}(C)}\pi_{n}(C)\in\mathbf{Q}^{\heart}(\mathrm{Spec}(\pi_{0}(B\otimes_{A}^{\mathbb{L}}C))$. Thus
$$(\mathbf{Aff}^{cn}_{\underline{\mathbf{C}}},\tau,\mathbf{P},\mathbf{A}^{\mathbf{Q}})$$
is a strong relative $(\infty,1)$-pre-geometry tuple.
Finally for $\mathrm{Spec}(A)$ let $\mathbf{Q}_{+}(\mathrm{Spec}(A))\subseteq\mathbf{QCoh}(\mathrm{Spec}(A))$ denote the class of connective objects $F$ such that each $\pi_{n}(F)\in\mathbf{Q}^{\heart}(\mathrm{Spec}(\pi_{0}(A))$.

We would like some conditions under which it is a  strong relative $(\infty,1)$-geometry tuple, i.e. the restricition of the topology to $\mathbf{A}$ is sub-canonical.

\begin{defn}
    A \textit{Cartan context} is a tuple 
    $$(\mathbf{Aff}^{cn}_{\underline{\mathbf{C}}},\tau,\mathbf{P},\mathbf{A}^{\heart},\mathbf{Q}^{\heart})$$
    where
$\mathbf{A}^{\heart}\subseteq\mathrm{Comm}(\mathbf{C}^{\heart})$ is a subcategory, $\mathbf{Q}^{\heart}\subseteq\mathbf{QCoh}^{\heart}|_{\mathbf{A}^{\heart}}$ a sub-presheaf, and  $(\infty,1)$-pre-geometry tuple $(\mathbf{Aff}^{cn}_{\underline{\mathbf{C}}},\tau,\mathbf{P})$ and $(\infty,1)$-pre-geometry tuple.
\begin{enumerate}
    \item
    $$(\mathbf{Aff}^{cn}_{\underline{\mathbf{C}}},\tau,\mathbf{P},\mathbf{A}^{\heart})$$
is a strong relative $(\infty,1)$-pre-geometry tuple
\item 
$\mathbf{Q}^{\heart}$ is thick, and closed under finite limits and colimits.
\item 
If $\mathrm{Spec}(B)\rightarrow\mathrm{Spec}(A)$ is a map in $\mathbf{P}$, and $\mathrm{Spec}(C)\rightarrow\mathrm{Spec}(A)$ is any map with $\mathrm{Spec}(C)\in\mathbf{A}^{\heart}$, then $\mathrm{Spec}(C\otimes_{A}^{\mathbb{L}}B)$ is in $\mathbf{A}^{\heart}$.
\item 
 $A\in\mathbf{Q}^{\heart}(\mathrm{Spec}(A))$ for any $\mathrm{Spec}(A)\in\mathbf{A}^{\heart}$. 
\item 
$\mathbf{P}$ consists of formally smooth maps,
\item
objects of $\mathbf{Q}^{\heart}$ are transverse to maps in $\mathbf{P}|_{\mathbf{A}^{\heart}}$
\item 
$\mathbf{Q}_{+}|_{\mathbf{A}^{\mathbf{Q}}}$ is a sub-presheaf of $\mathbf{QCoh}$.
\item
$\mathbf{Q}^{\heart}|_{\mathbf{A}^{\heart}}$ satisfies \v{C}ech descent for $\tau|_{\mathbf{A}^{\heart}}$.
\end{enumerate}
A Cartan context 
$$(\mathbf{Aff}^{cn}_{\underline{\mathbf{C}}},\tau,\mathbf{P},\mathbf{A}^{\heart},\mathbf{Q}^{\heart})$$
is said to be \textit{an \'{e}tale Cartan context} if in addition
\begin{enumerate}
    \item 
    $\mathbf{P}$ consists of formally \'{e}tale maps.
    \item 
    for each $\mathrm{Spec}(A)\in\mathbf{A}^{\heart}$ there exists a $\mathrm{Spec}(B)$ in $\mathbf{A}^{\heart}$ such that
    \begin{enumerate}
        \item 
        $\mathbb{L}_{\mathrm{Spec}(B)}\cong B^{\oplus n}$ for some $n$.
        \item 
        there exists a map of algebras $B\rightarrow A$ which is a surjection in $\mathbf{C}^{\heart}$.
    \end{enumerate}
\end{enumerate}
\end{defn}

Related conditions to the above are explored in \cite{RhiannonRepresentability}, comprising \textit{representability contexts} in the terminology of loc. cit. Definition 5.3.2. Here they are used to establish representability theorems for derived stacks.

By Lemma \ref{lem:hyperdescentheart}, and Proposition \ref{prop:descentsubcan}, we immediately get the following.

\begin{prop}
    Let 
      $$(\mathbf{Aff}^{cn}_{\underline{\mathbf{C}}},\tau,\mathbf{P},\mathbf{A}^{\heart},\mathbf{Q}^{\heart})$$
      be a Cartan context. Then
     $$(\mathbf{Aff}^{cn}_{\underline{\mathbf{C}}},\tau,\mathbf{P},\mathbf{A}^{\heart})$$
     and
         $$(\mathbf{Aff}^{cn}_{\underline{\mathbf{C}}},\tau,\mathbf{P},\mathbf{A}^{\mathbf{Q}})$$
         are  strong relative hyper-$(\infty,1)$-geometry tuples. Moreover $\mathbf{Q}_{\ge n}|_{\mathbf{A^{\mathbf{Q}}}}$ and $\mathbf{Perf}_{\ge n}|_{\mathbf{A}^{\mathbf{Q}}}$ satisfy hyperdescent for $\tau|_{\mathbf{A^{\mathbf{Q}}}}$ for each $n\ge0$. If covers in $\tau|_{\mathbf{A^{\mathbf{Q}}}}$ have finite refinements, then $\mathbf{Q}_{+}$ and $\mathbf{Perf}|_{\mathbf{A}^{\mathbf{Q}}}$  also satisfy hyperdescent. 

         Moreover if 
               $$(\mathbf{Aff}^{cn}_{\underline{\mathbf{C}}},\tau,\mathbf{P},\mathbf{A},\mathbf{Q})$$
               is an \'{e}tale Cartan context then 
                          $$(\mathbf{Aff}^{cn}_{\underline{\mathbf{C}}},\tau,\mathbf{P},\mathbf{A},\mathbf{Q})$$
                          is a weak relative \'{e}tale $(\infty,1)$-AG context.
\end{prop}

\begin{proof}
The claim about strong relative hyper-$(\infty,1)$-geometry tuples follows immediately from  Lemma \ref{lem:hyperdescentheart} and Proposition \ref{prop:descentsubcan}.

For the \'{e}tale context claim first we observe that for $M\in{}_{A}\mathbf{Mod}_{\ge1}$ we have $\pi_{0}(A\oplus M)\cong \pi_{0}(A)$ and for $n>0$ $\pi_{n}(A\oplus M)\cong\pi_{n}(A)\oplus\pi_{n}(M)$ which is clearly in $\mathbf{Q}^{\heart}(\mathrm{Spec}(\pi_{0}(A)))$. Moreover for $\mathrm{Spec}(A)\in\mathbf{A}$ we have $\mathbb{L}_{A}\in\mathbf{Q}(\mathrm{Spec}(A))$ by Lemma \ref{lem:cotangentsubcat} .
\end{proof}

\begin{thm}\label{thm:ABclassfiyingcovers}
Let 
$$(\mathbf{Aff}_{\mathbf{C}}^{cn},\tau,\mathbf{P},\mathbf{A}^{\heart},\mathbf{Q}^{\heart})$$
be an \'{e}tale context. Suppose that all covers in $\tau|_{\mathbf{A}^{\heart}}$ have a $\kappa$-small refinement for some $\kappa\ge\aleph_{1}$. Let $\{\mathrm{Spec}(B_{i})\rightarrow\mathrm{Spec}(A)\}$ be a cover in $\tau|_{\mathbf{A}}$. Let $M\in\mathbf{QCoh}(\mathrm{Spec}(A))$ be such that each $\pi_{m}(M)$ is transverse to maps in $\mathbf{P}|_{\mathbf{A}}$, and each $\pi_{m}(A)$ is formally $\kappa$-filtered. Then the natural map
$$M\rightarrow\mathbf{lim}_{n\in\mathbb{N}}\prod_{(i_{1},\ldots,i_{n})\in\mathcal{I}^{n}}M\otimes_{A}^{\mathbb{L}}B_{i_{1}}\otimes_{A}^{\mathbb{L}}B_{i_{2}}\otimes_{A}^{\mathbb{L}}\ldots\otimes_{A}^{\mathbb{L}}B_{i_{n}}$$
    is an equivalence  
    
\end{thm}

\begin{proof}
    Note that each $A\rightarrow B_{i}$ is derived strong. Then 
    $$\{\mathrm{Spec}(\pi_{0}(B_{i}))\rightarrow\mathrm{Spec}(\pi_{0}(A))\}$$
    is a cover. Thus we have that the map 
    $$\pi_{0}(A)\rightarrow\mathbf{lim}_{n\in\mathbb{N}}\prod_{(i_{1},\ldots,i_{n})\in\mathcal{I}^{n}}\pi_{0}(B_{i_{1}})\otimes_{\pi_{0}(A)}^{\mathbb{L}}\pi_{0}(B_{i_{2}})\otimes_{\pi_{0}(A)}^{\mathbb{L}}\ldots\otimes_{\pi_{0}(A)}^{\mathbb{L}}\pi_{0}(B_{i_{n}})$$
    is an equivalence. Tensoring with $\pi_{m}(M)$ we get 
    \begin{align*}
       \pi_{m}(M)&\cong\rightarrow\mathbf{lim}_{n\in\mathbb{N}}\prod_{(i_{1},\ldots,i_{n})\in\mathcal{I}^{n}}\pi_{m}(M)\otimes_{\pi_{0}(A)}^{\mathbb{L}}\pi_{0}(B_{i_{1}})\otimes_{\pi_{0}(A)}^{\mathbb{L}}\pi_{0}(B_{i_{2}})\otimes_{\pi_{0}(A)}^{\mathbb{L}}\ldots\otimes_{\pi_{0}(A)}^{\mathbb{L}}\pi_{0}(B_{i_{n}})\\
       &\cong\mathbf{lim}_{n\in\mathbb{N}}\prod_{(i_{1},\ldots,i_{n})\in\mathcal{I}^{n}}\pi_{m}(M\otimes_{A}^{\mathbb{L}}B_{i_{1}}\otimes_{A}^{\mathbb{L}}B_{i_{2}}\otimes_{A}^{\mathbb{L}}\ldots\otimes_{A}^{\mathbb{L}}B_{i_{n}})
       \end{align*}
       using \cite{ben2020fr} Lemma 7.8. A spectral sequence argument then gives that 
       $$M\rightarrow\mathbf{lim}_{n\in\mathbb{N}}\prod_{(i_{1},\ldots,i_{n})\in\mathcal{I}^{n}}M\otimes_{A}^{\mathbb{L}}B_{i_{1}}\otimes_{A}^{\mathbb{L}}B_{i_{2}}\otimes_{A}^{\mathbb{L}}\ldots\otimes_{A}^{\mathbb{L}}B_{i_{n}}$$
    is an equivalence  
\end{proof}

Rhiannon Savage is working on a vast generalisation of the representability theorem of \cite{toen2008homotopical} which will show that a stack $\mathcal{X}\in\mathbf{Stk}(\mathbf{S},\tau|_{\mathbf{S}})$ is $n$-geometric precisely if the classical truncation $t_{\le 0}(\mathcal{X})$ is $n$-geometric and $\mathcal{X}$ has an obstruction theory.

\section{Geometry Relative to Lawvere Theories}\label{sec:geomlawv}

Let $(\mathbf{C},\mathbf{C}_{\ge0},\mathbf{C}_{\le0},\mathbf{C}^{0})$ be a derived algebraic context. In this section we discuss geometry relative to a Lawvere theory $\mathrm{T}$ which is of homotopy $\mathbf{C}^{\heart}$-polynomial type. The Lawvere theories allow us to pick out' different kinds of geometries inside of the universal' geometric context $\mathbf{Aff}_{\mathbf{C}_{\ge0}}$.

Let $\underline{\mathbf{C}}=(\mathbf{C},\mathbf{C}_{\ge0},\mathbf{C}_{\le0},\mathbf{C}_{0},\mathbb{L}\mathrm{Sym},\theta)$ be a derived algebraic context,  and $\mathrm{T}$ a $\Lambda$-sorted Lawvere theory of homotopy $\mathbf{C}^{\heart}$-polynomial type.

Let is introduce some classes of maps. Recall Subsubsection \ref{subsubsection:standardT} for some of the terminology below. 

\begin{defn}
Let $\tau$ be a pre-topology on $\mathbf{Aff}^{cn}_{\underline{\mathbf{C}}}$. A map $f:\mathrm{Spec}(B)\rightarrow\mathrm{Spec}(A)$ is said to be 
\begin{enumerate}
\item
$(\mathrm{T},\tau)$-smooth if there is a cover $\{\mathrm{Spec}(B_{i})\rightarrow\mathrm{Spec}(B)\}$ in $\tau$ such that $B_{i}\rightarrow A$ is standard $\mathrm{T}$-smooth. The class of all $(\mathrm{T},\tau)$- smooth maps is denoted $\textbf{sm}^{\mathrm{T},\tau}$.
\item
$(\mathrm{T},\tau,k)$-smooth if
for a $k\ge0$ if there is a cover $\{\mathrm{Spec}(B_{i})\rightarrow\mathrm{Spec}(B)\}$ in $\tau$ such that $B_{i}\rightarrow A$ is standard $\mathrm{T}$-smooth of relative dimension $k$. The class of all $(\mathrm{T},\tau)$- smooth maps is denoted $\textbf{sm}^{\mathrm{T},\tau,k}$.
\item
$\mathrm{T}$-\textit{\'{e}tale} if it is standard $\mathrm{T}$-\'{e}tale. We denote the class of $\mathrm{T}$-\'{e}tale maps by $\textbf{\'{e}t}^{\mathrm{T}}$.
\end{enumerate}
\end{defn}

Note that for \'{e}tale maps we assume they are always globally' standard. This is because this is the case in both algebraic and rigid analytic geometry. Recall also that $\mathbf{rat}^{\mathrm{T}}$ is the class of $\mathrm{T}$-rational localisations (definition \ref{defn:Trat}). 
\begin{defn}
We define 
$${}_{A}\mathbf{Mod}^{\aleph_{1},\mathrm{T},RR}$$
to be the class of $A$-modules such that
\begin{enumerate}
    \item 
    each $\pi_{n}(M)$ is $(\mathrm{T},RR)$-quasicoherent (Definition \ref{defn:Tqcoh}).
    \item 
    each $\pi_{n}(M)$ is $\aleph_{1}$-metrisable in ${}_{\pi_{0}(A)}\mathrm{Mod}$ (Definition \ref{defn:metrisable}).
\end{enumerate}
\end{defn}

\begin{defn}
\begin{enumerate}
\item
We define the \textit{countable} $\mathrm{G-T}$-\textit{ pre-topology} $\mathrm{G}_{\mathrm{T}}^{\aleph_{1}}$ to be the pre-topology consisting of collections of maps $\{\mathrm{Spec}(A_{i})\rightarrow\mathrm{Spec}(A)\}_{i\in I}$ such that
\begin{enumerate}
\item
each map $f_{i}:A\rightarrow A_{i}$ is a $\mathrm{T}$-rational localisation.
\item
there is a countable subset $J\subset I$ such that for any $A$-module $M\in{}_{A}\mathbf{Mod}^{\aleph_{1},\mathrm{T},RR}$ which satisfies $A_{j}\otimes^{\mathbb{L}}_{A}M\cong 0$ for all $j\in J$ then $M\cong 0$.
\end{enumerate}
\item 
We define the \textit{countable pre-}$\mathrm{G-T}$-\textit{pre-topology}, $G^{pre,\aleph_{1}}_{\mathrm{T}}$ to be the subpre-topology of $\mathrm{Op}_{\mathrm{T}}^{\aleph_{1}}$ generated by covers $\{\mathrm{Spec}(B_{i})\rightarrow\mathrm{Spec}(A)\}$ in $\mathrm{Op}_{\mathrm{T}}^{\aleph_{1}}$ such that each $A\rightarrow B_{i}$ is a $\mathrm{T}$-rational localisation.
\item 
We define the \textit{countable }$\mathrm{G-T}$-\textit{topology}, $G^{\aleph_{1}}_{\mathrm{T}}$ to be the topology generated by $G^{pre,\aleph_{1}}_{\mathrm{T}}$.
\end{enumerate}
\end{defn}

The reason for introducing the pre-$G$ pre-topology is that, as we shall see, it often has good properties with respect to maps in such covers being derived strong. It also behaves well with respect to base-change. It is not clear to us that the $G$-pre-topology does so.

\begin{defn}
\begin{enumerate}
\item
We define the \textit{finite} $\mathrm{G-T}$ \textit{pre-topology}, $G^{pre}_{\mathrm{T}}$ to be the pre-topology consisting of collections of maps $\{\mathrm{Spec}(A_{i}),\mathrm{Spec}(A)\}$ such that
\begin{enumerate}
\item
each map $f_{i}:A\rightarrow A_{i}$ is a $\mathrm{T}$-rational localisation.
\item
there is a finite subset $J\subset I$ such that whenever $M\in{}_{A}\mathbf{Mod}$ is such that $A_{j}\otimes^{\mathbb{L}}_{A}M\cong 0$ for all $j\in J$ then $M\cong 0$.
\end{enumerate}
\item 
We define the \textit{finite} $\mathrm{G-T}$-\textit{topology}, $G_{\mathrm{T}}$ to be the topology generated by $\mathrm{G_{T}}^{pre}$.
\end{enumerate}
\end{defn}

Finally we define our general classes of smooth maps and subdomain embeddings relative to $\mathrm{T}$. These are just smooth maps which locally look like some class of standard smooth maps after passing to an open cover.

\begin{defn}
    \begin{enumerate}
    \item 
    The class of $\mathrm{T}$-\textit{smooth maps}
 is
 $$\mathbf{sm}^{\mathrm{T}}\defeq\mathbf{sm}^{\mathrm{T}-std,\mathrm{G_{T}}}$$
 \item 
 The class of $\mathrm{T}$-\textit{open smooth maps} is 
 $$\mathbf{sm}_{o}^{\mathrm{T}}\defeq\mathbf{sm}_{o}^{\mathrm{T}-std,\mathrm{G_{T}}}$$
 \item 
 The class of $\mathrm{T}$-\textit{subdomain embeddings} is 
  $$\mathbf{sub}^{\mathrm{T}}\defeq(\mathbf{rat}^{\mathrm{T}})^{\mathrm{G_{T}}}$$
  \item 
    The class of \textit{Lindel\"{o}f} $\mathrm{T}$-\textit{smooth maps}
 is
 $$\mathbf{sm}^{\mathrm{T},\aleph_{1}}\defeq\mathbf{sm}^{\mathrm{T}-std,\mathrm{G^{\aleph_{1}}_{T}}}$$
   \item 
    The class of \textit{Lindel\"{o}f} $\mathrm{T}$-\textit{open smooth maps}
 is
 $$\mathbf{sm}_{o}^{\mathrm{T},\aleph_{1}}\defeq\mathbf{sm}_{o}^{\mathrm{T}-std,\mathrm{G^{\aleph_{1}}_{T}}}$$
  \item 
 The class of \textit{Lindel\"{o}f} $\mathrm{T}$-\textit{subdomain embeddings} is 
  $$\mathbf{sub}^{\mathrm{T},\aleph_{1}}\defeq(\mathbf{rat}^{\mathrm{T}})^{\mathrm{G^{\aleph_{1}}_{T}}}$$
 \end{enumerate}
\end{defn}

There is one more class of maps which we will need later - smooth maps which are local for the \'{e}tale topology
$$\mathbf{sm}^{\mathrm{T},\textrm{\'{e}t}}\defeq\mathbf{sm}^{\mathrm{T},\tau_{\mathrm{T}-\textrm{\'{e}t}}}$$

\begin{defn}
We define the \textit{smooth (resp. \'{e}tale pre-topology)}, denoted $\tau_{\mathrm{T}-sm}$ (resp. $\tau_{\mathrm{T}-\textrm{\'{e}t}}$) to be the pre-topology consisting of finite collections of maps $\{\mathrm{Spec}(B_{i})\rightarrow\mathrm{Spec}(A)\}$ such that 
\begin{enumerate}
\item
$A\rightarrow B_{i}$ is in $\mathbf{sm}^{\mathrm{T}}$ (resp. is in $\textbf{\'{e}t}^{\mathrm{T}})$.
\item
$A\rightarrow \prod_{i}B_{i}$ is descendable. 
\end{enumerate}
\end{defn}

\begin{rem}
    By Example \ref{Example:htpymonodescendable} covers in $G_{\mathrm{T}}$ are descendable. Thus $G_{\mathrm{T}}\subseteq\tau_{\mathrm{T}-\textrm{\'{e}t}}$. 
\end{rem}

By Corollary \ref{cor:Postnikovdecsendablecontext} we have the following.

\begin{cor}
$$(\mathbf{Aff}_{\mathbf{C}_{\ge0}},\tau_{\mathrm{T}-
\textrm{\'{e}t}},\mathbf{sm}^{\mathrm{T},\textrm{\'{e}t}})$$
$$(\mathbf{Aff}_{\mathbf{C}_{\ge0}},\tau_{\mathrm{T}-\textrm{sm}},\mathbf{sm}^{\mathrm{T}})$$
$$(\mathbf{Aff}_{\mathbf{C}_{\ge0}},\mathrm{G^{pre}_{T}},\mathbf{sm}^{\mathrm{T}})$$
$$(\mathbf{Aff}_{\mathbf{C}_{\ge0}},\mathrm{G_{T}},\mathbf{sm}^{\mathrm{T}})$$
    are $(\infty,1)$-geometry tuples.
\end{cor}

\subsection{Forcing Strong Maps}

We would also like geometry tuples relative to discrete $\mathrm{T}$-finitely presented algebras, coherent $\mathrm{T}$-algebras, and (almost) $\mathrm{T}$-finitely presented algebras. Let us introduce some notation. Here Recall Definitions \ref{defn:T-cohalg} and\ref{defn:Tfinemb}. Let 
\begin{enumerate}
\item
$\mathbf{Aff}^{cn,\mathrm{T}-coh}\defeq(\mathbf{DAlg}(\mathbf{C})^{cn,\mathrm{T}-coh})^{op}$
\item
$\mathbf{Aff}^{\heart,\mathrm{T}-fp}$ denote the category opposite to the category of discrete finitely $\mathrm{T}$-presented algebras in$\mathbf{C}^{\heart}$.
\item 
$\mathbf{Aff}^{cn,\mathrm{T}-afp}$ denote the category opposite to almost $\mathrm{T}$-finitely presented algebras.
\item 
$\mathbf{Aff}^{\mathrm{T}-fp}$ denote the category opposite to $\mathrm{T}$-finitely presented algebras.
\item 
$\mathbf{Aff}^{\heart,\mathrm{T}-f}$ denote the category opposite to discrete $\mathrm{T}$-finitely embeddable algebras.
\item 
$\mathbf{Aff}^{\heart,\mathrm{T}-gf}$ denote the category opposite to discrete $\mathrm{T}$-globally finitely embeddable algebras.
\end{enumerate}

We will need some `strong' versions of topologies relative to discrete objects. For this we in introduce some more notation. Recall the definition of $\mathbf{Sm}^{\mathrm{T}}$ from Definition \ref{defn:Topensmooth}.

\begin{defn}\label{defn:Tstrongtop}
    Given a pre-topology $\tau$, we define the \textit{corresponding strong pre-topology} $\overline{\tau}$ to consist of covers $\{\mathrm{Spec}(A_{i})\rightarrow\mathrm{Spec}(A)\}$ in $\tau$ such that each $A\rightarrow A_{i}$ is in $\mathbf{Sm}^{\mathrm{T}}$. We also define the \textit{strong} $\mathrm{T}$-\textit{open pre-topology}, $\overline{\tau}_{\mathrm{T}-open}$, to consist of finite covers $\{\mathrm{Spec}(A_{i})\rightarrow\mathrm{Spec}(A)\}_{i\in\mathcal{I}}$ such that 
    \begin{enumerate}
        \item 
        $A\rightarrow\prod_{i\in\mathcal{I}}A_{i}$ is descendable.
        \item 
        each $A\rightarrow A_{i}$ is in $\mathbf{open}^{\mathrm{T}}$ (Definition \ref{defn:Topensmooth}).
    \end{enumerate}
\end{defn}

\begin{defn}\label{defn:Tstrongmaps}

We define the classes of 

\begin{enumerate}
\item
\textit{strong} $\mathrm{T}$-\textit{ rational maps by} 
$$\overline{\mathbf{rat}}^{\mathrm{T}}\defeq\mathbf{rat}^{\mathrm{T}}\cap\mathbf{Sm}^{\mathrm{T}}$$

\item
\textit{strong} $\mathrm{T}$-\textit{subdomain embeddings} by $\overline{\mathbf{sub}}^{\mathrm{T}}\defeq(\mathbf{rat}^{\mathrm{T}})^{\overline{\mathrm{G}}^{pre}_{\mathrm{T}}}\cap\mathbf{Sm}^{\mathrm{T}}$
and \textit{Lindel\"{o}f strong} $\mathrm{T}$-\textit{subdomain embeddings} by $\overline{\mathbf{sub}}^{\mathrm{T}}\defeq(\mathbf{rat}^{\mathrm{T}})^{\overline{\mathrm{G}}^{pre,\aleph_{1}}_{\mathrm{T}}}\cap\mathbf{Sm}^{\mathrm{T}}$
    \item 
    \textit{strong} $\mathrm{T}$-\textit{ \'{e}tale maps by} $$\overline{\textbf{\'{e}t}}^{\mathrm{T}}\defeq\textbf{\'{e}t}^{\mathrm{T}}\cap\mathbf{Sm}^{\mathrm{T}}$$
    \item 
    \textit{strong} $\mathrm{T}$-\textit{ smooth maps by}$$\overline{\mathbf{sm}}^{\mathrm{T}}\defeq\mathbf{sm}^{\mathrm{T},\overline{\mathrm{G}}^{pre}_{\mathrm{T}}}\cap\mathbf{Sm}^{\mathrm{T}}$$
and \textit{Lindel\"{o}f strong smooth maps} by
$$\overline{\mathbf{sm}}^{\mathrm{T},\aleph_{1}}\defeq\mathbf{sm}^{\mathrm{T},\overline{\mathrm{G}}^{pre,\aleph_{1}}_{\mathrm{T}}}\cap\mathbf{Sm}^{\mathrm{T}}$$
    \item 
    \textit{strong} $\mathrm{T}$-\textit{open smooth maps by}$$\overline{\mathbf{sm}}_{o}^{\mathrm{T}}\defeq\mathbf{sm}_{o}^{\mathrm{T},\overline{\mathrm{G}}^{pre}_{\mathrm{T}}}\cap\mathbf{Sm}^{\mathrm{T}}$$
and \textit{open Lindel\"{o}f strong smooth maps} by
$$\overline{\mathbf{sm}}_{o}^{\mathrm{T},\aleph_{1}}\defeq\mathbf{sm}_{o}^{\mathrm{T},\overline{\mathrm{G}}^{pre,\aleph_{1}}_{\mathrm{T}}}\cap\mathbf{Sm}^{\mathrm{T}}$$
\end{enumerate}
\item 
$$\overline{\mathbf{sm}}^{\mathrm{T},\textrm{\'{e}t}}\defeq\mathbf{sm}^{\mathrm{T},\overline{\tau}_{\mathrm{T}-\textrm{\'{e}t}}}\cap\mathbf{Sm}^{\mathrm{T}}$$
\end{defn}

Let us spell out what $\overline{\mathbf{P}}^{\mathrm{T}}$ is, for $\mathbf{P}$ the various classes above. Let $A\rightarrow B$ be a map in $\overline{\mathbf{P}}^{\mathrm{T}}$, and suppose that $\pi_{0}(A)$ is $\mathrm{T}$-finitely embeddable. Then by assumption $\pi_{0}(B)$ is also $\mathrm{T}$-finitely embeddable. There is a cover $\{\mathrm{Spec}(B_{i})\rightarrow\mathrm{Spec}(B)\}$ in $\overline{\mathrm{G}}_{\mathrm{T}}^{pre}$ such that each $A\rightarrow B_{i}$ is in in $\mathbf{P}$. But again by assumption $\{\mathrm{Spec}(\pi_{0}(B_{i}))\rightarrow\mathrm{Spec}(\pi_{0}(B)\}$ is a cover in $\mathrm{G}_{\mathrm{T}}^{pre}$. Moreover each map $\mathrm{Spec}(\pi_{0}(B_{i}))\rightarrow\mathrm{Spec}(\pi_{0}(A))$ is also in $\mathbf{P}$. Thus $\mathrm{Spec}(B)\rightarrow\mathrm{Spec}(A)$ is also in $\overline{\mathbf{P}}^{\mathrm{T}}$. Under some transversality assumptions on the $\pi_{n}(A)$ and $\pi_{n}(A_{i})$, there is in fact an equivalence between maps in  $\overline{\mathbf{P}}^{\mathrm{T}}$  and derived strong maps such that the induced map on $\pi_{0}$ is in $\overline{\mathbf{P}}^{\mathrm{T}}$. Thus, in lieu of flatness, we are imposing that what happens at the derived level also happens on the level of $0$-truncations.

In particular under transversality assumptions on $\pi_{n}(A)$,

\begin{prop}\label{prop:etaleregstrongimp}
\begin{enumerate}
    \item 
Let $A\in\mathbf{DAlg}^{\heart,\mathrm{T}-f}(\mathbf{C})$ be a discrete $\mathrm{T}$-finitely embeddable algebra such that whenever $A\otimes\mathrm{F}(\mathrm{Free_{T}}(\underline{\lambda}))\rightarrow B$ is $\mathrm{T}$-\'{e}tale (resp. a $\mathrm{T}$-open immersion) then $B$ is also in $\mathbf{DAlg}^{\heart,\mathrm{T}-f}(\mathbf{C})$.
Then any standard  $\mathrm{T}$-smooth map  (resp. standard open $\mathrm{T}$ -smooth map) $A\rightarrow B$ of relative dimension $k$ (Definition \ref{defn:Tstandardsmooth}) is $\mathrm{T}$-$\mathbb{I}^{\oplus k}$-smooth (Definition \ref{defn:Topensmooth}). In particular if $f:A\rightarrow B$ is in $\textbf{\'{e}t}^{\mathrm{T}}$ then it is in $\overline{\textbf{\'{e}t}}^{\mathrm{T}}$.
\item 
Suppose that $\mathrm{T}$ is $\mathbf{C}$-localisation coherent (Definition \ref{defn:loccoherent}). Then $\overline{\mathbf{rat}}^{\mathrm{T}}=\mathbf{rat}^{\mathrm{T}}$ 
\end{enumerate}
\end{prop}

\begin{proof}
\begin{enumerate}
    \item 
    Write $A\cong \mathrm{F}(\mathrm{Free_{T}}(\underline{\lambda})\big\slash I)$. It suffices to prove that if $\mathrm{F}(\mathrm{Free_{T}}(\underline{\lambda}))\big\slash I\rightarrow B$ is standard $\mathrm{T}$-smooth, then $B$ is of the form $\mathrm{F}(\mathrm{Free_{T}}(\underline{\gamma}))\big\slash J$. Now $B\cong(\mathrm{F}(\mathrm{Free_{T}}(\underline{\lambda}))\big\slash I)\otimes^{\mathbb{L}}\mathrm{F}(\mathrm{Free_{T}}(\delta_{1},\ldots,\delta_{k+c}))\big\slash\big\slash (f_{1},\ldots,f_{c})$. By assumption the quotient is underived and this suffices to prove the claim.
    
    \item 
    This is similar to the first part.
\end{enumerate}
\end{proof}

   \begin{prop}\label{prop:cartanT}
       Let $\mathbf{A}\subseteq\mathbf{Aff}^{cn,\mathrm{T}-coh}$ be a subcategory such that
       \begin{enumerate}
       \item 
       if $\mathrm{Spec}(A)\in\mathbf{A}$ then $\pi_{0}(A)$ is coherent.
           \item 
                  $\mathrm{Spec}(A)\in\mathbf{A}$ if and only if $\mathrm{Spec}(\pi_{0}(A))\in\mathbf{A}$ and for $n\ge0$ $\pi_{n}(A)$ is a finitely presented $\pi_{0}(A)$-module.
                  \item 
                  if $\mathrm{Spec}(A)\in\mathbf{A}$ then $\mathrm{Spec}(A\otimes^{\mathbb{L}}\mathrm{F}(\mathrm{Free_{T}}(\lambda_{1},\ldots,\lambda_{n})))\in\mathbf{A}$.
                  \item 
                  if $\mathrm{Spec}(A)\in\mathbf{A}^{\heart}$ then  $A$ is $\mathrm{T}$-\'{e}tale regular and $\mathrm{T}$-\'{e}tale coherent - Definition \ref{defn:Tetaleregular} (resp is $\mathrm{T}$-localisation regular and $\mathrm{T}$-localisation coherent - Definition \ref{defn:loccoherent}).
                  \item 
                  $\mathrm{Coh}|_{\mathbf{A}^{\heart}}$ satisfies descent. 
       \end{enumerate}
Then
$$(\mathbf{Aff}_{\mathbf{C}}^{cn},\tau_{\mathrm{T}-\textrm{\'{e}t}},\overline{\mathbf{sm}}^{\mathrm{T}},\mathbf{A}^{\heart},\mathrm{Coh})$$
(resp.
$$(\mathbf{Aff}_{\mathbf{C}}^{cn},\mathrm{G_{T}}^{pre},\overline{\mathbf{sm}}_{o}^{\mathrm{T}},\mathbf{A}^{\heart},\mathrm{Coh})$$
)
is a Cartan context. Moreover
$$(\mathbf{Aff}_{\mathbf{C}}^{cn},\tau_{\mathrm{T}-\textrm{\'{e}t}},\textbf{\'{e}t}^{\mathrm{T}},\mathbf{A},\mathbf{Coh}_{+})$$
(resp.
$$(\mathbf{Aff}_{\mathbf{C}}^{cn},\mathrm{G_{T}}^{pre},\textbf{\'{e}t}^{\mathrm{T}},\mathbf{A},\mathbf{Coh}_{+})$$
)
is an \'{e}tale Cartan context.
   \end{prop}

\begin{proof}

We prove that \'{e}tale case, the open case being similar. First we observe that these are well-defined contexts, i.e. if $\{\mathrm{Spec}(B_{i})\rightarrow\mathrm{Spec}(B)\}$ is a cover in $\tau_{\mathrm{T}-\textrm{\'{e}t}}$ then each $A\rightarrow B_{i}$ is in $\overline{\mathbf{sm}}^{\mathrm{T}}$. Indeed let $A\rightarrow C$ be a map with $C\in\mathbf{A}^{\heart}$. By the assumption and Proposition \ref{prop:etaleregstrongimp}, $B\otimes_{A}^{\mathbb{L}}C$ is also discrete, and is an object of $\mathbf{A}^{\heart}$.

\begin{enumerate}
    \item 
    Now we verify the Cartan context conditions.

The above shows that 

$$(\mathbf{Aff}_{\mathbf{C}}^{cn},\tau_{\mathrm{T}-\textrm{\'{e}t}},\overline{\mathbf{sm}}^{\mathrm{T}},\mathbf{A}^{\heart},\mathrm{Coh})$$
is a strong relative $(\infty,1)$-geometry tuple.

 For $\mathrm{Spec}(A)\in\mathbf{A}^{\heart}$ we have that $\mathrm{Coh}_{+}(\mathrm{Spec}(A))$ is thick and closed under finite limits and colimits by Proposition \ref{prop:cohthick}.

Again by Proposition \ref{prop:etaleregstrongimp} $\mathrm{Spec}(B)\rightarrow\mathrm{Spec}(A)$ is a map in $\overline{\mathbf{sm}}^{\mathrm{T}}$ and $\mathrm{Spec}(C)\rightarrow\mathrm{Spec}(A)$ with $A\in\mathbf{A}^{\heart}$ then $\mathrm{Spec}(C\otimes_{A}^{\mathbb{L}}B)$ is in $\mathbf{A}^{\heart}$.

 By definition $\mathrm{Spec}(A)\in\mathrm{Coh}_{+}(\mathrm{Spec}(A))$ for any $\mathrm{Spec}(A)\in\mathbf{A}$.

 By definition $\overline{\mathbf{sm}}^{\mathrm{T}}$ consists of formally smooth maps.
 
Let $f:\mathrm{Spec}(B)\rightarrow\mathrm{Spec}(A)$ be a map in $\mathbf{A}^{\heart}$. We prove transversality of finitely presented modules to smooth maps. By $\mathrm{T}$-\'{e}tale coherence and descent for transversality we may assume that $f$ is of the form $A\rightarrow A\otimes\mathrm{F}(\mathrm{Free_{T}}(\lambda_{1},\ldots,\lambda_{n}))$ which is even a flat map.

$\mathbf{Coh}_{+}$ is a subpresheaf of $\mathbf{QCoh}|_{\mathbf{A}}$ by Corollary \ref{cor:sthngleftcoh}

$\mathrm{Coh}_{+}|_{\mathbf{A}^{\heart}}$ satisfies \v{C}ech descent by assumption.
\item
For the \'{e}tale restriction, first note that if $\mathrm{Spec}(A)\in\mathbf{A}$ and $M\in\mathbf{Coh}_{\ge1}(\mathrm{Spec}(A))$ then $\mathrm{Spec}(A\oplus M)\in\mathbf{A}$. Indeed $\pi_{0}(A\oplus M)\cong\pi_{0}(A)$ and for $n>0$ $\pi_{n}(A\oplus M)\cong\pi_{n}(A)\oplus\pi_{n}(M)$ is clearly coherent. All that remains to prove is that if $\mathrm{Spec}(A)$ is in $\mathbf{A}$ then $\mathbb{L}_{\mathrm{Spec}(A)}$ is in $\mathbf{Coh}_{\ge0}(\mathrm{Spec}(A))$. This is Proposition \ref{prop:altLcoh}.
\end{enumerate}
\end{proof}

\begin{example}
Let $R$ be a unital commutative ring ring and consider the derived algebraic context $\mathbf{C}_{R}$ (Subsubsection \ref{subsec:embeddingalgebra}).We let $\mathrm{T}=\mathbb{A}_{\mathbf{C}_{R}}$ be the $R$-polynomial Lawvere theory. Then $\mathrm{G}_{\mathbb{A}_{\mathbb{C}_{R}}}$ is the usual Zariski topology. Any smooth (resp. \'{e}tale) morphism of affine schemes is finitely presented, so is automatically in $\mathbf{Sm}^{\mathbb{A}_{\mathbb{C}_{R}}}$ (resp. $\textbf{\'{E}t}^{\mathrm{T}}$). Thus the classes $\textbf{sm}^{\mathbb{A}_{\mathbf{C}_{R}}}=\overline{\textbf{sm}}^{\mathbb{A}_{\mathbf{C}_{R}}}$ and $\textbf{\'{e}t}^{\mathbb{A}_{\mathbf{C}_{R}}}=\overline{\textbf{\'{e}t}}^{\mathbb{A}_{\mathbf{C}_{R}}}$ are just the usual notions of smooth/ \'{e}tale maps of $R$-modules. Moreover, covers in the classical \'{e}tale and smooth topologies are descendable, since they are finitely presented and faithfully flat. Thus $\tau_{\mathbb{A}_{\mathbf{C}_{R}}-\textrm{\'{e}t}}$ and $\tau_{\mathbb{A}_{\mathbf{C}_{R}}-\textrm{sm}}$ are the usual \'{e}tale and smooth topologies. Moreover $\overline{\mathrm{G}}_{\mathrm{T}}=\mathrm{G}_{\mathrm{T}}$ is again the usual Zariski pre-topology.

\end{example}

For strongly coherent Lawvere theories we get some further $(\infty,1)$-geometry tuples.

\begin{thm}\label{thm:varioustuples}
Suppose that $\mathrm{T}$ is strongly $\mathbf{C}^{\heart}$-coherent and concretely of $\mathbf{C}^{\heart}$-polynomial type. Then
\begin{enumerate}
    \item
 The following are all strong $(\infty,1)$-geometry tuples.
$$(\mathbf{Aff}^{cn}_{\underline{\mathbf{C}}},\overline{\mathbf{sm}}^{\mathrm{T}},\overline{\mathrm{G}}_{\mathrm{T}},\mathbf{Aff}^{cn,\mathrm{T}-coh})$$
$$(\mathbf{Aff}^{cn}_{\underline{\mathbf{C}}},\overline{\mathbf{sm}}^{\mathrm{T},\textrm{\'{e}t}},\overline{\tau}_{\mathrm{T}-\textrm{\'{e}t}},\mathbf{Aff}^{cn,\mathrm{T}-coh})$$
$$(\mathbf{Aff}^{cn}_{\underline{\mathbf{C}}},\overline{\mathbf{sm}}^{\mathrm{T}},\overline{\mathrm{G}}_{\mathrm{T}},\mathbf{Aff}^{cn,\mathrm{T}-afp})$$
$$(\mathbf{Aff}^{cn}_{\underline{\mathbf{C}}},\overline{\mathbf{sm}}^{\mathrm{T},\textrm{\'{e}t}},\overline{\tau}_{\mathrm{T}-\textrm{\'{e}t}},\mathbf{Aff}^{cn,\mathrm{T}-afp})$$
$$(\mathbf{Aff}^{cn}_{\underline{\mathbf{C}}},\overline{\mathbf{sm}}^{\mathrm{T}},\overline{\mathrm{G}}_{\mathrm{T}},,\mathbf{Aff}^{cn,\mathrm{T}-fp})$$
$$(\mathbf{Aff}^{cn}_{\underline{\mathbf{C}}},\overline{\mathbf{sm}}^{\mathrm{T},\textrm{\'{e}t}},\overline{\tau}_{\mathrm{T}-\textrm{\'{e}t}},\mathbf{Aff}^{cn,\mathrm{T}-fp})$$

\item 
 The following are all strong $(\infty,1)$-geometry tuples.
 $$(\mathbf{Aff}^{cn}_{\underline{\mathbf{C}}},\overline{\mathbf{sm}}^{\mathrm{T}},\overline{\mathrm{G}}_{\mathrm{T}},\mathbf{Aff}^{\heart,\mathrm{T}-fp})$$
$$(\mathbf{Aff}^{cn}_{\underline{\mathbf{C}}},\mathbf{Aff}^{\heart,\mathrm{T}-fp},\overline{\mathbf{sm}}^{\mathrm{T},\textrm{\'{e}t}}\overline{\tau}_{\mathrm{T}-\textrm{\'{e}t}},\mathbf{Aff}^{\heart,\mathrm{T}-fp})$$
\item 
If in addition maps in $\textbf{\'{e}t}^{\mathrm{T}}$ are transverse to finitely presented $A$-modules, where $A$ is a discrete finitely $\mathrm{T}$-presented algebra, then 
$$(\mathbf{Aff}^{cn}_{\underline{\mathbf{C}}},\mathbf{Aff}^{cn,coh},\overline{\mathbf{sm}}^{\mathrm{T}},\overline{\tau}_{\mathrm{T}-\textrm{sm}})$$
is also a strong $(\infty,1)$-geometry tuple. 
\end{enumerate}
All of the topologies above satisfy \v{C}ech descent.
\end{thm}

\begin{proof}
    \begin{enumerate}
        \item 
        By $\mathbf{Aff}^{cn,coh}$ Proposition \ref{prop:Tcohtens}, $\mathbf{Aff}^{cn,coh}$ is closed under fibre products.  Again by  Proposition \ref{prop:Tcohtens}, Proposition \ref{prop:Textend}, if $A\rightarrow B$ is $\mathrm{T}$-standard smooth with $A\in\mathbf{DAlg}^{cn,\mathrm{T}-coh}(\mathbf{C})$, then $B\in\mathbf{DAlg}^{cn,\mathrm{T}-coh}(\mathbf{C})$. All affines are \v{C}ech stacks by Proposition \ref{prop:descentsubcan} and the fact that we have assumed descendability. 
        The claims for finitely $\mathrm{T}$-presented and almost finitely $\mathrm{T}$-presented objects is similar, using Corollary \ref{cor:tensprodafp}.
          \item 
        This is similar to the previous part, noting that since maps in $\mathbf{sm}^{\mathrm{T}}$ are in $\mathbf{Sm}^{\mathrm{T}}$, and if $\mathrm{Spec}(B)\rightarrow\mathrm{Spec}(A)$ is a map in $\mathbf{Sm}^{\mathrm{T}}$ with $A$ a discrete finitely $\mathrm{T}$-presented algebra, then $B$ is also a disrete finitely $\mathrm{T}$-presented algebra.
        \item
        All that remains to prove is that the tuple is strong. Let $f:A\rightarrow B$ be a map in $\mathbf{sm}^{\mathrm{T}}$ with $A\in\mathbf{DAlg}^{cn,\mathrm{T}-coh}$. Since $f\in\mathbf{Sm}^{\mathrm{T}}$, $\pi_{0}(B)$ is also finitely $\mathrm{T}$-presented. Now by the factorisation of standard $\mathrm{T}$-smooth morphisms, by the assumption on transversality to \'{e}tale morphisms, and by  and by Lemma \ref{lem:decsenttransv}, $\pi_{n}(A)$ is transverse to $\pi_{0}(A)$ over $\pi_{0}(B)$. Also since $f\in\mathbf{Sm}^{\mathrm{T}}$, $\pi_{0}(A)\rightarrow\pi_{0}(B)$ is also formally smooth. Thus $f$ is derived strong, and the result follows easily.
      
    \end{enumerate}
    The claim about \v{C}ech descent simply follows from the fact that we are working with the descendable topologies.
\end{proof}

\subsection{$\mathrm{T}$-Base Change}

Let 
$$\mathbf{G}:(\mathbf{C},\mathbf{C}_{\ge0},\mathbf{C}_{\le0},\mathbf{C}^{0})\rightarrow(\mathbf{D},\mathbf{D}_{\ge0},\mathbf{D}_{\le0},\mathbf{D}^{0})$$
be a transformation of derived algebraic contexts. Let $\mathrm{T}$ be a $\Gamma$-sorted Lawvere theory of homotopy $\mathbf{C}^{\heart}$-polynomial type, with corresponding functor
$$\mathbf{F}:\mathbf{sAlg}_{\mathrm{T}}\rightarrow\mathbf{DAlg}^{cn}(\underline{\mathbf{C}})$$
Then 
$$\mathbf{G}\circ\mathbf{F}:\mathbf{sAlg}_{\mathrm{T}}\rightarrow\mathbf{DAlg}^{cn}(\underline{\mathbf{D}})$$
realises $\mathrm{T}$ as a Lawvere theory of homotopy $\mathbf{D}^{\heart}$-polynomial type. Moreover $\mathbf{G}$ is left adjoint and strongly monoidal, so the induced functor
$$\mathbf{G}:\mathbf{DAlg}(\underline{\mathbf{C}})\rightarrow\mathbf{DAlg}(\underline{\mathbf{D}})$$
is also a left adjoint. Since it is strongly monoidal it sends descendable maps to descendable maps, and since it is commutes with colimits it sends $\mathrm{T}$-standard \'{e}tale maps and $\mathrm{T}$-rational localisations to $\mathrm{T}$-standard \'{e}tale maps and $\mathrm{T}$-rational localisations respectively. Thus we get functors of $(\infty,1)$-geometry contexts

$$\mathbf{G}:(\mathbf{Aff}_{\mathbf{C}_{\ge0}},\tau_{\mathrm{T}-
\textrm{\'{e}t}},\mathbf{sm}^{\mathrm{T}})\rightarrow (\mathbf{Aff}_{\mathbf{D}_{\ge0}},\tau_{\mathrm{T}-
\textrm{\'{e}t}},\mathbf{sm}^{\mathrm{T}}) $$
$$\mathbf{G}:(\mathbf{Aff}_{\mathbf{C}_{\ge0}},\tau_{\mathrm{T}-\textrm{sm}},\mathbf{sm}^{\mathrm{T}})\rightarrow (\mathbf{Aff}_{\mathbf{C}_{\ge0}},\tau_{\mathrm{T}-\textrm{sm}},\mathbf{sm}^{\mathrm{T}})$$
$$\mathbf{G}(\mathbf{Aff}_{\mathbf{C}_{\ge0}},\mathrm{G^{pre}_{T}},\mathbf{sm}^{\mathrm{T}})\rightarrow (\mathbf{Aff}_{\mathbf{D}_{\ge0}},\mathrm{G^{pre}_{T}},\mathbf{sm}^{\mathrm{T}})$$
$$\mathbf{G}:(\mathbf{Aff}_{\mathbf{C}_{\ge0}},\mathrm{G_{T}},\mathbf{sm}^{\mathrm{T}})\rightarrow (\mathbf{Aff}_{\mathbf{D}_{\ge0}},\mathrm{G_{T}},\mathbf{sm}^{\mathrm{T}})$$

\subsection{$\mathrm{T}$-Analytification}\label{subsec:tan}

Let $(\mathbf{C},\mathbf{C}_{\ge0},\mathbf{C}_{\le0},\mathbf{C}^{0})$ be derived algebraic context, let $\mathrm{T}$ be a Lawvere theory, $F:\mathrm{T}^{op}\rightarrow\mathrm{C}$ a fully faithful functor, and  
$$p:\mathbb{A}_{\mathbf{C}^{\heart}}\rightarrow\mathrm{T}$$
be a map of Lawvere theories, realising $\mathrm{T}$ as a Lawvere theory of homotopy $\mathbf{C}^{\heart}$ polynomial type. There is a straightforward way to define analytification relative to $\mathrm{T}$. Denote by 
$$\mathrm{An}^{T}:\mathbf{DAlg}^{cn}(\mathbf{C}_{alg})\rightarrow\mathbf{DAlg}^{cn}(\mathbf{C})$$
the functor given by the composition
\begin{displaymath}
\xymatrix{
\mathbf{DAlg}^{cn}(\mathbf{C}_{R})\cong\mathbf{sAlg}_{\mathbb{A}^{\mathbf{C}}}\ar[r]^{p} &\mathbf{sAlg}_{\mathrm{T}}\ar[r]^{F} & \mathbf{DAlg}^{cn}(\mathbf{C})
}
\end{displaymath}

Let $R=\mathbf{Map}(\mathbb{I},\mathbb{I})$. We may regard $\mathbf{C}_{R}\subseteq\mathbf{C}$ as a full symmetric monoidal subcategory closed under colimits, so that we also have 
$$\mathbf{DAlg}^{cn}(\mathbf{C}_{R})\rightarrow\mathbf{DAlg}^{cn}(\mathbf{C})$$
closed under colimits.
There is then a natural map 
$$A\rightarrow\mathrm{An}^{T}(A)$$
By Lemma \ref{lem:htpyepiT} this map is a homotopy epimorphism. This proves the following.

\begin{prop}[following \cite{ben2021analytification}]\label{prop:ancotan}
Let $A\rightarrow B$ be a map in $\mathbf{DAlg}^{cn}(\mathbf{C}_{R})$. We have $$\mathbb{L}_{\mathrm{An^{T}}(B)\big\slash\mathrm{An^{T}}(A)}\cong \mathrm{An_{T}}(B)\otimes^{\mathbb{L}}_{\mathrm{An_{T}}(A)}\mathbb{L}_{B\big\slash A}$$    
\end{prop}

\begin{lem}
Let $R$ be a Banach ring and $f:A\rightarrow B$ be a Zariski open immersion/ \'{e}tale/ standard smooth morphism in $\mathbf{DAlg}^{cn}(\mathbf{C}_{R_{alg}})$. Then $\mathrm{An}^{\mathrm{T}}(f)$ is $\mathrm{T}$-Zariski open immersion/ $\mathrm{T}$-\'{e}tale morphism/ $\mathrm{T}$-standard smooth morphism.
\end{lem}

\begin{proof}
We have
$$\mathrm{An}^{\mathrm{T}}(A\otimes R[x]\big\slash\big\slash(1-xf))\cong\mathrm{An}^{\mathrm{T}}(A)\otimes\mathrm{F}(\mathrm{Free_{T}}(1))\big\slash\big\slash(1-xf)$$
$$\mathrm{An}^{\mathrm{T}}(A\otimes R[x_{1},\ldots,x_{k+c}]\big\slash\big\slash(f_{1},\ldots,f_{c}))\cong\mathrm{An}^{\mathrm{T}}(A)\otimes^{\mathbb{L}}\mathrm{F}(\mathrm{Free_{T}}(n))\big\slash\big\slash(f_{1},\ldots,f_{c})$$
Now if $k=0$ then the map $\mathrm{An^{T}}(A)\rightarrow\mathrm{F}(\mathrm{Free_{T}}(c))\big\slash\big\slash(f_{1},\ldots,f_{c})$ is formally \'{e}tale by Proposition \ref{prop:ancotan}. This proves the claim.
\end{proof}

In general it is not clear that analytification sends general smooth maps to $\mathrm{T}$-\'{e}tale / $\mathrm{T}$-smooth maps. 


We will consider analytification of covers, and therefore of stacks, in future work.

\section{General Constructions: Group Stacks and Mapping Stacks}

Let us fix a relative $(\infty,1)$-geometry tuple 
$$(\mathbf{M},\tau,\mathbf{S},\mathbf{P})$$ 
 In this section we look at general constructions such as group stacks and mapping stacks.


\subsection{Group Stacks}

As $\mathbf{Stk}(\mathbf{S},\tau|_{\mathbf{S}})$ is a symmetric monoidal category, we may consider the $(\infty,1)$-category of group objects internal to it. A \textit{group stack} is an object of this category. Moreover we say that a group stack $\mathcal{G}$ acts on a stack $\mathcal{X}$ if there is a structure of a left $\mathcal{G}$-module on $\mathcal{X}$. 

\begin{defn}
Let $\mathcal{G}$ be a group stack acting on a stack $\mathcal{X}$. We define the quotient prestacks $\Bigr[\mathcal{X}\big\slash\mathcal{G}\Bigr]^{pre}$ to be the colimit of the diagram
$$n\mapsto\mathcal{G}^{n}\times\mathcal{X}$$
\end{defn}

In particular we define $B^{pre}\mathcal{G}\defeq\Bigr[\bullet\big\slash\mathcal{G}]^{pre}$, where $\mathcal{G}$ acts trivially on the point.

\begin{defn}
For $\mathcal{G}$ a group stack acting on a stack $\mathcal{X}$ we define $\Bigr[\mathcal{X}\big\slash\mathcal{G}\Bigr]$ to be the stackification of $\Bigr[\mathcal{X}\big\slash\mathcal{G}\Bigr]^{pre}$. In particular $B\mathcal{G}$ is the stackification of $B^{pre}\mathcal{G}$. 
\end{defn}


  \begin{rem}\label{rem:quotprops}
Now let $\mathcal{G}$ be a group stack acting on a stack $\mathcal{X}$, and consider the quotient stack $\Bigr[\mathcal{X}\big\slash G\Bigr]$. We have
  $$\mathcal{X}\cong\Bigr[\mathcal{X}\big\slash G\Bigr]\times_{BG}\bullet$$
  Let $\mathbf{P}$ be a class of maps in $\mathbf{Stk}(\mathbf{S},\tau|_{\mathbf{S}})$ which is closed under fibre products and is \textit{local}, i.e. if $f:\mathcal{X}\rightarrow\mathcal{Y}$ is a map and $U\rightarrow\mathcal{Y}$ an epimorphisms of stacks with $U\in\mathbf{S}$ such that $U\times_{\mathcal{Y}}\mathcal{X}\rightarrow U$ is in $\mathbf{P}$, then $f$ is in $\mathbf{P}$. 
Suppose that $\mathcal{G}\rightarrow\bullet$ is in $\mathbf{P}$. $\bullet\rightarrow B\mathcal{G}$ is an epimorphism, and $\bullet\times_{B\mathcal{G}}\bullet\cong\mathcal{G}$. Thus $\bullet\rightarrow B\mathcal{G}$ is also in $\mathbf{P}$. Hence $\mathcal{X}\rightarrow\Bigr[\mathcal{X}\big\slash G\Bigr]$ also satisfies property $\mathbf{P}$.

    \end{rem}

$\pi_{0}(\Bigr[\mathcal{X}\big\slash\mathcal{G}\Bigr])(U)$ may be identified with the set of pairs $(P\rightarrow U,\rho)$ where $P\rightarrow U$ is a map in $\mathrm{Ho}(\mathbf{Stk}(\mathbf{S},\tau|_{\mathbf{S}}))$ which is a principal $\mathcal{G}$-bundle, and $\rho:P\rightarrow\mathcal{X}$ is a map in $\mathrm{Ho}({}_{\mathcal{G}}\mathrm{Mod}(\mathbf{Stk}(\mathbf{S},\tau|_{\mathbf{S}})))$. Indeed $\mathcal{X}\rightarrow\Bigr[\mathcal{X}\big\slash\mathcal{G}\Bigr]$ is an epimorphism. Thus given a map $f$ in $\pi_{0}\mathbf{Map}(U,\Bigr[\mathcal{X}\big\slash\mathcal{G}\Bigr])$ there is a cover $U_{i}\rightarrow U$ such that $U_{i}\rightarrow U\rightarrow \Bigr[\mathcal{X}\big\slash\mathcal{G}\Bigr]$ factors through  $\mathcal{X}\rightarrow \Bigr[\mathcal{X}\big\slash\mathcal{G}\Bigr]$. Then $U_{i}\times_{\Bigr[\mathcal{X}\big\slash\mathcal{G}\Bigr]}\mathcal{X}\cong U_{i}\times \mathcal{G}$. Thus $U\times_{\Bigr[\mathcal{X}\big\slash\mathcal{G}\Bigr]}\mathcal{X}\rightarrow U$ is a principal $\mathcal{G}$-bundle. Moreover it is equipped with a $\mathcal{G}$-equivariant map to $\mathcal{X}\rightarrow\Bigr[\mathcal{X}\big\slash\mathcal{G}\Bigr]$. Conversely given a principal $\mathcal{G}$-bundle $P\rightarrow U$ equipped with a $\mathcal{G}$-equivariant map $P\rightarrow\mathcal{X}$ we get a map $U\rightarrow\Bigr[\mathcal{X}\big\slash\mathcal{G}\Bigr]$ by passing to quotients. $\pi_{1}(\Bigr[\mathcal{X}\big\slash\mathcal{G}\Bigr])(U)$ can also be described in terms of equivalence classes of isomorphisms of principal bundles over $U$, with compatible maps to $\mathcal{X}$. 

\subsection{Moduli of Representations}

Let $\mathcal{X}$ be any $(\infty,1)$-topos and $G,H$ group objects in $\mathcal{X}$.

\begin{defn}
The \textit{moduli space of representations of }$G$\textit{ in }$H$ is the mapping stack
$$\mathbf{Map}(BG,BH)$$
\end{defn}

\begin{lem}
There is a natural equivalence
$$\mathbf{Map}(BG,BH)\cong[\mathbf{Hom}_{\mathrm{Grp}}(G,H)\big\slash H]$$
\end{lem}

\begin{proof}
The suspension-loop adjunction
$$\adj{\Sigma}{\mathrm{Grp}(\mathcal{X})}{\mathcal{X}_{\bullet}}{\Omega}$$
realises $\mathrm{Grp}(\mathcal{X})$ as a coreflective subcategory of the category $\mathcal{X}_{\bullet}$ of pointed objects of $\mathcal{X}$. Moreover this is a $\mathcal{X}$-enriched adjunction. Now for any $X,Y\in\mathcal{X}_{\bullet}$ we have a fibre product diagram
\begin{displaymath}
\xymatrix{
\underline{\mathbf{Map}}_{\mathcal{X}_{\bullet}}(X,Y)\ar[d]\ar[r] & \underline{\mathbf{Map}}_{\mathcal{X}}(X,Y)\ar[d]\\
\bullet\ar[r] & \mathbf{Map}_{\mathcal{X}}(\bullet,Y)\cong Y
}
\end{displaymath}
Taking $Y=BH$ we find that $\underline{\mathbf{Map}}_{\mathcal{X}_{\bullet}}(X,BH)$ is an $H$-torsor over $\underline{\mathbf{Map}}_{\mathcal{X}}(X,BH)$, and thus $\underline{\mathbf{Map}}_{\mathcal{X}}(X,BH)\cong[\underline{\mathbf{Map}}_{\mathcal{X}_{\bullet}}(X,BH)\big\slash H]$. If $X=BG$, then this gives $\mathbf{Map}(BG,BH)\cong[\mathbf{Hom}_{\mathrm{Grp}}(G,H)\big\slash H]$.
\end{proof}

\subsection{Quotient and Mapping Stacks Relative to Monoidal Categories}

In this section we specialise to the case that $\mathbf{M}=\mathbf{Aff_{C}}$ where $\mathbf{C}$ is a presentably symmetric monoidal $(\infty,1)$-category. Let $H\in\mathbf{S}$ be a Hopf algebra \cite{ergusthesis}. Then $G=\mathrm{Spec}(H)$ is a group object of $\mathbf{Stk}(\mathbf{S},\tau|_{\mathbf{S}})$. 

\begin{defn}
The category $\mathbf{Rep(G)}$ of $G$-representations is the category of $H$-comodules in $\mathbf{C}$. 
\end{defn}

\begin{prop}
Let $G=\mathrm{Spec}(H)$ where $H$ is a commutative Hopf algebra. Then there is an equivalence
$$\mathbf{QCoh}(B^{pre}\mathcal{G})\cong\mathbf{Rep}(G)$$
\end{prop}

\begin{proof}
We have
\begin{align*}
\mathbf{QCoh}(B^{pre}\mathcal{G})&\cong \underset{n}\lim \mathbf{Mod}(H^{\otimes n})
\end{align*}
Thus $\mathbf{QCoh}(B^{pre}\mathcal{G})$ is precisely the category of coalgebras for the comonad $H\otimes(-):\mathbf{C}\rightarrow\mathbf{C}$, i.e. an $H$-comodule.
\end{proof}
%
%

\subsubsection{Cotangent Complexes of Groups and Quotients}

For simplicity let us assume in this subsubsection that we are working in the connective part $\mathbf{C}_{\ge0}$ of a derived algebraic context $\mathbf{C}$ which is enriched over $\mathbb{Q}$.
Let $\mathcal{G}$ be a group stack. By functoriality, the cotangent space $\mathbb{L}_{G,e}$ at the identity of $\mathcal{G}$ will have the structure of a Lie coalgebra in $\mathbf{C}_{\ge0}$. Let us make this precise in the discrete case. Let $\epsilon:A\rightarrow R$ be an augmented algebra in $\mathbf{C}^{\heart}$ with $A$ cofibrant, and let $\tilde{A}\rightarrow A$ be a map of augmented algebras which is a homotopy epimorphism. Then by \cite{kelly2019koszul} Proposition 3.24, we have that
$$\mathbb{L}_{R\big\slash \tilde{A}}\cong\mathrm{Ker}(\epsilon)\big\slash\mathrm{Ker}(\epsilon)^{2}$$
Now let $\tilde{H},H$ be Hopf algebras and $\tilde{H}\rightarrow H$ a map of Hopf algebras which is a homotopy epimorphism of commutative algebras. We regard $H$ as an augmented algebra via the counit $\epsilon:H\rightarrow R$. In this case the coalgebra structure on $\tilde{H}$ induces a Lie coalgebra structure on $\mathbb{L}_{R\big\slash \tilde{H}}$, which we denote by $\mathrm{coLie}(H)$. We then define
$$\mathrm{Lie}(H)\defeq\mathrm{coLie}(H)^{\vee}$$

\begin{lem}
Let $G$ be a discrete group stack over $R$ which has a cotangent complex. Let $\pi:G\to \spec (R)$ be the projection and $e:\spec (R)\to G$ the unit. Then 
$\mathbb{L}_{G/\spec (R)}\cong \pi^*e^*\mathbb{L}_{G/\spec (R)}$.
\end{lem}

\begin{proof} Let $m,p_1,p_2:G\times G\to G$ denote the multiplication, first projection and second projection respectively. Let $\tau:G\times G\to G\times G$ be defined by $\tau=(m,p_2)$ (the universal translation). $\tau$ is an automorphism of $G\times G$ over $G$, where $G\times G$ is viewed as a $G$ scheme via the second projection. So $\tau^*\mathbb{L}_{G\times G/G}\to \mathbb{L}_{G\times G/G}$ is an isomorphism. As $G\times G\to G$ is the pullback of $G\to \spec (R)$ via $p_1$ we have that $\mathbb{L}_{G\times G/G}\cong p_1^*\mathbb{L}_{G/\spec (R)}$. We get that $m^*\mathbb{L}_{G/\spec (R)}\cong p_1^*\mathbb{L}_{G/\spec (R)}$. Pulling back this isomorphism via $(e\circ\pi, id_G):G\to G\times G$ gives us $\mathbb{L}_{G/\spec (R)}\cong \pi^*e^*\mathbb{L}_{G/\spec (R)}$ as required. Note that all the pullbacks in the proof are derived. 
\end{proof}

Suppose we are in a situation where we know that $m$-geometric stacks have cotangent complexes, and that such cotangent complexes which are themselves in $\mathbf{N}$ which satisfies hyperdescent - we will deal with situations in which this is always true in future work. In particular by \cite{kelly2021analytic} it holds if $\mathbf{P}$ consists of formally \'{e}tale maps. Then the cotangent complex can be computed as follows, as in \cite{raskin2009cotangent} Section 2. Fix such a stack $\mathcal{X}$, and let $x:U\rightarrow\mathcal{X}$ be a map with $U\in\mathbf{S}$, and let $\coprod V_{i}\rightarrow\mathcal{X}$ be an atlas.
\begin{displaymath}
\xymatrix{
\underset{i}\coprod U\times_{\mathcal{X}}V_{i}\ar[r]\ar[d] & U\ar[d]\\
\underset{i}\coprod V_{i}\ar[r] & \mathcal{X}
}
\end{displaymath}
The top horizontal map is a formal $\mathbf{P}$-cover, and thus satisfies descent for $\mathbf{N}$. Now $U\times_{\mathcal{X}}V_{i}$ is alsot geometric, so $\mathbb{L}_{\coprod_{i}U\times_{\mathcal{X}}V_{i}\big\slash U}$ exists. Moreover by general nonsense this is the pullback of $\mathbb{L}_{U\big\slash\mathcal{X}}$ to $V_{i}$.


Now consider a stack of the form $\Bigr[\mathcal{X}\big\slash \mathcal{G}\Bigr]$ where $\mathcal{X}$ and $\mathcal{G}$ are geometric, and $\mathcal{G}\times\mathcal{X}\rightarrow\mathcal{X}$ is a map in $\mathbf{P}$ defining an action of $\mathcal{G}$ on $\mathcal{X}$. Let $\pi:\mathcal{X}\rightarrow\Bigr[\mathcal{X}\big\slash \mathcal{G}\Bigr]$ be the canonical projection, so that we have a pullback square
\begin{displaymath}
\xymatrix{
\mathcal{X}\times\mathcal{G}\ar[d]\ar[r] & \mathcal{X}\ar[d]\\
\mathcal{X}\ar[r]&\Bigr[\mathcal{X}\big\slash\mathcal{G}\Bigr]
}
\end{displaymath}
By descent we have $\mathbb{L}_{\mathcal{X}\big\slash\Bigr[\mathcal{X}\big\slash\mathcal{G}\Bigr]}\cong\mathcal{O}_{\mathcal{X}}\otimes\mathfrak{g}^{\vee}$. Thus we have
$$\mathbb{L}_{\Bigr[\mathcal{X}\big\slash\mathcal{G}\Bigr],\pi}\cong\mathrm{cone}(\mathbb{L}_{\mathcal{X}}\rightarrow\mathcal{O}_{\mathcal{X}}\otimes\mathfrak{g}^{\vee})$$

In particular we have $\mathbb{L}_{B\mathcal{G}}\cong\mathfrak{g}^{\vee}[-1]$ where $\mathfrak{g}$ is the Lie algebra of $\mathcal{G}$.

\subsubsection{Mapping Stacks and Their Cotangent Complexes}
The category $\mathbf{Stk}(\mathbf{S},\tau|_{\mathbf{S}})$ is a closed symmetric monoidal category when equipped with the Cartesian monoidal category. The internal hom $\underline{\mathbf{Map}}$ may be defined as 
$$\underline{\mathbf{Map}}(\mathcal{X},\mathcal{Y})(U)=\mathbf{Map}(\mathcal{X}\times U,\mathcal{Y})$$
where $\mathbf{Map}$ is the simplicial mapping space functor. The following fact will be useful when dealing with moduli of representations of groups.

Suppose now that $\mathbf{C}$ is a flat spectral algebraic pre-context and $\mathbf{M}=\mathbf{Aff}^{cn}_{\underline{\mathbf{C}}}$. Let $\mathcal{X}$ and $\mathcal{Y}$ be stacks, and consider the mapping stack
$\underline{\mathbf{Map}}(\mathcal{X},\mathcal{Y})$.

\begin{lem}[c.f. \cite{porta2018derivedhom} Lemma 8.4]\label{lem:cotangentmapping}
Let $f:\mathcal{X}\rightarrow S$ be a map of stacks with $S\in\mathbf{S}$, and let $Y\rightarrow S$ be a map of stacks whose global relative contangent complex exists and is perfect. Suppose that $\mathcal{X}\rightarrow S$ is strongly $\mathbf{QCoh}$-cohomologically proper (Definition \ref{defn:cohprop}) and $\mathbf{QCoh}$-perfect base change (Definition \ref{defn:perfbase}). Then
$$\underline{\mathbf{Map}}_{S}(\mathcal{X},\mathcal{Y})\rightarrow S$$
 has a connective cotangent complex.
\end{lem}

\begin{proof}
Consider the diagram
\begin{displaymath}
\xymatrix{
\mathcal{X}\times_{S}\underline{\mathbf{Map}}_{S}(\mathcal{X},\mathcal{Y})\ar[d]^{\pi}\ar[rr]^{ev} & & \mathcal{Y}\\
\underline{\mathbf{Map}}_{S}(\mathcal{X},\mathcal{Y})
}
\end{displaymath}
Recall the plus construction from Subsubsection \ref{subsubsec:plus}, and define $L\defeq\pi_{+}(ev^{*}(\mathbb{L}_{\mathcal{Y}\big\slash S}))$. Let $u:U\rightarrow\underline{\mathbf{Map}}_{S}(\mathcal{X},\mathcal{Y})$ be a map with $U$ affine. This is equivalent to the data of a map $f_{u}:\mathcal{X}\times_{S}U\rightarrow\mathcal{Y}$. Moreover we have a pullback diagram
\begin{displaymath}
\xymatrix{
\mathcal{X}\times_{S}U\ar[d]^{\pi_{u}}\ar[r]^{q} & \mathcal{X}\times_{S}\underline{\mathbf{Map}}(\mathcal{X},\mathcal{Y})\ar[d]^{\pi}\\
U\ar[r]^{u} & \underline{\mathbf{Map}}(\mathcal{X},\mathcal{Y})
}
\end{displaymath}
Let $F\in\mathbf{QCoh}(U)$. Then we have
\begin{equation}
\begin{split}
\mathbf{Der}_{\mathbf{Map}_{S}(\mathcal{X},\mathcal{Y})\big\slash S}(U,F)&\cong\mathbf{Map}_{\mathbf{QCoh}(\mathcal{X}\times_{S}U)}(f_{u}^{*}\mathbb{L}_{\mathcal{Y}\big\slash S}.\pi_{u}^{*}(F)) \\ & \cong\mathbf{Map}((\pi_{u})_{+}f^{*}_{u}\mathbb{L}_{\mathcal{Y}\big\slash S},F)
\end{split}
\end{equation}
Base change for the $+$-construction (Lemma \ref{lem:Perfbasechange}) then gives
$$(\pi_{u})_{+}f_{u}^{*}\mathbb{L}_{\mathcal{Y}\big\slash S}\cong u^{*}(\pi_{+}ev^{*}(\mathbb{L}_{\mathcal{Y}\big\slash S}))$$
By the assumption this is connective, and therefore defines the cotangent complex of the map 
$$\underline{\mathbf{Map}}_{S}(\mathcal{X},\mathcal{Y})\rightarrow S$$
\end{proof}

Let $\mathcal{X}$ be a stack such that $\mathcal{X}\rightarrow\bullet$ is a base change morphism and is strongly cohomologically proper. Let $\mathcal{G}$ be a geometric group stack such that $\mathcal{G}\rightarrow\bullet$ is in $n-\mathbf{P}$ for some $n$, and such that $\mathfrak{g}^{\vee}=\mathrm{coLie}(\mathcal{G})$ is reflexive and nuclear relative to $\mathbf{S}$. Consider the moduli space of principal $\mathcal{G}$-bundles on $\mathcal{X}$,
$$\mathrm{Bun}_{\mathcal{G}}(\mathcal{X})\defeq\underline{\mathbf{Map}}(\mathcal{X},B\mathcal{G})$$
Let $U\rightarrow\mathrm{Bun}_{\mathcal{G}}(\mathcal{X})$ be a map, classified by a principal $G$-bundle $P$ on $\mathcal{X}\times U$, or equivalently a map $f_{U}:\mathcal{X}\times U\rightarrow B\mathcal{G}$. The pullback of $\mathbb{L}_{B\mathcal{G}}\cong\mathfrak{g}^{\vee}[-1]$ is the $\mathcal{G}$-bundle $P\times_{\mathcal{G}}\mathfrak{g}^{\vee}[-1]$. This is a vector bundle which we may equivalently regard as a (locally free) quasi-coherent sheaf. Then \[\mathbb{L}_{\mathrm{Bun}_{\mathcal{G}}(\mathcal{X}),U}\cong(\pi_{U})_{+}(P\times_{\mathcal{G}}\mathfrak{g}^{\vee}[-1]).\] In particular if $U=\mathrm{Spec}(R)$, then \[\mathbb{L}_{\mathrm{Bun}_{\mathcal{G}}(\mathcal{X}),\mathrm{Spec}(R)}\cong\mathbb{R}\Gamma(\mathcal{X},P\times_{\mathcal{G}}\mathfrak{g}^{\vee}[-1]).\] 

Let us interpret this in the context of the moduli of representations of $G$ in $G'$, where $G$ and $G'$ are homotopicaly discrete group schemes, $G$ is strongly $\mathbf{QCoh}$-proper, and $\mathrm{coLie}(G')$ is free of finite rank. Usually $G'$ will be an open inside $\mathrm{GL_{n}}$.
Consider
$$\underline{\mathbf{Map}}(BG,BG')$$
Let $\pi:\mathrm{Spec}(A)\rightarrow\underline{\mathbf{Map}}(BG,BG')$. Perhaps after passing to a smaller neighbourhood, we may assume that $\pi$ is represented by a map $f:\mathrm{Spec}(A)\rightarrow\underline{\mathbf{Map}}_{\mathrm{Grp}}(G,G')$. The pullback of \[\mathbb{L}_{BG'}\cong\mathfrak{g'}^{\vee}[-1]\cong\mathrm{coLie}(G')[-1]\] is $A\otimes\mathrm{coLie}(G)[-1]$. We may regard this as a representation of the group object $\mathrm{Spec}(A)\times G\rightarrow\mathrm{Spec}(A)$ over $\mathrm{Spec}(A)$, with the action being induced from the restriction of the $G'$ action on $\mathrm{coLie}(G')[-1]$ along the map $\mathrm{Spec}(A)\times G\rightarrow G'$. The dual of this representation is $A\otimes\mathrm{Lie}(G')[1]$. We write this representation as $\mathrm{Ad}_{G'}(\mathrm{Spec}(A)\times G)$. For a representation $M\in{}_{A}\mathbf{Mod}$ of $\mathrm{Spec}(A)\times G$ denote by $\mathbb{R}_{A}(M)^{G}$ the group cohomology of $M$. The functor $\mathbb{R}_{A}(-)^{G}$ coincides with the pushofrward functor $\mathrm{Spec}(A)\times BG\rightarrow\mathrm{Spec}(A)$. By the proof of Lemma \ref{lem:cotangentmapping}, the cotangent complex of $\underline{\mathbf{Map}}(BG,BG')$ at $A$ is 
$$\mathbb{R}_{A}(\mathrm{Ad}_{G'}(\mathrm{Spec}(A)\times G))^{\vee}[-1].$$

As mentioned above most applications, $G'$ will be an `open' inside $\mathrm{GL_{n}}^{alg}\defeq\mathrm{Spec}(\mathrm{Sym}(\mathbb{I}^{\oplus n^{2}}[\mathrm{det}^{-1}])$, i.e. there will be a homotophy monomorphism $G\rightarrow\mathrm{GL_{n}}^{alg}$ - in practice this will be some $\mathrm{T}$-analytification of  $\mathrm{GL_{n}}^{alg}$.

\chapter{Bornological Analytic Geometry}\label{BAG}

In this chapter we apply the very general machinery we have developed up until now in order to give a model for analytic geometry, both Archimedean and non-Archimedean. We also explain how, working over Banach rings of integers, we can obtain a `universal analytic geometry'.

\section{Preliminary Observations}

In this section we collate some basic results connecting analytic geometry and bornological algebra.

\subsection{The Spectrum}\label{subsec:Spectr}

In Berkovich analytic geometry, the basic object of study is the \textit{Berkovich spectrum}. Such an object can be defined for bornological algebras as well.

\begin{defn}[\cite{bambozzi2015stein} Definition 4.1]
Let $A\in\mathrm{Comm}(\mathrm{CBorn}_{\mathbb{Z}_{an}}\mathrm{)}$, and let $K$ and $L$ be valued fields. We say that maps $k:A\rightarrow K$ and $l:A\rightarrow L$ in $\mathrm{Comm}(\mathrm{CBor}_{\mathbb{Z}_{an}}\mathrm{)}$ are \textit{equivalent} if there is a map $A\rightarrow M$ in $\mathrm{Comm}(\mathrm{CBor}_{\mathbb{Z}_{an}}\mathrm{)}$ with $M$ a valued field and maps $M\rightarrow K$, $M\rightarrow L$ of Banach fields, such that the compositions $A\rightarrow M\rightarrow K$ and $A\rightarrow M\rightarrow L$ are equal to $k$ and $l$ respectively. We denote by $\mathcal{M}(A)$ the set of all such equivalence classes of maps from $A$ to a valued field, and call it the \textit{spectrum of }$A$. It is equipped with the weak topology, i.e. the weakest topology such that for each $a\in A$, the `evalution map' $|ev_{f}|:A\rightarrow\mathbb{R}_{\ge 0}$, $\chi\mapsto |\chi(f)|$ is continuous. 
\end{defn}

\begin{rem}
An important subspace of $\mathcal{M}(A)$ is the \textit{maximal spectrum} $\mathrm{Max}(A)$.
Let $A\in\mathrm{Comm}(\mathrm{CBorn}_{\mathbb{Z}_{an}}\mathrm{)})$. The \textit{maximal spectrum} of $A$, denoted $\mathrm{Max}(A)$, is the set of admissible epimorphisms $A\rightarrow K$ where $K$ is a valued field, with finitely generated kernel. 
\end{rem}

\begin{prop}
Let $A_{1},\ldots,A_{n}$ be a finite collection of objects of $\mathrm{Comm}(\mathrm{CBorn}_{R})$. Then there is a natural bijection $\coprod_{i=1}^{n}\mathcal{M}(A_{i})\cong\mathcal{M}(A_{1}\times\cdots \times A_{n})$
\end{prop}

\begin{proof}
The map $\coprod_{i=1}^{n}\mathcal{M}(A_{i})\cong\mathcal{M}(A_{1}\times\cdots \times A_{n})$ sends the equivalence class of $A_{j}\rightarrow K$ to the equivalence class of the composition $A_{1}\times\cdots \times A_{n}\rightarrow A_{j}\rightarrow K$. This map is clearly injective. Let us prove that it is surjective. Let $A_{1}\times\cdots \times A_{n}\rightarrow K$ be a map to a field. The kernel of this map is an admissible prime ideal of $A_{1}\times\cdots \times A_{n}$. It must have the form $A_{1}\times\cdots \times A_{j-1}\times\mathfrak{p}\times A_{j+1}\times\cdots \times A_{n}$ where $\mathfrak{p}$ is a prime ideal of $A_{j}$ for some $j$. Thus the map $A_{1}\times\cdots \times A_{n}\rightarrow K$ factors through the projection $A_{1}\times\cdots \times A_{n}\rightarrow A_{j}$. 
\end{proof}

Let us establish some important instances when $\mathcal{M}(A)$ is non-empty.

 Let $R$ be a non-Archimedean Banach ring. We say that $R$ is \textit{strongly Noetherian} if for any $\lambda_{1},\ldots,\lambda_{n}$, $R<\lambda_{1},\ldots,\lambda_{n}>$ is Noetherian. We say that $R$ is a \textit{Tate ring} if it has a topologically nilpotent unit (see for example \cite{bambozzi2020sheafyness} Section 4).

\begin{prop}\label{prop:afndnonemptyspec}
Let $R$ be a strongly Noetherian Tate ring (e.g. a non-trivially valued field), and let $I\subset R\Bigr<\lambda_{1}^{-1}x_{1},\ldots,\lambda_{n}^{-1}x_{n}\Bigr>$ where $I$ is a proper ideal. Then there is a map $ R\Bigr<\lambda_{1}^{-1}x_{1},\ldots,\lambda_{n}^{-1}x_{n}\Bigr>\rightarrow L$ which vanishes on $I$, where $L$ is a Banach field.
\end{prop}

\begin{proof}
 All ideals of $R\Bigr<\rho_{1}^{-1}x_{1},\ldots,\rho_{n}^{-1}x_{n}\Bigr>$ are closed. Pick a maximal ideal of $R\Bigr<\rho_{1}^{-1}x_{1},\ldots,\rho_{n}^{-1}x_{n}\Bigr>$ containing $I$, and set $L=R\Bigr<\rho_{1}^{-1}x_{1},\ldots,\rho_{n}^{-1}x_{n}\Bigr>\big\slash I$.
\end{proof}

\begin{cor}\label{lem:nonzeromaptofield}
Let $k$ be a non-trivially valued non-Archimedean Banach field, and let $A$ be a non-zero $\mathrm{Afnd}_{k}$-finitely presented $k$-algebra. Then there is a non-trivially valued Banach field $F$, and a map $A\rightarrow F$.
\end{cor}

\begin{prop}\label{nonzeroaffinoid}
Let $k$ be a non-trivially valued non-Archimedean Banach field, and let $A$ be a non-zero $\mathrm{Afnd}^{\dagger}_{k}$-finitely presented $k$-algebra, or a non-zero $\mathrm{Afnd}^{nA,\dagger}_{k}$-finitely presented $k$-algebra in the non-Archimedean case. Then there is a non-trivially valued Banach field $F$, and a map $A\rightarrow F$.
\end{prop}

\begin{proof}
This is an immediate consequence of \cite{MR3448274} Lemma 4.3.
\end{proof}

\begin{prop}\label{prop:daggernonemptyspec}
Let $F$ be a  non-trivially valued Banach field, and let $I\subset \mathcal{O}(D^{n}_{<\lambda,F})$ where $I=(f_{1},\ldots, f_{m})$ is a closed, finitely generated ideal. Then there is a map $ \mathcal{O}(D^{n}_{<\lambda,F})\rightarrow L$ which vanishes on $I$, where $L$ is a Banach field.
\end{prop}

\begin{proof}
We write $ \mathcal{O}(D^{n}_{<\lambda,F})\cong\lim_{\underline{\rho}<\underline{\lambda}} F\Bigr<\rho_{1}^{-1}x_{1},\ldots,\rho_{n}^{-1}x_{n}\Bigr>^{\dagger}$. By \cite{MR3448274} all ideals of $F\Bigr<\rho_{1}^{-1}x_{1},\ldots,\rho_{n}^{-1}x_{n}\Bigr>^{\dagger}$ are closed. There is some $\underline{\rho}$ such that the image of $I$ in  $F\Bigr<\rho_{1}^{-1}x_{1},\ldots,\rho_{n}^{-1}x_{n}\Bigr>^{\dagger}$ is proper. Now \[F\Bigr<\rho_{1}^{-1}x_{1},\ldots,\rho_{n}^{-1}x_{n}\Bigr>^{\dagger}\big\slash(f_{1},\ldots, f_{m})\] is a non-zero object of $\mathrm{Ind}(\mathrm{Comm}(\mathrm{Ban}_{k}))$ so by \cite{MR3448274} its spectrum is non-empty.
\end{proof}
%
%
%

\subsection{Conservativity and Surjectivity}

Here we relate joint surjectivity of collections of maps of Berkovich spectra of bornological algebras to conservativity of base change functors between categories of modules.

\begin{lem}
Let $k$ be a non-trivially valued Banach field, and $\{A\rightarrow A_{i}\}$ be a collection of maps of $\mathrm{T}$-finitely presented algebras, where $\mathrm{T}\in\{\mathrm{Disc}_{k},\mathrm{Afnd}^{\dagger}_{k}\}$ in the Archimedean case or $\mathrm{T}\in\{\mathrm{Disc}^{nA}_{k},\mathrm{Afnd}^{nA,\dagger}_{k},\mathrm{Tate}_{k}\}$ in the non-Archimedean case. If the collection of functors
$$\{A_{i}\hat{\otimes}_{A}(-):{}_{A}\mathbf{Mod}^{\heart}\rightarrow{}_{A_{i}}\mathbf{Mod}^{\heart}\}$$
is jointly conservative when restricted to Banach $A$-modules, then $\{\mathcal{M}(A_{i})\rightarrow\mathcal{M}(A)\}$ is surjective.
\end{lem}

\begin{proof}
Let $A\rightarrow k$ be a map with $k$ a Banach field. By conservativity some $A\hat{\otimes}_{A_{i}}k\neq0$. By Corollary \ref{lem:nonzeromaptofield}, Proposition \ref{nonzeroaffinoid}, and Proposition \ref{prop:daggernonemptyspec}, there is a map $A\hat{\otimes}_{A_{i}}k\rightarrow l$ with $l$ a Banach field. Thus the class of $k$ in $\mathcal{M}(A)$ is in the image of $\mathcal{M}(A_{i})\rightarrow\mathcal{M}(A)$.
\end{proof}

\begin{cor}\label{cor:daggerconservativesurj}
Let $k$ be a non-trivially valued Banach field, let $\mathrm{T}\in\{\mathrm{Disc}_{k},\mathrm{Afnd}^{\dagger}_{k}\}$ in the Archimedean case or $\mathrm{T}\in\{\mathrm{Disc}^{nA}_{k},\mathrm{Afnd}^{nA,\dagger}_{k},\mathrm{Tate}_{k}\}$ in the non-Archimedean case, and let $\{A\rightarrow B_{i}\}$ be a collection of maps of finitely $\mathrm{T}$-presented complete bornological $k$-algebras. Suppose that $B_{i}\cong A\hat{\otimes}\mathrm{Free_{T}}(\lambda_{1},\ldots,\lambda_{n})\big\slash(f_{1}^{(i)},\ldots, f_{m_{i}}^{(i)})$ where $(f_{1}^{(i)},\ldots, f_{m_{i}}^{(i)})$ is a closed ideal and for any map $A\rightarrow K$ with $K$ a non-trivially valued field, the image of $f_{1}^{(i)},\ldots, f_{m_{i}}^{(i)}$ in $K\cong A\hat{\otimes}\mathrm{Free_{T}}(\lambda_{1},\ldots,\lambda_{n})\big\slash(f_{1}^{(i)},\ldots, f_{m_{i}}^{(i)})$ is a regular sequence. If the collection of functors \[\{B_{i}\hat{\otimes}_{A}(-):{}_{A}\mathrm{Mod}(\mathrm{LH(Ind(Ban}_{k}))\rightarrow{}_{B_{i}}\mathrm{Mod}(\mathrm{LH(Ind(Ban}_{k}))\}\] is jointly conservative then the induced map $\coprod\mathcal{M}(B_{i})\rightarrow\mathcal{M}(A)$ is surjective.
\end{cor}

\begin{proof}
Let $A\rightarrow F$ be a map with $F$ being a complete extension of $k$. By conservativity we must have that some $B_{i}\hat{\otimes}_{A}F$ is non-zero. Now $B_{i}\hat{\otimes}_{A}F\cong B_{i}\hat{\otimes}_{A}^{\mathbb{L}}F$ by Subsection \ref{subsubsec:derquotom}. Thus the tensor product $B_{i}\hat{\otimes}_{A}F$ may in fact be computed in $\mathrm{CBorn}_{F}$. Also by the various results of the previous section there is a map $B_{i}\hat{\otimes}_{A}F\rightarrow L$ with $L$ a complete extension of $k$. This gives a point of $\mathcal{M}(B_{i})$ which maps to the equivalence class of $A\rightarrow F$ in $\mathcal{M}(A)$.  
\end{proof}

\begin{prop}\label{lem:productloc}
Let $k$ be a non-trivially valued non-Archimedean Banach field. Let $\mathrm{T}\in\{\mathrm{Disc_{k}},\mathrm{Afnd}^{\dagger}_{k}\}$ in the Archimedean case, and $\mathrm{T}\in\{\mathrm{Disc}^{nA}_{k},\mathrm{Afnd}^{nA,\dagger}_{k},\mathrm{Tate}_{k}\}$ in the non-Archimedean case. Let $A_{1},\ldots,A_{n}\in\mathbf{DAlg}^{cn,\mathrm{T}-coh}$. Then 
$$\{\mathrm{Spec}(A_{i})\rightarrow\mathrm{Spec}(\bigoplus A_{j})\}_{i=1}^{n}$$
is a cover in $\overline{\tau}_{\mathrm{T}-open}$. 

\end{prop}

\begin{proof}
    Each map $\bigoplus A_{j}\rightarrow A_{i}$ is a derived strong homotopy monomorphism by Proposition \ref{prop:projhtpy}. We use Example \ref{Example:htpymonodescendable} - we just need to observe that ${}_{\prod_{i=1}^{n}A_{i}}\mathbf{Mod}\cong\prod_{i=1}^{n}{}_{A_{i}}\mathbf{Mod}$, and the functor $A_{j}\hat{\otimes}^{\mathbb{L}}_{\prod_{i=1}^{n}A_{i}}(-)$ is just the projection onto the $j$th factor.

    Now $\pi_{0}(\prod_{i=1}^{n}A_{i})\cong\prod_{i=1}^{n}\pi_{0}(A_{i})$. The direct sum of disrete $\mathrm{T}$-finitely presented algebras is discrete $\mathrm{T}$-finitely presented. If $C\cong \mathrm{F}(\mathrm{Free_{T}}(\lambda_{1},\ldots\lambda_{m}))\big\slash I$ and $D\cong  \mathrm{F}(\mathrm{Free_{T}}(\gamma_{1}\ldots\gamma_{n})\big\slash J)$ then $C\oplus D$ is given by $\mathrm{F}(\mathrm{Free_{T}}(\lambda_{1},\ldots,\lambda_{m},\gamma_{1},\ldots,\gamma_{n}))\big\slash (I+J+\{x_{i}y_{j}\}_{1\le i\le m,1\le j\le n})$ where $x_{i}$ is the $i$th coordinate of $\mathrm{F}(\mathrm{Free_{T}}(\lambda_{1},\ldots,\lambda_{m}))$, and $y_{j}$ is the $j$th coordinate of $\mathrm{F}(\mathrm{Free_{T}}(\gamma_{1},\ldots,\gamma_{n}))$.

    If $C$ is discretely $\mathrm{T}$-finitely presented and $\prod_{i=1}^{n}A_{i}\rightarrow C$ is a map, then it factors through $\prod_{i=1}^{n}\pi_{0}(A_{i})$. Moreover we have 
    $$C\hat{\otimes}^{\mathbb{L}}_{\prod_{i=1}^{n}A_{i}}A_{j}\cong C\hat{\otimes}_{\prod_{i=1}^{n}\pi_{0}(A_{i})}^{\mathbb{L}}\hat{\otimes}^{\mathbb{L}}_{\prod_{i=1}^{n}A_{i}}A_{j}\cong C\hat{\otimes}_{\prod_{i=1}^{n}\pi_{0}(A_{i})}^{\mathbb{L}}\cong\pi_{0}(C)\hat{\otimes}\pi_{0}(A_{j})$$
\end{proof}

\section{Analytic Geometries}

 To avoid repetition, throughout if $k$ is non-Archimedean then the '$nA$' decorations will be implicit.

\subsection{Stein Geometry}

Let $R$ be any Banach ring. We beging by discussing geometry relative to Stein spaces. In this context we will use the topology of open embeddings.

\begin{defn}
The \textit{ derived pre-}$G$-\textit{ analytic context over} $R$ is
$$(\mathbf{Aff}^{cn}_{\mathrm{Ind(Ban}_{R}\mathrm{)}},\mathbf{sm}_{o}^{\mathrm{Disc}_{R}},\mathrm{G_{\mathrm{Disc_{R}}}^{pre}})$$
The \textit{ derived }$G$-\textit{ analytic context over} $R$  is
$$(\mathbf{Aff}^{cn}_{\mathrm{Ind(Ban}_{R}\mathrm{)}},\mathbf{sm}_{o}^{\mathrm{Disc}_{R}},\mathrm{G_{\mathrm{Disc_{R}}}})$$
The \textit{ derived \'{e}tale analytic context over} $R$  is
$$(\mathbf{Aff}^{cn}_{\mathrm{Ind(Ban}_{R}\mathrm{)}},\mathbf{sm}^{\mathrm{Disc}_{R}},\tau_{\mathrm{Disc_{R}}-\textrm{\'{e}t}})$$
\end{defn}

\begin{rem}
    Over $\mathbb{C}$ a $\mathrm{Disc_{\mathbb{C}}}$-\'{e}tale map of discrete reduced dagger Stein algebras corresponds locally in the $\overline{\mathrm{G}}^{\aleph_{1}}_{\mathrm{Disc_{\mathbb{C}}}}$-topology to a $\mathrm{Disc_{\mathbb{C}}}$-rational embedding.
\end{rem}

\begin{defn}
The \textit{ countable strong derived pre-}$G$-\textit{ analytic context over} $R$ is
$$(\mathbf{Aff}^{cn}_{\mathrm{LH(Ind(Ban}_{R}\mathrm{)}},\overline{\mathbf{sm}}_{o}^{\mathrm{Disc}_{R},\aleph_{1}},\overline{\mathrm{G}}_{\mathrm{Disc_{R}}}^{pre,\aleph_{1}})$$
The \textit{ countable strong derived }$G$-\textit{ analytic context over} $R$ is
$$(\mathbf{Aff}^{cn}_{\mathrm{LH(Ind(Ban}_{R}\mathrm{)}},\overline{\mathbf{sm}}_{o}^{\mathrm{Disc}_{R},\aleph_{1}},\overline{\mathrm{G}}_{\mathrm{Disc_{R}}}^{\aleph_{1}})$$
\end{defn}

Now fix a non-trivially valued Banach field $k$. Let $\mathbf{dStn}_{k}^{\dagger}\defeq\mathbf{dStnAlg}^{\dagger}$ be the category opposite to the category of derived dagger Stein algebras. following \cite{bambozzi2015stein} we call this the category of derived dagger Stein spaces. Define $\mathbf{dStn}_{k}^{\dagger,f}$ to be the category of those derived dagger Stein spaces $\mathrm{Spec}(A))$ such that $\pi_{0}(A)\cong\mathcal{O}(D^{n}_{k,<\rho})\big\slash I$ where $I$ is a closed ideal. We call this the category of derived finitely embeddable Stein spaces. We also denote by $\mathbf{dStn}_{k}^{\dagger,gf}\subseteq\mathbf{dStn}_{k}^{\dagger,f}$ the full subcategory consisting of those $\mathrm{Spec}(A)$ where $\pi_{0}(A)\cong\mathcal{O}(D^{n}_{k,<\rho})\big\slash I$ and $I$ is closed and finitely generated. Denote by $\mathbf{Coad}_{+}\subseteq\mathbf{QCoh}|_{\mathbf{dStn}_{k}^{\dagger}}$ the sub-presheaf sending $\mathrm{Spec}(A)$ to $\mathbf{Coad}_{+}(A)$. Indeed this is a sub-presheaf precisely by Theorem \ref{thm:coadpresheaf}.

We introduce a slight variant on the class $\mathbf{Sm}^{\mathrm{T}}$, and denote it by $\mathbf{Sm}^{\mathrm{Stein^{\dagger}_{k}}}$, It is the class of formally smooth maps
$$A\rightarrow B$$
such that whenever $A\rightarrow C$ is a map with $C\cong\mathcal{O}(D^{n}_{k<\rho})\big\slash I$ for $I$ a \textit{closed ideal}, then $B\otimes_{A}^{\mathbb{L}}C$ is also of this form. This slight variant is necessary because not all ideals of $\mathcal{O}(D^{n}_{k<\rho})$ need be closed. All of the classes of maps from Definition \ref{defn:Tstrongtop} and Definition \ref{defn:Tstrongmaps} are modified in the obvious way.

\begin{thm}
Let $k$ be a non-trivially valued Banach field.
        
        \begin{enumerate}
            \item 
            $$(\mathbf{Aff}^{cn}_{\mathrm{LH(Ind(Ban}_{R}\mathrm{)}},\mathrm{G}_{\mathrm{Disc}_{k}}^{pre,\aleph_{1}},\overline{\mathbf{sm}}_{o}^{\mathrm{Disc}_{k},\aleph_{1}},(\mathbf{dStn}^{\dagger,f})^{\heart},\mathbf{Coad}^{\heart})$$
        $$(\mathbf{Aff}^{cn}_{\mathrm{LH(Ind(Ban}_{R}\mathrm{)}},\mathrm{G}^{pre}_{\mathrm{Disc}_{k}},\overline{\mathbf{sm}}_{o}^{\mathrm{Disc}_{k}},(\mathbf{dStn}^{\dagger,f})^{\heart},\mathbf{Coad}^{\heart})$$
        $$(\mathbf{Aff}^{cn}_{\mathrm{LH(Ind(Ban}_{R}\mathrm{)}},\mathrm{G}_{\mathrm{Disc}_{k}}^{pre,\aleph_{1}},\overline{\mathbf{sm}}_{o}^{\mathrm{Disc}_{k},\aleph_{1}},(\mathbf{dStn}^{\dagger,gf})^{\heart},\mathbf{Coad}^{\heart})$$
        $$(\mathbf{Aff}^{cn}_{\mathrm{LH(Ind(Ban}_{R}\mathrm{)}},\mathrm{G}^{pre}_{\mathrm{Disc}_{k}},\overline{\mathbf{sm}}_{o}^{\mathrm{Disc}_{k}},(\mathbf{dStn}^{\dagger,gf})^{\heart},\mathbf{Coad}^{\heart})$$
        are Cartan contexts. In particular $\mathbf{Coad}_{\ge n}$ and $\mathbf{Perf}_{\ge n}|_{\mathbf{dStn}^{\dagger,f}}$ satisfy hyperdescent in the first case, and $\mathbf{Coad}_{+}$ and $\mathbf{Perf}|_{\mathbf{dStn}^{\dagger,f}}$ satisfy descent in the second case. 
        \item 
    Restricting to \'{e}tale maps we get the following.
    $$(\mathbf{Aff}^{cn}_{\mathrm{LH(Ind(Ban}_{R}\mathrm{)}},\mathrm{G}_{\mathrm{Disc}_{k}}^{\aleph_{1},pre},\overline{\textbf{\'{et}}}^{\mathrm{Disc}_{k},\aleph_{1}},(\mathbf{dStn}^{\dagger,f})^{\heart},\mathbf{Coad}^{\heart})$$
        $$(\mathbf{Aff}^{cn}_{\mathrm{LH(Ind(Ban}_{R}\mathrm{)}},\mathrm{G}^{pre}_{\mathrm{Disc}_{k}},\overline{\textbf{\'{et}}}^{\mathrm{Disc}_{k}},(\mathbf{dStn}^{\dagger,f})^{\heart},\mathbf{Coad}^{\heart})$$
        $$(\mathbf{Aff}^{cn}_{\mathrm{LH(Ind(Ban}_{R}\mathrm{)}},\mathrm{G}_{\mathrm{Disc}_{k}}^{\aleph_{1},pre},\overline{\textbf{\'{et}}}^{\mathrm{Disc}_{k},\aleph_{1}},(\mathbf{dStn}^{\dagger,gf})^{\heart},\mathbf{Coad}^{\heart})$$
        $$(\mathbf{Aff}^{cn}_{\mathrm{LH(Ind(Ban}_{R}\mathrm{)}},\mathrm{G}^{pre}_{\mathrm{Disc}_{k}},\overline{\textbf{\'{et}}}^{\mathrm{Disc}_{k}},(\mathbf{dStn}^{\dagger,gf})^{\heart},\mathbf{Coad}^{\heart})$$
        are \'{e}tale Cartan contexts.
        \end{enumerate}
\end{thm}

\begin{proof}

First we observe that these are well-defined contexts, i.e. if $\mathrm{Spec}(B_{i})\rightarrow\mathrm{Spec}(B)\}$ is a cover in $\mathrm{G_{}}^{pre,\aleph_{1}}$ then each $A\rightarrow B_{i}$ is in $\overline{\mathbf{sm}}_{o}^{\mathrm{Disc_{k}},\aleph_{1}}$. Indeed let $A\rightarrow C$ be a map with $C$ a discrete dagger Stein algebra. By \cite{bambozzi2015stein} Theorem 5.5 the map $C\rightarrow\pi_{0}(B\otimes_{A}^{\mathbb{L}}C)$ is a homotopy monomorphism. Thus the map $C\rightarrow B\otimes_{A}^{\mathbb{L}}C$ is derived strong, so that $B\otimes_{A}^{\mathbb{L}}C$ is in fact a discrete dagger Stein algebra. 

\begin{enumerate}
    \item 
    Now we verify the Cartan context conditions.

The above also shows that

$$(\mathbf{Aff}^{cn}_{\mathrm{LH(Ind(Ban}_{R}\mathrm{)}},\mathrm{G}_{T}^{\aleph_{1}},\overline{\mathbf{sm}}_{o}^{\mathrm{Disc}_{k},\aleph_{1}},(\mathbf{dStn}^{\dagger,f})^{\heart})$$
$$(\mathbf{Aff}^{cn}_{\mathrm{LH(Ind(Ban}_{R}\mathrm{)}},\mathrm{G}_{T}^{\aleph_{1}},\overline{\mathbf{sm}}_{o}^{\mathrm{Disc}_{k},\aleph_{1}},(\mathbf{dStn}^{\dagger,f})^{\heart})$$
are strong relative $(\infty,1)$-pre-geometry tuples.

 For $A\in\mathrm{Stn}_{k}^{\dagger,f}$ we have that $\mathrm{Coad}(\mathrm{Spec}(A))$ is thick and closed under finite limits and colimits by Proposition \ref{prop:coad}.

 By construction if $\mathrm{Spec}(B)\rightarrow\mathrm{Spec}(A)$ is a map in $\overline{\mathbf{sm}}_{o}^{\mathrm{Disc}_{k},\aleph_{1}}$ and $\mathrm{Spec}(C)\rightarrow\mathrm{Spec}(A)$ with $A\in\mathrm{Stn}_{k}^{\dagger,f}$ then $\mathrm{Spec}(C\otimes_{A}^{\mathbb{L}}B)$ is in $\mathrm{Stn}_{k}^{\dagger,f}$.

 By definition $\mathrm{Spec}(A)\in\mathrm{Coad}(\mathrm{Spec}(A))$ for any $\mathrm{Spec}(A)\in\mathrm{Stn}_{k}^{\dagger,f}$.

 By definition $\overline{\mathbf{sm}}_{o}^{\mathrm{Disc}_{k},\aleph_{1}}$ consists of formally smooth maps.
 
Let $f:\mathrm{Spec}(B)\rightarrow\mathrm{Spec}(A)$ be a map  $\mathrm{Stn}_{k}^{\dagger,f}$. By descent for transversality we may assume that $f$ is of the form $A\rightarrow A\hat{\otimes}\mathcal{O}(D^{n}_{k,<\rho})$ which is even a flat map.

$\mathbf{Coad}_{+}$ is a subpresheaf of $\mathbf{QCoh}|_{\mathbf{dStn}^{\dagger}_{k}}$ by Theorem \ref{thm:coadpresheaf}.

That $\mathrm{Coad}|_{\mathrm{Stn_{k}}^{\dagger}}$ satisfies descent for $\mathrm{G_{T}}^{\aleph_{1}}|_{\mathrm{Stn^{\dagger}}_{k}}$ follows from descent for coherent sheaves on Stein spaces, Cartan's/ Kiehl's Theorems A and B, and the open mapping theorem. 
\item
For the \'{e}tale restriction, all that reamins to prove is that if $\mathrm{Spec}(A)$ is in $\mathbf{dStn}^{\dagger}_{k}$ then $\mathbb{L}_{\mathrm{Spec}(A)}$ is in $\mathbf{Coad}_{\ge0}(\mathrm{Spec}(A))$. This is Lemma \ref{lem:cotangentfembed}.
\end{enumerate}
\end{proof}

\begin{cor}
\begin{enumerate}
    \item 
    $$(\mathbf{Aff}^{cn}_{\mathrm{LH(Ind(Ban}_{R}\mathrm{)}},\mathrm{G}_{\mathrm{Disc}_{k}}^{\aleph_{1}},\overline{\mathbf{sm}}_{o}^{\mathrm{Disc}_{k},\aleph_{1}},\mathbf{dStn}^{\dagger,f})$$
        $$(\mathbf{Aff}^{cn}_{\mathrm{LH(Ind(Ban}_{R}\mathrm{)}},\mathrm{G}_{\mathrm{Disc}_{k}},\overline{\mathbf{sm}}_{o}^{\mathrm{Disc}_{k}},\mathbf{dStn}^{\dagger,f}))$$
        are strong, hyper relative $(\infty,1)$-geometric contexts.
        \item 
              $$(\mathbf{Aff}^{cn}_{\mathrm{LH(Ind(Ban}_{R}\mathrm{)}},\mathrm{G}_{\mathrm{Disc}_{k}}^{\aleph_{1}},\overline{\textbf{\'{e}t}}^{\mathrm{Disc}_{k},\aleph_{1}},\mathbf{Coad}_{+})$$
        $$(\mathbf{Aff}^{cn}_{\mathrm{LH(Ind(Ban}_{R}\mathrm{)}},\mathrm{G}_{\mathrm{Disc}_{k}},\overline{\textbf{\'{e}t}}^{\mathrm{Disc}_{k}},\mathbf{Coad}_{+}))$$
        are weak relative \'{e}tale $(\infty,1)$-AG contexts.

\end{enumerate}    
  
\end{cor}

. 

\begin{prop}
    Let $k$ be a non-trivially valued Banach field. Let $A\rightarrow B$ be a map of discrete finitely $\mathrm{Disc_{k}}$-presented algebras. Then the following are equivalent
    \begin{enumerate}
        \item 
        $f$ is in $\overline{\mathbf{open}}_{\mathrm{Stein^{\dagger}_{k}}}$
        \item 
        $f$ is a homotopy monomorphism
        \item 
        $\mathrm{Max}(f):\mathrm{Max}(B)\rightarrow\mathrm{Max}(A)$
        is an open immersion of Stein spaces.
    
    \end{enumerate}
\end{prop}

\begin{proof}
    This is a consequence of \cite{bambozzi2015stein} Theorem 5.5, Theorem 5.7
\end{proof}

\begin{lem}\label{lem:discratrefine}
    If $k$ is a non-trivially valued Banach field then the topologies $\overline{\mathrm{G}}^{pre,\aleph_{1}}_{\mathrm{Disc_{k}}-rat}|_{\mathbf{Aff}^{cn,\mathrm{T}-coh}}$ and $\overline{\mathrm{G}}^{\aleph_{1}}_{\mathrm{Disc_{k}}-open}|_{\mathbf{Aff}^{cn,\mathrm{Disc_{k}}-coh}}$ are equivalent.
\end{lem}

\begin{proof}
    Let $\{\mathrm{Spec}(A_{i})\rightarrow\mathrm{Spec}(A)\}$ be a cover in $\overline{\tau}^{\aleph_{1}}_{\mathrm{Disc_{R}}-rat}$. Then it is clearly a cover in $\overline{\tau}^{\aleph_{1}}_{\mathrm{Disc_{R}}-open}$. Conversely suppose that it is a cover in $\overline{\tau}^{\aleph_{1}}_{\mathrm{Disc_{R}}-open}|_{\mathbf{Aff}^{cn,\mathrm{T}-coh}}$. Then $\{\mathrm{Spec}(\pi_{0}(A_{i}))\rightarrow\mathrm{Spec}(\pi_{0}(A))\}$ is a cover of the finitely embeddable Stein space corresponding to $\mathrm{Spec}(\pi_{0}(A))$ by open Stein subdomains $\mathrm{Spec}(\pi_{0}(A_{i}))$. There is a surjective $\mathcal{O}(k^{n})\rightarrow\pi_{0}(A)$, and there is an open subset $V$ of $\mathbb{C}^{n}$ such that $\mathrm{Spec}(\pi_{0}(A_{i}))$ is the intersection of $V$ with $\mathrm{Spec}(\pi_{0}(A))$. Now $V$ may be covered by polydiscs 
    $\{U_{j_{i}}\}_{j_{i}\in\mathcal{J}_{i}}$ in $k^{n}$. 
    The maps $\mathcal{O}(k^{n})\rightarrow\mathcal{O}(U_{i_{j}})$ are $\mathrm{Disc}_{k}$-rational localisations. Then the maps $\pi_{0}(A)\rightarrow\pi_{0}(A_)\otimes_{\mathcal{O}(\mathbb{C}^{n})}^{\mathbb{L}}\mathcal{O}(U_{i_{j}})$ 
    are $\mathrm{Disc}_{k}$-rational localisations, and running over all $i$ and all $j_{i}$ we get a $\mathrm{Disc}_{k}$-rational cover of $\mathrm{Spec}(\pi_{0}(A))$ refining the original one. 
    Then $\{\mathrm{Spec}(A\otimes^{\mathbb{L}}_{\mathcal{O}(k^{n})}\mathcal{O}(U_{i_{j}})\rightarrow\mathrm{Spec}(A)\}$ is a $\mathrm{Disc}_{k}$-rational cover refining the original one, as required. 
\end{proof}

\begin{thm}\label{thm:steintopchar}
    Let $k$ be a non-trivially valued Banach field and let $A\in\mathbf{dStn}^{\dagger,f}$. Let $\{\mathrm{Spec}(A_{i})\rightarrow\mathrm{Spec}(A)\}$ be a collection of maps where $\pi_{0}(A)$ is a finitely embeddable dagger Stein algebra and each $\pi_{n}(A)$ is transverse to open immersions. The following are equivalent
    \begin{enumerate}
        \item 
        A collection of maps $\{\mathrm{Spec}(A_{i})\rightarrow\mathrm{Spec}(A)\}$ is a cover in $\overline{G}^{\aleph_{1}}_{\mathrm{Disc_{k}}}$ (Resp. $\overline{G}^{pre,\aleph_{1}}_{\mathrm{Disc_{k}}}$).
        \item 
               \begin{enumerate}
            \item 
            Each $A\rightarrow A_{i}$ is derived strong.
            \item 
                   $\{\mathrm{Spec}(\pi_{0}(A_{i})\rightarrow\mathrm{Spec}(\pi_{0}(A))\}$ is a cover in $\overline{G}^{\aleph_{1}}_{\mathrm{Disc_{k}}}$ (resp. $\overline{G}^{pre,\aleph_{1}}_{\mathrm{Disc_{k}}}$)
                   
        \end{enumerate}
        \item
        \begin{enumerate}
            \item 
            Each $A\rightarrow A_{i}$ is derived strong.
            \item 
                   $\{\mathrm{Max}(\pi_{0}(A_{i}))\rightarrow\mathrm{Max}(\pi_{0}(A))\}$ is a usual cover of Stein spaces by open Stein subdomains (resp. rational Stein subdomains)
                   
        \end{enumerate}
        \item 
             \begin{enumerate}
            \item 
            Each $A\rightarrow A_{i}$ is derived strong.
            \item 
            each $\pi_{0}(A)\rightarrow\pi_{0}(A_{i})$ is a homotopy monomorphism (resp. a $\mathrm{Disc_{k}}$-rational localisation).
            \item 
                   $\coprod\mathcal{M}(\pi_{0}(A_{i}))\rightarrow\mathcal{M}(\pi_{0}(A))$ is surjective.
                   
        \end{enumerate}
%
                   
        \end{enumerate}
 
\end{thm}

\begin{proof}
    $(1)$ and $(2)$ are equivalent by Theorem \ref{thm:ABclassfiyingcovers}. 
    
    The equivalence between $(2)$ and $(3)$ is then a question about covers of discrete objects. But this is then just \cite{bambozzi2015stein} Corollary 5.15. $(3)\Rightarrow (4)$ by Corollary \ref{cor:daggerconservativesurj}. Now suppose $(4)$ holds. Let $\mathfrak{m}$ be a finitely generated, and hence closed, maximal ideal of $\pi_{0}(A)$, and let $K=\pi_{0}(A)\big\slash\mathfrak{m}$. By the covering assumption there is some $\pi_{0}(A_{i})$ and a factorisation
    $$\pi_{0}(A)\rightarrow\pi_{0}(A_{i})\rightarrow K$$
    The map $\pi_{0}(A_{i})\rightarrow K$ is therefore an admissible epimorphism with kernel given by a finitely generated maximal ideal $\mathfrak{m}'$. Moreover the preimage of this ideal is $\mathfrak{m}$, as required.
\end{proof}



\subsubsection{A Comparison With Porta and Yue Yu}

Let us compare our formulation of derived complex analytic geometry with the formulation of Porta and Yue Yu. In future work we will also compare our non-Archimedean geometry with theirs. Porta and Yue Yu develop derived complex analytic geometry using the language of pregeometries introduced by Lurie in \cite{lurie2009derived}. Consider the pregeometry $\mathcal{T}_{an}(k)$ defined in \cite{porta2018derived} as follows. 

\begin{enumerate}
\item
The underlying category of $\mathcal{T}_{an}(k)$ is the category of smooth $k$-analytic spaces (rigid in the case that $k$ is a non-Archimedean field).
\item
A morphism in $\mathcal{T}_{an}(k)$ is admissible precisely if it is \'{e}tale.
\item
The topology on $\mathcal{T}_{an}(k)$ is the \'{e}tale topology.
\end{enumerate}

Porta and Yue Yu define the category $\mathbf{dAn}_{k}$ as the full subcategory of the category $\mathbf{Str}^{loc}(\mathcal{T}_{an}(k))$ of $\mathcal{T}_{an}(k)$-structured spaces consisting of those $X=(\mathcal{X},\mathcal{O}_{X})$ such that
\begin{enumerate}
    \item 
    $t_{\le0}(X)$ is the $\mathcal{T}_{an}(k)$-structured space corresponding to a $k$-analytic space in the usual sense.
    \item 
    each $\pi_{n}(\mathcal{O}_{X})$ is a coherent $\pi_{0}(\mathcal{O}_{X})$-module.
\end{enumerate}

There is a topology $\tau_{an}$ on $\mathbf{dAn}_{k}$ whereby a collection of maps $\{f_{i}:U_{i}\rightarrow U\}$ is a cover precisely if
\begin{enumerate}
    \item 
    $\{t_{\le0}(U_{i})\rightarrow t_{\le0}(U)\}$ is a cover by open immersions of analytic spaces in the complex case, or \'{e}tale maps in the rigid analytic case.
    \item 
    the natural map
    $$t_{\le0}(f_{i})^{*}\pi_{n}(\mathcal{O}_{U})\rightarrow\pi_{n}(\mathcal{O}_{U_{i}})$$
    is an isomorphism for each $n$.
\end{enumerate}

 There is a full subcategory $\mathbf{dStn}_{k}\subset\mathbf{dAn}_{k}$ consisting of those $X$ such that $t_{\le0}(X)$ is a Stein space. Any derived analytic space $X$ admits an effective epimorphism $\coprod U_{i}\rightarrow X$ with each $U_{i}$ a derived Stein space and each $U_{i}\rightarrow X$ being an open immersion. In fact all derived analytic spaces admit pseudo-representable hypercovers by Steins. Inspired by \cite{MR4036665}, we let $\mathbf{dStn}_{\mathbb{C}}^{gf}$ denote the full subcategory of $\mathbf{dStn}_{\mathbb{C}}$ on those objects $X$ such that $t_{\le0}X$ is \textit{globally finitely embeddable}, i.e. it admits a closed immersion $t_{\le0}X\rightarrow\mathbb{C}^{n}$ for some $n$ and the ideal of definition is finitely generated. Finally we let $\mathbf{dAn}_{k}^{gf,<\infty}$ denote the full subcategory of $\mathbf{dAn}_{k}$ consisting of those derived analytic spaces which for some $n$ admit a $n$-coskeletal hypercover by finitely embeddable derived Steins $X$. In particular they are geometric stacks for $(\mathbf{dStn}_{k}^{gf},\tau_{an}|_{\mathbf{dStn}^{gf}},\mathbf{open}|_{\mathbf{dStn}_{k}^{gf}})$ where $\mathbf{open}$ is the class of open immersions as defined above. This is not a major restriction. As pointed out in \cite{MR4036665} the globally finitely embeddable restriction is not a major one, as all Stein spaces are locally of his form.

The following is shown in \cite{MR4036665}.

\begin{prop}[ \cite{MR4036665}, Proposition 3.25]
There is a fully faithful functor
$$\mathbb{R}\Gamma:(\mathbf{dStn}^{f}_{\mathbb{C}})^{op}\rightarrow\mathbf{sAlg}_{\mathbf{EFC}_{\mathbb{C}}}$$
 Moreover the essential image is precisely the class of $\mathbf{EFC}_{k}$-algebras $A$ such that $\pi_{0}(A)$ is globally finitely embeddable and $\pi_{n}(A)$ is the module of global sections of a coherent sheaf as a $\pi_{0}(A)$-module. Finally we have $\pi_{n}\mathbb{R}\Gamma(X,\mathcal{O}_{X})\cong\Gamma(\pi_{n}(\mathcal{O}_{X})^{alg})$.
\end{prop}

In particular under the equivalence 
$$\mathbb{R}\Gamma:(\mathbf{dStn}^{f}_{k})^{op}\rightarrow\mathbf{sAlg}_{\mathbf{EFC}_{k}}$$
a collection of maps $\{f_{i}:U_{i}\rightarrow U\}$ is a cover in $\mathbf{dStn}^{f}$ if and only if 
\begin{enumerate}
    \item 
    $\{t_{\le0}(U_{i})\rightarrow t_{\le0}(U)\}$
    is a usual cover of Steins by Stein subodmains.
    \item 
    the map 
    $$\pi_{n}(\mathbb{R}\Gamma(U))\hat{\otimes}_{\pi_{0}(\mathbb{R}\Gamma(U))}\pi_{0}(\mathbb{R}\Gamma(U_{i}))\rightarrow\pi_{n}(\mathbb{R}\Gamma(U_{i}))$$
    is an isomorphism for all $n$.
\end{enumerate}

Now by Lemma \ref{lem:embedEFC} we have a fully faithful functor 
$$\mathbf{F}:\mathbf{sAlg}_{\mathbf{EFC}_{k}}\rightarrow\mathbf{DAlg}^{cn}\mathrm{(}\underline{\mathbf{Ch}}(\mathrm{Ind(Ban}_{k}\mathrm{))}$$
By composition we get a fully faithful functor
$$\mathbf{F}\circ\Gamma^{\mathrm{LT}^{\mathrm{EFC}_{k}}}:(\mathbf{dStn}^{PYY}_{k})^{op}\rightarrow\mathbf{DAlg}^{cn}\mathrm{(}\underline{\mathbf{Ch}}(\mathrm{Ind(Ban}_{k}\mathrm{))}$$
For technical reasons it is convenient to restrict to $\mathbf{dStn}^{gf}$ 


For such $X$, $\pi_{0}(\mathbf{F}\circ\Gamma (X))$ is a Fr\'{e}chet space. By Lemma \ref{lem:Tconstructsequence} the essential image of $\mathbf{F}\circ\Gamma|_{\mathbf{dStn}^{gf}}$  is the class of algebras $A$ such that $\pi_{0}(A)$ is a  globally finitely embeddable dagger Stein algebra and $\pi_{n}(A)$ is an \textit{algebraically coherent} $\pi_{0}(A)$-module. By algebraically coherent here we mean it is the essential image under 
$$\mathrm{F}_{\pi_{0}(A)}:{}_{\pi_{0}(A)}\mathrm{Mod}^{alg}\rightarrow{}_{\pi_{0}(\mathbf{F}(A)}\mathrm{Mod}(\mathrm{Ind(Ban}_{\mathbb{C}\mathrm{)}})$$
of a module of the form $\Gamma(\mathcal{F})$, where $\mathcal{F}$ is a coherent sheaf on the Stein space corresponding to $\pi_{0}(A)$. Note that, because ${}_{\pi_{0}(A)}\mathrm{Mod}^{alg}$ is not closed under limits in ${}_{\pi_{0}(\mathbf{F}(A)}\mathrm{Mod}(\mathrm{Ind(Ban}_{\mathbb{C}\mathrm{)}})$, such a module will be coadmissible in general only if $\mathcal{F}$ were a globally finitely generated coherent sheaf.

Now for $\mathrm{Spec}(B)\in\mathbf{dStn}^{gf}$ we have that the map
$${}_{B}\mathbf{Mod}^{alg}\rightarrow{}_{\mathbf{F}(B)}\mathbf{Mod}(\mathrm{Ind(Ban}_{\mathbb{C}}\mathrm{)})$$
is exact as a consequence of the open mapping theorem for Fr\'{e}chet spaces.. It is moreover monoidal. Let $A\rightarrow B$ be a map in $\mathbf{DAlg}^{cn}(\mathrm{Ind(Ban}_{\mathbb{C}}\mathrm{)})$ where both $\pi_{0}(A)$ and $\pi_{0}(B)$ are globally finitely embeddable dagger Steins and each $\pi_{n}(A)$ is an algebraically coherent $\pi_{0}(A)$-module. We denote the class of such Stein spaces by $\mathbf{dStn}^{\dagger,gf,alg}$.

Now 
$$(\mathbf{Aff}^{cn}_{\mathrm{Ind(Ban_{k})})},\mathbf{dStn}^{\dagger,gf,alg}_{k},\overline{\mathbf{sm}}_{o}^{\mathrm{Disc_{k}}},\overline{\tau}_{\mathrm{Disc_{k}}-open}^{\aleph_{1}})$$
is a strong relative $(\infty,1)$-geometry tuple. Since each $\pi_{n}(A)$ is algebraic it is transverse to open immersions $\pi_{0}(A)\rightarrow C$ with $C$ a finitely embeddable dagger Stein algebra. Thus a collection of maps
$$\{\mathrm{Spec}(B_{i})\rightarrow\mathrm{Spec}(A)\}$$
between objects of $\mathbf{dStn}^{\dagger,gf,alg}_{k}$ is in $\tau_{an}$ if and only if 
$$\{\mathrm{Spec}(\mathbf{F}(B_{i}))\rightarrow\mathrm{Spec}(\mathbf{F}(A))\}$$
is such that 
$$\{\mathrm{Spec}(\pi_{0}(\mathbf{F}(B_{i})))\rightarrow\mathrm{Spec}(\pi_{0}(\mathbf{F}(A))\}$$
is a cover in $\overline{\mathrm{G}}^{\aleph_{1}}_{\mathrm{Disc}_{k}-open}$ and each $A\rightarrow B_{i}$
 is derived strong. In particular the functor $\mathbf{F}\circ\Gamma^{\mathrm{LT}^{\mathrm{EFC}_{k}}}$ also identifies the topologies.



This proves the following.

\begin{thm}\label{thm:pyycompare}
Let $\{A\rightarrow B_{i}\}_{i\in\mathcal{I}}$ be a collection of maps in $\mathbf{sAlg}_{\mathrm{EFC}_\mathbb{C}}$ between $\mathrm{EFC}_{\mathbb{C}}$-algebras such that each $\pi_{0}(A)$ and each $\pi_{0}(B_{i})$ correspond to globally finitely embeddable Stein spaces, and each $\pi_{n}(A)$, $\pi_{n}(B_{i})$ are global sections of coherent sheaves. Then this corresponds to a cover $\tau_{an}$ if and only if $\{\mathbf{F}(A)\rightarrow \mathbf{F}(B)_{i}\}_{i\in\mathcal{I}}$ corresponds to a cover in $\overline{\tau}_{\mathrm{T}-open}$. In particular there is a fully faithful functor
$$\mathbf{dAn}^{f,<\infty}\rightarrow\mathbf{Stk}_{geom}   (\mathbf{Aff}^{cn}_{\mathrm{Ind(Ban_{k})})},\mathbf{dStn}^{\dagger,gf}_{k},\overline{\mathbf{sm}}_{\mathrm{Disc_{k}}},\overline{\tau}_{\mathrm{Disc_{k}}-open}^{\aleph_{1}})$$
\end{thm}

\subsection{Overconvergent Geometry}

Let $R$ be either an Archimedean or a non-Archimedean Banach ring. For simplicity, in what follows $\mathrm{Afnd}^{\dagger}$ will be used to denote both $\mathrm{Afnd}^{\dagger}$ and $\mathrm{Afnd}^{nA,\dagger}$.

\begin{defn}
The \textit{derived overconvergent }$G$-\textit{context over }$R$ is 
$$(\mathbf{Aff}^{cn}_{\mathrm{Ind(Ban}_{R}\mathrm{)}},\mathbf{sm}_{o}^{\mathrm{Afnd}^{\dagger}_{R}},\mathrm{G}_{\mathrm{Afnd}^{\dagger}_{R}-\textrm{rat}})$$
The \textit{derived overconvergent }\'{e}tale context over $R$ is 
$$(\mathbf{Aff}^{cn}_{\mathrm{Ind(Ban}_{R}\mathrm{)}},\mathbf{sm}^{\mathrm{Afnd}^{\dagger}_{R}},\tau_{\mathrm{Afnd}^{\dagger}_{R}-\textrm{\'{e}t}})$$
\end{defn}

\begin{prop}
    Let $k$ be a non-trivially valued Banach field. Let $A\rightarrow B$ be a map of discrete finitely $\mathrm{Afnd^{\dagger}_{k}}$-presented algebras. Then the following are equivalent
    \begin{enumerate}
        \item 
        $f$ is in $\overline{\mathbf{open}}_{\mathrm{Afnd^{\dagger}_{k}}}$
        \item 
        $f$ is a homotopy monomorphism
        \item 
        $\mathcal{M}(f):\mathcal{M}(B)\rightarrow\mathcal{M}(A)$
        is an open immersion of Berkovich dagger affinoids spaces.
    
    \end{enumerate}
\end{prop}

\begin{proof}
    This is a consequence of \cite{MR3448274} Lemma 5.4 and Theorem 5.11.
\end{proof}

The following can be proven exactly as in Lemma \ref{lem:discratrefine}

\begin{lem}
    If $k$ is a non-trivially valued Banach field then the topologies $\mathrm{G_{T}}|_{\mathbf{Aff}^{cn,\mathrm{Afnd^{\dagger}_{k}-coh}}}$ and $\overline{\tau}_{\mathrm{Afnd^{\dagger}_{k}}-open}|_{\mathbf{Aff}^{cn,\mathrm{Afnd^{\dagger}_{k}-coh}}}$ are equivalent.
\end{lem}

\begin{prop}
    Let $k$ be a non-trivially valued Banach field. Then
    
    \begin{enumerate}
        \item 
        $$(\mathbf{Aff}^{cn}_{\mathrm{Ind(Ban}_{k}\mathrm{)}},\mathrm{G}^{pre}_{\mathrm{Afnd}_{k}^{\dagger}},\overline{\mathbf{sm}}^{\mathrm{Afnd}_{k}^{\dagger}}_{o},\mathbf{Aff}^{\heart,\mathrm{Afnd}_{k}^{\dagger}-coh},\mathrm{Coh})$$
    is a Cartan context. In particular $\mathbf{Coh}_{\ge n}|_{\mathbf{Aff^{cn,\mathrm{Afnd}^{\dagger}_{k}-coh}}}$,$\mathbf{Perf}_{\ge n}|_{\mathbf{Aff^{cn,\mathrm{Afnd}^{\dagger}_{k}-coh}}}$, $\mathbf{Coh}_{+}|_{\mathbf{Aff^{cn,\mathrm{Afnd}^{\dagger}_{k}-coh}}}$, and $\mathbf{Perf}|_{\mathbf{Aff^{cn,\mathrm{Afnd}^{\dagger}_{k}-coh}}}$ satisfy hyperdescent.
    \item
    $$(\mathbf{Aff}^{cn}_{\mathrm{Ind(Ban}_{k}\mathrm{)}},\mathrm{G}^{pre}_{\mathrm{Afnd}_{k}^{\dagger}},\overline{\textbf{\'{e}t}}^{\mathrm{Afnd}_{k}^{\dagger}},\mathbf{Aff}^{\heart,\mathrm{Afnd}_{k}^{\dagger}-coh},\mathrm{Coh})$$
       is an \'{e}tale Cartan context.
    \end{enumerate}
       
\end{prop}

\begin{proof}
We use Proposition \ref{prop:cartanT}. By Proposition \ref{prop:NoetherianStrong} discrete finitely presented dagger affinoid algebras are in fact strongly Noetherian. Clearly if $A$ is a discrete finitely presented dagger affinoid then so is $A\hat{\otimes} k<r_{1}x_{1},\ldots,r_{n}x_{n}>^{\dagger}$. Any such $A$ is also $\mathrm{Afnd}^{\dagger}_{k}$-localisation regular by \cite{MR3448274} Theorem 5.4. $A$ is also $\mathrm{Afnd}^{\dagger}_{k}$-localisation coherent by Lemma \ref{lem:transverse}. It remains to prove \v{C}ech descent for discrete finitely presented modules. This follows from \cite{bambozzi2014generalization} Theorem 6.1.21.
\end{proof}

\begin{rem}
    The proof of descent for $\mathbf{Coh}_{+}$ in the dagger analytic non-Archimedean was also established in \cite{soor2023quasicoherent} Theorem 1.1.
\end{rem}

\begin{cor}
   Let $k$ be a non-trivially valued Banach field. Then
\begin{enumerate}
    \item 
$$(\mathbf{Aff}^{cn}_{\mathrm{Ind(Ban_{k})}},\mathbf{Aff}^{cn,\mathrm{Afnd^{\dagger}_{k}}-coh},\overline{\mathbf{sm}}_{o}^{\mathrm{Afnd^{\dagger}_{k}}},\overline{G}^{pre}_{\mathrm{Afnd^{\dagger}_{k}}})$$
        $$(\mathbf{Aff}^{cn}_{\mathrm{Ind(Ban_{k})}},\mathbf{Aff}^{cn,\mathrm{Afnd^{\dagger}_{k}}-\heart},\overline{\mathbf{sm}}_{o}^{\mathrm{Afnd^{\dagger}_{k}}},\overline{G}_{\mathrm{Afnd^{\dagger}_{k}}})$$
        are strong, relative hyper-$(\infty,1)$-geometry tuples.
        \item 
        $$(\mathbf{Aff}^{cn}_{\mathrm{Ind(Ban_{k})}},\mathbf{Aff}^{cn,\mathrm{Afnd^{\dagger}_{k}}-coh},\overline{\textbf{\'{e}t}}^{\mathrm{Afnd^{\dagger}_{k}}},\overline{G}^{pre}_{\mathrm{Afnd^{\dagger}_{k}}})$$
        $$(\mathbf{Aff}^{cn}_{\mathrm{Ind(Ban_{k})}},\mathbf{Aff}^{cn,\mathrm{Afnd^{\dagger}_{k}}-\heart},\overline{\textbf{\'{e}t}}^{\mathrm{Afnd^{\dagger}_{k}}},\overline{G}_{\mathrm{Afnd^{\dagger}_{k}}})$$
        are weak relative \'{e}tale $(\infty,1)$-AG contexts.
\end{enumerate}

    \end{cor}

\begin{thm}\label{lem:afndoverchar}
    Let $k$ be a non-trivially valued non-Archimedean Banach field and let $A\in\mathbf{DAlg}^{cn,\mathrm{Afnd^{\dagger}_{k}-coh}}$. The following are equivalent
    \begin{enumerate}
        \item 
        A collection of maps $\{\mathrm{Spec}(A_{i})\rightarrow\mathrm{Spec}(A)\}$ is a cover in $\mathrm{G^{pre}_{\mathrm{Afnd^{\dagger}}_{k}}}$ (resp. $\overline{\mathrm{G}}_{\mathrm{Afnd^{\dagger}}_{k}}$).
        \item 
               \begin{enumerate}
            \item 
            Each $A\rightarrow A_{i}$ is derived strong.
            \item 
                   $\{\mathrm{Spec}(\pi_{0}(A_{i})\rightarrow\mathrm{Spec}(\pi_{0}(A))\}$ is a cover in $\mathrm{G^{pre}_{Afnd^{\dagger}_{k}}}$ (resp. $\overline{\mathrm{G}}_{Afnd^{\dagger}_{k}}$)
                   
        \end{enumerate}
        \item
        \begin{enumerate}
            \item 
            Each $A\rightarrow A_{i}$ is derived strong.
            \item 
                   $\{\mathcal{M}(\pi_{0}(A_{i}))\rightarrow\mathcal{M}(\pi_{0}(A))\}$ is a usual cover of Berkovich dagger affinoid spaces by rational Berkovich dagger affinoid subdomains (resp. Berkovich dagger affinoid subdomains).
                   
        \end{enumerate}
            \item
        \begin{enumerate}
            \item 
            Each $A\rightarrow A_{i}$ is derived strong.
            \item 
                   $\{\mathrm{Max}(\pi_{0}(A_{i}))\rightarrow\mathrm{Max}(\pi_{0}(A))\}$ is a usual cover of rigid dagger affinoid domains by rational rigid subdomains (resp. rigid subdomains).
                   
        \end{enumerate}
%
                   
        \end{enumerate}

\end{thm}

\begin{proof}
    $(1)$ and $(2)$ are equivalent by Theorem \ref{thm:ABclassfiyingcovers}. The equivalence between $(2)$ and $(3)$ is then a question about covers of discrete objects. But this is then just \cite{MR3448274} Section 5. The equivalence of $(3)$ and $(4)$ follows from  \cite{soor2023quasicoherent} Proposition 3.3. 
\end{proof}


Let
$$(\mathrm{Afnd}^{\dagger}_{k},\mathrm{G})$$
denote the usual dagger affinoid $G$-site, and let $\mathrm{sm}$ denote the class of maps of dagger affinoids which locally in the $G$-topology are standard smooth.

\begin{cor}
For each $n\ge -1$ there is a fully faithful functor
$$\mathbf{Stk}_{n}(\mathrm{Afnd}^{\dagger}_{k},\mathbf{sm},\tau_{\textrm{\'{e}t}})\rightarrow\mathbf{Stk}_{n}(\mathbf{Aff}^{cn}_{\mathrm{Ind(Ban_{k})}},\mathbf{Aff}^{cn,\mathrm{Afnd^{\dagger}_{k}}-coh},\overline{\mathbf{sm}}_{o}^{\mathrm{Afnd^{\dagger}_{k}}},\overline{G}^{pre}_{\mathrm{Afnd^{\dagger}_{k}}})$$
\end{cor}

\subsection{Affinoid Geometry}

Let $R$ be a non-Archimedean Banach ring. In this setting it is natural to consider two contexts: the $G$-context and the \'{e}tale context separately in detail.

\subsubsection{The $G$-Context}
\begin{defn}
The \textit{derived affinoid }$G$\textit{ context over }$R$ is
$$(\mathbf{Aff}^{cn}_{\mathrm{LH(Ind(Ban^{nA}}_{R}\mathrm{)}},\mathbf{sm}^{\mathrm{Tate}_{R}},\mathrm{G}^{pre}_{\mathrm{Tate_{R}}})$$
\end{defn}

\begin{prop}
    Let $k$ be a non-trivially valued Banach field. Let $A\rightarrow B$ be a map of discrete finitely $\mathrm{Tate}_{k}$-presented algebras. Then the following are equivalent
    \begin{enumerate}
        \item 
        $f$ is in $\overline{\mathbf{open}}_{\mathrm{Tate}_{k}}$
        \item 
        $f$ is a homotopy monomorphism
        \item 
        $\mathcal{M}(f):\mathcal{M}(B)\rightarrow\mathcal{M}(A)$
        is an open immersion of Berkovich dagger affinoids spaces.
    
    \end{enumerate}
\end{prop}

\begin{proof}
    This is \cite{koren} Theorem 5.16 and Theorem 5.31.
\end{proof}

The following can be proven exactly as in Lemma \ref{lem:discratrefine}

\begin{lem}
    If $k$ is a non-trivially valued Banach field then the pre-topologies $\mathrm{G_{T}}|_{\mathbf{Aff}^{cn,\mathrm{Afnd^{\dagger}_{k}-coh}}}$ and $\mathrm{G^{pre}_{T}}|_{\mathbf{Aff}^{cn,\mathrm{Afnd^{\dagger}_{k}-coh}}}$ are equivalent.
\end{lem}

\begin{prop}
    Let $k$ be a non-trivially valued non-Archimedean Banach ring. Then 
    \begin{enumerate}
        \item 
    $$(\mathbf{Aff}^{cn}_{\mathrm{Ind(Ban}_{k}\mathrm{)}},\mathrm{G}^{pre}_{\mathrm{Tate}_{k}},\overline{\mathbf{sm}}^{\mathrm{Tate}_{k}}_{o},\mathbf{Aff}^{\heart,\mathrm{Tate}_{k}-coh},\mathrm{Coh})$$
    is a Cartan context. In particular $\mathbf{Coh}_{\ge n}|_{\mathbf{Aff^{cn,\mathrm{Tate_k}-coh}}}$, $\mathbf{Perf}_{\ge n}|_{\mathbf{Aff^{cn,\mathrm{Tate_k}-coh}}}$, $\mathbf{Coh}_{+}|_{\mathbf{Aff^{cn,\mathrm{Tate_k}-coh}}}$, and $\mathbf{Perf}|_{\mathbf{Aff^{cn,\mathrm{Tate_k}-coh}}}$ satisfy hyperdescent. 
    \item
    $$(\mathbf{Aff}^{cn}_{\mathrm{Ind(Ban}_{k}\mathrm{)}},\mathrm{G}^{pre}_{\mathrm{Tate}_{k}},\overline{\textbf{\'{e}t}}^{\mathrm{Tate}_{k}}_{o},\mathbf{Aff}^{\heart,\mathrm{Tate}_{k}-coh},\mathrm{Coh})$$
    is an \'{e}tale Cartan context.
    \end{enumerate}
    
\end{prop}

\begin{proof}
We use Proposition \ref{prop:cartanT}. By Proposition \ref{prop:NoetherianStrong} discrete finitely presented Tate algebras are in fact strongly Noetherian. Clearly if $A$ is a discrete finitely presented Tate algebra then so is $A\hat{\otimes} k<r_{1}x_{1},\ldots,r_{n}x_{n}>$. Any such $A$ is also $\mathrm{Tate}_{k}$-localisation coherent by \cite{koren} Lemma 5.13. $A$ is also $\mathrm{T}$-localisation coherent by Lemma \ref{lem:transverse}.  It remains to prove \v{C}ech descent for discrete finitely presented modules. This follows from e.g. \cite{MR4332074} Theorem 1.4.
\end{proof}

\begin{cor}
   Let $k$ be a non-trivially non-Archimedean valued Banach field. Then
\begin{enumerate}
    \item 
$$(\mathbf{Aff}^{cn}_{\mathrm{Ind(Ban_{k})}},\mathbf{Aff}^{cn,\mathrm{Tate}_{k}-coh},\overline{\mathbf{sm}}_{o}^{\mathrm{Tate}_{k}},\overline{G}^{pre}_{\mathrm{Tate}_{k}})$$
        $$(\mathbf{Aff}^{cn}_{\mathrm{Ind(Ban_{k})}},\mathbf{Aff}^{cn,\mathrm{Tate}_{k}-\heart},\overline{\mathbf{sm}}_{o}^{\mathrm{Tate}_{k}},\overline{G}_{\mathrm{Tate}_{k}})$$
        are strong, relative hyper-$(\infty,1)$-geometry tuples.

    \item 
$$(\mathbf{Aff}^{cn}_{\mathrm{Ind(Ban_{k})}},\mathbf{Aff}^{cn,\mathrm{Tate}_{k}-coh},\overline{\textbf{\'{e}t}}^{\mathrm{Tate}_{k}},\overline{G}^{pre}_{\mathrm{Tate}_{k}})$$
        $$(\mathbf{Aff}^{cn}_{\mathrm{Ind(Ban_{k})}},\mathbf{Aff}^{cn,\mathrm{Tate}_{k}-\heart},\overline{\textbf{\'{e}t}},\overline{G}_{\mathrm{Tate}_{k}})$$
        are weak relative \'{e}tale $(\infty,1)$-AG contexts.
\end{enumerate}

    \end{cor}

\begin{thm}\label{lem:affinoidcoverchar}
    Let $k$ be a non-trivially valued non-Archimedean Banach field and let $A\in\mathbf{DAlg}^{cn,\mathrm{Tate_{k}-coh}}$. The following are equivalent
    \begin{enumerate}
        \item 
        A collection of maps $\{\mathrm{Spec}(A_{i})\rightarrow\mathrm{Spec}(A)\}$ is a cover in $\mathrm{G^{pre}_{\mathrm{Tate}_{k}}}$ (resp. $\overline{\mathrm{G}}_{\mathrm{Tate}_{k}}$).
        \item 
               \begin{enumerate}
            \item 
            Each $A\rightarrow A_{i}$ is derived strong.
            \item 
                   $\{\mathrm{Spec}(\pi_{0}(A_{i})\rightarrow\mathrm{Spec}(\pi_{0}(A))\}$ is a cover in $\mathrm{G^{pre}_{\mathrm{Tate}_{k}}}$ (resp. $\overline{\mathrm{G}}_{\mathrm{Tate}_{k}}$)
                   
        \end{enumerate}
        \item
        \begin{enumerate}
            \item 
            Each $A\rightarrow A_{i}$ is derived strong.
            \item 
                   $\{\mathcal{M}(\pi_{0}(A_{i}))\rightarrow\mathcal{M}(\pi_{0}(A))\}$ is a usual cover of Berkovich affinoid spaces by rational Berkovich affinoid subdomains (resp. rational Berkovich affinoid subdomains).
                   
        \end{enumerate}
            \item
        \begin{enumerate}
            \item 
            Each $A\rightarrow A_{i}$ is derived strong.
            \item 
                   $\{\mathrm{Max}(\pi_{0}(A_{i}))\rightarrow\mathrm{Max}(\pi_{0}(A))\}$ is a usual cover of rigid domains by rational rigid subdomains (resp. rational rigid subdomains).
                   
        \end{enumerate}
%
                   
        \end{enumerate}
 
\end{thm}

\begin{proof}
    $(1)$ and $(2)$ are equivalent by Theorem \ref{thm:ABclassfiyingcovers}. The equivalence between $(2)$ and $(3)$ is then a question about covers of discrete objects. But this is then just \cite{koren} Theorem 5.37. Now suppose $(4)$ holds. Let $\mathfrak{m}$ be a closed maximal ideal of $\pi_{0}(A)$, and let $K=\pi_{0}(A)\big\slash\mathfrak{m}$. By the covering assumption there is some $\pi_{0}(A_{i})$ and a factorisation
    $$\pi_{0}(A)\rightarrow\pi_{0}(A_{i})\rightarrow K$$
    The map $\pi_{0}(A_{i})\rightarrow K$ is therefore an admissible epimorphism with kernel given by a closed maximal ideal $\mathfrak{m}'$. Moreover the preimage of this ideal is $\mathfrak{m}$, as required. Conversely suppose $(4)$ holds. $\mathrm{Max}(\pi_{0}(A))\subseteq\mathcal{M}(\pi_{0}(A))$ is dense. Moreover the image of $\coprod_{i}\mathcal{M}(\pi_{0}(A_{i}))\rightarrow\mathcal{M}(\pi_{0}(A))$ is closed by the proof of \cite{koren} Lemma 5.32. Thus it must be all of $\mathcal{M}(\pi_{0}(A))$. 
\end{proof}

Let
$$(\mathrm{Afnd}_{k},\mathrm{G})$$
denote the usual affinoid $G$-site, and let $\mathrm{sm}$ denote the class of maps of affinoids which locally in the $G$-topology are standard smooth.

\begin{cor}
For each $n\ge -1$ there is a fully faithful functor
$$\mathbf{Stk}_{n}(\mathrm{Afnd}_{k},\mathbf{sm},\tau_{\textrm{\'{e}t}})\rightarrow\mathbf{Stk}_{n}(\mathbf{Aff}^{cn}_{\mathrm{Ind(Ban_{k})}},\mathbf{Aff}^{cn,\mathrm{Tate}_{k}-coh},\overline{\mathbf{sm}}_{o}^{\mathrm{Tate}_{k}},\overline{G}^{pre}_{\mathrm{Tate}_{k}})$$
\end{cor}

\subsubsection{The \'{E}tale Context}

Over a general non-Archimedean ring $R$ we define the \'{e}tale context as follows.

\begin{defn}
The \textit{derived affinoid \'{e}tale context over }$R$ is
$$(\mathbf{Aff}^{cn}_{\mathrm{Ind(Ban^{nA}}_{R}\mathrm{)}},\mathbf{sm}^{\mathrm{Tate}_{R}},\tau_{\mathrm{Tate}_{R}-\textrm{\'{e}t}})$$
\end{defn}

Let us now specialise to a non-trivially valued Banach field $k$. We anaylse the structure of \'{e}tale maps.






\begin{thm}[\cite{MR2014891} Proposition 8.1.2]\label{thm:tateetalefact}
Let $A$ be a finitely $\mathrm{Tate}_{k}$-presented algebra and $A\rightarrow B=A\hat{\otimes}k<r_{1}x_{1},\ldots,r_{n}x_{n}>\big\slash(f_{1},\ldots,f_{n})$ a discrete standard $\mathrm{Tate}_{k}$-\'{e}tale map, i.e. the Jacobian of $(f_{1},\ldots,f_{n})$ is a unit in $B$. There exists a rational covering $\mathcal{M}(A_{j})\rightarrow\mathcal{M}(A)$, free finite rank $\mathrm{Tate}_{k}$-\'{e}tale morphisms $g_{j}:\mathcal{M}(S_{j})\rightarrow\mathcal{M}(A_{j})$ and $\mathrm{Tate}_{k}$-rational localiations $h_{j}:\mathcal{M}(B\hat{\otimes}_{A}A_{j})\rightarrow\mathcal{M}(S_{j})$ such that the following diagram commutes
\begin{displaymath}
\xymatrix{
\mathcal{M}(B\hat{\otimes}_{A}A_{j})\ar[rr]\ar[d] & & \mathcal{M}(B)\ar[d]\\
\mathcal{M}(S_{j})\ar[r] & \mathcal{M}(A_{j})\ar[r] & \mathcal{M}(A)
}
\end{displaymath}
\end{thm}

\begin{cor}\label{cor:fptransetale}
Finitely presented modules are transverse to discrete $\mathrm{Tate}_{k}$-\'{e}tale maps
\end{cor}

\begin{proof}
Let $f:A\rightarrow B$ be an \'{e}tale map. By descent for transversality, we may assume that $f$ factors as $A\rightarrow \tilde{B}\rightarrow B$ where $\tilde{B}$ is, as an $A$-module, free of finite rank, and $\tilde{B}\rightarrow B$ is a $\mathrm{Tate}_{k}$-rational localisation. Finite maps are flat, and coherent sheaves are transverse to $\mathrm{Tate}_{k}$-localisations.
\end{proof}

There is also a derived version of the factorisation above.

\begin{lem}\label{lem:etfact}
Let $f:\mathrm{Spec}(B)\rightarrow\mathrm{Spec}(A)$ be a derived strong $\mathrm{Tate}_{k}$-\'{e}tale map with $A$ in $\mathbf{DAlg}^{cn,\mathrm{Tate_{k}}-coh}$. There exists a finite rational covering $\{\mathrm{Spec}(A_{j})\rightarrow\mathrm{Spec}(A)\}$, free finite rank \'{e}tale morphisms $g_{j}:\mathrm{Spec}(S_{j})\rightarrow\mathrm{Spec}(A_{j})$, and rational localisations $h_{j}:\mathrm{Spec}(B\hat{\otimes}^{\mathbb{L}}_{A}A_{j})\rightarrow\mathrm{Spec}(S_{j})$ such that the following diagram commutes
\begin{displaymath}
\xymatrix{
\mathrm{Spec}(B\hat{\otimes}^{\mathbb{L}}_{A}A_{j})\ar[d]\ar[rr] & & \mathrm{Spec}(B)\ar[d]\\
\mathrm{Spec}(S_{j})\ar[r] &\mathrm{Spec}(A_{j})\ar[r] & \mathrm{Spec}(A)
}
\end{displaymath}
\end{lem}

\begin{proof}
The proof in the case that $A$ and hence $B$ are concentrated in homological degree $0$ is Proposition 8.1.2 in \cite{MR2014891}, and its proof. Now $\pi_{0}(A)\rightarrow\pi_{0}(B)$ is an \'{e}tale map, so we get a factorisation
\begin{displaymath}
\xymatrix{
\mathrm{Spec}(\pi_{0}(B)\hat{\otimes}^{\mathbb{L}}_{\pi_{0}(A)}\overline{A}_{j})\ar[d]\ar[rr] & & \mathrm{Spec}(\pi_{0}(B))\ar[d]\\
\mathrm{Spec}(\overline{S}_{j})\ar[r] &\mathrm{Spec}(\overline{A}_{j})\ar[r] & \mathrm{Spec}(\pi_{0}(A))
}
\end{displaymath}
Now each map  $\pi_{0}(A)\rightarrow\overline{A}_{j}$ may be written as \[\pi_{0}(A)\hat{\otimes}^{\mathbb{L}}k\Bigr<\lambda_{1}^{-1}x_{1},\ldots,\lambda_{n_{j}}^{-1}x_{n_{j}}\Bigr>\big\slash\big\slash(p_{1},\dots,p_{n_{j}}).\] Define $A_{j}\defeq A\hat{\otimes}^{\mathbb{L}}\Bigr<\lambda_{1}^{-1}x_{1},\ldots,\lambda_{n_{j}}^{-1}x_{n_{j}}\Bigr>\big\slash\big\slash(p_{1},\dots,p_{n_{j}})$. Then 
$$\{\mathrm{Spec}(A_{j})\rightarrow\mathrm{Spec}(A)\}$$
is a cover in $\mathrm{G_{Tate_{k}}}^{pre}$ by Theorem \ref{lem:affinoidcoverchar}

We also have \[\overline{S}_{j}\cong\overline{A}_{j}\hat{\otimes}^{\mathbb{L}}k\Bigr<\rho_{j,1}^{-1}x_{1},\ldots,\rho_{j,m_{j}}^{-1}x_{m_{j}}\Bigr>\big\slash\big\slash(q_{j,1},\ldots,q_{j,m_{j}}).\] Define $S_{j}\defeq A_{j}\hat{\otimes}^{\mathbb{L}}k\Bigr<\rho_{j,1}^{-1}x_{1},\ldots,\rho_{j,m_{j}}^{-1}x_{m_{j}}\Bigr>\big\slash\big\slash(q_{j,1},\ldots,q_{j,m_{j}})$. Then $S_{j}\rightarrow A_{j}$ is \'{e}tale. In particular it is derived strong, and therefore is free of finite rank. Now consider the algebra $B\hat{\otimes_{A}}^{\mathbb{L}}A_{j}$. Each map $A_{j}\rightarrow S_{j}$ and $A_{j}\rightarrow B\otimes_{A}^{\mathbb{L}} A_{j}$ is formally \'{e}tale, so $S_{j}\rightarrow B\otimes^{\mathbb{L}}_{A}A_{j}$ is formally \'{e}tale. Thus there is a unique lift of the localisation $\overline{h}_{j}:\overline{S}_{j}\rightarrow\pi_{0}(B)\hat{\otimes}^{\mathbb{L}}_{\pi_{0}(A)}\overline{A}_{j}$ to a map $h_{j}:S_{j}\rightarrow B\hat{\otimes}^{\mathbb{L}}_{A}A_{j}$ which must itself be formally \'{e}tale. Since $\pi_{0}(h_{j})=\overline{{h}}_{j}$, and $S_{j}\in\mathbf{DAlg}^{cn,\mathrm{Tate_{k}}-coh}$, $h_{j}$ is derived strong and hence is itself a homotopy epimorphism.
\end{proof}

   Let $A$ be a discrete finitely $\mathrm{Tate}_{k}$-presented algebra and $f:A\rightarrow B= A\hat{\otimes}k<r_{1}x_{1},\ldots,r_{n}x_{n}>\big\slash(f_{1},\ldots,f_{n})$ a discrete standard $\mathrm{Tate}_{k}$-\'{e}tale map, i.e. the Jacobian of $(f_{1},\ldots,f_{n})$ is a unit in $B$. Recall that we say $f$ is \textit{regular} if the map
    $$A\hat{\otimes}^{\mathbb{L}}k<r_{1}x_{1},\ldots,r_{n}x_{n}>\big\slash\big\slash(f_{1},\ldots,f_{n})\rightarrow A\hat{\otimes}k<r_{1}x_{1},\ldots,r_{n}x_{n}>\big\slash(f_{1},\ldots,f_{n})$$
    is an equivalence.



\begin{lem}\label{lem:etalediagramequiv}
Let $f_{i}:A\rightarrow B_{i}$ be discrete $\mathrm{Tate}_{k}$-\'{e}tale map of $\mathrm{Tate}_{k}$-finitely presented algebras such that $\mathcal{M}(\prod_{i}B_{i})\rightarrow\mathcal{M}(A)$ is an \'{e}tale cover. Suppose that for each $n$ the map $B^{\hat{\otimes}^{\mathbb{L}}_{A}n+1}\rightarrow B^{\hat{\otimes}_{A}n+1}$ is an equivalence. Then the map
$$A\rightarrow \underset{n}{\mathbf{lim}}B^{\hat{\otimes}^{\mathbb{L}}_{A}n+1}$$
is an equivalence.
\end{lem}

\begin{proof}
By \cite{MR4332074} Theorem 1.4 (or indeed \cite{conrad2003descent}) the map
$$A\rightarrow\underset{n}{\mathbf{lim}} B^{\hat{\otimes}_{A}n+1}$$
is an algebraic equivalence. Now by assumption the map $B^{\hat{\otimes}^{\mathbb{L}}_{A}n+1}\rightarrow B^{\hat{\otimes}_{A}n+1}$ is an equivalence. Thus it suffices to prove that
$$A\rightarrow \underset{n}{\mathbf{lim}} B^{\hat{\otimes}_{A}n+1}$$
is a quasi-isomorphism of Banach spaces. However by the open mapping theorem this follows from it being an algebraic equivalence.
\end{proof}

\begin{cor}\label{thm:etaledescendable}
Let $f:A\rightarrow B\defeq\prod_{i=1}^{n}B_{i}$ be a map with each $A\rightarrow B_{i}$ regular $\mathrm{Tate}_{k}$-\'{e}tale map of $\mathrm{Tate}_{k}$-finitely presented algebras. Suppose that for each $n$ the map $B^{\hat{\otimes}^{\mathbb{L}}_{A}n+1}\rightarrow B^{\hat{\otimes}_{A}n+1}$ is an equivalence. Then $\mathcal{M}(B)\rightarrow\mathcal{M}(A)$ is an \'{e}tale cover if and only if $f$ is descendable.
\end{cor}

\begin{proof}
If $f$ is descendable then by Proposition \ref{lem:productloc} $\mathcal{M}(B)\rightarrow\mathcal{M}(A)$ is surjective.

For the converse note that since each $A\rightarrow B_{i}$ is regular, the Koszul complex is acyclic. Thus $B_{i}$ has a finite resolution by direct sums of tensor product of $A$-modules of the form $A<\lambda_{1}x_{1},\ldots,\lambda_{k}x_{k}>$. These are compact projectives and hence formally $\aleph_{1}$-filtered. Now the claim follows from Lemma \ref{lem:etalediagramequiv} and Lemma \ref{lem:strongdesc}.


$$T^{aug}_{n}(f)\cong \underset{m}{\mathbf{lim}} B^{\hat{\otimes}_{A}m}\hat{\otimes}T^{aug}_{m}(f)$$
\end{proof}

\begin{defn}
Let $\mathbf{DAlg}^{cn,sp}\subseteq\mathbf{DAlg}^{cn,\mathrm{Tate_{k}}-coh}$ (resp. $\mathbf{DAlg}^{\heart,sp}\subseteq\mathbf{DAlg}^{\heart}$) denote the full subcategory consitying of those $A$ such that $\pi_{0}(A)$ is a sous-perfectoid ring. We also write $\mathbf{Aff}^{cn,sp}\defeq(\mathbf{DAlg}^{cn,sp})^{op}$ (resp. $\mathbf{Aff}^{\heart,sp}\defeq(\mathbf{DAlg}^{\heart,sp})^{op}$).
\end{defn}

\begin{lem}[\cite{camargo2024analytic} Lemma 3.5.13]\label{lem:camargo}
    Let $A\in\mathbf{DAlg}^{\heart}(\mathrm{Ind(Ban}_{k}\mathrm{)})$ be a sous-perfectoid algebra and let $B\in\mathbf{DAlg}^{\heart}(\mathrm{Ind(Ban}_{k}\mathrm{)})$ be an affinoid algebra. Then any $\mathrm{Tate}_{k}$-\'{e}tale map $A\rightarrow B$ is regular. 
\end{lem}

\begin{proof}
    Write $B=A\hat{\otimes}^{\mathbb{L}}k<\lambda_{1},\ldots,\lambda_{n}>\big\slash\big\slash(f_{1},\ldots,f_{n})$. Write $\tilde{B}=A\hat{\otimes}k<\lambda_{1},\ldots,\lambda_{n}>\big\slash(f_{1},\ldots,f_{n})\cong\pi_{0}(B)$. The claim is that the map $B\rightarrow\pi_{0}(B)$ is an equivalence. Now it suffices to show that $(f_{1},\ldots,f_{n})$ is an algebraically regular sequence by the open mapping property. This is what is proven in \cite{camargo2024analytic} Lemma 3.5.13.
\end{proof}

\begin{cor}[\cite{camargo2024analytic} Lemma 3.5.13]
If $A\rightarrow B$ is a $\mathrm{Tate_{k}}$-\'{e}tale map with $A$ in $\mathbf{DAlg}^{\heart,sp}$ then $B$ is in $\mathbf{DAlg}^{\heart,sp}$
\end{cor}

\begin{prop}
    Let $k$ be a non-trivially valued non-Archimedean Banach ring. Then 
    \begin{enumerate}
        \item 
    $$(\mathbf{Aff}^{cn}_{\mathrm{Ind(Ban}_{k}\mathrm{)}},\tau_{\mathrm{Tate_{k}}-\textrm{\'{e}t}},\overline{\mathbf{sm}}^{\mathrm{T},\textrm{\'{e}t}},\mathbf{Aff}^{\heart,sp},\mathrm{Coh})$$
    is a Cartan context. In particular $\mathbf{Coh}_{\ge n}|_{\mathbf{Aff^{cn,\mathrm{Tate_k}-coh}}}$, $\mathbf{Perf}_{\ge n}||_{\mathbf{Aff}^{cn,sp}}$, $\mathbf{Coh}_{+}|_{\mathbf{Aff}^{cn,sp}}$, and $\mathbf{Perf}|_{\mathbf{Aff}^{cn,sp}}$ satisfy hyperdescent. 
    \item
      $$(\mathbf{Aff}^{cn}_{\mathrm{Ind(Ban}_{k}\mathrm{)}},\tau_{\mathrm{Tate_{k}}-\textrm{\'{e}t}},\textbf{\'{e}t}^{\mathrm{Tate_{k}}},\mathbf{Aff}^{\heart,sp},\mathrm{Coh})$$
    is an \'{e}tale Cartan context.
    \end{enumerate}
    
\end{prop}

\begin{proof}
We again use Proposition \ref{prop:cartanT}. By Proposition \ref{prop:NoetherianStrong} discrete finitely presented Tate algebras are in fact strongly Noetherian. Clearly if $A$ is a discrete finitely presented Tate algebra then so is $A\hat{\otimes} k<r_{1}x_{1},\ldots,r_{n}x_{n}>$. Moreover it is sous-perfectoid when $A$. By Lemma \ref{lem:camargo} if $A\rightarrow B$ is a $\mathrm{Tate}_{k}$-\'{e}tale map with $A$ discrete and sous-perfectoid then $B$ is also discrete and sous-perfectoid. \'{E}tale maps are transverse to finitely presented modules by Corollary \ref{cor:fptransetale}. It remains to prove \v{C}ech descent for discrete finitely presented modules. This is a weaker version of the main result of \cite{conrad2003descent}.
\end{proof}

\begin{cor}
   Let $k$ be a non-trivially non-Archimedean valued Banach field. Then
\begin{enumerate}
    \item 
  $$(\mathbf{Aff}^{cn}_{\mathrm{Ind(Ban}_{k}\mathrm{)}},\mathbf{Aff}^{cn,sp},\tau_{\mathrm{Tate_{k}}-\textrm{\'{e}t}},\overline{\mathbf{sm}}^{\mathrm{T},\textrm{\'{e}t}})$$
  $$(\mathbf{Aff}^{cn}_{\mathrm{Ind(Ban}_{k}\mathrm{)}},\mathbf{Aff}^{\heart,sp},\tau_{\mathrm{Tate_{k}}-\textrm{\'{e}t}},\overline{\mathbf{sm}}^{\mathrm{T},\textrm{\'{e}t}})$$
        are strong, relative hyper-$(\infty,1)$-geometry tuples.

    \item 
  $$(\mathbf{Aff}^{cn}_{\mathrm{Ind(Ban}_{k}\mathrm{)}},\mathbf{Aff}^{cn,sp},\tau_{\mathrm{Tate_{k}}-\textrm{\'{e}t}},\textbf{\'{e}t}^{\mathrm{Tate}_{k}})$$
  $$(\mathbf{Aff}^{cn}_{\mathrm{Ind(Ban}_{k}\mathrm{)}},\mathbf{Aff}^{\heart,sp},\tau_{\mathrm{Tate_{k}}-\textrm{\'{e}t}},\textbf{\'{e}t}^{\mathrm{Tate}_{k}})$$
        are weak relative \'{e}tale $(\infty,1)$-AG contexts.
\end{enumerate}

    \end{cor}

\begin{thm}\label{thm:sousperfetale}
    Let $A\rightarrow B$ be a map in $\mathbf{Dalg}^{cn,sp}$. The following are equivalent
    \begin{enumerate}
        \item 
        $\mathrm{Spec}(B)\rightarrow\mathrm{Spec}(A)$ is a cover in $\tau_{\mathrm{Tate_{k}}-\textit{\'{e}t}}$.
        \item 
        $A\rightarrow B$ is derived strong and 
        $$\mathcal{M}(\pi_{0}(B))\rightarrow\mathcal{M}(\pi_{0}(A))$$ 
        is an \'{e}tale cover of Berkovich affinoid domains in the usual sense.
    \end{enumerate}
\end{thm}

\begin{proof}
     $(1)\Rightarrow(2)$ by Lemma \ref{lem:strongdesc}. Conversely suppose $(2)$. Again Lemma \ref{lem:strongdesc}, implies that we have an equivalence
     $$B\cong\mathbf{lim}_{n}B^{\hat{\otimes}_{A}^{\mathbb{L}}n}$$
     Using the Koszul complex, $B$ is a finite colimit of $A$-modules of the form $A\hat{\otimes}P$ for $P$ a compact projective. Thus $B$ is formally $\aleph_{1}$-filtered as an $A$-module. Thus Lemma \ref{lem:strongdesc} also gives $(2)\Rightarrow(1)$
\end{proof}




\subsection{Universal Analytic Geometry}

Consider the Banach rings of integers $\mathbb{Z}_{an}$ and $\mathbb{Z}_{triv}$. The former is equipped with the absolute value norm, and the latter with the trivial norm such that $|n|_{triv}=1$ for $n\neq0$. As mentioned previously, $\mathbb{Z}_{an}$ is the initial Banach ring, and $\mathbb{Z}_{triv}$ is the initial non-Archimedean Banach ring. In particular if $R$ is a Banach ring (either Archimedean or non-Archimedean) there is a unique map of Banach rings $\mathbb{Z}_{an}\rightarrow R$.

We gave definitions for the various sites above over arbitrary Banach rings $R$, In particular over $\mathbb{Z}_{an}$ we can define the various Stein, overconvergent, and dagger affinoid sites. Over $\mathbb{Z}_{triv}$ we can define the non-Archimedean versions of these, and also  the affinoid site.





It follows from Subsection \ref{subsec:changeofbase} and Proposition \ref{prop:basechangedescendable}, that 
if $R$ is a Banach ring (either Archimedean or non-Archimedean), the base change functor
$$R\hat{\otimes}_{\mathbb{Z}_{an}}(-)$$
(or $\pi_{R}\circ R\hat{\otimes}_{\mathbb{Z}_{an}}(-)$ in the case that $R$ is non-Archimedean)
sends covers in the Stein and dagger affinoid topologies over $\mathbb{Z}_{an}$ to covers in the corresponding Stein and dagger affinoid topologies over $R$. If $R$ is non-Archimedean then this is also true for the affinoid topologies and the base-change functor
$$R\hat{\otimes}^{nA}_{\mathbb{Z}_{triv}}(-)$$
Thus we get base-change functors on the level of geometric stacks. In particular we can consider $\mathbb{Z}$-forms of complex analytic and rigid analytic spaces.

\subsubsection{A Short Note on Analytification}

Using the process of Subsection \ref{subsec:tan}, we get various analytification functors. We will not go into detail here, but let us list some features - in particular when $\mathrm{An_{T}}$ sends covers to covers, and thus induces a functor on the level of stacks.

\begin{prop}
    Let $k$ be a non-trivially valued Banach field, let $R\cong\mathrm{Sym}(k^{\oplus m})\big\slash I$ be a finitely presented $k$ algebra, and let $f_{1},\ldots,f_{n},a_{1},\ldots,a_{n}$ be elements of $R$ such that $\sum_{i=1}^{n}a_{i}f_{i}=1$. Then 
    $$\{\mathrm{Spec}(\mathcal{O}(D_{k}^{m,<\infty})\big\slash I\hat{\otimes}\mathcal{O}(D^{1,<\infty}_{k})\big\slash(1-xf_{i})\rightarrow\mathrm{Spec}(\mathcal{O}(D_{k}^{m,<\infty})\big\slash I\}$$
    is cover in $\mathrm{G}^{pre}_{\mathrm{Disc_{k}}}$.
\end{prop}

\begin{proof}
    The maps are localisations and therefore homotopy monomorphisms. It remains to prove that the induced map 
    $$\coprod\mathrm{Max}(\mathcal{O}(D_{k}^{m,<\infty})\big\slash I\hat{\otimes}\mathcal{O}(D^{1,<\infty}_{k})\big\slash(1-f_{i}(x)))\rightarrow\mathrm{Max}(\mathcal{O}(D_{k}^{m,<\infty})\big\slash I)$$
    is surjective. As a set $\mathrm{Max}(\mathcal{O}(D_{k}^{m,<\infty})\big\slash I\hat{\otimes}\mathcal{O}(D^{1,<\infty}_{k})\big\slash (1-f_{i}x))$ is the set of points $y\in\mathrm{Max}(\mathcal{O}(D_{k}^{m,<\infty})\big\slash I)$ such that $f_{i}(y)\neq0$. Now $1=\sum_{i=1}^{n}a_{i}f_{i}$ so for each $y\in\mathrm{Max}(\mathcal{O}(D_{k}^{m,<\infty})\big\slash I)$ there is some $j$ such that $f_{j}(y)\neq0$. Hence \[y\in\mathrm{Max}(\mathcal{O}(D_{k}^{m,<\infty})\big\slash I\hat{\otimes}\mathcal{O}(D^{1,<\infty}_{k})\big\slash (1-f_{j}x)).\]
\end{proof}

There is also a non-Archimedean result.

\begin{lem}
Let $R$ be a strongly Noetherian Tate ring, and let $f_{1},\ldots,f_{n},a_{1},\ldots,a_{n}$ be elements of $R$ such that $\sum_{i=1}^{n}a_{i}f_{i}=1$. Then
$$\{\mathrm{Spec}(R\Big<|a_{i}|x\Bigr>\big\slash(1-f_{i}x))\rightarrow\mathrm{Spec}(R)\}_{i=1}^{n}$$
is a cover in $\mathrm{G}^{pre}_{\mathrm{Tate_{k}}}$.
\end{lem}

\begin{proof}
We use \cite{bambozzi2020sheafyness} Theorem 4.15.
The maps are localisations and therefore homotopy monomorphisms.
First assume that each $a_{i}=1$. It suffices to prove that the open sets $\mathrm{Spa}(R\Big<x\Bigr>\big\slash(1-f_{i}x))$ cover $\mathrm{Spa}(R)$. But 
$$\mathrm{Spa}(R\Big<x\Bigr>\big\slash(1-f_{i}x))=\{v:v(1)\le v(f_{i})\}$$
Now 
$$v(1)=v(f_{1}+\ldots+f_{n})\le\mathrm{max}_{1\le i\le n}\{v(f_{i})\}$$
Thus $v(1)\le v(f_{i})$ for some $i$, and the sets $\mathrm{Spa}(R\Big<x\Bigr>\big\slash(1-f_{i}x))$ cover $\mathrm{Spa}(R)$.

Now in the general case we have bounded maps 
$$R\Big<|a_{i}|x\Bigr>\big\slash(1-f_{i}x)\rightarrow R\Big<x\Bigr>\big\slash(1-a_{i}f_{i}x)$$
Since the collection of maps 
$$\{\mathrm{Spec}(R\Big<x\Bigr>\big\slash(1-a_{i}f_{i}x))\rightarrow\mathrm{Spec}(R)\}_{i=1}^{n}$$
is conservative, the collection of maps 
$$\{\mathrm{Spec}(R\Big<|a_{i}|x\Bigr>\big\slash(1-f_{i}x))\rightarrow\mathrm{Spec}(R)\}_{i=1}^{n}$$
is also conservative.
\end{proof}

\subsection{$C^{\infty}$-Geometry}\label{subsec:cinftygeom}

Consider the Lawvere theory $\mathrm{CartSp}_{\mathrm{smooth}}$ from Definition \ref{defn:lawvcart} governing the category $\mathrm{C}^{\infty}\mathrm{Ring}$ of $C^{\infty}$-rings. This is a Fermat theory concretely of homotopy $\mathrm{Ind(Ban_{\mathbb{R}})}$-polynomial type and concretely of $\mathrm{Ind(Ban_{\mathbb{R}})}$-polynomial type. We therefore get fully faithful coproduct preserving functors
$$\mathbf{F}:\mathbf{s}C^{\infty}\mathbf{Ring}\rightarrow\mathbf{DAlg}^{cn}(\mathbf{Ch}(\mathrm{Ind(Ban_{\mathbb{R}})})$$
$$\mathrm{F}:C^{\infty}\mathrm{Ring}\rightarrow\mathrm{Comm}(\mathrm{Ind(Ban_{\mathbb{R}})})$$

\begin{defn}
   A map $f:A\rightarrow B$ in $\mathbf{DAlg}^{cn}$ is said to be a  $C^{\infty}$-\textit{localisation} if it is a (derived) base change of the map 
   $$C^{\infty}(\mathbb{R})\rightarrow C^{\infty}(\mathbb{R}\setminus\{0\})$$
   $f$ is said to be a $C^{\infty}$-\textit{localisation} if it is a composition of basic $C^{\infty}$-localisations.
\end{defn}

The class of all $C^{\infty}$-localisation is denoted $\mathbf{loc}^{C^{\infty}}$

$C^{\infty}$-open immersions are closed under composition and pullback, and contain isomorphisms. By Theorem \ref{thm:localisationflat} we have the following.

\begin{prop}
    A $C^{\infty}$-open immersion is flat. 
\end{prop}

Corollary \ref{cor:flatderstrong} then gives the following.

\begin{cor}
    A map $A\rightarrow B$ in $\mathbf{DAlg}^{cn}$ is a $C^{\infty}$-localisation if and only if it is derived strong and $\pi_{0}(A)\rightarrow\pi_{0}(B)$ is a $C^{\infty}$-localisation.
\end{cor}

Following \cite{carchedi2023derived} let $\mathbf{s}C^{\infty}\mathbf{Ring}^{fp}\subseteq\mathbf{DAlg}^{cn}$ denote the category of compact objects in $\mathbf{s}C^{\infty}\mathbf{Ring}$. By 
\cite{carchedi2023derived} Corollary 3.26 there is an equivalence between $(\mathbf{s}C^{\infty}\mathbf{Ring}^{fp})^{op}$ and the category $\mathbf{DMfld}$ of \textit{derived manifolds}.

\begin{defn}
Say that a collection of maps
$$\{\mathrm{Spec}(A_{i})\rightarrow\mathrm{Spec}(A)\}$$
in $\mathbf{Dalg}^{cn}$ \textit{ is a }$C^{\infty}$-\textit{cover} if
\begin{enumerate}
    \item 
    each $A\rightarrow A_{i}$ is a $C^{\infty}$-localisation.
    \item 
    for any $\mathrm{Spec}(B)\in\mathbf{DMfld}$, 
    $$\{\mathrm{Spec}(A\hat{\otimes}_{A_{i}}^{\mathbb{L}}B)\rightarrow\mathrm{Spec}(B)\}$$
    corresponds to an open cover of $\mathbf{DMfld}$ in the usual sense.
\end{enumerate}
The pre-topology consisting of all such covers is denote $\tau_{C^{\infty}}^{\aleph_{1}}$
\end{defn}

Using the open mapping theorem, it follows from \cite{carchedi2023derived} Lemma 5.3 that if $\{\mathrm{Spec}(A_{i})\rightarrow\mathrm{Spec}(A)\}$ is a cover in $\mathbf{DMfld}$, then there is an equivalence
$$A\cong \underset{n}{\mathbf{lim}}\prod_{(i_{1},\ldots,i_{m})\in\mathcal{I}^{n}}A_{i_{1}}\hat{\otimes}_{A}^{\mathbb{L}}\cdots\hat{\otimes}_{A}^{\mathbb{L}}A_{i_{n}}$$
Thus the pre-topology $\tau_{C^{\infty}}^{\aleph_{1}}$ restricted to $\mathbf{DMfld}$ is \v{C}ech sub-canonical. Finally, let $\tilde{\mathbf{open}}^{C^{\infty},\aleph_{1}}$ denote $(\mathbf{loc}^{C^{\infty}})^{\tau_{C^{\infty}}^{\aleph_{1}}}$ be the class of maps which locally for the topology $\tau_{C^{\infty}}^{\aleph_{1}}$ are $C^{\infty}$-localisations, and let $\mathbf{open}^{C^{\infty}}$ denote the subclass of $(\mathbf{loc}^{C^{\infty}})^{\tau_{C^{\infty}}^{\aleph_{1}}}$ consisting of those maps $A\rightarrow B$ such that whenever $A\rightarrow C$ is a map with $\mathrm{Spec}(C)\in\mathbf{DMfld}$, then $B\otimes_{A}^{\mathbb{L}}C\in\mathbf{DMfld}$. Then 
$$(\mathbf{Aff}^{cn},\mathbf{DMfld},\mathbf{open}^{C^{\infty}},\tau_{C^{\infty}})$$
is relative geometry tuple.

\subsubsection{Discrete Smooth Manifolds}

Specialising to smooth manifolds is more satisfying. Let $\underline{\mathbf{C}}$ be a derived algebraic context. Let ${}_{A}\mathrm{Mod}^{\mathrm{\aleph_{1}}}$ denote the class of $A$-modules whose underlying object in $\mathbf{C}^{\heart}$ formally $\aleph_{1}$-filtered. For $B\in\mathbf{DAlg}^{cn}(\underline{\mathbf{C}})$ denote by ${}_{B}\mathbf{Mod}^{htpy-\aleph_{1}}$ the full subcategory of ${}_{B}\mathbf{Mod}$ such that each $\pi_{m}(M)$ is in ${}_{\pi_{0}(A)}\mathrm{Mod}^{\mathrm{\aleph_{1}}}$. In this instance we get a more categorical description of covers in $\tau^{C^{\infty},\aleph_{1}}$,
\begin{thm}
    Let $U\subset\mathbb{R}^{n}$ be open, and $\{U_{i}\rightarrow U\}_{i\in\mathcal{I}}$ a collection of smooth maps. 
  $\{U_{i}\rightarrow U\}_{i\in\mathcal{I}}$  is a cover of $U$ by smooth open immersions if and only if 
    \begin{enumerate}
        \item 
        each map $C^{\infty}(U)\rightarrow C^{\infty}(U_{i})$ is a flat, Zariski open immersion.
        \item 
        there is a countable subset $J\subset I$ such that for the  collection of functors
        $$C^{\infty}(U_{i})\hat{\otimes}_{C^{\infty}(U)}(-):{}_{C^{\infty}(U)}\mathrm{Mod}^{\mathrm{Fr}}\rightarrow{}_{C^{\infty}(U_{i})}\mathrm{Mod}^{\mathrm{Fr}}$$
        is jointly conservative.
    \end{enumerate}
    \end{thm}

    \begin{proof}
Suppose $\{U_{i}\rightarrow U\}$ is a cover. 
     Let $\{U_{i}\rightarrow U\}_{i\in\mathcal{I}}$ be a cover of $U$ by smooth open immersions. Let $V_{i}=U\setminus U_{i}$, and let $\phi_{i}$ be a function such that $V_{i}$ is the vanishing locus of $\phi_{i}$. Then the map $C^{\infty}(U)\rightarrow C^{\infty}(U_{i})$ factors through the map $C^{\infty}(U)\otimes_{C^{\infty}(\mathbb{R})}\mathbb{R}\rightarrow C^{\infty}(U_{i})$ from the localisation of $C^{\infty}(U)$ at $\phi_{i}$. Moreover this map is an isomorphism. Now consider the \v{C}ech complex
        $$0\rightarrow C^{\infty}(U)\rightarrow\prod_{i\in\mathcal{I}}C^{\infty}(U_{i})\rightarrow\cdots\rightarrow\prod_{(i_{1},\ldots,i_{n})\in\mathcal{I}^{n}}C^{\infty}(U_{i_{1}})\hat{\otimes}_{C^{\infty}(U)} \cdots \hat{\otimes}_{C^{\infty}(U)}C^{\infty}(U_{i_{n}})\rightarrow\cdots$$
        We have $C^{\infty}(U_{i_{1}})\hat{\otimes}_{C^{\infty}(U)}\cdots\hat{\otimes}_{C^{\infty}(U)} C^{\infty}(U_{i_{n}})\cong C^{\infty}(U_{i_{1}}\times_{U} \cdots \times_{U}U_{i_{n}})$. 
        
        By descent and flabbiness for the sheaf of smooth functions this sequence is algebraically exact. Moreover it is a complex of Fr\'{e}chet spaces, so it is also exact in $\mathrm{Ind(Ban_{\mathbb{R}})}$
       This means we have an equivalence
      $$C^{\infty}(U)\rightarrow \underset{k}{\Rlim} \prod_{(n_{1},\ldots,n_{k})\in\mathbb{N}^{k}}C^{\infty}(U_{i_{1}}\times\cdots\times U_{i_{n}})$$
      Now let $M\in {}_{C^{\infty}(U)}\mathrm{Mod}^{\aleph_{1}}$. Then by \cite{ben2020fr} Lemma 7.8 and since the limit is countable we have
      
      $$M\cong \underset{k}{\Rlim} \prod_{(n_{1},\ldots,n_{k})\in\mathbb{N}^{k}}C^{\infty}(U_{i_{1}})\hat{\otimes}_{C^{\infty}(U)}\cdots\hat{\otimes}_{C^{\infty}(U)} C^{\infty}(U_{i_{n}})\hat{\otimes}_{C^{\infty}(U)}M$$
So clearly the functors $C^{\infty}(U_{i})\hat{\otimes}_{C^{\infty}(U)}(-)$ are jointly conservative.

Conversely suppose each map $C^{\infty}(U)\rightarrow C^{\infty}(U_{i})$ is a flat, Zariski open immersion and  there is a countable subset $J\subset I$ such that for the  collection of functors
        $$C^{\infty}(U_{i})\hat{\otimes}_{C^{\infty}(U)}(-):{}_{C^{\infty}(U)}\mathrm{Mod}^{\aleph_{1}}\rightarrow{}_{C^{\infty}(U_{i})}\mathrm{Mod}^{\aleph_{1}}$$
        is jointly convservative. We just need to show that the map $\coprod U_{i}\rightarrow U$ is surjective. Let $x\in U$. The maximal ideal $\mathfrak{m}_{x}=\{f\in C^{\infty}(U):f(x)\neq 0\}$ is closed. Moreover $C^{\infty}(U)\big\slash\mathfrak{m}_{x}$ is a bornological Fr\'{e}chet space. Some $C^{\infty}(U_{i})\hat{\otimes}_{C^{\infty}(U)}C^{\infty}(U)\big\slash\mathfrak{m}_{x}$ is non-zero, and then $x\in U_{i}$. 
    
    \end{proof}

\subsubsection{The Finite Site}

If we are willing to restrict to finite covers things become far more elegant.

\begin{defn}
    Let $\tau^{C^{\infty}}$ denote the pre-topology consisting of finite covers
$$\{\mathrm{Spec}(A_{i})\rightarrow\mathrm{Spec}(A)\}_{i\in\mathcal{I}}$$
such that 
\begin{enumerate}
    \item 
    $A\rightarrow A_{i}$ a $C^{\infty}$-localisation.
    \item 
    $A\rightarrow \underset{i\in\mathcal{I}}{\prod} A_{i}$ is faithful.
\end{enumerate}
\end{defn}

In particular such covers are covers in the faithfully flat topology, and so satisfy hyperdescent.

\begin{prop}
    A finite collection $\{\mathrm{Spec}(A_{i})\rightarrow\mathrm{Spec}(A)\}$ of maps in $\mathbf{DMfld}$ is a cover by open submanifolds if and only if it is in $\tau^{C^{\infty}}$
\end{prop}

\begin{proof}
Each map is a 
    A cover $\{\mathrm{Spec}(A_{i})\rightarrow\mathrm{Spec}(A)\}_{i\in\mathcal{I}}$ in $\tau^{C^{\infty}}$  is a cover in the finite homotopy monomorphism topology and is therefore descendable. In particular we have 
$$A\rightarrow\underset{n}{\mathbf{lim}}\prod_{i_{n}\in\mathcal{I}_{n}}A_{i_{n}}$$
is an equivalence. This implies that $\{\mathrm{Spec}(A_{i})\rightarrow\mathrm{Spec}(A)\}_{i\in\mathcal{I}}$ is a cover by open submanifolds. 

Conversely if $\{\mathrm{Spec}(A_{i})\rightarrow\mathrm{Spec}(A)\}_{i\in\mathcal{I}}$ is a cover by open submanifolds then we have 
$$A\rightarrow \underset{n}{\mathbf{lim}}\prod_{i_{n}\in\mathcal{I}_{n}}A_{i_{n}}$$
is an equivalence
by \cite{carchedi2023derived} Lemma 5.3 This is a finite limit and so clearly $A\rightarrow\underset{i\in\mathcal{I}}{\prod} A_{i}$ is faithful.
\end{proof}

Let $\widetilde{\mathbf{open}}^{C^{\infty}}$ denote $(\mathbf{loc}^{C^{\infty}})^{\tau_{C^{\infty}}}$ be the class of maps which locally for the topology $\tau_{C^{\infty}}$ are $C^{\infty}$-localisations, and let $\mathbf{open}^{C^{\infty}}$ denote the subclass of $(\mathbf{loc}^{C^{\infty}})^{\tau_{C^{\infty}}}$ consisting of those maps $A\rightarrow B$ such that whenever $A\rightarrow C$ is a map with $\mathrm{Spec}(C)\in\mathbf{DMfld}$, then $B\otimes_{A}^{\mathbb{L}}C\in\mathbf{DMfld}$. 

\begin{cor}
    $$(\mathbf{Aff}^{cn},\mathbf{open}^{C^{\infty}},\tau^{C^{\infty}})$$
    is a hyper geometry context, and 
     $$(\mathbf{Aff}^{cn},\mathbf{DMfld},\mathbf{open}^{C^{\infty}},\tau^{C^{\infty}})$$
     is a strong hyper relative geometry context. Moreover $\mathbf{QCoh}$ satisfies hyperdescent.
\end{cor}

\section{Sheaves}\label{sec:classicalsheaves}

Recall that if $\mathpzc{E}$ is a complete and cocomplete quasi-abelian category, and $(X,\tau)$ is a ($1$-categorical) Grothendieck site, then a presheaf $\mathcal{F}\in\mathpzc{Fun}(X^{op},\mathpzc{E})$ is said to be a sheaf if for any covering $\{U_{i}\rightarrow U\}_{i\in\mathcal{I}}$ in $X$  the sequence
$$0\rightarrow\mathcal{F}(U)\rightarrow\prod_{i\in\mathcal{I}}\mathcal{F}(U_{i})\rightarrow\prod_{i,j\in\mathcal{I}}\mathcal{F}(U_{i}\cap U_{j})$$
is strictly exact. 

Throughout this section, by an \textit{analytic algebra} we will mean either a dagger affinoid algebra, a dagger Stein algebra, or a non-Archimedean affinoid algebra. A morphism of analytic algebras will mean a morphism between two dagger affinoid algebras, two dagger Stein algebras, or two affinoid algebras. 

Let $X=\mathcal{M}(A)$ for $A$ a dagger affinoid algebra, dagger Stein algebra, or affinoid algebra (in the case that $k$ is non-Archimedean). Its structure sheaf $\mathcal{O}_{\mathcal{M}(A)}$ has global sections $\mathcal{O}_{\mathcal{M}(A)}(\mathcal{M}(A))=A$. We consider the category
$${}_{\mathcal{O}_{\mathcal{M}(A)}}\mathrm{Mod}\defeq{}_{\mathcal{O}_{\mathcal{M}(A)}}\mathrm{Mod}(\mathpzc{Shv}(\mathcal{M}(A),\Born_{k}))$$
This is a quasi-abelian category. By \cite{qacs} the category ${}_{\mathcal{O}_{\mathcal{M}(A)}}\mathrm{Mod}$ is closed monoidal. Consider the functor $\mathcal{O}_{\mathcal{M}(A)}\hat{\otimes}(-):{}_{A}\mathrm{Mod}\rightarrow{}_{\mathcal{O}_{\mathcal{M}(A)}}\mathrm{Mod}$ sending an $A$-module $M$ to the sheafification of the assignment $U\mapsto\mathcal{O}_{\mathcal{M}(A)}(U)\hat{\otimes}_{A}M$. 

It is left adjoint to the global sections functor $\Gamma_{A}(\mathcal{M}(A),-)$. We therefore get an adjunction.

$$\adj{\mathcal{O}_{\mathcal{M}(A)}\hat{\otimes}_{A}(-)}{{}_{A}\mathrm{Mod}(\mathrm{Ch}(\mathrm{Ind(Ban_{k})}))}{{}_{\mathcal{O}_{\mathcal{M}(A)}}\mathrm{Mod}(\mathrm{Ch}(\mathrm{Shv}(\mathcal{M}(A),\mathrm{Ind(Ban_{k})}}{\Gamma_{A}(\mathcal{M}(A),-)}$$
The left hand side may be equipped with the projective model structure, and the right-hand side with the flat model structure by \cite{kelly2024flat} Theorem 1.4. This renders the adjunction aboe a Quillen adjunction
$$\adj{\mathcal{O}_{\mathcal{M}(A)}\hat{\otimes}_{A}(-)}{Ch({}_{A}\mathrm{Mod})}{Ch({}_{\mathcal{O}_{\mathcal{M}(A)}}\mathrm{Mod})}{\Gamma_{A}(\mathcal{M}(A),-)}$$
which presents an adjunction of $(\infty,1)$-categories
$$\adj{\mathcal{O}_{\mathcal{M}(A)}\hat{\otimes}^{\mathbb{L}}_{A}(-)}{\QCoh_{s}(\mathrm{Spec}(A))}{\Ch({}_{\mathcal{O}_{\mathcal{M}(A)}}\mathrm{Mod})}{\mathbb{R}\Gamma_{A}(\mathcal{M}(A),-)}$$

Our goal is to find a full subcategory of $\Ch({}_{\mathcal{O}_{X}}\mathrm{Mod})$ consisting of objects $\mathcal{F}$ such that the counit 
$$\mathbb{R}\Gamma_{A}(\mathcal{M}(A),\mathcal{O}_{\mathcal{M}(A)}\hat{\otimes}^{\mathbb{L}}_{A}\mathcal{F})\rightarrow \mathcal{F}$$
is an equivalence. The restriction of $\mathbb{R}\Gamma_{A}(\mathcal{M}(A),-)$ to this subcategory in then a fully faithful embedding. 

Much of the theory and proofs in the following is inspired by similar work of \cite{MR1420618}, particularly Chapter 4.

\begin{prop}\label{prop:derivedglobalsectionscompactcolimit}
Let $X$ be a Grothendieck site such that any cover can be refined to a finite subcover. Let $\mathpzc{E}$ be an elementary quasi-abelian category, $\mathcal{I}$ a filtered diagram, and $\mathcal{F}:\mathcal{I}\rightarrow\mathpzc{Shv}(X,\mathpzc{E})$ a diagram. Then for any object $U$ of $X$ there is an equivalence
$$\underset{i\in\mathcal{I}}\colim \mathbb{R}\Gamma(U,\mathcal{F}_{i})\cong\mathbb{R}\Gamma(U,\underset{i\in\mathcal{I}}\colim \mathcal{F}_{i})$$
\end{prop}

\begin{proof}
We may assume that $\mathpzc{E}$ is a the abelian category of modules over a ring. Then this follows from Lemma 21.16.1 in \cite{stacks-project}.
\end{proof}

\begin{prop}
Let  $E$ be a flat $\aleph_{1}$-metrisable bornological $k$-vector space, and let $X$ be an analytic space equipped with the countable $G$-topology. The assignment $U\mapsto\mathcal{O}_{X}(U)\hat{\otimes} E$ for $U\subset X$ a subdomain defines an acyclic sheaf on $X$. If $X$ is equipped with the finite $G$-topology then this is true for all bornological spaces $E$, not necessarily flat. 
\end{prop}

\begin{proof}
For $U\subset X$ the algebra $\mathcal{O}_{X}(U)$ is flat and $\aleph_{1}$-metrisable. Thus for $\{U_{i}\rightarrow U\}$ a countable cover by , $\mathcal{O}_{X}(U)$ and $\underset{i\in\mathcal{I}}\prod \mathcal{O}_{X}(U_{i})$ are flat and $\aleph_{1}$-metrisable.. Since $E$ is $\aleph_{1}$ filtered we have that the sequence
$$0\rightarrow\mathcal{O}_{X}(U)\hat{\otimes} E\rightarrow\Bigr(\prod_{i\in\mathcal{I}}\mathcal{O}_{X}(U_{i})\Bigr)\hat{\otimes} E\rightarrow\Bigr(\prod_{i,j\in\mathcal{I}}\mathcal{O}_{X}(U_{i}\cap U_{j})\Bigr)\hat{\otimes} E\rightarrow\ldots$$
is isomorphic to the complex
$$0\rightarrow\mathcal{O}_{X}(U)\hat{\otimes} E\rightarrow\Bigr(\prod_{i\in\mathcal{I}}\mathcal{O}_{X}(U_{i})\hat{\otimes} E\Bigr)\rightarrow\Bigr(\prod_{i,j\in\mathcal{I}}\mathcal{O}_{X}(U_{i}\cap U_{j})\hat{\otimes} E\Bigr)\rightarrow\ldots$$
Since $E$ is flat the second-complex is acyclic. For the finite topology the complex
$$0\rightarrow\mathcal{O}_{X}(U)\rightarrow\Bigr(\prod_{i\in\mathcal{I}}\mathcal{O}_{X}(U_{i})\Bigr)\rightarrow\Bigr(\prod_{i,j\in\mathcal{I}}\mathcal{O}_{X}(U_{i}\cap U_{j})\Bigr)\rightarrow\ldots$$
is a bounded acyclic complex of flats so it is stable after tensoring with any $E$.
\end{proof}

\begin{defn}
A sheaf of the form $U\mapsto\mathcal{O}_{X}(U)\hat{\otimes} E$ is called a \textit{bornologically free} $\mathcal{O}_{X}$-\textit{module}. If $E$ is $\aleph_{1}$-metrisable, then it is called an $\aleph_{1}$-bornologically free $\mathcal{O}_{X}$-module.
\end{defn}

\begin{defn}
Let $X$ be an analytic space. An object $\mathcal{F}\in{}_{\mathcal{O}_{X}}\mathrm{Mod}$ is said to be \textit{analytically quasi-coherent} if for every subdomain $\mathcal{M}(B)\subset\mathcal{M}(A)$, the map
$$\mathcal{O}_{X}(U)\hat{\otimes}^{\mathbb{L}}_{\mathcal{O}_{X}(X)}\mathcal{F}(X)\rightarrow\mathcal{F}(U)$$
is an equivalence. The full subcategory of ${}_{\mathcal{O}_{X}}\mathrm{Mod}$ consisting of analytic quasi-coherent sheaves is denoted $\mathrm{QCoh}(X)$. We also denote by $\mathrm{QCoh}^{\aleph_{1}}(X)$ the full subcategory of $\mathrm{QCoh}(X)$ consisting of those sheaves $\mathcal{F}$ such that for any  subdomain $U\subset X$, $\mathcal{F}(U)$ is an $\aleph_{1}$-metrisable bornological space. 
\end{defn}

Let us give a more algebraic definition, following \cite{RR}, \cite{koren}.

\begin{defn}
Let $A$ be an analytic algebra. An $A$-module $M$ is said to be an $RR$-\textit{quasi-coherent} if for any map $A\rightarrow B$ corresponding to a subdomain embedding, $B\hat{\otimes}^{\mathbb{L}}_{A}M\rightarrow B\hat{\otimes}_{A}M$ is an equivalence. The full subcategory of $RR$-quasi-coherent modules is denoted ${}_{A}\mathrm{Mod}^{RR}$. The full subcategory of  ${}_{A}\mathrm{Mod}^{RR}$ consisting of modules whose underlying bornological space is $\aleph_{1}$-metrisable is denoted ${}_{A}\mathrm{Mod}^{RR, \aleph_{1}}$
\end{defn}

It is essentially tautological that the adjunction
$$\adj{\mathcal{O}_{\mathcal{M}(A)}\hat{\otimes}_{A}(-)}{{}_{A}\mathrm{Mod}}{{}_{\mathcal{O}_{\mathcal{M}(A)}}\mathrm{Mod}}{\Gamma_{A}(\mathcal{M}(A),-)}$$
restricts to equivalences
$$\adj{\mathcal{O}_{\mathcal{M}(A)}\hat{\otimes}_{A}(-)}{{}_{A}\mathrm{Mod}^{RR}}{\mathrm{QCoh}(X)}{\Gamma_{A}(\mathcal{M}(A),-)}$$
$$\adj{\mathcal{O}_{\mathcal{M}(A)}\hat{\otimes}_{A}(-)}{{}_{A}\mathrm{Mod}^{RR,\aleph_{1}}}{\mathrm{QCoh}^{\aleph_{1}}(X)}{\Gamma_{A}(\mathcal{M}(A),-)}$$
For posterity, we will record this fact as part of a large theorem towards the end of this section.

\begin{rem}
A sheaf $\mathcal{F}$ is analytically quasi-coherent precisely if the bar complex
$$B_{\bullet}^{\mathcal{O}_{X}(X)}(\mathcal{O}_{X}(U),\mathcal{F}(X))\rightarrow\mathcal{F}(U)$$
is a resolution of $\mathcal{F}(U)$ for each subdomain $U$. 
\end{rem}

Using a spectral sequence argument, the following fact is obvious, but useful. 

\begin{prop}\label{prop:flatRR+}
Let $A$ be an analytic algebra, and let $M_{\bullet}\in Ch_{+}({}_{A}\mathrm{Mod}^{RR})$. Then for any subdomain $\mathcal{M}(B)\subset\mathcal{M}(A)$, the map
$$B\hat{\otimes}^{\mathbb{L}}_{A}M_{\bullet}\rightarrow B\hat{\otimes}_{A}M_{\bullet}$$
is an equivalence.
\end{prop}

\begin{lem}\label{Lem:L1flatacyclic}
In the dagger Stein, dagger affinoid, and affinoid cases (the latter for $k$ non-Archimedean) with the finite $G$-topology, a bornologically free $\mathcal{O}_{X}$-module on an analytic space $X$ is analytically quasi-coherent and acyclic.

In the dagger Stein case with the countable $G$-topology a flat $\aleph_{1}$-metrisable bornologically free $\mathcal{O}_{X}$-module on an analytic space $X$ is analytically quasi-coherent and acyclic. 
\end{lem}

\begin{proof}
The bar complex
$$B_{\bullet}^{\mathcal{O}_{X}(X)}(\mathcal{O}_{X}(U),\mathcal{O}_{X}(X))\rightarrow\mathcal{O}_{X}(U)$$
 is an acyclic sequence of flat modules, and so after tensoring with $E$,
$$B_{\bullet}^{\mathcal{O}_{X}(X)}(\mathcal{O}_{X}(U),\mathcal{O}_{X}(X)\hat{\otimes }E\rightarrow\mathcal{O}_{X}(U)\hat{\otimes} E$$
is still acyclic. 

\end{proof}

\begin{rem}
Note that for dagger affinoids we have shown that for any flat bornological algebra, $M\hat{\otimes} \mathcal{O}_{X}$ is acyclic.
\end{rem}

\begin{prop}
    Let $R$ be a Banach ring, and let $E$ be a $\aleph_{1}$-metrisable object of $\mathrm{Ind(Ban_{R})}$. Then there is a strongly flat resolution $P_{\bullet}\rightarrow E$ with each $P_{n}$ being $\aleph_{1}$-metrisable.  
\end{prop}

\begin{proof}
For $V\in\mathrm{Ban_{R}}$ we can functorially construct a resolution
$$P_{\bullet}(V)\rightarrow V$$
where each $P_{n}(V)$ is a projective Banach space. Write $E=\colim_{\mathcal{I}}E_{i}$ with $\mathcal{I}$ $\aleph_{1}$-filtered and each $E_{i}$ Banach. Then define 
$$P_{\bullet}(E)\defeq \underset{\mathcal{I}}{\colim} P_{\bullet}(E_{i})$$
\end{proof}

\begin{cor}
Let $X$ be a dagger Stein space, and $\mathcal{F}$ an $\aleph_{1}$-metrisable analytic quasi-coherent sheaf on $X$. Then there is a resolution of $\mathcal{F}$ by flat, acyclic $\aleph_{1}$-bornologically free $\mathcal{O}_{X}$-modules. If $X$ is a dagger Stein, affinoid, or dagger affinoid and we are working with the finite $G$-topology, then any analytic quasi-coherent sheaf admits a resolution by flat, acyclic, bornologically free $\mathcal{O}_{X}$-modules.
\end{cor}

\begin{proof}
Choose a flat resolution
$$P_{\bullet}\rightarrow\mathcal{F}(X)$$
with $P_{n}$ being flat, and being $\aleph_{1}$-metrisable in the Stein case. Then
$$\mathcal{O}_{X}\hat{\otimes}P_{\bullet}\rightarrow\mathcal{F}$$
is a resolution.

\end{proof}

\begin{lem}\label{lem:steinamplitude}
Let $X$ be an analytic space of finite dimension $n$, $\mathpzc{E}$ an elementary exact category, and $\mathcal{F}\in\mathpzc{Shv}(X,\mathpzc{E})$ a sheaf. Then there is an $n$ such that for $m<n$, $LH_{m}\mathbb{R}\Gamma(X,\mathcal{F})=0$.
\end{lem}

\begin{proof}
Since $\mathpzc{E}$ is elementary, we may in fact assume that it is the (abelian) category of modules over a ring. For $k=\mathbb{C}$ the claim is true for any finite dimensional Stein space since such spaces have finite covering dimension. For $k$ non-Archimedean and $X$ an affinoid or dagger affinoid, it follows from \cite{MR1386046}. Now suppose $X=\bigcup_{n\in\mathbb{N}}U_{n}$ is a dagger Stein space over non-Archimedean $k$. Again we write
$$\mathbb{R}\mathrm{lim}_{n}\mathbb{R}\Gamma(U_{n},(i_{n})_{*}(i_{n})^{*}\mathcal{F})$$
where $i_{n}:U_{n}\rightarrow X$ is the inclusion.  Now $\mathbb{R}\Gamma(U_{n},(i_{n})_{*}(i_{n})^{*}\mathcal{F})$ vanishes in degrees $< n$, and $\underset{n}\Rlim$ has amplitude $[0,-1]$, so $\mathbb{R}\Gamma(X,\mathcal{F})$ vanishes in degrees $< n-1$.
 
\end{proof}

\begin{cor}\label{cor:qcohacyclic}
Let $X$ be a dagger Stein space, let $U\subset X$ be a rational subdomain, and $\mathcal{F}$ an $\aleph_{1}$-metrisable analytic quasi-coherent sheaf. Then $\mathcal{F}|_{U}$ is acyclic as a sheaf. In the affinoid and dagger affinoid cases this is also true for any analytic quasi-coherent sheaf.
\end{cor}

\begin{proof}
Pick an acyclic topologically free resolution $F_{\bullet}\rightarrow\mathcal{F}$. Consider the exact sequences
$$0\rightarrow \ker d_{n}\rightarrow F_{n}\rightarrow \ker d_{n-1}\rightarrow 0$$
By the long exact sequence on cohomology, and using that $RH_{m}\Gamma_{A}(U,F_{m})=0$ for all $m\ge 1$, we find that 
$$LH_{m+k}\mathbb{R}\Gamma(U, \ker d_{n+k})\cong LH_{m}\mathbb{R}\Gamma(U, \ker d_{n})$$
By Lemma \ref{lem:steinamplitude}, we therefore find that all sheaves $\ker d_{n}$, in particular, $\mathcal{F}$, are acyclic on $U$. 
\end{proof}

An easy computation using the long exact $Tor$ sequence gives the following. 

\begin{lem}
Let 
$$0\rightarrow\mathcal{F}\rightarrow\mathcal{G}\rightarrow\mathcal{H}\rightarrow 0$$
be an exact sequence of $\mathcal{O}_{X}$-modules. If 
$$0\rightarrow\mathcal{F}(X)\rightarrow\mathcal{G}(X)\rightarrow\mathcal{H}(X)\rightarrow0$$
is exact, and any two of the sheaves are quasi-coherent, then so is the third.
\end{lem}

Using the long exact sequence for $\mathbb{R}\Gamma(U,-)$, and Corollary \ref{cor:qcohacyclic}, we get the following.

\begin{cor}
The category $\mathrm{QCoh}^{\aleph_{1}}$ is an extension closed subcategory of ${}_{\mathcal{O}_{X}}\mathrm{Mod}$. For the dagger Stein, affinoid, and dagger affinoid cases with the fintie $G$-topology, $\mathrm{QCoh}$ is an extension closed subcategory of ${}_{\mathcal{O}_{X}}\mathrm{Mod}$. 
\end{cor}

\begin{thm}\label{thm:equivcatsKKR}
Let $X$ be an analytic space. 
\begin{enumerate}
\item
The functor $\Gamma_{A}:\mathrm{QCoh}\rightarrow{}_{\mathcal{O}_{X}(X)}\mathrm{Mod}$ is fully faithful and reflects exactness, Its essential image is the category of RR-quasi-coherent modules. 
\item
If $X$ is a dagger Stein space the functor $\Gamma_{\mathcal{O}_{X}(X)}:\mathrm{QCoh}^{\aleph_{1}}\rightarrow{}_{\mathcal{O}_{X}(X)}\mathrm{Mod}$ is exact, and its essential image is the category of $\aleph_{1}$-RR quasicoherent modules. 
\item
If $X$ is a dagger Stein, affinoid or dagger affinoid space and we are working with the finite $G$-topology, the functor $\Gamma_{\mathcal{O}_{X}(X)}:\mathrm{QCoh}\rightarrow{}_{\mathcal{O}_{X}(X)}\mathrm{Mod}$ is exact.
\item
Let $X$ be a dagger Stein space. For $\mathcal{F}_{\bullet}\in Ch_{b}(\mathrm{QCoh}^{\aleph_{1}})$ the maps
$$\Gamma_{\mathcal{O}_{X}(X)}(X,\mathcal{F}_{\bullet})\rightarrow\mathbb{R}\Gamma_{\mathcal{O}_{X}(X)}(X,\mathcal{F}_{\bullet})$$
$$\mathcal{O}_{X}\hat{\otimes}^{\mathbb{L}}\Gamma_{\mathcal{O}_{X}(X)}(X,\mathcal{F}_{\bullet})\rightarrow\mathcal{O}_{X}\hat{\otimes}\Gamma_{\mathcal{O}_{X}(X)}(X,\mathcal{F}_{\bullet})$$
are equivalences.
\item
Let $X$ be a dagger Stein, affinoid, or dagger affinoid space equipped with the finite $G$-topology. For $\mathcal{F}_{\bullet}\in\Ch_{b}(\mathrm{QCoh})$ the maps
$$\Gamma_{\mathcal{O}_{X}(X)}(X,\mathcal{F}_{\bullet})\rightarrow\mathcal{F}_{\bullet}\rightarrow\mathbb{R}\Gamma_{\mathcal{O}_{X}(X)}(X,\mathcal{F}_{\bullet})$$
$$\mathcal{O}_{X}\hat{\otimes}^{\mathbb{L}}\Gamma_{\mathcal{O}_{X}(X)}(X,\mathcal{F}_{\bullet})\rightarrow\mathcal{O}_{X}\hat{\otimes}\Gamma_{\mathcal{O}_{X}(X)}(X,\mathcal{F}_{\bullet})$$
are equivalences.
\end{enumerate}
\end{thm}

\begin{proof}
\begin{enumerate}
\item
The fact that the functor is fully faithful, and that the essential image is the category of $RR$-quasicoherent modules is an immediate consequence of the fact that, by definition, for any Stein subdomain $U\subset X$, $\mathcal{O}_{X}(U)\hat{\otimes}^{\mathbb{L}}_{\mathcal{O}_{X}(X)}\mathcal{F}(X)\rightarrow\mathcal{F}(U)$ is an isomorphism. The fact that \[\Gamma_{\mathcal{O}_{X}(X)}(X,-):\mathrm{QCoh}^{\aleph_{1}}\rightarrow{}_{\mathcal{O}_{X}(X)}\mathrm{Mod}\] is exact  is an immediate consequence of Corollary \ref{cor:qcohacyclic}.
\item
Again this follows from Corollary \ref{cor:qcohacyclic}.
\item
First suppose the sequence is bounded. By an easy spectral sequence argument, and Part (1), the map \[\Gamma_{\mathcal{O}_{X}(X)}(X,\mathcal{F}_{\bullet})\rightarrow\mathbb{R}\Gamma_{\mathcal{O}_{X}(X)}(X,\mathcal{F}_{\bullet})\] is an equivalence. Since $\Gamma_{\mathcal{O}_{X}(X)}(X,\mathcal{F}_{\bullet})$ is a bounded sequence, by Proposition \ref{prop:flatRR+} we have 
\begin{equation*}
    \begin{split}
& \mathcal{O}_{X}(U)\hat{\otimes}^{\mathbb{L}}_{\mathcal{O}_{X}(X)}\mathbb{R}\Gamma_{\mathcal{O}_{X}(X)}(X,\mathcal{F}_{\bullet}) \\ & \cong \mathcal{O}_{X}(U)\hat{\otimes}_{\mathcal{O}_{X}(X)}\Gamma_{\mathcal{O}_{X}(X)}(X,\mathcal{F}_{\bullet}) \\ & \cong\mathcal{F}_{\bullet}(U)
    \end{split}
\end{equation*}
\item
This is identical to Part (3).
\end{enumerate}
\end{proof}

\begin{defn}
Let $X$ be an analytic space
\begin{enumerate}
\item
An object $M_{\bullet}\in Ch({}_{\mathcal{O}_{X}(X)}\mathrm{Mod})$ is said to be $K-RR$-\textit{quasi-coherent} if the homotopical unit map
$$M_{\bullet}\rightarrow\mathbb{R}\Gamma(X,\mathcal{O}_{X}(-)\hat{\otimes}^{\mathbb{L}}_{\mathcal{O}_{X}(X)}M_{\bullet})$$
is an equivalence. Denote the category of all $K-RR$-quasi-coherent modules by $Ch({}_{\mathcal{O}_{X}(X)}\mathrm{Mod})^{K-RR}$.
\item
An object $\mathcal{F}_{\bullet}\in Ch({}_{\mathcal{O}_{X}}\mathrm{Mod})$ is said to be $K$-\textit{quasi-coherent} if the homotopical counit
$$\mathcal{O}_{X}\hat{\otimes}^{\mathbb{L}}\mathbb{R}\Gamma_{\mathcal{O}_{X}(X)}(X,\mathcal{F}_{\bullet})\rightarrow\mathcal{F}_{\bullet}$$
is an equivalence. Denote the category of all $K-RR$-quasi-coherent modules by $Ch({}_{\mathcal{O}_{X}}\mathrm{Mod})^{K,\mathrm{QCoh}}$ .
\end{enumerate}
\end{defn}

Recall the adjunction

$$\adj{\mathcal{O}_{X}\hat{\otimes}(-)}{\Ch({}_{\mathcal{O}_{X}(X)}\Mod)}{\Ch({}_{\mathcal{O}_{X}}\Mod)}{\Gamma(X,-)}$$
of $(\infty,1)$-categories. The categories $\Ch({}_{\mathcal{O}_{X}(X)}\mathrm{Mod})^{K-RR}$ and $\Ch({}_{\mathcal{O}_{X}}\mathrm{Mod})^{K,\mathrm{QCoh}}$ present the $(\infty,1)$-category of fixed points of this adjunction in $\Ch({}_{\mathcal{O}_{X}(X)}\Mod)$ and $\Ch({}_{\mathcal{O}_{X}}\Mod)$ respectively. In particular, they present equivalent $(\infty,1)$-categories.

By the above we have the following.

\begin{prop}\label{prop:KKRbounded}
Let $X$ be a dagger Stein space, and let $F$ be a bounded above complex of $\aleph_{1}-RR$-quasicoherent $\mathcal{O}_{X}(X)$-modules. Then $F$ is $K-RR$-quasicoherent. If $X$ is  dagger Stein, dagger affinoid or an affinoid and is equipped with the finite $G$-topology, then this is true for any bounded above complex of bornological spaces.
\end{prop}



\begin{lem}\label{lem:glueRR}
The category $\Ch({}_{\mathcal{O}_{X}(X)}\mathrm{Mod})^{K-RR}$ is closed under finite homotopy limits, and countable homotopy limits in the affinoid and dagger Stein cases.
\end{lem}

\begin{proof}
In the affinoid and dagger affinoid cases this follows since both $\mathbb{R}\Gamma(X,-)$ and $\mathcal{O}_{X}(U)\hat{\otimes}^{\mathbb{L}}_{\mathcal{O}_{X}(X)}M_{\bullet}$ commute with finite limits. In the dagger Stein case the claim by similar reasoning, now using that fact that any localisation $\mathcal{O}_{X}(X)\rightarrow\mathcal{O}_{X}(U)$ is $\kappa$-filtered. 
\end{proof}

\subsection{Classical Quasi-coherent Sheaves, Descent, and an Example of Gabber}

Let $X$ be an analytic space. One might naïvely define a quasi-coherent sheaf on $X$ to be an $\mathcal{O}_{X}$-module $\mathcal{F}$ which locally has a presentation. 
$$\bigoplus_{i\in\mathcal{I}}\mathcal{O}_{X}|_{U}\rightarrow \bigoplus_{i\in\mathcal{J}}\mathcal{O}_{X}|_{U}\rightarrow\mathcal{F}|_{U}\rightarrow0$$

The following famous example of Gabber \cite{2840594} shows that in the context of analytic geometry over $\mathbb{C}$, with the usual topology, stable descent is the best that we can hope for.  Let $X$ be the open unit disc, and let $x',x''$ be two points of $X$.
 Write 
$$U'=X\setminus\{x'\}, U''=X\setminus\{x''\},\;\;U=U'\cap U''$$Define
$$\mathcal{F'}=\bigoplus_{n\in\mathbb{Z}}\mathcal{O}_{U'}e'_{n},\;\; \mathcal{F''}\defeq\bigoplus_{n\in\mathbb{Z}}\mathcal{O}_{U''}e''_{n}$$\newline
\\
Let $h\in\mathcal{O}_{X}(U)$ be a function with essential singularities at $x',x''$.  Identify $\mathcal{F}'|_{U}$ with $\mathcal{F}''|_{U}$ with $\underset{n\in\mathbb{Z}}{\bigoplus}\mathcal{O}_{X}(U)e_{n}$ by
$$e_{2m}=e'_{2m}|_{U}=e''_{2m}|_{U}+he''_{2m+1}|_{U}$$
$$e_{2m+1}=e''_{2m}|_{U}=e'_{2m}|_{U}+he'_{2m+1}|_{U}$$
These identifications induce isomorphisms of complete bornological $\mathcal{O}_{X}(U)$-modules
$$\mathcal{F}'|_{U}\cong\bigoplus_{n\in\mathbb{Z}}\mathcal{O}_{X}(U)e_{n}\cong\mathcal{F}''|_{U}$$
Let $f\in\mathcal{F}(X)$. $f|_{U'}=\sum_{\textrm{finite}}f'_{n}e'_{n}$, and $f|_{U''}=\sum_{\textrm{finite}}f''_{n}e''_{n}$. Noting that the restriction maps are injective and comparing coefficients one finds that 
$$f|_{U}=\sum_{\textrm{finite}}f_{n}e_{n}$$
with $f_{n}$ analytic at $x'$ for even $n$ and $x''$ for odd $n$, and $f_{n}+hf_{n-1}$ analytic at $x'$ for odd $n$ and $x''$ for even $n$. Let $n_{0}$ be the maximal $n$ with $f_{n}\neq0$. Then $f_{n_{0}}$ and $hf_{n_{0}}$ are both analytic at one of the $x',x''$.  Thus $\frac{hf_{n_{0}}}{f_{n_{0}}}$ is analytic at one of these, which contradicts the fact that $h$ has an essential singularity here. Thus we have $\mathcal{F}(X)\cong 0$. 

Now consider $X$ as a dagger Stein space equipped with the finite $G$-topology. Then $\mathcal{F}'$ and $\mathcal{F}''$ are bornological analytic quasi-coherent sheaves. We denote their global sections by $\mathrm{F}'$ and $\mathrm{F}''$ respectively.
By \v{C}ech descent in the finite cover classical topology, there is a unique object $\mathrm{F}$ of 
$$\mathbf{QCoh}(\mathrm{Spec}(\mathcal{O}_{X}(X)))$$
such that 
$$\mathcal{O}_{X}(U')\otimes^{\mathbb{L}}_{\mathcal{O}_{X}(X)}\mathrm{F}\cong\mathrm{F}'$$
$$\mathcal{O}_{X}(U'')\otimes^{\mathbb{L}}_{\mathcal{O}_{X}(X)}\mathrm{F}\cong\mathrm{F}''$$
$\mathrm{F}$ is moreover equivalent to a bounded complex of $\mathcal{O}_{X}(X)$-modules,
Moreover, this object $F$ lives in the category
$$\Ch({}_{\mathcal{O}_{X}(X)}\mathrm{Mod})^{K-RR}$$
and so may be regarded as a complex of sheaves on $X$ for the finite the finite $G$-topology. In the derived bornological sense Gabber's example also does not pose a problem: quasi-coherent sheaves on analytic spaces glue as long as one allows for the module of global sections to be a complex. Recently, Haohao Liu has put up some preprints on his website: (1) Quasi-coherent GAGA, 
(2) Sheaves with connection on complex tori, 
(3) Fourier-Mukai transform on complex tori, revisited,
(4) Quasi-coherent sheaves on complex analytic spaces. In these articles, he defines and  compares various notions of quasi-coherence in complex analytic geometry including the naive one above. He also points out some mistakes in work of the first author in his thesis work \cite{OrenThesis} and \cite{OrenNC} on quasi-coherent Fourier-Mukai transforms in complex-analytic geometry.

\subsection{Future Application: Bounded Cohomology and Moduli of Representations}

Let $G$ be a compact topological group and $V$ a Banach $R$-module. Denote by $C(G,V)$ the Banach space of continuous functions $G\rightarrow V$. This defines a commutative Hopf monad on $\mathrm{Ban_{R}}$, which extends to a Hopf monad on $\mathrm{Ind(Ban_{R})}$. If $R$ is non-Archimedean then this Hopf monad arises from the commutative Hopf algebra $C(G,R)$.  Note that there is always a map of commutative monads \[C(G,R)\hat{\otimes}_{R}(-)\rightarrow C(G,-).\]

\begin{defn}
Say that a compact group $G$ is $R$-\textit{compact} if the map
$$C(G,R)\hat{\otimes}_{R} C(G,R)\rightarrow C(G,C(G,R))$$
is an isomorphism in $\mathrm{Ban}_{R}$.
\end{defn}

\begin{example}
\begin{enumerate}
\item
If $R$ is a non-Archimedean Banach field then all compact groups are $R$-compact.
\item
If $R$ is any Banach ring, then finite groups (with the discrete topology) are $R$-compact. Moreover any pro-finite group is $R$-compact. 
\end{enumerate}
\end{example}

In the situation that $G$ is $R$-compact, we have that $C(G,R)$ is a commutative Hopf monoid. Thus $G^{R}\defeq \mathrm{Spec}(C(G,R))$ is an affine group stack. We can construct the classifying stack $BG^{R}$. We are interested in the moduli stack of dimension $n$ representations
$$\mathbf{Rep}_{n}(G^{R})\defeq\underline{\mathrm{Map}}(BG^{R},BGL_{n})$$

In future work we will show the following.

\begin{lem}
The stack $\mathbf{Rep}_{n}(G^{R})$ has a (potentially unbounded) cotangent complex. If $G$ is profinite and $R$ contains $\mathbb{Q}$ then the cotangent complex is concentrated in non-positive degrees.  
\end{lem}

This will involve in particular that the map $BG^{R}\rightarrow\mathrm{Spec}(R)$ is strongly $\mathbf{QCoh}$-cohomologically proper (Definition \ref{defn:cohprop}). We will then show the following.

\begin{lem}
Let $R$ be a non-Archimedean commutative Banach $\mathbb{Z}_{p}$-algebra. Let $K$ be a non-Archimedean valued field, and let $G_{K}$ be its absolute Galois group. Then for any $n\in\mathbb{N}_{0}$ the stack $\mathbf{Rep}_{n}(G_{K}^{R})$ has a bounded above cotangent complex.
\end{lem}

Then, using the forthoming representability theorem of Savage, we will establish the geometricity of the moduli of Galois representations.

\appendix

\chapter{Some Technical Details}\label{STD}

\section{Simplicial Objects and Sifted Cocompletions}

In this section we record some technical facts regarding model structures on categories of simplicial objects, functors between them, and the relation to free sifted cocompletions of $(\infty,1)$-categories.

\subsection{Model Structures on Simplicial Objects}

 Let $\mathrm{C}$ be a complete and cocomplete, category and $\mathcal{Z}=\{Z_{i}\}_{i\in\mathcal{I}}$ a set of tiny objects of $\mathrm{C}$ closed under finite coproducts. Consider the simplicial model structure on $\mathrm{sC}$ whereby a map $f:X\rightarrow Y$ is a weak equivalence (resp. a fibration) if for any $i\in\mathcal{I}$, the map $\mathrm{Hom}(Z_{i},f):\mathrm{Hom}(Z_{i},X)\rightarrow\mathrm{Hom}(Z_{i},Y)$ is a weak equivalence (resp. a fibration) of simplicial sets. By \cite{Goerss-Jardine} Section II.5 With this model structure $\mathrm{sC}$ becomes a combinatorial simplicial model category. $\mathrm{sC}$ is tensored over simplicial sets as a model category. A set of generating cofibrations is given by 
$$\{\partial\Delta[n]\otimes Z_{i}\rightarrow\Delta[n]\otimes Z_{i}\}_{n\in\mathbb{N}_{0},i\in\mathcal{I}}$$
A set of generating acyclic cofibrations is given by
$$\{\Lambda^{i}[n]\otimes Z_{i}\rightarrow\Delta[n]\otimes Z_{i}\}$$

The small object argument implies that there is a cofibration-acyclic fibration factorisation functor , where cofibrant replacement functor $Q$ is such that $Q(X)_{n}$ is a coproduct of objects in $\{Z_{i}\}_{i\in\mathcal{I}}$.

 Let $\mathcal{G}$ be a small category with finite coproducts and consider the free sifted cocompletion $\mathcal{P}_{\Sigma}(\mathcal{G})$. Consider the category $\mathrm{s}\mathcal{P}_{\Sigma}(\mathcal{G})$, and equip it with the model structure corresponding to taking $\mathcal{Z}=\mathcal{G}$.

In particular, there is a natural functor 
$$L:\mathcal{P}_{\Sigma}(\mathcal{Z})\rightarrow\mathrm{C}$$
and, extending to simplicial objects, a functor, which we also denote by $L$
$$L:\mathrm{s}\mathcal{P}_{\Sigma}(\mathcal{Z})\rightarrow\mathrm{sC}$$
This functor commutes with colimits and is therefore a left adjoint with right adjoint $R$. 

\begin{lem}\label{lem:rigidification}
Suppose every cofibrant object of $\mathrm{sC}$ is fibrant. Then the adjunction
$$\adj{L}{\mathrm{s}\mathcal{P}_{\Sigma}(\mathcal{Z})}{\mathrm{sC}}{R}$$
is a Quillen equivalence.
\end{lem}

\begin{proof}
By the description of the generating (acyclic) cofibrations in each category, $L$ is left Quillen. Each $Z_{i}\in\mathcal{Z}$ is tiny in both $\mathcal{P}_{\Sigma}(\mathcal{Z})$ and $\mathrm{sC}$. Thus $L$ is fully faithful when restricted to coproducts of objects of $\mathcal{Z}$. Since any object of $\mathrm{sC}$ is equivalent to an object which in each degree is a coproduct of objects of $\mathcal{Z}$, the derived functor $\mathbb{L}L$ is homotopy essentially surjective. We claim that it also induces equivalences of mapping spaces. Let $X_{\bullet}, Y_{\bullet}$ be fibrant-cofibrant objects of $\mathrm{s}\mathcal{P}_{\Sigma}(\mathcal{Z})$. Without loss of generality we may assume that each $X_{n}$ and each $Y_{n}$ is a disjoint union of objects of $\mathcal{Z}$. The mapping space in $\mathrm{s}\mathcal{P}_{\Sigma}(\mathcal{Z})$ is then given by 
$$\underline{\mathrm{Hom}}(X_{\bullet},Y_{\bullet})$$
On the other hand since $L(X_{\bullet})$ is cofibrant, and $L(Y_{\bullet})$ is cofibrant and therefore also fibrant by assumption,  the mapping space in $\mathrm{sC}$ is also given by
$$\underline{\mathrm{Hom}}(L(X_{\bullet}),L(Y_{\bullet}))$$
Again because $L$ is fully faithful when restricted to coproducts of objects of $\mathcal{Z}$ we have that 
$$L:\underline{\mathrm{Hom}}(X_{\bullet},Y_{\bullet})\rightarrow\underline{\mathrm{Hom}}(L(X_{\bullet}),L(Y_{\bullet}))$$
is an isomorphism.
\end{proof}

\subsection{Functors Between Categories of Simplicial Objects}

Many functors between categories of simplicial objects we consider in this work are not Quillen. However we still need to derived them. Let $\mathcal{G}$ be as in previous section, let $\mathpzc{C}$ be a combinatorial model category, and 
$$F:\mathcal{G}\rightarrow\mathrm{C}$$ 
be a functor. This extends by sifted colimits to a functor, which we also denote by $F$,
$$F:\mathcal{P}_{\Sigma}(\mathcal{G})\rightarrow\mathrm{C}$$
and then to a functor of simplicial objects, which we again denote by $F$,
$$F:\mathrm{s}\mathcal{P}_{\Sigma}(\mathcal{G})\rightarrow\mathrm{s}\mathrm{C}$$
and then a functor
$$\overline{F}:\mathrm{s}\mathcal{P}_{\Sigma}(\mathcal{G})\rightarrow\mathpzc{C}$$
by taking geometric realisations. 
Suppose that $F$ preserves weak equivalences between fibrant-cofibrant objects. Then $F$ may be derived to a functor
$$\mathrm{L}F:\mathrm{L^{H}}(\mathrm{s}\mathcal{P}_{\Sigma}(\mathcal{G}))\rightarrow\mathrm{L^{H}}(\mathrm{s}\mathrm{C})$$

Let us give some checkable conditions under which a functor $F:\mathcal{M}\rightarrow\mathcal{N}$ between simplicial model categories preserves equivalences between fibrant-cofibrant objects. 

\begin{lem}[c.f. \cite{borisov2017quasi} Proposition 27]\label{lem:funcsimp}
Let $F:\mathcal{M}\rightarrow\mathcal{N}$ be a functor between model categories with $\mathcal{M}$ simplicial. Let $\otimes:\mathrm{sSet}\times\mathcal{M}\rightarrow\mathcal{M}$ denote the tensoring of $\mathcal{M}$ over $\mathrm{sSet}$. Suppose that for any fibrant-cofibrant object $A$ the map
$$F(\sigma):F(A\otimes_{\mathcal{N}}\Delta[1])\rightarrow F(A)$$
is an equivalence. Then $F$ preserves weak equivalences between fibrant-cofibrant objects.
\end{lem}

\begin{proof}
%
%
More generally let $f,f':A\rightarrow A'$ be homotopic maps between fibrant-cofibrant objects. There is $\phi:A\otimes\Delta[1]\rightarrow A'$ such that $f=\phi\circ i_{1},f'=\phi\circ i_{2}$, where $i_{1},i_{2}:A\rightarrow A\otimes\Delta[1]$ are the canonical inclusions. We have $F(\sigma)\circ F(i_{1})=\mathrm{Id}=F(\sigma)\circ F(i_{2})$. Since $F(\sigma)$ is an equivalence by assumption $F(i_{1})=F(i_{2})$ in the homotopy category. Thus $F(f)=F(f')$. Now if $f:A\rightarrow B$ is an equivalence between fibrant-cofibrant objects in $\mathcal{M}$ then it has a homotopy inverse $f':B\rightarrow A$. By the above argument $F(f)$ and $F(f')$ are inverse to each other in $\mathrm{Ho}(\mathcal{N})$.
%
\end{proof}

Let us briefly specialise to the case that $\mathpzc{C}=\mathrm{sC}$ for some complete and cocomplete category $\mathrm{C}$ with model structure determined by  a set $\mathcal{Z}=\{\mathcal{Z}_{i}\}_{i\in\mathcal{I}}$ of tiny objects of $\mathrm{C}$ closed under finite coproducts. We wish to determine in this case when the functor $\mathbb{L}F$ is fully faithful. 
Let $\mathcal{G}^{\coprod}\subset\mathcal{P}_{\Sigma}(\mathcal{G})$ denote the full subcategory consisting of aribtrary coproducts of objects of $\mathcal{G}$.

\begin{prop}
 $F:\mathcal{G}\rightarrow\mathrm{C}$ be a functor such that the induced functor
$$F:\mathcal{G}^{\coprod}\rightarrow\mathrm{C}$$
is fully faithful. Then the functor
$$\mathrm{s}\mathcal{P}_{\Sigma}(\mathcal{G})\rightarrow\mathrm{sC}$$
is fully faithful on cofibrant objects.
\end{prop}

\begin{proof}
This is clear since cofibrant objects are degree-wise free (c.f. Proposition 27 in \cite{borisov2017quasi}).
\end{proof}

\begin{lem}
Suppose that 
\begin{enumerate}
\item
$\mathcal{G}^{\coprod}\rightarrow\mathrm{C}$ is fully faithful.
\item
$F:\mathrm{s}\mathcal{P}_{\Sigma}(\mathcal{G})\rightarrow\mathrm{sC}$ preserves weak equivalences between fibrant-cofibrant objects.
\item
For $X$ and $Y$ fibrant-cofibrant objects of $\mathrm{s}\mathcal{P}_{\Sigma}(\mathcal{G})$ the map \[\mathrm{\underline{\Hom}}_{\mathrm{sC}}(F(X),F(Y))\rightarrow\Map(F(X),F(Y))\] is an equivalence.
\end{enumerate}
Then the functor
$$\mathbf{F}:\mathrm{L}^{H}(\mathrm{s}\mathcal{P}_{\Sigma}(\mathcal{G}))\rightarrow\mathrm{L}^{H}(\mathrm{sC})$$
is fully faithful.
\end{lem}

\begin{proof}
 By assumption, for $X$ and $Y$ fibrant-cofibrant objects of $\mathrm{s}\mathcal{P}_{\Sigma}(\mathcal{G})$ the map $\mathrm{\underline{\Hom}}_{\mathrm{sC}}(F(X),F(Y))\rightarrow\Map(F(X),F(Y))$ is an equivalence. Thus it suffices to prove that the map
$$\mathrm{\underline{\Hom}}_{\mathrm{s}\mathcal{P}_{\Sigma}(\mathcal{G})}(X,Y)\rightarrow\mathrm{\underline{\Hom}}_{\mathrm{sC}}(F(X),F(Y))$$
is an isomorphism of simplicial sets for any fibrant-cofibrant objects. Every (fibrant) object of $\mathrm{s}\mathcal{P}_{\Sigma}(\mathcal{G})$ has a cofibrant (fibrant-cofibrant) resolution by an object $X_{\bullet}$ such that each $X_{n}$ is in $\mathcal{G}^{\coprod}$. Now the result is clear by assumption.
\end{proof}

\subsection{Rigidification of Sifted Cocompletions}
Let $\mathpzc{C}$ be a combinatorial simplicial model category and let $\mathcal{G}$ be a small category with finite products. Suppose further that geometric realisation of simplicial objects in $\mathpzc{C}$ commutes with finite products, and that filtered colimits of weak equivalences are weak equivalences in $\mathpzc{C}$. By \cite{MR1852091} Proposition A.5 we may pick a set $\mathcal{C}$ of cofibrant objects such that a map $f:C\rightarrow D$ in $\mathpzc{C}$ is a weak equivalence precisely if 
$$\mathrm{R}\underline{\mathrm{Hom}}(c,f):\mathrm{R}\underline{\mathrm{Hom}}(c,C)\rightarrow\mathrm{R}\underline{\mathrm{Hom}}(c,D)$$
is a weak equivalence for any $c\in\mathcal{C}$. Consider the category
$$\mathbf{Fun}(\mathcal{G},\mathpzc{C})$$
of functors equipped with the projective model structure, whereby a natural transformation $\eta:E\rightarrow F$ is a fibration (resp. an acyclic fibration) precisely if for any $g\in\mathcal{G}$, $\eta_{g}:E(g)\rightarrow F(g)$ is a fibration (resp. an acyclic fibration) in $\mathpzc{C}$. This is also a combinatorial model category, and presents the category
$$\mathbf{Fun}(\mathrm{N}(\mathcal{G}),\mathrm{L^{H}}(\mathpzc{C}))$$
of $(\infty,1)$-functors. Consider the left Bousfield localisation of 
$\mathbf{Fun}(\mathcal{G}^{op},\mathpzc{C})$ at those maps of the form
$$\coprod_{i=1}^{n}\mathrm{Hom}(g_{i},-)\otimes c\rightarrow\mathrm{Hom}(\prod_{i=1}^{n}g_{i},-)\otimes c$$
for any $c\in\mathcal{C}$ and any finite collection $\{g_{1},\ldots,g_{n}\}$ of objects of $\mathcal{G}$. 
We denote this model category by
$$\mathrm{Fun}^{h\times}(\mathcal{G}^{op},\mathpzc{C})$$
Let $F:\mathcal{G}^{op}\rightarrow\mathpzc{C}$ be fibrant in $\mathrm{Fun}^{h\times}(\mathcal{G}^{op},\mathpzc{C})$. Then $F$ is fibrant in the projective model structure on functors. Moreover the map
$$\mathrm{R}\underline{\mathrm{Hom}}(\mathrm{Hom}(\prod_{i=1}^{n}g_{i},-)\otimes c,F)\rightarrow\mathrm{R}\underline{\mathrm{Hom}}(\coprod_{i=1}^{n}\mathrm{Hom}(g_{i},-)\otimes c,F)\cong\prod_{i=1}^{n}\mathrm{R}\underline{\mathrm{Hom}}(\mathrm{Hom}(g_{i},-)\otimes c,F)$$
is an equivalence. Now both source objects are cofibrant, so this is just the simplicial mapping space
$$\underline{\mathrm{Hom}}(\prod_{i=1}^{n}g_{i},-)\otimes c,F)\rightarrow\prod_{i=1}^{n}\underline{\mathrm{Hom}}(\mathrm{Hom}(g_{i},-)\otimes c,F)$$
By adjunction this is just 
$$\underline{\mathrm{Hom}}(c,F(\prod_{i=1}^{n}g_{i}))\rightarrow\prod_{i=1}^{n}\underline{\mathrm{Hom}}(c,F(g_{i}))$$
Now again each $F(\prod_{i=1}^{n}g_{i})$ and each $F(g_{i})$ is fibrant, so this is
$$\mathrm{R}\underline{\mathrm{Hom}}(c,F(\prod_{i=1}^{n}g_{i}))\rightarrow\prod_{i=1}^{n}\mathrm{R}\underline{\mathrm{Hom}}(c,F(g_{i}))$$
Since $\mathcal{C}$ is a set of cofibrant generators, we have that $F$ is fibrant precisely if 
$$F(\prod_{i=1}^{n}g_{i})\rightarrow\prod_{i=1}^{n}F(g_{i})$$
is an equivalence for any finite collection $\{g_{1},\ldots,g_{n}\}$ of objects of $\mathcal{G}$

Now since $\mathpzc{C}$ is locally-presentable, the category $\mathrm{Fun}^{\times}(\mathcal{G},\mathpzc{C})$ of strict product-preserving functors is also locally-presentable. Indeed let $\mathcal{D}$ be a set of generators of $\mathpzc{C}$ under filtered colimits. Then  $\mathrm{Fun}^{\times}(\mathcal{G},\mathpzc{C})$ is the strict, $1$-categorical, localisation at maps of the form 
$$\coprod_{i=1}^{n}\mathrm{Hom}(g_{i},-)\otimes d\rightarrow\mathrm{Hom}(\prod_{i=1}^{n}g_{i},-)\otimes d$$
for any $d\in\mathcal{D}$ and any finite collection $\{g_{1},\ldots,g_{n}\}$ of objects of $\mathcal{G}$. There is an inclusion
$$\tilde{J}_{\mathcal{G}}\mathrm{Fun}^{\times}(\mathcal{G},\mathpzc{C})\rightarrow\mathrm{Fun}(\mathcal{G},\mathpzc{C})$$
This functor has a left adjoint $\tilde{K}_{\mathcal{G}}$ by the adjoint functor theorem. Indeed both categories are locally presentable, and the inclusion commutes with all limits and all sifted colimits.

For each $g\in\mathcal{G}$ let $U_{g}:\mathrm{Fun}^{\times}(\mathcal{G},\mathpzc{C})\rightarrow\mathpzc{C}$ denote functor of evaluation at $g$. This functor also has left adjoint $F_{g}:\mathpzc{C}\rightarrow\mathrm{Fun}^{\times}(\mathcal{G},\mathpzc{C})$ by the adjoint functor theorem. Exactly as in \cite{MR2263055} Theorem 4.7, we have the following

\begin{thm}
    There exists a cobminatorial model category structure on $\mathrm{Fun}^{\times}(\mathcal{G},\mathpzc{C})$ such that a map $\eta:E\rightarrow F$ is a fibration (resp. an acyclic fibration) precisely if $\eta_{g}:E(g)\rightarrow F(g)$ is a fibration (resp. an acyclic fibration) for any $g\in\mathcal{G}$.
\end{thm}

By construction the adjunction 
$$\adj{K_{\mathcal{G}}}{\mathrm{Fun}(\mathcal{G},\mathpzc{C})}{\mathrm{Fun}^{\times}(\mathcal{G},\mathpzc{C})}{J_{\mathcal{G}}}$$
is clearly Quillen. Moreover $K_{\mathcal{G}}$ sends maps of the form
$$\coprod_{i=1}^{n}\mathrm{Hom}(g_{i},-)\otimes c\rightarrow\mathrm{Hom}(\prod_{i=1}^{n}g_{i},-)\otimes c$$
to isomorphisms in $\mathrm{Fun}^{\times}(\mathcal{G},\mathpzc{C})$. It follows that under the assumptions stated above we have a Quillen adjunction
$$\adj{K_{\mathcal{G}}}{\mathrm{Fun}^{\h\times}(\mathcal{G},\mathpzc{C})}{\mathrm{Fun}^{\times}(\mathcal{G},\mathpzc{C})}{J_{\mathcal{G}}}$$

\begin{thm}\label{thm:rigidification}
    $$\adj{K_{\mathcal{G}}}{\mathrm{Fun}^{h\times}(\mathcal{G},\mathpzc{C})}{\mathrm{Fun}^{\times}(\mathcal{G},\mathpzc{C})}{J_{\mathcal{G}}}$$
\end{thm}

\begin{proof}[Proof Sketch]
    The proof works very similarly to \cite{MR2263055} Theorem 5.1/ \cite{MR1923968} Theorem 6.5. Crucially as in \cite{MR1923968} Theorem 6.5 it suffices to establish the analogue of \cite{MR1923968} Lemma 6.4, that for any cofibrant object $F\in\mathrm{Fun}^{h\times}(\mathcal{G},\mathpzc{C})$, the unit
    $$\eta_{X}:X\rightarrow J_{\mathcal{G}}K_{\mathcal{G}}X$$
    is a local equivalence. Now thanks to the assumption that geometric realisations of simplicial objects in $\mathpzc{C}$ commute with finite products, Lemma 6.1 of \cite{MR1923968} works in this greater generality, and Lemma 6.2 also follows. Now $K_{\mathcal{G}}$ preserves equivalences between cofibrant objects of $J_{\mathcal{G}}$ preserves all equivalences. Thus we may replace $X$ by any convenient cofibrant object. As a combinatorial model category $\mathpzc{C}$ has a presentations as a Bousfield localisation of a category of presheaves on some small category by \cite{MR1870516}.Thus we may in fact assume that $X$ is a geometric realisation of coproducts of objects of the form $\coprod_{i\in\mathcal{I}}\mathrm{Hom}(g_{i},-)\otimes c_{i}$ for some possibly infinite set $\mathcal{I}$. Now the analogue of Lemma 6.4 can be proved exactly as in \cite{MR1923968}.
\end{proof}

\subsection{Rigidification of Functors}

Let $\mathcal{G}$ be a small $1$-category. Consider the $(\infty,1)$-categorical free-sifted cocompletion
$$\mathbf{P}_{\Sigma}(\mathrm{N}(\mathcal{G}))\defeq\mathbf{Fun}^{\times}(\mathrm{N}(\mathcal{G})^{op},\mathbf{sSet})$$
By the previous subsection, this is presented by the $(\infty,1)$-category $\mathrm{s}\mathcal{P}_{\Sigma}(\mathcal{G})$ of simplicial objects in the $1$-categorical free sifted cocompletion of $\mathcal{G}$. 

Now let $\mathpzc{C}$ be a combinatorial model category, and let $F:\mathcal{G}\rightarrow\mathpzc{C}$ be a functor. Suppose that the induced functor

$$\overline{F}:\mathrm{s}\mathcal{P}_{\Sigma}(\mathcal{G})\rightarrow\mathpzc{C}$$
preserves weak equivalences between fibrant-cofibrant objects, so that we get a functor
$$\overline{\mathbf{F}}:\mathbf{P}_{\Sigma}(\mathrm{N}(\mathcal{G}))\rightarrow\mathrm{L^{H}}(\mathpzc{C})$$
Consider the restriction of this functor
$$\overline{\mathbf{F}}|_{\mathrm{N}(\mathcal{G})}:\mathrm{N}(\mathcal{G})\rightarrow\mathrm{L^{H}}(\mathpzc{C})$$
This functor extends by sifted colimits to a functor
$$\overline{\mathbf{F}}':\mathbf{P}_{\Sigma}(\mathrm{N}(\mathcal{G}))\rightarrow\mathrm{L^{H}}(\mathpzc{C})$$
We wish to understand when there is a natural equivalence 
$$\overline{\mathbf{F}}\cong\overline{\mathbf{F}}'$$
Since $\overline{\mathbf{F}}$ and $\overline{\mathbf{F}}'$ agree on $\mathrm{N}(\mathcal{G})$ it suffices to understand when $\overline{\mathbf{F}}$ commutes with sifted colimits. 

\begin{lem}\label{lem:rigidfunc}
Suppose that 
\begin{enumerate}
\item
$$\overline{F}:\mathrm{s}\mathcal{P}_{\Sigma}(\mathcal{G})\rightarrow\mathpzc{C}$$
preserves weak equivalences between cofibrant objects. 
\item 
filtered colimits of maps of the form $F(i)$ for $i$ a cofibration in $\mathrm{s}\mathcal{P}_{\Sigma}(\mathcal{G})$ presented homotopy colimits
\item
The geometric realisation of a simplicial object in $\mathrm{s}\mathpzc{C}$ of the form $F(X_{\bullet})$ for $X_{\bullet}$ a Reedy cofibrant simplicial object in $\mathrm{s}\mathcal{P}_{\Sigma}(\mathcal{G})$ presents the homotopy colimit of the diagram in $\mathpzc{C}$
\end{enumerate}
Denote by $\mathbb{L}F$ the associated derived functor
$$\overline{\mathbf{F}}:\mathrm{L}^{H}(\mathrm{s}\mathcal{P}_{\Sigma}(\mathcal{G}))\rightarrow\mathrm{L}^{H}(\mathrm{sC})$$
Then $\overline{\mathbf{F}}$ commutes with sifted colimits.
\end{lem}

\begin{proof}
The assumption on filtered colimits in $\mathpzc{C}$ implies that $\overline{\mathbf{F}}$ commutes with $(\infty,1)$-categorical filtered colimits, and the assumption on realisation of simplicial objects in $\mathpzc{C}$ implies that $\overline{\mathbf{F}}$ commutes with $(\infty,1)$-categorical geometric realisations. Thus $\overline{\mathbf{F}}$ commutes with all sifted colimits.
\end{proof}

\subsubsection{Example: Rigidification of Functors Between Categories of Simplicial Objects}

Let $\mathrm{C}$ be a locally presentable $1$-category and $\mathcal{Z}=\{Z_{i}\}_{i\in\mathcal{I}}$ a set of tiny objects of $\mathrm{C}$ closed under finite coproducts. Consider the category $\mathrm{sC}$ equipped with the model structure induced by $\mathcal{Z}$. Let $\mathcal{G}$ be a $1$-category and $F:\mathcal{G}\rightarrow\mathrm{C}$ a functor. Consider the induced functor
$$F:\mathrm{s}\mathcal{P}_{\Sigma}(\mathcal{G})\rightarrow\mathrm{sC}$$
Suppose that $F$ preserves weak equivalences between cofibrant objects. Now since $\mathcal{Z}$ consists of tiny objects, filtered colimits of weak equivalences are weak equivalences in $\mathrm{sC}$. Thus filtered colimits always present homotopy filtered colimits. Suppose further that the diagonal of a bisimplicial object $\Delta^{op}\rightarrow\mathrm{sC}$ is equivalent to its geometric realisation which happens, for example, if for each $Z_{i}\in\mathcal{Z}$, $\mathrm{Hom}(Z_{i},-)$ commutes with reflexive coequalisers. Then the conditions of Lemma \ref{lem:rigidfunc} hold. 

\section{Adjunctions of Stable $(\infty,1)$-Categories}

We have the following result regarding adjunctions between stable $(\infty,1)$-categories with $t$-structures and enough projectives which will be useful later for comparing modules in different contexts.

\begin{lem}
Let
$$\adj{L}{\mathbf{C}}{\mathbf{D}}{R}$$
be an adjunction of stable $(\infty,1)$-categories with $t$-structures, with $L$ being right exact (and $R$ being left exact) for the $t$-structures. Let $\overline{\mathbf{C}}\subset\mathbf{C}^{\heart}$ and $\overline{\mathbf{D}}\subset\mathbf{D}^{\heart}$ be full subcategories closed under finite limits and colimits such that 
\begin{enumerate}
\item
the adjunction 
$$\adj{L^{\heart}}{\mathbf{C}^{\heart}}{\mathbf{D}^{\heart}}{R^{\heart}}$$
restricts to an equivalence of categories
$$\adj{L^{\heart}}{\overline{\mathbf{C}}}{\overline{\mathbf{D}}}{R^{\heart}}$$
\item
$R(d)\in\mathbf{C}^{\heart}$ for any $d\in\overline{\mathbf{D}}$
\item
$L(c)\in\mathbf{D}^{\heart}$ for any $c\in\overline{\mathbf{C}}$
\end{enumerate}
 Let $\mathbf{C}^{\overline{\mathbf{C}}}$ (resp. $\mathbf{D}^{\overline{\mathbf{D}}}$) denote the full subcategory of $\mathbf{C}$ (resp. $\mathbf{D}$) consisting of objects $F$ (resp. $G$) such that $\pi_{n}(F)\in\overline{\mathbf{C}}$ (resp. $\pi_{n}(G)\in\overline{\mathbf{D}}$) for all $n\in\mathbb{Z}$. Then there is an induced equivalence
$$\adj{L}{\mathbf{C}^{\overline{C}}\cap\mathbf{C}^{b}}{\mathbf{D}^{\overline{D}}\cap\mathbf{D}^{b}}{R}$$
Suppose further that $\mathbf{C}$ and $\mathbf{D}$ are left complete (resp. right complete). If Mittag-Leffler sequences have no higher derived limits in $\mathbf{C}^{\heart}$  (co-Mittag-Leffler sequences have no higher derived colimits in $\mathbf{D}$) then there is also an induced equivalence
$$\adj{L}{\mathbf{C}^{\overline{C}}\cap\mathbf{C}_{+}}{\mathbf{D}^{\overline{D}}\cap\mathbf{D}_{+}}{R}$$
(resp.
$$\adj{L}{\mathbf{C}^{\overline{C}}\cap\mathbf{C}_{-}}{\mathbf{D}^{\overline{D}}\cap\mathbf{D}_{-}}{R}$$).
\end{lem}

\begin{proof}
Firstly let us show that indeed there is a well-defined adjunction
$$\adj{L}{\mathbf{C}^{\overline{C}}\cap\mathbf{C}^{b}}{\mathbf{D}^{\overline{D}}\cap\mathbf{D}^{b}}{R}$$
Let $d\in\mathbf{D}^{\overline{\mathbf{D}}}\cap\mathbf{D}_{\ge0}\cap\mathbf{D}_{\le n}$. There is a cofibre-fibre sequence
$$\pi_{n}(d)[n]\rightarrow d\rightarrow\tau_{\ge n-1}d$$
Applying $R$ we get a fibre sequence
$$R(\pi_{n}(d))[n]\rightarrow R(d)\rightarrow R(\tau_{\ge n-1}d)$$
A simple induction shows that $R(d)\in\mathbf{C}^{\overline{D}}\cap\mathbf{C}_{\ge0}\cap\mathbf{C}_{\le n}$ and that $\pi_{n}(R(d))\cong R(\pi_{n}(d))$
A similar argument shows that for any bounded object $c$
$\mathbf{C}^{\overline{\mathbf{C}}}$ we have 
$$L(\pi_{n}(c))\cong \pi_{n}(L(c))$$
Note also (see for example \cite{perverset} page 10) since $L$ is right $t$-exact we have for any $X\in\mathbf{C}$,
$$\tau_{\ge n}L(\tau_{\ge n}X)\cong\tau_{\ge n}L(X)$$
Thus we have for any $c\in\mathbf{C}^{\overline{C}}\cap\mathbf{C}_{+}$
$$L(\pi_{n}(c))\cong \pi_{n}(L(c))$$

Now we have

$$R(d)\cong \underset{n}{\Rlim}R(\tau_{\ge -n}d)$$

There is a fibre sequence

$$R(d)\rightarrow\prod_{n}R(\tau_{\ge -n}d)\rightarrow\prod_{n}R(\tau_{\ge -n}d)$$
Now we have $\pi_{k}(R(\tau_{\ge -n}d))\cong R(\pi_{k}(\tau_{\ge -n}d))$ which is $R\pi_{k}(d)$ for $k\le -n$ and $0$ otherwise. 

$$\pi_{p}(R(d))\cong\lim \pi_{p}R(\tau_{\ge -n}d)\cong R(\pi_{p}(d))\cong\mathrm{R}\lim R\pi_{p}(\tau_{\ge -n}(d))\cong R(\pi_{p}(d))$$

Each sequence of abelian groups
$$\cdots\rightarrow \pi_{p}(\tau_{\ge -n-1}R(d))\rightarrow \pi_{p}(\tau_{\ge -n}R(d))\rightarrow\cdots$$
is Mittag-Leffler, so the sequence
$$0\rightarrow\lim_{n}\pi_{k}(R(\tau_{\ge -n}d))\rightarrow\prod_{n}\pi_{k}(R(\tau_{\ge -n}d))\rightarrow\prod_{n}\pi_{k}(R(\tau_{\ge -n}d))\rightarrow0$$
is exact for each $k$.

Going back to the fibre sequence 
$$R(d)\rightarrow\prod_{n}R(\tau_{\ge -n}d)\rightarrow\prod_{n}R(\tau_{\ge -n}d)$$
gives that the sequences 
$$0\rightarrow\pi_{k}R(d)\rightarrow\prod_{n}\pi_{k}(R(\tau_{\ge -n}d))\rightarrow\prod_{n}\pi_{k}(R(\tau_{\ge -n}d))\rightarrow0$$
are exact, so that $\pi_{k}(R(d))\cong R(\pi_{k}(d))$

So we see that the maps
$$c\rightarrow RL(c)$$
$$LR(d)\rightarrow d$$
are homology isomorphisms, which proves the claim. The claim for bounded above sequences works similarly, as do both of the unbounded claims. 
\end{proof}

%
%
%

\end{document}